\documentclass{cmfne}
\usepackage[cp1251]{inputenc}
\usepackage{cite}
\usepackage[mathscr]{eucal}

\usepackage{amsfonts, amsbsy, amssymb,amsmath,graphicx}
\usepackage{graphicx,float}
\usepackage{epigraph}
\usepackage{caption}
\usepackage{subcaption}
\usepackage[all]{xy}

\def\currentyear{2018}

\firstpage{1} 
 % Номер УДК
\title{Transmutation Method and Boundary-Value Problems\\for Singular Elliptic Equations} % Название статьи или обзора
\author[Valeriy V. Katrakhov (1949--2010)]{V.\,V.~Katrakhov} % Инициалы и фамилия автора
\address{Voronezh and Vladivostok, Russia.\\~}
\author[Sergey M. Sitnik]{S.\,M.~Sitnik}
\address{Belgorod National Research University ``Belgorod State University,''\\
Institute of Engineering Technology and Natural Science, Belgorod, Russia.} % Адрес автора
\email{sitnik@bsu.edu.ru} % Электронный адрес автора

\theoremstyle{plain} % Далее вводятся окружения типа "теорема"
\newtheorem{theorem}{Theorem}[section]
\newtheorem*{theorem*}{Theorem}
\newtheorem{lemma}{Lemma}[section]
\newtheorem*{lemma*}{Lemma}
\newtheorem{corollary}{Corollary}[section]
\newtheorem*{corollary*}{Corollary}
\newtheorem{property}{Property}[section]
\theoremstyle{definition} % Далее вводятся окружения типа "определение"
\newtheorem{definition}{Definition}[section]
\newtheorem{remark}{Remark}[section]
\numberwithin{equation}{section}
\numberwithin{section}{chapter}

\renewcommand{\sectionmark}[1]{\markright{\sc #1}}
\newcommand{\sectionmarknum}[1]{\markright{\sc\thesection. #1}}

\renewcommand{\geqslant}{\ge}

\makeatletter
\renewcommand*\l@subsection{\@dottedtocline{2}{2em}{3.3em}}
\makeatother

\def\p{{,}}
\def\Re{\mathop{\rm Re}\nolimits}
\def\supp{\mathop{\rm supp}\nolimits}
\def\arg{\mathop{\rm arg}\nolimits}
\def\const{\mathop{\rm const}\nolimits}
\def\ds{\displaystyle}

\newcommand{\lr}[1]{\left(#1\right)}
\newcommand{\sq}{\sqrt}
\newcommand{\ov}{\overline}
\newcommand{\pr}{\partial}
\newcommand{\pd}{\partial}
\newcommand{\lrs}[1]{\left[#1\right]}
\newcommand{\N}{\mathbb{N}}

\newcommand{\R}{\mathbb{R}}

\begin{document}

\begin{abstract}
The main content of this text is composed from the two doctoral
dissertations of V.\,V.~Katrakhov (1989) and  S.\,M.~Sitnik
(2016). In our work, we systematically expound the theory of
transmutation operators and their applications to differential
equations with singularities in coefficients (in particular, with
Bessel operators);  for the first time, this is done in the
monograph format. Apart from a  detailed review and bibliography
on this theory, the book contains original results of the authors.
Significant part of these results is published with detailed
proofs for the first time. In the first chapter, we provide
    the historical background, necessary notation, definitions, and
auxiliary facts. In the second chapter, we provide the detailed
theory of Sonin and Poisson transmutations. In the third chapter,
we describe an important special class of the Buschman--Erd\'elyi
transmutations and their applications. In the fourth chapter, we
consider new weighted boundary-value
 problems with Sonin and
Poisson transmutations. In the fifth chapter, we consider
applications of a special form of the Buschman--Erd\'elyi
 transmutations to new boundary-value
  problems for elliptic equations with
    essential singularities of solutions. In the sixth chapter, we
describe a universal compositional method to construct
transmutations and provide its applications. In the seventh
chapter, we consider applications of the transmutation theory to
differential equations with variable coefficients: namely, to the
problem to construct  a new class of transmutations with sharp
estimates of kernels for perturbed differential equations with the
Bessel operator as well as to special cases of the well-known
Landis problem on exponential estimates of the growth rate  for
solutions of the stationary Schr\"odinger equation. The book is
concluded with  brief biographic essay of Valeriy Vyacheslavovich Katrakhov
 and a detailed bibliography containing~\ref{BibEnd}
references.
\end{abstract}

\maketitle

\vspace{-3ex}
\tableofcontents

\chapter*{Preface of Professor A.B. Muravnik}
\addcontentsline{toc}{chapter}{Preface of Professor A.B. Muravnik}

The theory of singular differential equations containing the
Bessel operator
$$\displaystyle{B_\nu = {1\over x^\nu}{d\over d x}
\left(x^\nu{d\over d x}\right)}$$
    and the theory of the corresponding function spaces
    (actually, these two research areas are     inextricably intertwined)
    are among mathematical directions such that their theoretical
    and applied value barely can be overestimated.
    As other actively developed directions of the forefront of the contemporary
    science, this area entirely confirms the known     aphorism
``founders do not compose text-books.''
    Indeed, during the last four decades of the previous century
 (actually, this theory in its contemporary form is created within the said period),
 a lot of papers with absolute breakthrough results are published,
 a  decent group of     highly skilled researchers is     nurtured
  (it suffices to note that at least five D.Sci. dissertations
  are successfully defended), but only one monograph is published;
  this is the fundamental work  \emph{Singular Elliptic Boundary-Value
 Problems}  (M.: Nauka, 1997) of the founder of this area  Ivan Aleksandrovich
Kipriyanov.

For the stage of the     explosive development of a new research
area, this is absolutely typical.
 However, during the initial 15 years of the current century, this science
 sustained very hard losses:  the founder I.\,A.~Kipriyanov, the
 classics
V.\,V.~Katrakhov and V.\,Z.~Meshkov, and other mathematicians, who
developed this direction persistently and finely, passed away. New
stage has come.
 Now, it is extremely important to sort out teachers' archives, look
 after their scientific heritage, and provide an access to it for
 the wide mathematical audience
(definitely, together with the further development of the
investigations).
 A decent example of such an activity is the monograph proposed to readers' attention,
 which is  the result of the herculean task (there is no exaggeration here) of S.\,M.~Sitnik
 to  systematize the heritage of his mentor
 V.\,V.~Katrakhov.

This work was begun by Katrakhov himself (more than twenty  years ago),
and, eventually, it is completed now. The very     cautious and
    meticulous work of the second author to preserve the heritage of the first one
    deserves the highest mark (frequently, even the     inimitable Katrakhov explanation
    manner, which was so much attractive for his readers and students, is preserved).
    However, the value of this work is beyond the pale of a simple sorting out of
     teacher's archive: apart from Katrakhov results (including the ones that were not
     published during his life), Sitnik included results obtained in their co-authorship
     and his own results obtained during those almost ten years that we live without
     Katrakhov.

The specified results form a consistent and reasonably structured
book devoted to one of the most efficient and (still)
 promising methods to study singular problems of the Kipriyanov direction;
 the method of transmutation (intertwining) operators is meant.
 Undoubtedly, the applications of this method are not restricted
 by singular differential equations.
 However, as it is shown during the last three decades (primarily, by the authors
 of the monograph presented to readers' attention),
   for equations with Bessel operators,  this method yields more than enough to treat
   transmutation operators as a natural part of the singular
   theory, comparable, e.\,g., with weight Kipriyanov function
   spaces.

Undoubtedly, the release  of this book will become an event of the
mathematical life.

\vspace{5mm}

\begin{flushright}
{\it A.\,B.~Muravnik}
\end{flushright}

\newpage

\chapter*{Author's  Preface}
\addcontentsline{toc}{chapter}{Authors' Preface}

 The theory of transmutation operators is a well-developed
 independent area of mathematics.
 A significant contribution to this theory and its applications to partial
 differential equations is provided by works of the Voronezh mathematician
 Valery Vyacheslavovich Katrakhov  (1949--2010),
 who is a disciple of  Ivan Aleksandrovich
Kipriyanov.

    Among important results of V.\,V.~Katrakhov is the
    investigation (by means of transmutation operator technique)
    of weight and spectral problems for differential
    equations and systems with Bessel operators.
  Together with I.\,A.~Kipriyanov, he introduced and studied
    equations with pseudodifferential operators defined via the
    Hankel transformation by means of the Sonin and Poisson transmutation operators.

The new class of boundary-value
 problems for the Poisson equation, introduced by  Katrakhov, deserves a special
 attention: they admit  solutions with essential singularities.
 He introduced a new class of transmutation operators obtained by compositions of the known
 Sonin and Poisson operators with fractional Riemann--Liouville
 integrals. Basing on this new class,  Katrakhov introduced
 special function spaces containing functions with essential
 singularities and proved embedding theorems for them as well as
 direct and inverse trace theorems.
 For functions without singularities, the specified spaces are reduced to Sobolev spaces.
 Thus, they are their direct generalizations.
 To ensure the well-posedness
  of the problems with essential singularities,  Katrakhov
  proposed a new natural boundary-value
  condition at an inner point of the domain: the limit of the
  convolution of the solution with a smoothing kernel of the
  Poisson kind is set.
  It is proposed to call this  new boundary-value
  condition the  $K$-\emph{trace}
  (in the honor of Katrakhov introduced this condition and studied  boundary-value
  problems with it in detail).
  Solutions of the Laplace equation, possessing singular points (including essential
  singularities) are completely characterized in terms of  $K$-traces.
   Katrakhov proved the well-posedness
   of the said problem (including the unique solvability and a priori estimates for solutions)
    in the specified function spaces.
    This result generalizes the theorems on solvability  in Sobolev spaces
     for elliptic equations
    with smooth solutions without singularities.
   In further works of  Katrakhov and his
   co-authors,   new boundary-value
   problems are generalized for equations with Bessel operators and singular potentials,
   for domains in Lobachevsky spaces, and for the case of angular
   points on the boundary of the domain.
   Main results of  V.\,V.~Katrakhov are briefly listed in \cite{S92,
S95} as well.

The present monograph consists of results from D.Sci.
dissertations of   V.\,V.~Katrakhov  (1989) and S.\,M.~Sitnik
(2016) (see, respectively, \cite{KatDis} and \cite{SitDis}).
    Results of the second author (who is the disciple) develop results of the first author
     (who is the mentor).
     Hopefully, the present book will promote a broader  knownness of Katrakhov's
     results, which are significantly interesting for the theory
     of degenerating and singular differential equations as well,
     as their development in the framework of ideas and methods of the theory of transmutation
     operators.
     Also, this book reflects the contribution into the development of the theory of
     differential equations and function theory, made by
      Ivan Aleksandrovich
Kipriyanov and the Voronezh mathematical school on singular and
degenerating differential equations, created by him.

\vspace{5mm}

\begin{flushright}
\it  S.\,M.~Sitnik. Belgorod--Voronezh,
 2018--2022.
\end{flushright}

\newpage

\chapter{Introduction}\label{ch1}
\section{Historical Data and Brief Book Content}\label{sec1}

Within several last decades, the interest to singular and
degenerate elliptic boundary-value
 problems goes up; this is caused by application needs.
  This direction is originated by the fundamental Keldysh work \cite{Kel}
  finding (on examples of second-order
  equations with power-like
  coefficients) main specific properties of the posing of boundary-value
   conditions for such equations.
  It is shown that that there relations between the coefficients on the characteristic
  part of the boundary such that the Dirichlet condition is to be posed
   (the $D$ problem), but it is to be changed for the boundedness condition for the solution
    (the $E$ problem) under other other relations. In the latter case,
    we say about the \emph{strong degeneration} (in the Nikol'skii--Lizorkin
    sense) of the corresponding problem. Analogs of the
  $D$ and  $E$ boundary-value
 problems is studied (by various authors) for quite general elliptic equations
  (see \cite{67,57, 58, 59, 68,48,39, 40, 41,15,72,9,Bitz2,76,Trib1} and references therein).
    If the degeneration is strong, then, apart from bounded solutions, the equation
    has solutions unbounded near the characteristic part of the boundary (i.\,e.,
    singular ones).
    In \cite{9, Bitz2}, Bitsadze proposes not to set the solution
    itself or its normal derivative on the characteristic part;
    instead, their products with weight functions selected in
    advance are set.
Such problems are called \emph{weight} ones. The weight Cauchy
problem for hyperbolic equations is studied in {Lio1,CSh} (see
references therein as well).  The weight Cauchy problem for
elliptic equations is studied in \cite{12,87,Yanu} and many other
works (see  \cite{SKM} and references therein as well).

On the other hand, strongly degenerating boundary-value
 problems arise in the theory of singular points of solutions (including regular ones)
 of elliptic equations.
    In this theory, classical results are obtained in \cite{14,37}.
    Newer approaches and results are provided in  \cite{54,53,85,69,18}.
    In all the cited papers, the main concern is to find conditions providing
    the singularities to be removable.
    For equations of mathematical physics, the corresponding facts are provided \cite{79}.
  The technique of degenerating equations is not used here because the problem posing leads
  to the elimination of the solution singularity, while it is impossible to set weight
  boundary-value problems without additional restrictions for the growth of the solution.

Again, we emphasize that the investigation of degenerating and
singular partial differential equations with variable coefficients
 is originated by \cite{Kel} (the $E$ and $D$ problems).
    We treat singular and degenerating equations in the sense of the classical monograph
    \cite{CSh}.
  Singular, degenerating, and  mixed-type
   equations (closely related to them) can be joined in the class of
    \textit{nonclassical equations of mathematical physics}; this term is proposed by Vragov
    (see \cite{Vrag}).
    The theory of equations of the specified types is developed
    by numerous mathematicians including Gellersted, Protter,
Radulescu,  Tricomi, Fichera, Holmgren, Cibrario, Aldashev,
A.~Andreev, Baranovskii, Bitsadze, Bubnov, Vasharin, Vekua,
Vishik, Volkodavov, Vragov, V.~Glushko, Jaiani, I.~Egorov,
 Zhegalov, Zarubin, V.\,A.~Il'in, ,A.~Il'in, Kal'menov, Kapilevich,
Kilbas, Klimentov, Kozhanov, Kudryavtsev, Lizorkin,
 Marichev, L.~Mikhailov, Mihlin,
E.~Moiseev, Nakhushev, N.~Nikolaev, S.~Nikol'skii, Oleinik,
Parasjuk, S.\,V.~Popov, Pul'kin, Pul'kina, Pyatkov, Rasulov,
Repin, Sabitov, Salakhitdinov, Skubachevskii, M.\,M.~Smirnov,
Soldatov, Rutkauskas, S.~Tersenov, Usmanov, V.\,Evg.~Fedorov,
V.\,Evs.~Fedorov, Frankl', Chibrikova, and Yanushauskas.
  Many papers are published in this area; we cite only
  \cite{Bitz2,EgFed,EPP,Glushko,Kip1,MaKiRe,Moi,Nah4,Pul2,Pyat,Rep,Smi,Sku1,Sku2,Sku3,Sku4,Ter1,
 Tricomi1,Yanu,CSh,Rad1}.

    Functional-differential equations form a special class (differential-difference
 equations or equations with retarded or deviated independent variables are their particular
 cases). Their theory is originated in works of Myshkis, Hale,
 G.~Kamenskii, and El'sgol'ts (see \cite{Mys,Hale,Els,ElNor,GKE,ZKNE}).
 Substantial results in this direction are obtained by Azbelev,
 Skubachevskii, V.~Maksimov, L.~Rakhmatullina,
Rossovskii, and Muravnik; we cite only
\cite{AMR,Sku1,Sku2,KaSk,Skub,Ross1,Ros3,VaRo,Mur}.
   V.\,A.~Rvachev and V.\,L.~Rvachev use functional-differential
 equations to define a special class of atomic functions applied in the approximation theory
 (see \cite{Rv1, Rv2, Rv3}) and numerous applied problems  (see \cite{Kra1, Kra2}).
    It is important that even prototype problems for  functional-differential
 equations might possess compactly supported solutions (such as Rvachev
  atomic functions); this phenomenon is possible nor for differential neither for difference
  equations. Also,  functional-differential
 equations include the important class of equations with Dunkl
 operators related to Bessel differential operators and appearing
    at the interface of the group theory and theory of group symmetries (Coxeter groups),
    differential equations, integral transformations,
    quantum physics, and crystallography (see \cite{Me1,DHS, Dun1, Dun2,
Dun3, Gal1,Gal2, Ros1,Ros2, Tri6}).  Also, the class of
 functional-differential  equations includes problems with  involute (Carleman) translations
  (see \cite{KaSa,Lit1}).

We select a special class of singular partial differential
equations such that its typical representative is the following
$B$-elliptic equation with Bessel operators with respect to each
variable:
\begin{equation}
\label{Bes2} \sum\limits_{k=1}^{n}B_{\nu,x_k}u(x_1,\dots, x_n)=f.
\end{equation}
  $B$-hyperbolic and $B$-parabolic
  equations are considered in the same way.
  This convenient terminology is introduced in \cite{Kip1}.
  The study of this class of equations is originated in works of Euler, Poisson, and Darboux.
 Then it is continued in the Weinstein generalized axially symmetric
potential theory (see \cite{Wei1, Wei2, Wei3}), in
 \cite{Bers1, Bers2, Bers3}, and in works of Aldashev, I.~Egorov, Zhitomirskii,
 Kapilevich,
Kilbas,  Kudryavtsev, Lizorkin, Sh.~Karimov, E.~Karimov, Marichev,
Matiychuk, L.~Mikhailov, Olevskii, Pul'kin,  M.\,M.~Smirnov,
S.~Tersenov, Khasanov, Kan Cher Khe, Yanushauskas, and others.
    The importance of equations from these classes is caused by their usage in applications
    as well: they are applied to problems of the theory of axially symmetric potential
    (see \cite{Wei1, Wei2, Wei3}),
    Euler--Poisson--Darboux
    equations (see \cite{VoNi,Dza1}, the Radon transformation and tomography (see
    \cite{Nat1, Rub3, Lud, Rub1, Rub2, Rub4}), gas dynamics and
    acoustics (see \cite{Bers1, Bers2, Bers3}), hydrodynamical stream theory
     (see \cite{Gur}), linearized Maxwell--Einstein
  equations (see \cite{Bitz1,Bitz12})), mechanics, elasticity and plasticity theory
  (see \cite{Dza2}), and many others.

One can say that the mentioned three classes (in the Kipriyanov
terminology) are considered in the following three monographs:
\cite{Kip1} (for
    $B$-elliptic equations), \cite{CSh}  (for
    $B$-hyperbolic ones), and \cite{Mat1} (for $B$-pa\-rabo\-lic
    ones).

The most comprehensive investigation of equations with Bessel
operators is performed by  Kipriyanov and his disciples
L.\,A.~Ivanov, Katrakhov, A.~Ryzhkov, Aziev, Arkhipov, Baydakov,
Bogachev, Brodskii, Vinogradova, V.~Zaytsev, Zasorin, G.~Kagan,
Katrakhova, Kipriyanova, Kononenko, Klyuchantsev, A.~Kulikov,
Larin, Leizin, Lyapin, Lyakhov, Muravnik, I.~Polovinkin,
A.~Sazonov, Sitnik, Ukrainskii, Shatskii, E.~Shishkina, and
Yaroslavtseva; main results of this direction are presented in
\cite{Kip1}.

\begin{figure}[H]
\centering
\begin{subfigure}[t]{0.45\textwidth}
\centering
\includegraphics[width=0.7\textwidth]{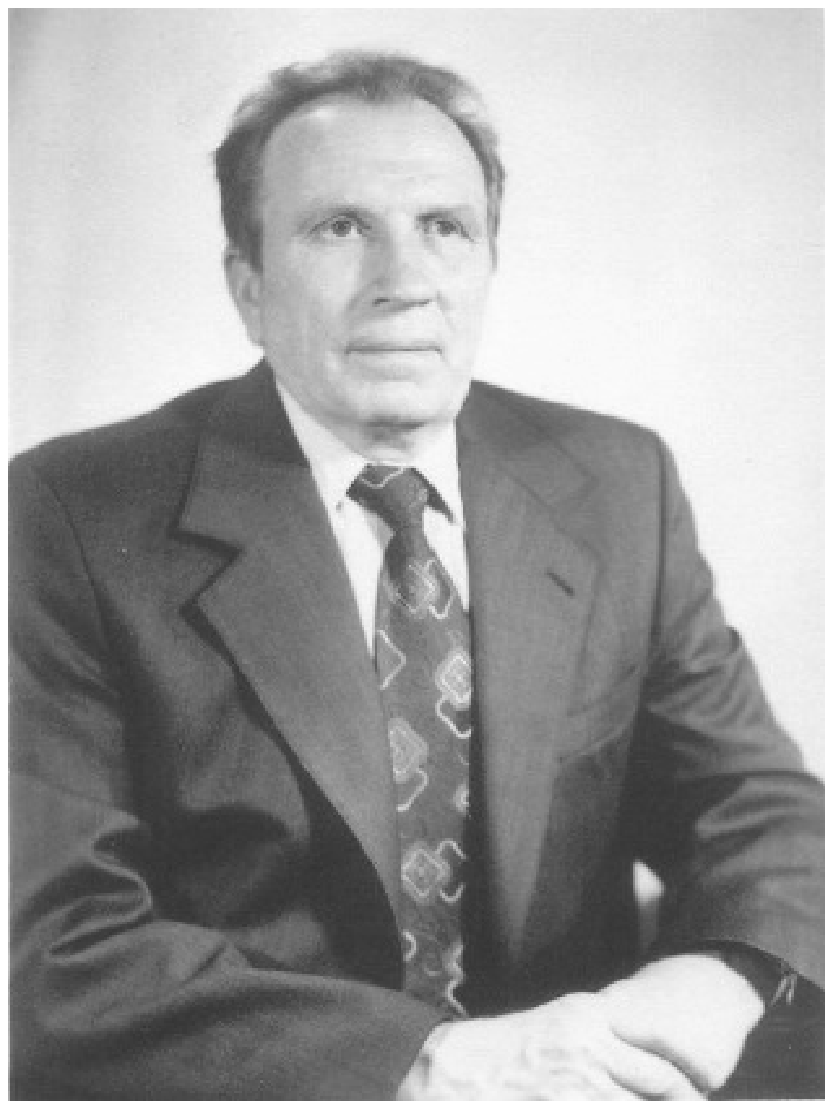}
\subcaption{\hfil Kipriyanov Ivan Alexandrovich\hfil\\
\hphantom{\hspace{2.9cm}}(1923--2001)} \label{pic1}
\end{subfigure}
\qquad
\begin{subfigure}[t]{0.45\textwidth}
\centering
\includegraphics[height=0.936\textwidth]{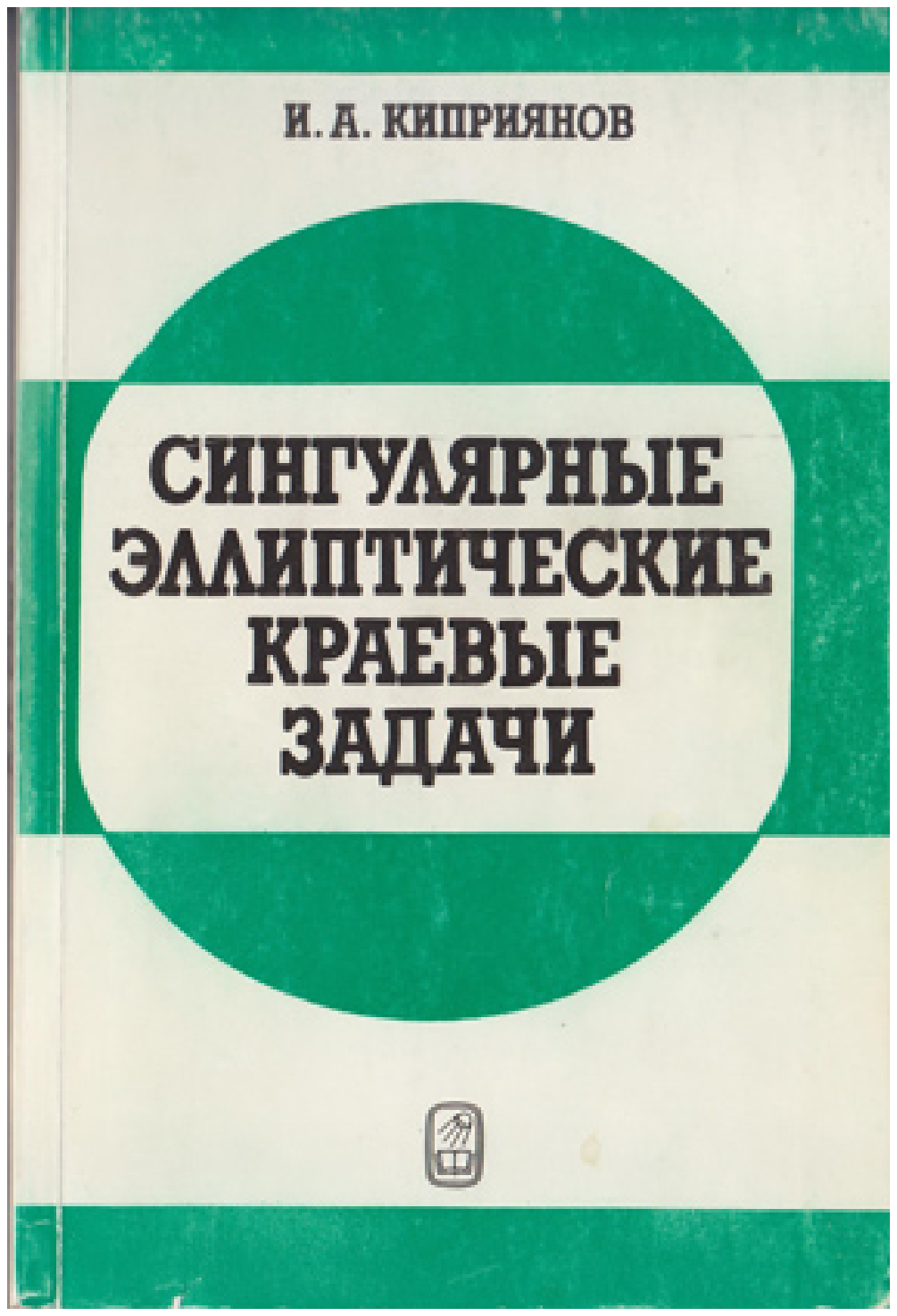}
\subcaption{\hfil Kipriyanov's book:\\ Singular elliptic boundary-value problems, 1997\hfil}\label{pic2}
\end{subfigure}
\end{figure}

To describe classes of solutions of the corresponding equations,
Kipriyanov introduced and investigated function spaces called in
his honor afterwards (see monographs
 \cite{Trib1,KN} containing special sections devoted to Kipriyanov  spaces).
 In \cite{KiKa1}, pseudodifferential operators are studied by means of
 So\-nin--Poisson--Delsarte
 transmutation operators. In a processed form,  these results are presented at a separate
 chapter of \cite{Car2}.
 In \cite{Iva1,Iva2,KipIv1}, the investigation of factorable  $B$-hyperbolic
 equations is originated.
 In \cite{42, KipIv2, KipIv3, KipIv4},  $B$-hyperbolic
 equations in Riemann and Lobachevskii spaces are investigated.
 In \cite{KipKo1, KipKo2, KipKo3, KipKu1}, fundamental solutions for  $B$-elliptic
 equations are constructed and studied.
    In \cite{Lyah1, Lyah2, Lyah3}, weight spherical harmonics,
    generalized potentials and  hypersingular  integrals, problems
    for pseudodifferential equations with Bessel operators, and
    generalizations of the Radon transformation are investigated.
In \cite{Baid}, a priori estimates of H\"older norms of solutions
of quasilinear $B$-elliptic
 equations are considered.
 In \cite{Zai1,Zai2}, factorable equations with singularities at
 coefficients are studied and an original method to investigate
 the Huygens  principle for them is proposed.
 This method allows one to reconsider the notions  ``weak lacuna,'' ``Huygens  principle,''
 and ``diffusion of waves,'' introducing the notion ``Huygens  principle of order $(p,q)$,''
 and to construct examples of equations, where the Huygens  principle is satisfied in
 even-dimensional spaces (earlier, this was treated impossible).
    In \cite{Sha1, Sha2, Sha3}, first-order
    singular and degenerating systems are considered.
 In \cite{Lar1} and \cite{Lar2, Lar3, Lar4, Lar5, Lar6}, various problems
 for the Fourier--Bessel
 transformation and singular differential equations (respectively) are considered.
 In \cite{Mur3, Mur4, Mur5, Mur6, Mur7}, properties of  the Fouri\-er--Bes\-sel
 transformation are considered. In \cite{Mur2,Mur8,Mur9}, singular differential equations
 and stabilization properties of their solutions are considered.
    In \cite{S1},
    bo\-undary-va\-lue problems with $K$-trace
    are considered for the stationary Schr\"odinger equation
    with singular potentials.
    In \cite{S46,14,SitDis,S42,FJSS}, the general composition method to construct
    various classes of transmutation operators is developed.
    In \cite{S3,S4,S66,S6,S46,S103,S400,SitDis,S42,S38,S401,S402}
    and \cite{S125,S46,S14,S42,S46}, new classes of transmutation operators    of
 Buschman--Erd\'elyi and
  Ve\-kua--Er\-d\'e\-lyi--Lowndes
  types (respectively) are introduced.
  Using the obtained norm estimates for Buschman--Erd\'elyi
  operators, we prove the equivalence of norms in Kipriyanov spaces and weight Sobolev
  spaces, generalize the Hardy inequality, and apply the method of transmutation operators
  to investigate special cases of the known Landis problem on the limit decrease rate
  of solutions of the stationary Schr\"odinger equation (see \cite{S75,S71,S3,S42}).
    Also, we modify the known method to construct transmutation
    operators for perturbed Bessel equations (see \cite{S63,S19,S43,S42,SitDis}),
 construct implicit integral forms of fractional powers of the
 Bessel and investigate their main properties
  (see \cite{S135,S133,S127,S123,S18,S700,SS,S42}),
  and consider applications of the method of transmutation operators to
    Euler--Poisson--Darboux
    equations  (see \cite{SSfiz,ShiE3}).

In works of Kipriyanov and mathematician of his school, the most
profound results for an important class of singular differential
equations with parameter are obtained.
 As we note above, in \cite{Kip1, 39}, weight function spaces are introduced by means
 of the closing of the set of functions such that they are smooth and even with respect
  to the
 direction normal to the boundary.
    Since not all traces are preserved after closing (this is a specific property
    of the theory of weight spaces), it follows that the boundary-value
    problem considered in these spaces (see \cite{39, 40, 41}) is similar to the problem
    of type $E$ with the homogeneous evennes  boundary-value
    conditions.
    Moreover, the mentioned parameter is real because the main method is the
    Fourier--Bessel transformation.
    We use another approach; it allows one to omit these restrictions
    and to consider the case of an arbitrary complex parameter;
    formally, this is a generalization neither from the viewpoint of
    boundary-value problems formulation nor from the viewpoint of
    physics applications.
    Unlike \cite{40,41}, we study general heterogeneous weight
    conditions and, therefore, consider singular solutions in new
    (broader) function spaces.
    If our method is applied under restrictions mentioned above, then the obtained
    results essentially coincide with Kipriyanov ones.
    Regarding this, we cite \cite{45} considering, for the case of a real
    parameter, general
     non-weight bo\-undary-value
     conditions set on the noncharacteristic part of the boundary;
 the evenness conditions are set otherwise.

For the first time, problems for operator-differential
 (abstract differential) equations with Bessel operators are studied in \cite{CSh}.
 In \cite{Glu1,Glu2, Glu3, Glu4, Glu55, Glu5, Glu6,
Glu7, Glu8, Glu9, Glu10, Glu11, Glu12, Glu13}, the following
 research direction for  abstract differential  equations are developed:
 the resolving of
 Euler--Poisson--Darboux
 and Bessel--Struve
   equations, operator co\-sine-functi\-ons,
  and Bessel and Legendre functions.

Thus, the current situation is as follows.
 Boundary-value problems of type $E$ are studied profoundly and
 comprehensively, while weight bo\-undary-value
 problems are considered only for quite simple (and only second-order)
  equations.
  Other settings of  bo\-undary-value
 conditions for strongly degenerating equations that can be used instead of weight settings
  if the latter ones are not possible are not studied sufficiently profoundly yet.
  In our opinion, this is caused by the following two reasons.
    First, traditional methods of the theory of partial
    differential equations do not work well for the specified
    problems (in particular, the variational method is not applicable at all for equations
    that are not strongly elliptic).
    Secondly, there was no theory function spaces such that the corresponding trace theorems
    hold for them. The known Kipriyanov weight function  spaces (see cite{66,48})
  usually applied for degenerating equations possess the following property: for functions
  from these spaces, the singularity order depends on the smoothness and decreases as the
      smoothness goes up. The latter circumstance is not acceptable in the considered case.
      This motivates the settings and investigation of new bo\-undary-value
 problems for strongly degenerating equations (in the corresponding function spaces) as well
 as the development of efficient methods to resolve them.

In brief, the structure of this book is as follows.

The first chapter is an introduction. It provides the work plan, a
brief history of the development of the theory of transmutation
operators, and collects main data about special functions,
integral transformations, and function spaces.

In the second chapter, the technique of transmutation operators
for the
    one-dimensional and multi-dimensional
    cases. This technique is used in further chapters (for detail, see the next section).
 The method of transmutation operators  appear in \cite{Lio1,Lio2,Lev7}.
 A known class of  transmutation operators is called \emph{Sonin--Poisson
 operators}.
 They are broadly applied in the theory of singular hyperbolic
   Euler--Poisson--Darboux
 equations (see the cited monographs of Lions,  Carroll, and Showalter).
 In \cite{42}, a number of important results in this direction is obtained for hyperbolic
 equations. Separate elements of this method are used (within the method of spherical means)
 for the study of the Cauchy problem for the classical wave equation
(see, e.\,g., \cite{49}).
 In \cite{76}, the corresponding transmutation operators are used
 in the theory of generalized analytic functions.
    In \cite{Kly1}, transmutation operators are constructed for other kinds of singular
  equations (these results are special cases of  transmutation operators obtained for
  the general case of hyper-Bessel
  operators,
see \cite{DK2,Dim,Kir1}).
 In \cite{Yar1} they are constructed for hyperbolic analogs of
 Sonin and Poisson operators in the trigonometric form.
 Also, the method of transmutation operators is used in the spectral theory of Schr\"odinger
 equations (see  \cite{Mar2,ShSa} and references therein).

In detail, results for the Sonin--Poisson
 class of transmutation operators are presented in the second chapter.
 We call them
 \emph{So\-nin--Poisson--Delsarte
 operators}: from the historical viewpoint, this term is more exact.
 Such transmutation operators, being integral ones, belong to the class of
 Erd\'elyi--Ko\-ber fractional integrals.

In \cite{47,86},
 Erd\'elyi--Kober transmutation operators are introduced. Their development
 is sufficiently comprehensively presented in \cite{SKM}.
 The main application  area for these (and related) operators is the constructing
 of explicit representations for classical solutions of  bo\-undary-value
 problems for se\-cond-or\-der
 equations of the
  Euler--Poisson--Darboux
  type.

    The present book is devoted to the further development of the
    method of transmutation operators. First of all, we use them to construct new
 Sobolev-type function spaces. It is impossible to that,
 applying transmutation operators known earlier. That is why we introduce
 new transmutation operators different from the integral operators studied
 (on finite intervals) by Buschman, Erd\'elyi, and others. For  Erd\'e\-lyi--Ko\-ber
 operators, new results (e.\,g., the corresponding assertions about traces)
  are obtained as well.

  Also, similar gaps were preserved in the multidimensional case (regarding the
 well-posedness of bo\-undary-va\-lue
 problems, constructing of transmutation operators, and selecting
 of suitable function spaces).
 The Radon operator (see \cite{17,82} for detail) transforms the multidimensional
 Laplace operator to the one-dimensional
 operator of the double differentiating with respect to the radial variable.
 Unfortunately, it is a smoothing operator in a way.
 This leads to the following situation: ``merits'' of the direct transformation
 are eliminated by qualitative ``demerits'' of the inverse one; actually,
  this is the main difficulty of the inversion theory for the
  Radon transformation and its numerical implementations.
    For a solution of this problem, we find a factor extinguishing
    this smoothability (it is a differential operator if the dimension of the space
  is odd and is a Liouville operator if it is even); the transmutation operator
  constructed this way is applied to construct new function space.
  Note that the used technique allows one to avoid the usage of projection spaces
  (usually, they arise in the Radon transformation theory).

In the second chapter, maps of known function spaces by means of
the introduced transmutation operators are studied in the
one-dimensional case, new spaces  $H_{\nu}^s \lr{E_{+}^1}$ are
constructed, and these spaces are generalized for the
multidimensional space.
 In the one-dimensional
 case, there are results with exact constants in estimates.
 Results of such kind are used in further chapters.

The book includes only results on transmutation operators and
boundary-value problems for differential equations.
    The same method of transmutation operators is applied to other
    problems. In particular, in \cite{25, 26, 33, KiKa1},
     the algebra of singular pseudodifferential operators
  and the theory of complex powers of singular elliptic operators
  are constructed; they generalize the known results of \cite{20,2,75}
  to the singular case.

    Regarding the theory of two classes of function spaces
    introduced in the book (i.\,e., Hilbert spaces and de\-numerably-normab\-le
ones), we note that spaces of such type satisfying theorems on
weight or nonlocal traces were not studied earlier.
    Our spaces cannot be treated as weight ones because we use integrodifferential
    operators to define their norms. Kipriyanov spaces (see \cite{Kip1}) are embedded
    into the former class specified above; Sobolev spaces (see \cite{77}) are embedded
    into the latter one.

In the third chapter, Buschman--Erd\'elyi
 transmutation operators and their modifications are considered.
    They contain all classes of transmutation operators known
    earlier and they have numerous applications.
 We introduce various classes of Buschman--Erd\'elyi
 transmutation operators: operators of the first and second kind,
 operators of zero order of smoothness, and unitary
  Sonin--Katrak\-hov and Poisson--Katra\-k\-hov
  operators.
  Among applications of these classes of transmutation operators,
  we consider unitary property and norm estimates in Lebesgue spaces on semiaxes,
 applications to the Copson lemma and the setting of Cauchy problems
 with nonlocal initial-value
 conditions, relations to the classical Hardy operators, the
  Bitsadze--Pashkovskii problem for the Maxwell--Einstein
  equation, and the proof of the norm equivalence in Kipriyanov
  spaces and Sobolev spaces with weight norms for the one-dimensional
 case.

In the fourth chapter, main results on weight boundary-value
 problems are presented.
 For higher-order
 equations, weight boundary-value
 conditions are set for the first time.
 Their main difference from results of Keldysh, Nikol'skii,
 Bitsadze, Kudryavtsev, Lizorkin, Kipriyanov,
    Vishik, Grushin, and other authors studied the $E$ problem and its analogs
  is the form of boundary-value
  conditions and their amount: we consider the number of boundary-value
  conditions equal to the half of the order of the equation.
  Using this, we are able (unlike the authors specified above) to consider
  regular and singular solutions simultaneously.
  If the coefficients of the equation have no singularities, then the weight function
  in the  boundary-value
  conditions disappears and our results coincide with classical results of the general
  theory of elliptic boundary-value
  problems (see \cite{1,83,62,51}).

In the first section, we set the weight boundary-value
  problem for homogeneous operators with constant coefficients.
  The method of transmutation operators reduces it to the boundary-value
  problem studied in the previous section.
  Note that, applying transmutation operators to equations with variable coefficients,
  we obtain pseudodifferential operators such that their symbols with respect to
  the dual variables have high-order
  singularities at the last coordinate hyperplane
  (see \cite{26,KiKa1}).
  So far, no theory of such pseudodifferential operators is developed.
  Therefore, in the case of variable coefficients, we use the Schauder technique
  to ``freeze'' coefficients and glue them afterwards, applying the partition of unity.
 According to this approach, we consider a boundary-value
 problem with slowly changing coefficients.
 To do that, we have to overcome a number of difficulties (both general ones
 and ones related to the weight character of   boundary-value
  conditions and to singularities of the coefficients).

Comparing the presented results with the corresponding results of
Bitsadze, Vasharin--Lizor\-kin,
 Yako\-v\-lev, Yanushauskas, and other results on weight
 boundary-value problems for   second-order
 elliptic equations, we easily see to see that  a lot of our results are new
 for this case as well. In particular, they are the main result about the
 Noetherian property, the form of the boundary-value
 conditions (and their order), the smoothness rising theorem, etc.
 Also, a case is found (seemingly, it was considered earlier) where
 the setting of the Dirichlet problem is impossible. Instead, a Neumann-type
 weight problem is well posed; in the physics examples provided below, it is an analog
 of the Sommerfeld radiation condition.

In the fifth chapter, a new boundary-value
 problem with singularities at isolated
 boundary points is studied for the Poisson equation.

 An arbitrary growth of the solution and the right-hand
 side is admitted at these points.
 In the case of the homogeneous (Laplace) equation,
 no a priori restrictions on the solution behavior is imposed;
 therefore, it is possible that it has an essential singularity (in the sense of
 the theory of analytic functions).
    In this general case, no known settings of boundary-value
    problems are applicable.
Indeed, in the cited works of Wiener, Keldysh, and other works of
this direction, only removable singularities are studied.
  In \cite{21}, isolated singularities of the solution are admitted, but these singular
  points are not isolated boundary points and at most power growth of the solution is
  admitted at these  points.
    Results in this direction are obtained, e.\,g., in \cite{19,35,Rut1, Rut2, Rut3}.

  Regarding this problem, we introduce a qualitatively new
    notion of the trace of a function such that its  singularities
    at isolated boundary points are admitted; in a way, such a trace is nonlocal.
    We construct new Frechet-type
    function spaces such that the corresponding trace theorems are
    valid for them.
 Here, the main result is the proof of the unique solvability of the formulated
 boundary-value problem and its Hadamard stability in the said spaces.
    Also, we classify singularities of harmonic functions in terms of the introduced
    trace. In \cite{S1, 44}, the same method is used to obtain similar results
    for other equations.

The considered boundary-value
 condition is convolutional; hence, this class of problems can be treated both as nonlocal
 ones and loaded ones (a special trace of solutions is given on a
 part of the boundary).

 In detail, let $\Omega \subset E^n$ be a bounded domain with a smooth boundary $\pr \Omega.$
    Without loss of generality, assume that the origin belongs to the domain.
    Let $\Omega_0 = \Omega \setminus 0.$ In $\Omega_0$, consider the Poisson equation
\begin{equation}
\Delta u = f (x), \  x \in \Omega_0,
\label{50}
\end{equation}
 with the following Dirichlet boundary-value
 condition on $\pr \Omega$:
\begin{equation}
 u|_{\pr \Omega} = g (x), \  x \in \pr \Omega.
\label{51}
\end{equation}
    The problem solved in the fifth chapter is as follows:
    one has to set a boundary-value
    condition at the origin and to select function spaces to ensure the unique solvability
    of the  boundary-value
    problem.
    In this case, no restrictions on the singularity order at the
    origin is imposed. This occurs in classical electrostatics
    problems. For example, let a spatial domain $\Omega_0$ be free of charges ($f=0$)
    and be surrounded by a grounded surface ($g=0$).
    Then the function $u$ is the potential of the electrostatic
    field created by a charged object placed to the  origin.
 If this object is the unit charge, then the function $u$ is called the \emph{Green function}.
 It has a singularity of the kind $|x|^{2-n}$ for $n \geq 3$ and of the kind $\ln |x|$
 for $n=2.$ If a dipole is placed to the origin, then the order of the singularity
 of the function $u$ becomes one unit greater, etc.
    In the general case, placing an infinite combination of multipoles of various orders
    to the origin, we obtain an infinite-order
 singularity of the function $u$. Moreover, the behavior of the function $u$ in a neighborhood
  of the origin is similar to the behavior of an analytic function in a neighborhood
  of an essential singularity.

    In the investigated case, known function spaces and known
    settings of boundary-value
    conditions are not applicable.
     For example, functions with essential singular points do not
     belong even to Schwartz distribution classes (see \cite{16}).
  In the spherical coordinates $r=|x|,$ $\vartheta = \dfrac{x}{|x|}$, Eq. \eqref{50}
  takes the form
\begin{equation}
\frac{\pr^2 u}{\pr r^2} +  \frac{n-1}{r} \frac{\pr u}{\pr r} +
 \frac{1}{r^2} \Delta_{\Theta} u = f \lr{r, \vartheta},
\label{52}
\end{equation}
    where $\Delta_{\Theta}$ is the Laplace--Beltrami
  operator on sphere.
  For such an equation, it is possible to set the $E$  boundary-value
    problem, but this leads to a removable singularity of the solution
     (provided that $f$ is a smooth function).
 It is easy to prove this by means of the classical technique (see \cite{79}).
 It is impossible to set weight boundary-value
 conditions: though Eq. \eqref{52} belongs to a class of known degenerating equations
 considered by Keldysh, the parameter $p$ is equal to $2$ in our case, while the inequality
  $p<2$ is a necessary condition of the possibility to set weight boundary-value
 conditions.
    Mapping the origin into infinity by means of the Kelvin transformation (see \cite{79}),
 one cannot reduce the problem to a known one because no external   boundary-value
    problems of this kind are studied earlier. Vise versa, this transformation allows one to
    transmit results obtained below to external boundary-value
    problems and no restrictions are imposed on the order of the growth of the solution.

    In the present paper, at the point $0 \in \pr \Omega_0$, the    following nonlocal
 boundary-value condition is set for the solution $u$:
\begin{equation}
\lim\limits_{r \to +0} r^{n-2} \int\limits_{\Theta} u \lr{r, \vartheta'} K_n
 \lr{r \vartheta, \vartheta'}  d \vartheta' = \Psi \lr{\vartheta'},
\label{53}
\end{equation}
    where $n \geq 3$ and $K_n (x, y)$ is the Poisson kernel associated with
    the unit circle, i.\,e.,
\begin{equation*}
K_n (x, y) = \frac{\Gamma \lr{\frac{n}{2}}}{2 \pi^{\frac{n}{2}}}
 \frac{1- |x|^2}{ |x-y|^n}, \  x,y \in E^n.
\end{equation*}
    For  $n=2$, in the polar coordinates $x_1 = r \cos \varphi$ and
     $x_2 = r\sin \varphi$, the corresponding boundary-va\-lue
      condition
 has the form
\begin{equation}
\lim\limits_{r \to +0} \frac{1}{2 \pi} \int\limits_{- \pi}^{\pi} u
\lr{r, \varphi'}  \lr{ \frac{2r \cos \lr{\varphi-\varphi'} -2
r^2}{1- 2r \cos \lr{\varphi-\varphi'} + r^2} +\frac{1}{\ln r}}  d
\varphi' = \Psi \lr{\varphi}. \label{54}
\end{equation}
 Conditions \eqref{53}-\eqref{54}
 mean that the function $u \lr{r, \vartheta}$ is averaged with a suitable kernel
 with respect to the angular variables $\vartheta$ for a fixed positive $r.$
 Once it is done, the passage to the limit as $r \to 0$ becomes possible.
 The
 left-hand side of \eqref{53} or \eqref{54} is called the \mbox{\it $\sigma$-trace}
 or the {\it $K$-trace}
  (in honor of Katrakhov introduced this condition) and denoted by
  the symbol  $\sigma u|_0$ (or $K u|_0$). The $K$-tra\-ce
  is different from zero only for functions with a singularity
   at the origin such that the order of this singularity is not exceeded
   by the singularity order of the fundamental solution of the Laplace equation.

It turns out that if a function is harmonic in a punctured
neighborhood of the origin, then its  $\sigma$-tra\-ce
 exists and it uniquely determines its singular part.
   To formulate this assertion strictly, we introduce the space
    $A\lr{\Theta}$ consisting of real functions   $\Psi$ such that they are
    analytic on the sphere $\Theta$ and the following norms are finite
    for each positive $h$:
\begin{equation}
\| \Psi \|_h = \lr{\sum\limits_{k=0}^{\infty}
\sum\limits_{l=1}^{d_k} |\Psi_{k, l}|^2 h^{-2k}}^{\frac{1}{2}},
\label{55}
\end{equation}
    where
\begin{equation*}
\Psi_{k, l} = \int\limits_{\Theta} \Psi \lr{\vartheta} Y_{k, l} \lr{\vartheta}
 d \vartheta
\end{equation*}
and $Y_{k,l}$ denotes the complete orthonormal system of spherical
 harmonics.

    Thus, the main results of the fifth chapter are as follows: we set
    boundary-value problems for the Poisson equation with isolated singular points
    such that its solution might have essential singularities at these points,
  introduce new nonlocal   boundary-value
    conditions in the form of
     $\sigma$-traces (or $K$-traces)
     at singular points,
 define the corresponding function spaces by means of transmutation operators,
 prove direct and inverse trace theorems, and prove that the considered
  boundary-value problems are well-posed
 in the introduced function spaces.

In the sixth chapter, the universal composition method for the
constructing of transmutation operators is explained and its
applications are provided.
 In this method, transmutation operators are constructed of ``blocks''
 that are classical integral
 transformations.
 This provides a unified possibility to construct all known explicit representations
 for  transmutation operators and obtain their numerous new classes.

In the seventh (final) chapter, applications of the method of
 transmutation operators to differential equations with variable
 coefficients are considered.
 First, we consider the problem to construct a new class of  transmutation operators
 for the perturbed Bessel operator, using a modified integral equation for the kernel.
 Also, we obtain exact estimates of kernels via special functions and extend
 the class of admissible potentials including Bargmann potentials,
 Yukawa potentials, and strongly singular potentials.
    For special cases of potentials, the found general estimates are refined.
    Then we consider the known Landis on the estimate of the rate of the exponential
    decay for solutions of the stationary Schr\"odinger equation (see \cite{Lan}).
    In this problem, the general answer is negative (see \cite{Mesh1,Mesh2}).
    However, for special-kind
    potentials, the Landis conjecture is confirmed.
    To obtain this result, we use the technique of transmutation operators with
  special-kind kernels.

\section[Theory of Transmutation Operators: History and Contemporary
State
 (Brief Outline)]{Theory of Transmutation Operators: History and Contemporary
State\\
 (Brief Outline)}\label{sec2}

    %   \vskip\baselineskip
    We start from the main definition.

\begin{definition}
 Let $(A,B)$ be a pair of operators.
 A nonzero operator $T$ is called a \textit{transmutation operator}
 if the relation
\begin{equation}
\label{1.1} {T\,A=B\,T}
\end{equation}
    holds.
    \end{definition}

Also, relation \eqref{1.1} is called the \textit{intertwining
property}; in this case, we say that the transmutation operator
$T$ \textit{intertwines} the operators $A$ and $B$. To convert
\eqref{1.1} into a strict definition, one has to specify function
spaces (sets) such that the operators $A,$ $B,$ and $T$ act in
them. Usually, the invertibility and continuity are included in
the definition of transmutation operators; however, these
requirements are not compulsory.
 In particular implementations, $A$ and $B$ are frequently (but  not compulsory)
 differential operators, while  $T$ is a linear integral operator; all the three ones act
 in standard spaces.

It is clear that the notion of transmutation operators is a
generalization of the linear-algebra
 notion of the similarity of matrices (see \cite{Gan, Hor, Tyr}).
 Note that no efficient methods to verify the similarity of two finite matrices exist
 because there are no efficient methods (apart from direct computations) to
 verify whether their     Jordan forms coincide each other.

 On the other hand, transmutation operators \textit{are not reduced to similar {\rm
 (}or equivalent{\rm )} operators} because, as a rule, the intertwined operators
 are unbounded in the natural spaces and, on the other hand, it is not guaranteed that
 the operator inverse to the transmutation one exists, acts in the
 same space, or is bounded.
 Thus, the spectra of  operators intertwined by a transmutation
operator are, as a rule, different from each other.
  Also, transmutation operators themselves might be bounded.
   For example, this occurs in the theory of Darboux transformations.
 This theory finds differential transmutation operators (substitution operators)
 between pairs of differential   operators; thus, in this case, all the three operators
 are unbounded in their natural spaces.
 The theory of Darboux transformations, treated as a section
 of the theory of differential equations, is coordinated with the general scheme
  of transmutation operators (in its extended version).
    If a transmutation operator is sought for a pair of operators,
    then it is not required the operators of the pair to be only differential ones.
    Problems for the following various types of operators occur in the theory
    of transmutation operators: integral operators,  integrodifferential operators,
differential-difference operators (e.\,g., Dunkl-type
 operators),
 infinite-order differential or  integrodifferential operators
(e.\,g., in relation to the Schur complementability lemma),
general linear operators in fixed function spaces, and
pseudodifferential and operator-differential
  (abstract differential) operators.

To provide an example, we briefly present the prototype scheme
illustrating the usage of transmutation operators to obtain
 binding relations between solutions of the perturbed and
 unperturbed equation provided that the intertwining property is
 proved for them.
 For example, suppose that rather complex operator $A$ is studied, while the needed properties
 are already known for a prototype (simpler) operator $B.$
 Then, if a transmutation operator of kind \eqref{1.1} exists, then it is frequently possible
 to extend properties of the prototype operator $B$ to $A.$
 This is a prototype scheme of a typical usage of transmutation operators
  in particular problems,
 explained in a few words.

In particular, if an equation  $Au=f$is considered, then, applying
the  transmutation operator $T$ with the intertwining property
\eqref{1.1} to it, we obtain an equation of the kind $Bv=g,$ where
  $v=Tu$ and $g=Tf.$ Therefore, if the equation with the operator  $B$ is simpler
  and relations for its solutions are already known, then we obtain the representations
  for solutions of the first equation as well: $u=T^{-1}v.$
  Certainly, in the framework of the explained approach, the existence of the inverse
  transmutation operator is needed and it has to act in the corresponding spaces;
  moreover, an explicit representation of this inverse operator is needed to obtain
  explicit representations of solutions.
  This is a simple example of applications of the technique of
transmutation operators in the theory of differential equations
(both ordinary and partial ones). Note that if a pair of linear
transmutation operators exists, then they provide the relation
between solutions and representations of solutions both for
homogeneous and inhomogeneous equations (as well as for equations
with a spectral parameter).

Let us clarify our terminology. In  Western publications, the term
  {\it transmutation} coming from Delsarte is conventional.
  According to Carrol, the similar term {\it
transformation} is used for the classical integral transformations
such as the Fourier transformation, Laplace transformation, Mellin
transformation, Hankel  transformation, etc.  Also, the term
 {\it transmutation} has an additional styling of a ``magical transformation,''
 rather vividly characterizing the action of transmutation operators.
 Let us cite \cite{Car3} literally:
``{\it Such operators are often called transformation operators by
the Russian school {\rm (}Levitan, Naimark, Marchenko et. al.{\rm
)}, but transformation seems too broad a term, and, since some of
the machinery seems {\rm ``}magical{\rm ''} at times, we have
followed Lions and Delsarte in using the word transmutation}.''
    In our opinion, this characterizes the case the best and the most exact.
    In Russian, the term \emph{transmutation operators} is proposed by Marchenko
    in 1940s (see \cite{Mar9}).

The necessity of the theory of transmutation operators is
confirmed by its numerous applications. Methods of transmutation
operators are especially important for the theory of differential
equations of various types.
    A lot of fundamental results for various classes of  differential
equations are proved by means of methods of transmutation
operators.

Nowadays, the theory of transmutation operators is a separate area
of mathematics, located at the intersection of differential,
integral, and integrodifferential equations, functional analysis,
function theory, complex analysis, theory of special functions and
fractional differential and integral calculus, harmonic analysis,
 theory of optimal control, inverse problems, and scattering problems.

Three main periods can be (conventionally) separated in the
development of the theory of transmutation operators.
    In the initial period, fundamental ideas and definitions are established:
    their source is the similarity theory for finite matrices (see \cite{Gan, Hor, Tyr}),
    separate results on the similarity of operators, and  separate results for basic
    differential equations.
    The idea of transmutation operators in the operator form is traditionally referred
    to Friedrichs (see \cite{Fri}).
 Actually, for the first time, the method of transmutation operators is developed and applied
    to obtain representations of solutions of differential equations by Letnikov in the XIXth
 century; this is the first real application of fractional differential and integral operators
  (treated  as transmutation operators) to problems for differential equations
  (see \cite{Shos, Koo1}).

The second period (conventionally) lasted from  1940s to 1980s; it
can be called the classical period.
 Numerous results in the theory of transmutation operators and their applications are obtained
 during this periods. The main directions and results of this period are as follows.

Methods of transmutation operators are successfully applied in the
theory of inverse problems: the generalized Fourier
transformation, spectral function, and solutions of the famous
Gel'fand--Levi\-tan
 equation are determined (see \cite{AM,
Mar1, Mar2, Mar3, Mar4, Mar5, Mar6, Mar7,Lev1, Lev2, Lev3, Lev4,
Lev5, Lev6, Lev7, Lev8}).
 In the scattering theory, the famous Marchenko equation is represented via
 transmutation operators
(see \cite{Lev1, Lev2, Lev3, Lev4, Lev5, Lev6, Lev7, Lev8,AM,
Mar1,Mar2,Fad1,Fad2}).
    For both classes of inverse problems, transmutation operators
    are the main tool because the mentioned classical equations
    are equations for kernels of transmutation operators, while,
    in inverse problems, the diagonal values of the kernels
    restore the desired potentials via the spectral function or
    the scattering data (see \cite{Laks1,Laks2, ShSa,
Kol1, Nizh1,Nizh2, New, Bloh}).
 For the Sturm--Liouville
 operators, the classical transmutation operators are constructed on the segment
  (see \cite{Povz}) and semiaxis
(see \cite{Levin2}).
 In the spectral theory, trace relations and the asymptotic behavior of the spectral function
 (see \cite{Lev1, Lev2, Lev3, Lev4, Lev5, Lev6, Lev7, Lev8}), estimates of kernels of
 transmutation operators, responsible for the stability of inverse problems
 and scattering problems (see \cite{AM, Mar1,Mar2}), and estimates of Jost solutions
 in the quantum  scattering theory (see \cite{Lev1,
Lev2, Lev3, Lev4, Lev5, Lev6, Lev7, Lev8, Sta1,Sta2, Soh1, Soh2,
Soh3, Soh4}) are obtained.
    One can say that the application of transmutation
operators to the  Sturm--Liouville
 equations with variable coefficients trivializes this equation  to the level of basic
 equations with trigonometric or exponential solutions.
 The Dirac system and other matrix systems of differential equations are investigated as well
  (see \cite{Lev6}).

    The theory of generalized analytic functions is developed;
    it can be treated as a section of the theory of transmutation
operators intertwining unperturbed and perturbed Cauchy--Riemann
 equations (see \cite{Bers1,Bers2, Berg, Vek3, Vek4, Boyar, Pol1, Pol2, Pol3})
 applied in problems of mechanics, elasticity theory, and gas dynamics.
  Methods  of transmutation
operators are used to create a new area of harmonic analysis,
studying various modifications of generalized translation
operators and generalized operator convolutions (see
\cite{Del1,Del2, Del7,Zhit,Lev2,Lev3}).
    A profound relation between transmutation operators and theorems of the Paley--Wiener
 type is found (see \cite{Sta1,Sta2,Ahi1,Che1, Che2, Che3,
Che4,Tri1,Tri2}).
 The theory of transmutation
operators allows one to create a new classification of special
functions and integral  operators with special functions in their
kernels  (see \cite{Car1, Car2, Car3,Koo1}).
 Note that the existence of Green and Riemann functions for various classes
 of differential equations as well as their explicit representations are used to find kernels
 of transmutation operators
 (see \cite{Sob, VNN, VZ}); this stimulates the search of such functions for various problems.

    In the theory of nonlinear differential equations, the Lax
    method is developed; it uses transmutation operators to prove
    the existence of solutions and construct solitons (see \cite{Zhu, ZMNP, AbSi, Car4}).
 Darboux transformations, which can be treated as transmutation operators,
  where the intertwined operators and the intertwining operator
  are differential ones, are broadly applied as well (see  \cite{MaSa});
  for the relation between the theory of the Darboux transformation and
   transmutation operators,  see the review \cite{BaSa}.
   In the quantum physics, to consider the Schr\"odinger equation
   and scattering-theory
   problems, the class of wave operators (a special class of transmutation operators)
   is studied.
    In the reviews \cite{Fad1,Fad2, Mar9}, general scattering and inverse problems
    are considered from the viewpoint of transmutation operators.
 In \cite{Kach1}, wave operators are constructed for
 scattering-theory problems with the Stark potential
 (though, in \cite{Lev1995}, the problem to construct the
 corresponding transmutation operator is stated to be unsolved).

In the transmutation-operator
 theory itself, restrictions related to the order of the differential operator are found.
 For differential operators with orders exceeding three, it is shown that classical
 Vol\-ter\-ra-kind transmutation operators exist only for the case of analytic coefficients
(see \cite{Sah1, Sah2, Sah3, Mal1, Mal2, Mal3, Mal4, Mal5});
 in the general case, transmutation operators have more complex structure and a passage
 to the complex plane is required even to construct real solutions
  (see \cite{Leo,Val,Hach1,Hach2,Mal1, Mal2, Mal3, Mal4, Mal5,Hrom1}).
 At  the same time,  in spaces of analytic functions, the equivalence of differential
 operators of equal orders is proved and a number of problems is studied
 (see \cite{FN, Fage1, Fage2, Fage3, Fage4, Fage5,Fage6, Mar5, Mar6, Mar7,Kor1, Kor2,Fish}).
 The solvability theory for the known Bianchi equation is constructed in order to
 be applied to the theory of transmutation operators (see \cite{FN}).

A separate area to apply transmutation operators is the theory
 of differential equations with singularities in coefficients.
 Primarily, this refers to equations with the Bessel operators
\begin{equation}\label{Bes1}{B_{\nu}u(x)=\frac{d^2
u}{dx^2}+\frac{2\nu+1}{x} \frac{du}{dx}.}\end{equation}
 At the initial stage, the pair of known Sonin and Poisson transmutation operators
  (they are defined in Chap. \ref{ch1} and are considered in detain in Chap. \ref{ch2})
 is applied to investigate equations of this class.
 As transmutation operators, they are introduced in \cite{Del1, Del2, Del3, Del4}.
 On the base of Delsarte ideas, their investigation is continued in
 \cite{Del5,Del6, Del7, Lio1, Lio2, Lio3}. That is why we use the term
\emph{ So\-nin--Poisson--Delsarte
  transmutation operators.}
  The known paper \cite{Lev7} about
So\-nin--Poisson--Delsarte
 operators is substantially based on classical Darboux results.

On the base of
 So\-nin--Poisson--Delsarte
 transmutation operators, Delsarte introduces the core notion of
 generalized translation operators.
 Numerous constructions of generalized harmonic analysis,  based
 on definitions of generalized translations and group structures
 introduced by means of them are developed.
 The research direction of generalized almost periodic functions, using  transmutation
operator of the
 So\-nin--Poisson--Delsarte
 type and generalized translation operators is founded in \cite{Del1,Del2} and continued
 in \cite{Del5,Del6, Lio1, Lio2, Lio3}.
 Independently this direction is exhaustively studied by Levitan
 in 1940 and 1947--1949;
 the results are included in  his classical monographs
  \cite{Lev2, Lev3, Lev4}.\footnote{Note that the  ``Almost periodic function'' entry
  of the Great Soviet Encyclopedia is
not exact: Delsarte is not mentioned and Levitan papers are dated
 by 1938. As far as the authors are aware, the first German-lingual
 paper of Levitan appeared in 1940.}
 For the first time, generalized Taylor expansions are constructed by Delsarte.
 Nowadays, they are reasonably called \emph{Tay\-lor--Delsar\-te
  series} (see \cite{Del2, Del7, Lev1, Lev2,Lev3};
  such series are studied in various papers up to now (see, e.\,g., \cite{Kam}).
 Note that the source and prototype of the majority of variants of the generalized
 harmonic analysis are Bessel operators and related differential equations.

Transmutation operators and generalized translation operators are
very important for the theory of partial differential equations
(see \cite{Mar9}). Using transmutation operators, one can
transform hard equations into simpler ones.
    In singular equations, generalized translation operators help to translate
    the singularity from the origin to an arbitrary point. Also,
 using them, one can construct fundamental solutions; then, using the generalized
 convolution, one can construct integral representations for solutions
 of the corresponding differential equations.

 For numerous generalizations of the Bessel operator, transmutation operators
 are considered as well. Transmutation operators for hyper-Bessel
 functions are important generalizations of
 So\-nin--Po\-is\-son--Del\-sarte
 operators. The theory of such functions is primarily founded in works of Kummer and Delerue.
 For the exhausted investigation of hyper-Bessel
 functions, the corresponding differential equations for them, and the corresponding
   transmutation
operators, see \cite{Dim, DK1,DK2}.
 The corresponding  transmutation
operators are reasonably called
 Sonin--Dimovski and Poisson--Dimovski
  transmutation
operators. They are studied in \cite{DK1,DK2, Kir1, Kir2, Kir3,
Kir4, Kir5} as well.
    In the theory of  hyper-Bessel
 functions and corresponding differential equations and
 transmutation operators, the famous Obreshkov integral transformation plays a key role.
 This transformation (in the general case, its kernel is represented through the Meyer
    $G$-function) generalizes the Laplace transformation, Mellin transformation, Fourier,
    Hankel, and Meyer
 sine-transformations and  cosine-transformations,
 and other classical integral
transformations.
 Various forms of hyper-Bessel
 functions and corresponding differential equations and
 transmutation operators as well as particular cases of the Obreshkov  transformation
 rediscovered many times afterwards.
 In our opinion,  the Obreshkov  transformation (together with the  Fourier, Mellin, and
 Laplace ones) belongs to a small set of fundamental ones: they are
 used (as blocks) to construct many other transformations as well
 as other objects and applications based on them.
 Historically, the  Obreshkov  transformation was the first known  integral transformation
 such that its kernel is expressed via the Meyer   $G$-function,
 but is not expressed via one (!) generalized hypergeometric function.
In the same way, the Stankovich integral transformation is
introduced:  its kernel is expressed via the Wright--Fox
 $H$-function, but is not expressed  via the Meyer $G$-function
 (which is simpler). The Stankovich  transformation is applied for the study of
  fractional-order differential equations of the fractional-diffusi\-on
  type (see \cite{Koch1,Koch2, EiIvKoch, Pshu1,Pshu2}).

Similar  transmutation operators are constructed for other
prototype operators such as
\begin{gather}
\label{510}
A=\frac{1}{v(x)}\frac{d}{dx}v(x)\frac{d}{dx},\\
v(x)=\sin^{2\nu+1}x,\, v(x)=\sinh^{2\nu+1}x, \,\textrm{or} \,
v(x)= (e^x-e^{-x})^{2\nu+1}(e^x+e^{-x})^{2\mu+1}
\end{gather}
 (see \cite{Car1, Car2, Car3, Yar1,Yar2}).
 For the theory, operators $A$ of kind \eqref{510} are important because,
 due to the famous Gel'fand relation, they represent the radial part of the Laplace operator
 on symmetric  spaces (see \cite{Hel1}). The Bessel operator is obtained if
 we assign $v(x)=x^{2\nu+1}$ in \eqref{510}.
 For the Airy operator $D^2+x$, transmutation operators are constructed as well.
 In \cite{Kach1}, its perturbed variant related to the quantum-mechanical
  Stark effect is considered. In \cite{Lev1995}, perturbed Hill operators
  with periodic potentials are considered. In \cite{Low1, Low2, Low3},
    Vekua--Erd\'elyi--Lowndes
  operators of the translation with respect to the spectral parameter are investigated.

The third (contemporary) period of the development of the theory
of transmutation operators includes the works since 1990s.
    A lot of important investigations refer to this period (see, e.\,g.,
 \cite{Mar9, CB, Lev1995, S42, S46, S38, S400, S401, SitDis}.
    In particular, the following research  directions are related to the application
    of methods of transmutation operators.
    The development of the theory of generalized analytic functions is continued
 (see \cite{Sol,Kli1, Kli2, Kli3, Kli5,Krav1}).
 Applications of transmutation operators to embeddings of function spaces
 and generalizations of Hardy operators (see \cite{S66, S6, S42, S46, S38, S400, S401})
 as well as to the
 constructing of various cases of the generalized translation and respective generalized
 variants of harmonic analysis (see
 \cite{Mar9,Gad1,Gul1,Gul2,Plat1, Plat2,
Plat3,Lyah1, Lyah2, Lyah3, LPSh1, LShFrac, LPSh2}) are found.
    In works of Muhli\-sov and his disciples, problems for $B$-potentials
    are considered.
    In works of Shishkina,
  Aisgersson-ty\-pe  me\-an-value
  theorems and new problems for the
   Euler--Poisson--Darboux
   equation are investigated, new classes of
   Riesz-type potentials with $B$-hyperbolic
   and ultrahyperbolic operators are constructed, and applications of these results
   to differential equations of the corresponding types are considered (see \cite{ShiR1,
ShiR2, ShiR3, ShiR4, ShiE1, ShiE2, SS, ShiE3, ShiE4, ShiE5,
ShiE6}).
   The applying of  transmutation operators and related methods in the theory
   of inverse problems and scattering theory is continued
    (see \cite{But, Ram1, Mar9, Yurko, ChCPR, PiSa}).
    For differential equations, the Darboux method and its modifications are developed
    (see \cite{Kap1, MaSa}).
    New classes of problems are considered for solutions with essential singularities
    on a part of the boundary at inner or angular points
    (see \cite{Kat1, 30,Kat2, 32, KatDis, Kat3, Kat4,KiKa3,KiKa4})
    and exact decay estimates are obtained for solutions of elliptic and ultraelliptic
   equations (see \cite{Mesh1,Mesh2,S3}).
    The usage of transmutation
operators for the investigation of fractional integral and
differential operators becomes a separate research direction (see
\cite{AKK, Dim, DK1,DK2, Kir1, Kir2, Kir3, Kir4, Kir5,ViRy,
Vir1,Vir2,Bac,SS,ShiE3,FJSS}).
    It is continued to use  transmutation-operator
    methods to study singular and degenerating
 boundary-value  problems and  pseudodifferential operators
(see \cite{Kat1, 30, Kat2, 32, KatDis, Kat3,
Kat4,KiKa3,KiKa4,Kan,Lyah1,Lyah2,Rep, Lar1, Lar2, Lar3, Lar4,
Lar5, Lar6}) as well as operator equations (see \cite{Glu1, Glu2,
Glu3, Glu4, Glu55, Glu5, Glu6, Glu7, Glu8, Glu9, Glu10, Glu11,
Glu12,FGP,FeIv}).
 Equations with Bessel operators and related areas are studied by Glushak
 (see \cite{Glu1, Glu2, Glu3, Glu4,
Glu55, Glu5, Glu6, Glu7, Glu8, Glu9, Glu10, Glu11, Glu12, Glu13}),
    Guliev (see \cite{Gul1,Gul2}), Lyakhov (see \cite{Lyah1, Lyah2,
Lyah3, LPSh1, LShFrac, LPSh2}), Pul'kina \cite{Pulkina},
    Sabitov (see \cite{SaIl}), Kravchenko (see \cite{CKT1, CKT2,
Krav1, Krav2, Krav3, Krav4, Krav5, Krav6, Krav7, Krav8, Krav9,
Krav10}) with his colleagues and disciples as well as Muravnik
(see \cite{Mur,Mur8,Mur9} for functional-differential
    equations and stabilization of solutions and \cite{Mur3,Mur4}
    for properties of the Fo\-urier--Bessel
    transformation), Volchkov (see \cite{Volch}), I.~Polovinkin
 (see \cite{LPSh1, LShFrac, LPSh2} for me\-an-value
 theorems for equations with Bessel operators), Shishkina
(see \cite{ShiR1, ShiR2, ShiR3, ShiR4, ShiE1, ShiE2, SS, ShiE3,
ShiE4, ShiE5, ShiE6, LPSh1, LShFrac, LPSh2} for $B$-hyperbo\-lic
    potentials and generalized averages), Sh. Karimov (see \cite{KarST,
KarST1, KarST2, KarST3, KarST4}), Hasanov (see  \cite{SaHa1},
E.~Karimov (see \cite{HaKa1}), Ergashev (see \cite{Erg1}),
Garipov, N.~Zaitseva, Mavlyaviev, Hushtova, and other researchers.

Recently, efficient numerical methods are developed to apply
transmutation  operators to compute solutions of differential
equations, their eigenfunctions, and spectral characteristics (see
\cite{CKT1, CKT2, Krav1, Krav2, Krav3, Krav4, Krav5, Krav6, Krav7,
Krav8, Krav9, Krav10}); these methods are based on the  spectral
parameter power series (developed by Kravchenko) for the
representation of solutions of Sturm--Liouville
    equations and their generalizations with Bessel or Dirac operators.
    Note that the idea to represent kernels of transmutation operators by series,
  being an alternative to the integral representation of transmutation operators,
 is natural because the kernels are obtained by the method of successive approximations,
  i.\,e., as Neumann series.
  Another series representation of transmutation operators, using generalized bases,
 is studied in \cite{FN}. For perturbed Bessel equations with
 variable potentials, series representations of kernels of transmutation operators
 is obtained in \cite{CFH,FH} in detail.

 Problems of
 Dirichlet--Neumann and Neumann--Dirichlet
 types, where a transmutation operator acts on
 boundary-value or initial-value
 conditions, preserving the differential expression, form a separate class.
 Such problems have important applications in mechanics (see \cite{BMYa, Yarem});
    also, this class of problems is closely related to the
    spectral theory and probability theory.
 Completed modifications of harmonic analysis for Bessel operators are constructed
 in \cite{Plat1, Plat2, Plat3}; for perturbed Bessel-type
  operators such modifications are constructed
 in \cite{Tri3, Tri4, Tri5}.
 Nowadays, harmonic analysis for
 Dunkl-type differential-difference
 operators is actively developed on the base of the corresponding generalizations of
   So\-nin--Poisson--Delsarte
   operators (see \cite{Me1, DHS, Dun1, Dun2, Dun3,
Gal1,Gal2, Rod, Ros1, Ros2, Tri6}).
 Once we have corresponding generalized translation
operators  (defined by means of transmutation operators), we are
able to define a generalized convolution and new algebraic and
group structures and to consider approximation problems for
functions (see \cite{AppHyp}). Fage's ideas developed for the
Bianchi in relation to the constructing of transmutation operators
of high-order
    differential equations are developed for more general equations in \cite{ZhM1,ZhM2}.
    Representations of solutions of fractional-order
    equations via solutions of integer-order
    equations (see the so-called
    subordination principle of Pruss, Kochubei, Bazhlekova, and Pskhu)
 can be treated as transmutation operators.
 Transmutation operators are applied in the theory of the Radon transformation
 and mathematical tomography (see \cite{Nat1, Hel1, Rub3, Rub1, Rub2, Rub4}
 and in expansions of functions in various series with respect to special functions
 (see \cite{Kam}). In \cite{Mar8,Mar9}, it is continued to apply
  transmutation operators to the quantum theory.

Also, the theory of transmutation operators is closely related to
the factorization  of differential operators (see \cite{Berk} and
works of A.~Shabat).
 For group properties of differential equations,
  the Schur complementability lemma is very important.
    It can be treated as the existence of a formal transmutation operator
    between a fractional integral operator and infinite-order
    differential operator.
    A problem of such kind is considered even in \cite{Burb}.
 In \cite{Car10}, it is attempted to construct  ``quantum'' transmutation
operators for $q$-differential
 operators.
 Various problems using ideas of transmutation operators or related methods
 are considered in \cite{Arsh, Bas, Han, GSPP, Kam, Matv, SPP} as well.

The possibility for the original and transmuted function
 to belong different spaces (traditionally emphasized by the usage of different
 notation for the variables) allows one to include all classical integral
 transformations (such as the Fourier one, Laplace one, Mellin one,
 Hankel one, Weierstrass one, Kantorovich--Lebedev
 one, Fock, Obreshkov one, Stankovich one, and others, see \cite{BE1,BE2, Ome})
  into the general scheme of
 transmutation operators.
 The Slater--Marichev
 theorem joining methods of the Mellin transformation with the theory of hypergeometric
 functions is applied for the computing of integrals needed to implement the method of
  transmutation operators (see \cite{Marich1, PBM}).
  Grinberg finite integral transformations (see \cite{Grin}) are included to the general
  scheme of transmutation operators as well.
  The investigation of the Green function and Riemann function used in the method of
  transmutation operators is continued in \cite{Ler, MaKiRe, Sob}.
  Transmutation operators are related to fractional (quadratic)
  Fourier (see \cite{OZK}) and Hankel integral transformations.
  The Chernyatin solution of the famous problem to find
  conditions for potentials such that they are necessary and sufficient
  for the justification of the Fourier method for the variable-coefficient
  wave equation opens new perspectives of the usage of this method to estimate
  kernels in the transmutation-operator
  theory (see \cite{Cher}).

Commuting operators of each nature satisfy the definition of
transmutation operators as well. Operators commuting with
derivatives are the most closed to the spirit and problems of the
transmutati\-on-ope\-rator theory.
    In this case, transmutation operators are frequently represented by formal series,
 pseudodifferential operators, or infinite-order
 differential operators. The description of commutants is directly related
  to the description of the whole family of transmutation operators (for a given pair)
  via its unique representative.
  In this class of problems, the theory of operator convolutions is applied; primarily,
  this refers to the Berg--Dimovski
  convolution  (see \cite{Dim,Bozh}).
  Also, applications in the theory of transmutation operators are
  found for results for commuting differential operators, coming from classical works
  of Burchnall and Chaundy.

An important area of the theory of transmutation operators is the
special class of Buschman--Erd\'elyi
    transmutation operators (see  Chap. \ref{ch3}).
    Under a suitable selection of parameters, this class of transmutation operators
    generalizes
 Sonin--Poisson--Delsarte
 transmutation operators and their inverse,
 Riemann--Liouville and   Erd\'elyi--Kober
 fractional integrodifferential operators, and
    Meler--Fock integral transformations. The term \emph{Buschman--Erd\'elyi
 operators} (being the most reasonable from the historical viewpoint) is introduced
  in \cite{S66, S6}; later, it is used by other authors as well.
  The importance of
Buschman--Erd\'elyi operators is mainly caused by their numerous
applications. For example, they occur in the following areas of
the theory of partial differential equations (see \cite{SKM}):
 the Dirichlet problem for the
  Euler--Poisson--Darboux
  in the quarter of plane
  and the setting of correspondence
  between values of solutions of the
   Euler--Poisson--Darboux
   equation on the initial-data
   manifold and the characteristic (see the Copson lemma above),
   the Radon transformation theory (because, due to the results of
    \cite{Lud, Nat1, Deans, Rub3, Rub1, Rub2, Rub4}, under the expansion
     of spherical harmonics,
 the action of the Radon transformation is reduced to the Buschman--Erd\'elyi
 with respect to the radial variable),
 the investigation of boundary-value
  problems for various equations with essential singularities inside the domain,
  the proof of the embedding of Kipriyanov spaces into weight Sobolev spaces,
  the finding of relations between transmutation operators and scattering-theory
  wave operators, and generalizations of the classical Sonin and Poisson
  integral representations and
  Sonin--Poisson--Delsarte
 transmutation operators.

For the most complete study of  Buschman--Erd\'elyi
    operators, see \cite{S1,S70,S72,S2,S73,S4,S66,S65,S6,S5,S7,S46,S14,S103,S400,SitDis,
SSfiz,S42,S94,S38,S401,S402}.
    Note that no role of Buschman--Erd\'elyi
    operators as transmutation operators is realized and considered before the cited
    works.

Essential parts of monographs
 \cite{CSh, Car1, Car2, Car3, Car4, Berg, GB, Lio1, BMYa, Yarem}
 are devoted to the theory of transmutation operators and their applications.
 Various aspects of transmutation operators are considered in
 \cite{SKM, Col1,Col2, Gil1,Gil2} and a number of other known monographs as well as
 in reviews \cite{CB,Lev1995,Mar9,S46,S42}.
  No books completely devoted to transmutation operators and comparable with
  \cite{CSh, Car1, Car2, Car3, Car4} are published in Russian yet.
  Therefore, the goal of the  present monograph is to fill this lacuna.
 Undoubtedly, \cite{FN} deserves to be mentioned.
 In this monograph, no results from the theory of transmutation operators are reflected.
 However, this is completely compensated by the presentation of the results
  of its authors about the constructing of such operators for
  high-order differential operators with variable coefficients, which is
   a very hard problem from the theory of transmutation operators.
  The cited monograph includes many other areas:
  the problem on operators commuting with derivatives in the space of analytic functions
  is solved (including the correction of erroneous results of Delsarte and Lions),
  the complete solvability theory for the Bianchi equation is created,
  the theory of operator-analytic functions (originated from \cite{Mar5, Mar6, Mar7})
  is created, and  differential, integral, and root-extraction
   operators are investigated in spaces of analytic functions.

In \cite{S41, S45, S140, S53}, a refinement method for
    Cauchy--Schwarz--Bunyakovskii
    inequalities, based on the method of means, is developed.
 Integral refinements obtained by this method  can be applied to norm estimates
 of various transmutation operators.

 Various estimates of special functions, providing a possibility to estimate kernels
  of transmutation operators (see, e.\,g., \cite{S61, S9, S12, S13, S15, S16,
S24, MS1, MS2, MS3}), are useful as well.

    Thus, methods of the theory of transmutation operators and related problems
    are applied  in works of many mathematicians such as Aliev,
 Begehr,  Bergman,  Betancor,  Boumenir, Braaksma,
 Bragg, R.~Carroll, Castillo-P\'erez,
 Chebli, Dimovski, Dunkl, Delsarte, Fitouhi,  Gasmi, R.~Gilbert,
 Hamza,  Holzleitner, Hriniv, Hristov,
 Hutson, Kalisch, Kalla, Koornwinder,
 Kiryakova, L\H{o}ffstr\H{o}m, J.~Lions,  Luchko, Moro,
 Mykytyuk,  Pym,  Rubin,  Santana-Bejarano,
 Santosa,  Siersma,  Sifi,  Sinclair,  de Snoo,
 Spiridonova,  Stempak,  Teschl,  Thyssen,
 Trim\`{e}che,  Tsankov,   Voit,  Tuan,
Z.~Agranovich,  Androshchuk, Bavrin, Baskakov, Britvina, Buterin,
Valitskii, Volk, Volchkov, Gadzhiev, Glushak, M.~Gorbachuk,
Gohberg, Guliev, Guseinov, Zhitomirskii, L.\,A.~Ivanov,
Kh.~Ishkin, Eliseeva, Eremin, Kapilevich, Sh.~Karimov, Karp,
Katrakhov, A.~Kachalov, Kilbas, Kipriyanov, Klyuchantsev,
Kononenko, Korobeinik, Kostenko, Kravchenko, M.~Krein, Kulish,
Kushnirchuk, G.\,I.~Laptev, B.~Levin, Levitan, Leontiev,
N.~Linchuk, S.~Linchuk, Lyatifova, Lyahov, Lyahovetskii, Malamud,
Marchenko, Matrosov, Matsaev, Muravnik, Nagnibida, Nizhnik,
Olevskii, Parfeniova, Prikarpatskii, S.~Platonov, Povzner, Rubin,
Rofe-Beketov, Sabitov, Samoilenko, Sakhnovich, Sitnik, Sokhin,
Stashevskaya, Torba, L.~Faddeev, Fage, Fishman, Khachatryan,
Khromov, Chudov, Shishkina, Shmulevich, Yurko, Yaremko, and
Yaroslavtseva.
    Definitely, this list is far from to be complete; it can be substantially extended.

It follows from the provided analysis that the method of
transmutation operators is efficient for the theory of
differential equations and related areas of mathematics, a lot of
works are devoted to it, and numerous classes of problems are
solved on its base.
 However, there are substantial lacunas and numerous unsolved problems.
 For instance, there is no detailed classification (with the description of main properties)
 for transmutation operators intertwining differential operators or solutions of differential
 equations with singularities at coefficients (including Bessel differential operators).
 Properties and applications of the basic class
of Sonin and Poisson transmutation operators are studied and
described in published papers in detail, but there is no
systematic explanations and proofs for many properties of
Buschman--Erd\'elyi operators, which are their important
generalizations.
    Before Sitnik works, it was not realized that Busch\-man--Erd\'elyi
    integral operators are transmutation operators for differential operators with Bessel
    operators. No general schemes to construct
transmutation operators of needed classes are developed to the
stage providing a possibility to construct (on their base)
explicit relations intertwining solutions of various differential
equations. Papers discovering the relation between transmutation
operators with main constructions of fractional calculus are
almost absent.
 No possibilities of the applying of
transmutation operators to proofs of the embedding of function
spaces such as Sobolev spaces and Kipriyanov spaces (including
energy spaces for singular partial differential equations with
Bessel operator with respect one or several variables)  are
considered. Imperfect methods are used to construct transmutation
operators for second-order
 differential operators with variable coefficients: they provide
 rough estimates for kernels of transmutation operators with indefinite constants,
 while inexact requirements for coefficients of differential equations lead to the restriction
 of their classes, e.\,g., classes of admissible potentials for Sturm--Liouville
  problems and their generalizations for differential equations with singularities at
  coefficients.
  Methods of transmutation operators are almost unknown to be applied to obtain
  exact estimates of solutions of differential equations, e.\,g., in problems like the known
  Landis problem.
  Also, the following paradoxical situation takes place:
  fractional powers of the Bessel operator, used in many papers, are defined only implicitly
  (in terms of the Fourier--Bessel
  or Hankel transformation); there are no relations to define them in the integral form
  explicitly though the theory of classical fractional Riemann--Liouville
  integrals is originated from the said representations. No explicit representations
  are obtained for various simple and natural constructions of transmutation operators
  for standard pairs of differential operators.
  General schemes for estimates of kernels of transmutation operators, requiring refinements
  and generalizations of classical inequalities,  are not introduced and
considered in broadly used function spaces. Recently obtained
exact inequalities for various special functions are not applied
for estimates of kernels of transmutation operators.

Solutions of a number of problems mentioned above are provided in
the present monograph.

\section{Main Definitions, Notation, and Properties:\\
 Special Functions, Function Spaces, and Integral Transformations}\label{sec3}
\sectionmarknum{Main Definitions, Notation, and Properties}

\subsection{Special functions}\label{sec3.1}

Here, we provide brief definitions and explanations for special
functions used in our monograph. We follow \cite{AS, BE1, BE2,
BE3, Luke3, Luke1,Luke2, Wit1, AAR, PBM123, NIST}.
 Also, a number of short historical comments is provided.

\subsubsection{Gamma-functions,
 beta-functions, psi-functions,
 Pochhammer symbols, and binomial coefficients}\label{sec3.1.1}

The gamma-function generalizes the factorial notion to the case of
numbers different from positive integers.
    In the general case, the  beta-function
    is defined via
 gamma-functions. The psi-function
    is the logarithmic derivative of the gamma-function.

    Let $z\in\mathbb{C}.$ According to Euler, the \textit{gamma-function}
     $\Gamma(z)$ is defined as the following limit:
$$
\Gamma(z)=\lim\limits_{N\rightarrow\infty}\frac{N!N^z}{z(z+1)(z+2)\ldots(z+N)},\,\,
z\in\mathbb{C}.
$$
    More frequently, the following definition via the Euler
    second-kind integral is used:
\begin{equation}\label{Gamma}
\Gamma(z)=\int\limits_{0}^\infty y^{z-1}e^{-y}dy,\qquad \Re z>0;
\end{equation}
    this integral converges for all complex $z$ such that $\Re x>0.$

Integrating relation \eqref{Gamma} by parts, we arrive at the
recurrent relation
\begin{equation}\label{Rec}
\Gamma(z+1)=z\Gamma(z).
\end{equation}
    Since $\Gamma(1)=1,$ it follows that, if $n$ is a positive integer,
    then the recurrent relation \eqref{Rec} leads to the relation
$$
\Gamma(n+1)=n\Gamma(n)=n(n-1)\Gamma(n-1)=\ldots=n(n-1)\cdot\ldots\cdot2\cdot
1\cdot\Gamma(1)
$$
 or
$$
\Gamma(n+1)=n!,
$$
    which allows one to treat the gamma-function
    as a generalization of the factorial notion.
  Presenting relation \eqref{Rec} in the form
\begin{equation}\label{ReGa}
\Gamma(z-1)=\frac{\Gamma(z)}{z-1},
\end{equation}
    we obtain an expression providing a possibility to define the gamma-function
    for negative value of the independent variable (the definition via
    \eqref{Gamma} is not acceptable in such a case).
    Relation \eqref{ReGa} shows that  $\Gamma(z)$ has jumps at the
    points  $z=0,-1,-2,-3,\ldots$

Applying relation \eqref{ReGa} many times, we obtain  the
following \emph{lowering and  raising} relations:
\begin{equation}\label{Povysh}
\Gamma(z+n)=z(z+1)\ldots(z+n-1)\Gamma(z),\qquad n=1,2,\ldots,
\end{equation}
    and
\begin{equation}\label{Ponizh}
\Gamma(z-n)=\frac{\Gamma(z)}{(z-n)(z-n+1)\ldots(z-1)},\qquad
n=1,2,\ldots
\end{equation}
    Note that
$$\Gamma\left(\frac{1}{2}\right)=\sqrt{\pi},\qquad
\Gamma\left(\frac{1}{2}+n\right)=\frac{(2n)!\sqrt{\pi}}{4^nn!},\qquad\textrm{and}\qquad
\Gamma\left(\frac{1}{2}-n\right)=\frac{(-4)^nn!\sqrt{\pi}}{(2n)!}.
$$
    The complement relation
\begin{equation}\label{Dopol}
\Gamma(z)\Gamma(1-z)=\frac{\pi}{\sin{z\pi}}
\end{equation}
 and the duplication (Legendre) formula
\begin{equation}\label{Lezh}
\Gamma(2z)=\frac{2^{2z-1}}{\sqrt{\pi}}\Gamma(z)\Gamma\left(z+\frac{1}{2}\right)
\end{equation}
    hold.

    The
 \emph{beta-function} $B(z,w)$ is closely related to the gamma-function.
 If the parameters  $z$ and $w$ satisfy the conditions $\Re z>0$
 and $\Re w>0,$ then the Euler beta-function
 is defined by the Euler first-kind
  integral
\begin{equation}\label{Beta}
B(z,w)=\int\limits_{0}^1t^{z-1}(1-t)^{w-1}dt.
\end{equation}
    If $\Re z\leq 0$ and $\Re w\leq 0$, then the
 beta-function is defined by the relation
\begin{equation}\label{BetaGamma}
B(z,w)=\frac{\Gamma(z)\Gamma(w)}{\Gamma(z+w)}.
\end{equation}

The
 \emph{psi-function} $\psi(z)$ is defined as the logarithmic derivative of the gamma-function:
$$
\psi(z)=\frac{d\,\ln\Gamma(z)}{dz}=\frac{\Gamma'(z)}{\Gamma(z)}.
$$
    The function $\psi(z)$ undergoes jumps at the points
$z=0,-1,-2,\ldots$
    The psi-function
    satisfies the representation
$$
\psi(z)=-\gamma+(z-1)\sum\limits_{n=0}^\infty\frac{1}{(n+1)(z+n)},
$$
    where
$$\gamma=\lim\limits_{m\rightarrow\infty}\left(\sum\limits_{n=1}^m\frac{1}{n}-\ln{m}\right)=0
\p5772156649\ldots$$
    denotes the Euler--Mascheroni
    constant (see \cite{AS}).
    It is obvious that $\psi(1)=-\gamma.$
    Note that the relation
\begin{equation}\label{Psi}
\int\limits_0^1\frac{t^x-t^y}{1-t}dt=\psi(y+1)-\psi(x+1)
\end{equation}
    holds.

 For integer values of $n$, the \emph{Pochhammer symbol} $(z)_n$
 is defined by the relation
$$
(z)_n=z(z+1)\ldots(z+n-1),\,\, n=1,2,\ldots,\,\,\textrm{and}\,\,
(z)_0\equiv 1.
$$
    The relations $ (z)_n=(-1)^n(1-n-z)_n,$ $ (1)_n=n! ,$ and
\begin{equation}\label{Poh}
(z)_n=\frac{\Gamma(z+n)}{\Gamma(z)}
\end{equation}
 hold.
    Relation \eqref{Poh} can be used to introduce the symbol
$(z)_n$ for real (complex) values of $n.$

    The \emph{{binomial coefficients}} are defined as follows:
$$
\biggl(\begin{array}{c}
$$\alpha$$ \\
$$n$$ \\
\end{array}\biggr)=\frac{(-1)^{n-1}\alpha\Gamma(n-\alpha)}{\Gamma(1-\alpha)\Gamma(n+1)}.
$$
    In particular, for $\alpha=m,$ $m=1,2,\ldots,$ we have the
    relations
$$
\biggl(\begin{array}{c}
$$m$$ \\
$$n$$ \\
\end{array}\biggr)=\frac{m!}{n!(m-n)!},\qquad m\geq n.
$$
    For arbitrary (complex) $\beta$ and $\alpha,$
$\alpha\neq -1,-2,\ldots,$ we assign
$$
\biggl(\begin{array}{c}
$$\alpha$$ \\
$$\beta$$ \\
\end{array}\biggr)=\frac{\Gamma(\alpha+1)}{\Gamma(\beta+1)\Gamma(\alpha-\beta+1)}.
$$
    For positive integer values of $k$, the relation
\begin{equation}\label{bin}
(-1)^k\biggl(\begin{array}{c}
$$\alpha$$ \\
$$k$$ \\
\end{array}\biggr)=\biggl(\begin{array}{c}
$$k-\alpha-1$$ \\
$$k$$ \\
\end{array}\biggr)=\frac{\Gamma(k-\beta)}{\Gamma(-\alpha)\Gamma(k+1)}.
\end{equation}
    holds.

%\begin{center}
%   \textbf{Функции Бесселя}
%\end{center}

\subsubsection{Bessel functions}\label{sec3.1.2}

The Bessel functions are defined as solutions of the Bessel
 differential equation
$$
x^2 \frac{d^2 y}{dx^2} + x \frac{dy}{dx} + (x^2 - \alpha^2)y = 0,
$$
    where the order $\alpha$ is an arbitrary complex number.

\emph{First-kind Bessel functions} denoted by $J_\alpha(x)$ are
solutions finite at the origin for integer or nonnegative
$\alpha.$
    These functions can be defined by means of the Taylor
    expansion around the origin or more general power series (if
    $\alpha$ is not integer):
$$
J_\alpha(x) = \sum\limits_{m=0}^\infty \frac{(-1)^m}{m!\, \Gamma(m+\alpha+1)}
 {\left({\frac{x}{2}}\right)}^{2m+\alpha}.
$$
    If $\alpha$ is not integer, then the functions $J_\alpha (x)$
    and $J_{-\alpha} (x)$ are linearly independent. Hence, they are solutions of the equation.
    If $\alpha$  is integer, then the following relation holds:
$$
J_{-\alpha}(x) = (-1)^{\alpha} J_{\alpha}(x).
$$
    This means that the functions are linearly  dependent (in this case).
    Then the second solution of the equation is the second-kind
    Bessel function, i.\,e., the
\textit{Neumann function}, which is the solution $N_\alpha(x)$ of
the Bessel equation, tending to infinity  as $x\to 0.$
    This function is related to $J_\alpha(x)$ as follows:
$$
N_\alpha(x) = \frac{J_\alpha(x) \cos(\alpha\pi) - J_{-\alpha}(x)}{\sin(\alpha\pi)}
$$
 (if $\alpha$ is integer, then the limit with respect to $\alpha,$
computed, e.\,g., by means of the L'Hospital rule, is taken). A
linear combination  of the
 first-kind and second-kind
 Bessel function is a complete solution of the Bessel equation:
$$
y(x) = C_1 J_\alpha(x) + C_2 N_\alpha(x).
$$
The notation $N_\alpha(x)=Y_\alpha(x)$ is frequently used as well.

\textit{Modified Bessel functions or Bessel functions of imaginary
independent  variables} are the function
$$
I_\nu(x)=i^{-\nu} J_\nu(ix)
$$
    and the  \textit{MacDonald function}
$$
K_\nu(x)= \frac{\pi}{2\sin(\pi\nu)} \left[I_{-\nu}(x) - I_\nu(x)\right], \nu\notin\mathbb{Z}.
$$
If $ \nu\in\mathbb{Z}$, then the MacDonald function is computed by
means of the limit passage with respect to the index  (the
L'Hospital rule is applied).

\textit{Hankel functions} or third-kind
 Bessel functions are linear combinations of
 first-kind and second-kind
 Bessel functions.
 Hence, they satisfy the Bessel equation as well.
    The first-kind
    Hankel function is
$$H_{\nu}^{(1)}(z)=J_{\nu}(z)+iN_{\nu}(z).$$
    The second-kind
    Hankel function is
$$H_{\nu}^{(2)}(z)=J_{\nu}(z)-iN_{\nu}(z).$$
The zero-index
    Hankel functions are fundamental solutions of the Helmholtz equation.
 They are represented by the  first-kind
 Bessel functions as follows:
$$
H_{\nu}^{(1)} (z) = \frac{J_{-\nu} (z) - e^{-\nu\pi i} J_{\nu}(z)}{i \sin (\nu\pi)}
    \qquad  \textrm{and}  \qquad
 H_{\nu}^{(2)} (z) = \frac{J_{-\nu} (z) -e^{\nu\pi i} J_{\nu} (z)}{- i \sin (\nu\pi)}.
$$

The \textit{normalized Bessel function} (the $j$-small
  Bessel function) $j_\nu$ is defined by the relation
\begin{equation}\label{FBess1}
j_\nu(x) ={2^\nu\Gamma(\nu+1)\over x^\nu}\,\,J_\nu(x),
\end{equation}
    where $J_\nu$ is  the  first-kind
 Bessel function
(see \cite[p.~10]{Kip1} and \cite{Lev1}).
    For \eqref{FBess1}, the following
relation holds:
\begin{equation}\label{RavDBes1}
T^y_x
j_{\frac{\gamma-1}{2}}(x)=j_{\frac{\gamma-1}{2}}(x)\,j_{\frac{\gamma-1}{2}}(y)
\end{equation}
(see, e.\,g.,\cite{Lev1}).

In \cite{Kuz1,Kuz2}, useful properties of Bessel functions with
applications to the analytic number theory are obtained.

\subsubsection{Gauss hypergeometric functions}\label{sec3.1.3}

The Gauss hypergeometric function is defined inside the disk
 $\{|z|{<}1\}$ as the sum of the hypergeometric  series
\begin{multline}\label{FG}
\,_2F_1(a,b;c;z)=F(a,b,c;z)=\sum\limits_{k=0}^\infty\frac{(a)_k(b)_k}{(c)_k}\frac{z^k}{k!}
\\
=\frac{\Gamma(c)}{\Gamma(b)\Gamma(c-b)}\int\limits_0^{1}t^{b-1}(1-t)^{c-b-1}(1-zt)^{-a}dt,
\,\,\Re c>\Re b>0
\end{multline}
  (see \cite[p. 373, relation 15.3.1]{AS}).
 For $|z|\geq1$, it is defined as the analytic continuation of
 this series (see \cite{AKdF, Bai}).
 In relation \eqref{FG}, the parameters $a,b,$ and $c$ and the
 variable  $z$ are allowed to be complex, $c $ is different from each nonpositive integer,
 and $(a)_k$ is the Pochhammer symbol defined by \eqref{Poh}.

The hypergeometric series \eqref{FG} converges only in the unit
disk of the complex plane. This causes the necessity to construct
the analytic continuation of the hypergeometric function to the
whole complex plane.
    One way to do this is to use the Euler integral representation
$$_2F_1(\alpha,\beta,\gamma;z) = { \Gamma(\gamma) \over
    \Gamma(\beta)\Gamma(\beta-\gamma) } \int\limits_{0}^{1}
t^{\beta-1} (1-t)^{\gamma-\beta-1} (1-tz)^{-\alpha} \,dt,
$$
$$
0<\Re \beta<\Re \gamma,\,\,\,\,|\arg(1-z)|<\pi;
$$
    the right-hand
    side is defined under the specified conditions providing the convergence of the integral.

    It is important that various special and elementary functions
    can be obtained from the hypergeometric function if we assign
    particular values to parameters and transform the independent
variable.

Examples for elementary functions are as follows:
$$
(1+x)^n = F(-n,\beta,\beta;-x),\qquad {1 \over x} \ln(1+x) =
F(1,1,2;-x),\qquad e^x = \lim\limits_{n \to \infty} F(1,n,1;{x
\over n}),
$$
$$ \cos x = \lim\limits_{\alpha,\;\beta \to \infty}
F\left(\alpha,\beta,\frac{1}{2}; -\frac{x^2}{4 \alpha
    \beta}\right),\; \textrm{and}\;
    \cosh x = \lim\limits_{\alpha,\;\beta \to \infty}
F\left(\alpha,\beta,\frac{1}{2};{ x^2 \over 4 \alpha
    \beta}\right).
$$
 The first-kind
 Bessel function and the  Gauss hypergeometric function are
 related as follows:
$$ J_\nu(z)= \lim\limits_{\alpha,\;\beta \to \infty} \left[
\frac{\left(\dfrac{z}{2}\right)^\nu}{\Gamma(\nu+1)}
F\left(\alpha,\beta,\nu+1; -\frac{z^2}{4 \alpha \beta}\right)
\right].
$$

\subsubsection{Legendre functions}\label{sec3.1.4}

The Legendre functions $P_{\nu}^\mu (x)$ and $Q_{\nu}^\mu (x)$ are
generalizations of Legendre polynomials to noninteger  degrees
(see \cite{BE1}).

 The Legendre functions $P_{\nu}^\mu (x)$ and $Q_{\nu}^\mu (x)$ are the solutions
 of the general Legendre equation
 $$
 (1-x^{2})\,y''-2xy'+\left[\lambda (\lambda +1)-{\frac {\mu ^{2}}{1-x^{2}}}\right]\,y=0,\,
$$
where the complex numbers  $\lambda$ and $\mu$ are called the
    degree and order (respectively) of the Legendre function.
    For complex values of the parameters and independent variable, these functions
    can be determined as follows:
$$
 P_{\lambda }^{\mu }(z)={\frac {1}{\Gamma (1-\mu )}}\left[{\frac {1+z}{1-z}}\right]^{\mu /2}
 \,_{2}F_{1}(-\lambda,\lambda +1;1-\mu;{\frac {1-z}{2}}),\qquad |1-z|<2,
$$
    and
$$ Q_{\lambda }^{\mu }(z){=}{\frac {{\sqrt {\pi }}\ \Gamma (\lambda {+}\mu {+}1)}
{2^{\lambda +1}\Gamma (\lambda {+}\frac{3}{2})}}{\frac {e^{i\mu \pi }
(z^{2}-1)^{\mu /2}}{z^{\lambda +\mu +1}}}\,_{2}F_{1}
\left({\frac {\lambda +\mu +1}{2}},{\frac {\lambda +\mu +2}{2}};\lambda +{\frac {3}{2}};
{\frac {1}{z^{2}}}\right),\quad |z|>1,
$$
    where $\Gamma$ is the gamma-function
    and $_{2}F_{1}$ is the hypergeometric function.
    Also, we use the direct values of the functions $P_{\nu}^\mu (x)$ and $ Q_{\nu}^\mu (x)$
    on the cut $x\in [-1;1],$ denoted by $\mathbb{P}_{\nu}^\mu (x)$ and
     $ \mathbb{Q}_{\nu}^\mu (x).$

 Legendre functions are special cases of the  Gauss hypergeometric function.
 Therefore, numerous series expansions, integral representations,
 and relations for the continuation as well as expressions via
 elementary functions for special values of parameters are known for them.

Legendre polynomials are Legendre functions of order $\mu=0$ for
nonnegative integers  $\lambda=n.$

\subsubsection{Mittag-Leffler functions}\label{sec3.1.5}

    The
 Mittag-Leffler function $E_{\alpha,\beta}(z)$ is the entire function (in $z\in\mathbb{C}$)
 of order $1/\alpha$ defined by the power series
\begin{equation}\label{ML}
E_{\alpha,
\beta}(z)=\sum\limits_{n=0}^{\infty}\frac{z^n}{\Gamma(\alpha
n+\beta)},\;\; z\in\mathbb{C},\;\alpha, \beta\in\mathbb{C},\;
\Re\alpha>0,\; \Re\beta>0
\end{equation}
 (see \cite{Dzh1, Djr2, BE3, GKMR, Kir5,Kir6, Jor, ViGa}).

 For $\alpha=1$, function \eqref{ML} is introduced by Mittag-Leffler
 in 1903. For the general case, it is introduced by Wiman in 1905.
 First, Mittag-Leffler
 and Wiman apply these functions in complex analysis (nontrivial examples of entire functions
 with noninteger growth orders and generalized summing methods).
 In USSR, these functions become broadly known from the famous
  Djrbashian monograph \cite{Dzh1}
 (also, see his monograph \cite{Djr2}).
 The most known application of
  Mittag-Leffler functions in the theory of integrodifferential equations
  and fractional calculus is the following fact:
  the resolvent of the Ri\-e\-mann--Liouville
  fractional integral is explicitly expressed via them according to the famous
  Hil\-le--Ta\-markin--Djrbashian
  relation  (see \cite{HT, SKM, Pshu1}).
  In \cite{GoMa}, due to numerous applications to frac\-ti\-o\-nal-order
  differential equations, this function is reasonably called the
 \textit{royal function of fractional calculus}.

The derivative of the Mittag-Leffler
 function is computed as follows:
$$
E_{\alpha,\beta}'(z)=\frac{d\,E_{\alpha,\beta}(z)}{dz}=
\sum\limits_{k=0}^\infty\frac{(1+k)z^k}{\Gamma(\beta+\alpha
    (1+k))}.
$$
Note that $E_{1,1}(z)=e^z.$

\subsubsection{Generalized hypergeometric-type
 functions}\label{sec3.1.6}

    Functions of the hypergeometric type are defined by a complex integral of the
  Mellin--Barnes type \cite{Marich1,PBM}.
 Generalized functions of the hypergeometric type are considered,
 e.\,g., in monographs \cite{Marich1, Sla, VN, Dwo1, Ext, KaSr, KiSa}.

For definite integrals, the computational method using generalized
 hypergeometric functions and based on the Slater--Marichev
 theorem (see \cite{Marich1,PBM}) is quite important.
 An original method to sum rather arbitrary series in terms of generalized
  hypergeometric functions and zeta-functions
   is developed in \cite{Ryko1,Ryko2}; however, it is almost out
 of use nowadays.
 Contemporary methods of theoretical and computer summing based on
   hypergeometric functions are considered in \cite{PWZ,Koepf}.

Numerous generalizations of  hypergeometric functions, e.\,g.,
$q$-hypergeometric functions (their theory comes from Heine),
Gel'fand--Graev $A$-hypergeometric
  systems, Spiridonov elliptic hypergeometric functions, and
  others.

 Lately, functions of the hypergeometric type (such as Mittag-Leffler
 functions, Wright functions, and Fox functions) are broadly applied
 in problems of probability theory and mathematical statistics. In particular important
 distribution functions, probability densities, and their characteristics are expressed
 via them.

\textit{Wright functions} are introduced in 1935 (see \cite{BE1,
ViGa, Jor, Kir1, Kir5,Kir6}). They intermediately generalize
 Mittag-Leffler functions and are defined by series similar to \eqref{ML};
 the difference is that the fraction for the common term of the series might contain an
 arbitrary finite amount of gamma-functions
 both and the numerator and denominator.

Let $p,q\in\mathbb{N}_0=\{0,1,2,\ldots\},$ $p^2+q^2\neq 0,$
$a_i,b_j\in\mathbb{C},$ and $\alpha_i,\beta_j\in\mathbb{R}$
($\alpha_i,\beta_j\neq0,$ $i=1,2,\ldots,p,$ $j=1,2,\ldots,q$).
    Then the Wright function is defined by the power series
$$
\,_p\Psi_q(z)=\sum\limits_{m=0}^\infty
\frac{\prod\limits_{i=1}^p\Gamma(a_i+\alpha_ik)}{\prod\limits_{j=1}^q\Gamma(b_j+\beta_jk)}
\frac{z^k}{k!},\qquad z\in\mathbb{C}.
$$
    Other terms, e.\,g., generalized Bessel functions (see \cite{Jor}), generalized
Mittag-Leffler functions (see \cite{GKMR}), multiparameter Mittag-Leffler
 functions (see \cite{Kir5,Kir6, ViGa}, etc  are used to call these functions as well.
 Wright functions represented by  Mellin images or
Mellin--Barnes integrals are Fox functions and vise versa (though
    such transformations are possible not for all values of parameters).
    Other terms are used for this class of functions as well: multiparameter
    hypergeometric functions, generalized or multiparameter Mittag-Leffler
     functions, multiparameter Bessel functions, etc.
  Note that the term  \emph{Bessel--Maitland
  functions} is crude: no mathematician Maitland ever exist and this is just the middle name
  of Wright; unfortunately, this inaccuracy is reproduced in the very respectable book
  \cite{BE2}.

An example of applications is as follows: the following variant
    of the Wright function is used to represent the resolvent
    for fractional powers of the Bessel operator:
\begin{equation}\label{GML}
J_{\gamma,\lambda}^\mu(z)=\sum\limits_{m=0}^\infty\frac{(-1)^m}
{\Gamma(\gamma+m\mu+\lambda+1)\Gamma(\lambda+m+1)}
\left(\frac{z}{2}\right)^{2m+\gamma+2\lambda}
\end{equation}
    (see \cite{S700,SS}).

The \textit{Meyer $G$-function}
 (see \cite{Dwo2, Sie}) is introduced in 1946.

The general definition of the  Meyer $G$-function
 is given by the following integral in the complex plane:
$$
G_{p,q}^{\,m,n} \!\left( \left. \begin{matrix} a_1, \dots, a_p \\
b_1, \dots, b_q \end{matrix} \; \right| \, z \right) =\frac{1}{2
\pi i} \int\limits_L \frac{\prod\limits_{j=1}^m \Gamma(b_j - s)
\prod\limits_{j=1}^n \Gamma(1 - a_j +s)} {\prod\limits_{j=m+1}^q
\Gamma(1 - b_j + s) \prod\limits_{j=n+1}^p \Gamma(a_j - s)} \,z^s
\,ds.
$$
    This is an integral of the so-called
    Mellin--Barnes type. It is well defined under the following assumptions:
     $ 0 \leq m' \leq q,$
$0\leq n\leq p,$ $m, n, p, q\in\mathbb{Z},$ and $a_k- b_j\neq 1,
2, 3,\ldots$ for $k = 1, 2,\ldots, n$ and $j = 1, 2,\ldots, m.$
    Under these assumptions, there no coinciding poles  $\Gamma(b_j-s),$ $j = 1,
2,\ldots, m,$ and $\Gamma(1-a_k+s),$ $k = 1, 2,\ldots, n,$
$z\neq0.$ Moreover, once the circuit over the region with poles is
over, the contour is continued up and down along the same vertical
line.

Is another contour is selected (such that the infinity is achieved
along horizontal lines), then new modifications of Meyer functions
are obtained (see \cite{KaPr3, KaLo1, KaLo2}).
 In such cases, it is possible that the Meyer function is defined, but its  Mellin transform
 dos not exist, or it exists, but the function itself is not restored by means of the inverse
 Mellin transformation (via the Mellin--Barnes
 integral).

In all cases, the   Meyer $G$-function
 is representable as a finite sum of generalized
    hypergeometric functions and satisfies a differential equation with polynomial
    coefficients (it can be presented explicitly).
     Meyer  functions are used, e.\,g., to represent the Obreshkov transformation
     (see \cite{Dim, Kir1}).

A special case of the Meyer function $G_{p,0}^{\,p,p}$ important
for applications is comprehensively investigated by  N\o rlund;
 this is the reason to propose the term \emph{Meyer--N\o rlund
 function} (see \cite{KaPr3, KaLo1, KaLo2}).

\textit{Fox functions} (see \cite{KiSa}). The $H$-function
 is introduced by Fox in 1961 (see \cite{KiSa, Bra}).

It is defined by means of the Mellin--Barnes
 integral as well:
$$
H(x)=H^{m,n}_{p,q}(z)=H^{m,n}_{p,q}\left[z\biggl|\left(
\begin{array}{c}
$$(a_p,A_p)$$ \\
$$(b_q,B_q)$$ \\
\end{array}
\right) \right]= H_{p,q}^{\,m,n} \!\left[ z \left| \begin{matrix}
(a_1, A_1) & (a_2, A_2) & \ldots & (a_p, A_p) \\
(b_1, B_1) & (b_2, B_2) & \ldots & (b_q, B_q) \end{matrix} \right.
\right]
$$
$$
=\frac{1}{2\pi i}\int\limits_L \frac
{\left(\prod\limits_{j=1}^m\Gamma(b_j+B_js)\right)\left(\prod\limits_{j=1}^n\Gamma(1-a_j-A_js)\right)}
{\left(\prod\limits_{j=m+1}^q\Gamma(1-b_j-B_js)\right)\left(\prod\limits_{j=n+1}^p\Gamma(a_j+A_js)\right)}
z^{-s} \, ds,
$$
    where each empty product is assigned to be equal to  one,
 $m,n,p,q\in \mathbb{N}_0,$ $0\leq n\leq p,$ $1\leq m\leq q,$
$A_i, B_j\in \mathbb{R}_+,$ $a_i$ and $b_j$ are either real or
complex, $i=1,\ldots,p,$ $j=1,\ldots,q,$ $L$ is a suitable contour
 separating the poles
$$
\zeta_{j\nu}=-\left(\frac{b_j+\nu}{B_j}\right),\qquad
j=1,\ldots,m,\qquad \nu=0,1,2,\ldots,
$$
 of the gamma-function
  $\Gamma(b_j+sB_j)$ from the poles
$$
\omega_{\lambda
k}=\left(\frac{1-a_\lambda+k}{A_\lambda}\right),\qquad
\lambda=1,\ldots,n,\qquad k=0,1,2,\ldots,
$$
of the  gamma-function
 $\Gamma(1-a_\lambda-sA_\lambda)$ such that
$$
A_\lambda(b_j+\nu)\neq B_j(a_\lambda-k-1),  j=1,\ldots,m,
\lambda=1,\ldots,n, \nu,k=0,1,2,\ldots,
$$
and, once the circuit over the region with poles is over, the
integration contour is continued up and down along the same
vertical line.

Is another contour is selected (such that the infinity is achieved
along horizontal lines), then new modifications of Fox functions
are obtained (see \cite{KaPr1,KaPr2}).
 In such cases, it is possible that the Fox function is defined, but its  Mellin transform
 dos not exist, or it exists, but the function itself is not restored by means of the inverse
 Mellin transformation (via the Mellin--Barnes
 integral).

If $A_j=B_k=C,$ $C>0,$ for $j =1,\ldots,p$ and $k=1,\ldots,q$ (a
special collection of parameters), then the Fox
 $H$-function is reduced to the Meyer $G$-function:
$$
H_{p,q}^{\,m,n} \!\left[ z \left| \begin{matrix}
(a_1, C) & (a_2, C) & \ldots & (a_p, C) \\
(b_1, C) & (b_2, C) & \ldots & (b_q, C) \end{matrix} \right.
\right] = \frac{1}{C} G_{p,q}^{\,m,n} \!\left( \left.
\begin{matrix} a_1, \dots, a_p \\
b_1, \dots, b_q \end{matrix} \; \right| \, z^{1/C} \right).
$$
  In the general case, the Fox
 $H$-function cannot be represented by a finite sum of generalized hypergeometric functions or
  $G$-functions (unlike the Meyer  $G$-function).
  If it is represented by a series, then it is called the\emph{Wright--Fox
  function}.
  Unlike the case of  Wright functions, no finite-order
  differential equation with polynomial coefficients is found such that the Fox function
  satisfies it (a conjecture that no such an equation exists is
  shared by the majority of experts, but is not proved yet).

In \cite{Koch1,Koch2, EiIvKoch}, Fox functions are applied to
construct the Green function for the fractional-diffusi\-on
 equation; the Stankovich transformation is  used for that.

Seemingly, Pincherle introduced integrals of the Mellin--Barnes
    type, used to define Meyer and Fox functions, for the first time (in 1888).
 Works of Barnes are referred to 1908--1910,
 while works of Mellin are referred to 1895 and 1909.

In our opinion, for functions introduced above, it is the most
correct (from the historical viewpoint) to use the general term
\emph{Wright--Fox  functions} and to specify the form they are
presented: either the Wright form (for the power series) or the
Fox form (for the complex integral).
 Nevertheless, the conventional terminology is as follows:
 functions defined by series  are called  Wright ones, functions defined by complex integrals
 are called  Meyer ones or Fox ones, and Fox functions represented by integrals of the Wright
  type are called the
   Wright--Fox ones.

Multivariable hypergeometric functions are used as well: they are
important for the theory of degenerating and singular differential
equations (see, e.\,g., \cite{VN, KarST2,KarST3, Erg1, CSh, HaKa1,
KarST, SaHa1, ShiE5}).

\subsection{Function spaces}\label{sec3.2}

Below, we provide notions and assertions (see \cite{KF,SKM,Bur1,
Bur2}) required to introduce transmutation operators and study
their properties.

As usual, $\mathbb{R}$ denotes the set of real numbers and
$\mathbb{C}$ denotes the set of complex numbers. Sets and
functions considered below are assumed to be measurable
 (unless otherwise stated).

\subsubsection{H\"older functions, absolutely continuous functions,
    and $AC^n$ class}\label{sec3.2.1}

In this section, we define the local and global H\"older
conditions and the function classes $AC$ and $AC^n$ (see
\cite{KF,SKM}).

Let $\Omega=[a,b],$ $-\infty{<}a{<}b{<}\infty$, denote a segment
of the real axis.

\begin{definition}
 We say that a function $f(x)$ satisfies the
    \textit{H\"older condition} of order $\lambda$ on  $\Omega$ if
    \begin{equation}\label{Ge}
    |f(x_1)-f(x_2)|\leq A|x_1-x_2|^\lambda
    \end{equation}
 for all $x_1,\,x_2\in\Omega,$
 where $A$ is a constant and $\lambda$ is called the  \textit{H\"older exponent}.
\end{definition}

\begin{definition}
 The set of all (in general, complex-valued)
 functions satisfying the H\"older condition of a fixed order $\lambda$ on $\Omega$
 is denoted by $H^\lambda=H^\lambda(\Omega).$
\end{definition}

If $\lambda> 1$, then the class $H^\lambda$ contains only
constants $f(x)\equiv \const$:
$$
|f'(x)|=\left|\lim\limits_{\Delta
x\to0}\frac{f(x_1)-f(x_2)}{x_1-x_2} \right| \le
A\lim\limits_{\Delta
    x\to0}\frac{|x_1-x_2|^\lambda}{|x_1-x_2|} | \le
A|x_1-x_2|^{\lambda-1}=0.
$$
    Therefore, the class $H^\lambda$ is interesting only for the case where $0<\lambda\leq
1.$

The class $H^1(\Omega)$ is called the \emph{Lipschitz class}.
    Define the class $AC(\Omega)$ of absolutely continuous functions,
    which is broader than $H^1(\Omega).$

\begin{definition}
  We say that a function $f(x)$ is   \textit{absolutely continuous} on a segment $\Omega$
  if for any positive $\varepsilon$ there exists a positive $\delta$ such that
  the inequality
     $$\sum\limits_{k=1}^n|f(b_k)-f(a_k)|<\varepsilon$$
  holds for any finite system of pairwise disjoint segments
   $[a_k,b_k]\subset\Omega,$ $k=1,2,\ldots,n,$ satisfying the
   inequality
   $$\sum\limits_{k=1}^n(b_k-a_k)<\delta.$$
 The class of all such functions is denoted by $AC(\Omega).$
\end{definition}

It is known (see \cite{KF}) that the class $AC(\Omega)$ coincides
with the class of antiderivatives of functions summable in the
Lebesgue sense, i.\,e., \begin{equation}\label{AC} f(x)\in
AC(\Omega)\Leftrightarrow
f(x)=c+\int\limits_a^x\varphi(t)dt,\,\,\,
\int\limits_a^b|\varphi(t)|dt<\infty,\,\,\varphi(t)=f'(t).
\end{equation}
    Therefore, each absolutely continuous function has an a.\,e.
    summable derivative $f'(x)$ (however, the existence of an a.\,e.
    summable derivative $f'(x)$ does not imply the absolute continuity).

\begin{definition}
  If $n=1,2,\ldots$ and $\Omega$ is a
    segment, then $AC^n(\Omega)$ is the class of functions $f(x)$ continuously differentiable
 on $\Omega$ up to order $n-1 $ and such that $f^{(n-1)}(x)\in AC(\Omega).$
\end{definition}

 It is obvious that $AC^1(\Omega)=AC(\Omega)$ and the class $AC^n(\Omega)$
  consists of functions representable by an $n$-multiple
 Lebesgue integral with a variable upper limit of a summable function,
  where the constant from \eqref{AC} is changed for the
  $(n-1)$-order polynomial:
\begin{equation}\label{ACn}
f(x)\in AC^n(\Omega)\Leftrightarrow
f(x)=\sum\limits_{k=0}^{n-1}c_k(x-a)^k+
    \underbrace{\int\limits_a^xdt\ldots\int\limits_a^{x}dt\int\limits_a^{x}}_n\varphi(t)dt,
\end{equation}
$$
\int\limits_a^b|\varphi(t)|dt<\infty,\,\,
c_k=\frac{f^{(k)}(a)}{k!},\,\,\varphi(t)=f^{(n)}(t).
$$
    Now, let $\Omega$ be an axis or semiaxis. In this case, defining the class
     $H^\lambda(\Omega)$, we add the requirement of the ``H\"older'' behavior at infinity.
 Namely, we say that a function $f(x)$ satisfies the H\"older condition in a neighborhood
 of infinity if
\begin{equation}\label{Geld}
\biggl|f\left(\frac{1}{x_1}\right)-f\left(\frac{1}{x_2}\right)\biggr|\leq
A\biggl|\frac{1}{x_1}-\frac{1}{x_2}\biggr|^\lambda
\end{equation}
 provided that the absolute values of $x_1$ and $x_2$ are sufficiently large.

\begin{definition}
 Let $\Omega$ be an axis or semiaxis. Then $H^\lambda=H^\lambda(\Omega)$ denotes
 the function class satisfying Condition \eqref{Ge} on each finite segment of $\Omega$
 and Condition \eqref{Geld} in a neighborhood of infinity.
\end{definition}

The union of Conditions \eqref{Ge} and \eqref{Geld} defining the
class  $H^\lambda(\Omega)$ for an infinite interval $\Omega$
    is equivalent to the unique condition
\begin{equation}\label{GlGe}
|f(x_1)-f(x_2)|\leq
A\frac{|x_1-x_2|^\lambda}{(1+|x_1|)^\lambda(1+|x_2|)^\lambda}
\end{equation}
 called the   {\it global H\"older condition}.

\subsubsection{$L_p$ class and its properties}\label{sec3.2.2}

Introduce the class of functions such that their $p$th degrees are
summable and provide several inequalities and theorems valid for
functions from this class (see \cite{KF,SKM}).
 Recently, generalizations of these spaces (introduced by F.~Riesz) to the case of variable
 exponents $p=p(x)$ are used for important applications (see,
e.\,g., \cite{KMRS,Rad2}).

Now, let $\Omega=[a,b],$ where $-\infty\leq a<b\leq\infty.$

\begin{definition}
 The set of all (in general, complex-valued)
  functions $f(x)$ measurable on  $\Omega$ and satisfying the inequality
  $$\int\limits_{\Omega}|f(x)|^pdx <    \infty,\quad 1 \leq p <\infty,$$
  is denoted by $L_p=L_p(\Omega)$.
\end{definition}

The norm in $L_p(\Omega)$ is defined by the relation
$$
||f||_{L_p(\Omega)}=\left(\int\limits_{\Omega}|f(x)|^pdx\right)^{1/p}.
$$
    If two elements of a space $L_p(\Omega)$ are functions
  different from each other on a zero-measure
  set, then they are treated as a same element of this space.

If $f\in L_p(\Omega)$ and $g\in L_p(\Omega)$, then the
\emph{Minkowski inequality}
$$
||f+g||_{L_p(\Omega)}\leq ||f||_{L_p(\Omega)}+||g||_{L_p(\Omega)}
$$
    holds; taking it into account, we see that   $L_p(\Omega)$ is a normed space.

If $f(x)\in L_p(\Omega),$ $g(x)\in L_{p'}(\Omega),$ and
$p'=\dfrac{p}{p-1},$ then the following \emph{H\"older inequality}
(more exactly,
 Rogers--H\"older--Riesz
 inequality) holds:
$$
\int\limits_{\Omega}|f(x)g(x)|dx\leq||f||_{L_p(\Omega)}||g||_{L_{p'}(\Omega)},\qquad
p'=\frac{p}{p-1}.
$$
    It is known that $L_p(\Omega)$ is a complete space.

Also, we need the following theorem allowing us to change the
order of the integrating in iterated integrals.

\begin{theorem*}[Fubini]
 Let $\Omega_1 =[a, b],$ $\Omega_2 = [c, d],$ $-\infty\leq a<b\leq\infty,$
  $-\infty\leq c<d\leq\infty,$ and
   $f(x, y)$ be a measurable function defined on $\Omega_1\times\Omega_2$.
   If at least one of the integrals
    $$
    \int\limits_{\Omega_1}dx\int\limits_{\Omega_2}f(x,y)dy,\qquad
    \int\limits_{\Omega_2}dy\int\limits_{\Omega_1}f(x,y)dx,\qquad
    \textrm{and}\qquad\iint\limits_{\Omega_1\times\Omega_2}f(x,y)dxdy
    $$
   {\rm (}absolutely{\rm )} converges, then they coincide each other.
\end{theorem*}

The following special case of the Fubini theorem holds:
\begin{equation}\label{Dirihle}
\int\limits_{a}^bdx\int\limits_{a}^xf(x,y)dy=\int\limits_{a}^bdy\int\limits_{y}^bf(x,y)dx
\end{equation}
    under the assumption that one of these integrals absolutely converges.
 Relation \eqref{Dirihle} is called the  \emph{Dirichlet formula}.

The generalized Minkowski inequality
$$
\left(\int\limits_{\Omega_1}dx\left|\int\limits_{\Omega_2}f(x,y)dy\right|^p\right)^{1/p}\leq
\int\limits_{\Omega_2}dy\left(\int\limits_{\Omega_1}|f(x,y)|^pdx\right)^{1/p}
$$
    holds as well.

\subsubsection{Spaces $C^m_{ev},$ $C^\infty_{ev},$ and $L_p^\gamma$}\label{sec3.2.3}

We say that  $\varphi$ is a \emph{test function} if it is
infinitely differentiable, even, and satisfies the estimates
$$
|D^q \varphi(x)|\leq\frac{C_{qr}}{(1+x^2)^r},
$$
 where $q$ and $r$ are arbitrary nonnegative integers, while
$C_{qr}$ are constants independent of $x.$ It is clear that if
$\varphi$ is a  test function, then $\dfrac{\varphi'(x)}{x}$ is a
 {test function} as well and, therefore, the Bessel operator is applicable to test functions
 arbitrarily many times and the estimates
$$
|B^q_\gamma \varphi(x)|\leq\frac{A_{qr}}{(1+x^2)^r}
$$
    hold for all arbitrary nonnegative integers $q$ and $r$.
 Denote the space of all test functions by  $S$.

Let $L_{2} (\mathbb{R}_{+}^1)$ denote the Hilbert space of
functions $f(y),$ $y>0,$ such that the norm
\begin{equation}
\| f \|_{L_{2} \lr{\mathbb{R}^1_+}} = \lr{\int\limits_0^{\infty}
|f(y)|^2 \, dy}^{\frac{1}{2}} \label{141}
\end{equation}
 is finite.

Let $L_{2, \nu} (\mathbb{R}_{+}^1),$ $\nu \geq - \dfrac{1}{2},$
 denote the weight Hilbert space of
functions  $f(y),$ $y>0,$  such that the norm
\begin{equation}
\| f \|_{L_{2, \nu} \lr{\mathbb{R}^1_+}} =
 \lr{\int\limits_0^{\infty} |f(y)|^2 y^{2 \nu +1} \, dy}^{\frac{1}{2}}
\label{1.4.1}
\end{equation}
is finite.

It is well known that the Fourier--Bessel
 transformation is unitary in $L_{2, \nu} (\mathbb{R}_{+}^1)$ and the
 Parseval identity
\begin{equation}
\|F_{\nu} f \|_{L_{2, \nu} \lr{\mathbb{R}^1_+}} =  \| f \|_{L_{2, \nu} \lr{\mathbb{R}^1_+}}.
\label{1.4.2}
\end{equation}
    holds.

The Kipriyanov function space $H^s_{\nu,+} (\mathbb{R}_{+}^1)$, $s
\geq 0,$ $\nu \geq - \dfrac{1}{2},$ introduced in \cite{Kip2} is
defined as the closure with respect to the norm
\begin{equation}\label{KipSp}
\| f \|_{H_{ \nu, +}^s \lr{\mathbb{R}^1_+}} =  \frac{1}{2^{\nu}\,
\Gamma (\nu +1)} \| (1+\eta^2)^{\frac{s}{2}} F_{\nu} f \|_{L_{2,
\nu} \lr{\mathbb{R}^1_+}}
\end{equation}
 of the set $\mathring{C}^{\infty}_{+} (\ov{\mathbb{R}_{+}^{1}})$ of functions even
in the Kipriyanov sense (see \cite{Kip2}). The evenness assumption
is essential because no finiteness of norm \eqref{KipSp} is
guaranteed otherwise. Other equivalent definitions of the norm are
possible for Kipriyanov  spaces as well.

The Sobolev space  $\mathring{H}^s (0, R),$ $s \geq 0,$
$0<R<\infty,$ is defined as the closure of the set
$\mathring{C}^{\infty} [0, R)$ with respect to the norm
$$
\| f \|_{\mathring{H}^s (0, R)} =  \|D^s  f \|_{L_{2, \nu} (0, R)}.
$$
    In the present monograph, we introduce a number of new functions (Banach and Fr\'echet)
     spaces  adjusted for the resolving of partial differential
     equations with singularities in coefficients of the equations
     and in the admitted solutions.

\subsection{Main integral transformations}\label{sec3.3}

Definitions and properties of the main integral transformations
can be found in \cite{BE1, BE2, BE3, PBM123, Ome, Ahi2, KGS, Pou,
Dzh1, KiSa, Marich1, KiSa, GSS}.

Unfortunate errors in inversions of  integral transformations,
occurred in \cite{BE1, BE2, BE3}, are corrected
in\cite{Kuz2,Kuz3}.

\subsubsection{Fourier transformation, sine transformation, cosine
transformation, and Hankel transformation}\label{sec3.3.1}

For the specified transformations, see \cite{Tit1, Dzh1, Pou}.
Important properties of  Hankel (Fourier--Bessel) transforms are
found in \cite{Lar1, Mur3, Mur4, Mur5, Mur6, Mur7}.

The \textit{Fourier transformation, sine transformation, cosine
transformation, and Hankel transformation} have the following
forms respectively:
\begin{eqnarray*}
    (Ff)(t)=\frac{1}{\sqrt{2\pi}}\int\limits_0^\infty \exp(-ity)f(y)\,dy,\\
    (F_c f)(t)=\sqrt{\frac{2}{\pi}}\int\limits_0^\infty \cos(ty)f(y)\,dy,\\
    (F_s f)(t)=\sqrt{\frac{2}{\pi}}\int\limits_0^\infty \sin(ty)f(y)\,dy.
\end{eqnarray*}
    Thy are unitary in $L_2(0,\infty)$ and coincide with their inverse transformations.

The Hankel (Fourier--Bessel)
 transformation has the form
$$
(H_\nu f)(t)=\frac{1}{t^\nu}\int\limits_{0}^{\infty}  J_\nu(ty)f(y)\,dy
$$
    or
$$(H_\nu f )(\xi)=\int\limits_{0}^{\infty} {j}_{\frac{\nu-1}{2}} (x\xi)\,
f(x)x^\nu\,dx.
$$
    Due to the property that ${j}_{\nu}(0)=1$, it is more convenient to use the
 Hankel transformation with the function ${j}_{\frac{\nu-1}{2}} (x\xi)$ in the kernel.

\subsubsection{Mellin transformation. Slater theorem}\label{sec3.3.2}

To compute the integral of a product of hypergeometric functions,
we use the method from \cite{Marich1, PBM}, based on the  Mellin
transformation. The Mellin transformation and Slater theorem are
studied in \cite{Marich1, PBM, PBM123}.

The  \emph{Mellin transformation} $f(x)$ is the function $g(s)$
    defined by the relation
\begin{equation}
\label{1710}
g(s)=M{f}(s)=\int\limits_0^\infty x^{s-1} f(x)\,dx.
\end{equation}
 Define the  Mellin convolution as follows:
\begin{equation}
\label{1711} (f_1*f_2)(x)=\int\limits_0^\infty
f_1\left(\frac{x}{y}\right) f_2(y)\,\frac{dy}{y}.
\end{equation}
  In  Mellin images, the convolution operator with kernel $K$ acts
  as a multiplier:
\begin{eqnarray}
\label{1712}
M{Af}(s)=\int\limits_0^\infty  K\left(\frac{x}{y}\right) f(y)\,\frac{dy}{y}=M{K*f}(s)
=m_A(s)M{f}(s),\\
    \nonumber m_A(s)=M{K}(s).\phantom{1111111111111111111}
\end{eqnarray}
    A convenient  algebraic approach for the investigation of operators of type
    \eqref{1712} is proposed in \cite{S6, S66}.
    It allows one to obtain needed estimates rapidly.
    Useful facts are gathered into the following assertion.

\begin{theorem}  \label{1tMel}
Let the convolution operator $A$ acts according to relation
\eqref{1712} as a multiplier in Mellin images.
    Then the following boundedness conditions for the direct and inverse operators
    and relations for their norms are satisfied{\rm:}
    \begin{enumerate}
    \item[(a)]
    it admits an extension till an
    operator bounded in $L_2(0,\infty)$ if and only if
    \begin{equation}
    \label{1713}
    \sup\limits_{\xi\in\mathbb{R}} \left|m_A\left( i\xi+\frac{1}{2}\right) \right|=M_2<\infty,
    \end{equation}
    where $\|A\|_{L_2}=M_2.$

    \item[(b)]
    under the additional condition that the kernel $K$ is nonnegative,
   it admits an extension till an
    operator bounded in  $L_p(0,\infty),$ ${p>1},$  if and only if
    \begin{equation}
    \label{1714}
    \sup\limits_{\xi\in\mathbb{R}}\left|m_A\left( i\xi+\frac{1}{p}\right) \right|
    =M_p<\infty,
    \end{equation}
    where $\|A\|_{L_p}=M_p.$

    \item[(c)]
    The inverse operator $A^{-1}$ acts according to relation
    \eqref{1712} as well and its multiplier is $\dfrac{1}{m_A}.$
     It admits an extension till an
    operator bounded in $L_2(0,\infty)$ if and only if
    \begin{equation}
    \label{1715}
    \inf\limits_{\xi\in\mathbb{R}}\left|m_A\left( i\xi+\frac{1}{2}\right)\right|=m_2>0,
    \end{equation}
    where $\|A^{-1}\|_{L_2}=\dfrac{1}{m_2}.$

    \item[(d)]
    Let operators  $A$ and $A^{-1}$ be defined and bounded in  $L_2(0,\infty).$
    Then they are unitary if and only if the relation
    \begin{equation}
    \label{1716}
    \left|m_A\left( i\xi+\frac{1}{2}\right) \right|=1
    \end{equation}
 holds for almost all $\xi.$
    \end{enumerate}
\end{theorem}

The last theorem summarizes results of many mathematicians: Schur,
Hardy, Littlewood, Polya, Kober, Mikhlin, and H\"ormander. It is
not known whether it is possible to omit the kernel nonnegativity
requirement in (b) under additional assumptions. In the general
case, there are no estimates in the range $\{0<p<1\}$ (an example
for Hardy operators is demonstrated by Burenkov).
  Below, we will see that Hardy operators are closely related to
 Buschman--Erd\'elyi transmutation operators.
 As far as we are aware, in \cite{Kob1}, the Mellin
 transformation technique is applied (for the first time) for norm estimates for
  Riemann--Liouville
 operators with pure imaginary degrees.
    This is why Part (b) of the above theorem is sometimes  called the Kober lemma.
    This is not completely exact because, actually, he proved the relation from Part (a)
 for the case of functions with alternating signs.

The  Mellin transformation is a generalized semiaxis Fourier
transformation with respect to the Haar measure $\dfrac{dy}{y}$
(see \cite{Hel1}). It is important for the theory of special
functions. For example, the gam\-ma-fun\-c\-ti\-on
   is the  Mellin transform of the exponential function.
   The following important achievement of 1970s is related to the Mellin transformation:
 the known Slater theorem allowing one to explicitly restore the majority of Mellin originals
    via their images is completely proved and adopted for the
    computing of integrals (a simple algorithm of hypergeometric functions is used for that,
    see \cite{Marich1, Sla, PBM}).
  Formally, this result can be easily obtained by means of the general Mellin--Barnes
   inversion relation (via residues). However, quite complicated and careful work is required
   to justify it strictly; one has to investigate the asymptotic behavior of hypergeometric
   functions near residues and at infinity, and a substantial diversity is typical for
   such a behavior.
   This work is only started by Slater, while it is completed by Marichev. The
 Slater--Marichev theorem yields a universal and powerful method to compute
 integrals. This method allows one to solve various problems from
 the theory of partial differential equations; also, it is the base
 of the forefront technologies of the symbol integrating of the \emph{MATHEMATICA}
 package.

The Slater theorem reads as follows (see \cite{Sla, Marich1,
PBM}). Let
\begin{equation}\label{4.13}
\Gamma\left[\begin{array}{cccc}
a_1, & a_2, &\ldots, & a_A \\
b_1, & b_2, &\ldots & b_B
\end{array}
\right]=\Gamma[(a),(b)]=\frac{\Gamma(a_1)\Gamma(a_2)\ldots\Gamma(a_A)}{\Gamma(b_1)\Gamma(b_2)
 \ldots\Gamma(b_B)},
\end{equation}
where each empty product is assigned to be equal to one,
\begin{equation*}
(a)+s=a_1+s,a_2+s,\ldots,a_A+s,
\end{equation*}
\begin{equation}\label{4.14}
(b)'-b_k=b_1-b_k,\ldots,b_{k-1}-b_k,b_{k+1}-b_k,\ldots,b_B-b_k,
\end{equation}
\begin{multline}\label{4.15}
\arraycolsep=2.5pt \Sigma_A(z)=\sum\limits_{j=1}^A
z^{a_j}\Gamma\left[\begin{array}{cc}
(a)'-a_j, & (b)+a_j \\
(c)-a_j & (d)+a_j
\end{array}
\right] \\
\times_{B+C}F_{A+D-1}\!\left(\begin{array}{ccc}
(b)+a_j, & 1+a_j-(c); & (-1)^{C-A}z \\
1+a_j-(a)', & (d)+a_j &
\end{array}
\right),
\end{multline}
\begin{multline}\label{4.16}
\arraycolsep=2pt \Sigma_B(1/z)=\sum\limits_{k=1}^B
z^{-b_k}\Gamma\left[\begin{array}{cc}
(b)'-b_k, & (a)+b_k \\
(d)-b_k & (c)+b_k
\end{array}
\right] \\
 \times_{A+D}F_{B+C-1}\!\left(\begin{array}{ccc}
(a)+b_k, & 1+a_k-(d); & \dfrac{(-1)^{D-B}}{z} \\
1+b_k-(b)', & (c)+b_k &
\end{array}
\right),
\end{multline}
    and
$$
|\arg{z}|<\pi.
$$
    If the series converge, then $\Sigma_A(z)$ and $\Sigma_B(1/z)$ are functions
    of the hypergeometric type and they pass into each other if we
    change places of the
 $A$-dimensional complex vector $(a)=a_1,a_2,\ldots,a_A$ and the analogous $B$-dimensional
 vector $(b),$   change places of the
$C$-dimensional  vector  $(c)$ and the $D$-dimensional
 vector $(d),$ and replace $z$ for $1/z.$
 These functions analytically depend on the complex parameters $(a),$ $(b),$ $(c),$ and $(d)$
 and on the variable $z.$ If the vectors $(a)$ (or $(b)$) have coinciding parameters
 or parameters such that the difference between them are integers,
 then the vectors $(a)'-a_j$ ($(b)'-b_k$) contain nonzero or negative integer
 components and, due to the property  $\Gamma(-n)=\infty,$ $n=0,1,2,\ldots,$
 the functions $\Sigma_A(z)$ and $\Sigma_B(1/z)$ might have ambiguous points of the kind
  $\infty-\infty.$ In such logarithmic cases,
  the values of $\Sigma_A(z)$ and $\Sigma_B(1/z)$ are understood as the corresponding limits
  of the  ``regular'' functions
$\Sigma_A(z)$  and $\Sigma_B(1/z)$ as their parameters
continuously tend to the considered special values.

Note that if the restriction $|\arg z|<\pi$ is not imposed, then
   $\Sigma_A(z)$ and $\Sigma_B(1/z)$ are not guaranteed to be one-valued
   functions.

\begin{theorem*}[Slater]
 Let
\begin{equation}\label{4.17}
K^*(s)=\Gamma\left[\begin{array}{cccc}
(a)+s & (b)-s \\
(c)+s, & (d)-s
\end{array}\right],
\end{equation}
 where $(a),$ $(b),$ $(c),$ and $(d)$ have $A,$
$B,$ $C,$ and $D$ components $a_j,$ $b_k,$ $c_l,$ and $d_m$
respectively.
 Let the following two groups of conditions be satisfied{\rm :}
\begin{equation}\label{4.181}
-\Re a_j<\Re s< \Re b_k,\,\, j=1,2,\ldots,A,\,\, k=1,2,\ldots,B,
\end{equation}
    and
\begin{equation}\label{4.182}
\left\{
\begin{array}{ll}
$$ A+B>C+D,$$ &  \\
$$ A+B=C+D,$$ & \hbox{$\Re s(A+D-B-C)<-\Re \nu,$} \\
$$ A=C,\, B=D,$$ & \hbox{$\Re \nu<0,$}
\end{array}
\right.
\end{equation}
 where
$$\nu=\sum\limits_{j=1}^Aa_j+\sum\limits_{k=1}^Bb_k-
\sum\limits_{l=1}^Cc_l-\sum\limits_{m=1}^Dd_m.
$$
 Then the following relations hold for the said values of $s${\rm:}
\begin{equation}\label{4.19}
K^*(s)=\left\{%
\begin{array}{ll}
$$\int\limits_0^\infty x^{s-1}\Sigma_A(x)dx,$$ & \hbox{$A+D>B+C,$} \\
$$\int\limits_0^1 x^{s-1} \Sigma_A(x)dx+\int\limits_1^\infty x^{s-1}\Sigma_B(1/x) dx,$$
 & \hbox{$A+D=B+C,$} \\
$$\int\limits_0^\infty x^{s-1} \Sigma_B(1/x)dx,$$ & \hbox{$A+D<B+C,$}, \\
\end{array}%
\right.
\end{equation}
    and
$\Sigma_A(1)=\Sigma_B(1)$ provided that $A+D=B+C,$ $\Re
\nu+C-A+1<0,$ and $A\geq C.$
\end{theorem*}

\begin{corollary*}
    If Conditions \eqref{4.181}-\eqref{4.182}
    are satisfied, then the preimage of the function
$$
K^*(s)=\Gamma\left[\begin{array}{cccc}
(a)+s & (b)-s \\
(c)+s, & (d)-s
\end{array}\right],
$$
    is a function $K(x)$ of the hypergeometric type such that it is equal to one of the
  following expressions{\rm:}
\begin{equation}\label{4.20}
K(x)=\left\{%
\begin{array}{ll}
$$\Sigma_A(x),$$ & \hbox{$x>0,A+D>B+C,$} \\
$$\Sigma_A(x),$$ & \hbox{$0<x<1,A+D=B+C,$} \\
$$\Sigma_B(1/x),$$ & \hbox{$x>1,A+D=B+C,$} \\
$$\Sigma_B(1/x),$$ & \hbox{$x>0,A+D<B+C,$} \\
\end{array}%
\right.
\end{equation}
 where $K(1)=\Sigma_A(1)-\Sigma_B(1),$ provided that $A+D=B+C,$ $\Re \nu+C-A+1<0,$ and
 $A\geq C.$
\end{corollary*}

%   \pagebreak

\begin{remark}
\label{r1}
 The corresponding Slater theorem (see \cite[Sec. 4.8, Th. 1]{Sla})
 contains the following  inaccuracies:
\begin{enumerate}
\item[(1)]
 only the condition $A+B\geq C+D$ is imposed instead of the
 conditions
 $$
\left\{
\begin{array}{ll}
$$A+B>C+D,$$ &  \\
$$A+B=C+D,$$ & \hbox{$\Re s(A+D-B-C)<-\Re \nu,$} \\
$$A=C,\, B=D,$$ & \hbox{$\Re \nu<0;$}
\end{array}
\right.
$$

\item[(2)]
 the condition $\Re \nu<0$  is imposed instead of the
 condition $\Re \nu+C-A+1<0$;

\item[(3)]
 it is stated that if  $A+D=B+C$, then the functions $\Sigma_A(z)$ and $\Sigma_B(1/z)$
 analytically continue each other (actually, this is true only under the assumption that
  $A+B>C+D.$)
\end{enumerate}
\end{remark}

\begin{remark}
\label{r2} If $|A+D-B-C|>1$ and $A+B=C+D$, then the restriction
imposed by \eqref{4.182} $\Re s,$ can be weakened till the
following condition:
\begin{equation}\label{4.183}
\Re s(A+D-B-C)<\frac{1}{2}-\Re \nu.
\end{equation}
\end{remark}

\begin{remark}
\label{r3}
    If there are parameters such that one or more conditions $a_j=c_l+n$  (or $a_j=-d_m-n$)
     [$b_k=d_m+n$ (or $b_k=-c_l-n$)], where $n=0,1, 2,\ldots,$ are satisfied
     and the vectors $(a)'-a_j$ and $(b)'-b_k$ contain no integer components,
     then the conditions referring to these parameters can be excluded from Conditions
     \eqref{4.181} or can be weakened.
 For $a_j=c_l+n$  [$b_k=d_m+n$], the corresponding conditions $\Re (s+a_j)>0$
  [$\Re (b_k-s)>0$] are removed; for $(a)'-a_j$ and $(b)'-b_k$, they are changed for the
  weakened requirements  $\Re (s+a_j)>-n-1$  [$\Re (b_k-s)>-n-1$].
If the vectors $(a)'-a_j$ and $(b)'-b_k$  contain integer
components, then the possibility to weaken restrictions
\eqref{4.181} is to be investigated separately.
\end{remark}

\begin{remark}
\label{r4}
 Restrictions \eqref{4.181}-\eqref{4.182}
 provide at least the conditional convergence of integrals \eqref{4.19}
 at $0$ and $\infty$ (and $1$). If these restrictions are broken,
 then, in general, the improper integrals \eqref{4.19} diverge;
 however, there are particular cases where they exist in principal value sense.
\end{remark}

\subsubsection{Various forms of fractional integrodifferentiating}
\label{sec3.3.3}

Operators of  fractional integrodifferentiating are studied,
e.\,g., in \cite{SKM, Nah1, Nah2, Nah3, KST, ViRy, GoMa, KT1,KT2}.
Useful properties of  fractional integrals and their inversions
for measures are provided in \cite{KaPr4} (see also
\cite{LShFrac}).

Operators of  fractional integrodifferentiating are important for
various areas of contemporary mathematics. For the theory of
special functions, the importance of the fractional
integrodifferentiating is reflected in the title of \cite{Kir4}
(according to Prof. Kilbas, the only exception is Fox functions).

   Consider fractional
Riemann--Liouville integrals and derivatives.

\begin{definition}
 Let $\varphi (x)\in L_{1} (a,b).$ Then the integrals
    \begin{equation} \label{RLI1}
    (I_{a+}^{\alpha } \varphi )(x)=\frac{1}{\Gamma (\alpha )} \int\limits
    _{a}^{x}\frac{\varphi(t)}{(x-t)^{1-\alpha } } dt,\qquad x>a,
    \end{equation}
    and
    \begin{equation} \label{RLI2}
    (I_{b-}^{\alpha } \varphi )(x)=\frac{1}{\Gamma (\alpha )} \int\limits
    _{x}^{b}\frac{\varphi(t)}{(t-x)^{1-\alpha } } dt,\qquad x<b,
    \end{equation}
   where $\alpha >0,$ are called the
   {\it left-side} and  {\it right-side}
   (respectively)
        {\it fractional
Riemann--Liouville integrals} of order $\alpha$.

 For a function $f(x),$ $x\in[a,b],$ the expressions
    \begin{equation} \label{RLD3} (D_{a+}^{\alpha }f)(x)=
    \frac{1}{\Gamma (n-\alpha )} \left(\frac{d}{dx} \right)^{n} \int\limits _{a}^{x}
    \frac{f(t)dt}{(x-t)^{\alpha -n+1} },   \end{equation}
    and
    \begin{equation} \label{RLD4} (D_{b-}^{\alpha }f)(x)=\frac{1}{\Gamma (n-\alpha )}
     \left(\frac{d}{dx} \right)^{n}
     \int\limits _{x}^{b}\frac{f(t)dt}{(t-x)^{\alpha -n+1} },  \end{equation}
 where $n=[\alpha ]+1,$ ${\alpha >0},$ are called the
  {\it left-side} and  {\it right-side}
   (respectively) {\it fractional
    Rie\-mann--Li\-o\-u\-ville derivative} of order $\alpha.$
\end{definition}

In special cases, the Riemann--Liouville
 operators most important for applications are defined (for positive values of $\alpha$)
  as follows:
\begin{equation}
\label{161} I_{0+,x}^{\alpha}f=\frac{1}{\Gamma(\alpha)}
\int\limits_0^x \left( x-t\right)^{\alpha-1}f(t)d\,t
 \quad\textrm{and} \quad
\nonumber
I_{-,x}^{\alpha}f=\frac{1}{\Gamma(\alpha)}\int\limits_x^\infty
\left( t-x\right)^{\alpha-1}f(t)d\,t.
\end{equation}
    For other  values of  $\alpha$, they are defined by means of the analytic continuation
      (regularization).

Note that there are numerous variants of fractional
integrodifferentiating operators:
 Gerasimov--Caputo ones introduced by Gerasimov in 1948 (see \cite{Nov}) and Caputo in 1968,
  Marchaud ones, Weyl ones, Riesz ones,
  Erd\'elyi--Kober ones, Hadamard ones, Gel'fond--Leontiev
  ones  (the double title is proposed by Korobeinik), Dzhrbashyan--Nersesyan
  ones, etc.
    Euler-type differential equations with main types of fractional
integrodifferentiating operators are considered, e.\,g.,
 in \cite{KiZhu,ZhuSi1}.

    Now, consider Erd\'elyi--Kober
fractional integrals and derivatives.

For $\alpha > 0$, \textit{Erd\'elyi--Kober operators} are defined
by the following relations:
\begin{eqnarray}
\label{162}
& & I_{0+;\, 2,\, y}^{\alpha} f = \frac{2}{\Gamma(\alpha)}x^{-2(\alpha+y)}
\int\limits_0^x (x^2-t^2)^{\alpha-1}t^{2y+1}f(t)\,dt, \\
& & I_{-;\, 2,\, y}^{\alpha} f = \frac{2}{\Gamma(\alpha)}x^{2y}
\int\limits_x^{\infty}
(t^2-x^2)^{\alpha-1}t^{2(1-\alpha-y)-1}f(t)\,dt.
\end{eqnarray}
    For $\alpha > -n,$ $n \in \mathbb{N},$ they  are defined
by the following relations:
\begin{eqnarray}
& & I_{0+;\, 2, y}^{\alpha} f =x^{-2(\alpha+y)} {\lr{\frac{d}{d x^2}}}^n  x^{2(\alpha
    +y+n)}I^{\alpha+n}_{0+; \, 2,\, y}f, \label{1.15} \\
& & I_{-;\, 2, y}^{\alpha} f = x^{2y} {\lr{-\frac{d}{d x^2}}}^n
x^{2(\alpha
    -y)}I^{\alpha+n}_{-; \, 2,\, y-n}f \label{1.16}.
\end{eqnarray}
    For other values of $\alpha$, they are defined by means of the analytic continuation
 (similarly to fractional Liouville integrodifferentiating operators).

In \cite{SKM}, the case of the integration limits $0$ and $\infty$
are not considered (unlike the present monograph). In
\cite{SaKiMar}, which is the next English version, these special
cases of the limits are allowed, but the definitions contain
inaccuracies; in particular, they lead to complex values under in
the integrand.

Now, consider the most general form of such operators, i.\,e.,
fractional integrals with respect to arbitrary functions.

The \textit{fractional integral  with respect to an arbitrary
function $g(x)$} is defined as follows:
\begin{gather}
\label{163}
I_{0+,g}^{\alpha}f=\frac{1}{\Gamma(\alpha)}
\int\limits_0^x \left( g(x)-g(t)\right)^{\alpha-1}g'(t)f(t)d\,t,\\
    %   \nonumber
I_{-,g}^{\alpha}f=\frac{1}{\Gamma(\alpha)}\int\limits_x^\infty
\left( g(t)-g(x)\right)^{\alpha-1}g'(t)f(t)d\,t,
\end{gather}
     where $\Re\alpha>0.$
     For the remaining values of $\alpha$, the relations are easily extended
      (see \cite{SKM}).
   Fractional integrals \eqref{161} are obtained if we assign $g(x)=x$ in \eqref{163}, the
    Erd\'elyi--Kober integrals \eqref{162} are obtained if we assign
$g(x)=x^2,$ and the Hadamard integrals are obtained if we assign
$g(x)=\ln x.$

The relation to transmutation operators is as follows: up to
factors, the
 Sonin--Poisson--Delsarte
  operators are the Erd\'elyi--Kober
  operators, i.\,e., the fractional powers
$$\left(\frac{d}{dg(x)}\right)^{-\alpha}=\left(\frac{d}{2xdx}\right)^{-\alpha},\quad g(x)
 =x^2.$$
    Therefore, the main properties of these transmutation operators can be obtained
    from the theory of fractional integrodifferentiating operators instead of to invent
    them again (as it is frequently done). A.~Dzhrbashyan
    attracted our attention to the fact that fractional
    integrating operators with respect to function \eqref{163} are special cases of more
    general operators introduced and studied by M.~Dzhrbashyan (see \cite{SKM}).

\subsubsection{Quadratic {\rm (}fractional{\rm )}
    Fourier and Hankel transformations}\label{sec3.3.4}

Briefly presenting main facts from the theory of the quadratic (fractional)
 Fourier--Fresnel transformation, we follow \cite{OZK,Kar1}, where
 more general information and references can be found.

Integer powers (the orbit) of the classical Fourier transformation
form a cyclic group of order four such that the fourth degree  of
this transformation yields the identical operator. Therefore, the
spectrum  of the classical Fourier transformation in
$L_2(-\infty,\infty)$
 consists of the following four points located on the unit circle: $1,$ $i,$ $-1,$
 and $-i.$ The idea to include this discrete group into a continuous one such that its
 spectrum completely fills the unit circle belongs to Wiener (in 1929,
  it is implemented by him).
  Nowadays, this group is called the fractional Fourier
transformation. For the multidimensional case, it is generalized
by Bargmann.

The fractional Fourier transformation is rediscovered by a number
of author many times.
%   В книге Антосика, Микусинского и
%Сикорского  под названием ``преобразование Fourier--Мелера''
%упоминается циклическая группа произвольного порядка, в которую
%можно включить преобразование Fourier. fractional Fourier
%transformation изучалось Гинандом  и Вольфом. Вавржинчик в
%приходит к fractional Fourier transformation рассматривая
%классическое преобразование Fourier в виде экспоненты от
%производящего оператора. В.\,Ф.~Осиповым независимо от предыдущих
%авторов построена теория fractional Fourier transformation на
%группах и введены соответствующие почти-периодические функции
%Бора--Френеля, изучены асимптотические свойства этого
%преобразования, рассмотрены приложения в гармоническом анализе и
%теории чисел~\cite{Os, AbOs}. Намиас переоткрыл дробное
%преобразования Fourier и использовал его для решения некоторых
%задач для уравнения Шр\"{e}\-дин\-гера. Керр исследовала дробное
%преобразование Fourier в пространстве $L_2$ и пространстве Шварца
%$S.$

A separate research direction is related to the fractional Fourier
transformation for purely imaginary values of the group parameter.
    In this case, the most propagated term is the \emph{Hermit semigroup}.
    This direction is originated by the Hille paper of  1926,
    where the operator of a fractional Fourier transformation with imaginary values
    of parameter arising in relation with the Abel summation of expansions with respect to
    Hermit polynomials. Later,  the Hermit semigroup is used by Babenko,
  Bekner, and Weisler to obtain inequalities in the theory of classic Fourier
  transformation.

Note that the quadratic Fourier transformation is one of the two
main components (the second one is
  inequalities for mean values in a special-kind
  complex plane) used by  Babenko and
  Bekner to prove the famous boundedness conditions (with exact constants) of the classical
   Fourier transformations in spaces $L_p$.
   Another interesting application of the quadratic Fourier transformation is
 related to the famous Pauli problem to determine the function via a spectral data collection.
 problems of this kind are usually  unsolvable for the classical transformation,
 but their elegant solutions are  found in the framework of the
 theory of  quadratic Fourier--Fresnel
  transformation.

Applications of the fractional Fourier transformation are numerous
as well.
 We mention only several ones. Applications in quantum mechanics are studied by
 Namias. Applications in optics and signal analysis are developed in \cite{OZK}.
 Analogs of the Heisenberg inequality,
 invariant with respect to the fractional Fourier transformation, are found by Mastard.
 It is shown that the multidimensional Wigner transformation is equal to the sixth-order
 root from the inverse Fourier transform; hence, it is a special case of the fractional
  Fourier transformation as well. New inversion relations for the fractional
  Fourier transformation are obtained by Bjun.  An analog of the fractional
   Fourier transformation for Hermit $q$-polynomials
   is introduced by Aski, Atakishev, and Suslov.
   Further generalizations to operators such that their kernels are bilinear generating
   functions of Aski--Wilson
   polynomials are considered in works of Aski and Rakhman.

The fractional Hankel transformation is much less studied.
    It is introduced by Kober and and studied by Ginand.
    Then it was rediscovered several times, e.\,g., by Namias.
 For real values of the group parameter, the fractional Hankel transformation is considered
  by Kerr in the space $L_2(0,\infty)$ and in Zemanyan spaces (see \cite{S60} as well).

In \cite{Kar1}, rather general approach to the constructing of
similar transformations by means of series expansions with respect
to known systems of orthogonal functions is considered.
 In particular, the classical Fourier transformation is obtained if the system of Hermit
 functions is selected, the quadratic Fo\-urier--Fresnel
  transformation  is obtained if the system of Laguerre functions is
  selected, and new semigroups of integral transformations are
  constructed if systems of Legendre, Chebyshev, or Gegenbauer
  functions  are selected.

    Following \cite{Kar1}, provide explicit integral relations for    the quadratic
    Fo\-urier--Fresnel transformation and quadratic Hankel transformation:
\begin{equation}\label{FrF}
(F^{\alpha}f)(y)=\frac{1}{\sqrt{\pi(1-e^{2i\alpha})}}\int\limits_{-\infty}^{\infty}\!\!
e^{-\frac{1}{2}i(x^2+y^2)\cot\alpha}e^{ixy\csc\alpha}f(x)dx
\end{equation}
    and
\begin{equation}\label{FrH}
(H_\nu^{\alpha}f)(y)=\frac{2\left(-e^{i\alpha}\right)^{-\frac{\nu}{2}}}{1-e^{i\alpha}}
\int\limits_0^{\infty}\!\!e^{-\frac{1}{2}i\cot\frac{\alpha}{2}\left(x^{2}+y^{2}\right)}
(xy)^{\frac{1}{2}}J_{\nu}\left(\frac{2xy\sqrt{-e^{i\alpha}}}{1-e^{i\alpha}}\right)f(x)dx.
\end{equation}
    In brief, consider relations of transformations forming the
    operation calculus for the quadratic Hankel transformation.

Introduce the differential operators
\begin{gather*}
A^{-}_\nu=x^{\nu +\frac{1}{2}}e^{-\frac{x^2}{2}}\frac{d}{dx}
x^{-\nu -\frac{1}{2}}e^{\frac{x^2}{2}}=-\frac{\nu +
\frac{1}{2}}{x}+x+\frac{d}{dx},
\\
A^{+}_\nu=x^{-\nu -\frac{1}{2}}e^{\frac{x^2}{2}}
\frac{d}{dx}x^{\nu
+\frac{1}{2}}e^{-\frac{x^2}{2}}=\frac{\nu+\frac{1}{2}}{x}-x+\frac{d}{dx},
\\
N_\nu=x^{\nu+\frac{1}{2}}\frac{d}{dx}x^{-\nu-\frac{1}{2}}=-\frac{\nu+\frac{1}{2}}{x}
+\frac{d}{dx},
\end{gather*}
    and
    $$
M_\nu=x^{-\nu-\frac{1}{2}}\frac{d}{dx}x^{\nu+\frac{1}{2}}=\frac{\nu+\frac{1}{2}}{x}
+\frac{d}{dx}.
    $$
    The following relations hold for these operators:
\begin{equation}\label{AMN}
L_{\nu}= -\frac{1}{4}D^2-\frac{\nu^2-1/4}{x^2} + \frac{1}{4}x^2 -
\frac{\nu+1}{2}
=-\frac{1}{4}A^{+}_{\nu}A^{-}_{\nu},~~~~A^{-}_{\nu}=N_{\nu}+x,
~~~~A^{+}_{\nu}=M_{\nu}-x,
\end{equation}
    and
\begin{equation}\label{LMN}
L_{\nu}=M_{\nu} N_{\nu},~~~~
x\frac{d}{dx}+\frac{d}{dx}x=N_{\nu}x+xM_{\nu}=M_{\nu}x+xN_{\nu}.
\end{equation}
  The multiplier by the independent variable is denoted by $X$.
Now, following \cite{Kar1}, one can represent the collection of
transmutation relations for operations for the quadratic Hankel
transformation:
\begin{align}\label{H+tD}
H_{\nu+1}^{\alpha}Xf&=\frac{1}{2}\left[(e^{-i\alpha}-1)N_{\nu}+
(e^{-i\alpha}+1)X\right]H_{\nu}^{\alpha}f,
\\
\label{H+ND}
H_{\nu+1}^{\alpha}N_{\nu}f&=\frac{1}{2}\left[(e^{-i\alpha}+1)N_{\nu}+
(e^{-i\alpha}-1)X\right]H_{\nu}^{\alpha}f,
\\
\label{NHD} N_{\nu}H_{\nu}^{\alpha}f&=
\frac{1}{2}H_{\nu+1}^{\alpha}\left[(e^{i\alpha}+1)N_{\nu}+(e^{i\alpha}-1)X\right]f,
\\
\label{xHD}
XH_{\nu}^{\alpha}f&=\frac{1}{2}H_{\nu+1}^{\alpha}\left[(e^{i\alpha}-1)N_{\nu}+
(e^{i\alpha}+1)X\right]f,
\\
\label{HtD} H_{\nu}^{\alpha}Xf&=
\frac{1}{2}\left[(1-e^{i\alpha})M_{\nu}+(1+e^{i\alpha})X\right]H_{\nu+1}^{\alpha}f,
\\
\label{HMD} H_{\nu}^{\alpha}M_{\nu}f&=
\frac{1}{2}\left[(1+e^{i\alpha})M_{\nu}+(1-e^{i\alpha})X\right]H_{\nu+1}^{\alpha}f,
\\
\label{MH+D} M_{\nu}H_{\nu+1}^{\alpha}f&=
\frac{1}{2}H_{\nu}^{\alpha}\left[(1+e^{-i\alpha})M_{\nu}+(1-e^{-i\alpha})X\right]f,
\\
\label{xH+D} XH_{\nu+1}^{\alpha}f&=
\frac{1}{2}H_{\nu}^{\alpha}\left[(1-e^{-i\alpha})M_{\nu}+(1+e^{-i\alpha})X\right]f,
\\
\label{A-H}
H_{\nu+1}^{\alpha}A^{-}_{\nu}f&=e^{-i\alpha}A^{-}_{\nu}H_{\nu}^{\alpha}f,
\qquad
A^{-}_{\nu}H_{\nu}^{\alpha}f=H_{\nu+1}^{\alpha}e^{i\alpha}A^{-}_{\nu}f,
\\
\label{A+H}
H_{\nu}^{\alpha}A^{+}_{\nu}f&=e^{i\alpha}A^{+}_{\nu}H_{\nu+1}^{\alpha}f,\qquad
A^{+}_{\nu}H_{\nu+1}^{\alpha}f=H_{\nu}^{\alpha}e^{-i\alpha}A^{+}_{\nu}f,
\\
\label{Ht^2D} H_{\nu}^{\alpha}X^{2}f&=
\left[X^{2}\cos^{2}{\frac{\alpha}{2}}-\frac{1}{2}i\sin{\alpha}\left[XD+DX\right]-
\sin^{2}{\frac{\alpha}{2}}L_{\nu}\right]H_{\nu}^{\alpha}f,
\\
\label{x^2HD} X^{2}H_{\nu}^{\alpha}f&=
H_{\nu}^{\alpha}\left[X^{2}\cos^{2}{\frac{\alpha}{2}}+\frac{1}{2}i\sin{\alpha}\left[XD+DX\right]-
\sin^{2}{\frac{\alpha}{2}}L_{\nu}\right]f,
\\
\label{HLD} H_{\nu}^{\alpha}L_{\nu}f&=
\left[-X^{2}\sin^{2}{\frac{\alpha}{2}}-\frac{1}{2}i\sin{\alpha}\left[XD+DX\right]+
\cos^{2}{\frac{\alpha}{2}}L_{\nu}\right]H_{\nu}^{\alpha}f,
\\
\label{LHD} L_{\nu}H_{\nu}^{\alpha}f&=
H_{\nu}^{\alpha}\left[-X^{2}\sin^{2}{\frac{\alpha}{2}}+\frac{1}{2}i\sin{\alpha}\left[XD+DX\right]+
\cos^{2}{\frac{\alpha}{2}}L_{\nu}\right]f,
\\
\label{HtD+DtD} H^\alpha_\nu\left[XD+DX\right]f&=
\left[-i\sin{\alpha}(X^2+L_\nu)+\cos\alpha\left[XD+DX\right]\right]H_{\nu}^{\alpha}f,
\end{align}
and
 \begin{equation}\label{xD+DxHD}
\left[XD+DX\right]H_{\nu}^{\alpha}f=H_{\nu}^{\alpha}\left[i\sin{\alpha}(X^2+L_\nu)+
\cos\alpha\left[XD+DX\right]\right]f.
    \end{equation}
 The above operation calculus relations for the quadratic
Hankel transformation allow one to use this integral
transformation in the composition method to construct
transmutation operators, explained in Chap. \ref{ch6}.

Thus, we have considered the main integral transformations used in
the present monograph. Regarding general theoretical aspects of
various classes of operators and function spaces,  see \cite{Nai,
SvFe, KF, KPS, Bas, Kus, Pas1, Sol, Fet}.

\subsubsection{Main classes of differential equations with Bessel
operators and related transmutation operators}\label{sec3.3.5}

To call differential equations with Bessel operators
\begin{equation}
\label{c1Bessel}
B_{\nu}u(x)=\frac{d^2 u}{dx^2}+\frac{2\nu+1}{x} \frac{du}{dx},
\end{equation}
 we use the following terminology (see \cite{Kip1}).

Equations with Bessel operators of the kind
\begin{equation}
\label{c1Bes1}
 \sum\limits_{k=1}^{n}B_{\nu,x_k}u(x_1,\dots,x_n)=f(t,x)
\end{equation}
    is called a {\it $B$-elliptic}
    equations.
  Sometimes, it is called the {\it
 Laplace--Bessel  equation}.

{\it $B$-hyperbolic} equations are equations with Bessel operators
of the kind
\begin{equation}
\label{c1Bes2} B_{\nu,t}u(t,x_1,\dots,
x_n)=\sum\limits_{k=1}^{n}B_{\nu,x_k}u(t,x_1,\dots, x_n) + f(t,x).
\end{equation}
 If the spatial variable is unique, then we obtain the {\it
Euler--Poisson--Darboux
    equation}.

{\it $B$-parabolic} equations are equations with Bessel operators
of the kind
\begin{equation}
\label{c1Bes3} \frac{\pr}{\pr t}u(t,x_1,\dots,
x_n)=\sum\limits_{k=1}^{n}B_{\nu,x_k}u(t,x_1,\dots, x_n) + f(t,x).
\end{equation}
  The same terms are preserved for incomplete equations,
  where one or several Bessel operators are reduced to the second
  derivatives and spectral parameters are added to the equations.
  The specified three classes of differential equations (according to the Kipriyanov
  classification) are considered in \cite{Kip1} ($B$-elliptic
  equations), \cite{CSh} ($B$-hyperbolic
  ones), and \cite{Mat1} ($B$-parabolic
  ones).
  We list several known classes of differential equations with Bessel
operators without any claims for this list to be complete; only
main references are provided. More complete information is
provided in \cite{SSfiz}, which is devoted to differential
equations with  Bessel operators.

$B$-elliptic equations of the kind
\begin{equation}\label{BEll}
\sum\limits_{i=1}^n\left( \frac{\partial^2 u}{\partial x_i^2}+
\frac{k_i}{x_i}\frac{\partial u}{\partial x_i}\right)=f(x)
\end{equation}
 form a type including the Weinstein equation of the
 generalized axially symmetric potential theory (see \cite{Wei1, Wei2,
Wei3}); equations of this class are profoundly and comprehensively
studied in works of Kipriyanov and his school (see \cite{Kip1}).

If the Bessel operator acts with respect to a spatial variable,
then we obtain the following singular variant of the wave equation
with an axially or central symmetry:
\begin{equation}\label{EPDWe}
\frac{\partial^2 u}{\partial t^2}= \frac{\partial^2 u}{\partial
x^2}+\frac{\nu}{x}\frac{\partial u}{\partial x}, \qquad
u=u(x,t),\qquad x>0,\qquad t\in\mathbb{R},\qquad \nu=\const.
\end{equation}
    Representations of solutions of Eq. \eqref{EPDWe} are obtained in \cite{Poisson}.
    The generalized radiation Weinstein problem (see \cite{CSh}) refers to this class as well.

If the Bessel operator acts with respect to  variable $t,$ then we
obtain the famous
  Euler--Pois\-son--Dar\-boux
  equation
\begin{equation}\label{EPD0}
\frac{\partial^2 u}{\partial t^2}+\frac{\nu}{t}\frac{\partial u}{\partial t}=
a^2\frac{\partial^2 u}{\partial x^2},\qquad u=u(x,t),\qquad t>0,\qquad x\in\mathbb{R},\qquad
a, \ \nu=\const.
\end{equation}
    For the first time, this equation appears in \cite[p. 227]{Euler}; then it
    is studied in \cite{Poisson,Riman,Darboux}.

For several spatial variables, the Cauchy problem
\begin{equation}\label{EPD1}
\frac{\partial^2 u}{\partial t^2}+\frac{\nu}{t}\frac{\partial u}{\partial t}=
\sum\limits_{i=1}^n\frac{\partial^2 u}{\partial x_i^2},\qquad u=u(x,t),\qquad t>0,\qquad x\in\mathbb{R}^n,\qquad
\nu=\const,
\end{equation}
\begin{equation}\label{UslEPD1}
u(x,0)=f(x),\qquad u_t(x,0)=0,
\end{equation}
    is considered in
and \cite{Diaz, Wei1, Wei2, Wei3} under the assumption that
$\nu>n-1$.

The generalized
 Euler--Poisson--Darboux
    equation
\begin{equation}\label{EPDWG}
\frac{\partial^2 u}{\partial t^2}+\frac{\nu}{t}\frac{\partial
u}{\partial t}= \frac{\partial^2 u}{\partial
x^2}+\frac{k}{x}\frac{\partial u}{\partial x},\qquad u=u(x,t),
\qquad t>0,\qquad x>0,\qquad \nu,k=\const,
\end{equation}
    and its
\begin{equation}\label{EPDM}
\frac{\partial^2 u}{\partial t^2}+\frac{\nu}{t}\frac{\partial u}{\partial t}=
\sum\limits_{i=1}^n\left( \frac{\partial^2 u}{\partial x_i^2}+\frac{k_i}{x_i}\frac{\partial u}{\partial x_i}\right),
\end{equation}
$$
 u=u(x_1,..,x_n,t),\qquad t>0,\qquad x_i>0,\qquad
 \nu,k_i=\const, \qquad i=1,..,n,
$$
 are studied in (\cite{CSh,Fox,LPSh1,LPSh2,ShiE3,Smi,76,Ter1}.
    The
     Euler--Poisson--Darboux
   equation  with the spectral potential
\begin{equation}\label{EPDSP}
\frac{\partial^2 u}{\partial t^2}+\frac{\nu}{t}\frac{\partial u}{\partial t}=
\frac{\partial^2 u}{\partial x^2}\pm\lambda^2 u,\qquad \lambda\in\mathbb{R}
\end{equation}
    is considered in \cite{Bresters2,Smi}.
     The generalized
     Euler--Poisson--Darboux
   equation  with the spectral potential
\begin{equation}\label{EPDSP1}
\frac{\partial^2 u}{\partial t^2}+\frac{\nu}{t}\frac{\partial u}{\partial t}=
\sum\limits_{i=1}^n\left( \frac{\partial^2 u}{\partial x_i^2}+\frac{k_i}{x_i}\frac{\partial u}{\partial x_i}\right) -\lambda^2 u,\qquad \lambda\in\mathbb{R}
\end{equation}
are studied in \cite{ShiE5}.

Further, generalized wave equations with variable potentials of
the kind
\begin{equation}\label{wpot}
\frac{\partial^2 u}{\partial t^2} = \frac{\partial^2 u}{\partial
x^2}+p(x)u\quad \textrm{and} \quad \frac{\partial^2 u}{\partial
t^2} +q(t)u= \frac{\partial^2 u}{\partial x^2}+p(x)u,
\end{equation}
 are generalized to $B$-hyperbolic
 equations with variable potentials of
the kind
\begin{equation}\label{Bpot1}
B_{\nu,t}u = B_{\nu,x}u + p(x)u\quad \textrm{and} \quad B_{\nu,t}u
+q(t)u= B_{\nu,x}u + p(x)u,
\end{equation}
 their multidimensional analogs
\begin{equation}\label{Bpot2}
B_{\nu,t}u +q(t)u= \sum\limits_{k=1}^n B_{\nu,x_k}u + p(x_k)u,
\end{equation}
 and the following $B$-ultrahyperbolic
  equations with variable potentials:
\begin{equation}\label{Bpot3}
\sum\limits_{j=1}^m \left(B_{\nu,t_j}u +q(t_j)u\right)=
\sum\limits_{k=1}^n \left(B_{\nu,x_k}u + p(x_k)u\right).
\end{equation}
    Usually, these equations are solved (explicitly or implicitly) by means
    of transmutation operators for Sturm--Liouville
    operators or perturbed Bessel operators.

Important class fractional-order
 integrodifferential operators with fractional powers of Bessel operators is investigated
 in recent years. This theory is founded in \cite{McB,Spr,Dim,Kir1, S140p, S135, S133}
 and is developed in \cite{S127, S123, 18, S700, SS, FJSS}.

Numerous applications of differential equations with  Bessel
operators are provided in cited papers as well.
    Thus, this class of equations is important both for the theory of partial
    differential equations and for practical applications.

Note that only linear    differential equations are considered in
the present book. Nonlinear ones require other approach; they form
 another research direction (see, e.\,g., \cite{Kap1,ZaiPol,Rad3,Rad4}).

The following useful Darboux--Weinstein
 relation (the term is introduced by J.-L.~Lions)
  can be verified intermediately:
\begin{equation}
B_{\nu} \lr{ y^{-2 \nu} f(y)} = y^{-2 \nu} B_{\nu} f(y);
 \label{9}
\end{equation}
 this easily reduces the case where $\Re \nu < 0$ to  the case where $\Re \nu > 0.$

Define the class of transmutation operators intertwining
 the Bessel differential operator with the second derivative:
\begin{equation}
\label{151} T\lr{B_\nu} f=\lr{D^2} Tf,\quad B_{\nu}=D^2+\frac{2\nu
+1}{x}D,\quad D^2=\frac{d^2}{dx^2},\quad \nu \in \mathbb{C}.
\end{equation}
    To construct  transmutation operators, one can establish relations between
    solutions of differential equations. Solutions of equations of the kind
     $B_\nu f=\lambda f$ are Bessel functions, while solutions of the equation
      $D^2f=\lambda f$ are trigonometric or exponential functions.
  Therefore transmutation operators of kind \eqref{151} grow from
  the following relations (Poisson and Sonin):
\begin{equation}
\label{152} J_{\nu}
(x)=\frac{1}{\sqrt{\pi}\Gamma({\nu+\frac{1}{2}})2^{\nu-1}x^\nu}
\int\limits_0^x \left(
x^2-t^2\right)^{\nu-\frac{1}{2}}\cos(t)\,dt,\quad
 \Re \nu> \frac{1}{2},
\end{equation}
    and
 \begin{equation}
\label{153} J_{\nu}
(x)=\frac{2^{\nu+1}x^{\nu}}{\sqrt{\pi}\Gamma({\frac{1}{2}-\nu})}
\int\limits_x^\infty \left(
t^2-x^2\right)^{-\nu-\frac{1}{2}}\sin(t)\,dt,\quad -\frac{1}{2}<
\Re \nu<\frac{1}{2}.
\end{equation}
    Integral \eqref{152} is studied since 1769 (Euler).
    For $\nu=0$, its value is computed by Parseval in 1805.
    For integer values of $\nu$, relation \eqref{152} is obtained by Plana in 1821.
    For half-integer
    values of $\nu$, Poisson derived it in 1823.  For the general case, relation \eqref{152}
    is established by Lommel in 1868. Relation \eqref{153} is derived
    by Sonin in 1880.

\begin{definition}
 The {\it Poisson transmutation operator} is the expression
\begin{equation}
\label{154} P_{\nu}f=\frac{1}{\Gamma(\nu+1)2^{\nu}x^{2\nu}}
\int\limits_0^x \left(
x^2-t^2\right)^{\nu-\frac{1}{2}}f(t)\,dt,\quad \Re \nu>
-\frac{1}{2}.
\end{equation}
The {\it Sonin transmutation operator} is the expression
\begin{equation}
\label{155}
S_{\nu}f=\frac{2^{\nu+\frac{1}{2}}}{\Gamma(\frac{1}{2}-\nu)}\frac{d}{dx}
\int\limits_0^x \left(
x^2-t^2\right)^{-\nu-\frac{1}{2}}t^{2\nu+1}f(t)\,dt,\quad \Re \nu<
\frac{1}{2}.
\end{equation}
 Operators \eqref{154}--\eqref{155}
 act as  transmutation operators according to the relations
\begin{equation}
\label{156} S_\nu B_\nu=D^2 S_\nu \quad \textrm{and}  \quad  P_\nu
D^2=B_\nu P_\nu.
\end{equation}
\end{definition}

They can be defined for all complex values of $\nu.$

The idea to study operators of kind \eqref{154}--\eqref{155}
 is attributed to Liouville.
  Sonin was the first to use them in the framework of the theory
  of  Bessel functions. Delsarte was the first to use them as transmutation
operators. Further, Delsarte ideas are used to continue the
investigation of these operators (Delsarte and Lions).
 That is why we call \eqref{154}--\eqref{155}
\emph{ So\-nin--Poisson--Delsarte
  transmutation operators}
  (see \cite{Lev7} as well).

It is possible to say that
  So\-nin--Poisson--Delsarte
  operators \eqref{154}--\eqref{155}
    are the most famous objects of the theory of transmutation
operators; hundreds of papers are devoted to them, their
applications, and their generalizations.

Apart from  So\-nin and Poisson  operators, we need similar ones
obtained as follows: the Rie\-mann--Li\-ouville
 operators $I^{\mu}$ are changed for the operators $I^{\mu}_e$
 defined by the relation
\begin{equation}
I^{\mu}_e =  \mathcal{E} I^{\mu}  \mathcal{E}^{-1},
\label{10}
\end{equation}
 where $\mathcal{E}$ is the multiplier by the function $e^x.$
 These operators are introduced in Chap. \ref{ch2}. In particular,
 the following relation is established in this chapter:
\begin{equation}
\mathcal{J}_{\mu, e} f(y) \equiv  I^{\mu}_e I^{-\mu} f(y) =  f(y)
- \mu \int\limits_y^{\infty} f(t) \, \Phi (\mu+1, 2; y-t) \, dt,
\label{11}
\end{equation}
    where  $\Phi (a, c; z)$ is the generalized hypergeometric function
 (other properties of the operators $I^{\mu}_e$ are provided in \cite{SKM}).
 Then the following relations hold:
\begin{equation}
P_{\nu, e} \equiv P_{\nu}^{\frac{1}{2}-\nu} I_e^{\nu-\frac{1}{2}}
=P_{\nu} \mathcal{J}_{\nu-\frac{1}{2}, e} \label{12}
\end{equation}
    and
\begin{equation}
S_{\nu, e} \equiv I_e^{\frac{1}{2}-\nu} S_{\nu}^{\nu-\frac{1}{2}}
= \mathcal{J}_{\frac{1}{2}-\nu, e} S_{\nu}. \label{13}
\end{equation}
  In Chap. \ref{ch2}, we study relations between operators of the Riemann--Liouville
  type with the following Fourier and Hankel transformations:
\begin{equation}
F f \lr{\eta} = \int\limits_{-\infty}^{\infty} f(y) e^{- i y \eta}
dy, \  F^{-1} g(y)  = \frac{1}{2 \pi}
\int\limits_{-\infty}^{\infty} g(\eta) e^{ i y \eta} d \eta,
\label{14}
\end{equation}
\begin{equation}
F_{-} f \lr{\eta} = \int\limits_{0}^{\infty} f(y) \sin \lr{y \eta}
dy, \  F^{-1}_{-}  g(y)  = \frac{2}{\pi} \int\limits_{0}^{\infty}
g(\eta) \sin \lr{y \eta} d \eta, \label{15}
\end{equation}
    and
\begin{equation}
F_{\nu} f \lr{\eta} = \int\limits_{0}^{\infty} f(y) j_{\nu} \lr{y
\eta} y^{2 \nu+1} dy,  \  F^{-1}_{\nu}  g(y) = \frac{2^{-2
\nu}}{\Gamma^2 (\nu+1)} \int\limits_{0}^{\infty} g(y) j_{\nu}
\lr{y \eta} \eta^{2 \nu+1} d \eta, \label{16}
\end{equation}
    where $j_{\nu} (t) = \dfrac{2^{\nu} \Gamma (\nu+1) J_{\nu}
(t)}{t^{\nu}}$ is the normalized (small)  Bessel function, while
$J_{\nu}$ is the classical first-kind
 Bessel function.

In the same way, a family of other transmutation operators
generalizing Bessel and Riesz one-di\-men\-sional
 potentials can be introduced. In \cite{Kat1}, applications of such operators can be found.

The core notion of generalized translations is introduced by
 Delsarte on the base of
So\-nin--Pois\-son--Del\-sarte
  transmutation operators.

\begin{definition}
    The {\it generalized translation operator} is the solution $u(x,y)=T_x^yf(x)$
    of the problem
\begin{equation}
\label{157}
\begin{gathered}
(B_\nu)_y u(x,y)=(\frac{\partial^2}{\partial y^2}+ \frac{2\nu
+1}{y}\frac{\partial}{\partial y})\  u(x,y)=
 \frac{\partial^2}{\partial x^2} \ u(x,y),
 \\
u(x,0)=f(x),\quad  u_y (x,0)=0.
\end{gathered}
\end{equation}
\end{definition}

To justify this term, we note that, in the special case
$\nu=-\dfrac{1}{2}$, the generalized translation operator is
reduced to the operator
$$T_x^yf(x)=\frac{1}{2}\left( f(x+y)+f(x-y)\right).$$
 For the generalized translation operator defined by \eqref{157}, the following
 explicit relation is obtained by Delsarte:
\begin{equation}
\label{158} T_x^y
f(x)=\frac{\Gamma(\nu+1)}{\sqrt{\pi}\Gamma(\nu+\frac{1}{2})}
\int\limits_0^{\pi}f(\sqrt{x^2+y^2-2xy \cos t }) \sin^{2\nu}t\,dt.
\end{equation}
    Arbitrary pairs of differential (or even arbitrary) operators can be considered
    in Definition  \eqref{157} as well. For example, the following definition yields
    the classical translation:
$$\frac{\partial u}{\partial x}=\frac{\partial u}{\partial y},\quad u(x,0)=f(x),
\quad  T_x^y f(x)=f(x+y).$$
    Note that the generalized translation operators given by \eqref{157}-\eqref{158}
    are explicitly expressed via the
So\-nin--Poisson--Delsarte
   transmutation operator given by \eqref{154}-\eqref{155}
(see \cite{Lev2,Lev3, Mar9}).

However, from the viewpoint of applications to solutions of
partial differential equations with singularities at coefficients,
investigated in the present monograph, the specified
    So\-nin--Poisson--Delsar\-te
 operators possess a number of disadvantages preventing their applying in various important
 cases. The disadvantages ar as follows.
  First,
 So\-nin--Poisson--Delsarte
 operators introduced above are  transmutation operators only on the set of even functions;
 therefore, it is impossible to consider functions with singularities at the origin.
 Secondly, they do not preserve the compactness of the support and the property to rapidly
 decrease at infinity.
    Thirdly, they change the smoothness of transmuted functions. For the first time,
    the attention to this fact is    drawn by
    J.-L.~Lions (see \cite{Lio1}).

This causes the necessity to introduce and study other classes of
transmutation operators for differential equations with Bessel
operators.

\subsubsection{Fractional powers of Bessel operators}\label{sec3.3.6}

Note that there exists a theory of fractional powers of Bessel
operators and their applications to fractional-order
 differential equations. This approach allows one to avoid the implicit
 definition of fractional powers of Bessel operators in images of the Hankel integral
 transformation; instead, they are explicitly defined as a particular integral operator
 with special functions in the kernels (see \cite{McB,
Spr, Dim, Kir1, S140p, S135, S133, S127, S123, 18, S700, SS,
FJSS}).

Similarly to classical derivatives, the above fractional powers
can be defined on the semiaxis by means of the following natural
property: in Hankel images, they act as multipliers by power
functions. Such an approach is justified (until a better one is
found) though no explicit  representations of fractional powers
can be obtained on this way.
    However, imagine that we are able to define usual  Riemann--Li\-o\-u\-vil\-le
    operators of fractional integrating only in Laplace or Mellin images,
    but no integral relations for these operators are known.
 It is clear that such a theory would be quite poor then (it would lose the majority
 of its
 most useful and nice results).
 The current theory fractional powers of Bessel
 operators is about at the explained state; that it is why it is important and
 interesting
 to construct their explicit integral forms.

 The main definitions and properties of fractional powers of Bessel operators are
  as follows.

    We consider real powers of the singular Bessel differential operator
\begin{equation}\label{Bess}
B_\nu= D^2+\frac{\nu}{x}D,\qquad \nu\geq 0
\end{equation}
    on the real semiaxis $(0,\infty).$

\begin{definition}
    Let $f(x)\in C^{2k}(0,b].$
    Under the assumption that $f^{(i)}(b)=0, 0\leq i \leq 2k-1,$ $k \in \N$,
 define the
  {\it right-hand Bessel fractional integrating operator} by the relation
\begin{multline*}
(B_{b-}^{{\nu},k}f)(x)=\frac{1}{{\Gamma (2k)}}\int\limits_{x}^{b}
\left(\frac{y^{2}-x^{2}}{2y}\right)^{2k-1}
{_2}F_1(k+\frac{{\nu}-1}{2},k;2k;1-\frac{x^{2}}{y^{2}})
f(y)\,dy\\
=\frac{\sqrt{{\pi}}}{2^{2k-1}{\Gamma (k)}}
\int\limits_{x}^{b}(y^{2}-x^{2})^{k-\frac{1}{2}}
\left(\frac{y}{x}\right)^{\frac{{\nu}}{2}}
P_{\frac{{\nu}}{2}-1}^{\frac{1}{2}-k}
\left(\frac{1}{2}\left(\frac{x}{y}+\frac{y}{x}\right)\right)
f(y)\,dy.
\end{multline*}
\end{definition}

\begin{definition}
   Under the assumption that $f^{(i)}(a)=0, 0\leq i \leq 2k-1, k \in N$,
 define the
 {\it right-hand Bessel fractional integrating operator} by the relation
\begin{multline*}
(B_{a+}^{{\nu},k}f)(x)=\frac{1}{{\Gamma (2k)}}\int\limits_{a}^{x}
\Bigl(\frac{x^{2}-y^{2}}{2x}\Bigr)^{2k-1}
{_2}F_1(k+\frac{{\nu}-1}{2},k;2k;;1-\frac{y^{2}}{x^{2}}) f(y)\,dy
\\
=\frac{\sqrt{{\pi}}}{2^{2k-1}{\Gamma (k)}}
\int\limits_{a}^{x}(x^{2}-y^{2})^{\left(k-\frac{1}{2}\right)}
\left(\frac{x}{y}\right)^{\frac{{\nu}}{2}}
P_{\frac{{\nu}}{2}-1}^{\frac{1}{2}-k}\left(\frac{1}{2}
\left(\frac{x}{y}+\frac{y}{x}\right)\right)
 f(y)\,dy,
\end{multline*}
    where ${_2}F_{1}$ is the Gauss hypergeometric function and
$P_{\nu}^{\mu}(z)$ is the Legendre function.
\end{definition}

It is useful to express Bessel fractional integrals via  Legendre
functions; it simplifies the original definition because the Gauss
hypergeometric function depends on three parameters, while the
Legendre function depends on two ones.

There exist versions of fractional Bessel integrals with arbitrary
integration limits as well as their further modifications (see
\cite{S140p, S135, S133, S127,S123, 18, S700, SS}).
 Most frequently, the operators $B_{0+}^{{\nu},k}$ and $B_{\infty-}^{{\nu},k}$ are
    used.

\begin{property}\label{p1}
 If $\nu=0,$ then the  fractional Bessel integral on the semiaxis $B_{0+}^{0,-\alpha}$
 is reduced to the  fractional Riemann--Liouville integral; namely, the relation
$$
(B_{\infty-}^{0,
\alpha}f)(x)=\frac{1}{\Gamma(2\alpha)}\int\limits_x^{\infty}(y-x)^{2\alpha-1}f(y)dy=
(I_{-}^{2\alpha}f)(x).
$$
    holds.
\end{property}

\begin{property}\label{p2}
    The relation
\begin{equation}\label{Prop2}
(B_{\infty-}^{{\nu},\alpha}f)(x)=
\frac{1}{2^{2\alpha}}J_{x^2}^{2\alpha,\frac{\nu{-}1}{2}-\alpha,-\alpha}
\left(x^{\frac{\nu{-}1}{2}}f(\sqrt{x})\right),
\end{equation}
    where
\begin{equation}\label{Saigo1}
J_{x}\,^{\gamma,\beta,\eta}f(x)=
\frac{1}{\Gamma(\gamma)}\int\limits_x^\infty(t-x)^{\gamma-1}t^{-\gamma-\beta}
\,_2F_1\left(\gamma+\beta,-\eta;\gamma;1-\frac{x}{t}\right)f(t)dt
\end{equation} is the Saigo  fractional integral   $($see
\rm{\cite{Rep}),}
   $\gamma>0,$ and $\beta$ and $\theta$ are real numbers, holds.
\end{property}

\begin{property}\label{p3}
    If $\lim\limits_{x\rightarrow +\infty}g(x)=0$ and
$\lim\limits_{x\rightarrow +\infty}g'(x)=0,$ then the relation
$$
(B_{\infty-}^{\nu,-1}B_\nu g)(x)=g(x)
$$
    holds.
\end{property}

\begin{property}\label{p4}
    If $x>0$ and $m+2\alpha+\nu<1,$ then the relation
\begin{equation}\label{Prop4}
B_{\infty-}^{\nu,\alpha}\,x^m=x^{2\alpha+m}\,2^{-2\alpha}\,\Gamma\left[
 \begin{array}{cc}
  $$-\alpha-\dfrac{m}{2},$$ & $$-\dfrac{\nu-1}{2}-\alpha-\dfrac{m}{2}$$ \\
  $$\dfrac{1-\nu-m}{2},$$ & $$-\dfrac{m}{2}$$  \\
        \end{array}
  \right]
\end{equation}
    holds.
\end{property}

\begin{property}\label{p5}
If $\alpha>0,$ then the  Mellin transform of the fractional Bessel
integral on the semiaxis has the form
\begin{equation}\label{Mellin2}
  M((B_{\infty-}^{\nu,\alpha}f)(x))(s)=\frac{1}{2^{2\alpha}}\,\,
  \Gamma\left[\begin{array}{cc}
$$\dfrac{s}{2},$$ & $$\dfrac{s}{2}-\dfrac{\nu-1}{2}$$ \\
$$\alpha+\dfrac{s}{2}-\dfrac{\nu-1}{2},$$ & $$\alpha+\dfrac{s}{2}$$  \\
                                                           \end{array}
                                                         \right] f^*(2\alpha+s).
\end{equation}
\end{property}

\begin{property}\label{p6}
    If $\alpha$ and $\beta$ are positive, then the fractional Bessel integral on
    the semiaxis possesses the semigroup property
\begin{equation}\label{SemiGroup3}
    B_{\infty-}^{\nu,\alpha} B_{\infty-}^{\nu,\beta}f=
     B_{\infty-}^{\nu, \alpha+\beta}f.
\end{equation}
\end{property}

 The resolvent expression for fractional powers of the Bessel operator, provided below,
 generalizes the famous relation for Riemann--Liouville
 fractional integrals, provided in \cite{HT} without a proof.
 For the first time, this relation is proved in \cite{Dzh1} by
 means of the     method of successive approximations, which is
 usual  for the theory of integral equations.
 Therefore, it seems to be reasonable (from the historical viewpoint)
 to call the resolvent relation for Riemann--Liouville
 fractional integrating operators the
\textit{Tamarkin--Hille--Dzhrbashyan}
    relation.
    Also, \cite{Dzh1} is the first book where properties of the Mittag-Lef\-f\-ler
 function are studied in details for the first time.

The
 Tamarkin--Hille--Dzhrbashyan
 relation is the most known application of Mittag-Leffler
  functions. Numerous papers about applications of fractional calculus
  to differential equations are based on it.
  In \cite{GoMa}, it is proposed to call the  Mittag-Leffler
  function the Queen function of the fractional
calculus.

\begin{property}\label{p7}
    The following integral representation holds {\rm(}on suitable functions{\rm)}
  for the resolvent of the Bessel fractional integrating operator
$B_{0+}^{\nu, \alpha},$ ${0\leq \nu <1:}$
\begin{equation}
R_\lambda f=-\frac{1}{\lambda}f-\frac{1}{\lambda}\int\limits_0^x
K(x,y)f(y)\,dy,
\end{equation}
    where
\begin{gather}\label{R}
K(x,y)=\frac{2y}{x^2-y^2}\int\limits_0^1 S_{\alpha,\nu}(z(t))
\frac{dt} {
\left( t\left(1-t\right) \right) ^{\frac{\nu+1}{2}} },\\
\nonumber
z(t)=\left(\frac{t(1-t)\left(x^2-y^2\right)^2}{\left(1-\left(1-\frac{x^2}{y^2}
\right)t\right)4y^2}\right)^\alpha,\quad\textrm{and}\quad
S_{\alpha,\nu}(z)=\sum\limits_{k=1}^{\infty}\frac{z^k}{\Gamma(\alpha
k+\frac{\nu-1}{2})\Gamma(\alpha k-\frac{\nu-1}{2})}
\end{gather}
 is a kind of  the  Wright--Fox
 hypergeometric function in the  Wright form $($see {\rm\cite{Kir4, Kir5, Kir6, KiSa})}.
\end{property}

Similar functions are applied in \cite{Pshu1} as well.
 They are a special case of more general Wright--Fox
 functions originally introduced as generalizations of Bessel functions
  (special Wright--Fox
 functions are described above). Further simplification of the representation
 of kernel \eqref{R} (if it is possible) is an interesting problem as well.

    The obtained resolvent relation for fractional powers of the
    Bessel operator allows one to consider problems with various boundary-value
    conditions for the
    ordinary integrodifferential equation of the kind
\begin{equation*}
B_\nu^\alpha u(x)-\lambda u(x)=f(x),
\end{equation*}
    where $B_\nu^\alpha$ is a fractional power of the Bessel operator.
  Also, it is possible to consider partial differential equations with fractional powers
  of Bessel operators and their Gerasimov--Caputo
  generalizations, which can be introduced on the base of the
  known generalized Taylor relations.
  Examples of such equations are the generalization of the $B$-elliptic
   fractional Laplace--Bessel
   equation
\begin{equation*}
\sum\limits_{k=1}^n B_{\nu_k}^{\alpha_k} u(x_1,x_2,\ldots
x_n)=f(x_1,x_2,\ldots x_n)
\end{equation*}
    and nonstationary equations of the kind
\begin{equation*}
B_{\nu,t}^\alpha u(x,t)=\Delta_x u(x,t)+f(x,t).
\end{equation*}
    Note that, to consider spectral properties of such equations,
    one has to study asymptotic properties of the function $K(x,y)$ from relation
    \eqref{R} in the complex plane as well root distributions of this function.

Riemann--Liouville fractional integrodifferentiating operators
have numerous applications because they are contained in the
remainder of the Taylor expansion.
    Therefore, once fractional powers of the  Bessel operator are defined,
    the problem to construct a generalized Taylor relation, where functions
  are expanded with respect to powers of the  Bessel operator, arises.
  This problem is not new. It has a certain history.

For the first time, expansions with respect to powers of the
Bessel operator (Taylor--Delsarte
 series) are obtained in \cite{Del2, Lev2, Lev3, Lev4}.
 In \cite{FN}, the general way to construct them is explained in terms of
 operator-analytic functions.
 However, Taylor--Delsarte
 series allow one to expand the generalized translation (instead of classical one)
 with respect to powers of the  Bessel operator.
 Actually, such expansions are operator variants of series for the Bessel functions
 in the same sense that classical Taylor series are operator version of series
 expansions
 of the exponential function.  Taylor--Delsar\-te
 series have their own application area. However, to solve partial
 differential equations numerically, one needs generalized Taylor
 expansions of another nature.
 Relations for the generalized translation are useless for the
 layer-to-layer
 recomputing of solutions, e.\,g., for the mesh method;
 instead, we need relations for the classical translation,
 expressing the solution on the computed layer via its
values at preceding layers. It turns out that, compared with the
generalized translation, it is much harder to construct such a
relation for the classical one; the reason is that, unlike the
case of the generalized translation,  they are not direct analogs
of identities known for special functions.

For the first time, the Taylor expansion of the needed type is
considered in \cite{KaKa}; the aim was to use it for the numerical
solving of equations with the Bessel operator by means of the
finite element method.
    However, the obtained result cannot be treated as desired relations in the explicit
    form because the coefficients are expressed by indefinite constants determined
    by a system of recurrent relations, while the remainder kernel is represented
    by a multiple integral.
    The reason is that it is impossible to guess the explicit form
    of kernels and remainders simultaneously until we know the particular
    remainder expression via fractional powers of the Bessel operator.

The above problem set by Katrakhov is solved in \cite{S140p} (see
also \cite{S135,S133, S127,S123, 18, S700, SS}): the final form of
the Taylor relation with powers of the  Bessel operator and
remainder in the form of a fractional power of the Bessel operator
is found.

\chapter[Sonin--Poisson--Delsarte
 Transmutation Operators and Their Modifications]
 {Sonin--Poisson--Delsarte
 Transmutation Operators\\
  and Their Modifications}\label{ch2}

In this chapter, we study properties of
 Sonin--Poisson--Delsarte
 and
Erd\'elyi--Kober transmutation operators and introduce new
transmutation operators. We establish their relations
    with the  Fourier and Hankel transformations and with
    Liouville-type fractional integrals. To do that, we use
    classical Liouville operators and operators located between them and Bessel potentials
    and inheriting good properties of both ones.
    The reason to introduce such operators is as follows:  Liouville operators are unbounded
    in Sobolev spaces if the domain is unbounded.
    We use transmutation operators to reduce Sobolev spaces to weight function spaces
     introduced by Kipriyanov. For the one-dimensional
     case, we introduce and investigate new function spaces as well.

    Further, we construct a class of transmutation operators
    transforming the multidimensional Laplace operator into the
    ordinary double-differentiation
    operator with respect to the radial
    variable. These operators are defined on functions admitting
    singularities at one point (for definiteness, we assume that it is the origin).

\section{One-Dimensional Transmutation Operators}\label{sec4}

\subsection{Main constructions of transmutation operators}\label{sec4.1}

 Recall the main Lions definition of transmutation operators (see \cite{58}).
 Let $A$ and $B$ be linear operators.
 Operators  $P$ and $S$ are called transmutation operators for $A$ and $B$ if
\begin{equation}
B = P A S, \   A = S B P, \label{1.1.1}
\end{equation}
    and the operators $P$ and $S$ are mutually invertible. This definition is not completely
    strict because domains and images if the participating operators are not given.
    This are described separately for each particular case.
    A basic example of transmutation operators is the Fourier transformation.
In this case,  $P$ and  $S$ are the direct and inverse Fourier
transformations respectively, $B$ is the differentiation operator,
and  $A$ is the multiplier by the dual variable. In this example,
 $A$ and $B$ are operators of different nature.
 In  the present book, we consider only the case where $A$ and $B$ are differential operators.
 In particular, in the current section, we consider the following operators:
\begin{equation*}
A=\frac{\pr^2}{\pr y^2}  \quad \textrm{and}  \quad B=B_{\nu}=
\frac{\pr^2}{\pr y^2} +\frac{2 \nu+1}{y}\frac{\pr}{\pr y}.
\end{equation*}
    The operator $B_{\nu}$ is called the \emph{Bessel operator} with parameter $\nu.$
 In this situation, the following  transmutation operators (see \cite{Lev7, Mar2}) are known:
\begin{equation}
P_{\nu, 0}^{\nu+\frac{1}{2}} f(y) =  \frac{2\, \Gamma \lr{\nu+1}
}{\sqrt{\pi}\, \Gamma \lr{\nu+\frac{1}{2}}} \, y^{-2\nu}
\int\limits_0^y \lr{y^2-t^2}^{\nu-\frac{1}{2}} f(t) \, dt
\label{1.1.2}
\end{equation}
    and
\begin{equation}
S_{\nu, 0}^{-\nu-\frac{1}{2}} f(y) =  \frac{\sqrt{\pi}}{ \Gamma
\lr{\nu+1} \, \Gamma \lr{\frac{1}{2} - \nu}} \frac{\pr}{\pr y}
\int\limits_0^y \lr{y^2-t^2}^{-\nu-\frac{1}{2}} t^{2 \nu+1} f(t)
\, dt, \label{1.1.3}
\end{equation}
 where $\Gamma(\mu)$ denoted the Euler gamma-function.
 The function $f(y)$ is defined on the semiaxis $E_{+}^1 = \{y>0\}$ and is smooth,
 compactly supported, and even. If an integral in the above relations diverges,
 one has to pass to its regularization. The operator $P_{\nu, 0}^{\nu+\frac{1}{2}}$ is called
 the \emph{Poisson operator}; the operator  $S_{\nu,0}^{-\nu-\frac{1}{2}}$ is called
 the \emph{Sonin  operator}. Also, we use the general term
   \emph{Sonin--Poisson--Delsarte
  transmutation operators}.

    From the viewpoint of applications to singular elliptic
    partial differential boundary-value
    problems studied in the present monograph, the specified operators possess a number
     of disadvantages
    preventing to apply them in this case.
     First, the operators $P_{\nu, 0}^{\nu+\frac{1}{2}}$ and
$S_{\nu, 0}^{-\nu-\frac{1}{2}}$ satisfy relations \eqref{1.1.1}
only on the set of even functions, while we have to investigate
functions increasing near the origin.
    Secondly, they do not preserve the compactness of the support and the property to rapidly
 decrease at infinity.
    Thirdly, they change the smoothness of transmuted functions. For the first time,
    the attention to this fact is    drawn by
    J.-L.~Lions (see \cite{Lio1, Lio2, Lio3}, but this can be easily found directly as well.

    Introduce transmutation operators free of all the disadvantage
    mentioned above. First, for operators \eqref{1.1.2}-\eqref{1.1.3},
      we find formally adjoint ones with respect to bilinear forms
      of the kind
\begin{equation*}
\lr{f, g}_{\nu} = \int\limits_0^{\infty}  y^{2 \nu+1} f(y) g(y) \, dy.
\end{equation*}
    We have the relation
\begin{equation*}
\lr{P_{\nu, 0}^{\nu+\frac{1}{2}} g, f}_{\nu} = \int\limits_0^{\infty} t^{2 \nu +1} f(t)
 \frac{2 \, \Gamma \lr{\nu+1}}{\sqrt{\pi}  \, \Gamma \lr{ \nu+\frac{1}{2} }} \,
  t^{- 2 \nu}  \int\limits_0^t \lr{t^2-y^2}^{\nu-\frac{1}{2}} g(y) \, dy dt
\end{equation*}
\begin{equation*}
 =  \frac{2 \, \Gamma \lr{\nu+1}}{\sqrt{\pi}  \, \Gamma \lr{ \nu+\frac{1}{2} }}
   \int\limits_0^{\infty} g(y)
      \int\limits_0^t t f(t) \lr{t^2-y^2}^{\nu-\frac{1}{2}}  \, dt dy =
       \lr{g,  S_{\nu}^{\nu+\frac{1}{2}} }_{-\frac{1}{2}},
\end{equation*}
    where the operator  $S_{\nu}^{\nu+\frac{1}{2}}$ is defined by
    the relation
\begin{equation}
S_{\nu}^{\nu+\frac{1}{2}} f(y)  =
\frac{2 \, \Gamma \lr{\nu+1}}{\sqrt{\pi}  \, \Gamma \lr{ \nu+\frac{1}{2} }}
  \int\limits_y^{\infty}  \lr{t^2-y^2}^{\nu-\frac{1}{2}} t f(t)  \, dt.
\label{1.1.4}
\end{equation}
    In the same way, we obtain the relation
\begin{equation*}
    \lr{S_{\nu, 0}^{-\nu-\frac{1}{2}} g, f}_{-\frac{1}{2}} =
     \frac{\sqrt{\pi} }{  \Gamma \lr{\nu+1} \, \Gamma \lr{ \frac{1}{2}-\nu }}
      \int\limits_0^{\infty} \frac{\pr}{\pr t} \int\limits_0^t
            \lr{t^2-y^2}^{-\nu-\frac{1}{2}} y^{2 \nu+1} g(y) \, dy  f(t) \, dt
\end{equation*}
\begin{equation*}
 =\frac{-\sqrt{\pi} }{  \Gamma \lr{\nu+1} \, \Gamma \lr{ \frac{1}{2}-\nu }}
  \int\limits_0^{\infty} y^{2 \nu+1} g(y)   \int\limits_y^{\infty}
    \lr{t^2-y^2}^{-\nu-\frac{1}{2}} \frac{\pr f(t)}{\pr t}   dt dy =
     \lr{g,  P_{\nu}^{-\nu-\frac{1}{2}} }_{\nu},
\end{equation*}
    where
\begin{equation}
    P_{\nu}^{-\nu-\frac{1}{2}} f(y)  =
     \frac{-\sqrt{\pi} }{  \Gamma \lr{\nu+1} \, \Gamma \lr{ \frac{1}{2}-\nu }}
      \int\limits_y^{\infty}  \lr{t^2-y^2}^{-\nu-\frac{1}{2}} \frac{\pr f(t)}{\pr t} \,
       dt.
    \label{1.1.5}
\end{equation}
    It is convenient to select other operators obtained from \eqref{1.1.4}-\eqref{1.1.5}
  by means of the moving of the differentiation operator from the second operator
  to the first one:
\begin{equation}
S_{\nu}^{\nu-\frac{1}{2}} f(y) = \frac{- 2 \, \Gamma \lr{\nu+1}}{
\sqrt{\pi} \, \Gamma \lr{ \nu+\frac{1}{2} }} \frac{\pr}{\pr y}
\int\limits_y^{\infty} \lr{t^2-y^2}^{\nu-\frac{1}{2}} t f(t) \, dt
\label{1.1.6}
\end{equation}
    and
\begin{equation}
P_{\nu}^{\frac{1}{2}-\nu} f(y) =  \frac{\sqrt{\pi} }{\Gamma
\lr{\nu+1}\, \Gamma \lr{\frac{1}{2}-\nu}} \int\limits_y^{\infty}
\lr{t^2-y^2}^{-\nu-\frac{1}{2}} f(t) \, dt, \label{1.1.7}
\end{equation}
    where $a^{\mu},$ $a>0,$ denotes (as above) the     principal branch
of the multifunction: $a^{\mu} = \exp (\mu \ln a).$

    The above reasoning are rather formal. They are provided to
    clarify the relation between the Poisson and Sonin operators and the
     Erd\'elyi--Kober.

Introduce the following notation.
    Let $R$ denote a positive number or the infinity.
    Let $C^{\infty} (0, R)$ denote the set of functions infinitely
    differentiable on the interval $(0, R)$.
  Let  $\mathring{C}^{\infty} (0, R)$  denote the subset of functions from
   $C^{\infty} (0, R)$ such that the support of each one is
   compact in $(0, R)$.
   Let $C^{\infty} [0, R)$   denote the subset of functions from $C^{\infty} (0, R),$
   such that all their derivatives are continuous up to the left-hand
   edge. Let $\mathring{C}^{\infty} [0, R)$   denote the subset of functions from
$C^{\infty} [0, R)$, vanishing in a neighborhood of the right-hand
   edge. Let $\mathring{C}^{\infty}_{\{0\}} (0, R)$
    denote the subset of functions from $C^{\infty} (0, R)$,
    vanishing in a neighborhood of the right-hand
   edge. Note that functions
$\mathring{C}^{\infty}_{\{0\}} (0, R)$ are allowed to increase
arbitrarily rapidly in a neighborhood of the point $y=0.$
    We say that a function $f \in C^{\infty} [0, R)$ is \emph{even}
    (\emph{odd}) if $D^k f(0) = 0$ for all odd (even) nonnegative values of $k.$
 Here, $D= \dfrac{\pr}{\pr y}$ and $D^k = D D^{k-1}.$
 A function $f \in C^{\infty} [0, R)$ is even if and only if the function
  $g(y)=f (\sqrt{y})$ belongs to the space $C^{\infty} [0, R^2).$
  This is easily proved by means of the Taylor expansion.
  The set of even (odd) functions is denoted by
  $C^{\infty}_{+} [0, R)$ ($C^{\infty}_{-} [0, R)$).
Assign $\mathring{C}^{\infty}_{\pm} [0, R]=C^{\infty}_{\pm} [0, R)
\cap \mathring{C}^{\infty} [0, R)$.
 Also, we use the notation
 $(0, \infty) = E^1_+$ and $[0, \infty) =\ov{E^1_+}.$

 For $\Re \nu \geq 0$, define the transmutation operators $S_{\nu}^{\nu - \frac{1}{2}}$
 on functions $f\in \mathring{C}^{\infty}_{\{0\}} \lr{E^1_+}$ by relation \eqref{1.1.6}.
 For $0 \leq \Re \nu < \dfrac{1}{2}$, define the operators $P_{\nu}^{\nu - \frac{1}{2}}$
 (on the same function class) by relation \eqref{1.1.7}.
    If $\Re \nu \geq \dfrac{1}{2}$, then the value of the function
 $P_{\nu}^{\frac{1}{2} - \nu } f(y)$ at each positive point $y$ is defined as
 the analytic continuation of integral \eqref{1.1.7} with respect to parameter $\nu.$
 If $\dfrac{1}{2} \leq \Re \nu <N+ \dfrac{1}{2},$ where $N$ is a positive integer,
 then this is equivalent to the defining of the operator $P_{\nu}^{ \frac{1}{2} - \nu }$
 by the relation
\begin{equation}
P_{\nu}^{\frac{1}{2}-\nu} f(y) =  \frac{(-1)^{N} 2^{-N} \sqrt{\pi}
} {\Gamma \lr{\nu+1}\, \Gamma \lr{N-\nu+\frac{1}{2}}}
\int\limits_y^{\infty} \lr{t^2-y^2}^{N-\nu-\frac{1}{2}}
\lr{\frac{\pr}{\pr t} \frac{1}{t}}^N f(t) \, dt, \label{1.1.8}
\end{equation}
    where the integral at the right-hand
    side converges.

    The introduced operators are transmutation ones because the
    following assertion holds.

\begin{theorem}\label{theorem:1_1_1}
    If $\Re \nu \geq 0,$ then the operators $P_{\nu}^{ \frac{1}{2} - \nu }$ and
  $S_{\nu}^{\nu -\frac{1}{2}}$ map the space
   $\mathring{C}_{\{0\}}^{\infty} \lr{E_{+}^{1}}$ into itself
  one-to-one
 and are mutually inverse. The relations
\begin{equation}
B_{\nu} P_{\nu}^{\frac{1}{2}-\nu}  = P_{\nu}^{\frac{1}{2}-\nu} D^2
    \quad \textrm{and}   \quad
D^2 S_{\nu}^{\nu-\frac{1}{2}}  = S_{\nu}^{\nu-\frac{1}{2}} B_{\nu}
\label{1.1.9}
\end{equation}
    hold.
    \end{theorem}

\begin{proof}
 Let $\Re \nu < N + \dfrac{1}{2},$ where $N$ is a positive integer.
 The substitution  $t \to ty$ reduces relation \eqref{1.1.8} to the form
\begin{multline}
P_{\nu}^{\frac{1}{2}-\nu} f(y) =
\frac{(-1)^{N} 2^{-N} \sqrt{\pi} \, y^{2(N-\nu)}} {\Gamma \lr{\nu+1}\,
 \Gamma \lr{N-\nu+\frac{1}{2}}}  \lr{\frac{\pr}{\pr y} \frac{1}{y}}^N
  \int\limits_1^{\infty} \lr{t^2-1}^{N-\nu-\frac{1}{2}} t^{-2N} f(t y) \, dt
  \\
 = \frac{(-1)^{N} 2^{-N} \sqrt{\pi}\, y^{2(N-\nu)}} {\Gamma \lr{\nu+1}\,
 \Gamma \lr{N-\nu+\frac{1}{2}}}  \lr{\frac{\pr}{\pr y} \frac{1}{y}}^N
  \int\limits_y^{\infty} y^{2 \nu} \lr{t^2-y^2}^{N-\nu-\frac{1}{2}} t^{-2N} f(t) \, dt.
  \label{1.1.10}
\end{multline}
    This immediately implies that the function $P_{\nu}^{ \frac{1}{2} - \nu }f(y)$
 is infinitely differentiable for positive $y$ and is compactly supported
 provided that $f \in \mathring{C}_{\{0\}}^{\infty} \lr{E_{+}^{1}}.$
 Moreover, the upper bound of the support is not enlarged under the passage
 from $f$  to  $P_{\nu}^{ \frac{1}{2} - \nu }f$.
 Thus, the operator $P_{\nu}^{ \frac{1}{2} - \nu }$ maps the space
 $\mathring{C}_{\{0\}}^{\infty} \lr{E_{+}^{1}}$ into itself.
 In the same way, this is proved for the operator $S_{\nu}^{\nu -\frac{1}{2}}.$
 Let us show that
  $P_{\nu}^{ \frac{1}{2} - \nu } =\lr{S_{\nu}^{\nu - \frac{1}{2}}}^{-1}.$
  By virtue of the Fubini theorem, the following  relation holds  provided that
$\Re \nu < N + \dfrac{1}{2}$:
$$
    \begin{gathered}
P_{\nu}^{\frac{1}{2}-\nu} S_{\nu}^{\nu-\frac{1}{2}} f(y) =
\frac{(-1)^{N+1} 2^{1-N} y^{2(N-\nu)}} {\Gamma
\lr{\nu+\frac{1}{2}}\, \Gamma \lr{N-\nu+\frac{1}{2}}}
\lr{\frac{\pr}{\pr y} \frac{1}{y}}^N  \int\limits_y^{\infty} y^{2
\nu} \tau^{1-2 \nu} \\
\times \frac{\pr}{\pr \tau} \lr{\tau^{2
\nu+1} f (\tau)} \int\limits_y^{\tau}
\lr{t^2-y^2}^{N-\nu-\frac{1}{2}} t^{-2N-1}
\lr{\tau^2-t^2}^{\nu-\frac{1}{2}} \, dt d \tau.
    \end{gathered}
$$
    In the internal integral, change the variables as follows:
$\dfrac{1}{t^2} = \dfrac{1}{y^2} + z
\lr{\dfrac{1}{\tau^2}-\dfrac{1}{y^2}}.$
 This yields the relation
$$
    \begin{gathered}
\int\limits_y^{\tau} \lr{t^2-y^2}^{N-\nu-\frac{1}{2}}  t^{-2N-1}
\lr{\tau^2-t^2}^{\nu-\frac{1}{2}} \, dt = \frac{1}{2} \, y^{2 N -
2 \nu -1 } \tau^{2 \nu -1}  \lr{\frac{1}{y^2}-\frac{1}{\tau^2}}^N
 \int\limits_0^1 z^{N-\nu-\frac{1}{2}} (1-z)^{\nu-\frac{1}{2}} dz
 \\
=\frac{(\tau^2-y^2)^{2 \nu - 2N-1 }\,\Gamma \lr{N-\nu+\frac{1}{2}}
\Gamma \lr{\nu+\frac{1}{2}}}{2 \, y^{2 \nu+1} \Gamma \lr{\nu+1}}.
    \end{gathered}
$$
    Hence, the following relation holds:
$$
P_{\nu}^{\frac{1}{2}-\nu} S_{\nu}^{\nu-\frac{1}{2}} f(y) =
\frac{(-1)^{N}  y^{2N-2\nu}}  {2^N N!}   \lr{\frac{\pr}{\pr y}
\frac{1}{y}}^N  \int\limits_y^{\infty} \frac{1}{y} (\tau^2-y^2)^N
\tau^{-2 N}  \frac{\pr}{\pr \tau} \lr{\tau^{2 \nu+1} f (\tau)} d
\tau.
$$
Integrate by parts. Then differentiate with respect to the
parameter. This yields the relation
$$
    \begin{gathered}
P_{\nu}^{\frac{1}{2}-\nu} S_{\nu}^{\nu-\frac{1}{2}} f(y) =
\frac{(-1)^{N} 2^N} {2^N N!} y^{2N-2\nu}  \frac{\pr}{\pr y}
\lr{\frac{\pr}{ y \pr y} }^{N-1}  \int\limits_y^{\infty}
(\tau^2-y^2)^{N-1} \tau^{2 \nu-2 N}   f (\tau) d \tau
    \\
=\frac{-N!} { N!} y^{2N-2\nu}  \frac{\pr}{\pr y}
\int\limits_y^{\infty}  \tau^{2 \nu-2 N}   f (\tau) d \tau=f(y).
    \end{gathered}
$$
The relation
 $S_{\nu}^{\nu - \frac{1}{2}}P_{\nu}^{\frac{1}{2}-\nu} f(y)=f(y)$
 is proved in the same way.
 Thus, it is proved that the operators
  $P_{\nu}^{\frac{1}{2}-\nu}$ and
   $S_{\nu}^{\nu -\frac{1}{2}}$
   map the space $\mathring{C}_{\{0\}}^{\infty} \lr{E_{+}^{1}}$ into itself
   one-to-one.

Let us proof relations \eqref{1.1.9}. It suffices to prove only
one of them (e.\,g., the first one) because the other one is its
corollary.
    Since it is possible to apply the analytic continuation principle, it suffices
    to consider the case where $\Re \nu <\dfrac{1}{2}.$
 For $f \in \mathring{C}_{\{0\}}^{\infty}\lr{E_{+}^{1}}$, it follows from \eqref{1.1.10}
 that
\begin{multline}
B_{\nu} P_{\nu}^{\frac{1}{2}-\nu} f(y) =
 \frac{\sqrt{\pi}} {\Gamma \lr{\nu+1}\,
 \Gamma \lr{\frac{1}{2}-\nu}}  \lr{\frac{\pr^2}{\pr y^2}+ \frac{2 \nu+1}{y}
 \frac{\pr}{\pr y}}  \int\limits_y^{\infty} \lr{t^2-y^2}^{-\nu-\frac{1}{2}}  f(t) \, dt
 \\
= \frac{\sqrt{\pi}} {\Gamma \lr{\nu+1}\, \Gamma
\lr{\frac{1}{2}-\nu}}    \int\limits_1^{\infty}
\lr{t^2-1}^{-\nu-\frac{1}{2}} \lr{\frac{\pr^2}{\pr y^2}+ \frac{2
\nu+1}{y} \frac{\pr}{\pr y}} \lr{y^{- 2 \nu} f(t y)} \, dt.
\label{1.1.11}
\end{multline}
    Further, since
$$
\lr{\frac{\pr^2}{\pr y^2}+ \frac{2 \nu+1}{y} \frac{\pr}{\pr y}}
\lr{y^{- 2 \nu} f(t y)}  = y^{- 2 \nu} t^2 \left. \frac{\pr^2 f
(\tau)}{\pr \tau^2} \right|_{\tau = t y} + (1-2 \nu) y^{- 2 \nu -
1 } t \left. \frac{\pr f (\tau)}{\pr \tau} \right|_{\tau = t y}
$$
    and
$$
 (1-2 \nu) y^{- 2 \nu } \int\limits_1^{\infty}  t
\lr{t^2-1}^{-\nu-\frac{1}{2}} \left. \frac{\pr f (\tau)}{\pr \tau}
\right|_{\tau = t y} dt = -  y^{- 2 \nu - 1 }
\int\limits_1^{\infty} \frac{\pr \lr{t^2-1}^{\frac{1}{2}-\nu}}{\pr
t} \left. \frac{\pr f (\tau)}{\pr \tau} \right|_{\tau = t y} dt =
$$
$$
=- y^{- 2 \nu} \int\limits_1^{\infty}
\lr{t^2-1}^{-\nu-\frac{1}{2}} (t^2-1) \left. \frac{\pr^2 f
(\tau)}{\pr \tau^2} \right|_{\tau = t y} dt,
$$
    it follows that the right-hand
    side of relation \eqref{1.1.11} is reduced to the form
$$
\frac{\sqrt{\pi}} {\Gamma \lr{\nu+1}\, \Gamma
\lr{\frac{1}{2}-\nu}}    y^{- 2 \nu } \int\limits_1^{\infty}
\lr{t^2-1}^{-\nu-\frac{1}{2}}  \left. \frac{\pr^2 f (\tau)}{\pr
\tau^2} \right|_{\tau = t y} dt.
$$
    This completes the proof of relation \eqref{1.1.9} and,
    therefore, of the theorem.
\end{proof}
    Theorem \ref{theorem:1_1_1} is not new. It is proved by many authors under various
    restrictions (see \cite{SKM}). Its full proof is presented because it is hard
    to provide a particular reference and this theorem plays a key role
    in the theory of transmutation operators.

    For $\Re \mu > 0$, the Liouville operator $I^{\mu}$ is defined on functions
     $f \in
\mathring{C}_{\{0\}}^{\infty} \lr{E_{+}^{1}}$ as follows:
\begin{equation}
I^{\mu} f (y) = \frac{1}{\Gamma (\mu)} \int\limits_y^{\infty} (t-y)^{\mu-1} f(t) \, dt.
\label{1.1.12}
\end{equation}
 For other values of the complex parameter $\mu$, the function $I^{\mu} f (y)$
 is defined by means of the analytic continuation with respect to
 the parameter $\mu.$ If $\Re \mu > -M,$ where $M$ is a nonnegative integer,
 then this is equivalent to the following definition of $I^{\mu} f (y)$:
\begin{equation}
I^{\mu} f (y) := \frac{(-1)^M}{\Gamma \lr{\mu+M}}
\int\limits_y^{\infty} (t-y)^{M+\mu-1} \frac{\pr^M f(t)}{\pr t^M}
\, dt. \label{1.1.13}
\end{equation}
    Hence, the operators $I^{\mu}$ are defined on
$\mathring{C}_{\{0\}}^{\infty} \lr{E_{+}^{1}}$ for all complex
values of $\mu.$ They map the space
 $\mathring{C}_{\{0\}}^{\infty}\lr{E_{+}^{1}}$
 into itself
 one-to-one.
 The group property
\begin{equation}
I^{\mu} I^{\nu} = I^{\nu} I^{\mu} = I^{\mu+\nu} \label{1.1.14}
\end{equation}
    holds.

Now, for $\Re \nu \geq 0$, define the operators $P_{\nu}$ and
$S_{\nu}$ on functions
  $f \in \mathring{C}_{\{0\}}^{\infty}\lr{E_{+}^{1}}$ as follows:
\begin{equation}
P_{\nu} f = P_{\nu}^{\frac{1}{2} - \nu} I^{\nu -\frac{1}{2}} f \,
 \  \textrm{and} \ \, S_{\nu} f = I^{\frac{1}{2}-\nu} S_{\nu}^{ \nu - \frac{1}{2} }  f.
\label{1.1.15}
\end{equation}
    The above properties of Liouville operators and Theorem \eqref{theorem:1_1_1}
    imply the following assertion.

\begin{theorem} \label{theorem:1_1_2}
    If $\Re \nu \geq 0,$ then the operators $P_{\nu}$ and $S_{\nu}$ map
    the space $\mathring{C}_{\{0\}}^{\infty} \lr{E_{+}^{1}}$ into itself
    one-to-one
    and are mutually inverse.
    If $f \in
\mathring{C}_{\{0\}}^{\infty} \lr{E_{+}^{1}},$ then the relations
\begin{equation}
B_{\nu} P_{\nu} f  = P_{\nu} D^2 f
    \,
 \  \textrm{and} \ \,
D^2 S_{\nu} f  = S_{\nu} B_{\nu} f \label{1.1.16}
\end{equation}
    hold.
    \end{theorem}

Hence,  $P_{\nu}$ and $S_{\nu}$ are transmutation operators
indeed. They admit the following representations that are derived
 on the standard way.

Let $f \in \mathring{C}_{\{0\}}^{\infty} \lr{E_{+}^{1}}.$ Then,
for $0 \leq \Re \nu < \dfrac{1}{2}$, relations \eqref{1.1.15} and
\eqref{1.1.13} imply that
$$
    \begin{gathered}
P_{\nu} f (y) =  P_{\nu}^{\frac{1}{2} - \nu} I^{\nu -\frac{1}{2}}
f (y)= - P_{\nu}^{\frac{1}{2} - \nu} I^{\nu +\frac{1}{2}} D f (y)
    \\
=\frac{- \sqrt{\pi}} {\Gamma \lr{\nu+1}\, \Gamma
\lr{\frac{1}{2}-\nu} \,  \Gamma \lr{\nu+\frac{1}{2}} }
\int\limits_y^{\infty} \lr{t^2-y^2}^{-\nu-\frac{1}{2}}
\int\limits_t^{\infty} \lr{\tau-t}^{\nu-\frac{1}{2}}  D f(\tau)
d \tau dt
\\
=\frac{- \sqrt{\pi}} {\Gamma \lr{\nu+1}\, \Gamma
\lr{\frac{1}{2}-\nu} \,  \Gamma \lr{\nu+\frac{1}{2}} }
\int\limits_y^{\infty} D f(\tau) \int\limits_y^{\tau}
\lr{t^2-y^2}^{-\nu-\frac{1}{2}} \lr{\tau-t}^{\nu-\frac{1}{2}}  dt
d \tau.
    \end{gathered}
$$
    To compute the internal integral, introduce the new variable
     $z = \dfrac{t-y}{\tau-y}.$ Then
\begin{multline}
\int\limits_y^{\tau} \lr{t^2-y^2}^{-\nu-\frac{1}{2}}
\lr{\tau-t}^{\nu-\frac{1}{2}}  dt = (2 y)^{-\nu-\frac{1}{2}}
\int\limits_0^1 z^{-\nu-\frac{1}{2}} (1-z)^{\nu-\frac{1}{2}}
\lr{1-\frac{y - \tau}{2y}z}^{-\nu-\frac{1}{2}} dz = \\
=\frac{\Gamma \lr{\frac{1}{2}-\nu} \, \Gamma \lr{\frac{1}{2}+\nu}
} {(2y)^{\nu+\frac{1}{2}}}\, {_2F_1} \lr{\nu+\frac{1}{2};
\frac{1}{2}-\nu; l; \frac{y - \tau}{2y}}, \label{1.1.17}
\end{multline}
    where ${_2F_1 (a, b; c; \zeta)}$ denotes the Gauss hypergeometric function
    and the Euler representation
\begin{equation}
{_2F_1} (a, b; c; \zeta) =  \frac{\Gamma (c)}{\Gamma (b) \, \Gamma
(c-b)} \int\limits_0^1 t^{b-1} (1-t)^{c-b-1} (1 - \zeta t)^{-a} dt
\label{1.1.18}
\end{equation}
 (see \cite[p. 72]{BE1}) is used.

For positive values of $\zeta,$ the first-kind
 Legendre function  $P_{\mu}^0 (\zeta)$ can be defined by the relation
\begin{equation}
P_{\mu}^0 (\zeta) =  \frac{1}{\pi} \int\limits_0^{\pi}
\lr{\zeta+\sqrt{\zeta^2-1} \cos t}^{\mu} dt \label{1.1.19}
\end{equation}
 (see  \cite[p. 156]{BE1}).

Note that $P_{\mu}^0 (\zeta)$ is an analytic function of the
complex parameter $\mu.$ The  hypergeometric function from
relation \eqref{1.1.18} is expressed via the Legendre function as
follows (see \cite[p. 127]{BE1}):
$$
 {_2F_1} \lr{\nu+\frac{1}{2}; \frac{1}{2}-\nu; l; \frac{1}{2} - \frac{1}{2}
  \frac{\tau}{y}} = P^0_{\nu-\frac{1}{2}} \lr{\frac{\tau}{y}}.
$$
    Thus, the following representation of the operator $P_{\nu}$ is obtained:
\begin{equation}
P_{\nu} f (y)  = \frac{ -\sqrt{\pi}} {2^{\nu+\frac{1}{2}} \,
\Gamma \lr{\nu+1}} \, y^{-\nu-\frac{1}{2}}\,
\int\limits_y^{\infty} \frac{\pr f(\tau)}{\pr \tau}\,
P^0_{\nu-\frac{1}{2}}  \lr{\frac{\tau}{y}} d \tau. \label{1.1.20}
\end{equation}
    This relation is proved under the following additional restriction:
 $0 \leq \Re \nu < \dfrac{1}{2}.$
 However, if the function  $f$ and a positive value of $y$ are fixed, then
 both the
 left-hand and right-hand
 sides are analytic functions of parameter $\nu$ provided that $\Re
\nu \geq 0.$
    Therefore, by virtue of the uniqueness of the analytic continuation,
 relation \eqref{1.1.20} holds for all positive $y,$ all
$f \in \mathring{C}_{\{0\}}^{\infty} \lr{E_{+}^{1}}$, and all
$\nu$ such that $\Re \nu \geq 0.$

    The associated Legendre function $P_{\mu}^{-1} (\zeta)$ is
    defined by the relation
$$
\frac{\pr }{\pr \zeta} \, P_{\mu}^0 (\zeta) = \frac{\mu
(\mu+1)}{\sqrt{\zeta^2 -1}} \, P_{\mu}^{-1} (\zeta), \quad   \zeta
\geq 1.
$$
    Then, integrating \eqref{1.1.20} by parts, we see that
\begin{equation}
P_{\nu} f (y) = \frac{ \sqrt{\pi} \, y^{-\nu-\frac{1}{2}}}
{2^{\nu+\frac{1}{2}}  \, \Gamma \lr{\nu+1}}  \lr{ f(y)+\lr{\nu^2 -
\frac{1}{4}} \int\limits_1^{\infty} f (t y)
\frac{P_{\nu-\frac{1}{2}}^{-1} (t)}{\sqrt{t^2 -1}}\, dt}.
\label{1.1.21}
\end{equation}
    This presentation of the transmutation operator  $P_{\nu}$
    shows that it does not change the smoothness of the
    transformed function because  the kernel of the
    integral operator from \eqref{1.1.21} has no nonsummable  singularities
    provided that $1 \leq t < \infty.$
 This assertion holds only for positive values of the independent variable.
 For $y=0$, it is possible that the function $P_{\nu} f (y)$ has a singularity
 though the function $f$  has no singularities. The behavior of the kernel at infinity
 is important in this case. For example, $P_{\frac{1}{2}} f (y)=
\frac{1}{y} f(y).$

    Let us derive a similar representation for the operator $S_{\nu}.$
 If  $0\leq \Re \nu < \dfrac{1}{2}$ and
 $f \in\mathring{C}_{\{0\}}^{\infty} \lr{E_{+}^{1}}$, then relations \eqref{1.1.12}
 and \eqref{1.1.16} imply that
\begin{multline}
S_{\nu} f (y) =  I^{\frac{1}{2} - \nu} S_{\nu}^{\nu -\frac{1}{2}}
f (y) = \frac{- 2 \, \Gamma \lr{\nu+1}} {\sqrt{\pi} \, \Gamma
\lr{\frac{1}{2}+\nu} \, \Gamma \lr{\frac{1}{2}-\nu} } \frac{\pr
}{\pr y} \int\limits_y^{\infty}   \lr{t-y}^{\nu-\frac{1}{2}}
\int\limits_t^{\infty} \lr{\tau^2-t^2}^{-\nu-\frac{1}{2}} \tau f
(\tau) \,     d \tau dt  \\
 =\frac{- 2 \, \Gamma \lr{\nu+1}} {\sqrt{\pi}\, \Gamma
\lr{\frac{1}{2}+\nu} \, \Gamma \lr{\frac{1}{2}-\nu} } \frac{\pr
}{\pr y}    \int\limits_y^{\infty}  \tau f (\tau)
\int\limits_y^{\tau} \lr{t-y}^{\nu-\frac{1}{2}}
\lr{\tau^2-t^2}^{-\nu-\frac{1}{2}} dt d \tau. \label{1.1.22}
\end{multline}
    To compute the internal integral,
    introduce the new variable  $z = \dfrac{t-y}{\tau-y}.$
    This yields the relation
$$
    \begin{gathered}
 \int\limits_y^{\tau} \lr{t-y}^{\nu-\frac{1}{2}}
 \lr{\tau^2-t^2}^{-\nu-\frac{1}{2}} dt  = (\tau+
 y)^{\nu-\frac{1}{2}} \int\limits_0^1 z^{-\nu-\frac{1}{2}}
 (1-z)^{\nu-\frac{1}{2}} \lr{1-\frac{y -
 \tau}{y+\tau}z}^{\nu-\frac{1}{2}} dz
 \\
=\frac{\Gamma
 \lr{\frac{1}{2}-\nu} \, \Gamma \lr{\frac{1}{2}+\nu} }
 {(\tau+y)^{\nu-\frac{1}{2}}}\, {_2F_1} \lr{\frac{1}{2}-\nu;
 \frac{1}{2}-\nu; 1; \frac{y - \tau}{y+\tau}}.
    \end{gathered}
$$
The last hypergeometric function is reduced to the Legendre
function as follows
$$
{_2F_1} \lr{\frac{1}{2}-\nu; \frac{1}{2}-\nu; 1; \frac{\zeta -
1}{\zeta+1}}  = 2^{\nu -\frac{1}{2}} (\zeta+1)^{\frac{1}{2}-\nu}
P^0_{\nu-\frac{1}{2}} \lr{\zeta}
$$
    (see \cite[p. 128]{BE1}).
 Taking this into account, we derive the following relation from \eqref{1.1.22}:
\begin{equation}
S_{\nu} f (y) = \frac {-2^{\nu+\frac{1}{2}} \, \Gamma \lr{\nu+1}}{
\sqrt{\pi} }   \frac{\pr}{\pr y} \int\limits_y^{\infty}
\tau^{\nu+\frac{1}{2}} f (\tau)  P^0_{\nu-\frac{1}{2}}
\lr{\frac{y}{\tau}} \, d \tau. \label{1.1.23}
\end{equation}
    By virtue of the analytic continuation principle,
    the obtained representation of the operator $S_{\nu}$ holds provided that
     $\Re \nu \geq 0.$
 Differentiating the integral with respect to the parameter, we
 obtain that
$$
S_{\nu} f (y) = \frac {2^{\nu+\frac{1}{2}} \,  \Gamma \lr{\nu+1}}{
\sqrt{\pi} }\, y^{\nu+\frac{1}{2}} \lr{ f(y)-\lr{\nu^2 -
\frac{1}{4}} \int\limits_1^{\infty} t^{\nu+\frac{1}{2}} f (t y)
\frac{P_{\nu-\frac{1}{2}}^{-1} \lr{\frac{1}{t}}}{\sqrt{t^2 -1}}\,
dt}.
$$
    We see the operator  $S_{\nu}$ does not change the smoothness of functions as well
    provided that the values of the independent variable are positive.

Also, we use other transmutation operators and, respectively, use
another class of fractional integrals to construct them.
 Properties of such integrals are studied in the next section.

\subsection{Fractional integrals of the Riemann--Liouville
 type}\label{sec4.2}

    For $\Re \mu > 0$, define the operator $I^{\mu}_e$ on functions from the class
     $\mathring{C}_{\{0\}}^{\infty} \lr{E_{+}^{1}}$ as follows:
\begin{equation}
    I^{\mu}_e f(y) = \frac{1}{ \Gamma \lr{\mu}} \int\limits_y^{\infty} \lr{t-y}^{\mu-1}
     e^{y-t} f(t) \, dt, \quad   y>0.
    \label{1.2.1}
\end{equation}
    If  $\Re \mu > -M,$ where $M$ is a nonnegative integer, then we assign
\begin{equation}
    I^{\mu}_e f(y) = \frac{(-1)^M}{ \Gamma \lr{\mu+M}}\,
      e^{y} \int\limits_y^{\infty} \lr{t-y}^{\mu+M-1} \frac{\pr^M}{\pr t^M} \lr{e^{-t} f(t)}
       \, dt.
    \label{1.2.2}
\end{equation}
    If $\mathcal{E}$ denotes the multiplier by the function $e^y$ and $\mathcal{E}^{-1}$ denotes
    its inverse, then
\begin{equation}
I^{\mu}_e =  \mathcal{E} I^{\mu}  \mathcal{E}^{-1}. \label{1.2.3}
\end{equation}
    This relation between the operators $I^{\mu}_e$ and $I^{\mu}$ holds for all complex
    values of $\mu.$
    This and the property of Liouville operators imply the group property
    of operators $I^{\mu}_e$:
\begin{equation}
I^{\mu}_e  I^{\nu}_e = I^{\nu}_e I^{\mu}_e =  I^{\mu+\nu}_e.
\label{1.2.4}
\end{equation}
Hence, each operator $I^{\mu}_e$ maps the space
$\mathring{C}_{\{0\}}^{\infty} \lr{E_{+}^{1}}$ onto itself and its
inverse operator is $I^{-\mu}_e.$

The operators  $I^{\mu}_e$ and $I^{\nu}$ commute each other.
    This follows from their convolutional nature.
    In particular, $I^{\mu}_e$ commutes with the differentiating operator.

 Let us find explicit expressions for operators of the kind $I^{\nu} I^{\mu}_e.$
 First, assume that  $\Re \nu > 0$ and $\Re \mu > 0.$
    Then, for each  function  $f \in \mathring{C}_{\{0\}}^{\infty}\lr{E_{+}^{1}}$,
    we have the relation
\begin{multline}
I^{\nu} I^{\mu}_e f(y) = \frac{1} {\Gamma \lr{\nu}\, \Gamma \lr{\mu}}
\int\limits_y^{\infty} \lr{t-y}^{\nu-1} \int\limits_t^{\infty} (\tau-t)^{\mu}  e^{t -\tau}
 f(\tau) \, d \tau dt  \\
 =\frac{1} {\Gamma \lr{\nu}\, \Gamma \lr{\mu}} \int\limits_y^{\infty}   f(\tau)
 \lr{\tau-y}^{\nu+\mu-1} \int\limits_0^1 z^{\mu-1} (1-z)^{\nu-1}  e^{(y -\tau)z}  \, d z
 d \tau. \label{1.2.5}
\end{multline}
    The last internal integral is expressed via the degenerate
    hypergeometric function $\Phi \lr{a, c; \zeta}$
 (its another standard notation is ${_1F_1}\lr{a, c; \zeta}$).
  For $\Re c > \Re a > 0$, the function $\Phi$ can be defined by
  the relation
\begin{equation}
\Phi \lr{a, c; \zeta} = \frac{\Gamma (c)}{\Gamma (a)  \,  \Gamma
(c-a)} \int\limits_0^1 z^{a-1} (1-z)^{c-a-1}  e^{\zeta z}  \, d z.
\label{1.2.6}
\end{equation}
    Its properties provided below are well known (see \cite[Ch. 6]{BE1}).
The function $\dfrac{1}{\Gamma (c)} \Phi \lr{a, c; \zeta} $ is
analytically continued up to an entire function of its parameters
$a$ and $c$ and the variable $\zeta.$
    At the removable singular points $c= -m=0,-1,-2, \dots,$ we
    assign
\begin{equation}
\frac{1}{\Gamma (c)} \Phi \lr{a, c; \zeta} = a (a+1) \dots (a+m)
\,  \frac{\zeta^{m+1}}{(m+1)!}  \, \Phi \lr{a+m+1, m+2; \zeta}.
\label{1.2.7}
\end{equation}
    The following relations hold:
\begin{equation}
\frac{\pr}{\pr \zeta} \lr{\zeta^{c} \, \Phi \lr{a, c+1; \zeta}} =
c \, \zeta^{c-1} \Phi \lr{a, c; \zeta}, \label{1.2.8}
\end{equation}
\begin{equation}
\frac{\pr}{\pr \zeta}\,  \Phi \lr{a, c; \zeta}= \frac{a}{c}\, \Phi \lr{a+1, c+1; \zeta},
\label{1.2.9}
\end{equation}
    and
\begin{equation}
 \Phi \lr{a, c; \zeta}= \frac{\Gamma(c)}{\Gamma(c-a)} (- \zeta)^{-a} \lr{1+O (|\zeta|^{-1})},
  \  \zeta \to - \infty.
\label{1.2.10}
\end{equation}
  Getting back to relation \eqref{1.2.5}, we take into account
  \eqref{1.2.6} and find that
\begin{equation}
I^{\nu} I^{\mu}_e f(y)  = \frac{1}{\Gamma (\nu+\mu)}
\int\limits_y^{\infty} f (\tau) (\tau -y)^{\nu+\mu-1} \Phi
\lr{\mu, \mu+\nu; y - \tau} d \tau. \label{1.2.11}
\end{equation}
A similar representation of the operator  $I^{\mu}_e I^{-\mu} =
I^{-\mu} I^{\mu}_e$ for complex values of $\mu$ is derived as
follows. Using the previous relation for $\Re \mu > 0$ and $m=[\Re
\mu]+1$, we obtain the following relation provided that $f \in
\mathring{C}_{\{0\}}^{\infty} \lr{E_{+}^{1}}$:
\begin{equation}
I^{-\mu} I^{\mu}_e f(y)  = (-1)^m \frac{1}{\Gamma (m)} \int\limits_y^{\infty} \frac{\pr^m f (\tau)}{\pr \tau^m} (\tau -y)^{m-1} \Phi \lr{\mu, m; y - \tau} d \tau.
\label{1.2.12}
\end{equation}
    Integrate \eqref{1.2.12}  by parts  $m-1$ times and use relation
    \eqref{1.2.8}. This yields the relation
\begin{equation}
I^{-\mu} I^{\mu}_e f(y)  =  - \int\limits_y^{\infty} \frac{\pr f
(\tau)}{\pr \tau} \Phi \lr{\mu, 1; y - \tau} d \tau.
\label{1.2.13}
\end{equation}
Integrating \eqref{1.2.13}  by parts, we find that
\begin{equation}
I^{\mu} I^{\mu}_e f(y)  = f(y) \Phi \lr{\mu, 1; 0} +
\int\limits_y^{\infty} f (\tau)  \frac{\pr }{\pr \tau} \Phi
\lr{\mu, 1; y - \tau} d \tau. \label{1.2.14}
\end{equation}
    Taking into account that $\Phi \lr{\mu, 1; 0}=1$ and using relation \eqref{1.2.9},
    we reduce relation  \eqref{1.2.14} to the following final
    form:
\begin{equation}
I^{-\mu} I^{\mu}_e f(y)  =  f(y) - \mu \int\limits_y^{\infty} f
(\tau) \Phi \lr{\mu, 2; y - \tau} d \tau. \label{1.2.15}
\end{equation}
    This relation is proved here for $\Re \mu > 0;$ actually, it holds for
    all functions $f \in \mathring{C}_{\{0\}}^{\infty}
\lr{E_{+}^{1}}$ for all complex  $\mu$
    because it is possible to apply the analytic continuation.

Let us investigate the action of operators $I^{\mu}_e$ in Sobolev
spaces over the semiaxis; denote these spaces by $H^s
\lr{E_{+}^{1}},$ where $s$ is integer and $s \geq 0.$ These spaces
are defined as the closure of the set $\mathring{C}^{\infty}
\lr{\overline{E_{+}^{1}}}$ with respect to the norm
\begin{equation}
\| f \|_{H^s \lr{E_{+}^1}} = \lr{\sum\limits_{k=0}^{s} \|D^k f
\|_{L_{2} \lr{E_{+}^1}}^2}^{\frac{1}{2}}, \label{1.2.16}
\end{equation}
    where
$$
\| f \|_{L_{2} \lr{E_{+}^1}} = \lr{\int\limits_0^{\infty} |f(y)|^2  \,dy}^{\frac{1}{2}}.
$$

\begin{lemma} \label{lem:1.2.1}
    The operator $I^s_e$ is extended up to a bounded operator
    isomorphically mapping the space $L_2 \lr{E_{+}^{1}}$ onto $H^s \lr{E_{+}^{1}}.$
\end{lemma}

\begin{proof}
If $f \in \mathring{C}^{\infty} (\ov{E_{+}^{1}})$, then, using the
Leibnitz rule for $0 \leq k \leq s,$ we obtain the relation
\begin{multline}
D^k I^s_e f = D^k \mathcal{E} I^s \mathcal{E}^{-1} f
 = \sum\limits_{m=0}^k {k \choose m} \mathcal{E} D^m I^s \mathcal{E}^{-1} f  \\
=\sum\limits_{m=0}^k (-1)^m {k \choose m} \mathcal{E} I^{s-m}
\mathcal{E}^{-1} f = \sum\limits_{m=0}^k (-1)^m {k \choose m}
I_e^{s-m} f. \label{1.2.17}
\end{multline}
    If $\Re \mu > 0$, then, using the generalized Minkowski inequality, we obtain
    the estimate
$$
    \begin{gathered}
| I_e^{\mu} f \|_{L_{2} \lr{E_{+}^1}} = \frac{1}{|\Gamma(\mu)|}
\left\| \int\limits_0^{\infty} t^{\mu-1} e^{-t} f(y+t)  \,dt
\right\|_{L_{2} \lr{E_{+}^1}}
    \\
\leq \| f \|_{L_{2} \lr{E_{+}^1}} \, \frac{1}{| \Gamma(\mu)|} \,
\int\limits_0^{\infty} t^{\Re \mu-1} e^{-t} dt =
 \| f \|_{L_{2} \lr{E_{+}^1}} \, \frac{\Gamma \lr{\Re \mu}}{| \Gamma(\mu)|}
    \end{gathered}
$$
 preserved in the trivial case $\mu=0$ as well. This and \eqref{1.2.17}
 implies that
\begin{equation}
\| I_e^s f \|_{H^s \lr{E_{+}^1}} \leq c \| f \|_{L_{2} \lr{E_{+}^1}}.
\label{1.2.18}
\end{equation}
    To prove the inverse estimate, we apply the Leibnitz rule and see that
$$
I^{-s}_e f = \mathcal{E} I^{-s} \mathcal{E}^{-1} f= (-1)^s
\mathcal{E} D^s \mathcal{E}^{-1} f =(-1)^s \sum\limits_{k=0}^s
(-1)^k {s \choose m} D^{s-k} f
$$
    provided that $f\in \mathring{C}^{\infty} (\ov{E_{+}^{1}})$.
 Then, for each $f \in \mathring{C}^{\infty} (\ov{E_{+}^{1}})$, we obtain that
\begin{equation}
 \| f \|_{L_{2} \lr{E_{+}^1}} =
  \left\|I_e^{-s} I_e^{s} f \right\|_{L_{2} \lr{E_{+}^1}} \leq
   c \| I_e^s f \|_{H^s \lr{E_{+}^1}}.
\label{1.2.19}
\end{equation}
    Since $I^s_e \lr{ \mathring{C}^{\infty} (\ov{E_{+}^{1}})}=
\mathring{C}^{\infty} (\ov{E_{+}^{1}}),$ it follows that estimates
\eqref{1.2.18}-\eqref{1.2.19} complete the proof of the lemma.
\end{proof}

\begin{corollary}\label{cor:1.2.1}
 If $s \geq 0,$ then the norm
  $\| \tilde{f} \|_{H^s \lr{E_{+}^1}} =  \|I_e^{-s} f \|_{L_{2} \lr{E_{+}^1}}$
 is equivalent to norm \eqref{1.2.16} on the space  $H^s \lr{E_{+}^1}$.
\end{corollary}

\begin{lemma}\label{lem:1.2.2}
 Let $s \geq 0, s' \geq 0,$ and $\mu$ is a complex number. Then the operator
 $I_e^{\mu}$ continuously maps the space $H^s \lr{E_{+}^1}$ into $H^{s'}\lr{E_{+}^1}$
 provided that $s-s'+ \Re \mu > 0$. If  $\mu$ is real,  then the same holds
 for $s-s'+\mu=0$ as well.
\end{lemma}

\begin{proof}
Using the norm $\widetilde{\|\cdot  \|}_{H^s}$ and the group
property of operators $I_e^{\nu},$ we obtain  the relation
$$
\widetilde{| I_e^{\mu} f \|}_{H^{s'} \lr{E_{+}^1}} =   \| I_e^{\mu
- s'} f \|_{L_{2} \lr{E_{+}^1}} = \| I_e^{s- s'+\mu } I_e^{-s} f
\|_{L_{2} \lr{E_{+}^1}}.
$$
    It remains to note that the operator $I_e^{\nu}$ maps $L_{2} \lr{E_{+}^1}$
    into itself continuously provided that  $\Re \nu > 0$ or $\nu=0.$
    This completes the proof of the lemma.
\end{proof}
    The next result is used in the second chapter.

\begin{lemma}\label{lem:1.2.3}
 Let a function $a$ belong to $C^{\infty} (\ov{E_{+}^1})$ and be
 bounded together with all its derivatives. Then, for each
 nonnegative $s$ and each complex $\mu,$
  the operator $I_e^{- \mu} a I_e^{\mu}$ continuously maps the
  space $H^{s} \lr{E_{+}^1}$ into itself and the estimate
\begin{equation}
\| I_e^{-\mu} a I_e^{\mu}  f \|_{H^{s} \lr{E_{+}^1}} \leq  c  \|
f \|_{H^{s} \lr{E_{+}^1}}  \max\limits_{k \leq |[\Re \mu]|+s+2}\,
\sup\limits_{y>0} \left| \frac{\pr^k a(y)}{\pr y^k} \right|,
\label{1.2.20}
\end{equation}
    where the positive constant $c$ depends neither on $f$ nor on $a,$ holds.
\end{lemma}

\begin{proof}
    Using Corollary \ref{cor:1.2.1}, reduce the general case where $s \geq 0$
 to the case where $s=0.$ First, assume that $\Re \mu >0.$ Using the definition
 of operators $I_e^{\mu}$ and assuming that $m = [\Re \mu]+1,$ we
 obtain the following relation for each arbitrary function
  $f \in \mathring{C}^{\infty} (\ov{E_{+}^{1}})$:
$$
    \begin{gathered}
I_e^{-\mu} \lr{a I_e^{\mu} f} (y) = \frac{(-1)^m e^y}{\Gamma
(m-\mu) \, \Gamma(\mu)} \, D^m_y \int\limits_y^{\infty}
(t-y)^{m-\mu-1} a(t)\int\limits_t^{\infty} (\tau-t)^{\mu-1}
e^{-\tau} f (\tau) d \tau dt
    \\
=\frac{(-1)^m e^y}{\Gamma (m-\mu) \,
\Gamma(\mu)} \, D^m_y \int\limits_y^{\infty} e^{- \tau} f (\tau)
 (\tau-t)^{m-1} \int\limits_0^1 z^{m-\mu-1} (1-z)^{\mu-1} a (y + z (\tau-y))\, dz
  d \tau.
    \end{gathered}
$$
    Introduce the notation $a^{(k)} (t) = D^k a(t).$ Then
$$
    \begin{gathered}
I_e^{-\mu} a I_e^{\mu}  f (y) = a(y) f(y) - \frac{e^y}{ \Gamma(m -
\mu) \, \Gamma(\mu)}  \int\limits_y^{\infty} e^{- \tau}
 f (\tau)  \sum\limits_{k=1}^m  {m \choose k}   \frac{(m-1)!}{(k-1)!}
 \\
\times (y - \tau)^{k-1}  \int\limits_0^1 z^{m-\mu-1}
(1-z)^{k+\mu-1} a^{(k)}  (y + z (\tau-y))\, dz d \tau.
    \end{gathered}
$$
    Hence, the following estimate holds:
$$
\left| I_e^{-\mu} a I_e^{\mu}  f (y)  \right| \leq c
\sum\limits_{k=0}^m \sup\limits_{t>0} |a^{(k)} (t)| I^k_e \lr{|f|}
(y).
$$
    Applying Lemma \ref{lem:1.2.1}, we complete the proof for the case where
    $\Re \mu>0.$ If $\Re \mu < 0$ and $m = [- \Re \mu]+1,$ then
$$
    \begin{gathered}
I_e^{-\mu} a I_e^{\mu}  f (y) =  \frac{(-1)^m e^y}{ \Gamma(m +
\mu) \, \Gamma(-\mu)} \int\limits_y^{\infty}  (t - y)^{-\mu-1}
a(t) \int\limits_t^{\infty} (\tau -t)^{m+\mu-1}  D^m \lr{e^{-
\tau} f (\tau)} \,  d \tau d t
  \\
  =  \frac{(-1)^m e^y}{ \Gamma(m +
\mu) \, \Gamma(-\mu)} \int\limits_y^{\infty}  D^m \lr{e^{- \tau} f
(\tau)}  (\tau - y)^{m-1}   \int\limits_0^1 z^{-\mu-1}
(1-z)^{m+\mu-1} a (y + z (\tau-y))\, dz d \tau.
    \end{gathered}
$$
    As above, integrating by parts the  $m$th time, we see that estimate
    \eqref{1.2.20} hods for $\Re \mu <0$ as well. Finally, let $\Re \mu = 0.$
    Then
$$
    \begin{gathered}
I_e^{-\mu} a I_e^{\mu}  f (y)  =  \frac{e^y}{ \Gamma(1 + \mu) \,
\Gamma(1-\mu)} \,  D_y \int\limits_y^{\infty}   (t - y)^{- \mu}
\int\limits_t^{\infty} (\tau -t)^{\mu} D \lr{e^{-\tau} f (\tau)}
\, d \tau dt
    \\
=  \frac{e^y}{ \Gamma(1 + \mu) \, \Gamma(1-\mu)}  \, D_y
\int\limits_y^{\infty} D_{\tau}  \lr{e^{- \tau} f (\tau)}  (\tau -
y)  \int\limits_0^1 z^{-\mu} (1-z)^{\mu} a (y + z (\tau-y))\, dz d
\tau.
    \end{gathered}
$$
    The remaining reasoning is similar to the previous one.
    This completes the proof of the lemma.
\end{proof}

\begin{corollary}\label{cor:1.2.2}
 Under the assumptions of Lemma \ref{lem:1.2.3}, the following estimate holds
 for $\Re (\nu+\mu) > 0:$
    \begin{equation}
    \| I_e^{\nu} a I_e^{\mu}  f \|_{H^{s} \lr{E_{+}^1}} \leq
    c  \|  f \|_{H^{s} \lr{E_{+}^1}}
     \max\limits_{k \leq \min\limits \lr{|[\Re \nu]|, |[\Re \mu]|} +s+2}
      \sup\limits_{y>0} \left|D^k a(y) \right|.
    \label{1.2.21}
    \end{equation}
\end{corollary}

\begin{proof}
 It suffices to apply the relation
$$
I_e^{\nu} a I_e^{\mu}  = I_e^{\nu+\mu} I_e^{-\mu} a I_e^{\mu}
=I_e^{\nu} a I_e^{-\nu} I_e^{\nu+\mu}
$$
    and to use the boundedness of the operator $I_e^{\lambda}$
    under the assumption that
$\Re \lambda > 0.$
\end{proof}

\subsection{The relation to the Fourier and Hankel transformations}\label{sec4.3}

    Define the direct and inverse Fourier transformations by the
    relations
$$
F f (\eta) = \int\limits_{-\infty}^{\infty} f (y) e^{-i y \eta}
dy\, \ \textrm{and} \, \ F^{-1} g(y)  = \frac{1}{2 \pi}
\int\limits_{-\infty}^{\infty} g (\eta) e^{i y \eta} d \eta.
$$
    The Fourier cosine-transformation
$$
F_+ f (\eta) = \int\limits_{0}^{\infty} f (y) \cos (y \eta) \, dy
\, \ \textrm{and} \, \ F^{-1}_{+} g(y)  = \frac{2}{\pi}
\int\limits_{0}^{\infty} g (\eta) \cos (y \eta) \, d \eta
$$
    and the Fourier sine-transformation
$$
F_- f (\eta) = \int\limits_{0}^{\infty} f (y) \sin (y \eta) \,
dy\, \ \textrm{and} \, \ F^{-1}_{-} g(y)  = \frac{2}{\pi}
\int\limits_{0}^{\infty} g (\eta) \sin (y \eta) \, d \eta
$$
    are used as well.

 The direct and inverse Hankel transformations are
\begin{eqnarray}
F_{\nu} f (\eta) = \int\limits_{0}^{\infty} f (y) j_{\nu} (y \eta)
y^{2 \nu +1}\, dy \label{1.3.1}
\end{eqnarray}
    and
$$
F^{-1}_{\nu} g(y) =  \frac{1}{2^{2 \nu}\, \Gamma^2 (\nu+1)}
\int\limits_0^{\infty} g(\eta) j_{\nu} (y \eta) \eta^{2 \nu+1}\, d
\eta,
$$
    where $\nu \geq - \dfrac{1}{2}.$
 The normalized Bessel function $j_{\nu} (\lambda t)$
 is the solution of the following problem
$$
B_{\nu} f = - \lambda^2 f, \  f(0)=1, \  f'(0)=0.
$$
    It is related to the first-kind
     Bessel function $J_{\nu} (t)$ as follows:
\begin{eqnarray}
t^{\nu} j_{\nu} (t) = 2^{\nu} \Gamma (\nu+1) J_{\nu} (t).
\label{1.3.2}
\end{eqnarray}
 Let $\mbox{Cont}$ be an operator of the continuation of functions,
 acting from  $\mathring{C}^{\infty} (\ov{E_{+}^{1}})$ to
$\mathring{C}^{\infty} \lr{E^1}.$
    An example of such an operator is constructed in \cite{74}.
 The following (simpler) continuation operator can be proposed as well:
$$
\mbox{Cont} f(y) = \left\{
\begin{array}{ll}
 \chi (y) \int\limits_0^{\infty} f(-\lambda y) \psi (\lambda) \, d \lambda, & y<0, \\
 f(y), & y \geq 0,
\end{array}
\right.
$$
    where $\psi (\lambda) = \dfrac{1}{\pi} \int\limits_0^{\infty} \sin
(x \sqrt{\lambda}) \psi_1 (x) \, dx,$ $\psi_1 (x) = \mbox{sh} (x)
\chi (x),$ and $\chi$ is an arbitrary even function from
$\mathring{C}^{\infty} \lr{E^1}$ equal to one in a neighborhood of
the origin. The introduced function  $\psi$ is bounded and rapidly
decreases at infinity and $\int\limits_0^{\infty} \lambda^n \psi
(\lambda) \, d \lambda = (-1)^n.$ This easily implies that the
operator $\mbox{Cont}$ maps the space $\mathring{C}^{\infty}
(\ov{E_{+}^{1}})$ into $\mathring{C}^{\infty} \lr{E^1}.$

 First, assume that $0 \leq \Re \nu < \dfrac{1}{2}.$
 Then, for each $f \in \mathring{C}^{\infty}(\ov{E_{+}^{1}})$,
 it follows from \eqref{1.1.7} that
$$
    \begin{gathered}
P_{\nu}^{\frac{1}{2}-\nu} f(y) =  \frac{\sqrt{\pi} }{\Gamma
\lr{\nu+1}\, \Gamma \lr{\frac{1}{2}-\nu}} y^{- 2 \nu}
\int\limits_1^{\infty} \lr{t^2-1}^{-\nu-\frac{1}{2}} f(yt) \, dt
    \\
= \frac{y^{- 2 \nu}}{\sqrt{\pi}  \, \Gamma \lr{\nu+1}\, \Gamma
\lr{\frac{1}{2}-\nu}}  \int\limits_1^{\infty}
\lr{t^2-1}^{-\nu-\frac{1}{2}} \int\limits_{-\infty}^{\infty} e^{i
t \eta y} F \mbox{Cont} f (\eta) \, d \eta dt.
    \end{gathered}
$$
    From the restriction imposed on $\nu$, it follows that the integrand function
    is summable. Now, applying the Fubini theorem, we obtain that
$$
P_{\nu}^{\frac{1}{2}-\nu} f(y) =  \frac{y^{- 2 \nu}}{\sqrt{\pi} \,
\Gamma \lr{\nu+1}\,  \Gamma \lr{\frac{1}{2}-\nu}}
\int\limits_{-\infty}^{\infty} F \mbox{Cont} f (\eta)
\int\limits_1^{\infty} \lr{t^2-1}^{-\nu-\frac{1}{2}}  e^{i t \eta
y}  \, dt d \eta.
$$
    The internal integral is expressed via the
 first-kind  Hankel function $H_{\nu}^{(1)}$ (the third-kind
   Bessel function) by the relation
$$
\int\limits_1^{\infty} \lr{t^2-1}^{-\nu-\frac{1}{2}}  e^{i t y
\eta}  \, dt  = i \sqrt{\pi} \, 2^{- \nu -1} (y \eta)^{\nu}\,
\Gamma \lr{\frac{1}{2}-\nu} H^{(1)}_{\nu} (y \eta)
$$
 (see \cite[p. 95]{BE2}).
    Substituting this relation in the previous one, we see that
\begin{equation}
P_{\nu}^{\frac{1}{2}-\nu} f(y) =  \frac{i y^{-  \nu}}{2^{\nu+2} \,
\Gamma \lr{\nu+1}} \int\limits_{-\infty}^{\infty} H^{(1)}_{\nu} (y
\eta) F \mbox{Cont} f (\eta) \, d \eta. \label{1.3.3}
\end{equation}
    The right-hand
 side of this relation is an analytic function for $\Re \nu \geq 0$
    (see the representations of Hankel functions on \cite[p. 95]{BE2}).
Hence, relation \eqref{1.3.3} holds provided that $\Re \nu \geq
0.$

    From the way of the proof, we see that the value of the integral from the right-hand
 side of \eqref{1.3.3} does not depend on the choice of the continuation operator
 for $y>0$.
 This fact can be proved directly as well: one has to use relation \eqref{1.3.3}
 and the generalization of the Paley--Wiener
 theorem given in  \cite{23}.

The following assertion (used below) is a direct corollary from
the obtained relation.

\begin{lemma}\label{lem:1.3.1}
If $f \in \mathring{C}^{\infty} (\ov{E_{+}^{1}}),$ then
$$
\lim\limits_{y \to +0} y^{2 \nu} P_{\nu}^{\frac{1}{2}-\nu} f(y)  =
\frac{1}{2 \nu} f (0)
$$
 provided that $\Re\nu> 0$  and
$$
\lim\limits_{y \to +0} \frac{1}{\ln \frac{1}{y}} P_{0}^{\frac{1}{2}} f(y) = f (0).
$$
\end{lemma}

\begin{proof}
    This assertion easily follows from the asymptotic relations
$$
z^{\nu}  H^{(1)}_{\nu}  (z) = - \frac{i 2^{\nu} \,
\Gamma(\nu)}{\pi} + o(1), \quad   \Re \nu > 0,\quad  z \to 0,
$$
    and
$$
H^{(1)}_0  (z) = \frac{2 i}{\pi} \ln z + O(1), \quad   z \to 0,
$$
 for Hankel functions (see \cite{BE2}) and from the boundedness of
 functions $H^{(1)}_{\nu}  (z)$ for $z \geq \varepsilon > 0.$
\end{proof}
    Consider the case where $\nu $ is real and nonnegative. Let $f \in
\mathring{C}^{\infty}_{-} (\ov{E_{+}^{1}}).$  Then its odd
continuation to the whole line belongs to the space
$\mathring{C}^{\infty} \lr{E^1}.$
 Then take the odd continuation as the continuation operator in \eqref{1.3.3}.
 This yields the relation
$$
 \begin{gathered}
P_{\nu}^{\frac{1}{2}-\nu} f(y) =   \frac{y^{-  \nu}}{2^{\nu+1} \,
\Gamma \lr{\nu+1}} \int\limits_{-\infty}^{\infty} H^{(1)}_{\nu} (y
\eta) \eta^{\nu} F_{-} f (\eta) \, d \eta= \frac{y^{-
\nu}}{2^{\nu+1} \, \Gamma \lr{\nu+1}}
\lr{\int\limits_{-\infty}^0+\int\limits_0^{\infty}}
    \\
= \frac{y^{-  \nu}}{2^{\nu+1} \, \Gamma \lr{\nu+1}}
\int\limits_0^{\infty} \lr{H^{(1)}_{\nu} (y \eta)- e^{i \pi \nu}
H^{(1)}_{\nu} (-y \eta)} \eta^{\nu} F_{-} f (\eta) \, d \eta.
\end{gathered}
$$
    Further, since
$- e^{i \pi \nu} H^{(1)}_{\nu} (-z) = H^{(2)}_{\nu} (z) $
 (see \cite[p. 91]{BE2}) and
  $J_{\nu} (z) =\frac{H^{(1)}_{\nu}(z)+H^{(2)}_{\nu}(z)}{2} $
 (see \cite[p. 12]{BE2}), it follows that
$$
P_{\nu}^{\frac{1}{2}-\nu} f(y) =   \frac{1}{2^{\nu+1} \, \Gamma
\lr{\nu+1}}\, y^{-  \nu} \int\limits_0^{\infty} J_{\nu} (y \eta)
\eta^{\nu} F_{-} f (\eta) \, d \eta.
$$
    Changing the function $J_{\nu}$ for $j_{\nu}$ according to relation \eqref{1.3.2},
 we obtain that
\begin{equation}
P_{\nu}^{\frac{1}{2}-\nu} f(y) =  \frac{1}{2^{2 \nu}  \, \Gamma^2
\lr{\nu+1}}  \int\limits_{0}^{\infty} j_{\nu} (y \eta)  \eta^{2
\nu} F_{-} f (\eta) \, d \eta. \label{1.3.4}
\end{equation}

\begin{lemma}\label{lem:1.3.2}
 If $f \in \mathring{C}^{\infty}_{-} (\ov{E_{+}^{1}})$ and $\nu
\geq 0,$ then
\begin{equation}
P_{\nu}^{\frac{1}{2}-\nu} f =  F_{\nu}^{-1} \lr{\frac{1}{\eta} F_{-} f}.
\label{1.3.5}
\end{equation}
    If $f \in \mathring{C}^{\infty}_{+} (\ov{E_{+}^{1}}),$ then
\begin{equation}
S_{\nu}^{\nu-\frac{1}{2}} f =  F_{\nu}^{-1} \lr{\eta F_{\nu} f}.
\label{1.3.6}
\end{equation}
    The operator $P_{\nu}^{\frac{1}{2}-\nu}$ maps the space
     $\mathring{C}^{\infty}_{-} (\ov{E_{+}^{1}})$ onto
$\mathring{C}^{\infty}_{+} (\ov{E_{+}^{1}})$
 one-to-one.
  The operator $S_{\nu}^{\nu-\frac{1}{2}}$ is its inverse.
\end{lemma}

\begin{proof}
    Relation \eqref{1.3.5} is a brief form of relation \eqref{1.3.4}.
    To prove relation \eqref{1.3.6}, we invert the previous relation.
  Let $f \in \mathring{C}^{\infty}_{-} (\ov{E_{+}^{1}}).$
  Then $F_{-} f (\eta)$ is an odd rapidly decreasing function.
  Hence, $\dfrac{1}{\eta} F_{-} f (\eta)$ is an even smooth rapidly decreasing function.
  The operator $F_{\nu}$ is an automorphism of this function class.
  Hence,
   $F_{\nu}\Big(\dfrac{1}{\eta} F_{-} f\Big) \in \mathring{C}^{\infty}_{+}
(\ov{E_{+}^{1}})$
    because the compactness of its support follows from the representation
    of the operator by relation \eqref{1.1.7}.
    In the same way, it is proved that the operator $S_{\nu}^{\nu-\frac{1}{2}}$
    maps the space $\mathring{C}^{\infty}_{+} (\ov{E_{+}^{1}})$
    into $\mathring{C}^{\infty}_{-} (\ov{E_{+}^{1}}).$
    Therefore, these maps are surjective, which completes the proof of the lemma.
\end{proof}
    Note that there is another proof of the lemma.
    It is based on the fact that the operator  $F_{\nu}$ maps
     $\mathring{C}^{\infty}_{+}(\ov{E_{+}^{1}})$ onto the set of even entire functions
     of the exponential type (see \cite{23}).

    Let us illustrate the constructing of a new class of
    transmutation operators by means of analogs of relations \eqref{1.3.5}-\eqref{1.3.6}
    (see \cite{Kat1}). Let
$\nu \geq -\dfrac{1}{2}.$
    Assign
\begin{equation}\label{preCM}
P_{\nu, \pm}^{(\varphi)}  =   F_{\nu}^{-1} \lr{ \varphi \lr{\eta}
F_{\pm} }\, \  \textrm{and}\, \   S_{\nu, \pm}^{(\varphi)}  =
F_{\pm}^{-1} \lr{ \frac{1}{\varphi \lr{\eta}} F_{\nu} },
\end{equation}
    where  $\varphi (\eta)$ is a function.
    One can easily see that
$P_{\nu, \pm}^{(\varphi)}$ and $S_{\nu, \pm}^{(\varphi)}$ are
transmutation operators on suitable domains. The most  important
special cases are as follows.
    If  $\varphi (\eta) = \eta^{2\nu+1},$ then the operators $P_{\nu, +}^{(\varphi)}$
    and  $S_{\nu,+}^{(\varphi)}$ coincide with the classical Poisson and Sonin
    operators. If $\varphi (\eta) \equiv 1,$ then $P_{\nu,
+}^{(\varphi)}= P_{\nu}^{-\frac{1}{2}-\nu}$ and  $S_{\nu,
+}^{(\varphi)}= S_{\nu}^{\nu+\frac{1}{2}}.$ If $\varphi (\eta) =
\eta^{\nu+\frac{1}{2}},$ then
\begin{multline}
P_{\nu, +}^{(\varphi)} f(y) = F^{-1}_{\nu}
\lr{\eta^{\nu+\frac{1}{2}} F_+ f} (y)
\\
=\frac{-1}{2^{\nu+2}\,\Gamma (\nu+1)}
\left( \frac{\Gamma \lr{\nu+\frac{1}{2}}\,\cos (\pi \nu)}{2^{\nu}\, \Gamma (\nu+1)}
   \int\limits_y^{\infty} {_2F_1} \lr{\frac{\nu}{2}+\frac{3}{4}, \frac{\nu}{2}+
   \frac{1}{4}; \nu+1; \frac{y^2}{t^2}} D f(t) \, dt + {} \right.\\
\left. {}+ \frac{\sqrt{2} \, \Gamma
\lr{\frac{\nu}{2}+\frac{3}{4}}}{\Gamma
\lr{\frac{\nu}{2}+\frac{1}{4}}} \frac{1}{y} \int\limits_0^y
{_2F_1} \lr{\frac{\nu}{2}+\frac{3}{4}, \frac{3}{4}-\frac{\nu}{2};
\frac{3}{2}; \frac{t^2}{y^2}} D f(t) \, dt \right)
 \label{1.3.8}
\end{multline}
    and
\begin{multline*}
S_{\nu, +}^{(\varphi)} f(y) = F^{-1}_{+}
\lr{\eta^{-\nu-\frac{1}{2}} F_{\nu} f} (y)
\\
=\frac{2^{\nu-1}\,\Gamma (\nu+1)}{\pi} \frac{\pr}{\pr y} \left(
\frac{\Gamma \lr{\nu+\frac{1}{2}}\,\cos (\pi \nu)}{2^{\nu}\,
\Gamma (\nu+1)}  y^{-\nu-\frac{1}{2}} \int\limits_0^y {_2F_1}
\lr{\frac{\nu}{2}+\frac{3}{4}, \frac{\nu}{2}+\frac{1}{4}; \nu+1;
\frac{t^2}{y^2}} t^{2 \nu +1} f(t) \, dt{}\right.
\\
\left. {}+ \frac{\sqrt{2} \, \Gamma
\lr{\frac{\nu}{2}+\frac{3}{4}}}{\Gamma
\lr{\frac{\nu}{2}+\frac{1}{4}}} y \int\limits_y^{\infty} {_2F_1}
\lr{\frac{\nu}{2}+\frac{3}{4}, \frac{3}{4}-\frac{\nu}{2};
\frac{3}{2}; \frac{y^2}{t^2}} t^{\nu - \frac{1}{2}} f(t) \, dt
\right).
\end{multline*}
 In \cite{Kat1}, operators \eqref{1.3.8} are called \emph{isometric} because they
 isometrically map spaces  $L_2$ into $L_{2, \nu}.$
 In the same paper, their applications to pseudodifferential operators and
 to spectral theory are considered. In Chap. \ref{ch3}, it is
 shown that one can simplify the obtained representations,
 expressing the kernels of the transmutation operators via the Legendre function.
 Such operators corrected by a power-function
  factor can be reduced to unitary operators in one Lebesgue space
  over the semiaxis.
  In terms of Chap. \ref{ch3} (see below), they are combinations of Buschman--Erd\'elyi
  first-kind and second-kind
  transmutation operators.
  Their unitary property is especially useful in applications,
e.\,g., they bind solutions of perturbed differential equations
(i.\,e., equations with Bessel operators) and unperturbed ones
 (i.\,e., equations with second derivatives)  such that the norm is preserved.
 In \cite{SitDis}, it is proposed to use the terms
  \emph{So\-nin--Kat\-rakhov operators} and \emph{Poisson--Katrakhov
  operators} (see Chap. \ref{ch3}) to call this important class of operators
 introduced by Katrakhov.
 Note that the above Katrakhov construction given by \eqref{preCM} served as an example
 for Sitnik to  develop the general  composition method to
 construct transmutation operators of various classes (see Chap. \ref{ch6}).

Now, introduce the  transmutation operators  $P_{\nu, e}$ and
$S_{\nu, e}$ by the relations
\begin{equation}
P_{\nu,e} = P_{\nu}^{\frac{1}{2}-\nu} I_e^{\nu - \frac{1}{2}}\, \
\textrm{and}\, \ S_{\nu,e} = I_e^{\frac{1}{2}-\nu}
S_{\nu}^{\nu-\frac{1}{2}}. \label{1.3.9}
\end{equation}
    Also, they can be expressed in the form
\begin{equation}
    P_{\nu,e} = P_{\nu} \mathcal{J}_{\nu - \frac{1}{2}, e}\, \
\textrm{and}\, \  S_{\nu,e} = \mathcal{J}_{\frac{1}{2}-\nu, e}
S_{\nu},
    \label{1.3.10}
\end{equation}
    where $\mathcal{J}_{\mu, e} = I_e^{\mu} I^{-\mu}.$
    The operators $P_{\nu}^{\frac{1}{2}-\nu},$ $S_{\nu}^{\nu-\frac{1}{2}},$ and
    $P_{\nu}$ used above are defined in Sec. \ref{sec4.1}; the operators  $I_e^{\mu}$
 are defined in Sec. \ref{sec4.2}.

From the results of the mentioned sections, it follows that the
operators $P_{\nu,e}$ and $S_{\nu,e}$ map the space
$\mathring{C}^{\infty}_{\{0\}} \lr{{E_{+}^{1}}}$ onto itself
    one-to-one,
 are mutually inverse, and satisfy the relations
\begin{equation}
    B_{\nu} P_{\nu,e} = P_{\nu,e} D^2\, \
\textrm{and}\, \   D^2 S_{\nu,e} = S_{\nu, e} B_{\nu}.
    \label{1.3.11}
\end{equation}
    Thus,  $P_{\nu,e}$ and $S_{\nu,e}$ are   transmutation operators  as well.

    To clarify the relation between the operator $P_{\nu,e}$ and Fourier transformation,
    assume that $f \in \mathring{C}^{\infty} \lr{E^{1}}$ and $\mu$
    is a complex number. Then, for $\Re \mu > 0$, it follows from the Fubini
    theorem that
$$
F I_e^{\mu} f (\eta) = \frac{1}{\Gamma (\mu)}
\int\limits_{-\infty}^{\infty} e^{- i y \eta}
\int\limits_y^{\infty} (t-y)^{\mu -1 } e^{y-t} f(t) \, dt dy =
 \frac{1}{\Gamma (\mu)} \int\limits_{-\infty}^{\infty} f(t)
  \int\limits_{-\infty}^t  e^{- i y \eta} (t-y)^{\mu -1 } e^{y-t} \, dy dt.
$$
    Since
$$
\int\limits_{-\infty}^t  e^{- i y \eta} (t-y)^{\mu -1 } e^{y-t} \,
dy  =  e^{- i t \eta} (1 - i \eta)^{- \mu} \Gamma(\mu),
$$
     it follows that
\begin{equation}
F I_e^{\mu} f (\eta) = (1 - i \eta)^{- \mu} F f (\eta).
\label{1.3.12}
\end{equation}
    Now, let $\Re \mu \leq 0$ and $m = [-\Re \mu]+1.$ Then, using the same method,
    we obtain that
$$
    \begin{gathered}
F I_e^{\mu} f (\eta) = \frac{(-1)^m}{\Gamma(m+\mu)}
\int\limits_{-\infty}^{\infty} e^{- i t \eta} e^y
\int\limits_y^{\infty} (t-y)^{m+\mu -1 } D^m_t \lr{e^{-t} f(t)} \,
dt dy
    \\
= \frac{(-1)^m}{\Gamma(m+\mu)} \int\limits_{-\infty}^{\infty}
D^m_t \lr{e^{-t} f(t)}  \int\limits_{-\infty}^t  e^{- i y \eta}
(t-y)^{m+\mu -1 } e^{y} \, dy dt
    \\
= \frac{(-1)^m}{\Gamma(m+\mu)} \int\limits_{-\infty}^{\infty}
D^m_t \lr{e^{-t} f(t)}   e^{- i t \eta + t}  (1 - i \eta)^{-
(m+\mu)}  \Gamma(m+\mu) \, dt.
    \end{gathered}
$$
   Now, integrating by parts, we obtain relation \eqref{1.3.12} again.
   Hence, the said relation is proved for all complex values of $\mu.$
   Note that the following branch of the power function is selected throughout
   the whole monograph:
$$
(1 - i \eta)^{- \mu} = \exp \lr{\frac{\mu}{2} \ln (1+ \eta^2) - i
\mu \, \mbox{atan}\, \eta }.
$$
    Taking into account relations \eqref{1.3.12}-\eqref{1.3.13},
    we obtain the following representation of the operator   $P_{\nu,
e}$, valid under the assumption that  $\Re \mu \geq 0$:
\begin{equation}
P_{\nu, e} f(y) = \frac{i y^{- \nu}}{ 2^{\nu+2}\,  \Gamma (\nu+1)}
\int\limits_{-\infty}^{\infty} H^1_{\nu} (y \eta) \eta^{\nu} (1-i
\eta)^{\frac{1}{2}-\nu} F \mbox{Cont} f(\eta) \, d \eta,
 \label{1.3.13}
\end{equation}
    where $f \in \mathring{C}^{\infty}_{\{0\}} \lr{{E_{+}^{1}}}.$

\subsection{Transmutation operators and function spaces
(one-dimensional theory)}\label{sec4.4}

In this section, norm estimates of the transmutation operators
introduced above are considered. Here, we follow \cite{KatDis}.
    The applied method substantially using the  Mellin
    transformation and its multipliers to estimate norms of
 transmutation operators is proposed by Sitnik;
 main results of the present section belong to him as well.
    The Sitnik original approach based on the theory of
 Buschman--Erd\'elyi  transmutation operators and explicit usage of  Mellin multipliers
 is presented in  Chap. \ref{ch3}.
 Another variant to present these results (see \cite{KatDis})
 is given afterwards.
 In fact, there is no sense to separate the specified results because they are
 discussed by both authors (many times and comprehensively) during a number of years,
 while they were looking for suitable approaches and proving technique.
    Applications of the specified estimates of norms of transmutation operators
    to weighted boundary-value
    problems  (see Chap. \ref{ch4}) and  boundary-value
    problems with isolated singularities  (see Chap. \ref{ch5})
    are obtained entirely by Katrakhov.

 Let $L_{2, \nu} \lr{E_{+}^1},$ $\nu \geq -
\dfrac{1}{2},$ be the Hilbert space of functions $f(y),$ $y>0,$
such that the following norm is finite for each one:
\begin{equation}
\| f \|_{L_{2, \nu} \lr{E^1_+}} = \lr{\int\limits_0^{\infty}
|f(y)|^2 y^{2 \nu +1} \, dy}^{\frac{1}{2}}. \label{1.4.1b}
\end{equation}
    It is well known that the Hankel transformation $F_{\nu}$ maps
    $L_{2, \nu} \lr{E_{+}^1}$ onto itself and the Parseval
    relation
\begin{equation}
\|F_{\nu} f \|_{L_{2, \nu} \lr{E^1_+}} = 2^{\nu} \, \Gamma (\nu
+1) \| f \|_{L_{2, \nu} \lr{E^1_+}} \label{1.4.2b}
\end{equation}
 holds.

The function space  $H^s_{\nu,+} \lr{E_{+}^1},$ $s \geq 0,$ $\nu
\geq - \dfrac{1}{2},$ introduced in \cite{Kip1}, is defined as the
 closure of the function set  $\mathring{C}^{\infty}_{+} (\ov{E_{+}^{1}})$
 with respect to the norm
\begin{equation}
    \| f \|_{H_{ \nu, +}^s \lr{E^1_+}} =
    \frac{1}{2^{\nu}\, \Gamma (\nu +1)} \| (1+\eta^2)^{\frac{s}{2}}
     F_{\nu} f \|_{L_{2, \nu} \lr{E^1_+}}.
    \label{1.4.3}
\end{equation}
    Here, the evenness assumption is essential because the finiteness of norm
    \eqref{1.4.3} is not guaranteed otherwise.

Closing the set $\mathring{C}^{\infty}_{+} [0, R),$ $0<R<\infty,$
with respect to the same norm, we obtain the space
$\mathring{H}^{s}_{\nu, +} (0, R)$ continuously embedded into the
space $H^{s}_{\nu, +} \lr{E_{+}^1}.$ In the space $\mathring{H}^{s}_{\nu, +} (0, R)$,
 the expression
\begin{equation}
\| f \|_{\mathring{H}_{ \nu, +}^s (0, R)} =  \frac{1}{2^{\nu} \,
\Gamma (\nu +1)} \| \eta^{s} F_{\nu} f \|_{L_{2, \nu} \lr{E^1_+}}
\label{1.4.4}
\end{equation}
    is a norm equivalent to norm \eqref{1.4.3}. From the Parseval equality
    given by \eqref{1.4.2b}, it easily follows that if  $s$ is
    even and nonnegative, then
$$
\| f \|_{\mathring{H}_{ \nu, +}^s (0, R)} =
\|B_{\nu}^{\frac{s}{2}}  f \|_{L_{2, \nu} (0, R)}.
$$
  Define the Sobolev space $\mathring{H}^s (0, R),$ $s
 \geq 0,$ $0<R<\infty,$ as the closure of the set $\mathring{C}^{\infty} [0, R)$
 with respect to the norm
$$
\| f \|_{\mathring{H}^s (0, R)} =  \|D^s  f \|_{L_{2, \nu} (0, R)}.
$$

\begin{lemma}\label{lem:1.4.1}
    If $\nu \geq 0,$ then the operator  $S_{\nu}$ is extended up to a continuous
    bounded operator mapping the space $\mathring{H}_{\nu, +}^s (0, R)$ into
    $\mathring{H}^s (0, R)$ such that the estimate
\begin{multline}
\frac{2^{\nu+1}\, \Gamma (\nu +1)}{\sqrt{\pi}} \min\limits \lr{
\frac{1}{\sqrt{2}}, | \cos \frac{\pi (\nu+s)}{2}| }
\| f \|_{\mathring{H}_{ \nu, +}^s (0, R)}  \\
\leq \|S_{\nu} f \|_{\mathring{H}^s (0, R)} \leq \frac{2^{\nu+1}\,
\Gamma (\nu +1)}{\sqrt{\pi}}  \min\limits \lr{ \frac{1}{\sqrt{2}},
| \cos \frac{\pi (\nu+s)}{2}| } \| f \|_{\mathring{H}_{ \nu, +}^s
(0, R)} \label{1.4.5}
\end{multline}
    holds.
\end{lemma}

\begin{proof}
We consider only the case where $s=0$ because the arguing for the
case where $s>0$ is the same.
 First, we prove the validity of the relation
\begin{equation}
S_{\nu} f(y) = \frac{2}{\pi} \int\limits_0^{\infty} \cos \lr{y
\eta - \frac{\pi \nu}{2} - \frac{\pi}{4}} \eta^{\nu+ \frac{1}{2}}
F_{\nu} f (\eta) \, d \eta  \label{1.4.6}
\end{equation}
 for $\nu \geq 0$ and $f \in \mathring{C}^{\infty}_{+}
\lr{E_{+}^{1}}.$
    The following representation of the operator $S_{\nu}$ on this set is obtained
    earlier:
\begin{equation}
S_{\nu} f = I^{\frac{1}{2}-\nu} F^{-1}_{-} \eta F_{\nu} f.
\label{1.4.7}
\end{equation}
    This implies the validity of \eqref{1.4.6} for half-integer
    values of $\nu$ because the operator $I^{\frac{1}{2}-\nu}$ is
    the differentiation operator
 $(-1)^{\nu-\frac{1}{2}}D^{\nu-\frac{1}{2}}$
 in this case.
  Now, let  $n$ be a nonnegative even number and $n-\dfrac{1}{2}<\nu<n+\dfrac{1}{2}.$
 Since the function $\eta^{1+n} F_{\nu} f (\eta)$ is smooth and odd and rapidly
 decreases at infinity, then, due to the Fubini theorem, we have
 the  relations
\begin{multline}
S_{\nu} f(y) = I^{\frac{1}{2}-\nu} F^{-1}_{-} \eta F_{\nu} f (y) =
(-1)^n I^{\frac{1}{2}-\nu+n} D^n F^{-1}_{-} \eta F_{\nu} f (y)
\\
=(-1)^{\frac{n}{2}} I^{\frac{1}{2}-\nu+n} F^{-1}_{-} \eta^{1+n}
F_{\nu} f (y)
\\
= \frac{2 (-1)^{\frac{n}{2}}}{\pi \, \Gamma
\lr{\frac{1}{2}-\nu+n}} \int\limits_y^{\infty} (t-y)^{n - \nu
-\frac{1}{2}} \int\limits_0^{\infty} \sin (t \eta) \eta^{1+n}
F_{\nu} f (\eta) \, d \eta dt
\\
= \frac{2 (-1)^{\frac{n}{2}}}{\pi \, \Gamma
\lr{\frac{1}{2}-\nu+n}} \lim\limits_{A \to \infty} \int\limits_y^A
(t-y)^{n - \nu -\frac{1}{2}} \int\limits_0^{\infty} \sin (t \eta)
\eta^{1+n} F_{\nu} f (\eta) \, d \eta dt
\\
= \frac{2 (-1)^{\frac{n}{2}}}{\pi \, \Gamma
\lr{\frac{1}{2}-\nu+n}} \lim\limits_{A \to \infty}
\int\limits_0^{\infty} \eta^{1+n} f (\eta) \int\limits_y^A
(t-y)^{n - \nu -\frac{1}{2}} \sin (t \eta)  \, dt d \eta
\\
= \frac{2 (-1)^{\frac{n}{2}}}{\pi \, \Gamma
\lr{\frac{1}{2}-\nu+n}} \lim\limits_{A \to \infty}
\int\limits_0^{\infty} \eta^{1+n} f (\eta)
\lr{\int\limits_y^{\infty}-\int\limits_A^{\infty}} \, d \eta,
\label{1.4.8}
\end{multline}
 where two internal integral of the last expression are treated as
  improper integrals. The former one is equal to
\begin{equation}
\int\limits_y^{\infty} (t-y)^{n - \nu -\frac{1}{2}} \sin (t \eta)
\, dt = \eta^{ \nu-n -\frac{1}{2}} \Gamma \lr{n-\nu+\frac{1}{2}}
\cos \lr{y \eta - \frac{\pi (\nu-n)}{2} - \frac{\pi}{4}}.
\label{1.4.9}
\end{equation}
 By means of the integrating by parts, the latter one is estimated as follows:
\begin{multline}
\left| \int\limits_A^{\infty} (t-y)^{n - \nu -\frac{1}{2}} \sin (t
\eta) \, dt \right| = \frac{1}{\eta} \left| \int\limits_A^{\infty}
(t-y)^{n - \nu -\frac{1}{2}} \frac{\pr}{\pr t}\cos (t \eta) \, dt
\right|
\\
\leq  \frac{1}{\eta} \lr{(A-y)^{n - \nu -\frac{1}{2}} |\cos (A
\eta)| + \lr{\frac{1}{2}+\nu-n} \left| \int\limits_A^{\infty}
(t-y)^{n - \nu -\frac{3}{2}} \cos (t \eta) \, dt \right|}
\\
\leq \frac{1}{\eta} \lr{(A-y)^{n - \nu -\frac{1}{2}} +
\lr{\frac{1}{2}+\nu-n}  \int\limits_A^{\infty} (t-y)^{n - \nu
-\frac{3}{2}} \, dt } = \frac{2}{\eta} (A-y)^{n - \nu
-\frac{1}{2}}. \label{1.4.10}
\end{multline}
    Substituting \eqref{1.4.9}-\eqref{1.4.10}
    in \eqref{1.4.8}, we obtain relation \eqref{1.4.6}.
    Now, let $n$ be an odd positive number such that
$n-\dfrac{1}{2}<\nu<n+\dfrac{1}{2}.$
 Then, as above, we obtain that
$$
    \begin{gathered}
S_{\nu} f(y) =  (-1)^n I^{\frac{1}{2}-\nu+n} D^n F^{-1}_{-} \eta
F_{\nu} f (y) =  (-1)^{\frac{n+1}{2}} I^{\frac{1}{2}-\nu+n}
F^{-1}_{+} \eta^{1+n} F_{\nu} f (y)
    \\
= \frac{2 (-1)^{\frac{n+1}{2}}}{\pi \, \Gamma
\lr{\frac{1}{2}-\nu+n}} \int\limits_y^{\infty} (t-y)^{n - \nu
-\frac{1}{2}} \int\limits_0^{\infty} \cos (t \eta) \eta^{1+n}
F_{\nu} f (\eta) \, d \eta dt
    \\
= \frac{2 (-1)^{\frac{n+1}{2}}}{\pi \, \Gamma
\lr{\frac{1}{2}-\nu+n}}  \int\limits_0^{\infty} \eta^{1+n} F_{\nu}
f (\eta) \int\limits_y^{\infty} (t-y)^{n - \nu -\frac{1}{2}} \cos
(t \eta)  \, dt d \eta.
    \end{gathered}
$$
    Since
\begin{equation}
\int\limits_y^{\infty} (t-y)^{n - \nu -\frac{1}{2}} \cos (t \eta)
\, dt = - \eta^{ \nu-n -\frac{1}{2}} \Gamma \lr{\frac{1}{2}-\nu+n}
\sin \lr{y \eta - \frac{\pi (\nu-n)}{2} - \frac{\pi}{4}},
    \label{1.4.11}
\end{equation}
    it follows that the validity of \eqref{1.4.6} is proved for this case as well.

Consider the operator $A_{\nu}$ of the kind
\begin{equation}
A_{\nu} f(y) = \int\limits_0^{\infty}  \cos \lr{y \eta - \frac{\pi
\nu}{2} - \frac{\pi}{4}} f (\eta)\, d \eta, \label{1.4.12}
\end{equation}
 occurring in \eqref{1.4.6}.
 It is a linear combination of the Fourier trigonometric transformations.
 Therefore, it is bounded in $L_{2} \lr{E_{+}^1}.$
 Introduce the following modified  Mellin operator $M$:
\begin{equation}
    M g (p) =
     \frac{1}{\sqrt{2 \pi}} \int\limits_0^{\infty} y^{ i p - \frac{1}{2}}  g(y) \, d y,
     \  p \in E^1.
    \label{1.4.13}
\end{equation}
    It is easy to see that the operator $M$ isometrically maps the
    space  $L_{2} \lr{E_{+}^1}$ onto $L_{2} \lr{E_{+}^1}.$

    For simplicity, assume that $f \in \mathring{C}^{\infty}\lr{E_{+}^{1}}.$
  Then, using relations \eqref{1.4.9} and \eqref{1.4.11} and the arguing used
  to derive relation \eqref{1.4.6}, we obtain that
$$
    \begin{gathered}
M A_{\nu} f (p) =  \frac{1}{\sqrt{2 \pi}} \int\limits_0^{\infty}
y^{ i p - \frac{1}{2}}  \int\limits_0^{\infty}  \cos \lr{y \eta -
\frac{\pi \nu}{2} - \frac{\pi}{4}} f (\eta)\, d \eta dy
\\
=  \frac{1}{\sqrt{2 \pi}} \int\limits_0^{\infty}  f (\eta)
\int\limits_0^{\infty} y^{ i p - \frac{1}{2}}   \cos \lr{y \eta -
\frac{\pi \nu}{2} - \frac{\pi}{4}} \, d \eta dy
 \\
 = \frac{-1}{\sqrt{2 \pi}} \int\limits_0^{\infty}  f (\eta) \eta^{ i
p - \frac{1}{2}}   \Gamma \lr{i p +\frac{1}{2}} \sin \frac{\pi (i
p - \nu -1)}{2} \, d \eta
    \\
    =  -   \Gamma \lr{i p +\frac{1}{2}} \sin
\frac{\pi (i p - \nu -1)}{2}  M f (-p) \equiv a_{\nu} (p) M f
(-p).
    \end{gathered}
$$
    This and the isometry property of the operator $M$ imply the    following
 two-side estimate for operator \eqref{1.4.12}:
\begin{equation}
\inf\limits_{p \in E^1} |a_{\nu} (p)| \| f \|_{L_2 \lr{E^1_+}}
\leq \|A_{\nu} f \|_{L_2 \lr{E^1_+}} \leq \sup\limits_{p \in E^1}
|a_{\nu} (p)| \| f \|_{L_2 \lr{E^1_+}}, \label{1.4.14}
\end{equation}
    where the constants are exact. From the known relations
$$
 \left| \Gamma \lr{ip + \dfrac{1}{2}} \right|=\sqrt{\pi} /
\sqrt{\cosh (\pi p)}$$
 and
 $$\left|  \sin  \dfrac{\pi (i p - \nu
-1)}{2} \right|^2 = \lr{\sinh \dfrac{\pi p}{2}}^2 + \lr{\sin
\dfrac{\pi (\nu+1)}{2}}^2,
    $$
    we find that
$$
    \begin{gathered}
 \sup\limits_{p \in E^1} |a_{\nu} (p)| = \sqrt{\pi}   \sup\limits_{p \in E^1}
  \frac{\sqrt{\lr{\sinh \frac{\pi p}{2}}^2 +
   \lr{\cos  \frac{\pi \nu}{2}}^2}}{\sqrt{\cosh (\pi p)}} =
    \sqrt{\pi} \lr{\sup\limits \frac{\lr{\sinh \frac{\pi p}{2}}^2 +
     \lr{\cos  \frac{\pi \nu}{2}}^2}{1+2\lr{\sinh \frac{\pi p}{2}}^2} }^{\frac{1}{2}}
        \\
 =  \sqrt{\pi} \lr{\sup\limits_{t \geq 0}
 \frac{t + \lr{\cos  \frac{\pi \nu}{2}}^2}{1+2t} }^{\frac{1}{2}} =
  \sqrt{\pi} \lr{\max\limits \left\{ \left.\frac{t + \lr{\cos  \frac{\pi \nu}{2}}^2}
  {1+2t}\right|_{t=0},
  \left.\frac{t + \lr{\cos  \frac{\pi \nu}{2}}^2}{1+2t}\right|_{t=+ \infty}
   \right\}}^{\frac{1}{2}}
    \\
= \sqrt{\pi} \lr{\max\limits \left\{ \frac{1}{2},  \lr{\cos
\frac{\pi \nu}{2}}^2 \right\}}^{\frac{1}{2}} = \sqrt{\pi}
\max\limits \left\{ \frac{1}{\sqrt{2}},  |\cos  \frac{\pi \nu}{2}|
\right\}.
    \end{gathered}
$$
    In the same way, we obtain the relation
$$
 \inf\limits_{p \in E^1} |a_{\nu} (p)| =
  \sqrt{\pi} \min\limits \left\{ \frac{1}{\sqrt{2}},
  |\cos  \frac{\pi \nu}{2}|  \right\}.
$$
 Thus, estimate \eqref{1.4.14} takes the form
\begin{equation}
\sqrt{\pi} \min\limits \left\{ \frac{1}{\sqrt{2}}, |\cos \frac{\pi
\nu}{2}|  \right\} \| f \|_{L_2 \lr{E^1_+}}\leq \|A_{\nu} f
\|_{L_2 \lr{E^1_+}} \leq  \sqrt{\pi} \max\limits \left\{
\frac{1}{\sqrt{2}}, |\cos  \frac{\pi \nu}{2}|  \right\} \| f
\|_{L_2 \lr{E^1_+}}. \label{1.4.15}
\end{equation}
    Note that the operator $\sqrt{\dfrac{2}{\pi}} A_{\nu}$
is unitary in  $L_2 $ only for half-integer
 values of
$\nu$ (in this case, $A_{\nu}$ is reduced to the
cosine-transformation or sine-transformation)
  and its bounded inverse operator exists in $L_2 \lr{E_{+}^1}$ only for $\nu \neq
\pm 1, \pm 3, \dots$

    To complete the proof, it suffices to combine relation \eqref{1.4.6},
    estimate \eqref{1.4.15}, and the Parseval equality given by \eqref{1.4.2b}.
 \end{proof}
    In Chap. \ref{ch3}, more clear proof of this important lemma is provided.
  The original Katrakhov proof contained the wrong assertion about the
  isometry property for all values of the parameter.
  In \cite{S66, S6}, it is corrected.
 For the first time, the idea to apply the Mellin transformation technique
 to estimate norms of transmutation operators (frequently, the Slater theorem is used
 together with this technique) is proposed by Sitnik (see \cite{S66, S6}).
 The remarkable Katrakhov idea to ``correct'' Sonin and Poisson operators to ensure
  they to be bounded is a same space is implemented in the general form in Chap. \ref{ch3};
  the
  Buschman--Erd\'elyi operators of the zero-order
 smoothness, introduced in the specified chapter, are used for that purpose.

\begin{corollary} \label{cor:1.4.1}
 If $s=0$ and $R=\infty$, then the constants in estimate \eqref{1.4.5} are exact.
\end{corollary}

This assertion follows from the exactness of the constants in
estimate \eqref{1.4.15}.

\begin{corollary} \label{cor:1.4.2}
  For each positive half-integer
  $\nu,$ the operator $\sqrt{\pi} 2^{- \nu - \frac{1}{2}} \, \Gamma^{-1} (\nu+1) S_{\nu}$
   isometrically maps the space $L_{2, \nu} \lr{E_{+}^1}$ onto $L_{2} \lr{E_{+}^1}.$
\end{corollary}

\begin{proof}
    The isometry property follows from relation \eqref{1.4.5},
    where one can assign $R=\infty$ provided that $s=0$.
It follows from Lemma \ref{lem:1.3.2} that the image of the
operator $\sqrt{\pi} 2^{- \nu - \frac{1}{2}} \, \Gamma^{-1}
(\nu+1) S_{\nu}$ consists of all functions $f$ of the kind
 $f =I^{\frac{1}{2} - \nu} g,$ $g \in \mathring{C}_{-}^{\infty}
(\ov{E_{+}^{1}}).$
    This set is dense in the space $L_2\lr{E^1_+},$ which completes the proof.
\end{proof}
  Introduce one more function space $\mathring{H}_{\nu}^s (0,R).$
  Let $\mathring{C}_{\nu}^{\infty} (0, R)$ denote the set of all functions $f$
  admitting the representation $f =P_{\nu} g,$ where
   $g \in \mathring{C}_{\nu}^{\infty} [0, R).$
    It is clear that $\mathring{C}_{\nu}^{\infty} (0, R) \subset
\mathring{C}_{\{0\}}^{\infty} (0, R) \subset
\mathring{C}_{\{0\}}^{\infty} \lr{E_{+}^1}.$
 On $\mathring{C}_{\nu}^{\infty} (0, R)$, introduce the norm
\begin{equation}
\| f \|_{\mathring{H}_{ \nu}^s (0, R)} =  \frac{\sqrt{\pi}}
{2^{\nu+\frac{1}{2}}\, \Gamma (\nu +1)} \| S_{\nu} f
\|_{\mathring{H}_{ \nu}^s (0, R)}. \label{1.4.16}
\end{equation}
    This norm is well defined because
    $S_{\nu} f \in \mathring{C}_{\nu}^{\infty} [0, R)$
    by assumption.

 Let $\mathring{H}^s_{loc} (0, R)$ denote the set of all functions $f$ vanishing for
  $y \geq R$ and such that the seminorms
$$
P_{\varepsilon, s} (f) = \|D^s f \|_{L_2 (\varepsilon, \infty)}
$$
    are finite for each positive $\varepsilon.$
 The set   $\mathring{H}^s_{loc} (0, R)$ is a Fr\'echet set with respect to the topology
 generated by these seminorms.

The closure of the lineal $\mathring{C}_{\nu}^{\infty} (0, R)$ in
$\mathring{H}^s_{loc} (0, R)$ with respect to norm \eqref{1.4.16}
    is denoted by $\mathring{H}^s_{\nu} (0, R).$
    In Sec. \ref{sec3}, the embedding
$\mathring{H}^s_{\nu} (0, R) \subset \mathring{H}^s_{loc} (0, R)$
    is proved for more general case.

    Since $\mathring{C}^{\infty}_{+} [0, R) \subset
\mathring{C}^{\infty}_{\nu} (0, R)$ for each $\nu,$ it follows
from Lemma \ref{lem:1.4.1} that the space  $\mathring{H}^s_{\nu,
+} (0, R)$ is continuously embedded in $\mathring{H}^s_{\nu} (0,
R)$ and the induced and original norms of $\mathring{H}^s_{\nu, +}
(0, R)$ are equivalent to each other provided that  $s+\nu \neq
1,3,5, \dots$

    Introduce the weight function $\sigma_{\nu}$ by the following relation:
\begin{equation*}
\sigma_{\nu} (y) = \left\{
\begin{array}{lll}
y^{2 \nu}  & \mbox{if} \  \Re \nu>0, \\
\dfrac{1}{\ln y} & \mbox{if} \   \nu=0.
\end{array}
\right.
\end{equation*}
     The limit
$$
\sigma_{\nu} f |_{y=0} := \lim\limits_{y \to +0} \sigma_{\nu} (y)
f(y)
$$
    is called the \emph{weight boundary value} (\emph{weight}
$\sigma_{\nu}$-\emph{trace}) of the function $f$ at the point
$y=0$.

\begin{theorem} \label{teo:1.4.1}
    Let $\nu \geq 0,$ $s \geq 1,$ and $s+\nu>1.$
    Let $0<R<\infty.$ Then each function $f$ {\rm (}perhaps, corrected at a zero-measure
  set{\rm )} from the space $\mathring{H}^s_{\nu} (0, R)$ has a weight $\sigma_{\nu}$-trace
  at the point $y=0,$ the  inequalities
    \begin{equation}
    |\sigma_{\nu} f |_{y=0} | \leq \left\{
    \begin{array}{ll}
    \dfrac{2^{\nu-1} R^{s+\nu-1} \, \Gamma (\nu)}{ \sqrt{\pi} \, \Gamma \lr{s+\nu - \frac{1}{2}} \sqrt{s+\nu-1}}
    \|  f \|_{\mathring{H}_{ \nu}^s (0, R)}  & \mbox{if} \   \nu>0,  \\
    \dfrac{2 R^{s-1}}{ \sqrt{\pi} \, \Gamma \lr{s - \frac{1}{2}} \sqrt{s-1}}
    \|  f \|_{\mathring{H}_{0}^s (0, R)}  & \mbox{if} \   \nu=0
    \end{array}
    \right.
    \label{1.4.17}
    \end{equation}
hold, and the constants in these inequalities are exact.
\end{theorem}

\begin{proof}
    Regarding inequality \eqref{1.4.17} itself, it suffices to prove it for $f
\in\mathring{C}^{\infty}_{\nu} (0, R).$
 Due to Lemma \ref{lem:1.3.1}, each such function has a $\sigma_{\nu}$-trace
    and
\begin{equation}
\lim\limits_{y \to +0} \sigma_{\nu} (y) f(y) = \lim\limits_{y \to
+0} \sigma_{\nu} (y) P_{\nu}^{\frac{1}{2} - \nu}
S_{\nu}^{\nu-\frac{1}{2}} f(y)= \left\{
\begin{array}{ll}
\dfrac{1}{2 \nu} S_{\nu}^{\nu-\frac{1}{2}} f|_{y=0} & \mbox{if} \ \nu > 0, \\
 S_{0}^{-\frac{1}{2}} f|_{y=0}   & \mbox{if} \ \nu = 0.
\end{array}
\right. \label{1.4.18}
\end{equation}
    Let $f \in\mathring{C}^{\infty} [0, R).$ Then, due to the Cauchy--Bunyakovsky
  inequality, the following relation holds for    $\alpha > \dfrac{1}{2}$:
$$
    \begin{gathered}
|g(0)|= \left| I^{\alpha} I^{-\alpha} g|_{y=0} \right| =
\frac{1}{\Gamma (\alpha)} \left| \int\limits_0^R t^{\alpha -1}
I^{-\alpha} g(t) \, dt  \right|
    \\
\leq \frac{1}{\Gamma (\alpha)} \left( \int\limits_0^R t^{2\alpha
-2}  \, dt  \right)^{\frac{1}{2}}
 \left( \int\limits_0^R \left| I^{-\alpha} g \right|^2 dt  \right)^{\frac{1}{2}}
  = \frac{R^{\alpha - \frac{1}{2}}}{ \Gamma (\alpha) \sqrt {2 \alpha-1}}
   \left\| I^{-\alpha} g \right\|_{L_2 (0, R)}.
    \end{gathered}
$$
    In the last inequality, substitute the function
$S_{\nu}^{\nu-\frac{1}{2}} f = I^{\nu - \frac{1}{2}} S_{\nu} f $
 (obviously, it belongs to the space $\mathring{C}^{\infty} [0,R)$) for
 $g$; this yields the inequality
$$
\left|S_{\nu}^{\nu-\frac{1}{2}} f(0) \right| \leq  \frac{R^{\alpha
- \frac{1}{2}}}{ \Gamma (\alpha) \sqrt {2 \alpha-1}} \left\|
I^{\nu-\alpha-\frac{1}{2}} S_{\nu} f \right\|_{L_2 (0, R)}.
$$
Here, changing  $\alpha$ for $s+\nu - \dfrac{1}{2}$ and taking
into account relation \eqref{1.4.18}, we arrive at the inequality
\eqref{1.4.17}.

    To prove that its constants are exact, consider the function
$$
f_0 (y) = \left\{
\begin{array}{ll}
\int\limits_y^R (t-y)^{s+\nu - \frac{3}{2}} t^{s+\nu-\frac{3}{2}}
\, dt  & \mbox{if}
 \  0<y<R, \\
0   & \mbox{if} \ y \geq R.
\end{array}
\right.
$$
    To prove that the function  $I^{\frac{1}{2}-\nu} f_0$ belongs
    to $\mathring{H}^s (0, R),$ it suffices to show that the function
     $I^{\frac{1}{2}-s-\nu} f_0$ is bounded in a neighborhood of the point $y=R.$
 For $y<R$, we have the relation
$$
    \begin{gathered}
I^{\frac{1}{2}-s-\nu} f_0 (y) =
\frac{(-1)^{[\nu-\frac{1}{2}]+s+1}}{\Gamma \lr{\frac{3}{2} -s-\nu
+[s+\nu-\frac{1}{2}]}} D_y^{[\nu-\frac{1}{2}]+s+1} \int\limits_y^R
(t-y)^{\frac{1}{2}-\nu+[\nu - \frac{1}{2}]} \int\limits_y^R (\tau-t)^{s+\nu - \frac{3}{2}}
\tau^{s+\nu-\frac{3}{2}} \, d \tau dt
    \\
=\frac{(-1)^{[\nu-\frac{1}{2}]+s+1} \, \Gamma \lr{s+\nu -
\frac{1}{2}}}{\Gamma \lr{[\nu- \frac{1}{2}] +s+1}}  D_y^{[\nu-
\frac{1}{2}] +s+1 }  \int\limits_y^R \tau^{s+\nu-\frac{3}{2}}
(\tau -y)^{s+[\nu-\frac{1}{2}]} d \tau = \Gamma \lr{s+\nu
-\frac{1}{2} } y^{s+\nu-\frac{3}{2}}.
    \end{gathered}
$$
    Hence, $I^{\frac{1}{2}-\nu} f_0 \in \mathring{H}^s (0, R)$ and
\begin{equation}
\left\| I^{\frac{1}{2} - \nu} f_0 \right\|_{\mathring{H}^s \lr{0,
R}} = \frac{\Gamma \lr{s+\nu -\frac{1}{2} } }{\sqrt{2s+2 \nu - 2}}
R^{s+\nu-1}. \label{1.4.19}
\end{equation}
Then the function  $f_1 (y) := P_{\nu} I^{\frac{1}{2} - \nu} f_0$
belongs to $\mathring{H}^s_{\nu} \lr{0, R}$ and the following
relation holds for each positive $\nu$:
$$
\sigma_{\nu} f_1 |_{y=0} = \frac{1}{2 \nu} S^{\nu -
\frac{1}{2}}_{\nu} f_1 (0) = \frac{1}{2 \nu}  f (0) =
\frac{R^{2s+2\nu-2}}{4 \nu (s+\nu-1)}.
$$
    From \eqref{1.4.19}, we have the relation
$$
\left\| f_1 \right\|_{\mathring{H}^s_{\nu} \lr{0, R}} =
\frac{\sqrt{\pi} } {2^{\nu+\frac{1}{2}} \, \Gamma (\nu+1)} \left\|
I^{\frac{1}{2} - \nu} f_0 \right\|_{\mathring{H}^s_{\nu} \lr{0,
R}}  =  \frac{\sqrt{\pi} \,  \Gamma \lr{s+\nu -\frac{1}{2}}
R^{s+\nu-1}} {2^{\nu+1} \, \Gamma (\nu+1) \sqrt{s+\nu-1}}.
$$
    The last two relations show that the first inequality
 of \eqref{1.4.17} is an identity provided that $f_1 \in
\mathring{H}^s_{\nu} (0, R)$. The same is valid for the second
inequality as well, which completes the proof.
\end{proof}
  From \cite{BE1}, it is known that
\begin{equation}
\lim\limits_{\nu \to +\infty}  \frac{\Gamma (\nu + \alpha)}{\Gamma
(\nu + \beta)} \, \nu^{\beta-\alpha} =1. \label{1.4.20}
\end{equation}
    Therefore, the following assertion holds.

\begin{corollary} \label{cor:1.4.3}
    Under the assumptions of Theorem \ref{teo:1.4.1}, the estimate
\begin{equation}
\left| \sigma_{\nu} f|_{y=0}  \right| \leq c (s,  R) 2^{\nu}
R^{\nu} (\nu+1)^{-s} \left\| f \right\|_{\mathring{H}^s_{\nu}
\lr{0, R}} \label{1.4.21}
\end{equation}
    holds, where $\nu \geq 0$ and the constant  $c (s, R)$ depends only on  $s$ and $R.$
\end{corollary}

\begin{corollary} \label{cor:1.4.4}
    Let $\nu \geq 0,$ $s-2k- \dfrac{1}{2}>0,$ $s-2k+\nu>1,$
$k=0,1,\dots,$ and $0<R<\infty.$
 Then, for  each  $f \in
\mathring{H}^s_{\nu} (0, R)$    {\rm (}perhaps, corrected on a
zero-measure set{\rm )} there exists a $\sigma_{\nu}$-trace
    of the function $B_{\nu}^k f$ and the estimate
\begin{equation}
\left| \sigma_{\nu} B^k_{\nu} f|_{y=0}  \right| \leq c (s, k, R)
2^{\nu} R^{\nu}  (\nu+1)^{2k-s} \left\| f
\right\|_{\mathring{H}^s_{\nu} \lr{0, R}}, \label{1.4.22}
\end{equation}
    where the constant depends only on $s,$ $k,$ and $R,$ holds.
\end{corollary}

\section{Multidimensional Transmutation Operators}\label{sec5}

\subsection{Properties of Sobolev spaces}\label{sec5.1}

First, we provide several known results.

 Let $E^n$ denote the Euclidean  $n$-dimensional
 space of points $x = (x_1, \dots, x_n).$
 Let $\Theta$ be the unit sphere in $E^n.$
Introduce the spherical coordinates  $r \geq 0$ and $\vartheta \in
\Theta,$ where $r = |x|$ and $\vartheta = \dfrac{x}{|x|}.$
 In the angular coordinates $\varphi_1, \dots, \varphi_{n-1},$
 the vector $\vartheta = (\vartheta_1, \dots, \vartheta_n)$
 is expressed by the relations
\begin{eqnarray*}
\vartheta_1 = \cos \varphi_1, \\
\vartheta_2 = \sin \varphi_1 \cos \varphi_2, \\
\dots \dots \dots \dots \dots \\
\vartheta_{n-1} = \sin \varphi_1 \sin \varphi_2 \dots \sin \varphi_{n-2}
 \cos \varphi_{n-1}, \\
\vartheta_{n} = \sin \varphi_1 \sin \varphi_2 \dots \sin
\varphi_{n-2}  \sin \varphi_{n-1}.
\end{eqnarray*}
    In the spherical coordinates, the Laplace operator
$$
\Delta = \sum\limits_{j=1}^n \frac{\pr^2}{\pr x^2_j}
$$
 has the form
$$
\Delta = B_{\frac{n}{2}-1} + \frac{1}{r^2} \Delta_{\Theta},
$$
    where the  Bessel operator
    $B_{\nu} = \dfrac{\pr^2}{\pr r^2} + \frac{2\nu +1}{r} \frac{\pr}{\pr r},$
     $\nu = \dfrac{n}{2}-1,$ is called the \emph{radial part} of the operator $\Delta,$
     and the operator
$$
    \begin{gathered}
\Delta_{\Theta} = \frac{1}{\sin^{n-2} \varphi_1} \, \frac{\pr}{\pr
\varphi_1} \lr{\sin^{n-2}  \varphi_1  \, \frac{\pr}{\pr \varphi_1}
} + \frac{1}{\sin^{2} \varphi_1 \sin^{n-3} \varphi_2}  \,
\frac{\pr}{\pr \varphi_2} \lr{\sin^{n-3} \varphi_2 \,
\frac{\pr}{\pr \varphi_2} }+
    \\
\dots + \frac{1}{\sin^{2} \varphi_1 \sin^{2} \varphi_2 \dots
\sin^{2} \varphi_{n-2}}   \frac{\pr^2}{\pr \varphi^2_{n-1}}
    \end{gathered}
$$
 is called the \emph{angular part} of the operator $\Delta.$

  A \emph{spherical harmonic} of order  $k=0,1,2, \dots$ is
 a function $Y_k (\vartheta)$ satisfying the equation
\begin{equation}
\Delta_{\Theta} Y_k + k (n+k-2) Y_k=0
\label{2.1.1}
\end{equation}
    on the sphere $\Theta.$
    This equation has
    $d_k = \dfrac{(n+2k-2)(k+n-3)!}{k! (n-2)!}$
  solutions $Y_{k, l},$ $l=1,\dots, d_k,$ linearly independent and orthonormal in the sense
  of the space $L_2 (\Theta)$  of functions square-summable
  on the sphere $\Theta$.
    The system of functions $Y_{k, l},$ $k=0,1, \dots,$ $l=1, \dots, d_k,$
    forms an orthonormal base in $L_2 (\Theta).$ For a function $f$,
    coefficients of its expansion in a series with respect to spherical harmonics
 are determined by the relation
\begin{equation}
f_{k, l} (r) = \int\limits_{\Theta} f (r, \vartheta) Y_{k, l} (\vartheta) d \vartheta.
\label{2.1.2}
\end{equation}
    If $f \in \mathring{C}^{\infty} \lr{E^n}$, then the series
\begin{equation}
f (r) = \sum\limits_{k=0}^{\infty} \sum\limits_{l=1}^{d_k}  f_{k,
l} (r) Y_{k, l} (\vartheta) \label{2.1.3}
\end{equation}
 with respect to spherical harmonics converges to $f$ absolutely and uniformly.

Let $L_2 \lr{E^n},$ denote the space of functions $f(x)$ with the
finite norm
$$
\|f\|_{L_2 \lr{E^n}} = \lr{\int\limits_{E^n} |f(x)|^2 dx}^{\frac{1}{2}}.
$$
    As above, the set of functions with  the
finite norm
$$
\|f\|_{L_{2, \nu} \lr{E^1_+}} =  \lr{\int\limits_0^{\infty}
|f(r)|^2 r^{2 \nu+1} dr}^{\frac{1}{2}}
$$
 is denoted by $L_{2, \nu} \lr{E^1_{+}}.$

\begin{theorem}\label{teo:2.1.1}
    Let $f \in L_2 \lr{E^n}.$ Then series \eqref{2.1.3}
    {\rm(}with respect to spherical harmonics{\rm)}
 converges to $f$ with respect to the norm of the space $L_2\lr{E^n}$
 and
\begin{equation}
\|f\|_{L_2 \lr{E^n}}^2 = \sum\limits_{k=0}^{\infty}
\sum\limits_{l=1}^{d_k}  \|r^{-k} f_{k, l} (r) \|^2_{L_{2,
\frac{n}{2}+k-1} (E^1_+)}. \label{2.1.4}
\end{equation}
    Conversely, let functions $f_{k, l},$ $k=0,1, \dots,$ $l=1, \dots,
d_k,$ be such that $r^{-k} f_{k, l} \in L_{2, \frac{n}{2}+k-1}
\lr{E_{+}^1}$ and the series from the right-hand
    side of relation \eqref{2.1.4} converges.
    Then the series from \eqref{2.1.3} converges with respect to the norm of the space
 $L_2\lr{E^n},$ the sum $f $ of this series belongs to $ L_2 \lr{E^n},$ and its norm is
 computed according to relation \eqref{2.1.4}.
\end{theorem}

    The $n$-dimensional
    Fourier transformation $F$ acts on functions $f = g (r)Y_k (\vartheta)$
    as follows:
$$
F f(\xi) = \int\limits_{E^n} f(x) e^{-i \langle x, \xi \rangle } dx,
$$
    where $\xi = (\xi_1, \dots, \xi_n) \in E^n$ and
    $\langle x, \xi
\rangle=x_1 \xi_1 + \dots + x_n \xi_n$.
 Let $\rho = |\xi|$ and $\theta = \dfrac{\xi}{|\xi|}$ be spherical coordinates
 in the dual space $E^n.$ If $r^{-k} g(r) \in L_{2,\frac{n}{2}+k-1} \lr{E_{+}^1},$
 then
\begin{equation}
F \lr{g Y_k} (\rho, \theta) =  (-i)^k (2 \pi )^{\frac{n}{2}}
\rho^{\frac{2-n}{2}} Y_k (\theta) \int\limits_0^{\infty}
J_{\frac{n}{2}+k-1} (r \rho) g(r) r^{\frac{n}{2}} \, dr.
 \label{2.1.5}
\end{equation}
    Changing the first-kind
 Bessel function $J_{\nu}$ for the normalized function by means of the relation
  $$2^{\nu}\Gamma ({\nu+1}) J_{\nu} (z) = z^{\nu}j_{\nu} (z)$$
  and using the one-dimensional
   Hankel transformation $F_{\nu}$ defined by relation \eqref{1.3.1},
 one can represent relation \eqref{2.1.5} as follows:
\begin{equation}
F \lr{g Y_k} (\rho, \theta) =  \frac{ (-i)^k (2 \pi
)^{\frac{n}{2}}}{2^{\frac{n}{2}+k-1} \, \Gamma \lr{\frac{n}{2}+k}}
\rho^{k} Y_k (\theta) F_{k+\frac{n}{2}-1} \lr{r^{-k} g} (\rho).
\label{2.1.6}
\end{equation}
 All the above facts are well known (see, e.\,g., \cite{77,78}).
 We just formulate them in terms convenient for  our consideration.
    Relations \eqref{2.1.5}--\eqref{2.1.6}
    are obtained by Bochner; this is a special case of the known Funk--Hecke
    theorem (see \cite{BE2,SKM}).

Let $U_R \subset E^n$ be an open ball of a finite radius $R$
centered at the origin. As usual, let $\mathring{C}^{\infty}
\lr{U_R} = \{f: f \in \mathring{C}^{\infty} \lr{E^n}, \  \supp f
\subset U_R \}.$
  For nonnegative integers  $s$, define the space $\mathring{H}^s \lr{U_R}$ as the closure
  of the  function set  $\mathring{C}^{\infty} \lr{U_R}$ with  respect to the norm
\begin{equation}
\| f \|_{\mathring{H}^s \lr{U_R}} =  \lr{2 \pi}^{-\frac{n}{2}}
\lr{\int\limits_{E^n} |F f (\xi)|^2 |\xi|^{2s} \, d
\xi}^{\frac{1}{2}}, \label{2.1.7}
\end{equation}
    where $|\xi| = \lr{\xi_1^2+\dots+\xi_n^2}^{1/2}.$

Let $\mathring{T}^{\infty}_{+} \lr{U_R}$ denote the set of
 functions of the kind
\begin{equation}
f (r, \vartheta) = \sum\limits_{k=0}^{\mathcal{K}}
\sum\limits_{l=1}^{d_k} f_{k, l} (r) Y_{k, l}(\vartheta),
\label{2.1.8}
\end{equation}
 where $r^{-k } f_{k, l} \in \mathring{C}^{\infty}_{+} [0,R)$ and the positive integer
  $\mathcal{K}$ depends on the function $f.$

\begin{lemma} \label{lem:2.1.1}
  The set $\mathring{T}^{\infty}_{+} \lr{U_R}$ is dense in the space
     $\mathring{H}^s \lr{U_R}$ provided the  $s \geq 0.$
\end{lemma}

\begin{proof}
 Since the set of linear combinations of functions of the kind
\begin{equation}
f = \chi (r) Q(x), \label{2.1.9}
\end{equation}
    where $\chi (r)$ is an arbitrary function from the space
    $\mathring{C}^{\infty}_{+} [0, R)$ and $Q (x)$ are homogeneous polynomials of
     $n$ variables, is dense in $\mathring{H}^s \lr{U_R},$
     it suffices to show that each function  $f$ defined by
     relation \eqref{2.1.9} belongs to $\mathring{T}^{\infty}_{+} \lr{U_R}.$
 Denote the degree  of the polynomial $Q$ by $q$.
 Then, from the Gauss representation for homogeneous polynomials, we obtain the expansion
$$
Q(x) = Q_0 (x) +|x|^2 Q_1 (x) + \dots + |x|^{2 l} Q_l (x),
$$
 where $2 l \leq q$ and $Q_j$ are homogeneous harmonic polynomials of
 degree  $q-2j,$ $j=0, \dots, l.$
 From the homogeneity property, we conclude that
$$
Q_j (x) = r^{q-2j} Q_j( \vartheta),
$$
 where the functions $Q_j (\vartheta) = Y_{q-2j} (\vartheta)$ are spherical harmonics
 of order  $q-2j.$ Combining the last two relations, we obtain
 that
$$
f(x) = \chi (r) Q(x) = \chi (r) \sum\limits_{j=0}^l r^{2j}
r^{q-2j} Y_{q-2j} (\vartheta).
$$
    This expansion is equivalent to expansion \eqref{2.1.8}
    because
$r^{2 j} \chi(r) \in \mathring{C}^{\infty}_{+} [0, R)$ and
functions $Y_{q-2j} (\vartheta)$ can be represented by linear
combinations of orthonormal harmonics
 $Y_{q-2j, l},$ $l=1, \dots,d_{q-2j},$
 which completes the proof of the lemma.
\end{proof}

\begin{lemma}\label{lem5.2}
    Let  $r^{-k} f_{k, l}  \in\mathring{H}^s_{\frac{n}{2}+k-1} (0, R).$
    Then the function
$f=\sum\limits_{k=0}^{\mathcal{K}} \sum\limits_{l=1}^{d_k} f_{k,
l} (r) Y_{k, l} (\vartheta),$ $\mathcal{K} < \infty,$
    belongs to the space $\mathring{H}^s \lr{U_R}$ and the
    relation
\begin{equation}
\|f\|^2_{\mathring{H}^s \lr{U_R}} =\sum\limits_{k=0}^{\mathcal{K}}
\sum\limits_{l=1}^{d_k} \|r^{-k} f_{k, l}
\|^2_{\mathring{H}^s_{\frac{n}{2}+k-1} (0, R)} \label{2.1.10}
\end{equation}
    holds.
\end{lemma}

\begin{proof}
    In the spherical coordinates, relation \eqref{2.1.7} has the
    form
$$
\|f\|^2_{\mathring{H}^s \lr{U_R}} = \frac{1}{(2 \pi)^n}
\int\limits_0^{\infty} \int\limits_{\Theta} |F f(\rho,
\vartheta)|^2 d \vartheta \rho^{2s+n-1} \, d \rho.
$$
    Using this and \eqref{2.1.6}, we obtain the relation
$$
    \begin{gathered}
\|f\|^2_{\mathring{H}^s \lr{U_R}} = \frac{1}{(2 \pi)^n}
\int\limits_0^{\infty} \int\limits_{\Theta} \left|
\sum\limits_{k=0}^{\mathcal{K}} \sum\limits_{l=1}^{d_k}
\frac{(-i)^k (2 \pi)^{\frac{n}{2}}}{2^{\frac{n}{2}+k-1} \, \Gamma
\lr{\frac{n}{2}+k}}  \rho^k Y_{k, l} (\vartheta)
F_{\frac{n}{2}+k-1} (r^{-k} f_{k, l}) \right|^2 d \vartheta
\rho^{2s+n-1} d \rho
    \\
= \sum\limits_{k=0}^{\mathcal{K}} \sum\limits_{l=1}^{d_k}
\frac{1}{2^{n+2k-2} \, \Gamma \lr{\frac{n}{2}+k}}
\int\limits_0^{\infty} \left| F_{\frac{n}{2}+k-1} (r^{-k} f_{k,
l}) \right|^2  \rho^{2k+2s+n-1} d \rho
    \end{gathered}
$$
    due to the orthonormality of the system of spherical harmonics $Y_{k, l}$
    in the space $L_2 (\Theta).$ Combining the obtained relation
    with the definition of the norm of the space $\mathring{H}^s_{\nu, +} (0, R)$
     (see Sec. \ref{sec4.4}), we arrive at relation \eqref{2.1.10}, completing
     the proof of the lemma.
\end{proof}

\begin{theorem} \label{teo:2.1.2}
A function $f$ belongs to the space $\mathring{H}^s \lr{U_R},$ $s
\geq 0,$ $0 < R< \infty,$ if and only if the functions
 $r^{-k}f_{k, l} (r)$ belong to the spaces
$\mathring{H}^s_{\frac{n}{2}+k-1, +} (0, R)$
 and the numerical series
\begin{equation}
\sum\limits_{k=0}^{\infty} \sum\limits_{l=1}^{d_k} \|r^{-k} f_{k,
l} \|^2_{\mathring{H}^s_{\frac{n}{2}+k-1, +} (0, R)}
\label{2.1.11}
\end{equation}
    converges.
    If this is satisfied, then the function series given by \eqref{2.1.3} converges to
    the function $f$ with respect to the norm of the space $\mathring{H}^s \lr{U_R}$
    and
\begin{equation}
\|f\|^2_{\mathring{H}^s \lr{U_R}} =\sum\limits_{k=0}^{\infty}
\sum\limits_{l=1}^{d_k} \|r^{-k} f_{k, l}
\|^2_{\mathring{H}^s_{\frac{n}{2}+k-1, +} (0, R)}. \label{2.1.12}
\end{equation}
\end{theorem}

\begin{proof}
    \

{\it Necessity}.
 Let $f \in \mathring{H}^s \lr{U_R}.$
    Then, due to Theorem \ref{teo:2.1.1}, the functions $r^{-k} f_{k,
    l}$ belong to the space $L_{2, \frac{n}{2}+k-1, +} (0, R)$ and series \eqref{2.1.3}
    converges to $f$ with respect to the norm of the space $L_2 \lr{U_R}.$
  Applying the Fourier transformation to \eqref{2.1.3} and taking into account
  relation \eqref{2.1.6}, we obtain that
$$
F f (\rho, \vartheta) = \sum\limits_{k=0}^{\infty}
\sum\limits_{l=1}^{d_k} \frac{ (-i)^k (2 \pi
)^{\frac{n}{2}}}{2^{\frac{n}{2}+k-1} \, \Gamma \lr{\frac{n}{2}+k}}
\,  \rho^{k} \, Y_{k, l} (\vartheta) F_{k+\frac{n}{2}-1}
\lr{r^{-k} f_{k, l}}.
$$
    The series at the right-hand
    side converges to the function $F f \in L_2 (E^n)$ in $L_2 \lr{E^n}$.
    Assign
$$
\widetilde{\sum\limits}_{\mathcal{K}}=
\sum\limits_{k=0}^{\mathcal{K}} \sum\limits_{l=1}^{d_k} \frac{
(-i)^k (2 \pi )^{\frac{n}{2}}}{2^{\frac{n}{2}+k-1} \, \Gamma
\lr{\frac{n}{2}+k}} \, \rho^{k} \, Y_{k, l} (\theta)
F_{k+\frac{n}{2}-1} \lr{r^{-k} f_{k, l}}.
$$
    Then, for almost all positive $\rho$, the functions $F f (\rho, \vartheta)$
    and
$\widetilde{\sum\limits}_{\mathcal{K}} (\rho, \vartheta)$
    belong to the spaces $L_2 (\Theta).$
    Taking into account the orthonormality of $Y_{k, l},$ we
    conclude that
$$
    \begin{gathered}
0 \leq \| F f(p, \cdot) - \widetilde{\sum\limits}_{\mathcal{K}}
(p, \cdot)   \|_{L_2 (\Theta)} = \| F f(p, \cdot)   \|_{L_2
(\Theta)}
    \\
- \sum\limits_{k=0}^{\mathcal{K}} \sum\limits_{l=1}^{d_k} \frac{(2
\pi)^n}{ 2^{n+2k-2} \, \Gamma^2 \lr{\frac{n}{2}+k}} |\rho^k
F_{\frac{n}{2}+k-1} (r^{-k} f_{k, l})|^2.
    \end{gathered}
$$
 Multiply both sides of this inequality by the function $\rho^{2s+n-1}$ and integrate.
 Then
$$
    \begin{gathered}
\sum\limits_{k=0}^{\mathcal{K}} \sum\limits_{l=1}^{d_k} \frac{1}{
2^{n+2k-2} \, \Gamma^2 \lr{\frac{n}{2}+k}} \int\limits_0^{\infty}
\rho^{2k+2s+n-1} | F_{\frac{n}{2}+k-1} (r^{-k} f_{k, l})|^2 \, d
\rho
    \\
\leq \frac{1}{(2 \pi)^n} \int\limits_0^{\infty}
\int\limits_{\Theta} |F f(\rho, \vartheta)|^2 \rho^{2s+n-1} \, d
\rho d \vartheta = \| f \|^2_{\mathring{H}^s \lr{U_R}}.
    \end{gathered}
$$
    The right-hand
    side is bounded. Hence, each term of the left-hand
    side is bounded as well and the numerical series converges.
 In  particular, this means that the functions
  $r^{-k} f_{k, l}$ belong to
$\mathring{H}^s_{\frac{n}{2}+k-1, +} (0, R)$ and the sequence
$\left\{\widetilde{\sum\limits}_{\mathcal{K}}\right\}$ is
fundamental in $F \mathring{H}^s.$
 Since it converges to the function $F f$ in the sense of $L_2 (E^n),$ it follows that it
 converges in the sense of the space  $F \mathring{H}^s$ as well.
 Then the series from  \eqref{2.1.3} converges to $f$ in the space $\mathring{H}^s
(U_R),$ which leads to relation \eqref{2.1.12}.

{\it Sufficiency}. Let the functions $r^{-k} f_{k, l}$ belong to
$\mathring{H}^s_{\frac{n}{2}+k-1, +} (0, R)$ and the numerical
series given by \eqref{2.1.11} converge. Then the sequence
$$
{\sum\limits}_{\mathcal{K}}= \sum\limits_{k=0}^{\mathcal{K}}
\sum\limits_{l=1}^{d_k} f_{k, l} (r) Y_{k,l} (\vartheta)
$$
    is fundamental in  $\mathring{H}^s (0, R).$
    Thus, there exist $f \in \mathring{H}^s (0, R)$ satisfying relation \eqref{2.1.3},
    where the series converges with respect to the norm of the space $\mathring{H}^s (U_R).$
    This completes the proof.
\end{proof}
    Define the space $H^s (E^n)$ as the closure of the function set
    $\mathring{C}^{\infty} (E^n)$ with respect to the norm
$$
\| f \|^2_{H^s \lr{E^n}} =  \frac{1}{(2 \pi)^n} \int\limits_{E^n}
\| F f (\xi) \|^2 (1+ |\xi|^2)^s \, d \xi.
$$
    Then, similarly to Theorem \ref{teo:2.1.2}, the following
    assertion is proved.
\begin{theorem} \label{teo:2.1.3}
    A function  $f$ belongs to the space $H^s (E^n)$ if and only if the functions
     $r^{-k} f_{k, l}$ belong to the space $\mathring{H}^s_{\frac{n}{2}+k-1} (E_{+}^1)$
     and the series
\begin{equation}
\sum\limits_{k=0}^{\infty} \sum\limits_{l=1}^{d_k} \|r^{-k} f_{k,
l} \|^2_{H^s_{\frac{n}{2}+k-1, +} (E_{+}^1)} \label{2.1.13}
\end{equation}
 converges. If this is satisfied, then the function series given by \eqref{2.1.3} converges
 to the function $f$ with respect to the norm of the space $H^s (E^n)$ and
  $\|f\|^2_{H^s (E^n)}$ is equal to the sum of series \eqref{2.1.12}.
\end{theorem}

\subsection{Multidimensional transmutation operators: definitions}\label{sec5.2}

 Let $\mathring{T}^{\infty}_{\{0\}} (U_{R, 0}),$ where $U_{R, 0}=U_R \setminus 0,$
 denote the set of functions of the kind
\begin{equation}
f (r, \vartheta)= \sum\limits_{k=0}^{\mathcal{K}}
\sum\limits_{l=1}^{d_k} f_{k, l} (r) Y_{k,l} (\vartheta),
\label{2.2.1}
\end{equation}
    where $ f_{k, l} (r) \in \mathring{C}^{\infty}_{\{0\}} (0, R)$
    and $\mathcal{K}=\mathcal{K}(f)$ is a positive integer.
    Arbitrary singularities at the origin are admitted for functions
$f_{k, l}$ (and, therefore, for the function $f$).

 On the space $\mathring{T}^{\infty}_{\{0\}} (U_{R, 0})$ define the operator $\mathfrak{G}_n$
 by the relation
$$
\mathfrak{G}_n f (r, \vartheta)= \sum\limits_{k=0}^{\mathcal{K}}
\sum\limits_{l=1}^{d_k} \frac{\sqrt{\pi}
r^{\frac{1-n}{2}}}{2^{\frac{n}{2}+k-\frac{1}{2}} \, \Gamma
\lr{\frac{n}{2}+k}} S_{\frac{n}{2}+k-1} (r^{-k} f_{k, l}) Y_{k, l}
(\vartheta),
$$
    where $S_{\nu}$ are the transmutation operators introduced in Sec. \ref{sec4.1},
    while $f_{k, l}$ is the coefficient of the expansion of the function $f$ in the series
    with respect to spherical harmonics $Y_{k, l}.$

    Let us derive another representation of the operator $\mathfrak{G}_n$ such that
    no expansions with respect to spherical harmonics is used.
    Substituting \eqref{2.1.2} in \eqref{2.2.1}, we find that
$$
\mathfrak{G}_n f (r, \vartheta) = \int\limits_{\Theta}
\sum\limits_{k=0}^{\mathcal{K}} \frac{\sqrt{\pi}
r^{\frac{1-n}{2}}}{2^{\frac{n}{2}+k-\frac{1}{2}} \, \Gamma
\lr{\frac{n}{2}+k}} S_{\frac{n}{2}+k-1} (r^{-k} f (r, \vartheta'))
\sum\limits_{l=1}^{d_k} Y_{k,l} (\vartheta) Y_{k,l} (\vartheta')
\, d \vartheta'.
$$
    Due to the addition theorem for spherical harmonics (see \cite[p. 235]{BE2}),
    we obtain the relation
$$
\sum\limits_{l=1}^{d_k} Y_{k,l} (\vartheta) Y_{k,l} (\vartheta') =
\frac{\Gamma \lr{\frac{n}{2}-1}}{2 \pi^{\frac{n}{2}}}
\lr{\frac{n}{2}+k-1} C_k^{\frac{n}{2}-1} (\gamma),
$$
    where $\gamma = \langle \vartheta, \vartheta' \rangle$ denotes
    the scalar product of the vectors $\vartheta, \vartheta' \in
E^n.$ Since $|\vartheta|=|\vartheta'|=1,$ it follows that $\gamma$
is the cosine of the angle between them.
 Let $C_k^{\lambda}$ denote the ultraspherical Gegenbauer
 polynomials.
    From the representation of the operator $S_{\nu}$ by relation \eqref{1.1.23},
    we see that
\begin{multline}\mathfrak{G}_n f (r, \vartheta) = \frac{\Gamma
\lr{\frac{n}{2}{-}1}}{4 \pi^{\frac{n}{2}}}  \int\limits_{\Theta}
\sum\limits_{k=0}^{\mathcal{K}} \frac{ 2 \sqrt{\pi}
}{2^{\frac{n}{2}{+}k{-}\frac{1}{2}} \, \Gamma
\lr{\frac{n}{2}{+}k}} r^{\frac{1{-}n}{2}} S_{\frac{n}{2}{+}k{-}1}
(r^{{-}k} f (r, \vartheta'))
    \\
    \times \lr{\frac{n}{2}{+}k-1}
C_k^{\frac{n}{2}-1} (\gamma) \, d \vartheta'
\\
= \frac{ - \Gamma \lr{\frac{n}{2}-1}}{4 \pi^{\frac{n}{2}}}
r^{\frac{1-n}{2}}
 \\
 \times \frac{\pr}{\pr r}
 \left( r^{\frac{1+n}{2}}
\int\limits_{\Theta} \sum\limits_{k=0}^{\mathcal{K}} 2
\lr{\frac{n}{2}+k-1} C_k^{\frac{n}{2}-1} (\gamma) \int\limits_0^1
t^{- \frac{n+3}{2}} \cdot P^0_{\frac{n}{2}+k-\frac{3}{2}} (t) f
\lr{\frac{r}{t}, \vartheta'} \, dt d \vartheta' \right),
\label{2.2.3}
\end{multline}
 where $P^0_{\nu}$ is the first-kind
 Legendre function.

    Further arguing depends on the evenness or oddness of the
    dimension $n.$
    First, we assume that it is odd and $n \geq 3$.
    The following known relations for Legendre and  Gegenbauer
 polynomials are proved, e.\,g., in \cite[Chap. 10]{BE2}.

For positive integers $p$, the relation
\begin{equation}
D_{\gamma}^p C_m^{\lambda} (\gamma) =  \frac{2^p \, \Gamma
(\lambda+p)}{\Gamma (\lambda)} C_{m-p}^{\lambda+p} (\gamma)
\label{2.2.4}
\end{equation}
    holds. Assigning $\lambda = \dfrac{1}{2},$
    $m=k+\dfrac{n-3}{2},$ and
$p=\dfrac{n-3}{2}$ in this relation, we obtain that
$$
D_{\gamma}^{\frac{n-3}{2}}
C_{\frac{n}{2}+k-\frac{3}{2}}^{\frac{1}{2}} (\gamma)  =
2^{\frac{n-3}{2}} \frac{\Gamma \lr{\frac{n}{2}-1}}{\Gamma
\lr{\frac{1}{2}}} C_{k}^{\frac{n}{2}-1} (\gamma).
$$
    Since $C_{\nu}^{\frac{1}{2}} \equiv P_{\nu}^0,$ it follows
    that
$$
D^{\frac{n-3}{2}} P_{\frac{n}{2}+k-\frac{3}{2}}^{0} (\gamma) = 2^{\frac{n-3}{2}}
\frac{\Gamma \lr{\frac{n}{2}-1}}{\sqrt{\pi}} C_{k}^{\frac{n}{2}-1} (\gamma).
$$
 Now, substituting the last relation in \eqref{2.2.3}, we find
 that
\begin{multline}
\mathfrak{G}_n f (r, \vartheta) = - 2^{- \frac{1+n}{2}}
\pi^{\frac{1-n}{2}} r^{\frac{1-n}{2}} D_r \int\limits_{\Theta}
r^{\frac{1+n}{2}} D_{\gamma}^{\frac{n-3}{2}} \sum\limits_k 2
\lr{\frac{n}{2}+k-1}    P_{\frac{n}{2}+k-\frac{3}{2}}^{0} (\gamma)
\\
{}\times\int\limits_0^1 t^{-\frac{n+3}{2}}  P_{\frac{n}{2}+k-\frac{3}{2}}^{0}  (t) f
 \lr{\frac{r}{t}, \vartheta'} \, dt d \vartheta'
  \\
= - 2^{- \frac{1+n}{2}} \pi^{\frac{1-n}{2}} r^{\frac{1-n}{2}} D_r
\int\limits_{\Theta} r^{\frac{1+n}{2}}
D_{\gamma}^{r^{\frac{n-3}{2}}}
\sum\limits_{k=\frac{n-3}{2}}^{\infty} (2k+1) P_k^0 (\gamma)
\int\limits_0^1 t^{-\frac{n+3}{2}} P_{k}^{0}  (t) f
\lr{\frac{r}{t}, \vartheta'} \, dt d \vartheta' \label{2.2.5}
\end{multline}
 (after the change of the summation index).
 Since the Legendre polynomial $P_k^0$ has degree $k,$ it follows
 that $D_{\gamma}^{(n-3)/2} P_k^0 (\gamma) \equiv 0$ provided that
  $k<\dfrac{n-3}{2}.$
  Hence, the lower summation index in the last sum can be assigned to be equal to zero
  and the sum is not changed after this operation. Further, the function
   $t^{- \frac{n+3}{2}} f\lr{\frac{r}{t}, \vartheta'}$ of variable $t$ is infinitely
   differentiable for $0<t \leq 1$ and is identically equal to zero in a neighborhood
   of the left-hand
   edge.
   Assigning it to be equal to zero on the segment $[-1, 0],$ we
   obtain a function from the class $C^{\infty} [-1, 1].$
   For each function $g \in C^{\infty} [-1, 1]$, the following expansion
 in the Fourier series with respect to Legendre polynomials hold (see \cite[Chap. 10]{BE2}):
$$
g (\gamma) = \sum\limits_{k=0}^{\infty} (2k+1) P_{k}^{0}  (\gamma)
\int\limits_{-1}^1 g(t) P_{k}^{0}  (t) \, dt.
$$
 In our case, this yields the expression
$$
\sum\limits_{k=0}^{\infty} (2k+1) P_{k}^{0}  (\gamma)
\int\limits_{0}^1t^{-\frac{n+3}{2}}  P_{k}^{0}  (t) f
\lr{\frac{r}{t}, \vartheta'} \, dt =
 \left\{
\begin{array}{ll}
\gamma^{-\frac{n+3}{2}}   f \lr{\frac{r}{\gamma}, \vartheta'}  & \mbox{if} \  \gamma>0, \\
0  & \mbox{if} \  \gamma \leq 0.
\end{array}
\right.
$$
    Substituting it in \eqref{2.2.5}, we find the following (final) representation
    of the operator $\mathfrak{G}_n$ for the case where $n$ is odd:
\begin{equation}
\mathfrak{G}_n f (r, \vartheta)  = - 2^{- \frac{1+n}{2}}
\pi^{\frac{1-n}{2}} r^{\frac{1-n}{2}} D_r \int\limits_{\langle
\vartheta, \vartheta'\rangle > 0} r^{\frac{1+n}{2}}   \left.
D_{\gamma}^{\frac{n-3}{2}} \lr{\gamma^{-\frac{n+3}{2}}   f
\lr{\frac{r}{\gamma}, \vartheta'}} \right|_{\gamma=\langle
\vartheta, \vartheta'\rangle} d \vartheta'. \label{2.2.6}
\end{equation}
    The similar arguing for the case where $n$ is even is as
    follows. First, assume that $n\geq 4.$
    Again, use relation \eqref{2.2.4}, but assign $\lambda=1,$ $p = \dfrac{n}{2}-2,$ and
    $m = \dfrac{n}{2}-2+k$ in it. Then
$$
D_{\gamma}^{\frac{n}{2}-2} C^1_{\frac{n}{2}+k-2} (\gamma) =
2^{\frac{n}{2}-2} \, \Gamma \lr{\frac{n}{2}-1} C_k^{\frac{n}{2}-1}
(\gamma).
$$
    Since  $C^1_{\frac{n}{2}+k-2}$ is a second-kind
    Chebyshev polynomial, it follows that
$$
C^1_{\frac{n}{2}+k-2} (\gamma) =   \frac{\sin
\lr{\lr{\frac{n}{2}+k-1}} \arccos \gamma}{\sin (\arccos \gamma)}.
$$
  (see \cite[p. 185]{BE2}). Then
$$
C^1_{\frac{n}{2}+k-2} (\gamma) =   \frac{1}{\frac{n}{2}+k-1}
D_{\gamma} \cos \lr{\lr{\frac{n}{2}+k-1} \arccos \gamma}.
$$
    Hence,
$$
D_{\gamma}^{\frac{n}{2}-1} \cos \lr{\lr{\frac{n}{2}+k-1} \arccos
\gamma}  = 2^{\frac{n}{2}-2} \, \Gamma \lr{\frac{n}{2}-1}
\lr{\frac{n}{2}+k-1} C_k^{\frac{n}{2}-1} (\gamma).
$$
    Substituting this relation in \eqref{2.2.3}, we find that
\begin{multline}
 \mathfrak{G}_n f (r, \vartheta) = - \frac{1}{\sqrt{2}} (2
\pi)^{\frac{3-n}{2}} r^{\frac{1-n}{2}} D_r
\int\limits_{\Theta} r^{\frac{1+n}{2}} D_{\gamma}^{\frac{n}{2}-1}
 \int\limits_0^1 t^{- \frac{n+3}{2}} \times \\
 \times f \lr{\frac{r}{t}, \vartheta'} \sum\limits_{k=0}^{\infty} \cos
\lr{\lr{\frac{n}{2}+k-1} \arccos \gamma}
P_{\frac{n}{2}+k-\frac{3}{2}}^{0} (t)  \, dt d \vartheta'
 \\
 = - 2^{\frac{2-n}{2}} \pi^{\frac{3-n}{2}} r^{\frac{1-n}{2}} D_r \int\limits_{\Theta}
  r^{\frac{1+n}{2}} D_{\gamma}^{\frac{n}{2}-1}
  \int\limits_0^1 t^{- \frac{n+3}{2}}  f \lr{\frac{r}{t}, \vartheta'}
    \sum\limits_{k=\frac{n}{2} -1}^{\infty} \cos \lr{k \arccos
\gamma} P_{k-\frac{1}{2}}^{0} (t) \, dt d \vartheta'.
\label{2.2.7}
\end{multline}
    The first-kind
    Chebyshev polynomial $\cos \lr{k \arccos \gamma}$ has degree $k.$
 Then $D_{\gamma}^{\frac{n}{2}-1} \lr{\cos\lr{k \arccos \gamma}}$ provided that
  $0 \leq k < \dfrac{n}{2}-1$ and the following relation holds (see \cite[p. 166]{BE2}):
\begin{equation}
P_{-\frac{1}{2}}^{0} (t) + 2 \sum\limits_{k=1}^{\infty} \cos \lr{k
\arccos \gamma} P_{k-\frac{1}{2}}^{0} (t) = \left\{
\begin{array}{ll}
\dfrac{1}{ \sqrt{2 (\gamma-t)}}  & \mbox{if} \  \gamma>t,  \\
0 &  \mbox{if} \  -1<\gamma<t.
\end{array}
\right.
\label{2.2.8}
\end{equation}
    Taking into account the above, we find the final representation
    \begin{equation}
\mathfrak{G}_n f (r, \vartheta) =  - 2^{\frac{1+n}{2}}
\pi^{-\frac{n}{2}} r^{\frac{1-n}{2}} D_r \int\limits_{\langle
\vartheta, \vartheta'\rangle > 0} r^{\frac{1+n}{2}}
D_{\gamma}^{\frac{n}{2}-1}  \int\limits_0^{\gamma} t^{-
\frac{n+3}{2}}    \left. \frac{f \lr{\frac{r}{t},
\vartheta'}}{\sqrt{\gamma-t}} dt \right|_{\gamma=\langle
\vartheta, \vartheta'\rangle} d \vartheta'. \label{2.2.9}
\end{equation}
  The last relation is proved for even values of $n$ such that $n \geq 4.$
  The proposed method of the proof does not work for the case where $n=2$.
  However, the relation itself holds for this case as well.
  Really, in the polar coordinates $x_1 = r \cos \varphi $ and $x_2 = r \sin \varphi,$
  we have the relation
$$
f (r, \varphi) = f_0 (r) + \sum\limits_{k=1}^{\infty} \lr{f_{k, 1}
(r) \cos (k \varphi) + f_{k, 2} (r) \sin (k \varphi)},
$$
 where
$$
f_{k, 1} (r) =  \frac{1}{\pi} \int\limits_{- \pi}^{\pi} f (r,
\varphi) \cos (k \varphi) \, d \varphi, \ f_{k, 2} (r) =
\frac{1}{\pi} \int\limits_{- \pi}^{\pi} f (r, \varphi) \sin (k
\varphi) \, d \varphi, \ \textrm{and} \ f_{0} (r) = \frac{1}{2
\pi} \int\limits_{- \pi}^{\pi} f (r, \varphi)  \, d \varphi.
$$
    Therefore,
$$
    \begin{gathered}
\mathfrak{G}_2 f (r, \varphi) = \sqrt{\frac{\pi}{2 r}} \left(
\frac{1}{2 \pi}  \int\limits_{- \pi}^{\pi} S_0 f (r, \varphi') \,
d \varphi'  \right.
 \\
\left. + \sum\limits_{k=1}^{\infty} \frac{1}{\pi 2^k \, k!}
\int\limits_{- \pi}^{\pi} \lr{\cos (k \varphi)\cos (k \varphi')
+\sin (k \varphi) \sin (k \varphi')} S_k\lr {r^{-k} f (r,
\varphi')}  \, d \varphi' \right)
    \\
=  \frac{-1}{2 \pi \sqrt{r}} D_r \int\limits_{- \pi}^{\pi}
r^{\frac{3}{2}}  \int\limits_0^1 f \lr{\frac{r}{t}, \varphi'}
r^{-\frac{5}{2}} \lr{P_{-\frac{1}{2}}^0 (t) + 2
\sum\limits_{k=1}^{\infty}  \cos (k (\varphi- \varphi'))
P_{k-\frac{1}{2}}^0 (t)} d t d \varphi'.
    \end{gathered}
$$
    Taking into account \eqref{2.2.8}, we obtain the relation
\begin{equation}
\mathfrak{G}_2 f (r, \varphi) =   \frac{-1}{2^{\frac{3}{2}} \pi
\sqrt{r}} D_r \int\limits_{|\varphi-\varphi'|<\frac{\pi}{2}}
r^{\frac{3}{2}}  \int\limits_0^{\cos  (\varphi- \varphi')}
\frac{t^{-\frac{5}{2}} f \lr{\frac{r}{t}, \varphi'}}{\sqrt{\cos
(\varphi- \varphi')-t}} \, d t d \varphi' \label{2.2.10}
\end{equation}
    completing the proof of relation \eqref{2.2.9}  for the case where $n=2.$

Note that the proof of relations \eqref{2.2.6}, \eqref{2.2.9}, and
\eqref{2.2.10} is not complete because the change of finite sums
by series is not justified. The said justification uses standard
properties of expansions of functions with respect to Legendre and
Chebyshev polynomials; here, we omit this part of the proof.

Define the operator $\mathfrak{B}_n$ on the set
$\mathring{T}^{\infty}_{\{0\}} \lr{U_{R, 0}}$ as follows:
\begin{equation}
\mathfrak{B}_n f \lr{r, \vartheta} =
\sum\limits_{k=0}^{\mathcal{K}} \sum\limits_{l=1}^{d_k}
\frac{2^{\frac{n}{2}+k-\frac{1}{2}} \, \Gamma \lr{\frac{n}{2}+k}
}{\sqrt{\pi}} r^k P_{\frac{n}{2}+k-\frac{1}{2}}
\lr{r^{\frac{n-1}{2}} f_{k, l}} Y_{k, l} \lr{\vartheta},
\label{2.2.11}
\end{equation}
    where $f_{k, l}$ are defined by \eqref{2.1.2}, while $P_{\nu}$
    is the transmutation operator from Sec. \ref{sec4.1}.
    No representation of kind \eqref{2.2.6} or \eqref{2.2.9} can
    be obtained for the operator $\mathfrak{B}_n$ because the operator $\mathfrak{G}_n$
removes singularities of functions at the origin, while the
operator $\mathfrak{B}_n$ generates them. For example, the
operator $\mathfrak{B}_n$ can map a function without singularities
into
 a function with an essential singularity at the origin.

\begin{theorem} \label{teo:2.2.1}
    The operators $\mathfrak{G}_n$ and $\mathfrak{B}_n$ map the
    space $\mathring{T}^{\infty}_{\{0\}} \lr{U_{R, 0}}$ into
itself and are mutually inverse. If $f \in
\mathring{T}^{\infty}_{\{0\}} \lr{U_{R, 0}},$ then
\begin{equation}
\mathfrak{G}_n  \Delta f = r^{\frac{1-n}{2}} D^2_r
\lr{r^{\frac{n-1}{2}} \mathfrak{G}_n f} \quad \textrm{and} \quad
\Delta \mathfrak{B}_n f = \mathfrak{B}_n \lr{ r^{\frac{1-n}{2}}
D^2_r \lr{r^{\frac{n-1}{2}}  f}}. \label{2.2.12}
\end{equation}
\end{theorem}

\begin{proof}
    Since the spherical harmonics $Y_{k, l}$ are eigenfunctions of the operator
     $\Delta_{\Theta},$ it follows that
$$
    \begin{gathered}
\Delta f =  \lr{D^2_r + \frac{n-1}{r} D_r + \frac{1}{r^2}
\Delta_{\Theta}} \sum\limits_{k=0}^{\mathcal{K}}
\sum\limits_{l=1}^{d_k}  f_{k, l} (r) Y_{k, l} \lr{\vartheta}
    \\
= \sum\limits_{k=0}^{\mathcal{K}} \sum\limits_{l=1}^{d_k}
\lr{D^2_r + \frac{n-1}{r} D_r - \frac{k(n+k-2)}{r^2} } f_{k, l}
(r) Y_{k, l} \lr{\vartheta} =
 \sum\limits_{k=0}^{\mathcal{K}} \sum\limits_{l=1}^{d_k} r^k
B_{\frac{n}{2}+k-1}  \lr{r^{-k} f_{k, l} } Y_{k, l}
\lr{\vartheta},
    \end{gathered}
$$
    where $B_{\nu}$ denotes the Bessel operator acting with respect to the radial variable.
 Substituting this expression in \eqref{2.2.1}, we obtain that
$$
\mathfrak{G} \Delta f = \sum\limits_{k=0}^{\mathcal{K}}
\sum\limits_{l=1}^{d_k} \frac{\sqrt{\pi}
r^{\frac{1-n}{2}}}{2^{\frac{n}{2}+k-1} \, \Gamma
\lr{\frac{n}{2}+k}} S_{\frac{n}{2}+k-1} B_{\frac{n}{2}+k-1}
\lr{r^{-k} f_{k, l} } Y_{k, l} \lr{\vartheta}.
$$
    Due to Theorem \ref{theorem:1_1_2}, we have the relation $S_{\nu} B_{\nu} = D_r^2
S_{\nu}$ leading to  relation  \eqref{2.2.12}.
 The second relation follows from the first one, which completes the proof of the theorem.
\end{proof}
  This theorem allows one to treat  $\mathfrak{G}_n$ and
  $\mathfrak{B}_n$ as transmutation operators transforming the
  multidimensional Laplace operator  $\Delta$ into an ordinary differential operator
  of the kind
$r^{\frac{1-n}{2}} D^2_r r^{\frac{n-1}{2}}.$
    The operators $r^{\frac{n-1}{2}} \mathfrak{G}_n$ and $\mathfrak{B}_n
r^{\frac{1-n}{2}}$ transform  $\Delta$ into $D^2_r.$

Note that integral representations for the multidimensional
transmutation operators $\mathfrak{G}_n$ and $\mathfrak{B}_n$ are
obtained in the Sitnik course work in 1981  (the problem is set by
Katrakhov). No special functions are used in this work. Instead,
combinatorial properties are used, which allows one to obtain a
strict and relatively simple proof of integral representations
 \eqref{2.2.6} and \eqref{2.2.9} and, therefore, of Theorem \ref{teo:2.2.1}
 (under the least assumptions about  the smoothness of the function).
 The Katrakhov proof presented here is more conceptual: it demonstrates relations
 to constructions of the theory of transmutation operators and special functions.
 However, as we note above, no completely strict proof is accessible and excessive
 requirements for the smoothness of functions are needed.

\subsection{$L_2$-theory of multidimensional transmutation operators}\label{sec5.3}

In this section, results about transmutation operators, related to
spaces $L_2,$ are provided. The main result is as follows.

\begin{theorem} \label{teo:2.3.1}
If $f \in \mathring{T}_{+}^{\infty} \lr{U_R},$
 then the relation
\begin{equation}
\|\mathfrak{G}_n f \|_{L_2 \lr{U_R}} =  \| f \|_{L_2 \lr{U_R}}
\label{2.3.1}
\end{equation}
 holds for odd $n$ and the estimate
\begin{equation}
\|\mathfrak{G}_n f \|_{L_2 \lr{U_R}} \leq \sqrt{2}  \| f \|_{L_2 \lr{U_R}},
\label{2.3.2}
\end{equation}
 where the constant is exact if  $R = \infty,$ holds for even $n$.
\end{theorem}

\begin{proof}
    Representing the $L_2$-norm
    in the spherical coordinates and using \eqref{2.2.1}, we  obtain the relation
$$
    \begin{gathered}
\|\mathfrak{G}_n f \|_{L_2 \lr{U_R}}^2 = \int\limits_0^R
\int\limits_{\Theta} | \mathfrak{G}_n f (r, \vartheta) |^2 d
\vartheta r^{n-1} dr
    \\
= \int\limits_0^R \int\limits_{\Theta}  \left|
\sum\limits_{k=0}^{\mathcal{K}} \sum\limits_{l=1}^{d_k}
\frac{\sqrt{\pi} r^{\frac{n-1}{2}}}{2^{\frac{n}{2}+k-1} \, \Gamma
\lr{\frac{n}{2}+k}} S_{\frac{n}{2}+k-1} \lr{r^{-k} f_{k, l} }
Y_{k, l} \lr{\vartheta} \right|^2 \, d \vartheta dr.
    \end{gathered}
$$
    Since the system of spherical harmonics $Y_{k, l}$ is
    orthonormal in $L_2 \lr{\Theta}$, it follows that
$$
\|\mathfrak{G}_n f \|_{L_2 \lr{U_R}}^2 =
\sum\limits_{k=0}^{\mathcal{K}} \sum\limits_{l=1}^{d_k}
\frac{\pi}{2^{n+2k-1} \, \Gamma \lr{\frac{n}{2}+k}} \left\|
S_{\frac{n}{2}+k-1} \lr{r^{-k} f_{k, l} }  \right\|^2_{L_2 (0,
R)}.
$$
    Due to Lemma \ref{lem:1.4.1}, this implies the following relation provided that $n$
    is odd:
$$
\|\mathfrak{G}_n f \|_{L_2 \lr{U_R}}^2 =
\sum\limits_{k=0}^{\mathcal{K}} \sum\limits_{l=1}^{d_k}  \left|
r^{-k} f_{k, l}   \right|^2_{L_{2, \frac{n}{2}+k-1} (0, R)}.
$$
    Then \eqref{2.3.1} follows from Theorem \ref{teo:1.4.1}.
    Estimate \eqref{2.3.2} is proved in the same way.
    This completes the proof of the theorem.
\end{proof}
  Note that the finiteness of $R$ plays no role here. Therefore, relations
  \eqref{2.3.1}-\eqref{2.3.2} remain to be valid for $R=\infty$ as well
  provided that $f \in \mathring{T}_{+}^{\infty}
\lr{U_{\infty}} =  \mathring{T}_{+}^{\infty} \lr{E^n}.$
  From Corollary \ref{cor:1.4.2}, we conclude that the domain
 $\mathring{T}_{+}^{\infty} \lr{E^n}$ of the operator $\mathfrak{G}_n$
  is dense in $L_2 (E^n).$
  Therefore, the operator $\mathfrak{G}_n$ can be extended as a continuous and bounded
  in $L_2 (E^n)$ operator. Denote such an extension by $\mathfrak{G}_n^{L_2}.$
  let  $f, g \in\mathring{T}_{+}^{\infty} \lr{E^n}.$
  Then \eqref{2.2.1} yields the relation
$$
    \begin{gathered}
\lr{\mathfrak{G}_n^{L_2} f, g }_{L_2 (E^n)} =
\int\limits_0^{\infty} r^{n-1}  \int\limits_{\Theta}
\mathfrak{G}_n^{L_2} (r, \vartheta) \cdot g (r, \vartheta) \, d
\vartheta dr
    \\
 =\int\limits_0^{\infty} r^{n-1} \int\limits_{\Theta}
  \sum\limits_{k=0}^{\mathcal{K}(f)} \sum\limits_{l=1}^{d_k}
   \lr{ \frac{\sqrt{\pi}}{2^{\frac{n}{2}+k-1} \, \Gamma \lr{\frac{n}{2}+k}}
    r^{\frac{1-n}{2}}  S_{\frac{n}{2}+k-1}
 \lr{r^{-k} f_{k, l} } Y_{k, l} \lr{\vartheta} }   \sum\limits_{k=0}^{\mathcal{K}(g)}
  \sum\limits_{l=1}^{d_k} g_{k, l} (r)  Y_{k, l} \lr{\vartheta} \, d \vartheta d r
  \\
= \sum\limits_{k=0}^{\min\limits \lr{\mathcal{K}(f),
\mathcal{K}(g)}} \sum\limits_{l=1}^{d_k}
\frac{\sqrt{\pi}}{2^{\frac{n}{2}+k-1} \, \Gamma
\lr{\frac{n}{2}+k}}   \int\limits_0^{\infty}  S_{\frac{n}{2}+k-1}
\lr{r^{-k} f_{k, l} } g_{k, l} (r) r^{\frac{n-1}{2}} \, d r.
    \end{gathered}
$$
    Using relation \eqref{1.1.23} and integrating by parts, we
    find that
$$
    \begin{gathered}
\lr{\mathfrak{G}_n^{L_2} f, g }_{L_2 (E^n)} =
\sum\limits_{k=0}^{\min\limits \lr{\mathcal{K}(f),
\mathcal{K}(g)}} \sum\limits_{l=1}^{d_k}  \int\limits_0^{\infty}
\int\limits_r^{\infty} \tau^{\frac{n-1}{2}} f_{k, l} (\tau)
P_{\frac{n}{2}+k-\frac{3}{2}}^0 \lr{\frac{r}{\tau}} \, d \tau \,
D_r \lr{ r^{\frac{n-1}{2}}  g_{k, l} (r)} dr
 \\
=  \sum\limits_{k=0}^{\min\limits \lr{\mathcal{K}(f),
\mathcal{K}(g)}} \sum\limits_{l=1}^{d_k}  \int\limits_0^{\infty}
\tau^{\frac{n-1}{2}} f_{k, l} (\tau) \int\limits_0^{\tau} D_r \lr{
r^{\frac{n-1}{2}}  g_{k, l} (r)}   P_{\frac{n}{2}+k-\frac{3}{2}}^0
\lr{\frac{r}{\tau}} \, d r d \tau = \lr{f,
\lr{\mathfrak{G}_n^{L_2} }^* g }_{L_2 (E^n)},
    \end{gathered}
$$
    where
$$
\lr{\mathfrak{G}_n^{L_2} }^* g (r, \vartheta) =
\sum\limits_{k=0}^{\mathcal{K}(g)} \sum\limits_{l=1}^{d_k}
r^{\frac{n-1}{2}} \int\limits_0^r  D_{\tau} \lr{
{\tau}^{\frac{n-1}{2}}  g_{k, l} (\tau)}
P_{\frac{n}{2}+k-\frac{3}{2}}^0 \lr{\frac{\tau}{r}} \, d \tau
Y_{k, l} \lr{\vartheta}.
$$
 One can sum the last series, arguing as in the derivation of the corresponding representation
 for the operator $\mathfrak{G}_n$:
\begin{equation}
\lr{\mathfrak{G}_n^{L_2} }^* g (r, \vartheta) =  2^{-
\frac{1+n}{2}} \pi^{\frac{1-n}{2}} \int\limits_{\langle \vartheta,
\vartheta' \rangle > 0} D^{\frac{n-1}{2}}_{\gamma} \left.
\lr{\gamma^{\frac{n-1}{2}} g(\gamma r, \vartheta')}
\right|_{\gamma=\langle \vartheta, \vartheta' \rangle} d
\vartheta' \label{2.3.3}
\end{equation}
 if $n$ is odd and
\begin{equation}
\lr{\mathfrak{G}_n^{L_2} }^* g (r, \vartheta)  = 2^{-
\frac{1+n}{2}} \pi^{-\frac{n}{2}} \int\limits_{\langle \vartheta,
\vartheta' \rangle > 0} D^{\frac{n-2}{2}}_{\gamma}
\int\limits_0^{\gamma}  \left. D_{\tau} \lr{\tau^{\frac{n-1}{2}}
g(\tau r, \vartheta')} \frac{d \tau}{\sqrt{\gamma-\tau}}
\right|_{\gamma=\langle \vartheta, \vartheta' \rangle} d
\vartheta' \label{2.3.4}
\end{equation}
 if $n$ is even.

Since the operator $\mathfrak{G}_n^{L_2}$ is unitary provided that
$n$ is odd, it follows that relation \eqref{2.3.3} expresses the
inverse operator for $\mathfrak{G}_n^{L_2}$ as well.

 The operators $\mathfrak{G}_n^{L_2}$ can be expressed by means of the Fourier
 transformation as well. Let $f \in\mathring{T}_{+}^{\infty} \lr{E^n}.$
 Then, according to Lemma \ref{lem:1.3.2}, we obtain the relation
$$
\mathfrak{G}_n^{L_2} f (r, \vartheta) =
\sum\limits_{k=0}^{\mathcal{K}}  \sum\limits_{l=1}^{d_k} \frac{\pi
r^{\frac{1-n}{2}} Y_{k, l}
\lr{\vartheta}}{2^{\frac{n}{2}+k-\frac{1}{2}} \, \Gamma
\lr{\frac{n}{2}+k}} I^{\frac{3}{2}-\frac{n}{2}-k} F^{-1}_{-}
\lr{\rho F_{\frac{n}{2}+k-1} (r^{-k} f_{k, l}) },
$$
 where $I^{\mu}$ is the Liouville operator, $F_{-}^{-1}$ is the inverse sine-transformation,
  and $F_{\nu}$ is the Hankel transformation.
  By assumption, the functions $\rho F_{\frac{n}{2}+k-1} (r^{-k} f_{k, l})$ are odd.
  Therefore, the
 sine-transformation can be changed (with the corresponding factor) to the one-dimensional
 Fourier transformation as follows: $F_{-}^{-1} g = - 2 i F^{-1} g.$
 Then
$$
\mathfrak{G}_n^{L_2} f (r, \vartheta) = (2 \pi)^{\frac{1-n}{2}} r^{\frac{1-n}{2}}   I^{\frac{1-n}{2}} F^{-1} \sum\limits_{k=0}^{\mathcal{K}} \sum\limits_{l=1}^{d_k}
\frac{(-i)^k (2 \pi)^{\frac{n}{2}} \rho^k}{2^{\frac{n}{2}+k-1} \, \Gamma \lr{\frac{n}{2}+k}}
 F_{\frac{n}{2}+k-1} (r^{-k} f_{k, l})  Y_{k, l} \lr{\vartheta}.
$$
    Comparing this expression with relation \eqref{2.1.6}, we find
    that
$$
\mathfrak{G}_{n}^{L_2} f (r, \vartheta)  = (2 \pi
r)^{\frac{1-n}{2}} I^{\frac{1-n}{2}} F^{-1}
\lr{\rho^{\frac{n-1}{2}} F f},
$$
    where $F$ is the multidimensional  Fourier transformation.
    The final form of this relation is as follows:
$$
\mathfrak{G}_{n}^{L_2} f (r, \vartheta) =  (-2 \pi i
r)^{\frac{1-n}{2}} F^{-1} \lr{\rho^{\frac{n-1}{2}} F f}.
$$
    The last relations implicitly include the operator of the continuation to negative values
    of the radial variable; its form is clear from the context.

\chapter[Theory of
Buschman--Erd\'elyi Transmutation Operators]{Theory of
Buschman--Erd\'elyi\\
 Transmutation Operators}\label{ch3}

In this chapter, the theory of Buschman--Erd\'elyi
 transmutation operators is systematically explained.
 Special cases of integral
 operators of this class are considered in the previous chapter.

Note, that, following the tradition of the theory of transmutation
operators and the corresponding literature, we frequently use the
term ``operators'' in cases where the term ``differential
expressions'' is more exact.
    In main theorems, the function spaces of their validity are provided.
 For results containing explicit relations, if no particular function classes are mentioned,
  then we assume that compactly supported on the positive semiaxis functions
  (infinitely differentiable functions different from zero on an
  interval $(a,b)$, where $a>0$ and $b<\infty$) are meant.
  Also, very useful classes of test functions, considered in \cite{McB}, fit.

\section[First-Kind Buschman--Erd\'elyi
 Integral Transmutation Operators of Zero-Order
  Smoothness]{First-Kind Buschman--Erd\'elyi
 Integral Transmutation Operators\\
  of Zero-Order
   Smoothness}\label{sec6}

The considered class of transmutation operators is very important:
under a suitable choice of parameters, it simultaneously
generalizes
 So\-nin--Poisson--Delsarte
  transmutation operator and their adjoint ones,
Riemann--Liouville and Erd\'elyi--Kober
 fractional integrodifferentiating operators,  and
    Meler--Fock integral transformations.
    For the first time, integral operators of the specified kind, containing Legendre
    functions in kernels, occur in \cite{Cop1, Cop2, Cop3}.
    The detailed investigation of their solvability and invertibility is started in
    \cite{Bus1,Bus2, Erd1, Erd2, Erd3, Erd4, Erd5}. Also,
Buschman--Erd\'elyi operators or their analogs are studied in
\cite{Hig1,Ta1,Ta2,Love1,Love2,Din,Vir1,KiSk2, KiSk1}, etc;
 the main studied problems are to solve integral equations with
 such operators, to factorize them, and to convert them.
 These results are partially mentioned in  \cite{SKM} though the case of the integration
 limits selected by us, is not considered in this monograph (apart from one collection of
 composition relations), being treated as a special one.
 Several results regarding to the special choice of limits are added to the English extended
 edition of the monograph (see \cite{SaKiMar}).

The term ``Buschman--Erd\'elyi
 operators'' is introduced in \cite{S66, S6}; later, it is used by other authors as well.
 In \cite{SKM}, the term ``Buschman operators'' occurs.
 In the theory of the Radon transformation and mathematical tomography, the term
 ``Chebyshev--Gegenbauer  operators'' is used (see \cite{Rub3, Rub1, Rub2,Rub4}).
 The most complete study of Buschman--Erd\'elyi
 operators is contained in
  \cite{S1,S70,S72,S2,S73,S4,S66,S65,S6,S5,S7,S46,S14,S103,S400,SitDis,
  SSfiz,S42,S94,S38,S401,S402}.
  Note that the role of
   Buschman--Erd\'elyi operators as transmutation operators is not considered in earlier
   publications.

Among relatively recent papers dealing with  Buschman--Erd\'elyi
  operators  treated as integral ones, we note papers of Virchenko
(see \cite{Vir1,Vir2}), Kilbas, Rubin, Glushak, and their
disciples.
    For instance, in \cite{KiSk1, KiSk2},  Buschman--Erd\'elyi
    operators in Lebesgue weight spaces and their multidimensional generalizations
 in the form of integrals over  pyramidal regions are considered.
 In \cite{Vir1}, generalizations of standard Legendre functions are
 introduced and integrals close to
Buschman--Erd\'elyi operators are considered on the whole positive
semiaxis (note that Buschman--Erd\'elyi
 operators are defined only on a part of the positive semiaxis);
  no special case of the considered operators are Busch\-man--Er\-d\'elyi
  ones. In Rubin papers, domain sets and images of integral  Buschman--Erd\'elyi
   (Cheby\-shev--Ge\-gen\-ba\-u\-er) operators in function spaces are described
   (see \cite{Rub3, Rub1, Rub2, Rub4}); these results are applied to the
   Radon-transformation theory and tomography.
 Buschman--Erd\'elyi
 transmutation operators are used in \cite{Glu2, Glu3}.

Let us define
 first-kind Buschman--Erd\'elyi
 operators.

\begin{definition}
    The following integral operators are called the
    {\it first-kind Buschman--Erd\'elyi
    operators}:
\begin{align}
\label{2BE1}
B_{0+}^{\nu,\mu}f&=\int\limits_0^x \left( x^2-t^2\right)^{-\frac{\mu}{2}}P_\nu^\mu
 \left(\frac{x}{t}\right)f(t)d\,t,\\
\label{2BE1a}
E_{0+}^{\nu,\mu}f&=\int\limits_0^x
\left( x^2-t^2\right)^{-\frac{\mu}{2}}\mathbb{P}_\nu^\mu \left(\frac{t}{x}\right)
f(t)d\,t,\\
\label{2BE2}
B_{-}^{\nu,\mu}f&=\int\limits_x^\infty
 \left( t^2-x^2\right)^{-\frac{\mu}{2}}P_\nu^\mu \left(\frac{t}{x}\right)f(t)d\,t,\\
\label{2BE2a} E_{-}^{\nu,\mu}f&=\int\limits_x^\infty \left(
t^2-x^2\right)^{-\frac{\mu}{2}}\mathbb{P}_\nu^\mu
\left(\frac{x}{t}\right)f(t)d\,t.
\end{align}
\end{definition}

Here,  $P_\nu^\mu(z)$ is the first-kind
 Legendre function (see \cite{BE1}), $\mathbb{P}_\nu^\mu(z)$
  is the same function on the cut $-1\leq  t \leq 1,$
   $f(x)$ is a locally summable function satisfying growth
   restrictions as $x\to 0$ and $x\to\infty.$
   The parameters $\mu$ and $\nu$ are complex numbers such that $\Re \mu <1$
   and $\Re \nu \geq -1/2.$

For the first time, integral operators of the specified kind with
Legendre functions in kernels arise in Copson papers of late 1950s
about the
 Euler--Poisson--Darboux
 equations. In \cite{Cop1,Cop2}, the following assertion is considered.

\textbf{Copson lemma.} {\it Consider the two-variable
 partial differential equation
\begin{equation} \label{C1}
\frac{\pd^2 u(x,y)}{\pd x^2}+\frac{2\alpha}{x}\frac{\pd u(x,y)}{\pd x}=
\frac{\pd^2 u(x,y)}{\pd y^2}+\frac{2\beta}{y}\frac{\pd u(x,y)}{\pd y}
\end{equation}
{\rm (}the generalized
 Euler--Poisson--Darboux
 equation or, in Kipriyanov terms,
 $B$-hyperbolic equation{\rm )} in the open quarter of plane $\{x>0,\ y>0\}$
with parameters $\beta>\alpha>0$ and the following boundary-value
 conditions on the coordinate axes {\rm (}characteristics{\rm
):}
\begin{equation} \label{C2}
u(x,0)=f(x), u(0,y)=g(y), f(0)=g(0).
\end{equation}
    It is assumed that solutions $ u(x,y)$ {\rm (}if they exist{\rm
)} are continuously differentiable
    in the closed first quadrant and have continuous second derivatives in the open
    quadrant and the boundary-value
    functions $f(x)$ and $g(y)$ are continuously differentiable.

Then, if there exists a solution of the set problem, then the
following relations hold for it{\rm:}
\begin{equation}
\label{2Cop1}
\frac{\pd u}{\pd y}=0, y=0,  \frac{\pd u}{\pd x}=0, x=0,
\end{equation}
\begin{equation}
\label{2Cop2}
    \begin{gathered}
     2^\beta \Gamma(\beta{+}\frac{1}{2})\int\limits_0^1
f(xt)t^{\alpha{+}\beta{+}1} \lr{1{-}t^2}^{\frac{\beta
{-}1}{2}}P_{{-}\alpha}^{1{-}\beta}(t) dt
    \\
     =2^\alpha
\Gamma(\alpha{+}\frac{1}{2})\int\limits_0^1
g(xt)t^{\alpha{+}\beta{+}1} \lr{1{-}t^2}^{\frac{\alpha
{-}1}{2}}P_{{-}\beta}^{1{-}\alpha}(t) dt,
    \end{gathered}
\end{equation}
 and, therefore,
\begin{equation}
\label{2Cop3}
g(y)=\frac{2\Gamma(\beta+\frac{1}{2})}{\Gamma(\alpha+\frac{1}{2})
\Gamma(\beta-\alpha)}y^{1-2\beta} \int\limits_0^y
x^{2\alpha-1}f(x) \lr{y^2-x^2}^{\beta-\alpha-1}x \,dx,
\end{equation}
    where   $P_\nu^\mu(z)$ is the firs-kind
    Legendre function from} \cite{BE1}.

Relations \eqref{2Cop1} are obvious.
 In the cited paper, \eqref{2Cop2} is derived nonstrictly, i.\,e., it is found
 the the boundary-value
 functions (or values of the solution characteristics) cannot be arbitrary;
 they are related by Buschman--Erd\'elyi
 operators (in contemporary terms). This is the main content of the Copson lemma.
 Also, it is additionally proved in the present chapter that if two functions
  are related by
Buschman--Erd\'elyi operators of the specified orders, then
\eqref{2Cop3} is satisfied, i.\,e., they are related by
Erd\'elyi--Kober operators, which are simpler.

However, this does not imply that one can immediately convert (at
least, formally) the corresponding Buschman--Erd\'elyi
 operator. To do that, one has to represent the arbitrary function
 at the right-hand
 side of the corresponding equation by the Buschman--Erd\'elyi
 operator of the corresponding order (in order to apply the Copson lemma then).
 However, to do that, one has to convert the Buschman--Erd\'elyi
 operator,  which is exactly the problem to be solved.
 Thus, it would not be correct to attribute the first result on the converting
 of  Buschman--Erd\'elyi
 operators to Copson.

Proofs in \cite{Cop1,Cop2} are not strict. They are rather notes
to compose a plan of the proof.
 Perhaps, this is the reason not to include this result \cite{Cop4}.
 Note that the citation of \cite{Cop1} in \cite{SKM} and other papers
 is not completely correct; here, we provide the correct one.
  This pioneering Copson paper is developed in \cite{Cop3}:
  the strict derivation is provided, suitable function classes are introduced, and
  the relation to fractional integrals and  Kober--Erd\'elyi
  operators is clarified.

Now, we present results for
 Buschman--Erd\'elyi transmutation operators and their applications to
 differential equations with singularities at coefficients.

    In the sequel, we deal with the semiaxis. Therefore, $L_2$ denotes the space
  $L_2(0, \infty)$ and $L_{2, k}$ denotes the weight  space $L_{2, k}(0, \infty).$

First, we extend Definitions \eqref{2BE1}--\eqref{2BE2a}
  to the case where $\mu =1$ (this case is important, but it is not investigated
  earlier).

\begin{definition}
 For $\mu =1,$ the {\it
 Buschman--Erd\'elyi operators of zero-order
 smoothness}
  are introduced by the relations
\begin{align}
\label{2BE01}
B_{0+}^{\nu,1}f&={_1 S^{\nu}_{0+}f}=\frac{d}{dx}\int\limits_0^x P_\nu
\left(\frac{x}{t}\right)f(t)\,dt,\\
\label{2BE02}
E_{0+}^{\nu,1}f&={_1 P^{\nu}_{-}}f=\int\limits_0^x P_\nu
 \left(\frac{t}{x}\right)\frac{df(t)}{dt}\,dt,\\
\label{2BE03}
B_{-}^{\nu,1}f&={_1 S^{\nu}_{-}}f=\int\limits_x^\infty P_\nu
 \left(\frac{t}{x}\right)(-\frac{df(t)}{dt})\,dt,\\
\label{2BE04} E_{-}^{\nu,1}f&={_1
P^{\nu}_{0+}}f=(-\frac{d}{dx})\int\limits_x^\infty P_\nu
\left(\frac{x}{t}\right)f(t)\,dt,
\end{align}
    where $P_\nu(z)=P_\nu^0(z)$ is the Legendre function.
\end{definition}
 Under obvious additional conditions on the functions from \eqref{2BE01}--\eqref{2BE04},
 those relations can be differentiated under the integral sign
  or integrated by parts.

\begin{theorem}\label{2fact1}
 On suitable functions, the Buschman--Erd\'elyi
 can be factorized via fractional
  Riemann--Liouville and  Buschman--Erd\'elyi
  integral of  zero-order
 smoothness as follows{\rm :}
\begin{equation}\label{1.9}{B_{0+}^{\nu,\,\mu} f=I_{0+}^{1-\mu}~
 {_1 S^{\nu}_{0+}f},~B_{-}^{\nu, \,\mu} f={_1 P^{\nu}_{-}}~ I_{-}^{1-\mu}f}
 \end{equation}
 and
\begin{equation}\label{1.10}{E_{0+}^{\nu,\,\mu} f={_1
P^{\nu}_{0+}}~I_{0+}^{1-\mu}f,~E_{-}^{\nu, \, \mu} f=
I_{-}^{1-\mu}~{_1 S^{\nu}_{-}}f.}\end{equation}
\end{theorem}

\begin{proof}
    It suffices to prove the first relation
     {\rm(}other ones are proved in the same way{\rm)}.
 Taking into account the definitions, the compactness of the support of the function
 $f(x),$ the agreement from the beginning of the chapter, and the semigroup property
 of the fractional
Riemann--Liouville integrals, we obtain that
 $$ B_{0+}^{\nu,\,\mu}
f=I_{0+}^{1-\mu}~ {_1 S^{\nu}_{0+}f}=
 I_{0+}^{-\mu}~\int\limits_0^t P_\nu \left(\frac{t}{y}\right)f(y)\,dy=
 \frac{1}{\Gamma(-\mu)}~\int\limits_0^x (x-t)^{-\mu-1} \left(\int\limits_0^t
P_\nu \left(\frac{t}{y}\right)f(y)\,dy \right)d\,t.
$$
 Now, we change the order of the integrating (this is possible because the support
  of the function is compact). To compute the internal integral, we imply
  \cite[Vol. 3, p. 163, (7)]{PBM123}.
  We obtain the required integral representation for the
  first-kind Buschman--Erd\'elyi
 operators, which completes the proof of the theorem.
\end{proof}
  These important relations allow one to ``separate'' the parameters $\nu$
  and $\mu.$
  We prove that operators \eqref{2BE01}--\eqref{2BE04}
    are isomorphisms of the spaces $L_2(0, \infty)$ unless $\nu$
    is equal to an exceptional value. Therefore, the action of operators
 \eqref{2BE1}--\eqref{2BE2a} in $L_2$-type
 spaces is similar (in a way) to fractional integrodifferentiating
 operators $I^{1-\mu}$ (for $\nu=0,$ they coincide each other).
 Further, Buschman--Erd\'elyi
 operators are defined for all values of $\mu.$
    Basing on this, we introduce the following definition.

\begin{definition}
 The number $\rho=1-\Re \mu $ is called the  {\it smoothness order} of the
  Buschman--Erd\'elyi operators defined by \eqref{2BE01}--\eqref{2BE04}.
\end{definition}

Thus, if $\rho > 0$ (i.\,e., if $\Re \mu > 1$), then Buschman--Erd\'elyi
 operators increase the smoothness in $L_2 (0, \infty)$-spaces;
 if $\rho < 0$ (i.\,e., if  $\Re \mu < 1$),  then they decrease the smoothness in
 these spaces.
 According to the specified definition, operators \eqref{2BE01}--\eqref{2BE04}
  with $\rho =0$ are operators of the zero-order
  smoothness.
  Here, we treat smoothing operators (operators increasing the smoothness) as operators
  representable in the form $A=D^k B,$ where $k>0,$ while $B$ is a bounded in
   $L_2 (0, \infty)$ operator.
   Operators decreasing the smoothness are operators acting from the  space
    $C^k(0, \infty), k>0,$ to the Lebesgue space $L_2 (0, \infty).$

Below, we list the main properties of
    first-kind Buschman--Erd\'elyi
    operators
    \eqref{2BE1}--\eqref{2BE2a} with the   first-kind
    Legendre function in the kernel.
    We provide no proofs of these properties because they all
    follow from the main properties of Legendre functions.
    Domains of operators are denoted by  $\mathfrak{D}(B_{0+}^{\nu, \,
\mu }),$ $\mathfrak{D}(E_{0+}^{\nu, \, \mu }),$ etc.

Basic  properties of Legendre functions lead to the following
identities expressing the symmetry with respect to parameters,
adjacency relations, and the translation property with respect to
  parameters of Buschman--Erd\'elyi
  operators:
\begin{equation}\label{2.1}
\begin{aligned}
& B_{0+}^{\nu, \, \mu} f = B_{0+}^{-\nu-1,\mu} f, & E_{0+}^{\nu,\,
\mu}  f=E_{0+}^{-\nu-1, \, \mu}f,
\\
& B _{-}^{\nu, \, \mu} f=B_{-}^{-\nu-1, \, \mu}f, & E_{-}^{\nu,
\,\mu} f = E_{-}^{-\nu-1, \, \mu} f,
\end{aligned}
\end{equation}
\begin{equation}\label{2.2}
\begin{aligned}
& & (2 \nu +1)x \, B_{0+}^{\nu,\, \mu} {\frac{1}{x}f}=
(\nu-\mu+1)B_{0+}^{\nu+1, \, \mu}
f +  (\nu+\mu)B_{0+}^{\nu-1, \, \mu}f, \\
& & (2 \nu +1) \frac{1}{x} \, B_{-}^{\nu,\mu} {x
f}=(\nu-\mu+1)B_{-}^{\nu+1,\mu} f +  (\nu+\mu)B_{-}^{\nu-1,\mu}f,
\end{aligned}
\end{equation}
    and
\begin{equation}\label{2.3}
\begin{aligned}
& & B_{0+}^{\nu-1, \, \mu}f - B_{0+}^{\nu+1, \, \mu}f =
 -(2 \nu +1)B_{0+}^{0, \, \mu-1} { \frac{1}{x}f},\\
& & B_{-}^{\nu-1, \, \mu}f - B_{-}^{\nu+1, \, \mu}f = -(2 \nu
+1)\frac{1}{x}B_{-}^{\nu, \, \mu-1} f.
\end{aligned}
\end{equation}
    From the factorizing relations for $L_{\nu}$ (see \eqref{275} below),
 we obtain the identities
\begin{equation}
 B_{0+}^{\nu, \, \mu -1}\left(\frac{d}{dx}-
\frac{\nu}{x}\right) f =
 B_{0+}^{\nu-1, \, \mu } f  \label{2.4}
\end{equation}
    and
\begin{equation}
  B_{0+}^{\nu, \, \mu -1}\left(\frac{d}{dx}+
\frac{\nu}{x}\right) f = B_{0+}^{\nu+1,  \mu } f \label{2.5}
\end{equation}
 valid under the following assumptions:
  $\Re \mu < 1,$ $\Re \nu > -\dfrac{1}{2},$
$$
\lim\limits_{y \to 0} f(y) / y^{\nu}=0, ~ f \in
\mathfrak{D}(B_{0+}^{\nu \pm 1, \, \mu }),\quad\textrm{and}\quad
\left(\frac{d}{dx}\pm \frac{\nu}{x}\right)f \in
\mathfrak{D}(B_{0+}^{\nu,  \, \mu -1}).
$$
    Relations \eqref{2.1} allow  one to restrict the consideration by the case where
    $\Re \nu \geq -\dfrac{1}{2}.$
    If an operator acts on a function, then the function is to belong to the domain
    of the operator.
    The operator $E_{0+}^{\nu,\, \mu}$ satisfies the same relations
    as the operator $B_{0+}^{\nu,\, \mu}$ does.

\begin{theorem}
    The
 Buschman--Erd\'elyi operators given by \eqref{2BE1}--\eqref{2BE2a}
 are defined if $\Re \mu < 1$ or $\mu \in \mathbb{N}$
 and the following assumptions are satisfied{\rm:}
\begin{enumerate}
\item[(a)] for the operator $B_{0+}^{\nu, \, \mu},$ we have
$$
\int\limits_0^x \sq y \,|f(y) \ln y|\, dy < \infty
$$
 provided that $\nu=-\dfrac{1}{2}$ and $\mu \neq \dfrac{1}{2},$
 while we have
$$
\int\limits_0^x y^{-\Re \nu} \,|f(y)|\, dy < \infty
$$
 for all other cases{\rm;}

\item[(b)] for the operator $E_{0+}^{\nu, \, \mu},$ no additional
assumptions are imposed{\rm;}

\item[(c)] for the operator $E_{-}^{\nu, \, \mu},$ we have
$$
\int\limits_x^{\infty} y^{-\Re \nu} \,|f(y)|\, dy < \infty;
$$

\item[(d)] for the operator $B_{-}^{\nu, \, \mu},$ we have
$$
\int\limits_x^{\infty} y^{-\frac{1}{2}-\Re \nu} \,|\ln y \cdot
f(y)|\, dy < \infty
$$
 provided that  $\nu=-\dfrac{1}{2}$ and $\mu \neq \dfrac{1}{2},$  while we have
$$
\int\limits_x^{\infty} y^{\Re(\nu-\mu)} \,|f(y)|\, dy < \infty
$$
 for all other cases.
\end{enumerate}
\end{theorem}

In this theorem, it is assumed that  $f(x)$ is locally summable
over $(0, \infty)$ and $x$ is an arbitrary positive number.
    The proof follows from estimates of the absolute values of the integrals
     and the usage of the known asymptotic behavior of Legendre functions
     at the origin and infinity (see  \cite{BE1,NIST}).

There exist special values of the parameters $\nu$ and $\mu$ such that
  Buschman--Erd\'elyi operators are reduced to simpler operators for them.
  For example, for $\mu=-\nu$ or $\mu=\nu+2$, they are
Erd\'elyi--Kober operators.
    For $\nu = 0$, they are fractional integrodifferentiating operators
     $I_{0+}^{1-\mu}$ and $I_{-}^{1-\mu}.$ For
$\nu=-\dfrac{1}{2}$ and $\mu=0$ or $\mu=1,$ the kernels are
expressed via elliptic integrals.
 For $\mu=0,$ $x=1,$ and $v=it-\dfrac{1}{2},$ the operator
  $B_{-}^{\nu, \, 0}$ is the     Meler--Fock
  transformation (up to a constant).
  Thus, the
  first-kind Buschman--Erd\'elyi
  operators are generalizations of all the specified classes of standard integral
  operators.

\begin{theorem}
    Assume that  $\Re \mu < 0$ or $\mu = m \in \mathbb{N},$
$-m \leq \nu \leq m-1,$ and $\nu \in \mathbb{Z}.$ Then the
identities
\begin{equation}
 \frac{d}{dx} B_{0+}^{\nu, \, \mu } f = B_{0+}^{\nu, \, \mu +1 }
 f,
E_{0+}^{\nu, \, \mu }\frac{d \, f}{dx} = E_{0+}^{\nu, \, \mu +1 }
f  \label{2.6}
\end{equation}
    and
\begin{equation}
  B_{-}^{\nu, \, \mu } \lr{-\frac{d \, f}{dx}} = B_{-}^{\nu, \,
\mu + 1}f,   \lr{-\frac{d}{dx}} E_{-}^{\nu, \, \mu } f =
E_{-}^{\nu, \, \mu +1 } f \label{2.7}
\end{equation}
    are satisfied provided that all the specified operators are defined.
\end{theorem}

This theorem allows one to define the Buschman--Erd\'elyi
operators for the case where $\Re \mu \geq 1,$ redefining
 them for positive integers $\mu.$

\begin{definition}
 Let $\sigma$ be such that $\Re \sigma \geq 1.$
 Let  $m$ denote the least positive integer such that
$\sigma= \mu +m$ and $\Re \mu <1.$ Then the definition of the
Buschman--Erd\'elyi operators is completed by the relations
\begin{eqnarray}
& & B_{0+}^{\nu, \, \sigma}=B_{0+}^{\nu, \, \mu + m}=\lr{\frac{d}{dx}}^m \,
B_{0+}^{\nu, \, \mu}, \nonumber\\
& & E_{0+}^{\nu, \, \sigma}=E_{0+}^{\nu, \, \mu + m}=
E_{0+}^{\nu, \, \mu} \, \lr{\frac{d}{dx}}^m, \label{2.9}\\
& & B_{-}^{\nu, \, \sigma}=B_{-}^{\nu, \, \mu + m}=
B_{-}^{\nu, \, \mu} \lr{-\frac{d}{dx}}^m, \nonumber\\
& & E_{-}^{\nu, \, \sigma}=E_{-}^{\nu, \, \mu +
m}=\lr{-\frac{d}{dx}}^m E_{-}^{\nu, \, \mu}. \nonumber
\end{eqnarray}
\end{definition}

Note that if $\mu$ is a positive integer, then the Buschman--Erd\'elyi
 operators are defined by relations \eqref{2BE1}--\eqref{2BE2a}
 as well.
 We use relation \eqref{2.9} to redefine them for these values of  $\mu$.
 Thus, in the sequel, the symbols
$B_{0+}^{\nu, \, \mu},$ $E_{0+}^{\nu, \, \mu},$ $B_{-}^{\nu, \,
\mu},$ and $E_{-}^{\nu, \, \mu}$ are used to denote the operators
defined by relations \eqref{2BE1}--\eqref{2BE2a}
 provided that $\Re \mu < 1$ and relations \eqref{2.9}  provided that $\Re \mu \geq 1.$

Apart from the Bessel operator, consider the differential operator
\begin{equation}
\label{275}
L_{\nu}=D^2-\frac{\nu(\nu+1)}{x^2}=\left(\frac{d}{dx}-\frac{\nu}{x}\right)
\left(\frac{d}{dx}+\frac{\nu}{x}\right)=\left(\frac{d}{dx}+\frac{\nu+1}{x}\right)
\left(\frac{d}{dx}-\frac{\nu+1}{x}\right);
\end{equation}
    for $\nu \in \mathbb{N},$ it is the angular momentum operator from quantum
    mechanics. They are linked by the following easy verifiable relations.

 Let a pair  $\{X_\nu, Y_\nu\}$ of transmutation operators intertwine $L_{\nu}$
 and the second derivative, i.\,e.,
\begin{equation}
\label{276} X_\nu L_{\nu}=D^2 X_\nu \, ~\textrm{and}\, ~ Y_\nu D^2
= L_{\nu} Y_\nu.
\end{equation}
    Introduce the new pair of transmutation operators as follows:
\begin{equation}
\label{277} S_\nu=X_{\nu-1/2} x^{\nu+1/2} \, ~\textrm{and}\, ~
P_\nu=x^{-(\nu+1/2)} Y_{\nu-1/2}.
\end{equation}
    This pair intertwines  the Bessel operator and the second
    derivative:
\begin{equation}
\label{278} S_\nu B_\nu = D^2 S_\nu \, ~\textrm{and}\, ~  P_\nu
D^2 = B_\nu P_\nu.
\end{equation}
 Conversely, the specified relations can be used to pass from transmutation operators
 for the Bessel operator to transmutation operators for the angular momentum operator.
 Namely, if a pair $\{S_\nu, P_\nu\}$ of transmutation operators possessing
 property \eqref{278} is given, then the pair of transmutation
 operators defined by the relations
\begin{equation}
\label{279} X_\nu=S_{\nu+1/2} x^{-(\nu+1)} \, ~\textrm{and}\, ~
Y_\nu=x^{\nu+1} P_{\nu+1/2}
\end{equation}
    satisfy relations \eqref{276}.

Transmutation operators intertwining the second derivative with
the angular momentum operator have the following advantage
compared with  transmutation operators intertwining the second
derivative with the Bessel operator: under suitable conditions,
they are bounded in a same space instead of a pair of different
spaces.
 Transmutation operators acting according to relations \eqref{276}
 are still called
    Sonin-type and Poisson-type
    transmutation operators respectively.

Let us find conditions guaranteeing that
 first-kind Buschman--Erd\'elyi
 operators are  transmutation operators.

Define the class $\Phi(B_{0+}^{\nu, \, \mu})$ as the set of
functions $f$ such that
\begin{enumerate}
\item[(1)] the function $f(x)$ belongs to the set
$\mathfrak{D}(B_{0+}^{\nu, \, \mu})\bigcap\limits C^2(0, \infty);$

\item[(2)] the relation  $\lim\limits_{y \to 0}\left| \dfrac{\ln
y}{\sq y} f(y)+\sq y \ln y \cdot f'(y) \right|=0$ holds if
$\nu=-\dfrac{1}{2}$ and
$ \mu \neq \dfrac{1}{2},$\\
 the relation  $ \lim\limits_{y \to 0}{(\nu+1) y^{\nu} f(y)- y^{\nu+1} f'(y) }=0$
  holds if $\mu= \nu +1$ and $ \Re \nu \neq - \dfrac{1}{2},$\\
and the relation  $ \lim\limits_{y \to 0}\left( \nu
\dfrac{f(y)}{y^{\nu+1}}+\dfrac{f'(y)}{y^{\nu}} \right)=0 $ holds
in all other cases.
\end{enumerate}

\begin{theorem}
\label{2OPB0}
 If $\Re\mu \leq 1,$ then  the operator
  $B_{0+}^{\nu, \, \mu}$ is a Sonin-type
  transmutation operator satisfying relation \eqref{276}
   on functions $f(x) \in \Phi(B_{0+}^{\nu, \, \mu})$.
\end{theorem}

    A similar result holds for other Buschman--Erd\'elyi
 operators. Also, $E_{-}^{\nu, \, \mu}$  is a Sonin-type
  operator, while $E_{0+}^{\nu, \, \mu}$ and $B_{-}^{\nu, \, \mu}$ are Poisson-type
  ones. The proof (see \cite{S66}) is based on general conditions  for kernels
    of intertwining transmutation operators and asymptotic properties
    of Legendre functions.

One can consider the case where the lower limit in the
corresponding integrals from \eqref{2BE1}--\eqref{2BE2a}
 is equal to an arbitrary positive number  $a$ or the upper limit  is equal
 to an arbitrary positive and finite number   $b.$
 In this case, all results of the present section are preserved,
 but their formulations are substantially simplified. In particular, all
  Buschman--Erd\'elyi are defined under the only condition that $\Re \mu < 1$ in
  \eqref{2BE1}--\eqref{2BE2a} and are transmutation operators on functions $f(x)$
  such that $f(a)=f'(a)=0$ or  $f(b)=f'(b)=0$.

From the obtained theorem, it follows that Buschman--Erd\'elyi
transmutation operators bind eigenfunctions of the Bessel operator
and of the second derivative.
    Thus, a half of Buschman--Erd\'elyi
     transmutation operators maps trigonometric or exponential functions
into reduced  Bessel functions, while the other half acts
conversely. These relations are not provided here; it is easy to
 obtain them explicitly. They all are generalizations of the original
 Sonin and Poisson relations given by \eqref{151}-\eqref{152}.

    Now, consider various factorizations of
Buschman--Erd\'elyi operators via  Erd\'elyi--Kober
  operators and fractional
  Riemann--Liouville integrals  (see definitions in Chap. \ref{ch1}).

\begin{theorem}
\label{2factBE}
 The following relations factorizing  Buschman--Erd\'elyi
 operators via fractional integrodifferentiating operators and Erd\'elyi--Kober
 operators hold{\rm:}
\begin{equation}
  B_{0+}^{\nu, \, \mu}= I_{0+}^{\nu+1-\mu} I_{0+; \, 2, \, \nu+
\frac{1}{2}}^{-(\nu+1)} {\lr{\frac{2}{x}}}^{\nu+1}\label{2.17},
\end{equation}
\begin{equation}
  E_{0+}^{\nu, \, \mu}= {\lr{\frac{x}{2}}}^{\nu+1} I_{0+; \, 2,  -
 \frac{1}{2}}^{\nu+1} I_{0+}^{-(\nu+\mu)}  \label{2.18},
\end{equation}
\begin{equation}
  B_{-}^{\nu, \, \mu}= {\lr{\frac{2}{x}}}^{\nu+1}I_{-; \, 2, \, \nu+ 1}^{-(\nu+1)}
 I_{-}^{\nu - \mu+2}  \label{2.19},
\end{equation}
    and
\begin{equation}
  E_{-}^{\nu, \, \mu}= I_{-}^{-(\nu+\mu)} I_{-; \, 2, \, 0}
^{\nu+1} {\lr{\frac{x}{2}}}^{\nu+1}  \label{2.20}.
\end{equation}
\end{theorem}

 The proof is provided in \cite{S66}; it is similar to the one in \cite{Kat2,
32, KatDis, Kat3}. The proof using Mellin multipliers (see below)
is even simpler.

Many main properties of Buschman--Erd\'elyi
 operators treated as integral operators (not as transmutation ones)
 can be derived from Theorem \ref{2factBE}.
 For instance, one can obtain that the formal inverse to the  Buschman--Erd\'elyi
 operator with parameters
($\nu, \mu$) is the same  operator with parameters ($\nu,  2 -\mu$).
 One of these operators (the original one and the inverse one) has the integral
 representation given by \eqref{2BE1}--\eqref{2BE2a}, while the other one is defined
 by relations \eqref{2.9}; one of them has a positive order of smoothness, while the other one
 has a negative one (the only exception is the case of operators of zero-order
    smoothness). Also, Theorem \ref{2factBE} allows one to
    complete the definition of Buschman--Erd\'elyi
    operators on the whole range of parameters.
    This completion is coordinated with  \eqref{2.9}.
    Note that factorizations \eqref{2.17}--\eqref{2.20}
    are new compared with ones provided in \cite{SKM}.

Consider the relation between Buschman--Erd\'elyi
 operators and
  Sonin--Poisson--Delsarte
 intertwining operators.
    We prefer to introduce new definitions for them in order to preserve the unified
    notation system in this section.

\begin{definition}
    Redefine {\it
    Sonin--Poisson--Delsarte
 intertwining operators} (see Chap. \ref{ch1}) by the relations
\begin{equation}
 {_0S_{0+}^{\nu}}=B^{\nu, \, \nu+2}_{0+}= I_{0+; \, 2, \, \nu+
 \frac{1}{2}}^{-(\nu+1)} {\lr{\frac{2}{x}}}^{\nu+1}, \label{2.21}
 \end{equation}
 \begin{equation}
  {_0P_{0+}^{\nu}}=E^{\nu,  - \nu}_{0+}=
  {\lr{\frac{x}{2}}}^{\nu+1} I_{0+; \, 2,  - \frac{1}{2}}^{\nu+1}, \label{2.22}
 \end{equation}
 \begin{equation}
  {_0P_{-}^{\nu}}=B^{\nu, \, \nu+2}_{-}=
 {\lr{\frac{2}{x}}}^{\nu+1}I_{-; \, 2, \, \nu+ 1}^{-(\nu+1)}, \label{2.23}
 \end{equation}
    and
 \begin{equation}
 {_0S_{-}^{\nu}}=E^{\nu,  - \nu}_{-}= I_{-; \, 2, \, 0}^{\nu+1} {\lr{\frac{x}{2}}}^{\nu+1}. \label{2.24}
\end{equation}
\end{definition}

 Combining this definition with \eqref{162} and \eqref{1.15}--\eqref{1.16},
 we obtain the following integral representations:
\begin{align*}
{_0S_{0+}^{\nu}} f &= \begin{cases}
\ds\frac{2^{\nu+2}}{\Gamma(-\nu-1)}x \int\limits_0^x
(x^2-t^2)^{-\nu-2}t^{\nu+1}f(t)\,dt, & \Re \nu < -1,
\\
\ds\frac{2^{\nu+1}}{\Gamma(-\nu)} \frac{d}{dx} \int\limits_0^x
(x^2-t^2)^{-\nu-1}t^{\nu+1}f(t)\,dt, & \Re \nu < 0,
\end{cases}
\\
{_0P_{0+}^{\nu}} f &= \begin{cases} \ds\frac{1}{2^{\nu}
\Gamma(\nu+1)}x^{-\nu} \int\limits_0^x (x^2-t^2)^{\nu}f(t)\,dt, &
\Re \nu > -1,
\\
\ds\frac{1}{2^{\nu} \Gamma(\nu+2)}\frac{1}{x^{\nu+1}} \frac{d}{dx}
\int\limits_0^x (x^2-t^2)^{\nu+1}f(t)\,dt, & \Re \nu > -2,
\end{cases}
\\
{_0P_{-}^{\nu}} f &= \begin{cases}
\ds\frac{2^{\nu+2}}{\Gamma(-\nu-1)}x^{\nu +1}
\int\limits_x^{\infty} (t^2-x^2)^{-\nu-2}tf(t)\,dt, & \Re \nu <
-1,
\\
\ds\frac{2^{\nu+1}}{\Gamma(-\nu)} x^{\nu} \lr{-\frac{d}{dx}}
\int\limits_x^{\infty} (t^2-x^2)^{-\nu-1}tf(t)\,dt, & \Re \nu < 0,
\end{cases}
\end{align*}
    and
    \begin{equation*}
{_0S_{-}^{\nu}} f = \begin{cases} \ds\frac{1}{2^{\nu}
\Gamma(\nu+1)}\int\limits_x^{\infty}
(t^2-x^2)^{\nu}t^{-\nu}f(t)\,dt, & \Re \nu > -1,
\\
\ds\frac{1}{2^{\nu+1} \Gamma(\nu+2)} \lr{-\frac{1}{x}
\frac{d}{dx}} \int\limits_x^{\infty}
(t^2-x^2)^{\nu+1}t^{-\nu}f(t)\,dt, & \Re \nu > -2.
\end{cases}
\end{equation*}
  These operators are
  Sonin-type or Poisson-type
  intertwining operators.
  If we construct new transmutation operators for the angular momentum operator
   (see above), then we obtain Sonin-type
    operators
$$
X_{\nu} f= {_0S_{0+}^{\nu- \frac{1}{2}}} x^{\nu} f=
\frac{2^{\nu+\frac{3}{2}}}{\Gamma(-\nu-\frac{1}{2})}x
\int\limits_0^x (x^2-t^2)^{-\nu-\frac{3}{2}}t^{2 \nu+1}f(t)\,dt
$$
 for the case where $\Re \nu < -1/2$
\begin{equation}\label{2.25}{X_{\nu} f= S_{\nu} f=
\frac{2^{\nu+\frac{1}{2}}}{\Gamma(\frac{1}{2}-\nu)}\frac{d}{dx}
\int\limits_0^x (x^2-t^2)^{-\nu-\frac{1}{2}}t^{2
\nu+1}f(t)\,dt}
 \end{equation}
 for the case where $\Re \nu < 1/2.$

  In the same way, we obtain the Poisson-type
  operator of the kind
\begin{equation}\label{2.26}{Y_{\nu} f= P_{\nu} f=
\frac{1}{2^{\nu}\Gamma(\nu+1)}\frac{1}{x^{2 \nu}} \int\limits_0^x
(x^2-t^2)^{\nu-\frac{1}{2}}f(t)\,dt}\end{equation}
 provided that $\Re \nu > -1/2.$

Note that relations \eqref{1.9}--\eqref{1.10}
    are derived from Theorem \ref{2factBE} as well.

Passing to operators \eqref{2BE01}--\eqref{2BE04},
 we note that if the function $f(x)$ admits the differentiating under the integral sign or
 integrating by parts, then operators \eqref{2BE01}--\eqref{2BE04}
 take the form
\begin{equation}
 {_1S_{0+}^{\nu}}f=
f(x)+\int\limits_0^x \frac{\partial}{\partial
x}P_{\nu}\lr{\frac{x}{y}}f(y)dy, \label{2.31}
\end{equation}
\begin{equation}
 {_1P_{0+}^{\nu}}f=
f(x)-\int\limits_0^x \frac{\partial}{\partial
y}P_{\nu}\lr{\frac{y}{x}}f(y)dy, \label{2.32}
\end{equation}
\begin{equation}
 {_1P_{-}^{\nu}}f=
f(x)+\int\limits_x^{\infty} \frac{\partial}{\partial
y}P_{\nu}\lr{\frac{y}{x}}f(y)dy, \label{2.33}
\end{equation}
    and
\begin{equation}
 {_1S_{-}^{\nu}}f=f(x)-\int\limits_x^{\infty}
\frac{\partial}{\partial x}P_{\nu}\lr{\frac{x}{y}}f(y)dy.
\label{2.34}
\end{equation}
 Relations \eqref{2.32} and \eqref{2.33} are satisfied under the
 additional necessary conditions
$$
\lim\limits_{x \to 0} P_{\nu}(0) f(x) =
0\,~\textrm{and}\,~\lim\limits_{x \to \infty} P_{\nu}(x) f(x) = 0
$$
    respectively.

 It is easy to prove that, under suitable conditions for the functions, operators
  \eqref{2BE01}--\eqref{2BE04} are transmutation operators.
 They intertwine the angular momentum operator and the second derivative.

Theorem \ref{2factBE} yields the following factorizations for
Buschman--Erd\'elyi operators of the zero-order
 smoothness.

\begin{corollary}
    The relations
\begin{equation}
 {_1S_{0+}^{\nu}}= I_{0+}^{\nu+1} I_{0+; \, 2, \, \nu+ \frac{1}{2}}^{-(\nu+1)}
  {\lr{\frac{2}{x}}}^{\nu+1}, \label{2.35}
\end{equation}
\begin{equation}
 {_1P_{0+}^{\nu}}= {\lr{\frac{x}{2}}}^{\nu+1} I_{0+; \, 2,  - \frac{1}{2}}^{\nu+1}
  I_{0+}^{-(\nu+1)}, \label{2.36}
\end{equation}
\begin{equation}
 {_1P_{-}^{\nu}}= {\lr{\frac{2}{x}}}^{\nu+1}I_{-; \, 2, \, \nu+ 1}^{-(\nu+1)}
  I_{-}^{\nu+1}, \label{2.3.7}
\end{equation}
    and
\begin{equation}
 {_1S_{-}^{\nu}}= I_{-}^{-(\nu+1)} I_{-; \, 2, \, 0} ^{\nu+1}
{\lr{\frac{x}{2}}}^{\nu+1}. \label{2.3.8}
\end{equation}
 hold.
\end{corollary}

Consider properties of operators \eqref{2BE01} in detail. In
\cite{Kat2, 32, KatDis, Kat3}, a similar operator is constructed
as follows: the standard Sonin transmutation operator is
multiplied by a usual fractional integral to mutually compensate
the smoothness of these two operators and to obtain a new operator
acting in the same $L_2(0,\infty)$-type
    space (see Chap. \ref{ch2}).
  Later, it is found that the same can be done by means of known tools because
  the Sonin transmutation operator is a special case of Erd\'elyi--Kober
  operators. There exists a remarkable Erd\'elyi theorem allowing one to select
  a standard Riemann--Liouville
  integral from the fractional integral with respect to each function (see \cite{SKM}).
  This yields the following assertion.

\begin{theorem}\label{2tErd}
  For  $\Re \alpha > 0,$ consider the  Erd\'elyi--Kober
 fractional integrodifferentiating operator with respect to the function $g(x)=x^2:$
$$
I_{0+; \, x^2}^{\alpha} f = \frac{1}{\Gamma(\alpha)}
\int\limits_0^x (x^2-t^2)^{\alpha-1} 2t f(t)\,dt.
$$
    The following
 representation of the
 Erd\'elyi--Kober operator via a fractional Riemann--Liouville
 integral and operator \eqref{2BE01} takes place for $0<\Re \alpha <\dfrac{1}{2}:$
\begin{equation}
\label{2799}
  I_{0+,x^2}^{\alpha}(f)(x)  =  I_{0+}^{\alpha} \left( \left( 2x\right)^\alpha f(x) +
\int\limits_0^x \lr{\frac{\pd}{\pd x} P_{-\alpha}\left(
\frac{x}{t} \right)} \left( 2t\right)^\alpha f(t)\,dt \right)=
B_{0+}^{\nu,1}\lr{\lr{2x}^\alpha f},
\end{equation}
    where $I_{0+}^{\alpha}$ is the classical fractional Riemann--Liouville
 integral.
\end{theorem}

\begin{proof}
    From the Erd\'elyi theorem of \cite{SKM}, we obtain a
    representation of the kind
$$
    (2x)^{\alpha} f(x) +\int\limits_0^x \frac{\partial}{\partial x}\Phi(x, s)f(s)\,ds.
$$
    For the kernel $\Phi$, the following representation from
    \cite{SKM} is valid:
$$
\Phi (x, s)= \frac{\sin \pi \alpha}{\pi} 2s \cdot \int\limits_s^x
(x-u)^{-\alpha}(u-s)^{\alpha-1}
\frac{(u-1)^{1-\alpha}}{(u^2-s^2)^{1-\alpha}} du= \frac{\sin \pi
\alpha}{\pi} \cdot 2s \cdot \int\limits_s^x
(x-u)^{-\alpha}{(u^2-s^2)^{\alpha-1}} du.
$$
Using \cite[p. 301, (1)]{PBM123} to compute the integral, we
obtain that
$$
\Phi (x, s)= \frac{\sin \pi \alpha}{\pi}\cdot 2s \cdot
(2s)^{\alpha-1}\frac{\pi}{\sin \pi \alpha}\,{_2F_1(\alpha,
1-\alpha; 1; \frac{1}{2}-\frac{1}{2}\frac{x}{s})}=
(2s)^{\alpha}{_2F_1(\alpha, 1-\alpha; 1;
\frac{1}{2}-\frac{1}{2}\frac{x}{s})}.
$$
It remains to use \cite[p. 129, (14)]{PBM123}.
\end{proof}
  Operators of the zero-order
  smoothness  have the following specific property:
  estimates in a \textit{same} $L_p(0,\infty)$-type
  space can be obtained only for them.
  Taking into account the structure of these operators, we conclude that
  it is convenient to use the Mellin-transformation
  technique and Slater theorem  (see Chap. \ref{ch1}).

\begin{theorem} \label{2tmult}~\par
\begin{enumerate}
\item
 Buschman--Erd\'elyi operators of the zero-order
  smoothness act according to \eqref{1712}, i.\,e., as multipliers in
   Mellin images. The following relations hold for their symbols{\rm :}
\begin{multline}
m_{{_1S_{0+}^{\nu}}}(s)=\frac{\Gamma(-\frac{s}{2}+\frac{\nu}{2}+1)
\Gamma(-\frac{s}{2}-\frac{\nu}{2}+\frac{1}{2})}{\Gamma(\frac{1}{2}-\frac{s}{2})
\Gamma(1-\frac{s}{2})}
 \\
=\frac{2^{-s}}{\sq{ \pi}}
\frac{\Gamma(-\frac{s}{2}-\frac{\nu}{2}+\frac{1}{2})
 \Gamma(-\frac{s}{2}+\frac{\nu}{2}+1)}{\Gamma(1-s)}, \Re s <
 \min\limits (2 + \Re \nu, 1- \Re \nu)\label{2.311},
   \end{multline}
\begin{equation}
m_{{_1P_{0+}^{\nu}}}(s)=\frac{\Gamma(\frac{1}{2}-\frac{s}{2})\Gamma(1-\frac{s}{2})}
{\Gamma(-\frac{s}{2}+\frac{\nu}{2}+1)
\Gamma(-\frac{s}{2}-\frac{\nu}{2}+\frac{1}{2})}, \Re s < 1, \label{2.312}
\end{equation}
\begin{equation}
m_{{_1P_{-}^{\nu}}}(s)=\frac{\Gamma(\frac{s}{2}+\frac{\nu}{2}+1)
\Gamma(\frac{s}{2}-\frac{\nu}{2})}{\Gamma(\frac{s}{2})
\Gamma(\frac{s}{2}+\frac{1}{2})}, \Re s > \max\limits(\Re \nu,
-1-\Re \nu),
 \label{2.313}
  \end{equation}
  and
\begin{equation}
m_{{_1S_{-}^{\nu}}}(s)=\frac{\Gamma(\frac{s}{2})\Gamma(\frac{s}{2}+\frac{1}{2})}
{\Gamma(\frac{s}{2}+\frac{\nu}{2}+\frac{1}{2})
\Gamma(\frac{s}{2}-\frac{\nu}{2})}, \Re s >0. \label{2.314}
\end{equation}

\item The symbols satisfy the relations
\begin{equation}
 m_{{_1P_{0+}^{\nu}}}(s)=1/m_{{_1S_{0+}^{\nu}}}(s),~
 m_{{_1P_{-}^{\nu}}}(s)=1/m_{{_1S_{-}^{\nu}}}(s) \label{2.315}
 \end{equation}
 and
\begin{equation}
 m_{{_1P_{-}^{\nu}}}(s)=m_{{_1S_{0+}^{\nu}}}(1-s),~
m_{{_1P_{0+}^{\nu}}}(s)=m_{{_1S_{-}^{\nu}}}(1-s). \label{2.316}
\end{equation}

\item
    For
 $L_2$-norms of Buschman--Erd\'elyi
  operators of the zero-order
  smoothness, the relations
\begin{equation}
 \| _1{S_{0+}^{\nu}} \| = \| _1{P_{-}^{\nu}}\|=
 1/ \min\limits(1, \sq{1- \sin \pi \nu}) \label{2.322}
 \end{equation}
 and
 \begin{equation}
 \| _1{P_{0+}^{\nu}}\| = \| _1{S_{-}^{\nu}}\|= \max\limits(1,
\sq{1- \sin \pi \nu}) \label{2.323}
\end{equation}
    are valid.

\item
    The norms of operators
 \eqref{2BE01}--\eqref{2BE04} are $2$-periodic
 with respect to $\nu,$ i.\,e.,  $\|x^{\nu}\|=\|x^{\nu+2}\|,$
 where $x^{\nu}$ is any of operators \eqref{2BE01}--\eqref{2BE04}.

\item
    The norms of the operators ${_1 S_{0+}^{\nu}}$ and ${_1 P_{-}^{\nu}}$
    are not bounded in total with respect to  $\nu$ and each this norm is not less
    than $1.$
    If  $\sin \pi \nu \leq 0,$ then these norms are equal to $1.$
    The specified operators are unbounded in $L_2$ if and only if
    $\sin \pi \nu
= 1$ {\rm (}or $\nu=(2k) + 1/2,~k \in \mathbb{Z}${\rm )}.

\item
    The norms of the operators  ${_1 P_{0+}^{\nu}}$ and ${_1 S_{-}^{\nu}}$
 are bounded in total with respect to $\nu$ and each this norm does not exceed $\sq{2}.$
 All these operators are bounded in $L_2$ for all values of $\nu.$
 If $\sin \pi \nu \geq 0,$ then their $L_2$-norm
 is equal to $1.$
    The greatest value of the norm, equal to $\sq 2,$ is achieved  if and only if
   $\sin \pi \nu = -1$ {\rm (}or $\nu=-1/2+(2k),~k \in \mathbb{Z}${\rm )}.
\end{enumerate}
\end{theorem}

\begin{proof}
 We prove the claimed assertions only for the first operator.
 For the other ones, the proofs are similar.

1. First, we prove that relation \eqref{1712} with symbol
\eqref{2.311} holds.
 Successively using relations \cite[p. 130, (7), p. 129, (2), and p. 130, (4)]{Marich1},
 we obtain the relation
$$
    \begin{gathered}
M\left[B_{0+}^{\nu,1}\right](s)=\frac{\Gamma(2-s)}{\Gamma(1-s)}
M\left[\int\limits_0^\infty \left\{H(\frac{x}{y}-1)P_\nu
(\frac{x}{y}) \right\} \left\{y f(y)\right\}\frac{dy}{y}
\right](s-1)
    \\
=\frac{\Gamma(2-s)}{\Gamma(1-s)}   M \left[(x^2-1)_+^0P_\nu^0 (x)
\right] (s-1) M\left[f\right](s),
    \end{gathered}
$$
    where the following notation from \cite{Marich1} is used for the
   truncated power function and Heaviside function respectively:
$$
x_+^\alpha=\left\{
\begin{array}{rl}
x^\alpha, & \mbox{if } x\geqslant 0, \\
0, & \mbox{if } x<0  \\
\end{array}\right.
\quad \textrm{and} \quad H(x)=x_+^0=\left\{
\begin{array}{rl}
1, & \mbox{if } x\geqslant 0, \\
0, & \mbox{if } x<0. \\
\end{array}\right.
$$
    Further, using relations \cite[p. 234, 14(1) and p. 130, (4)]{Marich1},
    we obtain that
$$
M\lrs{(x-1)_+^0 P_\nu^0 (\sqrt x)}(s)=
\frac{\Gamma(\frac{1}{2}+\frac{\nu}{2}-s)\Gamma(-\frac{\nu}{2}-s)}
{\Gamma(1-s)\Gamma(\frac{1}{2}-s)}
$$
    and
$$
M\left[(x^2-1)_+^0 P_\nu^0 (x) \right](s-1)=\frac{1}{2} \frac {
\Gamma(\frac{1}{2}+\frac{\nu}{2}-\frac{s-1}{2})
\Gamma(-\frac{\nu}{2}-\frac{s-1}{2}) } {
\Gamma(1-\frac{s-1}{2})\Gamma(\frac{1}{2}-\frac{s-1}{2}) }=
\frac{1}{2}\frac { \Gamma(-\frac{s}{2}+\frac{\nu}{2}+1)
\Gamma(-\frac{s}{2}-\frac{\nu}{2}+\frac{1}{2}) }
{\Gamma(-\frac{s}{2}+\frac{3}{2})\Gamma(-\frac{s}{2}+1)}
$$
 provided that $\Re s<\min\limits(2+\Re\nu, \,1-\Re\nu).$
 This yields the following relation for the symbol:
$$
M\lrs{B_{0+}^{\nu,1}}(s)=\frac{1}{2}
\frac{\Gamma(2-s)}{\Gamma(1-s)}
\Gamma\lr{-\frac{s}{2}+\frac{3}{2}}\Gamma\lr{-\frac{s}{2}+1}.
$$
    Applying the Legendre duplication formula for the independent variable of the
 gamma-function to $\Gamma\lr{2-s}$ (see, e.\,g., \cite{BE1}), we
 obtain that
$$
M\lrs{B_{0+}^{\nu,1}}(s)=\frac{2^{-s}}{\sqrt\pi}  \frac {
\Gamma(-\frac{s}{2}+\frac{\nu}{2}+1)
\Gamma(-\frac{s}{2}-\frac{\nu}{2}+\frac{1}{2}) } {\Gamma(1-s)}.
$$
    Applying the Legendre duplication formula to $\Gamma(1-s)$ again,
 we arrive at the claimed relation for symbol \eqref{2.311}.

If $0<\Re s<1$, then this relation holds for all values of the
parameter $\nu,$ which is verified intermediately.

2. Now, we prove the relation for norm \eqref{2.322}.
 From the relation found for the symbol, by virtue of Theorem \ref{1tMel},
 we obtain that the following relations holds on the line $\{\Re s=1/2, s=i u+1/2\}$:
$$
|M\lrs{B_{0+}^{\nu,1}}(i u+1/2)|=\frac{1}{\sqrt{2\pi}}\left|\frac
{ \Gamma(-i\frac{u}{2}-\frac{\nu}{2}+\frac{1}{4})
\Gamma(-i\frac{u}{2}+\frac{\nu}{2}+\frac{3}{4}) }
{\Gamma(\frac{1}{2}-iu)}\right|.
$$
    In the sequel, we omit the generating operator in the notation of the symbol.
    Express the modulus of a complex number by the relation $|z|=\sqrt{z\bar{z}}$
 and use the identity $\ov{\Gamma(z)}=\Gamma(\bar z)$ following from the integral
 definition of the gamma-function.
 The last identity holds for the class of
  so-called real-analytic
  functions including the gamma-function.
  This yields the relation
$$
|M\lrs{B_{0+}^{\nu,1}}(i u+1/2)|=
\frac{1}{\sqrt{2\pi}}\left|\frac {
\Gamma(-i\frac{u}{2}-\frac{\nu}{2}+\frac{1}{4})
\Gamma(i\frac{u}{2}-\frac{\nu}{2}+\frac{1}{4})
\Gamma(-i\frac{u}{2}+\frac{\nu}{2}+\frac{3}{4})
\Gamma(i\frac{u}{2}+\frac{\nu}{2}+\frac{3}{4}) }
{\Gamma(\frac{1}{2}-iu)\Gamma(\frac{1}{2}+iu)}\right|.
$$
   In the numerator, combine the extreme factors and the mid-factors.
 Transform the obtained three pairs of gamma-functions
 according to the relation
$$
    \Gamma(\frac{1}{2}+z)\  \Gamma(\frac{1}{2}-z)=\frac{\pi}{\cos \pi z}.
$$
(see \cite{BE1}). This yields the relation
$$
|M\lrs{B_{0+}^{\nu,1}}(i u+1/2)|= \sqrt{ \frac{\cos(\pi i u)}
{2\cos\pi(\frac{\nu}{2}+\frac{1}{4}+i\frac{u}{2})
\cos\pi(\frac{\nu}{2}+\frac{1}{4}-i\frac{u}{2})} }= \sqrt{
\frac{\cosh(\pi  u)}{\cosh\pi u-\sin\pi\nu} }.
$$
    Introduce the notation  $t=\cosh\pi u, 1\le t <\infty,$ and apply the assumption
    of Theorem \ref{1tMel}. We obtain that
$$
\sup\limits_{u\in\R} |m(i u+\frac{1}{2})|=\sup\limits_{1\le t
<\infty} \sqrt{ \frac{t}{t-\sin\pi\nu} }.
$$
Therefore, if $\sin\pi\nu\ge 0,$ then the supremum is achieved for
$t=1$ and the claimed relation
$$
\|B_{0+}^{\nu,1}\|_{L_2}=\frac{1}{\sqrt{1-\sin\pi\nu}}
$$
    holds for the norm.
    If $\sin\pi\nu\le 0,$ then the supremum is achieved as $t\to\infty$
    and the relation
$$
\|B_{0+}^{\nu,1}\|_{L_2}=1
$$
    holds.

Now, assertions 3--6
 of the theorem directly follow from the obtained relation for the norm and assumptions
 of Theorem \ref{1tMel}, which completes the proof of the theorem.
\end{proof}
    It is important that
 Buschman--Erd\'elyi operators of the zero-order
 smoothness are unitary for integer values of $\nu.$
    Note that if  $L_{\nu}$ is interpreted as the quantum-mechanical
    operator of the angular momentum, then values of the parameter $\nu$
     are nonnegative integers.
     The next assertion is one of the main results of the present
     chapter.

\begin{theorem}\label{2tunit}
Operators \eqref{2BE01}--\eqref{2BE04}
 are unitary in $L_2$ if and only if $\nu$ is an integer.
 If this is satisfied, then
  $({_1 S_{0+}^{\nu}},$ ${_1 P_{-}^{\nu}})$ and
$({_1 S_{-}^{\nu}},$ ${_1 P_{0+}^{\nu}})$ are mutually inverse.
\end{theorem}

\begin{proof}
    For $\nu \in \mathbb{Z},$ we obtain that $\sin\pi\nu=0$ and
    the modulus of the corresponding symbol
 in relation \eqref{2.311} is identically equal to zero on the line $\Re s =\dfrac{1}{2}.$
Therefore, due to property (d) of Theorem \ref{1tMel}, the
specified operator and its inverse are unitary in $L_2(0,\infty).$
    Now, since the corresponding operator pairs are adjoint  in $L_2(0,\infty),$
it follows that they are mutually inverse, which completes the
proof of the theorem.
\end{proof}
    For the first time, this theorem was formulated in \cite{Kat1, 30, Kat2},
 but it was stated that the unitary property  takes place for all $\nu$.
 In \cite{S1, S70, S72, S2, S73}, it is corrected (see \cite{S66, S6, S46, S14, S400, S42,
S94, S38, S401, S402, SitDis} as well).

To formulate a special case as a corollary, assume that operators
 \eqref{2BE01}--\eqref{2BE04} are defined on functions $f(x)$ such that representations
 \eqref{2.31}--\eqref{2.34} hold (to achieve this, it suffices to assume that $x f(x) \to 0$
 as $x \to 0$). Then, for $\nu=1,$ we have the relations
\begin{equation}\label{2.325}{_1{P_{0+}^{1}}f=(I-H_1)f\,~\textrm{and}
\,~_1{S_{-}^{1}}f=(I-H_2)f,}
\end{equation}
    where
$H_1$ and $H_2$ are the Hardy operators,
\begin{equation}\label{2.326}{H_1 f = \frac{1}{x} \int\limits_0^x f(y) dy\,~\textrm{and}
\,~H_2 f = \int\limits_x^{\infty}  \frac{f(y)}{y}
dy}\end{equation}
   (see Chap. \ref{ch1}), and $I$ is the identity operator.

\begin{corollary}
 Operators \eqref{2.325} are unitary and mutually inverse in $L_2$.
    They intertwine differential expressions $d^2 / d x^2$ and $d^2 / d x^2 - 2/ x^2.$
\end{corollary}

Also, one can show that operators \eqref{2.325} are Cayley
transforms of symmetric operators $\pm 2 i (x f(x))$ under a
suitable choice of domains.

 In the unitary case,
  Buschman--Erd\'elyi operators of the zero-order
  smoothness form a pair of Watson biorthogonal transformations,
  while their kernels form pairs of nonsymmetric  Fourier kernels (see \cite{Dzh1}).

For the theory of integral equations, it is important to study the
unitary property.
    In this case, the inverse operator is to be sought in the form of an integral
    such that its integration limits are different from the ones of the original integral.

Also, consider the case where $\nu = i \alpha - \dfrac{1}{2},$
$\alpha \in \mathbb{R},$ related to the  Meler--Fock
  transformation.

\begin{theorem}
 Let $\nu = i \alpha - \dfrac{1}{2},$ $\alpha \in \mathbb{R}.$
 Then all operators \eqref{2BE01}--\eqref{2BE04}
 are bounded in $L_2$ and their norms satisfy the relations
$$
\| _1{S_{0+}^{i \alpha - \frac{1}{2}}}\|=\| _1{P_{-}^{i \alpha -
\frac{1}{2}}}\|=1.
$$
\end{theorem}

The proof is the same as in the case where $\nu$ is real.

Further, we list several general properties of operators such that
they act acting according to rule \eqref{1712}, i.\,e., as a
multiplier in Mellin images, and intertwine the second derivative
and the angular momentum operator.

\begin{theorem} \label{2tOPmult}
 Let an operator  $S_{\nu}$ act according to relations \eqref{1712} and \eqref{276}.
 Then
\begin{enumerate}
\item[(a)] its symbol satisfies the functional equation
\begin{equation}\label{2.5.1}{m(s)=m(s-2)\frac{(s-1)(s-2)}{(s-1)(s-2)-\nu(\nu+1)};}
\end{equation}

\item[(b)] if a function $p(s)$ is $2$-periodic, i.\,e.,
$p(s)=p(s-2),$  then the function $p(s)m(s)$
 is the symbol of another transmutation operator $S_2^{\nu}$
 intertwining  $L_{\nu}$ and the second derivative according to rule \eqref{276}.
\end{enumerate}
\end{theorem}

\begin{proof}
The latter assertion follows from the former one.
    To obtain Eq. \eqref{2.5.1} from \eqref{276}, we apply the Mellin transformation and use
 transmutation relations for basic operations (see \cite{Marich1}).
\end{proof}
The last theorem shows (again) that it is  convenient to study
transmutation operators in terms of Mellin multiplier.

The  Stieltjes transformation is defined by the following relation
(see, e.\,g., \cite{SKM}):
$$
(S f)(x)= \int\limits_0^{\infty} \frac{f(t)}{x+t} dt.
$$
    This operator is expressed by \eqref{1712} with symbol $p(s)=
\pi /sin (\pi s)$ and is bounded in $L_2.$ It is obvious
that $p(s)=p(s-2).$
    Therefore, it follows from Theorem \ref{2tOPmult} that the composition
    of the Stieltjes transformation with bounded intertwining operators given by
    \eqref{2BE01}--\eqref{2BE04} is a transmutation operator of the same type,
    bounded in $L_2.$

From the above reasoning, it follows that
$$
\|S\|_{L_r}=|\pi / \sin \frac{\pi}{r}|,\quad r>1.
$$
    On the other hand, the following relation holds:
$$
\|S\|_{L_{2,\, k}}=|\pi / \sin \pi k |,\quad k \notin  \mathbb{Z}.
$$
 Estimates in the weight spaces $L_{r, \, k}, r>0,$ are obtained in the same way.

    Now, consider the operator $H^{\nu}$ of kind \eqref{1712} with  symbol
\begin{equation}\label{2009}
 m(s)=\sq{\frac{\sin \pi s - \sin \pi \nu}{\sin \pi s}}.
\end{equation}
    From Theorem \ref{2tmult}, we obtain that the modulus of $m(s)$ is equal to one
    on the line $\Re s = \dfrac{1}{2}$. Then this theorem implies that
\begin{equation}\label{2010}
 {\|H^{\nu}\|_{L_2}=\|_1{P_{0+}^{\nu}}\|_{L_2}=\|_1{S_{-}^{\nu}}\|_{L_2}.}
 \end{equation}
 Therefore, the conclusion of Theorem \ref{2tmult} is valid for the operator $H^{\nu}$.
  In particular, the operator $H^{\nu}$ is bounded in $L_2$ for each $\nu.$

    Note that the formal relation between this operator and the Stieltjes
 transformation is as follows:
\begin{equation}\label{2.53}
    {H^{\nu}=(1-\frac{\sin \pi \nu}{\pi}S)^{\frac{1}{2}}.}
    \end{equation}
    Apart from  $H^{\nu}$, introduce the operator  $\mathfrak{D}^{\nu}$
    with symbol
$$
m_{\mathfrak{D}^{\nu}}(s)=\sq{\frac{\sin \pi s}{\sin \pi s - \sin \pi \nu}}.
$$
    This implies that
\begin{equation}\label{2.54}
{\|\mathfrak{D}^{\nu}\|_{L_2}=\|_1{S_{0+}^{\nu}}\|_{L_2}=\|_1{P_{-}^{\nu}}\|_{L_2}}
\end{equation}
 and the operator $\mathfrak{D}^{\nu}$ is bounded provided that
  $\sin \pi \nu \neq 1.$

\begin{theorem}
    Consider the operator compositions
\begin{equation}
 _3{S^{\nu}_{0+}}={_1S^{\nu}_{0+}}H^{\nu},~ _3{S^{\nu}_{-}}
=\mathfrak{D}^{\nu} {_1S^{\nu}_{-}} \label{2.55}
\end{equation}
    and
\begin{equation}
   _3{P^{\nu}_{0+}}=\mathfrak{D}^{\nu} {_1P^{\nu}_{0+}},~
_3{P^{\nu}_{-}}={_1P^{\nu}_{-}}H^{\nu}. \label{2.56}
\end{equation}
    The operators  ${_3 S_{0+}^{\nu}},$ ${_3 S_{-}^{\nu}}$ are
    Sonin-type transmutation operators, while ${_3
P_{0+}^{\nu}},$ ${_3 P_{-}^{\nu}}$ are Poisson-type transmutation
operators. All  these operators are unitary in $L_2.$ If $\sin \pi
\nu \neq 1,$ then compositions \eqref{2.55}--\eqref{2.56}
    can be computed in each order.
\end{theorem}

The proof of this theorem is obvious; it follows from the passage
to multipliers.
    In \cite{Lud}, a similar idea is applied to change the Radon transformation such
    that it becomes isometric.

Below, we obtain an explicit integral transformation of
transmutation operators intertwining  $L_{\nu}$ and $d^2/dx^2$
such that they are unitary  for all $\nu \in \mathbb{R}$.

Investigate the relations between
 left-side and right-side
  Buschman--Erd\'elyi operators. These relations are similar to the ones binding
   left-side and right-side
   fractional Riemann--Liouville
   integrals (see \cite[p. 163--171]{SKM}).
  Introduce the operator
   \begin{equation}\label{2.57}
   {C^{\nu} f=f(x)-\frac{\sin \pi \nu}{\pi} S f,}
   \end{equation}
 where $S$ is the Stieltjes transformation.
 We provide the following properties of $C^{\nu}$, omitting the proofs:

\begin{enumerate}
\item[(1)] $\|C^{\nu}\|_{L_2}=\min\limits (1, 1 - \sin \pi \nu)
\leq 1, ~ \nu \in \mathbb{R};$

\item[(2)] $\|C^{\nu}\|_{L_2}=1+\cosh \pi \alpha, ~  \nu=i \alpha
- \dfrac{1}{2},~ \alpha \in \mathbb{R}.$
\end{enumerate}

\begin{theorem}
 If  $\nu \in \mathbb{R},$ then the following identities for the
 compositions hold{\rm:}
\begin{equation}
 C^{\nu}={_1S_{-}^{\nu}} \ {_1P_{0+}^{\nu}}=
 {_1P_{0+}^{\nu}} \  {_1S_{-}^{\nu}}, \label{2.58}
\end{equation}
\begin{equation}
 {_1S_{-}^{\nu}}={_1S_{0+}^{\nu}} \ C^{\nu}, ~   {_1P_{0+}^{\nu}}=
  {_1P_{-}^{\nu}} \ C^{\nu}, \label{2.59}
\end{equation}
    and
\begin{equation}
 {_1S_{-}^{\nu}}=C^{\nu} \  {_1S_{0+}^{\nu}}, ~
{_1P_{0+}^{\nu}}=C^{\nu} \  {_1P_{-}^{\nu}}, ~ \sin \pi \nu \neq
1. \label{2.510}
\end{equation}
\end{theorem}

\section[Second-Kind Buschman--Erd\'elyi
 Integral Transmutation Operators\\
  and
  Sonin--Katrakhov and Poisson--Katrakhov
 Unitary Transmutation Operators]{Second-Kind Buschman--Erd\'elyi
 Integral Transmutation Operators\\
  and
  Sonin--Katrakhov and Poisson--Katrakhov
 Unitary Transmutation Operators}\label{sec7}
  %  \sectionmarknum{Операторы преобразования Buschman--Erd\'elyi второго рода,
%  Сонина--Катрахова и
%Пуассона--Катрахова}

Define and study
 second-kind Buschman--Erd\'elyi
   operators.
   In this section, several proofs are omitted (for brevity) because they mainly
   repeat the proofs from the previous section.

\begin{definition}
    Introduce another pair of
    {\it Buschman--Erd\'elyi operators} with second-kind
    Legendre functions (see \cite{BE1}) in the kernel:
\begin{equation}\label{2.61}{{_2S^{\nu}}f=\frac{2}{\pi}
 \left( - \int\limits_0^x (x^2-y^2)^{-\frac{1}{2}}Q_{\nu}^1 (\frac{x}{y}) f(y) dy  +
\int\limits_x^{\infty} (y^2-x^2)^{-\frac{1}{2}}\mathbb{Q}_{\nu}^1
(\frac{x}{y}) f(y) dy\right),}\end{equation}
\begin{equation}\label{2.62}{{_2P^{\nu}}f=\frac{2}{\pi} \left( - \int\limits_0^x
(x^2-y^2)^{-\frac{1}{2}}\mathbb{Q}_{\nu}^1 (\frac{y}{x}) f(y) dy -
\int\limits_x^{\infty} (y^2-x^2)^{-\frac{1}{2}}Q_{\nu}^1
(\frac{y}{x}) f(y) dy\right).}\end{equation}
\end{definition}

If $y \to x \pm 0$, then the integrals are understood in the sense
of the principal value. Under suitable conditions for the function
$f(x)$, these operators are defined and are intertwining, operator
\eqref{2.61} has the Sonin type, and operator \eqref{2.62} has the
Poisson type (here, we omit the proofs of all these facts).

\begin{theorem}\label{2rod}
  Operators \eqref{2.61}-\eqref{2.62}
 are representable by \eqref{1712} with symbols
\begin{eqnarray}
& & m_{_2S^{\nu}}(s)=p(s) \ m_{_1S_{-}^{\nu}}(s), \label{2.63}\\
& & m_{_2P^{\nu}}(s)=\frac{1}{p(s)} \ m_{_1P_{-}^{\nu}}(s),
\label{2.64}
\end{eqnarray}
    where the symbols of the operators  ${_1S_-^{\nu}}$ and
${_1P_-^{\nu}}$ are defined by relations \eqref{2.313} and
\eqref{2.314} respectively, while the $2$-periodic
 function $p(s)$ is expressed as follows{\rm:}
\begin{equation}\label{2.65}{p(s)=\frac{\sin \pi \nu+ \cos \pi s}{\sin \pi \nu -
 \sin \pi s}.}\end{equation}
\end{theorem}

First, we prove the following assertion.
\begin{lemma} \label{2lem2rod}
For $\Re \nu < 1,$ consider the integral operator
\begin{equation}\label{2.66}
{_3S^{\nu,\mu}}f=\frac{2}{\pi}
 \left( \int\limits_0^x (x^2+y^2)^{-\frac{\mu}{2}} e^{-\mu \pi i} Q_{\nu}^{\mu}
 ( \frac{x}{y}) f(y)\, dy +
 \int\limits_x^{\infty} (y^2+x^2)^{-\frac{\mu}{2}}\mathbb{Q}_{\nu}^{\mu} (\frac{x}{y})
  f(y)\, dy\right),
\end{equation}
    where $Q_{\nu}^{\mu}(z)$ is the
 second-kind Legendre function and $\mathbb{Q}_{\nu}^{\mu}(z)$ is its value on the cut
 {\rm(}note that this operator is more general than operator \eqref{2.61}{\rm)}.

On functions from  $C_0^{\infty}(0, \infty),$
 operator \eqref{2.66} is defined and acts as follows{\rm:}
\begin{equation}\label{2.67}
    \begin{gathered}
M{{_3S^{\nu,\mu}}}(s)=m(s)  M{x^{1-\mu} f}(s),
    \\
 m(s)=2^{\mu-1} \left( \frac{ \cos \pi(\mu-s) - \cos \pi \nu}{ \sin \pi(\mu-s)
  - \sin \pi \nu}  \right)  \left( \frac{\Gamma(\frac{s}{2})\Gamma(\frac{s}{2}+
  \frac{1}{2}))}{\Gamma(\frac{s}{2}+\frac{1-\nu-\mu}{2})
   \Gamma(\frac{s}{2}+1+\frac{\nu-\mu}{2})} \right).
       \end{gathered}
\end{equation}
\end{lemma}

\begin{proof}
    From the asymptotic behavior of the Legendre functions (see \cite{BE1}),
    it follows that operator \eqref{2.66} is defined.
    Representing it as the  Mellin convolution and successively applying relations
 \cite[p. 31, (2.50), p. 283, (10), p. 251, 40(1), and p. 130, (5)]{Marich1},
 we obtain relation \eqref{2.67}, where the symbol is equal to
$$
          \begin{gathered}
\frac{2^{\mu{-}2}}{\pi^2} \frac{\sin \pi(\nu{-}\mu)}{\sin \pi
\mu}\Gamma(\frac{s}{2})\Gamma(\frac{s}{2}{+}\frac{1}{2})
\Gamma(\frac{1{+}\mu{+}\nu}{2}{-}\frac{s}{2})
\Gamma(\frac{\mu{-}\nu}{2}{-}\frac{s}{2}){+}
\frac{2^{\mu{-}1}}{\sin \pi \mu}
\frac{\Gamma(\frac{1{+}\mu{+}\nu}{2}{-}\frac{s}{2})
\Gamma(\frac{\mu{-}\nu}{2}{-}\frac{s}{2})}
{\Gamma(1{-}\frac{s}{2})\Gamma(\frac{1}{2}{-}\frac{s}{2})}{}
    \\
{}+\frac{2^{\mu-1} \cos \pi \mu }{\sin \pi \mu} \frac{
\Gamma(\frac{s}{2})\Gamma(\frac{s}{2}+\frac{1}{2})}
{\Gamma(\frac{s}{2}+\frac{1-\nu-\mu}{2})
\Gamma(\frac{s}{2}+1+\frac{\nu-\mu}{2})}.
\end{gathered}
$$
    In the last two cases, ``lift'' the gamma-functions
    from the denominator to the numerator,  passing from $\Gamma (z)$
    to  $\Gamma (1-z).$
    We obtain that
$$
m(s)=\frac{2^{\mu-1}}{\pi^2 \sin \pi \mu} A(s)
\Gamma(\frac{s}{2})\Gamma(\frac{s}{2}+\frac{1}{2})
\Gamma(\frac{\mu+\nu+1-s}{2})\Gamma(\frac{\mu-\nu-s}{2}).
$$
    Successively transforming the expression  $A(s)$ by means of basic trigonometric
    rules, we see that
$$
          \begin{gathered}
A(s)= \sin \pi (\nu - \mu) +
 2 \sin \frac{\pi s}{2} \cos \frac{\pi s}{2} - 2 \cos \pi \mu
  \cos \frac{\pi}{2}(s-\nu-\mu) \sin \frac{\pi}{2}(s+\nu-\mu)
  \\
=\sin \pi (\nu - \mu) + \sin \pi s - \cos \pi \mu (\sin \pi \nu  +
\sin \pi (s-\mu))=
 \sin \pi \mu (\cos \pi (s- \mu) - \cos \pi \nu).
          \end{gathered}
$$
Substitute this expression in $m(s).$
 Moving the last two gamma-functions
 to the  denominator, we obtain \eqref{2.67}. Applying the specified relations from
 \cite{Marich1}, we arrive at the following restrictions for values of the variable:
\begin{equation}\label{2.68}{0 < \Re (\nu+\mu)< \min\limits
(1+\Re(\nu+ \mu),~ \Re (\mu-\nu)).}\end{equation}
    In our case, they can be weakened.
     \end{proof}
 Pass to the proof of Theorem \ref{2rod}. From relation \eqref{2.67}, it follows that
 one can pass to the limit as $\mu\to 1-0.$
 This yields \eqref{2.63} and \eqref{2.65}.
 A similar arguing is used to prove \eqref{2.64}.

Note that \eqref{2.66} is another family of Sonin-type
    transmutation operators.

Proofs of the following results are entirely the same as the ones
for the corresponding assertions regarding the first-kind
 Buschman--Erd\'elyi operators.

\begin{theorem}
    The following relations for norms hold{\rm:}
\begin{equation}
   \| {_2S^{\nu}} \|_{L_2}= \max\limits (1, \sq{1+\sin \pi \nu}) \label{2.69}
\end{equation}
    and
\begin{equation}
   \| {_2P^{\nu}} \|_{L_2}= 1 / {\min\limits (1, \sq{1+\sin \pi
\nu})}. \label{2.610}
\end{equation}
\end{theorem}

From this theorem, it follows that the operator ${_2S^{\nu}}$ is
bounded for all $\nu$ and the operator  ${_2P^{\nu}}$ is not
continuous if and only if $\sin \pi \nu=-1.$

\begin{theorem}
  The operators ${_2S^{\nu}}$ and ${_2P^{\nu}}$ are unitary in $L_2$
 if and only if $\nu$ is integer.
\end{theorem}

\begin{theorem}
    If $\nu=i \beta+1/2,~\beta \in \mathbb{R},$ then
\begin{equation}\label{2.6.11}{\| {_2S^{\nu}} \|_{L_2}=
\sq{1+\cosh \pi \beta},~\| {_2P^{\nu}} \|_{L_2}=1.}\end{equation}
\end{theorem}

\begin{theorem}
    The representations
\begin{equation}
  {_2S^0} f = \frac{2}{\pi} \int\limits_0^{\infty} \frac{y}{x^2-y^2}f(y)\,dy
 \label{2.612}
 \end{equation}
 and
 \begin{equation}
  {_2S^{-1}} f = \frac{2}{\pi} \int\limits_0^{\infty}
\frac{x}{x^2-y^2}f(y)\,dy \label{2.613}
\end{equation}
    are valid.
\end{theorem}

Thus, in this case, the operator  ${_2S^{\nu}}$ is reduced to a
    pair of known Hilbert transformations
    on the semiaxis (see \cite{SKM}).

Now, we can solve the important problem to construct transmutation
operators unitary for all values of $\nu.$ Such operators are
defined by the relations
\begin{equation}
  S_U^{\nu} f = - \sin \frac{\pi \nu}{2}\  {_2S^{\nu}}f+
\cos \frac{\pi \nu}{2}\  {_1S_-^{\nu}}f  \label{2.614}
\end{equation}
 and
\begin{equation}
  P_U^{\nu} f = - \sin \frac{\pi \nu}{2}\  {_2P^{\nu}}f+ \cos
\frac{\pi \nu}{2}\  {_1P_-^{\nu}}f. \label{2.615}
\end{equation}
 For all real values of $\nu,$ they are linear combinations of the
 first-kind and second-kind
  Buschman--Erd\'e\-lyi transmutation operators of the zero-order
  smoothness. They can be called the
   third-kind Busch\-man--Er\-d\'elyi
   operators.
   Their integral form is as follows:
\begin{multline} \label{2.616}
S_U^{\nu} f = \cos \frac{\pi \nu}{2} \left(- \frac{d}{dx} \right)
 \int\limits_x^{\infty} P_{\nu}\lr{\frac{x}{y}} f(y)\,dy   {}\\
+ \frac{2}{\pi} \sin \frac{\pi \nu}{2} \left(  \int\limits_0^x
(x^2-y^2)^{-\frac{1}{2}}Q_{\nu}^1 \lr{\frac{x}{y}} f(y)\,dy
\right. -  \int\limits_x^{\infty}
(y^2-x^2)^{-\frac{1}{2}}\mathbb{Q}_{\nu}^1 \lr{\frac{x}{y}}
f(y)\,dy \Biggl. \Biggr)
\end{multline}
    and
\begin{multline}\label{2.617}
P_U^{\nu} f = \cos \frac{\pi \nu}{2}  \int\limits_0^{x} P_{\nu}\lr{\frac{y}{x}}
 \left( \frac{d}{dy} \right) f(y)\,dy {}  \\
  {}-\frac{2}{\pi} \sin \frac{\pi \nu}{2}
  \left( - \int\limits_0^x (x^2-y^2)^{-\frac{1}{2}}\mathbb{Q}_{\nu}^1\lr{\frac{y}{x}}
  f(y)\,dy   \right.
 - \int\limits_x^{\infty} (y^2-x^2)^{-\frac{1}{2}} Q_{\nu}^1 \lr{\frac{y}{x}} f(y)\,dy
  \Biggl. \Biggr).
\end{multline}

\begin{theorem}\label{2unit}
 For all $\nu,$ operators
 \eqref{2.614}-\eqref{2.615} or \eqref{2.616}-\eqref{2.617}
 are unitary, mutually adjoint, and mutually inverse in $L_2.$
 They are intertwining operators and they act according to relations \eqref{275}.
 The operator $S_U^{\nu}$ has the Sonin {\rm
(}Sonin--Katrakhov{\rm )} type, while the operator  $P_U^{\nu}$
has the Poisson {\rm (}Poisson--Katrakhov{\rm )} type.
\end{theorem}

\begin{proof}
For one of the symbols, verify that relation \eqref{1716} is
satisfied. Similarly to Theorem \ref{2tmult}, we use the Mellin
convolution and relations for the Mellin transforms of special
functions to obtain the following relation for the symbol:
$$
    \begin{gathered}
M\lrs{S_U^\nu}(s)= -\sin(\frac{\pi\nu}{2}) \frac {-\cos\pi
s-\cos\pi\nu}{\sin\pi s-\sin\pi\nu}
\frac{\Gamma(\frac{s}{2})\Gamma(\frac{s}{2}+\frac{1}{2})}
{\Gamma(\frac{s}{2}-\frac{\nu}{2})\Gamma(\frac{s}{2}+\frac{\nu}{2}+\frac{1}{2})}
+\cos(\frac{\pi\nu}{2})
\frac{\Gamma(\frac{s}{2})\Gamma(\frac{s}{2}+\frac{1}{2})}
{\Gamma(\frac{s}{2}-\frac{\nu}{2})
\Gamma(\frac{s}{2}+\frac{\nu}{2}+\frac{1}{2})}
    \\
=\left( -\sin(\frac{\pi\nu}{2})  \frac {-\cos\pi
s-\cos\pi\nu}{\sin\pi s-\sin\pi\nu} + \cos(\frac{\pi\nu}{2})
\right) \frac{\Gamma(\frac{s}{2})\Gamma(\frac{s}{2}+\frac{1}{2})}
{\Gamma(\frac{s}{2}-\frac{\nu}{2})
\Gamma(\frac{s}{2}+\frac{\nu}{2}+\frac{1}{2})}
    \end{gathered}
$$
 (for the complete computation, see \cite{S66}).
 Further, according to Theorem \ref{1tMel}, taking into account the arguing
 from \ref{2tmult}, consider
$$
    \begin{gathered}
\left|M\lrs{S_U^\nu}(i u+\frac{1}{2})\right|=
\left|-\sin(\frac{\pi\nu}{2})  \frac {-\cos\pi (i
u+\frac{1}{2})-\cos\pi\nu}{\sin\pi (i u+\frac{1}{2})-\sin\pi\nu} +
\cos(\frac{\pi\nu}{2}) \right| \sqrt{\frac{\cos\pi u -\sin\pi
\nu}{\cos\pi u}}
    \\
=\left|
\frac{\sin(\frac{\pi\nu}{2})-\cos\left(\frac{\pi\nu}{2}+\pi i
u\right)} {\sin\pi\nu-\cos\pi i u} \right| \sqrt{\frac{\cos\pi u
-\sin\pi \nu}{\cos\pi u}}= \left|
\frac{\sin\frac{\pi\nu}{2}-\cos\frac{\pi\nu}{2}\cosh\pi u +
i\sin\frac{\pi\nu}{2}\sinh\pi u} {\sqrt{\cosh\pi
u(\cosh\pi-\sin\pi \nu)}} \right|.
 \end{gathered}
$$
Compute the absolute value. Then changes the (trigonometric and
hyperbolic) sines for the  (trigonometric and hyperbolic) cosines,
    we obtain the relation
$$
\left|M\lrs{S_U^\nu}(i u+\frac{1}{2})\right|=
\sqrt{\frac{\left(\sin\frac{\pi\nu}{2}-\cos\frac{\pi\nu}{2}\cosh\pi
u\right)^2+ \left(\sin\frac{\pi\nu}{2}\sinh\pi u\right)^2}
{\cosh\pi u(\cosh\pi-\sin\pi \nu)}}= \sqrt{\frac{\cosh^2\pi
u-\sin\pi\nu\cosh\pi u}{\cosh\pi u(\cosh\pi-\sin\pi \nu)}} =1.
$$
This proves the unitary property.
   The mutual adjointness follows from Definitions \eqref{2.614}-\eqref{2.615}
   if the operators are treated as extensions from the set of compactly supported
   functions.
   Hence, these unitary operators are mutually inverse as well.
  Verifying the assumptions of Theorem \ref{1tMel}, imposed in the symbols, we
  confirm that  these unitary operators have the Sonin and Poisson types,
  which completes the proof of the theorem.
\end{proof}
    This result completes the story about the constructing of unitary
 Sonin-type and Poisson-type
 transmutation operators. Unitary transmutation operators are
 closely related to the unitary property of the scattering
 operator in problems of quantum mechanics (see \cite{AM, Fad1,Fad2, ShSa}).

    Consider the special case following from Theorem \ref{2unit} for $\nu=1.$
    We obtain the pair
\begin{equation}
\label{282} B_{0+}^{1,1}f=f(x)-\frac{1}{x}\int\limits_0^x
f(y)\,dy,\  B_{-}^{1,1}f=f(x)-\int\limits_x^\infty
\frac{f(y)}{y}\,dy
\end{equation}
    of simple operators related to the famous Hardy operators
\begin{equation}
\label{283} H_1f=\frac{1}{x}\int\limits_0^x f(y)\,dy,
H_2f=\int\limits_x^\infty \frac{f(y)}{y}\,dy.
\end{equation}
(see \cite{OpKu} for Hardy inequalities).
    Our results yield the following assertion.

\begin{theorem} \label{2tHar}
    Operators \eqref{282} are mutually inverse and unitary in $L_2(0,\infty)$.
    They intertwine $\dfrac{d^2}{dx^2}$ and $\dfrac{d^2}{dx^2}-\dfrac{2}{x^2}.$
\end{theorem}

It follows from \eqref{283} that Buschman--Erd\'elyi
 operators can be treated as  generalizations of Hardy operators and inequalities for their
 norms generalize Hardy inequalities; this allow us to look at this class from a new
 viewpoint. Also, one can show that operators \eqref{282} are Cayley transformations
 of the symmetric operators $\pm 2i (xf(x))$ under a suitable choice of the domains.
 In \cite{S6, S66}, this is considered for spaces with power weights as well.

The result of Theorem \ref{2tHar} about the unitary property is
considered in \cite{KMP}: its elementary proof and applications
are provided.
  Theorem \ref{2unit} provides several other pairs of simple operators
   such that they are unitary in $L_2(0,\infty)$ and are special
   cases of Buschman--Erd\'elyi
    operators for integer values of $\nu$ (see \cite{S6, S66}):
\begin{align*}
\label{284} U_3f&= f+\int\limits_0^x f(y)\,\frac{dy}{y}, &  U_4f&=
f+\frac{1}{x}\int\limits_x^\infty f(y)\,dy,\\\nonumber U_5f&=
f+3x\int\limits_0^x f(y)\,\frac{dy}{y^2}, &  U_6f&=
f-\frac{3}{x^2}\int\limits_0^x y f(y)\,dy,\\\nonumber
U_7f&=f+\frac{3}{x^2}\int\limits_x^\infty y f(y)\,dy, & U_8f&=f-3x
\int\limits_x^\infty f(y)\frac{dy}{y^2},\\\nonumber
U_9f&=f+\frac{1}{2}\int\limits_0^x
\left(\frac{15x^2}{y^3}-\frac{3}{y}\right)f(y)\,dy, &
U_{10}f&=f+\frac{1}{2}\int\limits_x^\infty
\left(\frac{15y^2}{x^3}-\frac{3}{x}\right)f(y)\,dy.\\\nonumber
\end{align*}
 This list can be continued. Note that the above examples contradict the following
 statement of Marchenko: ``among Volterra operators, only the
 identity operator is unitary.''
  The cited statement refers to the second-kind
  integral Volterra operators; the majority of classical transmutation operators are
  of this form.

In  1980s, Katrakhov constructed (for the first time)
transmutation operators such that their form is similar to
\eqref{2.616}-\eqref{2.617}, but their kernels are expressed via
the general Gauss hypergeometric function.
    This is the reason to call them the
 Sonin--Katrakhov and Poisson--Katrakhov
 transmutation operators.
 We succeed to express them via
 first-kind and second-kind
 Legendre functions and to include them in the general scheme to construct transmutation
 operators by means of the composition method (see Chap. \ref{ch6} below).
 The simplest factorization relations of kind \eqref{2.614}-\eqref{2.615}
 are the main elements of this scheme.
 Now, the constructing of such operators is not a special skilful trick anymore;
 it is a component of a general technique of the constructing of classes
  of such transmutation operators by means of the composition method.

\section[Applications of
    Buschman--Erd\'elyi, Sonin--Katrakhov,
 and Poisson--Katrakhov\\
 Transmutation Operators to Differential Equations with
 Singularities at Coefficients]{Applications of
    Buschman--Erd\'elyi, Sonin--Katrakhov,
 and Poisson--Katrakhov
 Transmutation Operators to Differential Equations
  with
 Singularities at Coefficients}\label{sec8}
%\sectionmarknum{Приложения операторов преобразования к дифференциальным уравнениям с особенностями в
%коэффициентах}

\subsection[{Applications of
   Buschman--Erd\'elyi transmutation operators to problems\\
    for the
   Euler--Poisson--Darboux
   equations and to the Copson lemma}]{Applications of
   Buschman--Erd\'elyi transmutation operators to problems for the
   Euler--Poisson--Darboux
   equations and to the Copson lemma}\label{sec8.1}

 Now, consider  applications of the
  first-kind Buschman--Erd\'elyi
   transmutation operators of the zero-order
   smoothness to generalizations and refinements of the Copson lemma.

\begin{theorem}
    Consider the Dirichlet problem in the quarter of plane for the
     Euler--Poisson--Darboux
  equation under Conditions \eqref{C1}-\eqref{C2}.
  Then the following relations between the data of the Dirichlet problem,  expressed by
   first-kind Buschman--Erd\'elyi
   transmutation operators of the zero-order
   smoothness, are satisfied{\rm :}
\begin{equation}
\frac{c_\beta}{x^\beta}
E_{0+}^{-\alpha,1-\beta}(y^{\alpha+\beta+1}f(y))=
\frac{c_\alpha}{x^\alpha}
E_{0+}^{-\beta,1-\alpha}(y^{\alpha+\beta+1}g(y)),
c_\beta=2\Gamma(\beta+1/2),
\end{equation}
    and
\begin{equation}
\frac{c_\beta}{x^\beta} {_1 P_{0+}^{-\alpha}}
I_{0+}^{\beta}(y^{\alpha+\beta+1}f(y))= \frac{c_\alpha}{x^\alpha}
{_1 P_{0+}^{-\beta}} I_{0+}^{\alpha}(y^{\alpha+\beta+1}g(y)).
\end{equation}
\end{theorem}

These relations directly follow from Theorem \ref{2fact1}.
 On the other hand, applying Theorem \ref{2factBE}, where
  Buschman--Erd\'elyi operators of the zero-order
   smoothness are factorized via
    Riemann--Liouville and Erd\'elyi--Kober
    fractional integrals, we obtain the following result.

\begin{theorem}
 Under the assumptions of the  Copson lemma, the following relations expressed via
 Erd\'elyi--Kober operators are satisfied for the data of the Dirichlet problem{\rm:}
\begin{equation}
x^{\alpha+\beta+1}f(x)=\frac{c_\alpha}{c_\beta}
I_{0+;2;-1/2}^{\alpha-\beta} (y^{\alpha+\beta+1}g(y)).
\end{equation}
\end{theorem}

The last relation refines the corresponding result from the
original Copson work.

The relation between the axial data, expressed via
 Buschman--Erd\'elyi operators of the zero-order
   smoothness, allows one to establish additional results.

\begin{theorem}
 Let $\alpha$ and $\beta$ be integers.
 Then the relations
\begin{equation}
c_\beta \|I_{0+}^{\beta}(y^{\alpha+\beta+1}f(y))\|= c_\alpha
\|I_{0+}^{\alpha}(y^{\alpha+\beta+1}g(y))\|
\end{equation}
    are satisfied for the weight norms of the Dirichlet data in
    the weight space $L_{2,x^k}(0,\infty).$
\end{theorem}

    This result follows from the assumption about the unitary property of the
    Buschman--Erd\'elyi operators of the  zero-order
   smoothness in Theorem \ref{2tunit}.
 In this case, the data can be expressed via the inverse
  Buschman--Erd\'elyi operator of the zero-order
  smoothness (instead of the Erd\'elyi--Kober
  operator), obtaining a new relation; the same theorem is used
  for that purpose.
  For arbitrary noninteger values of the parameters, Theorem \ref{2tmult}
  yields a relation between weight norms of the Dirichlet data as well.

Consider applications of transmutation operators from the
considered classes to settings of the Cauchy problem for the
    Euler--Poisson--Darboux
  equation.

Consider  the
    Euler--Poisson--Darboux
  equation in the half-space
$$
B_{\alpha,\, t} u(t,x)= \frac{ \pd^2 \, u }{\pd t^2} +  \frac{2
\alpha+1}{t} \frac{\pd u}{\pd t}=\Delta_x u+F(t, x),
$$
    where $t>0$ and $x \in \mathbb{R}^n.$
 The following procedure allows one to obtain various settings of
 boundary-value conditions for $t=0$ by a unified method.
 Use relations \eqref{275} to form the  transmutation operators $X_{\alpha, \, t}$
 and $Y_{\alpha, \, t}.$
    Assume that $X_{\alpha, \, t} u:=v(t,x)$ and $X_{\alpha, \, t} F:=G(t,x)$ exist.
  Let the regular Cauchy problem
   \begin{equation}\label{2.7.28}{\frac{ \pd^2 \, v }{\pd t^2}
=\Delta_x v+G,~ v|_{t=0}=\varphi (x),~ v'_t|_{t=0}=\psi
(x)}\end{equation}
 be correctly solvable  in the half-space.
 Then we obtain the following boundary-value
 conditions for the
Euler--Poisson--Darboux
 equation:
\begin{equation}\label{2.7.29}{X_{\alpha} u|_{t=0}=a(x),~(X_{\alpha}
u)'|_{t=0}=b(x).}\end{equation}
 Boundary-value conditions depend on the choice of transmutation operators
  $X_{\alpha, t}$ (Sonin operators, Busch\-man--Erd\'elyi
 operators,
  first-kind  or  second-kind
  Buschman--Erd\'elyi operators of the zero-order
  smoothness, third-kind
  unitary operators defined by \eqref{2.616}, or generalized Buschman--Erd\'elyi
  operators).
  In each particular case, following the procedure explained above,
  one can reduce them to simpler analytic relations.
  Using integral transmutation operators of various types,
  we obtain nonlocal  boundary-value
 conditions (including conditions suitable for the considering of solutions with
 singularities) for each particular transmutation operator.

    This scheme is generalized for differential with many variables such that
  Bessel operators with different parameters act with respect to
  them and to equations of other types.
  Applying transmutation operators, one can reduce singular (degenerating)
  equations with Bessel operators with respect to one or several variables
   (Euler--Poisson--Darboux
   equations, singular heat equations, Kipriyanov $B$-el\-liptic
    equations, equations of the Weinstein generalized  axially  symmetric potential
      theory, etc.) to regular ones.
 In this case, a priori estimates for the singular case are obtained as corollaries
 from the corresponding  a priori estimates for regular case provided that
 transmutation operators themselves are estimated in suitable
 spaces.
    A lot of such estimates are provided above.

\subsection{Buschman--Erd\'elyi, Sonin--Katrakhov,
 and Poisson--Katrakhov
 transmutation operators and relations between solutions of differential equations}
 \label{sec8.2}

Each pair of transmutation operators from the present chapter,
    intertwining Bessel operators and the second derivative, can
    be used to set relations between solutions of  differential equations
    with Bessel operators (i.\,e., with singularities at coefficients)
    of the kind
$$
\sum\limits a_k B_{\nu_k,x_k} u(x) =f(x)
$$
    and solutions of the unperturbed constant-coefficient
    equation
$$
\sum\limits a_k \frac{\pr^2 v(x)}{\pr x_k^2} =g(x).
$$
    If the pairs of mutually inverse transmutation operators act
    by the relations
\begin{equation}
S_\nu B_\nu=D^2 S_\nu \quad \textrm{and}\quad  P_\nu D^2=B_\nu
P_\nu
\end{equation}
    (with respect to each variable), then solutions of the perturbed and unperturbed
    equation are related as follows:
     $$ u(x)=\prod\limits_k S_{\nu_k} v(x) \ \textrm{and} \
v(x)=\prod\limits_k P_{\nu_k} u(x)
$$
    (operators of the Sonin and Poisson types).
 Results on the boundedness, estimates of norms, and unitary properties  of
 the transmutation operators automatically lead to the corresponding assertions about
 pairs of solutions of the differential equations.
 We provide only this scheme, omitting the corresponding relations and estimates
 for solutions of  differential equations with singularities at coefficients.

    Consider the following example allowing us to use the constructed
 transmutation operators of various classes and their estimates to construct solutions
 of a nonlinear equation. It is known from \cite{Bitz1,Bitz12} that a number of problems
 of mathematical physics are reduced to the Maxwell--Einstein
  equations
\begin{equation}\label{ur1}
\Delta u +\frac{1}{x}u_x-\frac{1}{u}(1-\frac{u^2}{A^2-u^2})(u_x^2+u_y^2)=0
\end{equation}
    and
\begin{equation}\label{ur2}
\Delta u +\frac{1}{x}u_x-\frac{1}{u}(u_x^2+u_y^2)=0.
\end{equation}
    Using the cited Bitsadze--Pashkovskii
     results and the technique of transmutation operators, developed above,
  we obtain the following application of the considered classes of transmutation
  operators to the nonlinear Maxwell--Einstein
  equations.

\begin{theorem}
    Let $P$ be an arbitrary Poisson-type
    transmutation operator possessing the intertwining property
\begin{equation}
P\ D^2=(\frac{d^2}{d x^2}+\frac{1}{x}\frac{d}{dx})\ P
\end{equation}
 on smooth functions.
 Let $g(x,y)$ be an arbitrary harmonic function.
 Then the function $u_1(x,y)=\dfrac{A}{\cosh(a P_x g(x,y))}$ satisfies Eq. \eqref{ur1}
 and the function   $u_2(x,y)=\exp(b P_x g(x,y))$ satisfies Eq. \eqref{ur2},
 where $a$ and $b$ are arbitrary constants.
\end{theorem}

    We constructed various classes of  Poisson-type
    transmutation operators:
    first-kind and second-kind
    Bus\-chman--Erd\'elyi operators, operators of the zero-order
    smoothness,
    Sonin--Katrakhov and Pois\-son--Kat\-rakhov
    operators, etc.
    Now, we can use them in the above theorem in order to obtain representations
    of solutions of nonlinear  Maxwell--Einstein
  equations via harmonic functions.

\subsection[Applications of
 Sonin--Katrakhov and Poisson--Katrakhov
transmutation operators\\
 to integrodifferential
equations]{Applications of
 Sonin--Katrakhov and Poisson--Katrakhov
transmutation operators to integrodifferential
equations}\label{sec8.3}

Apply unitary
    Sonin--Katrakhov and Poisson--Katrakhov
  operators to the corresponding integrodifferential equations.

\begin{theorem}
Let $f(x)$ and $g(x)$ belong to $L_2(0,\infty)$ and be
continuously differentiable on the semiaxis.
 Then the integrodifferential equations
\begin{multline}
g(x) = \cos \frac{\pi \nu}{2} \left(- \frac{d}{dx} \right)
 \int\limits_x^{\infty} P_{\nu}\lr{\frac{x}{y}} f(y)\,dy  {}\\
{}+ \frac{2}{\pi} \sin \frac{\pi \nu}{2} \left(  \int\limits_0^x
(x^2-y^2)^{-\frac{1}{2}}Q_{\nu}^1 \lr{\frac{x}{y}} f(y)\,dy
\right. -  \int\limits_x^{\infty} (y^2-x^2)^{-\frac{1}{2}}
\mathbb{Q}_{\nu}^1 \lr{\frac{x}{y}} f(y)\,dy \Biggl. \Biggr)
\end{multline}
    and
\begin{multline}
f(x) = \cos \frac{\pi \nu}{2}  \int\limits_0^{x} P_{\nu}\lr{\frac{y}{x}}
\left( \frac{d}{dy} \right) g(y)\,dy  {} \\
{}-\frac{2}{\pi} \sin \frac{\pi \nu}{2} \left( - \int\limits_0^x
(x^2-y^2)^{-\frac{1}{2}}\mathbb{Q}_{\nu}^1\lr{\frac{y}{x}}
g(y)\,dy   \right. - \int\limits_x^{\infty}
(y^2-x^2)^{-\frac{1}{2}} Q_{\nu}^1 \lr{\frac{y}{x}} g(y)\,dy
\Biggl. \Biggr)
\end{multline}
    are mutually inverse and are solved by the above relations.
    In the specified space,
    the norms of the solutions and right-hand
    parts are equal to each other.
\end{theorem}

This is an application of the theorem on the unitary property of
 Sonin--Katrakhov and Poisson--Katrakhov
  transmutation operators.
  For special values of the parameter $\nu$ such that Legendre functions are expressed
  via simpler functions for those values, we obtain a list of particular
  integrodifferential equations such that explicit solutions with estimates
   of their norms are obtained. Here, we omit this list.

\section[Buschman--Erd\'elyi Transmutation Operators\\
 and Norm
Equivalence Between Kipriyanov Spaces and Weight Sobolev
Spaces]{Buschman--Erd\'elyi Transmutation Operators\\
 and Norm
Equivalence  Between Kipriyanov Spaces
 and Weight Sobolev Spaces}\label{sec9}
%\sectionmarknum{\sП\sр\sи\sл\sо\sж\sе\sн\sи\sя\s{ }\sо\sп\sе\sр\sа\sт\sо\sр\sо\sв\s{ }\sБ\sу\sш\sм\sа\sн\sа\s--\sЭ\sр\sд\sе\sй\sи\s{ }\sк\s{
%}\sэ\sк\sв\sи\sв\sа\sл\sе\sн\sт\sн\sо\sс\sт\sи\s{ }\sн\sо\sр\sм\s{
%}\sп\sр\sо\sс\sт\sр\sа\sн\sс\sт\sв\s{
%}\sК\sи\sп\sр\sи\sя\sн\sо\sв\sа\s{ }\sи\s{
%}\sв\sе\sс\sо\sв\sы\sх\s{ }\sп\sр\sо\sс\sт\sр\sа\sн\sс\sт\sв{
%}\sС\sо\sб\sо\sл\sе\sв\sа}

 Kipriyanov spaces are the best suitable for the study of $B$-elliptic
 partial differential equations (see \cite{Kip1,Kip2}).
 Therefore, properties of these spaces are important for the study
 of differential equations with singularities in coefficients.
 In the present section, it is proved that, for the one-dimensional
 case, the norms in Kipriyanov spaces and weight
 Sobolev spaces are equivalent to each other.
 For several prototype regions in multidimensional spaces, these results hold as well
 and they are, in fact, restatements of the results
 on the boundedness and unitary property of
 Buschman--Erd\'elyi transmutation operators of the zero-order
 smoothness in Lebesgue spaces on the semiaxis
  (see the previous section).
  Thus, the mentioned results about conditions of the boundedness and unitary property
   in Lebesgue spaces on the semiaxis for
   Buschman--Erd\'elyi transmutation operators of the zero-order
 smoothness are important for the theory of partial differential equations
 with Bessel operators: in this area, they find both standard and unexpected
 applications.
 Results for energy spaces for the corresponding differential equations
 are obtained as well.

In this section, the list of main results is provided without
proofs.
 In essence, the proofs follow from the properties of Buschman--Erd\'elyi
 operators, provided above.

 Define  $\mathfrak{D}(0, \infty)$ as the set of functions from
 $\mathfrak{D}(0, \infty),$ compactly supported at infinity.
 On this function set, introduce the seminorms
\begin{equation}
  \|f\|_{h_2^{\alpha}}=\|I_-^{\alpha}f\|_{L_2(0, \infty)}  \label{2.71}
\end{equation}
    and
\begin{equation}
   \|f\|_{\widehat{h}_2^{\alpha}}=\|x^{\alpha}
(-\frac{1}{x}\frac{d}{dx})^{\alpha}f\|_{L_2(0, \infty)},
\label{2.72}
\end{equation}
    where $I_-^{\alpha}$ is the Riemann--Liouville
    fractional derivative, the operator in \eqref{2.72} is
    defined by the relation
\begin{equation}\label{2.73}
 {(-\frac{1}{x}\frac{d}{dx})^{\beta}=2^{\beta}I_{-; \, 2,\,0}^{-\beta}x^{-2 \beta},}
    \end{equation}
  $I_{-; 2, \,0}^{-\beta}$ is the Erd\'elyi--Kober
  operator, and $\alpha$ is an arbitrary real number.
 For $\beta = n \in\mathbb{N}_0,$ expression \eqref{2.73}
 is treated in the classical sense, which is coordinated with Definitions
 \eqref{162}--\eqref{1.16}.

\begin{lemma}
 If $f(x) \in \mathfrak{D}(0, \infty),$ then the identities
\begin{equation}
  I_-^{\alpha}f={_1S_-^{\alpha-1}} {x^{\alpha} (-\frac{1}{x}\frac{d}{dx})^{\alpha}} f
   \label{2.74}
\end{equation}
    and
\begin{equation}
   x^{\alpha}
(-\frac{1}{x}\frac{d}{dx})^{\alpha}f={_1P_-^{\alpha-1}}
I_-^{\alpha}f \label{2.75}
\end{equation}
    hold.
\end{lemma}

    Thus, first-kind
    Buschman--Erd\'elyi operators of the zero-order
    smoothness bind the operators from the definitions of seminorms
    \eqref{2.71} and \eqref{2.72} provided that $\alpha \in \mathbb{N}$.

\begin{lemma}
 If $f(x) \in \mathfrak{D}(0, \infty),$ then the inequalities
\begin{equation}
   \|f\|_{h_2^{\alpha}} \leq \max\limits (1, \sq{1+\sin \pi \alpha})
 \|f\|_{\widehat{h}_2^{\alpha}}  \label{2.77}
\end{equation}
    and
\begin{equation}
   \|f\|_{\widehat{h}_2^{\alpha}} \leq \frac{1}{\min\limits (1,
\sq{1+\sin \pi \alpha})} \|f\|_{h_2^{\alpha}}, \label{2.78}
\end{equation}
 where $\alpha$ is each real number different from
     $-\dfrac{1}{2}+2k, k \in \mathbb{Z},$
     hold.
\end{lemma}

Note that the constants in inequalities \eqref{2.77}-\eqref{2.78}
    are not exceeded by one.
    If $\sin\pi \alpha = -1 $
 (i.\,e., $\alpha = -\dfrac{1}{2}+2k, k \in \mathbb{Z}$), then estimate \eqref{2.78}
 does not hold.

The above lemmas are restatements of results about factorization
relations and norm estimates for
    Buschman--Erd\'elyi  operators of the zero-order
    smoothness.

On  $\mathfrak{D} (0, \infty )$, introduce the Sobolev norm
\begin{equation}\label{2.79}
 {\|f\|_{W_2^{\alpha}}=\|f\|_{L_2 (0, \infty)}+\|f\|_{h_2^{\alpha}}}
 \end{equation}
    and the norm
\begin{equation}\label{2.710}
 {\|f\|_{\widehat{W}_2^{\alpha}}=\|f\|_{L_2 (0,\infty)}+\|f\|_{\widehat{h}_2^{\alpha}}}.
 \end{equation}
    Denote the spaces  $W_2^{\alpha}$ and $\widehat{W}_2^{\alpha}$ as the closures
    of $\mathfrak{D}(0, \infty)$ with respect to norms \eqref{2.79} and \eqref{2.710}
    respectively.

\begin{theorem}~\par
\begin{enumerate}
\item[(a)]
    For each $\alpha \in \mathbb{R},$ then space $\widehat{W}_2^{\alpha}$
    is continuously embedded in $W_2^{\alpha}$ and
\begin{equation}\label{2.711}
    {\|f\|_{W_2^{\alpha}}\leq A_1\|f\|_{\widehat{W}_2^{\alpha}},}
    \end{equation}
    where
$A_1=\max\limits (1, \sq{1+\sin \pi \alpha}).$

\item[(b)]
 Let $\sin \pi \alpha \neq -1,$ i.\,e., $\alpha \neq
-\dfrac{1}{2} + 2k, ~ k \in \mathbb{Z}.  $
    Then the inverse embedding of $W_2^{\alpha}$ in $\widehat{W}_2^{\alpha}$
 takes place and
\begin{equation}\label{2.712}{\|f\|_{\widehat{W}_2^{\alpha}}\leq A_2
\|f\|_{W_2^{\alpha}},}\end{equation}
    where
$A_2 =1/  \min\limits (1, \sq{1+\sin \pi \alpha}).$

\item[(c)]
 If $\sin \pi \alpha \neq -1,$ then the spaces $W_2^{\alpha}$ and
  $\widehat{W}_2^{\alpha}$ are isomorphic to each other and their norms are
  equivalent to each other.

\item[(d)]
    The constants in the embedding inequalities
    \eqref{2.711}-\eqref{2.712} are exact.
\end{enumerate}
\end{theorem}

In fact, this theorem follows from the boundedness result of
 Buschman--Erd\'elyi operators of the zero-order
 smoothness in $L_2,$ namely, from Theorem \ref{2tmult}.
 On the other hand,  the unitary property of these operators,
  proved in Theorem \ref{2tunit}, implies the following result about the equivalent norms
 in the considered variants of Sobolev spaces.

\begin{theorem} For each $s \in \mathbb{Z},$ the norms
\begin{equation}
  \|f\|_{W_2^{\alpha}} = \sum\limits_{j=0}^s \| \mathfrak{D}_-^j f\|_{L_2}  \label{2.713}
 \end{equation}
    and
\begin{equation}
  \|f\|_{\widehat{W}_2^{\alpha}}=\sum\limits_{j=0}^s \|
x^j(-\frac{1}{x}\frac{d}{dx})^j f \|_{L_2} \label{2.714}
\end{equation}
    are equivalent to each other in the Sobolev space.
Moreover, each term of \eqref{2.713} is identically equal to the
corresponding term of \eqref{2.714} with the same subscript $j.$
\end{theorem}

As we note above, the scale of spaces, introduced in
\cite{Kip1,Kip2}, have a substantial impact on the theory of
partial differential equations with the Bessel operator(s) with
respect to one or several variables. These spaces can be defined
as follows. Consider the subspace of $\mathfrak{D}(0, \infty),$
    consisting even functions such that all their derivatives of odd orders are equal to zero
    for $x=0.$ Denote this set by $\mathfrak{D}_c (0, \infty)$ and introduce the following
    norm on it:
\begin{equation}\label{2.715}
 {\|f\|_{\widetilde{W}_{2, k}^s} =\|f\|_{L_{2, k}}+\|B_k^{\frac{s}{2}}\|_{L_{2, k}}},
 \end{equation}
    where $s$ is an even number and $B^{s/2}_k$ is the corresponding iteration
    of the Bessel operator.  For even values of    $s$, the Kipriyanov space
 is defined as the closure of $\mathfrak{D}_c (0, \infty)$ with respect to norm \eqref{2.715}.
 It is known from \cite{Kip2} that a norm equivalent to \eqref{2.715}
 can be defined by the relation
\begin{equation}\label{2.716}{\|f\|_{\widetilde{W}_{2, k}^s} =
\|f\|_{L_{2, k}}+\|x^s(-\frac{1}{x}\frac{d}{dx})^s f\|_{L_{2,
k}}}.\end{equation}
    This allows one to define the norm in $\widetilde{W}_{2, \, k}^s$ for all values of $s.$
    In fact, this approach coincide with a approach from \cite{Kip2}; the other approach
    is based on the usage of the Fourier--Bessel
    transformation.
    In the sequel, we assume that $\widetilde{W}_{2, k}^s$ is normed according to
    relation \eqref{2.716}.

Introduce the weight Sobolev norm
\begin{equation}\label{2.717}{\|f\|_{W_{2, k}^s} =
\|f\|_{L_{2, k}}+\|\mathfrak{D}_-^s f\|_{L_{2, k}}}\end{equation}
    and define the space $W_{2, \, k}^s$ as the closure of
 $\mathfrak{D}_c (0,\infty)$ with respect to this norm.

\begin{theorem} \label{2tvloz1}~\par
\begin{enumerate}
\item[(a)] Let $k \neq -n$ and $n \in \mathbb{N}.$ Then the space
$\widetilde{W}_{2, \, k}^s$ is continuously embedded in $W_{2, \,
k}^s$ and there exists a positive constant $A_3$ such that
\begin{equation}\label{2.718}{\|f\|_{W_{2, k}^s}\leq A_3 \|f\|_{\widetilde{W}_{2,
k}^s}.}\end{equation}

\item[(b)] Let $k+s \neq -2m_1-1, ~ k-s \neq -2m_2-2, ~ m_1 \in
\mathbb{N}_0, $ and $ m_2 \in \mathbb{N}_0.$
    Then the inverse embedding of  $W_{2, \, k}^s$ in $\widetilde{W}_{2, \, k}^s$
    takes place and there exists  a positive constant $A_4$ such
    that
\begin{equation}\label{2.719}
    {\|f\|_{\widetilde{W}_{2, k}^s}\leq A_4 \|f\|_{W_{2,k}^s}.}
    \end{equation}

\item[(c)]
 If the specified conditions are not satisfied, then the corresponding embeddings do not
   take place.
\end{enumerate}
\end{theorem}

\begin{corollary}
    Let the following conditions be satisfied{\rm :} $k \neq -n,$ $n \in \mathbb{N},$
$k+s \neq -2m_1-1,$ $m_1 \in \mathbb{N}_0,$ and $k-s \neq
-2m_2-2,$ $m_2 \in \mathbb{N}_0.$
    Then the Kipriyanov space can be defined as the closure of
 $\mathfrak{D}_c(0, \infty)$ with respect to the weight Sobolev norm given by \eqref{2.717}.
\end{corollary}

\begin{corollary}
    The exact values of the constants in inequalities
\eqref{2.718}-\eqref{2.719} are as follows{\rm:}
$$
A_3 = \max\limits (1, \|{_1S_-^{s-1}} \| _ {L_{2, k}})\, ~
\textrm{and}\, ~ A_4=\max\limits(1, \|{_1P_-^{s-1}}\|_{L_{2, k}}).
$$
\end{corollary}

It is obvious that the presented theorem and its corollary
    follow from the results for the Buschman--Erd\'elyi
    operators, presented above.
Note that the $L_{2, \, k}$-norms
    of
 Buschman--Erd\'elyi operators of the zero-order
 smoothness yield the values of the exact constants in inequalities
 \eqref{2.718}-\eqref{2.719}.
 Norm estimates for Buschman--Erd\'elyi
 operators in Banach spaces $L_{p, \alpha}$ allow one to consider
 embeddings for the corresponding function spaces.

The inequality for seminorms, leading to embedding \eqref{2.718}
($s$ is an integer), is obtained in \cite{Lis}.
    Embeddings similar to the ones obtained in Theorem \ref{2tvloz1} are studied
    in \cite{Lei1,Lei2}. In Theorem \ref{2tvloz1}, we refine conditions providing
    the absence of embeddings from \cite{Lei1,Lei2}.
    In fact,  Theorem \ref{2tvloz1} establishes more exact (compared with \cite{Lei1,Lei2})
    inequalities between the corresponding seminorms.
    This is possible because Buschman--Erd\'elyi
    transmutation operators studied above in detail are applied.

  Pass to
  right-side intertwining operators defined by  \eqref{2BE01}--\eqref{2BE04}.
  We show that, in the general case, they implement an isomorphism between the Sobolev
  and Kipriyanov spaces.
  Define the spaces $H^{2s},$ $H_{\alpha}^{2s},$ and $K_{\alpha}^{2s}$ as the closures
  of the function set $\mathfrak{D} (0, \infty)$ with respect to the norms
\begin{eqnarray}
& & \|f\|_{H^{2s}} = \|f\|_{L_2}+\|I_-^{2s} f\|_{L_2}, \label{2.720} \\
& & \|f\|_{H^{2s}_{\alpha}} =
 \|f\|_{L_{2, {\alpha}}}+\|I_-^{2s} f\|_{L_{2, {\alpha}}}, \label{2.721}
\end{eqnarray}
 and
 \begin{equation}
  \|f\|_{K^{2s}_{\alpha}} = \|f\|_{L_{2,
{\alpha}}}+\|B_{\alpha}^s f\|_{L_{2, {\alpha}}}, \label{2.722}
    \end{equation}
    where $s$ is an integer and $\alpha \in \mathbb{R}.$
 Also, define the following pair of operators of type  \eqref{276}:
\begin{equation}\label{2.723}{{_1X_-^{\alpha}}={_1S_-^{\alpha-\frac{1}{2}}}
x^{\alpha+\frac{1}{2}}, ~ {_1Y_-^{\alpha}}=
x^{-(\alpha+\frac{1}{2})} {_1P_-^{\alpha-\frac{1}{2}}}.
}\end{equation}

\begin{theorem} \label{2tvloz2}
    Let $\alpha \in \mathbb{R}$ and $s \in \mathbb{N}.$ Then the operator
    ${_1X^2_-}$ continuously acts from $H^{2s}_{\alpha}$ to $K^{2s}_{\alpha}$
    and
\begin{equation}\label{2.724}
    {\|{_1X_-^{\alpha}}f\|_{H^{2s}_{\alpha}} \leq A_5\|f\|_{K^{2s}_{\alpha}},}
    \end{equation}
 where
$A_5=\|{_1X_-^{\alpha}}\|_{H^{2s}_{\alpha} \to
K^{2s}_{\alpha}}=\|{_1S_-^{\alpha-\frac{1}{2}}}\|_{L_2}=
\max\limits (1, \, \sq{1+\cos \pi \alpha}).$

    Let $s \in \mathbb{N}, \alpha \neq 2k+1,$ and $ k \in \mathbb{Z}$
{\rm (}or $\cos \pi \alpha \neq -1${\rm )}.
 Then the operator ${_1Y_-^{\alpha}}$ continuously acts from $K^{2s}_{\alpha}$
 to $H^{2s}_{\alpha}$ and the inequality
\begin{equation}\label{2.725}
    {\|{_1Y_-^{\alpha}}f\|_{K^{2s}_{\alpha}} \leq A_6\|f\|_{H^{2s}_{\alpha}}}
    \end{equation}
    holds with the constant
$$
A_6=\|{_1Y_-^{\alpha}}\|_{K^{2s}_{\alpha} \to
H^{2s}_{\alpha}}=\|{_1P_-^{\alpha-\frac{1}{2}}}\|_{L_2}=1/
\max\limits (1, \, \sq{1+\cos \pi \alpha}).
$$
\end{theorem}
    All assertions of the theorem follow from the above results for properties of
    norms of
 Buschman--Erd\'elyi transmutation operators of the zero-order
 smoothness.

    The  Bessel operator is the radial part of the Laplacian in $\mathbb{R}^n.$
  Under this interpretation of this operator,  Condition
  \ref{2tvloz2} is satisfied for the case of odd-dimensional
    spaces.

\begin{theorem} \label{2tvloz3}
    Assume that $\alpha \neq 2k+1,$ $k \in
\mathbb{Z},$ $\alpha \neq -n,$ $n \in \mathbb{N},$ $\alpha+2s \neq
-2m_1-1,$ $m_1 \in \mathbb{N}_0,$ $\alpha-2s \neq -2m_2-2,$ and
$m_2 \in \mathbb{N}_0.$
    Then \eqref{2.723} are operators of the topological isomorphism
    between
    the Sobolev space $H^{2s}$ and   the weight Sobolev space $H^{2s}_{\alpha}.$
\end{theorem}

    It is obvious that all assumptions of Theorem \ref{2tvloz3} are satisfied for all
    half-integer values of $\alpha.$ Therefore, the following assertion holds.

\begin{theorem}
 Let $s \in \mathbb{N}$ and $\alpha - \dfrac{1}{2} \in \mathbb{Z}.$
    Then \eqref{2.723} are operators of the topological isomorphism
    between
    the Sobolev space $H^{2s}$ and   the weight Sobolev space $H^{2s}_{\alpha}.$
\end{theorem}

In the same way, relations of kind \eqref{2.723} can be used to
introduce the operators ${_1X^{\alpha}}_{0+}$ and
${_1Y^{\alpha}}_{0+}.$
    Another application of the above results is the action of operators \eqref{2.723}
    in spaces with norms \eqref{2.720}--\eqref{2.722}
    for arbitrary weights not related to the constant $\alpha$ in the
    Bessel operator $B_{\alpha}.$

Results for the one-dimensional
 case, obtained in this section, are obviously extended to the multidimensional case
 for the region consisting of the Cartesian product of positive semiaxis or segments
    with respect to each variable.
    In particular, in the two-dimensional
 case, the obtained results can be immediately applied to estimate norms and prove
 embeddings in the first positive quadrant or a rectangle contained in it.

Thus, in this section, we use
 Buschman--Erd\'elyi transmutation operators of the zero-order
 smoothness to confirm that Kipriyanov spaces are isomorphic to weight
 Sobolev spaces.
 The presented results in now way diminish the importance
 of Kipriyanov spaces or the necessity to use them for suitable
 problems in the theory of partial differential equations with Bessel operators.
  The core importance of Kipriyanov spaces for the theory of partial differential
   equations (of various types) with Bessel operators is reflected
   by the general methodological approach formulated by Kudryavtsev as follows:
    ``{\it each differential equation is to be studied in its own space}''
    (a fragment of a plenary lecture listened by the author at a
    conference).

Also, the obtained results confirm the utility and efficiency of
 Buschman--Erd\'elyi  transmutation operators of the zero-order
 smoothness (in fact, this special class of operators is introduced
 by Katrakhov in 1980s) for the theory of differential
equations.

\chapter[General Weight Boundary-Value
 Problems for Singular Elliptic Equations]{General Weight Boundary-Value
 Problems\\
  for Singular Elliptic Equations}\label{ch4}

In a natural way, this chapter is decomposed into two parts.
    In the first part, new function spaces are introduced and studied for the half-space
    case; then they are introduced and studied for bounded regions.
    In the former case, we use transmutation operators introduced in the first section.
 In the latter one, the spaces are introduced by means of local charts.
 Also, we provide all the data to be used below, including
 internal embedding theorems, theorems on weight traces, and the
 corresponding theorems about multipliers.

    In the second part, we construct the regularizer for a regular
    equation with constant coefficients and nonlocal boundary-value
    conditions containing a Liouville-type
    operator. This is an auxiliary result.
    The main results of the chapter are concentrated at the final section:
    we set and study the general weight boundary-value
  problem for singular elliptic high-order
  equations. For the case of ``constant'' coefficients, we apply the method of
  transmutation operators, reducing the weight problem to the regular one
   specified above.
 For the case of variable coefficients, we apply a classical elliptic technique
  going back to Schauder.
  We construct a two-side
  regularizer in the spaces introduced above; this easily implies
  all main assertions regarding the  Noetherian property of the
  studied weight boundary-value
  problem.

\section{Function Spaces $H_{\nu}^s(E_{+}^{n+1})$}\label{sec10}

\subsection{Definitions and internal embedding theorems}\label{sec10.1}

 Let $E^{n+1}$ be the Euclidean $(n+1)$-di\-men\-sional
  space of points $x= (x', y)= (x_1, \dots, x_n, y),$ where $x' \in E^n,$ $y \in E^1.$
Let $E^{n+1}_{+}$ denote the half-space
 $\{y>0\}$ and $\ov{E^{n+1}_{+}}$  denote its closure.
 Throughout the present chapter, we use the transmutation operators $P_{\nu, e}$
 and $S_{\nu, e}$ (see Sec. \ref{sec4.3}) with respect to the last variable, i.\,e.,
 with respect to $y.$

Introduce the following notation. Let $\mathring{C}^{\infty}
(\ov{E^{n+1}_{+}})$ be the set of infinitely differentiable in
$\ov{E^{n+1}_{+}}$ functions such that their supports are compact
in  $\ov{E^{n+1}_{+}}$. Let $\mathring{C}^{\infty}_{\nu}
(E^{n+1}_{+})$ be the set of functions admitting the
representation  $f = P_{\nu, e}\, g,$ where $g \in
\mathring{C}^{\infty} (\ov{E^{n+1}_{+}}).$
    Denote this as follows:
$\mathring{C}^{\infty}_{\nu} (E^{n+1}_{+}) = P_{\nu, e}\,
\mathring{C}^{\infty} (\ov{E^{n+1}_{+}}).$

 For positive integer values of $s,$ introduce the following norm
 on the space $\mathring{C}^{\infty}_{\nu} (E^{n+1}_{+})$:
\begin{equation}
\| f \|_{H^s_{\nu} \lr{E^{n+1}_{+}}} = \| S_{\nu, e} f \|_{H^s
\lr{E^{n+1}_{+}}},
    \label{3.1.1}
\end{equation}
 where $H^s (\Omega) \equiv W^s_2 (\Omega)$ denote the Sobolev spaces.
 Norm \eqref{3.1.1} is well defined because $S_{\nu, e} = \lr{P_{\nu, e}}^{-1}$
  and, therefore, $S_{\nu, e}\, f \in \mathring{C}^{\infty} (\ov{E^{n+1}_{+}}).$
 Let $H^s_{loc} (E^{n+1}_{+})$ denote the set of functions defined
 in $E^{n+1}_{+}$ and belonging to the space $H^s (E^{n+1}_{a, b})$
 for each layer
${E^{n+1}_{a, b} = \{x=(x', y): x' \in E^n, a<y<b \}},$
$0<a<b<\infty$.
    This set equipped by the system of seminorms  $p_{s, a, b} (f) = \| f\|_{H^s
(E^{n+1}_{a, b})}$ is a Frechet space.

\begin{lemma} \label{lem: 3.1.1}
Let  $s$ be even and nonnegative and $\Re \nu \geq 0.$ Then the
estimate
\begin{equation}
p_{s, a, b} (f) \leq c \,  \|  f \|_{H^s_{\nu} \lr{E^{n+1}_{+}}},
\label{3.1.2}
\end{equation}
    where $0<a<b<\infty$ and $c$ is a positive constant independent of $f,$
    holds for any function
     $ f \in\mathring{C}^{\infty}_{\nu} \lr{E^{n+1}_{+}}.$
\end{lemma}

\begin{proof}
    By assumption, the function $g=S_{\nu} f$ belongs to the space
$\mathring{C}^{\infty} (\ov{E^{n+1}_{+}})$ and ${P_{\nu, e} g=f}.$
    Apply the following relation to this function (see \eqref{1.3.13}):
\begin{equation}
P_{\nu, e} \, g (x', y) = c_{\nu}\, y^{-\nu}
\int\limits_{-\infty}^{\infty} H^{(1)}_{\nu} (y \eta) \eta^{\nu}
(1 - i \eta)^{\frac{1}{2}- \nu} F \mbox{Cont} g (x', \eta)   \, d
\eta, \label{3.1.3}
\end{equation}
    where $F$ is the Fourier transformation with respect to the last
    variable and $\mbox{Cont}$ is the Whitney continuation operator acting from
 $E^{n+1}_{+}$ to $E^{n+1}.$ First, assign $s=0.$
 From the Cauchy--Bunyakovsky
  inequality, we obtain the estimate
\begin{multline}
\| P_{\nu, e}\, g \|^2_{L_2 \lr{E^{n+1}_{a, b}}} \leq c
\int\limits_{E^{n+1}_{a, b}} \int\limits_{|\eta| < 1}
|H^{(1)}_{\nu} (y \eta) \eta^{\nu} (1 - i \eta)^{\frac{1}{2}- \nu}
|^2   \, d \eta  \int\limits_{|\eta| < 1} | F \mbox{Cont} g |^2 \,
d \eta   \, d x  {}
\\
{} + c  \int\limits_{E^{n+1}_{a, b}}  | \int\limits_{|\eta| > 1}
H^{(1)}_{\nu} (y \eta) \, \eta^{\nu} (1 - i \eta)^{\frac{1}{2}-
\nu} F \mbox{Cont} g (x', \eta) \, d \eta|^2 \, dx. \label{3.1.4}
\end{multline}
The first internal integral of the first term of the right-hand
    part of \eqref{3.1.4} is bounded because, according to
    asymptotic relations for Hankel functions (see \cite{BE2}), the integrand function
    is either bounded or has a logarithmic singularity at the point $\eta=0.$
    Therefore, extending the integration limits and taking into
    account the Parseval equality, we estimate the first term by  $c\,
\|\mbox{Cont} g \|_{L_2 (E^{n+1}_{+})}.$ To estimate the second
term, one has to use the asymptotic behavior of the Hankel
function $H_{\nu}^{(1)} (z)$ as $z \to \infty.$ Its behaves as
follows (see \cite[p. 98]{BE2}):
\begin{equation}
H_{\nu}^{(1)} (z) = \frac{e^{i z}}{\sqrt{z}} \lr{c+ O(|z|^{-1})},
\label{3.1.5}
\end{equation}
 where $c$ is a constant. Using  the Parseval equality,
 we obtain that
$$
\int\limits_a^b \left| \int\limits_{|\eta| > 1} e^{i z}  (y
\eta)^{-\frac{1}{2}} \, \eta^{\nu} (1 - i \eta)^{\frac{1}{2}- \nu}
F \mbox{Cont} g (x', \eta)   \, d \eta\right|^2 \, dy \leq
 c \int\limits_{-\infty}^{\infty} | F \mbox{Cont} g |^2 \, d \eta
\leq c \int\limits_{- \infty}^{\infty} |  \mbox{Cont} g (x', y)
|^2 \, d y.
$$
    Thus, we obtain the desired estimate of the expression
    corresponding to the first term of \eqref{3.1.5}.
  To estimate the expression
    corresponding to the  remainder, we take into account  the Cauchy--Bunyakovsky
  inequality and obtain that
$$
    \begin{gathered}
\int\limits_a^b | \int\limits_{|\eta| > 1}   (y
\eta)^{-\frac{1}{2}} \eta^{\nu}  (1 - i \eta)^{\frac{1}{2}- \nu}
O(|y \eta|^{-1}) F \mbox{Cont} g (x', \eta)   \, d \eta|^2 \, dy
    \\
\leq c \int\limits_a^b  \int\limits_{|\eta| > 1} | y \eta |^{-2}
\, d \eta \int\limits_{|\eta| > 1}  | F \mbox{Cont} g |^2 \, d
\eta \,  dy \leq c \int\limits_{- \infty}^{\infty} |  \mbox{Cont}
g (x', y) |^2 \, d y.
    \end{gathered}
$$
    Hence, since $\|\mbox{Cont} g\|_{L_2 (E^n)} \leq c \,
\|g\|_{L_2 (E^{n+1}_{+})},$ it follows that inequality
\eqref{3.1.2} is proved for $s=0$.

Consider the case where $s$ is even and positive. First, we note
that
$$
p^2_{s, a, b} (f) \leq \sum\limits_{|\alpha'|+2 \alpha_{n+1} \leq
s} \|D^{\alpha'}_{x'} \, D^{2 \alpha_{n+1}}_y \, f \|^2_{L_2
(E^{n+1}_{a, b})} \leq
 c \sum\limits_{|\alpha'|+2 \alpha_{n+1} \leq s} \|D^{\alpha'}_{x'}
\, B^{\alpha_{n+1}}_y \, f \|^2_{L_2 (E^{n+1}_{a, b})},
$$
    where
$$
D_{x'}^{\alpha'} = \frac{\pr^{|\alpha'|}}{\pr x_1^{\alpha_1} \pr
x_2^{\alpha_2} \dots \pr x_n^{\alpha_n}},~ |\alpha'| = \alpha_1 +
\dots + \alpha_n,~\textrm{and}~B_y = \frac{\pr^2}{\pr y^2} +
\frac{2 \nu +1}{y} \frac{\pr }{\pr y}.
$$
  Since $S_{\nu, e} B = D^2 S_{\nu, e},$ it follows that
$$
p^2_{s, a, b} (f) \leq c \sum\limits_{|\alpha'|+2 \alpha_{n+1}
\leq s} \|S_{\nu}\, D^{\alpha'}_{x'} \, B^{ \alpha_{n+1}}_y \, f
\|^2_{L_2 (E^{n+1}_{+})} \leq
 c\,  \| S_{\nu, e} \, f \|_{H^s (E^{n+1}_{+})} = c \,
\|f\|^2_{H^s_{\nu} (E^{n+1}_{+})},
$$
    which completes the proof of the lemma.
\end{proof}
  By virtue of this lemma, the space $\mathring{C}^{\infty}_{\nu} ({E}^{n+1}_{+})$
    endowed with norm \eqref{3.1.1} is continuously embedded in
    the complete space $H^s_{loc} ({E}^{n+1}_{+}).$ Closing
$\mathring{C}^{\infty} ({E}^{n+1}_{+})$ with respect to norm
\eqref{3.1.1} and assuming that $s \geq 0$ and $\Re \nu \geq 0,$
 we obtain a space denoted by $H^s_{\nu}({E}^{n+1}_{+}).$

 Lemma \ref{lem: 3.1.1} intermediately implies the following assertion.

\begin{corollary} {\label{cor: 3.1.2}}
  If $s \geq 0$ and $\Re \nu \geq 0,$ then the Hilbert space $H^s_{\nu} ({E}^{n+1}_{+})$
  is continuously embedded in the Frechet space $H^s_{loc} ({E}^{n+1}_{+}).$
\end{corollary}

\begin{lemma}{\label{lem:3.1.2}}
    Let $s$ be even and nonnegative, $\Re \nu \geq 0,$ and $0<a<b<\infty.$
  Let $f$ belong to $H^s_{\nu} ({E}^{n+1}_{+})$ and its support is contained in the layer
   $E^{n+1}_{a, b}.$ Then $f \in H^s({E}^{n+1}_{+})$ and
\begin{equation}
c' \, \| f \|_{H^s (E^{n+1}_{+})} \leq  \| f \|_{H^s_{\nu}
(E^{n+1}_{+})} \leq c'' \, \| f \|_{H^s (E^{n+1}_{+})},
\label{3.1.6}
\end{equation}
    where $c'$ and $c''$ are positive constants independent of $f.$
\end{lemma}

\begin{proof}
    The fact that $f \in H^s ({E}^{n+1}_{+})$ and the left-hand
 inequality of \eqref{3.1.6} are proved within the proof of Lemma \ref{lem: 3.1.1}.
 To prove the right-hand
 inequality, we recall that
\begin{equation}
S_{\nu, e} = \mathcal{J}_{\nu - \frac{1}{2}, e}\, S_{\nu},
\label{3.1.7}
\end{equation}
    where
\begin{equation}
 \mathcal{J}_{\nu - \frac{1}{2}, e} \,f (x', y) =
 f (x', y) + \lr{\frac{1}{2}-\nu} \int\limits_y^{\infty}   f (x', t)
  \Phi \lr{\nu + \frac{1}{2}, 2; y-t} dt,
\label{3.1.8}
\end{equation}
\begin{equation}
S_{\nu} f (x', y) =  c_{\nu} \lr{y^{\nu + \frac{1}{2}} f (x', y) -
\int\limits_y^{\infty} t^{\nu + \frac{1}{2}}  f (x', t) \frac{\pr
}{\pr y} P_{\nu - \frac{1}{2}}^0 \lr{\frac{y}{t}} dt},
\label{3.1.9}
\end{equation}
  $\Phi$ is the degenerate hypergeometric function, and $P^0_{\mu}$ of the first-kind
  Legendre function.
  These functions, being kernels of integral operators, are smooth.
    Therefore, since $f(x', y)=0$ provided that $y<a$ or $y>b,$
    it follows that the norm of the function  $S_{\nu, e} f$ in the space $H^s (E^{n+1}_{+})$
    is estimated via the norm of the function $f$ in the same space,
    which completes the proof of the lemma.
\end{proof}
 Introduce the space $H^s_{\nu} (E^{n+1}_{0, b})$ for the layer
$E^{n+1}_{0, b},$ $0<b<\infty,$ as the closure with
    respect to norm \eqref{3.1.1} of the set of all
functions from $\mathring{C}^{\infty}_{\nu} ({E}^{n+1}_{+})$
    such that the support of each one is contained in $\ov{{E}^{n+1}_{0, b}}$.

\begin{lemma}{\label{lem: 3.1.2}}
 Let $s$ be even and nonnegative, $\Re \nu \geq 0,$ and $0<b< \infty.$
 Then norm \eqref{3.1.1} and the norm
\begin{equation}
 \| f \|_{H^s (E^{n+1}_{0, b})} = \| S_{\nu} f \|_{H^s (E^{n+1}_{0, b})}
\label{3.1.10}
\end{equation}
    are equivalent on the space $H^s_{\nu} (E^{n+1}_{0, b})$.
\end{lemma}

\begin{proof}
    Let $f \in \mathring{C}^{\infty}_{\nu} ({E}^{n+1}_{+})$
 and $\supp f \subset \ov{E^{n+1}_{+}}.$ By definition, this means that
  $S_{\nu, e}\, f \in \mathring{C}^{\infty}_{\nu}(\ov{E^{n+1}_{+}}).$
  Then $S_{\nu} f = \mathcal{J}_{\frac{1}{2}-\nu, e}\, S_{\nu, e}\, f \in
\mathring{C}^{\infty}_{\nu} (\ov{E^{n+1}_{+}})$ as well; this is
easily derived from \eqref{3.1.7}-\eqref{3.1.8}.
    Note that relation \eqref{3.1.8} holds for all complex values of $\nu.$
 The operators $\mathcal{J}_{\mu, e}$ and $\mathcal{J}_{-\mu, e}$ are inverse to each other
 on the space $\mathring{C}^{\infty}_{\nu}(\ov{E^{n+1}_{+}}).$
 Further, from the smoothness of the function $\Phi$, we easily
 obtain the estimate
$$
c' \, \|  \mathcal{J}_{\mu, e} g \|_{H^s (E^{n+1}_{+})} \leq \| g
\|_{H^s (\ov{E^{n+1}_{+}})} \leq c'' \, \| \mathcal{J}_{\mu, e} g
\|_{H^s (E^{n+1}_{+})}
$$
    valid for all $g \in \mathring{C}^{\infty}(\ov{E^{n+1}_{+}})$ such that
 $\supp g \subset \ov{{E}^{n+1}_{0,b}}.$
 The constants $c'$ and $c''$ depend on $b,$ but do not depend on
 the function $g.$ Thus, the equivalence of norms \eqref{3.1.1} and \eqref{3.1.10}
 is proved on a dense set of $f \in
\mathring{C}^{\infty}_{\nu} ({E}^{n+1}_{+}),$ $\supp f \subset
\ov{{E}^{n+1}_{0, b}},$ which completes the proof of the lemma.
\end{proof}
    Note that the equivalence of norms \eqref{3.1.1} and \eqref{3.1.10} in the spaces
     $H^s_{\nu}(E^{n+1}_{0, b})$ is proved under the assumption that $b$ is finite.
 For  $b=\infty$, this is not valid.
 One can show that, for $b=\infty$, norm \eqref{3.1.10} is subordinated to norm \eqref{3.1.1}.
  Therefore, if norm \eqref{3.1.1} is changed for norm \eqref{3.1.10} in the definition
  of the space $H^s (E^{n+1}_{+})$, then we obtain a wider space and it turns out that
  it is not embedded in $H^s_{loc} (E^{n+1}_{+}).$
  This is the reason to use norm \eqref{3.1.1} generated by the transmutation operator
  $S_{\nu, e}$ (instead of $S_{\nu}$).

    Let us find the relation between the introduced spaces and
    Kipriyanov function spaces from \cite{Kip1}. Let $F$ be the $(n+1)$-dimensional
  Fourier transformation
$$
F f(\xi)  = \int\limits_{E^{n+1}} f (x', y) e^{- i \langle \xi', x' \rangle - i y \eta} dx,
$$
    where $\xi = (\xi', \eta) \in E^{n+1},$ $\xi' \in E^n,$ $\langle
\xi', x' \rangle = \xi_1 x_1 + \dots +  \xi_n x_n.$
 Let $F_{\nu}$ be the
 Fourier--Bessel (or Fourier--Hankel)
 transformation
$$
F_{\nu} f(\xi)  =  \int\limits_{E^{n+1}_{+}} f (x', y) e^{- i
\langle \xi', x' \rangle} j_{\nu} (y \eta) y^{2 \nu +1 } \, dx.
$$
    The spaces $H^s_{\nu, +} (E^{n+1}_{0, b}),$ where $s \geq 0$ and
 $\nu\geq -\dfrac{1}{2},$ are  defined the closure of the set of even functions
  $f \in \mathring{C}^{\infty} (\ov{E^{n+1}_{+}})$ such that
   $\supp f \subset \ov{E^{n+1}_{0, b}}$ (the set of such functions is denoted by
    $\mathring{C}^{\infty}_{+}(\ov{E^{n+1}_{0, b}})$) with respect
    to the norm
\begin{equation}
\| f \|_{H^s_{\nu, +} (E^{n+1}_{+})} =   \lr{ ~
\int\limits_{E^{n+1}_{+}} |F_{\nu} f (\xi) |^2 (1+|\xi|^2)^s
{\eta}^{2 \nu +1 } \, d \xi}^{\frac{1}{2}}. \label{3.1.11}
\end{equation}
 Multiply applying Lemma \ref{lem: 3.1.2}, we arrive at the following assertion.

\begin{lemma}{\label{lem: 3.1.4}}
    Let $s$ be even and nonnegative and $\nu$ be real and nonnegative.
    Then the space $H^s_{\nu, +} (E^{n+1}_{0, b})$ is continuously
    embedded in $H^s_{\nu} (E^{n+1}_{0, b}).$
 For $\nu \neq 1,3,5,\dots,$ the set $H^s_{\nu, +} (E^{n+1}_{0, b})$ forms a proper subspace
 of the space $H^s_{\nu} (E^{n+1}_{0, b})$ and the induced norm on
  $H^s_{\nu, +} (E^{n+1}_{0,b})$ is equivalent to its proper one.
\end{lemma}

    Now, consider internal embedding theorems of the introduced spaces $H^s_{\nu}.$

\begin{theorem} \label{teo: 3.1.1}
 Let $\Re \nu \geq 0$ and $s'>s \geq 0.$
 Then the space $H^{s'}_{\nu} (E^{n+1}_{+})$ is continuously embedded in
  $H^s_{\nu} (E^{n+1}_{+})$ and the corresponding inequality for
  the norms holds.
\end{theorem}

This assertion is the direct corollary from the definition of the
spaces $H^s_{\nu}$ and from embedding theorems for the spaces
$H^s.$

\begin{theorem}{\label{teo: 3.1.2}}
 Let $\Re \nu \geq 0$ and $s'>s \geq 0.$ Let $f_k,$ $k=1,2,\dots,$
 be a bounded sequence of functions from $H^{s'}_{\nu} (E^{n+1}_{+})$
 such that $\supp f_k \subset \mathcal{Y},$ where $\mathcal{Y}$ is a bounded set
 of $E^{n+1}_{+}.$ Then there exists a subsequence converging with respect to the norm
 of the space $H^s_{\nu} (E^{n+1}_{+}).$
\end{theorem}

\begin{proof}
One can construct a  sequence of functions  $\varphi_k (x) \in
\mathring{C}^{\infty}_{\nu} ({E}^{n+1}_{+})$ such that $\supp
\varphi_k \subset \mathcal{Y}'$ and $\|f_k - \varphi_k
\|_{H^{s'}_{\nu} (E^{n+1}_{+})} \to 0$ as $k \to \infty,$ where
$\mathcal{Y}'$ is a compactum in $\ov{E^{n+1}_{+}}.$
    Then the sequence of norms $\|\varphi_k \|_{H^{s'}_{\nu}
(E^{n+1}_{+})} = \|S_{\nu, e} \varphi_k \|_{H^{s'}_{\nu}
(E^{n+1}_{+})}$ is bounded.
 Due to the complete continuity of the embedding of the spaces $H^s$
  (see \cite{66}), there exists a function $g \in H^s (E^{n+1}_{+})$ such that
  a subsequence $S_{\nu, e} \varphi_{k_m}$ converges to it.
Hence, it is fundamental in  $H^s (E^{n+1}_{+}).$
    Then the sequence $\{f_{k_m}\}$  is fundamental in  $H^s_{\nu}
(E^{n+1}_{+}).$ By virtue of the completeness of the space
$H^s_{\nu},$ this subsequence is convergent, which completes the
proof.
\end{proof}

\subsection{Results about multipliers}\label{sec10.2}

    In this section, we find sufficient conditions on the function
$a(x),$ providing the continuity of the map $f \to a f$ in spaces
$H^s_{\nu} (E^{n+1}_{+}).$ We provides estimates of the norm of
this map as well.

\begin{lemma}{\label{lem: 3.2.1}}
 Let $a(x) \in C^{\infty}  (\ov{E^{n+1}_{+}})$ and
\begin{equation}
D^k_y a(x) = 0  \  \mbox{for} \  y = 0, \  k =1,2, \dots
\label{3.2.1}
\end{equation}
    Then, for $\Re \nu \geq 0,$ the function $a f$ belongs to the
    set
$\mathring{C}^{\infty}_{\nu} ({E}^{n+1}_{+})$ if $f \in
\mathring{C}^{\infty}_{\nu} ({E}^{n+1}_{+}).$
\end{lemma}

\begin{proof}
    Let us show that if $g \in \mathring{C}^{\infty}
(\ov{E^{n+1}_{+}}),$ then the function $S_{\nu, e}\, (a P_{\nu,
e}\, g)$ belongs to the same class. Since $S_{\nu, e} =
I^{\frac{1}{2}-\nu}_{e} S^{\nu- \frac{1}{2}}_{\nu}$, $P_{\nu, e} =
P^{\frac{1}{2}-\nu}_{\nu} I^{\nu- \frac{1}{2}}_{e}$, and the
operators $I^{\mu}_e$ map $\mathring{C}^{\infty}
(\ov{E^{n+1}_{+}})$ onto itself for each complex  $\mu,$ it
follows that it suffices to show that $S^{\nu- \frac{1}{2}}_{\nu}
(a P^{\frac{1}{2}-\nu}_{\nu} g) \in \mathring{C}^{\infty}
(\ov{E^{n+1}_{+}}).$
    From the definition of the transmutation operators
$S^{\nu- \frac{1}{2}}_{\nu}$ and $P^{\frac{1}{2}-\nu}_{\nu},$
where $N$ is a positive integer, we obtain that
$$
    \begin{gathered}
S_{\nu}^{\nu - \frac{1}{2}} \lr{a P_{\nu}^{ \frac{1}{2} - \nu} f}
(x', y) = \frac{(-1)^{N+1} \, 2^{1-N}}{\Gamma \lr{\nu +
\frac{1}{2}} \,  \Gamma \lr{N-\nu + \frac{1}{2}}}\, \frac{\pr}{\pr
y} \int\limits_y^{\infty} (t^2-y^2)^{\nu - \frac{1}{2}} a (x',
t)\, t  {}
 \\
{}\times \int\limits_t^{\infty} (\tau^2-t^2)^{N-\nu - \frac{1}{2}}
\,  \lr{\frac{\pr}{\pr \tau} \frac{1}{\tau}}^N g(x', \tau) \, d
\tau d t
    \\
= \frac{(-1)^{N+1} \, 2^{1-N}}{\Gamma \lr{\nu + \frac{1}{2}} \,
\Gamma \lr{N-\nu + \frac{1}{2}}}\, \frac{\pr}{\pr y}
\int\limits_y^{\infty} \lr{\lr{\frac{\pr}{\pr \tau}
\frac{1}{\tau}}^N g(x', \tau)} \,
 \int\limits_y^{\tau} (\tau^2-y^2)^{\nu - \frac{1}{2}} \, t \,
 (\tau^2-t^2)^{N-\nu - \frac{1}{2}}   a (x', t)\, d t d \tau.
    \end{gathered}
$$
    In the last internal integral, change the variable as follows:
     $t = \sqrt{y^2 + z (\tau^2-y^2)}.$ Then
$$
\int\limits_y^{\tau} (\tau^2-y^2)^{\nu - \frac{1}{2}} \, t \,
(\tau^2-t^2)^{N-\nu - \frac{1}{2}}   a (x', t)\, d t =
\frac{(\tau^2-y^2)^N }{2}
 \int\limits_0^1 z^{\nu - \frac{1}{2}}  z^{N-\nu -
\frac{1}{2}} a (x', \sqrt{y^2 + z (\tau^2-y^2)}) \, d z.
$$
    Integrating this by parts, we find that
$$
     \begin{gathered}
S_{\nu}^{- \frac{1}{2} + \nu} \lr{a P_{\nu}^{ \frac{1}{2} - \nu}
g} (x', y) = \frac{-1}{\Gamma \lr{\nu + \frac{1}{2}} \,  \Gamma
\lr{N-\nu + \frac{1}{2}}}\, \frac{\pr}{\pr y}
\int\limits_y^{\infty} g(x', \tau)
    \\
\times \frac{\pr^N}{(\pr \tau^2)^N} \left(  (\tau^2-y^2)^N
\int\limits_0^1 z^{\nu - \frac{1}{2}}  z^{N-\nu - \frac{1}{2}} a
(x', \sqrt{y^2 + z (\tau^2-y^2)}) \, d z \right)  d \tau
    \\
= \frac{-1}{\Gamma \lr{\nu + \frac{1}{2}} \,  \Gamma \lr{N-\nu +
\frac{1}{2}}}\, \frac{\pr}{\pr y} \int\limits_y^{\infty} g(x',
\tau) \int\limits_0^1 \sum\limits_{k=0}^{N} 2^{k-N} \binom{N}k
\frac{N!}{(N-k)!}
    \\
\times (\tau^2-y^2)^{N-k} z^{\nu + N - k - \frac{1}{2}}
(1-z)^{N-\nu - \frac{1}{2}}\left. \lr{\frac{\pr}{\lambda \pr
\lambda}}^{N-k} a(x', \lambda)\right|_{\lambda= \sqrt{y^2 + z
(\tau^2-y^2)}} dz d \tau,
     \end{gathered}
$$
    where $\binom{N}k$ denote the binomial coefficients.
    Differentiating this with respect to the parameter $y,$ we obtain the relation
\begin{equation}
S_{\nu}^{\nu- \frac{1}{2}} \lr{a P_{\nu}^{ \frac{1}{2} - \nu} g}
(x', y) =  a (x', y) g (x', y) + \int\limits_y^{\infty} g(x',
\tau) a_{\nu} (x', y, \tau) \, d \tau, \label{3.2.2}
\end{equation}
    where
\begin{multline}
    a_{\nu} (x', y, \tau) = \frac{-1}{\Gamma \lr{\nu +\frac{1}{2}}
    \, \Gamma \lr{N-\nu +\frac{1}{2}}} \frac{\pr}{\pr y}
    \sum\limits_{k=0}^{N}   \binom{N}k \frac{2^{k-N} \, N!}{(N-k)!}
    (\tau^2-y^2)^{N-k}  \\
\left. \times \int\limits_0^1
    z^{\nu + N - k - \frac{1}{2}} (1-z)^{N-\nu - \frac{1}{2}}
    \lr{\frac{\pr}{\lambda \pr \lambda}}^{N-k} a(x',
    \lambda)\right|_{\lambda= \sqrt{y^2 + z (\tau^2-y^2)}} dz.
    \label{3.2.3}
\end{multline}
    It follows from \eqref{3.2.3} and Condition \eqref{3.2.1} that  the function
 $a_{\nu} (x', y, \tau)$ is infinitely differentiable for $x' \in
E^n,$ $y  \geq 0,$ and $ \tau \geq 0.$
    Therefore, the left-hand
    part of relation \eqref{3.2.2} is infinitely differentiable as well.
    The compactness of its support is obvious, which completes the proof of the
    lemma.
\end{proof}
    Now, impose additional restrictions of the function $a(x)$.
 Assume that there exists a finite $R$ such that
\begin{equation}
D_y \, a (x', y) = 0, \  y \geq  R,
\label{3.2.4}
\end{equation}
    and
\begin{equation}
\sup\limits_{x \in E_{+}^{n+1}} \left| D^{\alpha'}_{x'}
\lr{\frac{1}{y} D_y}^{\alpha_{n+1}} a(x)\right| := M^{\alpha',
\alpha_{n+1} } := M^{\alpha}  < \infty  \label{3.2.5}
\end{equation}
    for all multi-indices
     $\alpha$.

Under these assumptions, properties of the function $a_{\nu} (x',
y, \tau)$ can be refined. The estimate
\begin{equation}
 \left| D^{\alpha'}_{x'}  \lr{\frac{1}{y} D_y}^{l} \lr{\frac{1}{\tau} D_{\tau}}^{m}
  a_{\nu} (x', y, \tau) \right| \leq c \, (1+\tau)^{- 2 \Re \nu -1-m}
   \max\limits_{k \leq N+l+m+1}  M^{\alpha', k},
\label{3.2.6}
\end{equation}
 where $c$ is a positive constant independent on $x',$ $y,$ $\tau,$
 and the function $a,$ holds.
 Indeed, the right-hand
 side of relation \eqref{3.2.3} contains only derivatives of  $a$
 with respect to the last variable, but does contain the function  $a$  itself.
  Then, it follows from \eqref{3.2.4} that
\begin{equation}
a_{\nu} (x, y, \tau) = 0, \  \tau \geq y \geq  R.
\label{3.2.7}
\end{equation}
    Hence, estimate \eqref{3.2.6} holds for $\tau \geq y\geq R$ a fortiori.
  Obviously, it holds for $0 \leq y \leq \tau \leq R$ as well.
  Now,  let $\tau>R>y$;   then, due to \eqref{3.2.4}, we have the
  relations
$$
 \begin{gathered}
\left| \left. \int\limits_0^1 z^{\nu + N - k - \frac{1}{2}}
(1-z)^{N-\nu - \frac{1}{2}}  \lr{\frac{\pr}{\lambda \pr
\lambda}}^{N-k} a(x', \lambda)\right|_{\lambda= \sqrt{y^2 + z
(\tau^2-y^2)}} dz \right|
    \\
= \left| \left. \int\limits_0^{\frac{R^2-y^2}{\tau^2-y^2}} z^{\nu
+ N - k - \frac{1}{2}}  (1-z)^{N-\nu - \frac{1}{2}}
\lr{\frac{\pr}{\lambda \pr \lambda}}^{N-k} a(x',
\lambda)\right|_{\lambda= \sqrt{y^2 + z (\tau^2-y^2)}} dz \right|
    \\
\leq \frac{1}{ |\nu + N - k + \frac{1}{2}|}
\lr{\frac{R^2-y^2}{\tau^2-y^2}}^{\Re \nu + N - k + \frac{1}{2}}
\sup\limits_{y \geq 0} \left|  \lr{\frac{\pr}{y \pr y}}^{N-k}
a(x', y) \right|
 \end{gathered}
$$
 leading to inequality \eqref{3.2.6}.

 Relation \eqref{3.2.3} yields a similar representation for the operator
  $S_{\nu, e}\, a\, P_{\nu, e}.$
  Since $S_{\nu, e}=I^{\frac{1}{2} - \nu} S^{\nu - \frac{1}{2}}_{\nu}$ and
  $P_{\nu, e}=P_{\nu}^{\frac{1}{2} - \nu} I^{\nu - \frac{1}{2}}_{e},$ it follows
  from \eqref{3.2.3} that
\begin{equation}
S_{\nu, e} \lr{a \, P_{\nu, e} \, g} (x', y) = I_{e}^{\frac{1}{2}
- \nu} \lr{a \, I_{e}^{\nu - \frac{1}{2} } g} (x', y) +
I_{e}^{\frac{1}{2} - \nu}  \int\limits_y^{\infty} a_{\nu} (x', y,
\tau) I_{e}^{\nu - \frac{1}{2} } g (x', \tau)\, d \tau,
\label{3.2.8}
\end{equation}
    where $g \in \mathring{C}^{\infty} (\ov{E^{n+1}_{+}}).$
    The first term of the right-hand
    side is considered in Lemma \ref{lem:1.2.3}.
    Therefore, it suffices to estimate the second term.
    For $\dfrac{1}{2}< \Re \nu < N+\dfrac{1}{2}$, we have the
    relation
\begin{multline}
I_{e}^{\frac{1}{2} - \nu}  \int\limits_y^{\infty} a_{\nu} (x', y,
\tau) I_{e}^{\nu - \frac{1}{2} } g (x', \tau)\, d \tau
= (-1)^N e^y I^{\frac{1}{2} - \nu+N} D^N_y \int\limits_y^{\infty}
 e^{-y} a_{\nu} (x', y, \tau) I_{e}^{\nu - \frac{1}{2} } g (x', \tau)\, d \tau  \\
\left. = \sum\limits_{k=0}^{N-1} (-1)^k  I^{\frac{3}{2} - \nu+k}
 \lr{e^y D^k_y \lr{e^{-\tau} a_{\nu} (x', y, \tau)}\right|_{\tau=y}
  I^{\nu - \frac{1}{2} } g (x', y)}  \\
+ (-1)^N  I_e^{\frac{1}{2} - \nu+N} \int\limits_y^{\infty} e^{y}
D^N_y \lr{e^{-y} a_{\nu} (x', y, \tau)} I_{e}^{\nu - \frac{1}{2} }
g (x', \tau)\, d \tau. \label{3.2.9}
\end{multline}
    From relation \eqref{3.2.3}, we see that  the function
$$
e^y D^k_y \left. \lr{e^{-y} a_{\nu} (x', y, \tau)}\right|_{\tau=y}, \  k = 0, \dots, N-1,
$$
    satisfies all conditions of Corollary \ref{cor:1.2.2}.
    Therefore,
\begin{multline}
\| \left. \sum\limits_{k=0}^{N-1} (-1)^k  I^{\frac{3}{2} - \nu+k}
 \lr{e^y D^k_y \lr{e^{-y} a_{\nu} (x', y, \tau)} \right|_{\tau=y}  I_e^{\nu - \frac{1}{2} }
  g (x', y)} \|_{H^s (E_{+}^{n+1})}  \\
\leq c \sum\limits_{\substack{|\alpha'| \leq  s \\
k \leq 3 N +s+1}} M^{\alpha', k} \|g\|_{H^s (E^{n+1}_{+})},
\label{3.2.10}
\end{multline}
 where $c$ is a positive constant independent of the functions $a$ and $g.$

To estimate the last integral of \eqref{3.2.9}, we note that
 the operators $I_e^{\frac{1}{2}-\nu+N}$ and $I_e^{\nu-\frac{1}{2}}$
 belong to the class $L (H^s (E^{n+1}_{+}), H^s  (E^{n+1}_{+}))$
 because $\Re \Big(\dfrac{1}{2}-\nu+N\Big)>0$ and $\Re \nu -\dfrac{1}{2}> 0$
   (see Lemma \ref{lem:1.2.2}).
 Estimate \eqref{3.2.6} shows that the integral operator
$$
w \to \int\limits_y^{\infty} w(x', \tau) \, e^y \, D_y^N
\lr{e^{-y} a_{\nu} (x', y, \tau)}  d \tau
$$
 belongs to the class ${L} \lr{ H^s (E^{n+1}_{+}), H^s
(E^{n+1}_{+})}$ for each $s$ and its norm does not exceed
$$
c \sum\limits_{  k \leq 3 N +s+1,~ |\alpha'| \leq  s} M^{\alpha', k},
$$
    where  $c$ is a positive constant independent of the function $a.$

 The same estimates hold for  $0 \leq \Re \nu \leq\dfrac{1}{2}.$
 To prove them, it suffices to estimate the norm of the second term in relation \eqref{3.2.8}.
 Integrating by parts, we find that
$$
I_e^{\frac{1}{2} {-} \nu}\!\! \int\limits_y^{\infty}   I_e^{\nu
{-} \frac{1}{2}} g (x'\!, \tau)  a_{\nu} (x'\!, y\!, \tau)  d \tau
=  I_e^{\frac{1}{2} {-} \nu}  a_{\nu} (x'\!, y\!, y) I_e^{\nu {+}
\frac{1}{2}} g {+}
 I_e^{\frac{1}{2} {-} \nu}\!\! \int\limits_y^{\infty}
  e^{{-} \tau} D_{\tau} \lr{ e^{ \tau} a_{\nu} (x'\!, y\!, \tau)} I_e^{\nu {+}
   \frac{1}{2}} g (x'\!, \tau)  d \tau.
$$
    Now, it easily follows from inequality \eqref{3.2.6} that the
 $ H^s (E^{n+1}_{+})$-norm of this expression does not exceed
$c \sum\limits_{k \leq 4+s,~|\alpha'|\leq s} M^{\alpha', k},$
    where  $c$ is a positive constant independent of the function $a.$

Thus, the following assertion is proved.

\begin{theorem}\label{teo: 3.2.1}
Let $a \in C^{\infty} (\ov{E^{n+1}_{+}})$ and Conditions
 \eqref{3.2.1}, \eqref{3.2.4}, and \eqref{3.2.5} be satisfied.
 Then, for each complex $\nu$ from the half-plane
  $\{\Re \nu \geq 0\},$ each  nonnegative
$s,$ and each positive ${R}$ such that it does not exceed a
positive number $R_0,$ the operator $S_{\nu, e} \, a \, P_{\nu,
e}$ continuously maps the space $ H^s (E^{n+1}_{+})$ into itself
and the estimate
\begin{equation}
\| S_{ \nu, e} \, a \, P_{ \nu, e} f \|_{H^s (E^{n+1}_{+})}
 \leq c \, \|f\|_{H^s (E^{n+1}_{+})} \sum\limits_{\substack{|\alpha'| \leq  s \\
\alpha_{n+1} \leq 3 N +s+1}} M^{\alpha', \alpha_{n+1} }
\label{3.2.11}
\end{equation}
 holds, where $N$ is the least positive integer satisfying the inequality
  $\Re \nu <N+1.$ The positive constant $c$ depends on $\nu, n, s,$ and $R_0,$ but does not
  depend on the functions  $a$ and $f$  and on $R.$
\end{theorem}

\begin{corollary}\label{cor: 3.2.1}
 Under the assumptions of the theorem, the estimate
    \begin{equation}
    \|  a  f \|_{H^s_{\nu} (E^{n+1}_{+})}
     \leq c  \sum\limits_{\substack{|\alpha'| \leq  s \\
\alpha_{n+1} \leq 3 N +s+1}} M^{\alpha', \alpha_{n+1} }
\|f\|_{H^s_{\nu} (E^{n+1}_{+})}
    \label{3.2.12}
    \end{equation}
    holds.
\end{corollary}

\begin{corollary}\label{cor: 3.2.2}
    Let  the assumptions of the theorem be satisfied and
$$
\sup\limits_{x \in E^{n+1}_{+}}  \left| D_{x'}^{\alpha'}
\lr{\frac{1}{y} D_y}^{\alpha_{n+1}} \frac{a(x', y)}{y} \right| =
\widetilde{M}^{\alpha', \alpha_{n+1} } < \infty.
$$
    Then  the estimate
\begin{equation}
\|  a \, D_y \,  f \|_{H^s_{\nu} (E^{n+1}_{+})} \leq c
  \sum\limits_{\substack{|\alpha'|
 \leq  s \\
\alpha_{n+1} \leq 3 N +s+1}} \widetilde{M}^{\alpha', \alpha_{n+1}
} \|f\|_{H^{s+1}_{\nu} (E^{n+1}_{+})} \label{3.2.13}
\end{equation}
    holds.
\end{corollary}

\begin{proof}
    The relation
$$
y \, D_y \, P_{\nu}^{\frac{1}{2} - \nu} = P_{\nu}^{\frac{1}{2} -
\nu} \lr{y D_y-2 \nu}
$$
    easily follows from the definition of the transmutation operator
    (see \eqref{1.1.7}-\eqref{1.1.8}).
    Then, similarly to \eqref{3.2.8}, we have the following relation
    for the function $S_{\nu, e} f$:
$$
S_{\nu, e} \lr{a \, D_y \, P_{\nu, e}\, g } = I_e^{\frac{1}{2} -
\nu } \lr{a (D_y - \frac{2 \nu}{y}) I^{\nu - \frac{1}{2}} } +
 I_e^{\frac{1}{2} - \nu} \int\limits_y^{\infty} \lr{\tau D_{\tau}
- 2 \nu}   \widetilde{a}_{\nu}  (x', y, \tau) I_e^{\nu -
\frac{1}{2}}  g (x', \tau)  \, d \tau.
$$
 Here, the function $\widetilde{a}_{\nu}$ is defined by relation \eqref{3.2.3}, where
 the function $a$ is changed for the function $\dfrac{1}{y} \, a (x', y).$
 Integrating by parts, we reduce the last relation to the form
$$
    \begin{gathered}
S_{\nu, e} \lr{a \, D_y \, P_{\nu, e} \, g } = I_e^{\frac{1}{2} -
\nu }  \lr{a \, I_e^{\nu - \frac{1}{2}} \, D_y \, g  } - 2 \nu
I_e^{\frac{1}{2} - \nu } \lr{ \frac{a}{y} I^{\nu - \frac{1}{2}} g}
-
 2 \nu \, I_e^{\frac{1}{2} - \nu} \int\limits_y^{\infty} I^{\nu -
\frac{1}{2}}  g (x', \tau)  \widetilde{a}_{\nu}  (x', y, \tau) \,
d \tau
    \\
{}- I_e^{\frac{1}{2} - \nu} \lr{ y \, \widetilde{a}_{\nu}  (x', y,
\tau) I_e^{\nu - \frac{1}{2}} g} -
  I_e^{\frac{1}{2} - \nu} \int\limits_y^{\infty} \tau
\widetilde{a}_{\nu}  (x', y, \tau) I^{\nu - \frac{1}{2}}  D_{\tau}
\, g (x', \tau) \, d \tau.
    \end{gathered}
$$
    It is easy to see that each term satisfies conditions providing a possibility to apply
    the arguing scheme from the proof of the previous theorem.
    This implies the validity of estimate \eqref{3.2.13}, which
    completes the proof of the corollary.
\end{proof}

\subsection{Weight traces}\label{sec10.3}

In this section, we introduce weight traces and prove direct
theorems about weight traces.

As above, define the weight function $\sigma_{\nu} (y)$ by the
relations
\begin{equation*}
\sigma_{\nu} (y) = \left\{
\begin{array}{lll}
y^{2 \nu}  & \mbox{if} \  \Re \nu>0, \\
 \dfrac{1}{ - \ln y}  & \mbox{if} \   \nu=0, \\
1  & \mbox{if} \  \Re \nu<0.
\end{array}
\right.
\end{equation*}
    The special case of imaginary values of the parameter $\nu$ is separately considered
    below.

 In the classical sense, the \emph{weight $\sigma_{\nu}$-trace}
 of a function is the limit
$$
\left. \sigma_{\nu} f \right|_{y=0} = \lim\limits_{y \to + 0 }
\sigma_{\nu} (y) f (x', y)  = \psi(x').
$$
 Let us show that functions $f \in \mathring{C}_{\nu}^{\infty}
(E_{+}^{n+1})$ possess $\sigma_{\nu}$-traces
 and they belong to the space $\mathring{C}^{\infty} (E^{n}).$
 If $f \in \mathring{C}_{\nu}^{\infty} (E_{+}^{n+1}),$ then, by the definition of this space,
 there exists  $g \in\mathring{C}_{\nu}^{\infty} (\ov{E_{+}^{n+1}})$
 such that $f=P_{\nu, e} g.$ Then $g = S_{\nu, e} f.$ In Sec.
 \ref{sec4.3}, for the one-dimensional
 case, the following relation is proved:
\begin{equation}
\lim\limits_{y \to + 0 } \sigma_{\nu} (y) P_{\nu}^{\frac{1}{2} - \nu } f_1 (y) = \left\{
\begin{array}{ll}
\dfrac{1}{2 \nu} f_1 (0)  & \mbox{if} \  \Re \nu>0, \\
f_1 (0) & \mbox{if} \   \nu=0.
\end{array}
\right.
\label{3.3.1}
\end{equation}
    Obviously, the same result holds in the multidimensional case as well.
    From this result, using the relation
 $P_{\nu, e}= P_{\nu}^{\frac{1}{2}-\nu} I_e^{\nu- \frac{1}{2}},$
 we derive that
\begin{equation}
\left. \sigma_{\nu} f \right|_{y=0} = \left\{
\begin{array}{ll}
\dfrac{1}{2 \nu}  I_e^{\nu - \frac{1}{2}} g (x', y) \Big|_{y=0}  & \mbox{if}
 \  \Re \nu>0, \\
 I_e^{ - \frac{1}{2}} g (x', y) \Big|_{y=0}  & \mbox{if} \
\nu=0.
\end{array}
\right.
\label{3.3.2}
\end{equation}
   Thus, $\left. \sigma_{\nu} f \right|_{y=0} \in
\mathring{C}^{\infty} (E^n)$ if $f \in \mathring{C}_{\nu}^{\infty}
(E_{+}^{n+1}).$

According to  relation \eqref{3.3.2}, we have to study traces of
fractional integrals $I^{\mu}_e$ in spaces $H^s.$
    Such results are used in the fifth section as well.

Let $g \in \mathring{C}^{\infty} (\ov{E^{n+1}_{+}})$ and
$\mbox{Cont} g$ be the Whitney continuation of this function to
the whole space $E^{n+1}.$
    Then, for each complex $\mu$ and each positive integer $k$,
    the relation
$$
\left. D_y^k \, I^{\mu}_e \, g (x', y) \right|_{y=0} =   \left.
D^k \, I^{\mu}_e \, \mbox{Cont}  g (x', y) \right|_{y=0} = \psi_k
(x'),
$$
    where $\psi_k \in \mathring{C}^{\infty} (E^n)$, holds.
 In Fourier images, this relation has the form
$$
\frac{1}{2 \pi} \int\limits_{- \infty}^{\infty} (1 - i \eta)^{-
\mu} i^k \eta^k F \mbox{Cont}   g (\xi', \eta) \, d \eta = F'
\psi_k (x'),
$$
    where $F'$ denotes the  Fourier transformation with respect to the initial $n$ variables.
 According to   the Cauchy--Bunyakovsky
  inequality, we obtain that
\begin{multline}
\| (1 + |\xi'|^2)^s F'  \psi_k (x')  \|_{L_2(E^n)}^2
    \\
     \leq  c
\int\limits_{E^n} (1 + |\xi'|^2)^s   \int\limits_{-
\infty}^{\infty} \frac{(1+\eta^2)^{- \Re \mu} \eta^{2 k}}{ (1 +
|\xi|^2)^{s'}} \, d \eta \int\limits_{- \infty}^{\infty} (1 +
|\xi'|^2)^{s'} \left| F \mbox{Cont}  g (\xi) \right|^2 \, d \eta d
\xi'. \label{3.3.3}
\end{multline}

\begin{lemma} \label{lem: 3.3.1.}
    If $s'>k- \Re \mu + \dfrac{1}{2} > 0,$ then
\begin{equation}
\int\limits_{- \infty}^{\infty} \frac{(1+\eta^2)^{- \Re \mu}
\eta^{2 k}}{ (1 + |\xi|^2)^{s'}}  \, d \eta \leq c \, (1 +
|\xi'|^2)^{k-s'+\frac{1}{2} - \Re \mu}, \label{3.3.4}
\end{equation}
    where $c$ is a positive constant independent of $\xi' \in E^n.$
\end{lemma}

\begin{proof}
    Decompose the integral from \eqref{3.3.4} into two ones. For
    the first one, we have the inequality
$$
 \int\limits_{|\eta|>1} \frac{(1+\eta^2)^{- \Re \mu} \eta^{2 k}}{ (1 + |\xi|^2)^{s'}} \,
 d \eta \leq c_1
  \int\limits_{- \infty}^{\infty} \frac{ \eta^{2 k - \Re \mu}}{ (1 + |\xi|^2)^{s'}} \, d \eta
   = c_2 \, (1 + |\xi'|^2)^{k-s' - \Re \mu +\frac{1}{2}}.
$$
    The second integral admits the simple estimate
$$
\int\limits_{|\eta|<1} \frac{(1+\eta^2)^{- \Re \mu} \eta^{2 k}}{
(1 + |\xi|^2)^{s'}} \, d \eta \leq c_3 \, (1 + |\xi'|^2)^{-s'},
$$
    which completes the proof of the lemma.
\end{proof}
    Combining inequalities \eqref{3.3.3} and \eqref{3.3.4}, we
    obtain that
$$
    \begin{gathered}
\| \psi_k \|^2_{H^s(E^n)} \leq  c \int\limits_{E^{n+1}}  (1 +
|\xi'|^2)^{s-s' + k - \Re \mu +\frac{1}{2}} (1 + |\xi|^2)^{-s'}
\left| F \mbox{Cont}  g (\xi) \right|^2 \, d \xi
    \\
\leq c \int\limits_{E^{n+1}}  (1 + |\xi|^2)^{s + k - \Re \mu
+\frac{1}{2}}   \left| F \mbox{Cont}  g (\xi) \right|^2 \, d \xi
\leq c \, \| g \|_{H^{s + k - \Re \mu +\frac{1}{2}}
(E^{n+1}_{+})}.
    \end{gathered}$$
Hence, the following assertion is proved.

\begin{theorem} \label{teo: 3.3.1}
    Let $s>k- \Re \mu + \dfrac{1}{2} > 0.$ Then, for each
$g \in H^s (E_{+}^{n+1}),$ there exists a trace of the function
$D^k_y I_e^{\mu} g$ on the hyperplane $\{y=0\}$ such that it
belongs to the space $H^{s - k + \Re \mu - \dfrac{1}{2}} (E^{n}).$
   The estimate
\begin{equation}
\|\left. D^k_y I^{\mu}_e g \right|_{y=0} \|_{H^{s - k + \Re \mu -
\frac{1}{2}} (E^{n})}  \leq c \, \| g \|_{H^{s} (E^{n+1}_{+})},
\label{3.3.5}
\end{equation}
    where $c$ is a positive constant independent on the function $g,$ holds.
\end{theorem}

Getting back to study weight traces, we assume that $f \in
\mathring{C}_{\nu}^{\infty} (E_{+}^{n+1}).$
 Applying relation \eqref{3.3.1} to the function $B^k f$
 belonging to $\mathring{C}_{\nu}^{\infty} (E_{+}^{n+1})$ as well,
 we see that
\begin{equation}\label{3.3.6}
\begin{aligned}
\left. \sigma_{\nu} B^k f \right|_{y=0} &= \left\{
\begin{array}{ll}
\dfrac{1}{2 \nu} I_e^{\nu - \frac{1}{2}} S_{\nu, e} \, B^k f \Big|_{y=0}  & \mbox{if}
 \  \Re \nu>0, \\
 I_e^{ - \frac{1}{2}} S_{0, e} B^k f \Big|_{y=0} &
\mbox{if} \   \nu=0
\end{array}
\right. \\
&= \left\{
\begin{array}{ll}
\dfrac{1}{2 \nu}  D_y^{2 k} I_e^{\nu - \frac{1}{2}} S_{\nu, e} f \Big|_{y=0}  & \mbox{if}
 \  \Re \nu>0, \\
 D_y^{2 k} I_e^{ - \frac{1}{2}} S_{0, e}  f \Big|_{y=0}
& \mbox{if} \   \nu=0.
\end{array}
\right.
\end{aligned}
\end{equation}
    The relation $S_{\nu, e} B^k= D_y^{2 k} S_{\nu,e}$ is used here.
    Combining \eqref{3.3.1} and \eqref{3.3.6} with Theorem \ref{teo: 3.3.1} and
    the definition of the norm of the space $H^s_{\nu}(E^{n+1}_{+}),$
    we arrive at the following theorem on weight traces.

\begin{theorem} \label{teo: 3.3.2}
    Let $\Re \nu > 0$ or $\nu = 0$. Let $s> 2 k - \Re \nu +1 >0.$
 Then the map $\left. f \to \sigma_{\nu} B^k f \right|_{y=0}$
    defined by relation \eqref{3.3.6} for  $f \in
\mathring{C}_{\nu}^{\infty} (E_{+}^{n+1})$ is extended up to a
linear bounded map of the space $H^s_{\nu} (E_{+}^{n+1})$ into the
space $H^{s-2 k + \Re \nu -1}_{\nu} (E^{n}).$ The inequality
$$
\|\left. \sigma_{\nu} B^k f \right|_{y=0} \|_{H^{s - 2k + \Re \nu
- 1} (E^{n})}  \leq c \, \| g \|_{H^{s}_{\nu} (E^{n+1}_{+})},
$$
    where $c$ is a positive constant independent on the function $f,$ holds.
\end{theorem}
  Consider the case of imaginary values of the parameter $\nu.$
  No relations of kind \eqref{3.3.1} hold in this case. Indeed, if $\Re\nu \geq 0$,
  then, due to \eqref{1.3.3}, any $g \in \mathring{C}^{\infty} (\ov{E_{+}^{n+1}})$
  satisfies the relation
\begin{equation}
P_{\nu}^{\frac{1}{2}-\nu} g (x', y) =  \frac{i
y^{-\nu}}{2^{\nu+2}\, \Gamma (\nu+1)}
\int\limits_{-\infty}^{\infty} H^{(1)}_{\nu} (y \eta) \eta^{\nu} F
\mbox{Cont} g (x', \eta) \, d \eta. \label{3.3.7}
\end{equation}
    If $\Re \nu =0,$ but $\nu \neq 0,$ then the Hankel function $ H^{(1)}_{\nu}(z)$
    has the following form in a neighborhood of the point $z=0$ (see
    \cite{BE2}):
$$
H^{(1)}_{\nu}  (z) = {c'}_{\nu} z^{- \nu} \lr{1+o(1)}+{c''}_{\nu} z^{- \nu} \lr{1+o(1)},
$$
    where $c'_{\nu}$ and $c''_{\nu}$ are known constants.
    The functions $z^{- \nu}$ and $z^{\nu}$ oscillate near the origin:
    $z^{- \nu} = \cos (\mu \ln z) - i \sin (\mu \ln z)$ and
$z^{\nu} = \cos (\mu \ln z) + i \sin (\mu \ln z).$
  Therefore, no power function $\varphi (z)$ satisfying the
  relation $\lim\limits_{z \to 0} H^{(1)}_{\nu}  (z) =1$ exist.
  Hence, in general, no weight traces of the function  $P_{\nu}^{\frac{1}{2}- \nu} g(x', y)$
  exist. That is why we consider traces of its derivative instead.

    Use the following recurrent relation for the function
$H^{(1)}_{\nu}$ (see \cite[p. 20]{BE2}):
$$
\frac{\pr}{z \pr z} \lr{z^{- \nu} H^{(1)}_{\nu}  (z)} = - z^{-\nu-1} H^{(1)}_{\nu+1}  (z).
$$
 Then, from \eqref{3.3.7}, we conclude that
$$
\frac{\pr}{y \pr y} P_{\nu}^{\frac{1}{2}-\nu} g (x', y)  =
\frac{-i y^{-\nu-1}}{2^{\nu+2} \Gamma (\nu+1)}
\int\limits_{-\infty}^{\infty} H^{(1)}_{\nu+1} (y \eta)
\eta^{\nu+1} F \mbox{Cont} g (x', \eta) \, d \eta.
$$
Hence,
$$
\frac{\pr}{y \pr y} P_{\nu}^{\frac{1}{2}-\nu} g (x', y)  = -2
(\nu+1) P_{\nu}^{-\frac{1}{2}-\nu} g (x', y).
$$
    This and \eqref{3.3.1} implies the relation
\begin{equation}
\lim\limits_{y \to +0}  y^{2 \nu +1} \frac{\pr}{ \pr y}
P_{\nu}^{\frac{1}{2}-\nu} g (x', y)  = - g(x', 0). \label{3.3.8}
\end{equation}

 Relation \eqref{3.3.8} valid for all $\nu$ from the half-plane
  $\Re \nu \geq 0$ (including imaginary values of $\nu$) is used to obtain
  the following theorem on weight traces.

\begin{theorem} \label{teo: 3.3.3}
 Let $\Re \nu \geq 0$ and $s> 2 k - \Re \nu +1 >0.$

 Then the map  $\left. f \to y^{2 \nu +1} D_y B^k f \right|_{y=0}$
 defined by the relation
    $$
    y^{2 \nu +1} \left. D_y B^k f \right|_{y=0} =
     - D^{2 k}_y I_e^{\nu - \frac{1}{2}} \left. S_{\nu, e} f \right|_{y=0}
    $$
  under the assumption that  $f \in \mathring{C}_{\nu}^{\infty} (E_{+}^{n+1})$
   is extended up to a continuous and bounded map of the space
      $H^s_{\nu} (E_{+}^{n+1})$ into the space $H^{s-2 k + \Re \nu -1} (E^{n}).$
 The inequality
    $$
    \|y^{2 \nu +1} \left. D_y B^k f \right|_{y=0} \|_{H^{s - 2k + \Re \nu - 1} (E^{n})}
    \leq \|f\|_{H_{\nu}^s (E^{n+1}_{+})},
    $$
  where $c$ is a positive constant independent of $f,$ holds.
\end{theorem}

\section{Function Spaces in $H^{s}_{\nu}
(\Omega)$}\label{sec11}

\subsection{Partition of unity and definitions of function spaces}\label{sec11.1}

Let $\Omega$ be a bounded domain of the half-space
 $E^{n+1}_{+}$ with boundary of the class $C^{\infty}.$
 Let a domain $\Omega_0 $ located strictly inside
  $\Omega$ (this means that $\ov{\Omega}_0 \subset \Omega$ or, which is the same,
$\ov{\Omega}_0 \cap \pr \Omega = \varnothing$) and domains
$\Omega_l,$ $l =1, 2, \dots, \ov{l},$ such that their
intersections with the boundary $\pr \Omega$ are nonempty form a
covering of $\ov{\Omega}.$
    Assign $\Omega_l^{+} = \Omega_l \cap\Omega,$ $l =1, 2, \dots, \ov{l}.$
  Let there exist diffeomorphisms $\varkappa_l$ of the class $C^{\infty},$
  mapping the domains $\Omega_l,$ $l =1, 2, \dots, \ov{l},$ into domains $\omega_l$
 located in the space $E^{n+1}.$ Let $\varkappa_0$ denote the identity transformation.
 Let the domain $\Omega_l^{+}$ be mapped into $\omega_l^{+} := \omega_l
\cap E^{n+1}_{+}$ and a part of the boundaries $\Omega_l \cap \pr
\Omega$ be mapped into a part  $\omega_l \cap \{ y=0 \}$ of the
hyperplane. Also, we assume that if  $\Omega_l \cap \Omega_{l'}
\neq \varnothing,$ then the natural map $\varkappa_{l'}
\varkappa_{l}^{-1}$ of the part $\varkappa_{l } \lr{\Omega_l^+
\cap \Omega_{l'}^+ }$ of the domain $\omega_l^+$ into the part
$\varkappa_{l '} \lr{\Omega_l^+ \cap \Omega_{l'}^+ }$ of the
domain $\omega_{l'}^+$ is a nondegenerate transformation with a
positive Jacobian determined by relations $x'=x' (\widetilde{x}')$
and $y = \widetilde{y}$ in a neighborhood of the boundary.

    In the sequel, the specified covering of the domain is
    assumed to be fixed; all other coverings used below are
    obtained by a refinement of the fixed one.

\begin{lemma}\label{lem: 4.1.1.}
 There exist functions $h_{l} \in \mathring{C}^{\infty} (E^{n+1})$ with the following
 properties{\rm :}
    \begin{enumerate}
    \item[(1)] $h_{l} (x)=0$ if $x \in E^{n+1} \setminus \Omega_{l},$ ${l}=0,
     \dots, \ov{{l}};$
    \item[(2)] $0 \leq h_{l} (x) \leq 1$ for all $x \in E^n;$
    \item[(3)] $h_0 (x)+ \ldots + h_{\ov{l}} (x)=1$ provided that $x \in \Omega;$
    \item[(4)] $D_y h_{l} (x) =0,$ ${l}=0, \dots, \ov{l},$
  in local coordinates in a neighborhood of the hyperplane $\{y=0\}.$
    \end{enumerate}

\end{lemma}

\begin{proof}
 Let $U_R (x)$ be an open ball of a positive radius $R,$ centered at a point $x \in E^{n+1}$.
 Let
$$
\omega_{l, \varepsilon} = \{x \in \omega_{l}: U_{\varepsilon}
\subset \omega_{l} \}.
$$
    It is easy to see that if the positive $\varepsilon$ is sufficiently small, then
    the domains $\varkappa_{l}^{-1} \omega_{ {l}, \varepsilon}$ form
    a covering of $\ov{\Omega}.$

Let a positive  $\varepsilon$ be such that the domains
$\varkappa_{l}^{-1} \omega_{ {l}, \varepsilon},$ $l = 0, \dots,
\ov{l},$ form a covering of the domain $\ov{\Omega}.$
    Consider a domain $\omega_{ {l}, \varepsilon}.$ Let $\mathcal{K}_{\delta} (x)$ denote
    an open $n+1$-dimensional
    cube centered at a point $x$ such that the length of its edge is equal to $2 \delta$
    and its faces are parallel the coordinate hyperplanes.
    Since each point $x \in\omega_{ {l}, \varepsilon} \cap \{y=0\}$ is contained
    in $\omega_l$ together with a ball of radius $\varepsilon,$
    it follows that there exists its finite covering by cubes of the kind
 $\mathcal{K}_{\varepsilon_p/2} (x^p),$ $p=1, \dots, \ov{p} < \infty,$ $\varepsilon_p > 0,$
 with centers at   points $x^p \in \omega_{ {l}, \varepsilon} \cap
\{y=0\}$ respectively such that $\mathcal{K}_{\varepsilon_p} (x^p)
\subset \omega_l.$ Cover the remaining part of the domain
    $\omega_{ l, \varepsilon}$ by cubes of the kind
$\mathcal{K}_{\varepsilon_p/2} (x^p),$ $p= \ov{p}+1, \dots,
\ov{\ov{p}} < \infty,$ such that $\mathcal{K}_{\varepsilon_p}
(x^p) \subset \omega_l$ and $\ov{\mathcal{K}}_{\varepsilon_p}
(x^p)$ has no common points with the hyperplane  $\{y=0\}.$

Let $\varphi_{\delta} (t)$ be a one-variable
 functions such that $\varphi_{\delta} \in \mathring{C}^{\infty}
(E^1),$ $\varphi_{\delta} (t) \geq 0$ for $t \in E^1,$
$\varphi_{\delta} (t)=1$ for $|t| \leq \delta,$ and
$\varphi_{\delta} (t)=0$ for $|t| \geq 2 \delta.$

 Now, introduce the function $\psi_{l}(x),$ $x= (x', y),$ by the
 relation
$$
\psi_{l}(x) = \sum\limits_{p=1}^{\ov{p}} \varphi_{\varepsilon_p}
(y) \prod\limits_{q=1}^{n} \varphi_{\varepsilon_p} (x_q-x_q^p) +
\sum\limits_{p=\ov{p}+1}^{\ov{\ov{p}}} \psi_{{l}, p} (x),
$$
    where $x_q^p$ denotes the $q$th coordinate of the point $x^p,$
    while $\psi_{l, p}(x)$ are functions possessing the following
    properties:
$\psi_{{l}, p} (x) \in \mathring{C}^{\infty} (E^{n+1}),$ $\psi_{l,
p}(x) \geq 0$ for $x \in E^{n+1},$ $\psi_{l, p}(x)=1$ for $x \in
\mathcal{K}_{\varepsilon_p/2} (x^p),$ and $\psi_{l, p}(x)=0$ for
$x \in E^{n+1} \setminus \mathcal{K}_{\varepsilon_p} (x^p),$ $p=
\ov{p}+1, \dots, \ov{\bar{p}}.$

The proof of the existence of functions $\varphi_{\delta}$ and
$\psi_{l, p}$ satisfying all listed conditions can be found,
e.\,g., in \cite{65}. The constructed function  $\psi_l (x)$
    possesses the following properties: $\psi_{l} (x) \in \mathring{C}^{\infty}
(E^{n+1}),$ $\psi_{l} (x) \geq 0$ for $x \in E^{n+1},$ $\psi_{l}
(x) \geq 1$ for $x \in \omega_{l, \varepsilon},$ and $\supp \psi_l
\subset \omega_l.$
 It is easy to see that $D_y \psi_l (x)
\equiv 0$ for $|y| < \delta,$ where $\delta$ is a positive number.

The desired functions $h_l$ define by the relation
$$
h_{l} (x) = \frac{\varkappa_{l}^{-1}
\psi_{l}}{\sum\limits_{m=0}^{\ov{l}} \varkappa_{l}^{-1} \psi_m}.
$$
It is easy to verify that Conditions (1)--(4)
 are satisfied, which completes the proof of the lemma.
\end{proof}
 In the sequel, we say that a function collection $\{h_l\}$ forms
 a partition   {\it subordinated to the covering} $\{\Omega_l\}$ if  Conditions (1)--(4)
 are satisfied.

Now, introduce the space $H^s_{\nu} (\Omega),$ $s \geq 0,$ $\Re
\nu \geq 0$ as the set of functions $f$ defined in the domain
$\Omega$ and such that the functions $\varkappa_l (h_l f)$ belong
to the space $H^s_{\nu} (E^{n+1}_{+}).$ This is a Hilbert space
with the norm
\begin{equation}
\| f \|_{H^s_{\nu} (\Omega)} =  \lr{ \sum\limits_{l=0}^{\ov{l}} \|
\varkappa_l (h_l f) \|^2_{H^s_{\nu} (E^{n+1}_{+})}}^{\frac{1}{2}},
\label{4.1.1}
\end{equation}
    where $\{h_l\}$ is a partition of unity, subordinated to the
    covering $\{\Omega_l\}.$

To show that norms \eqref{4.1.1} are equivalent for different
choices of the partition of unity, obtain more general result.
    Let $\{{\Omega'}_{l'}\}$ be another covering of $\ov{\Omega}.$
  Assume that for each $l'=0, \dots, \ov{l'}$ there exists   $l=0, \dots, \ov{l}$
  such that ${\Omega'}_{l'}\subset {\Omega}_{l}.$
  If this condition is satisfied, then the covering
$\{{\Omega'}_{l'}\}$ is called a \emph{refinement} of the covering
$\{{\Omega}_{l}\}.$ Let $\{h'_{l'}\}$ be a partition of unity,
subordinated to  the covering $\{{\Omega'}_{l'}\}.$ Then  the
following norm can be introduced:
\begin{equation}
\| f \|_{H^s_{\nu} (\Omega)} =  \lr{ \sum\limits_{l'=0}^{\ov{l'}}
\| \varkappa'_{l'} (h'_{l'} f) \|^2_{H^s_{\nu}
(E^{n+1}_{+})}}^{\frac{1}{2}}, \label{4.1.2}
\end{equation}
 where the map is such that $\varkappa'_{l'}=\varkappa_l,$ $l=l(l'),$
 and $l$ is selected such that $\supp h'_{l'} \subset
\Omega_l.$
 If there are more than one such $l,$ then any of them is
 satisfied
(it follows from properties of  $\varkappa_l$ that the
corresponding terms are equivalent norms).

    Estimate a term from the right-hand
    side of \eqref{4.1.2}. We have the inequalities
$$
    \begin{gathered}
\| \varkappa'_{l'} (h'_{l'} f) \|_{H^s_{\nu} (E^{n+1}_{+})} =  \|
\sum\limits_{l=0}^{\ov{l}} \varkappa'_{l'} (h'_{l'} h_{l} f)
\|_{H^s_{\nu} (E^{n+1}_{+})} \leq
  \sum\limits_{l=0}^{\ov{l}}
   \| \varkappa'_{l'} (h'_{l'} h_{l} f) \|_{H^s_{\nu} (E^{n+1}_{+})}
   \\
\leq c_1 \sum\limits_{l=0}^{\ov{l}} \| \varkappa_{l} (h'_{l'}
h_{l} f)  \|_{H^s_{\nu} (E^{n+1}_{+})} \leq c_2
\sum\limits_{l=0}^{\ov{l}} \| \varkappa_{l} (h_{l} f)
\|_{H^s_{\nu} (E^{n+1}_{+})}.
    \end{gathered}
$$
    In the next to the last inequality, properties of the map
  $\varkappa'_{l'} \varkappa^{-1}_{l}$ are used (see the beginning of the section);
  in the last inequality, results from Sec. \ref{sec3.2} about multipliers are used.
  Property  4 of functions ${h'}_{l'}$ is used here as well.
  The inverse estimate is proved in the same way. Thus, the equivalence of norms \eqref{4.1.1}
   and \eqref{4.1.2} is established.

Consider another structure property of the space $H^s_{\nu}
(\Omega).$ Let $f \in H^s_{\nu} (\Omega)$ and $a \in
\mathring{C}^{\infty} (\Omega).$ Then results of Sec. \ref{sec3.1}
imply that the function $a f$ belongs to the space $H^s (\Omega)$
    and the estimate
\begin{equation}
c_1  \| a f \|_{H^s (\Omega)} \leq   \| a f \|_{H^s_{\nu}
(\Omega)} \leq c_2   \| a f \|_{H^s (\Omega)}, \label{4.1.3}
\end{equation}
    where $c_1$ and $ c_2 $ are positive constant independent of the function $f,$ holds.

\begin{lemma} \label{lem: 4.1.2}
  Let $\widetilde{\Omega}$ be a strictly internal subdomain of the domain $\Omega.$
  Then the restriction of each $f \in H^s_{\nu} (\Omega)$ to
    $\widetilde{\Omega}$ belongs to the space $H^s (\widetilde{\Omega})$    and the estimate
    \begin{equation}
        \| f \|_{H^s (\widetilde{\Omega})} \leq c \,  \|  f \|_{H^s_{\nu} (\Omega)},
        \label{4.1.4}
    \end{equation}
 where $c$ is a positive constant independent of the function $f,$ holds.
\end{lemma}

\begin{proof}
Introduce $a \in \mathring{C}^{\infty} (\Omega)$ such that $a (x)
= 1$ provided that $x \in \widetilde{\Omega}.$ From inequality
\eqref{4.1.3}, we conclude that
$$
\| f \|_{H^s (\widetilde{\Omega})} \leq   \| a f \|_{H^s (\Omega)}
\leq  c \, \| a f \|_{H^s_{\nu} (\Omega)}.
$$
 From Theorem \ref{teo: 3.2.1}, we obtain the estimate
$$
\| a f \|_{H^s_{\nu} (\Omega)} \leq  c \, \|  f \|_{H^s_{\nu}
(\Omega)},
$$
where $c$ is a positive constant independent of the function $f,$
    which completes the proof of the lemma.
\end{proof}
 The Kipriyanov space $H^s_{\nu, +} (\Omega),$ $s \geq 0,$ $\nu \geq 0,$
 is defined as the closure with respect to the
norm
\begin{equation}
\| f \|_{H^s_{\nu, +} (\Omega)}  =\lr{ \sum\limits_{l=0}^{\ov{l}}
\| \varkappa_{l} (h_{l} f) \|^2_{H^s_{\nu, +}
(E^{n+1}_{+})}}^{\frac{1}{2}}, \label{4.1.5}
\end{equation}
 where the norms $H^s_{\nu, +} (E_{+}^{n+1})$ are defined in Sec.
 \ref{sec3.1},
 of the set of all
functions $f$ from $C^{\infty} (\ov{\Omega})$ such that  $D_y^{2 k
+1} f = 0,$ $y=0,$ $k=0, 1, \dots,$ in each local coordinate
system (see \cite{Kip1}).

\begin{lemma} \label{lem: 4.1.3}
 Let $s \geq 0,$  $\nu \geq 0,$ and $\nu \neq 1, 3, 5, \dots$
 Then the space $H^s_{\nu, +} (\Omega)$ is a proper space of the space $H^s_{\nu} (\Omega)$
 and norm  \eqref{4.1.1} is equivalent to norm \eqref{4.1.5} on $H^s_{\nu, +} (\Omega)$.
 For $\nu = 1, 3,  \dots,$ the space
       $H^s_{\nu, +} (\Omega)$ is embedded in $H^s_{\nu} (\Omega).$
\end{lemma}

\begin{proof}
Taking into account that the functions  $h_l f$ (see the
definition of the norm) have compact supports, one can apply Lemma
\ref{lem: 3.1.4}, which completes the proof.
\end{proof}

\subsection{Embedding theorems}\label{sec11.2}

We start from internal embedding theorems. The next assertion
easily follows from Theorem \ref{teo: 3.1.1}.

\begin{theorem} \label{teo: 4.2.1}
 Let $s'>s \geq 0$ and $\Re \nu \geq 0.$ Then the space $H^{s'}_{\nu} (\Omega)$
 is  embedded in $H^{s}_{\nu} (\Omega)$.
\end{theorem}

\begin{theorem} \label{teo: 4.2.2}
    If $s'>s \geq 0$ and $\Re \nu \geq 0,$ then the embedding
    operator from $H^{s'}_{\nu} (\Omega)$ to $H^{s}_{\nu} (\Omega)$
    is completely continuous.
\end{theorem}

\begin{proof}
Let a sequence of functions $f_j,$ $j=1,2, \dots,$ be bounded with
respect to the norm of the space $H^{s}_{\nu} (\Omega).$
    Let $\mathsf{y}_j,$ $j=1,2, \dots,$ denote an expanding system
    of strictly internal subdomains of the domain $\Omega$ such that their union
     coincides with  $\Omega.$
  Due to Lemma \ref{lem: 4.1.2}, the sequence of $f_j$ is bounded with respect to the norm
 of each space $H^{s'}(\mathsf{y}_j);$ therefore, it contains a subsequence $\{f_{1j}\}$
 converging with respect to the norm of the subspace $H^{s} (\mathsf{y}_1)$.
 Denote the limit of this subsequence by $f^1.$
 From $f_{1j},$ select a subsequence $f_{2j}$ converging to $f^2$
 with respect to the norm $H^{s}(\mathsf{y}_2).$ It is clear that $f^2|_\mathsf{y} = f^1.$
 Repeating this procedure, we find a function $f$ such that
 $f|_{\mathsf{y}_j} \in H^s (\mathsf{y}_j)$ for each $j$.
 Form a diagonal sequence $\{f_{jj}\}.$ It satisfies the relation
$$
\| \left. (f-f_{jj}) \right|_{\mathsf{y}_p} \|_{H^s
(\mathsf{y}_p)} \to 0, ~ j \to \infty, ~  p  = 1,2, \dots
$$
    Without loss of generality, we assume that $\{fj\}$ itself possesses this property.

Further, the sequence $\{\varkappa_l (h_l f_j) \}$ is a subset of
$H^{s'}_{\nu} (\omega_l^{+})$ and  is bounded with respect to the
norm
 of this space. Then, due to Theorem \ref{teo: 3.1.2}, there exist
  $g_l \in H^{s}_{\nu} (\omega_l^{+})$ such that a subsequence of
  the sequence $\{\varkappa_l (h_l f_j)\}$ converges to it.
 Again, without loss of generality, we assume that the sequence $\{fj\}$ itself
 possesses this property for all $l=0,1,\dots, \ov{l}$.
  Let $\varepsilon>0.$ Let $\omega^{+}_{l, \varepsilon}$ denote the part of the domain
   $\omega_l^{+},$ located in the subspace $\{y> \varepsilon\}.$
 Then, by Lemma \ref{lem: 4.1.2}, the sequence
$\{\varkappa_l (h_l f_j)|_{\{y > \varepsilon\}}\}$ converges to
the function $\varkappa_l (h_l f)|_{\{y > \varepsilon\}}$  with
respect to the norm
 of the space $H^s_{\nu} (\omega^{+}_{l, \varepsilon}).$
    Since  $\varepsilon$ is selected arbitrarily, it follows that $g_l= \varkappa_l
(h_l f)$ for all $l.$ Then $f \in H^s_{\nu} (\Omega),$ which
completes the proof of the theorem.
\end{proof}
  To consider the assertion about weight traces, we need boundary spaces $H^s (\pr \Omega)$
  (see \cite{66}). An equivalent norm of the space $H^s (\pr \Omega)$ can be defined
  by the relation
$$
\| \varphi \|^2_{H^s(\pr \Omega)} =  \sum\limits_{l=0}^{\ov{l}} \|
\varkappa_{l} (h_{l} f) \|_{H^s (E^{n})}.
$$
Let $\sigma_{\nu}$ denote the function coinciding with the weight
 function $\sigma_{\nu} (y)$ (see above) in each local coordinate
system.

We say that a function  $f\in H^s_{\nu} (\Omega)$ has a weight
trace $\left. \sigma_{\nu} f \right|_{\pr \Omega} = \varphi$ on
the boundary $\pr \Omega$ if the functions $h_l f,$ $l=0,1,\dots,
\ov{l},$ have a weight $\sigma_{\nu}$-trace
 in each local coordinate system.

Let $B_{\nu}$ be a differential operator coinciding with the
Bessel operator in each local coordinate system.

\begin{theorem}\label{teo: 4.2.3}
  Let $\Re \nu > 0$ or $\nu=0.$ Let $s>2k+1- \Re \nu > 0.$
  Then each $f \in H^s_{\nu} (\Omega)$ has a weight trace $B_{\nu}^k f$ and the inequality
$$
\| \left. \sigma_{\nu} B^k_{\nu} f \right|_{\pr \Omega}
\|_{H^{s-2k-1+ \Re \nu} (\pr \Omega)}  \leq c \, \| f
\|_{H^s_{\nu} (\Omega)},
$$
 where $c$ is a positive constant independent of the function $f,$ holds.

 If $\Re \nu \geq 0$ and $s>2k+1- \Re \nu > 0,$ then each $f \in H^s_{\nu} (\Omega)$
 has weight traces
$\left. \sigma_{\nu + \frac{1}{2}} D B^k_{\nu} f \right|_{\pr
\Omega}$ of functions $D B^k f$ and the inequalities
$$
\| \left. \sigma_{\nu +  \frac{1}{2}} D B^k_{\nu} f \right|_{\pr
\Omega} \|_{H^{s-2k-1+ \Re \nu} (\pr \Omega)} \leq c \, \| f
\|_{H^s_{\nu} (\Omega)}
$$
    hold, where  $c$ is a positive constant independent of the function $f$ and $D$
denotes an operator coinciding with $\dfrac{\pr }{\pr y}$ in each
local coordinate system near the boundary.
\end{theorem}

This assertion follows from Theorems \ref{teo: 3.3.2}-\ref{teo:
3.3.3}.

\section[Half-Space Elliptic Boundary-Value
 Problems\\
   with Nonlocal
 Fractional-Order Boundary-Value
 Conditions]{Half-Space Elliptic Boundary-Value
 Problems\\
  with Nonlocal
 Fractional-Order Boundary-Value
 Conditions}\label{sec12}

\subsection{Boundary-value
 problem: setting}\label{sec12.1}

    Consider the following homogeneous partial differential
    operator with constant coefficients:
$$
A \lr{\frac{1}{i} D} = \sum\limits_{|\alpha|=2 m} a_{\alpha}
\lr{\frac{1}{i} D}^{\alpha},
$$
    where $\alpha = (\alpha_1, \dots, \alpha_{n+1}) = (\alpha',
\alpha_{n+1}),$ $\alpha_j$ are nonnegative integers, $|\alpha| =
\alpha_1+ \dots+ \alpha_{n+1},$ and
$$
\lr{\frac{1}{i} D}^{\alpha} = \frac{\pr^{|\alpha|}}{(i \pr
x_1)^{\alpha_1} \dots (i \pr x_{n+1})^{\alpha_{n+1}}}.
$$
 We say that the operator $A$ is \emph{elliptic} if
$$
A(\xi) = \sum\limits_{|\alpha|=2 m} a_{\alpha}  \xi^{\alpha} \neq
0
$$
    for all $\xi \in E^{n+1},$ $\xi \neq 0,$ where $\xi^{\alpha} =
\xi^{\alpha_1}_1 \dots \xi^{\alpha_{n+1}}_{n+1}.$

Assume that the following \emph{proper ellipticity} condition is
satisfied: for each $\xi' \in E^n,$ $\xi \neq 0,$ the
characteristic polynomial $A (\xi) = A (\xi', \eta)$ of the
complex variable $\eta$ has $m$ roots (counted with their
multiplicities) $\eta_j^{+},$ $j=1, \dots, m,$ with positive
imaginary parts.

Let $G_j \Big(\dfrac{1}{i} D\Big),$ $j=0, \dots, m-1,$ be boundary
operators with constant coefficients such that  $G_j (\xi)$ is a
homogeneous polynomial of degree $m_j \geq 0.$

Assume that the following Shapiro--Lopatinskii
 condition is satisfied: for each $0 \neq \xi' \in
E^n,$ the polynomials $G_j (\xi', \eta)$ of the variable $\eta$
are linearly independent modulo polynomial $A^{+} (\xi', \eta) =
(\eta - \eta_1^{+} (\xi ')) \dots (\eta - \eta_m^{+} (\xi ')).$

If the above conditions are satisfied, then we say that the
operators  $\{A, G_j\} $ as well as their corresponding
polynomials  (symbols) form an \emph{elliptic collection}.

Consider the boundary-value
    problem
\begin{equation}
\begin{cases}
A \lr{\dfrac{1}{i} D} u(x) = f (x), &  x \in E^{n+1}_{+}, \\
G_j \lr{\dfrac{1}{i} D} \left.  D^k_y I^{\mu}_e u \right|_{y=0} =
g_j (x'), &  x' \in E^n, \  j=0,\dots, m-1,
\end{cases}
\label{5.1.1}
\end{equation}
    where $k$ is a nonnegative integer and  $I^{\mu}_e$ is the  Liouville-type
 operator introduced in Sec. \ref{sec4.2} and acting  with respect to the variable $y.$
 Its order $\mu$ might be a complex number.

Once the Fourier transformation with respect to  $n$ initial
variables, defined by the relation
$$
\widetilde{f} (\xi', y) = F' f (\xi', y) =  \int\limits_{E^n} e^{-
i \langle x', \xi' \rangle} f (x', y) \, d x'
$$
    is applied, problem \eqref{5.1.1} takes the form
\begin{equation}
\begin{cases}
A \lr{\xi', \dfrac{1}{i} D_y}  \widetilde{u} (\xi', y) = \widetilde{f} (\xi', y), &  \xi' \in E^{n}, \  y>0, \\
G_j \lr{\xi', \dfrac{1}{i} D_y}  \left.  D^k I^{\mu}_e
\widetilde{u} \right|_{y=0} = g_j (\xi'), &  \xi' \in E^n, \
j=0,\dots, m-1.
\end{cases}
\label{5.1.2}
\end{equation}

\begin{lemma} \label{lem:5.1.1}
 Let operators $A$ and $G_j,$ $j=0,\dots, m-1,$ form an elliptic collection.
  Let the following relations be satisfied{\rm:}
\begin{equation}
k - \Re \mu + \frac{1}{2} > 0, \ 2  m +s -k + \Re \mu -
\frac{1}{2} - \max\limits_j m_j >0. \label{5.1.3}
\end{equation}
Let the function $\widetilde{f} (\xi', y)$ belong to the space
$H^s (E^1_{+})$ with respect to the last variable.
    Then, for each nonzero $\xi' \in E^n,$  problem \eqref{5.1.2}
    has a unique solution in the space $H^{2m+s} (E_{+}^1)$ and
    the estimate
\begin{multline}
c \, \sum\limits_{l=0}^{2 m + s} |\xi'|^{2(2m+s-l)}
 \| D^{l}_y \widetilde{u} (\xi', y) \|_{L_2 (E^1_+)}  \\
\leq  \sum\limits_{l=0}^{s} |\xi'|^{2(s-l)} \| D^{l}_y
\widetilde{f} (\xi', y) \|_{L_2 (E^1_+)} + \sum\limits_{j=0}^{ m -
1} |\xi'|^{2(2m+s-m_j-k + \Re \mu - \frac{1}{2})}
|\widetilde{g}_j (\xi')|^2, \label{5.1.4}
\end{multline}
 where $c$ is a positive constant independent of $f$ and $g_j,$ holds.
 If $|\xi'| \geq 1,$ then this constant does not depend on $\xi'.$
\end{lemma}

    This lemma can be treated as known though the authors are not
    able to provide a particular reference.
    Its proof is standard; it can be easily restored according to
    the following scheme. Once a particular solution of the
    equation is found by means of the continuation and applying
    the Fourier transformation with respect to the last variable,
    it suffices to consider the boundary-value
    problem with the homogeneous equation.
    Applying the operator $D_y^k I_e^{\mu}$ for the function $D_y^k
I_e^{\mu} \widetilde{u}$ to the equation, we obtain a usual
    boundary-value problem such that its solution is represented by a contour integral.
  Then the inverse operator is applied (the operator $D_y^k I_e^{\mu}$ is invertible
  on the space of rapidly decreasing functions) and a direct estimate of the obtained
  expression is provided.

\subsection{Regularizer and a priori estimates}\label{sec12.2}

    Problem \eqref{5.1.1} generates the operator
\begin{equation*}
\mathfrak{U}: u \to \mathfrak{U} u =  \left\{ A  u, G_0  D_y^k
I_e^{\mu} u |_{y=0}, \dots, G_{m-1}  D_y^k I_e^{\mu} u |_{y=0}
\right\}.
\end{equation*}
    Introduce
\begin{equation*}
\mathcal{H}^s  \lr{E_{+}^{n+1}, E^n, m} = H^s (E_{+}^{n+1})
\times \prod\limits_{j=0}^{m-1} H^{s+2m-m_j-k+\Re \mu -
\frac{1}{2}} \lr{E^n}.
\end{equation*}
 In $\mathcal{H}^s$, introduce the direct-product
  topology.

It follows from Theorem \ref{3.3.1} that  the operator
$\mathfrak{U}$ continuously maps the space $H^s (E_{+}^{n+1})$
into $\mathcal{H}^s \lr{E_{+}^{n+1}, E^n, m}.$

Let $\Phi = \{f(x), g_0 (x'), \dots,  g_{m-1} (x')\}$ be an
element of the space $\mathcal{H}^s.$
    An operator $\mathfrak{R}_{\mbox{l}}:
\mathcal{H}^s \to  H^{s+m}$ such that
$$
\mathfrak{R}_{\mbox{l}} \mathfrak{U} u =  u + T_{\mbox{l}} u, \  u
\in H^{s+2 m } (E^{n+1}_+),
$$
 where ${T}_{\mbox{l}}$ is a smoothing operator,
${T}_{\mbox{l}}: H^{s+2m} \to  H^{s+2m+1},$ is called a
\emph{left-hand regularizer} for $\mathfrak{U}$.

    An operator  $\mathfrak{R}_{\mbox{r}}: \mathcal{H}^s \to  H^{s+2m}$
    is called a \emph{right-hand
     regularizer} if
$$
  \mathfrak{U}  \mathfrak{R}_{\mbox{r}} \Phi = \Phi + T_{\mbox{r}} \Phi,
$$
    where ${T}_{\mbox{r}}: \mathcal{H}^s \to  \mathcal{H}^{s+1}$
    is a smoothing operator.

 An operator $\mathfrak{R}$ is called a two-side
 regularizer (or regularizer) if it is a
 left-hand regularizer and a right-hand
 one.

The main goal of the present section is to construct a regularizer
for the operator $\mathfrak{U}$ introduced above.
    We change the H\"ormander scheme to construct a  regularizer (see \cite{83}).
 We construct the desired regularizer as follows. Assume that, for $|\xi'| \geq 1$,
  a function $U (\xi', y)$ is a unique solution of the
  boundary-value problem
\begin{equation}
    \begin{cases}
        A \lr{\xi', \dfrac{1}{i} D_y}  {U} (\xi', y)
         = \widetilde{f} (\xi', y),  &  y>0, \\
        G_j \lr{\xi', \dfrac{1}{i} D_y}  \left.  D^k I^{\mu}_e {U} \right|_{y=0}
         = \widetilde{g}_j (\xi'), &  j=0,\dots, m-1, \\
        \widetilde{f} (\xi', y) = F' f, \quad  \widetilde{g}_j (\xi') = F' g_j,
    \end{cases}
    \label{5.2.1}
\end{equation}
    where $F'$ denotes the Fourier operator with respect the initial $n$
variables.
    For $|\xi'|<1 $, the function $U (\xi', y)$ satisfies the
    problem
\begin{equation}
\begin{cases}
A \lr{\xi'_0, \dfrac{1}{i} D_y}  {U} (\xi', y) = \widetilde{f} (\xi', y),  &  y>0, \\
G_j \lr{\xi'_0, \dfrac{1}{i} D_y}  \left.  D^k I^{\mu}_e {U}
\right|_{y=0} = \widetilde{g}_j (\xi'), &  j=0,\dots, m-1,
\end{cases}
\label{5.2.2}
\end{equation}
    where $\xi'_0 \in E^n$ is an arbitrary fixed point of the unit sphere, $|\xi'_0|=1.$
  Note that the difference from the  H\"ormander scheme is the way
   to define the function $U (\xi', y)$ for $|\xi'| < 1.$

Now, define the operator $\mathfrak{R}$ by the relation
\begin{equation}
\mathfrak{R} \Phi = \lr{F'}^{-1} {U},
\label{5.2.3}
\end{equation}
 where $\Phi = \{f, g_0, \dots, g_{m-1}\}.$

Let us show that the introduced operator $\mathfrak{R}$
continuously maps the space  $\mathcal{H}^s (E_{+}^{n+1}, E^n, m)$
into $H^{s+2m} (E_{+}^{n+1}).$ We have the inequalities
$$
    \begin{gathered}
\| \mathfrak{R} \Phi  \|_{H^{s+2m} (E^{n+1}_{+})} =  \| (F')^{-1}
U  \|^2_{H^{s+2m} (E^{n+1}_{+})}
    \\
\leq c \lr{ \| (F')^{-1}  U  \|_{L_2 (E^{n+1}_{+})} +
\sum\limits_{l=0}^{2m+s} \sum\limits_{|\alpha'|=2m+s-l}
\|D_{x'}^{\alpha'} D_y^{l} (F')^{-1}  U  \|_{L_2 (E^{n+1}_{+})}^2}
    \\
\leq c \lr{ \|  U (\xi', y)  \|^2_{L_2 (E^{n+1}_{+})} +
\sum\limits_{l=0}^{2m+s} \sum\limits_{|\alpha'|=2m+s-l}
\|{\xi'}^{\alpha'} D_y^{l}   U (\xi', y)  \|_{L_2
(E^{n+1}_{+})}^2}.
    \end{gathered}
$$
    taking into account that definitions of the function $U (\xi', y)$ are different
    for the cases where $|\xi'|
\geq 1$ and  $|\xi'| < 1,$ we start from the former case.
 We have the inequality
$$
    \begin{gathered}
I_1 = \int\limits_{|\xi'|>1} \int\limits_0^{\infty} \lr{ |U (\xi',
y)|^2+ \sum\limits_{l=0}^{2m+s} \sum\limits_{|\alpha'|=2m+s-l}
|\xi_1^{\alpha_1} \dots \xi_n^{\alpha_n}  D_y^{l}   U (\xi',
y)|^2} \, d \xi' dy
    \\
\leq c \sum\limits_{l=0}^{2m+s} \int\limits_{|\xi'|>1}
\int\limits_0^{\infty} (1+|\xi'|^2)^{2m+s-l}  | D_y^{l}   U (\xi',
y)|^2 \, dy d \xi'.
    \end{gathered}
$$
    Applying Lemma \ref{lem:5.1.1}, we obtain that
$$
    \begin{gathered}
I_1 \leq c \int\limits_{|\xi'|>1} \left( \sum\limits_{l=0}^{s}
(1+|\xi'|^2)^{s-l}  \int\limits_0^{\infty} | D_y^{l}   F' f (\xi',
y)|^2 \, dy  \right.
    \\
\left. + \sum\limits_{j=0}^{m-1} (1+|\xi'|^2)^{2m+s-m_j-k + \Re
\mu - \frac{1}{2}} |F' g_j (\xi')|^2  \right) \, d \xi' \leq c \,
\|\Phi \|^2_{\mathcal{H}^s},
    \end{gathered}
$$
    where $c$ is a positive constant independent of $\Phi.$

Consider the case where $|\xi'| < 1.$ From the definition of the
function  (see \eqref{5.2.2}), we obtain that
$$
I_2  = \int\limits_{|\xi'|<1} \int\limits_0^{\infty} \lr{ |U
(\xi', y)|^2+ \sum\limits_{l=0}^{2m+s}
\sum\limits_{|\alpha'|=2m+s-l}\hspace{-3ex}  |{\xi'}^{\alpha'}
D_y^{l} U (\xi', y)|^2} \, dy d \xi'  \leq
 c \sum\limits_{l=0}^{2m+s} \int\limits_{|\xi'|<1} \int\limits_0^{\infty} | D_y^{l}
    U (\xi', y)|^2 \, dy d \xi'.
$$
Applying Lemma \ref{lem:5.1.1} again, we obtain that
$$
I_2 \leq c \sum\limits_{l=0}^{s} \int\limits_{|\xi'|<1} \lr{
\int\limits_0^{\infty}  | D_y^{l}   F' (\xi', y) |^2 \, dy +
\sum\limits_{l=0}^{m-1} |F' g_j (\xi') |^2 } d \xi'.
$$
    This implies the estimate $I_2 \leq c \, \| \Phi
\|^2_{\mathcal{H}^s}.$
 Combining the estimates for $I_1$ and $I_2,$ we see that
$$
\| \mathfrak{R} \Phi \|_{H^{s+2m} (E_{+}^{n+1})}  \leq c \, \|
\Phi \|_{\mathcal{H}^s (E_{+}^{n+1}, E^n, m)}.
$$
Thus, the continuity of the operator  $\mathfrak{R}$ is proved.

Let us verify that the operator $\mathfrak{R}$  is a
 left-hand regularizer. We have the relation
\begin{equation}
\mathfrak{R} \mathfrak{U} u =u +  \frac{1}{(2 \pi)^n}
\int\limits_{|\xi'|<1} e^{i \langle x', \xi' \rangle} V_{(\xi',
y)} \, d \xi' = u + T_{\mbox{l}} u, \label{5.2.4}
\end{equation}
where the function $V (\xi', y)$ satisfies the boundary-value
 problem
\begin{equation}
    \begin{cases}
        A \lr{\xi'_0, \dfrac{1}{i} D_y}  V (\xi', y) {=} \left[ A \lr{\xi'_0, \dfrac{1}{i} D_y} -A \lr{\xi', \dfrac{1}{i} D_y}   \right] F'u(\xi', y), &  y>0, \\
        G_j \lr{\xi'_0, \dfrac{1}{i} D_y}
         \left.  D^k_y I^{\mu}_e  V \right|_{y=0} {=}\\
         \left[ G_j \lr{\xi'_0, \dfrac{1}{i} D_y} {-} G_j \lr{\xi', \dfrac{1}{i} D_y}
           \right]  D^k_y I^{\mu}_e   \left. F'u(\xi', y) \right|_{y=0},
           & j=\overline{0,\ m-1}.
    \end{cases}
    \label{5.2.5}
\end{equation}
  $T_{\mbox{l}} u = \lr{F'}^{-1} V_1,$   $V_1 (\xi', y) =
V(\xi', y) $ provided that $|\xi'|<1,$ and $V_1 (\xi', y)=0$
provided that $|\xi'| \geq 1.$ Then
$$
\| T_{\mbox{l}} u \|_{H^{s+2m} (E^{n+1}_{+})} \leq c
\sum\limits_{l=0}^{s+2m} \int\limits_{|\xi'|<1}
\int\limits_0^{\infty} | D^{l}_y V(\xi', y)|^2 \, dy d\xi'.
$$
Applying the estimate of Lemma \ref{lem:5.1.1} to the function $V$
and taking into account that the operators in the square brackets
in \eqref{5.2.5} has the same principal part, we obtain the
estimate
$$
    \begin{gathered}
\| T_{\mbox{l}} u \|_{H^{s+2m} (E^{n+1}_{+})} \leq c \left(
\sum\limits_{l=0}^{s} \int\limits_{|\xi'|<1}  \| D^{l+2m-1}_y  F'u
(\xi', y) \|^2_{L_2 (E^1_{+})} \,  d\xi' \right.
    \\
+ \left. \sum\limits_{l=0}^{m_j-1}  \sum\limits_{j=0}^{s-1}
\int\limits_{|\xi'|<1} |\left.  D^{l+k}_y I^{\mu}_e  F' u (\xi',
y) \right|_{y=0}|^2  \, d \xi' \right).
    \end{gathered}
$$
    It is obvious that the first sum is estimated by the norm
 $\| u \|_{H^{s+2m+1}(E^{n+1}_{+})}.$
 From Theorem \ref{teo: 3.3.1}, it follows that the second sum can be estimated by the same
 norm as well. Thus, the  estimate
$$
\| T_{\mbox{l}} u \|_{H^{s+2m} (E^{n+1}_{+})} \leq c \, \| u
\|_{H^{s+2m+1} (E^{n+1}_{+})},
$$
    where the constant does not depend on the function $u,$ is proved.
Hence, $T_{\mbox{l}}$ is a smoothing operator.

    In the same way, it is proved that the operator $\mathfrak{R}$ is a right-hand
 regularizer.

 Thus, the following assertion is proved.

\begin{theorem} \label{teo: 5.2.1}
 Let operators $A$ and $G_j,$ $j=0, \dots, m-1,$ form an elliptic collection.
  Let the relations
 $$
 k - \Re \mu + \frac{1}{2}>0  \quad \textrm{and} \quad
 s+2m-k+ \Re \mu -  \frac{1}{2} - \max\limits_{j} m_j > 0
 $$
 be satisfied.
 Then the operator $\mathfrak{U}$ has a regularizer $\mathfrak{R},$
 belonging to $L \lr{\mathcal{H}^s(E^{n+1}_{+}, E^n, m), H^{s+2m} (E^{n+1}_{+})}.$
 If $u\in H^{2m+s} (E^{n+1}_{+})$ and there exists a positive $p$ such that
  $\mathfrak{U} u \in H^{s+p},$ then $u \in H^{2m+s+p} (E^{n+1}_{+})$
  and the a priori estimate
$$
c \, \| u \|_{H^{s+2m} (E^{n+1}_{+})} \leq \| A u \|_{H^{s}
(E^{n+1}_{+})}  + \sum\limits_{j=0}^{m-1} \| G_j D_y^k I^{\mu}_e
u|_{y=0}  \|_{H^{s_j} (E^n)} +
 \| u \|_{H^{s+2m-1} (E^{n+1}_{+})},
$$
 where $c$ is a positive constant independent of the function $u$ and
  $s_j = s + 2m-m_j - k+ \Re \mu - \dfrac{1}{2},$ holds.
\end{theorem}

\section{General Weight Boundary-Value
  Problems for Singular Elliptic Equations}\label{sec13}

\subsection{Half-space boundary-value
  problems; constant coefficients}\label{sec13.1}

 In the half-space
  $E_{+}^{n+1}$, consider the equation
\begin{equation}
    A \lr {\frac{1}{i} D_{x'}, \frac{1}{i^2} B_y}  u \equiv \sum\limits_{|\alpha'|
     + 2 \alpha_{n+1} = 2 m}  a_{\alpha} \lr{\frac{1}{i} D_{x'}}^{\alpha'}
     \lr{\frac{1}{i^2} B_y}^{\alpha_{n+1}} u(x) = f(x),
    \label{6.1.1}
\end{equation}
    where $B$ is the Bessel operator with a complex parameter $\nu.$
    On the hyperplane $\{y = 0\},$ consider the two types of boundary-value
    conditions:
\begin{equation}
\sigma_{\nu} G_j \lr {\frac{1}{i} D_{x'}, \frac{1}{i^2} B_y} B_y^k
\left. u \right|_{y=0}  = g_j (x'), \  j=0, \dots, m-1,
\label{6.1.2}
\end{equation}
    under the assumption that $\Re \nu > 0$ or $\nu = 0,$ and
\begin{equation}
y^{2 \nu +1} D_y G_j \lr {\frac{1}{i} D_{x'}, \frac{1}{i^2} B_y}
B_y^k  \left. u \right|_{y=0} = g_j (x'), \  j=0, \dots, m-1,
\label{6.1.3}
\end{equation}
   under the assumption that $\Re \nu \geq 0.$
   The operators $G_j,$ $j=0, \dots, m-1,$ have the form
$$
G_j \lr {\frac{1}{i} D_{x'}, \frac{1}{i^2} B_y}  =
\sum\limits_{|\alpha'| + 2 \alpha_{n+1} = m_j}  g_{\alpha j}
\lr{\frac{1}{i} D_{x'}}^{\alpha'} \lr{\frac{1}{i^2}
B_y}^{\alpha_{n+1}}.
$$
    The specified  boundary-value
  problems generate the operators $\mathfrak{U}_{\nu}$
 and $\mathfrak{U'}_{\nu}$ defined by the relations
$$
\mathfrak{U}_{\nu} u = \left\{ Au, \sigma_{\nu} \left. G_0 B^k_y u
\right|_{y=0}, \dots,   \sigma_{\nu} \left. G_{m-1} B^k_y u
\right|_{y=0}   \right\}
$$
    and
$$
\mathfrak{U'}_{\nu} u =  \left\{ Au, y^{2 \nu +1} D_y \sigma_{\nu}
\left. G_0 B^k_y u \right|_{y=0}, \dots,  y^{2 \nu +1} D_y  \left.
G_{m-1} B^k_y u \right|_{y=0}   \right\}.
$$
    Our goal is to construct regularizers for
$\mathfrak{U}_{\nu}$ and $\mathfrak{U'}_{\nu}.$
 To do that, we use the method of transmutation operators.

 Introduce  the space
$$
\mathcal{H}^s_{\nu} \lr{E^{n+1}_{+}, E^n, m} = H^s_{\nu}
\lr{E^{n+1}_{+}} \times \prod\limits_{j=0}^{m-1} H^{2m+s-m_j - 2 k
- 1 + \Re \nu} (E^n)
$$
    and endow it with the direct-product
  topology. Assume that
\begin{equation}
2 k + 1 - \Re \nu > 0 \quad \textrm{and}  \quad s+2 m -2 k- 1 +
\Re \nu - \max\limits_j m_j > 0. \label{6.1.4}
\end{equation}
  Then, from results of Sec. \ref{sec3}, it follows that the operators
$\mathfrak{U}_{\nu}$ and $\mathfrak{U'}_{\nu}$ continuously map
the spaces $H^{s+2 m}_{\nu} \lr{E^{n+1}_{+}}$ into the spaces
$H^{s}_{\nu} \lr{E^{n+1}_{+}, E^n, m}.$

Below, we study the operator $\mathfrak{U}_{\nu}$ in detail.
 The operator $\mathfrak{U'}_{\nu}$ is studied in the same way.

Apart from the operator $\mathfrak{U}_{\nu}$, consider the
operator $\mathfrak{U}$ defined by the relation
$$
    \begin{gathered}
\mathfrak{U} u = \left\{ A \lr{ \frac{1}{i} D_{x'}, \frac{1}{i^2}
D^2_{y}} u, G_0  \lr{ \frac{1}{i} D_{x'}, \frac{1}{i^2} D^2_{y}}
D_y^{2 k } \left. I^{\nu -\frac{1}{2}}_e u \right|_{y=0}, \dots,
\right.
 \\
\left. G_{m-1}  \lr{ \frac{1}{i} D_{x'},  \frac{1}{i^2} D^2_{y}}
D_y^{2 k }  \left. I^{\nu -\frac{1}{2}}_e u \right|_{y=0}
\right\}.
    \end{gathered}
$$
    In the previous section, it is shown that this operator continuously maps the space
     $H^{s+2 m}\lr{E^{n+1}_{+}}$ into the space
$$
\mathcal{H}^s \lr{E^{n+1}_{+}, E^n, m} =   H^{s} \lr{E^{n+1}_{+}}
\times \prod\limits_{j=0}^{m-1} H^{s+2m-m_j - 2 k + \Re \nu -1}
\lr{E^{n}}
$$
    endowed with the direct-product
    topology.

Introduce transmutation operators $\mathfrak{B}_{\nu, e}$ and
$\mathfrak{G}_{\nu, e}$ related to the studied  boundary-value
  problem.
  This merely extends the domains of the transmutation operators  $P_{\nu, e}$
  and
$S_{\nu, e}$ introduced and studied above.
    Let $\Phi$ denote the function collection
$\Phi=\{f, g_0, \dots, g_{m-1}\}.$
    Assign
$\mathfrak{B}_{\nu, e} \Phi=\{P_{\nu, e} f, c_{\nu} g_0, \dots,
c_{\nu} g_{m-1}\} $
    and
$ \mathfrak{G}_{\nu, e} \Phi=\{S_{\nu, e} f, \frac{1}{c_{\nu}}
g_0, \dots, \frac{1}{c_{\nu}} g_{m-1}\}, $ where $c_{\nu} = 2 \nu$
 provided that $\Re \nu > 0,$   $c_0=1,$ and the operators $P_{\nu,
 e}$ and
$S_{\nu, e}$ are defined in Sec. \ref{sec1}.
    The operator $\mathfrak{B}_{\nu, e}$ isomorphically maps the
    space $\mathcal{H}^s$ into $\mathcal{H}^s_{\nu},$ while  $\mathfrak{G}_{\nu, e}$
   realizes the inverse map.

 Note that the operators $\mathfrak{U}$ and
$\mathfrak{U}_{\nu}$  are bound by the simple relations
$$
\mathfrak{U}_{\nu} =  \mathfrak{B}_{\nu, e} \mathfrak{U} S_{\nu,
e}\quad \textrm{and} \quad  \mathfrak{U}  =  \mathfrak{G}_{\nu, e}
\mathfrak{U}_{\nu} P_{\nu, e}.
$$
    Once we have such relations, we can construct a regularizer of
    the singular operator $\mathfrak{U}_{\nu}$ by means of the relation
$\mathfrak{R}_{\nu} =  P_{\nu, e} \mathfrak{R} \mathfrak{G}_{\nu,
e},$ where $\mathfrak{R}$ is the regularizer of
    the operator $\mathfrak{U},$ constructed in Sec. \ref{sec5} under the assumption
    that the operators $A \Big( \dfrac{1}{i} D_{x'},  \dfrac{1}{i^2}
D^2_{y}\Big)$ and $G_j \Big( \dfrac{1}{i} D_{x'}, \dfrac{1}{i^2}
D^2_{y}\Big)$ form an elliptic collection.
    Let us show that the operator $\mathfrak{R}_{\nu}$ is a regularizer indeed.

    We have the relation
$$
\mathfrak{U}_{\nu} \mathfrak{R}_{\nu} =  \mathfrak{B}_{\nu, e}
\mathfrak{U} S_{\nu, e}  P_{\nu, e}  \mathfrak{R}
\mathfrak{G}_{\nu, e} =\mathfrak{B}_{\nu, e} \mathfrak{U}
\mathfrak{R} \mathfrak{G}_{\nu, e} = \mathfrak{B}_{\nu, e}  (I +
T_{\mbox{r}}) \mathfrak{G}_{\nu, e} = I + T_{\mbox{r}, \nu},
$$
    where $T_{\mbox{r}, \nu } = \mathfrak{B}_{\nu, e} T_{\mbox{r}}
\mathfrak{G}_{\nu, e}.$ Since $ T_{\mbox{r}} \in L(\mathcal{H}^s,
\mathcal{H}^{s+1}),$ it follows that
 $T_{\mbox{r}, \nu } \in L(\mathcal{H}^s_{\nu}, \mathcal{H}^{s+1}_{\nu}).$

 In the same way, we obtain that
$$
 \mathfrak{R}_{\nu} \mathfrak{U}_{\nu}
  = I + T_{\mbox{l}, \nu} ~\, \textrm{and}  ~\,   T_{\mbox{l}, \nu}
 = P_{\nu, e}  T_{\mbox{l}} S_{\nu, e}.
$$
    Since $T_{\mbox{l}} \in L(H^{s+2 m} (E^{n+1}_{+}), H^{s+2m+1}
(E^{n+1}_{+})),$ it follows that $T_{\mbox{l}, \nu} \in  L (H^{s+2
m}_{\nu} (E^{n+1}_{+}), H^{s+2m+1}_{\nu} (E^{n+1}_{+})).$

Thus, the following assertion is proved.

\begin{lemma} \label{teo: 6.1.1}
Let operators $A \Big( \dfrac{1}{i} D_{x'},  \dfrac{1}{i^2}
D^2_{y}\Big)$ and $G_j  \Big( \dfrac{1}{i} D_{x'},  \dfrac{1}{i^2}
D^2_{y}\Big)$ form an elliptic collection and relations
\eqref{6.1.4} be satisfied.
    Then the operators $\mathfrak{U}_{\nu}$ and $\mathfrak{U}_{\nu}'$ have
    regularizers continuously mapping the spaces $\mathcal{H}^s_{\nu}$
    into $H^{s+2 m}_{\nu}(E^{n+1}_{+})$ for all admissible $\nu.$
\end{lemma}

\subsection{Half-space boundary-value
  problems; low-variation
   coefficients}\label{sec13.2}

 In
  $E^{n+1}_{+}$, consider the  singular elliptic partial differential
  equation
\begin{equation}
    A \lr {x, \frac{1}{i} D_{x'}, \frac{1}{i^2} B_y}
     u \equiv \sum\limits_{|\alpha'| + 2 \alpha_{n+1} \leq 2 m}
      a_{\alpha} (x) \lr{\frac{1}{i} D_{x'}}^{\alpha'}
      \lr{\frac{1}{i^2} B_y}^{\alpha_{n+1}} u = f(x).
    \label{6.2.1}
\end{equation}
    If $\Re \nu >0$ and $\nu=0$, then  boundary-value
    conditions of the following kind are added to Eq. \eqref{6.2.1}:
\begin{multline}
\sigma_{\nu} (y) G_j \lr {\frac{1}{i} D_{x'}, \frac{1}{i^2} B_y}
B_y^k \left. u \right|_{y=0}
\equiv \sigma_{\nu}(y) \hspace{-3ex}\sum\limits_{|\alpha'| + 2 \alpha_{n+1} \leq m_j}\hspace{-3ex}  g_{\alpha j} (x') \lr{\frac{1}{i} D_{x'}}^{\alpha'} \lr{\frac{1}{i^2} B_y}^{\alpha_{n+1}}  B_y^k \left. u \right|_{y=0}  = g_j (x'), \\
 j=0, \dots, m-1.
\label{6.2.2}
\end{multline}
    Also,   boundary-value
    conditions of the following kind are considered in the case where
     $\Re \nu \geq 0$:
\begin{equation}
y^{2 \nu +1} D_y G_j \lr {x', \frac{1}{i} D_{x'}, \frac{1}{i^2}
B_y} B_y^k  \left. u \right|_{y=0} = g_j (x'), \  j=0, \dots, m-1.
\label{6.2.3}
\end{equation}
 It is assumed that the coefficients  $a_{\alpha}$ and $g_{\alpha j}$ are
 infinitely differentiable and $\left. D^p_y a_{\alpha}\right|_{y=0} = 0$
 for $p= 1, 2, \dots$
 For consistency, we do not compute $\dfrac{1}{i^2}=-1$ in the above notation.

Assume that the following conditions are satisfied:
\begin{enumerate}
\item[(1)]
    the polynomials $A (0, \xi', \eta^2)$ and $G_j (0, \xi',\eta^2)$ form
    an elliptic collection (i.\,e., satisfy the Shapiro--Lopatinskii
     condition);

\item[(2)]
 if $|\alpha'| + 2 \alpha_{n+1} = 2 m$ and $x \in
\ov{E_{+}^{n+1}},$ then the inequalities $|a_{\alpha} (x) -
a_{\alpha} (0)| < \varepsilon,$ $\Big|\Big(\dfrac{1}{y}
D_{y}\Big)^p a_{\alpha} (x) \Big| < \varepsilon,$ and $1 \leq p
\leq 3 \Re \nu +s+1$ hold;

\item[(3)]
  there exists a positive $R$ such that  $a_{\alpha} (x) = a_{\alpha}(0)$
  provided that $|x| \geq R$ and $|\alpha'| + 2 \alpha_{n+1} =2 m;$

\item[(4)]
 if $|\alpha'| + 2 \alpha_{n+1} < 2 m,$ then
$a_{\alpha} (x) \in \mathring{C}^{\infty} (\ov{E^{n+1}_{+}});$

\item[(5)]
 if $j = 0, \dots, m-1$ and $|\alpha'| + 2 \alpha_{n+1} =
m_j,$ then $|g_{\alpha_j } (x) - g_{\alpha_j} (0)| < \varepsilon$
for all $x' \in E^n$;

\item[(6)]
  if $j = 0, \dots, m-1$ and $|\alpha'| + 2 \alpha_{n+1} =
m_j,$  then $g_{\alpha_j } (x')=g_{\alpha_j} (0)$ provided that
$|x'| \geq R;$

\item[(7)]
  if $j = 0, \dots, m-1$ and $|\alpha'| + 2 \alpha_{n+1} =
m_j,$  then  $g_{\alpha_j } (x') \in \mathring{C}^{\infty} (E^n).$
\end{enumerate}
  Define the operators $\mathfrak{U}_{\nu}$ and $\mathfrak{U}_{\nu}'$
 by the relations
$$
\mathfrak{U}_{\nu} u = \left\{ Au, \sigma_{\nu}  G_0 \left. B^k_y
u  \right|_{y=0},  \dots, \sigma_{\nu}  G_{m-1} \left. B^k u
\right|_{y=0}  \right\}
$$
    and
$$
\mathfrak{U}_{\nu}' u = \left\{ Au, y^{2 \nu + 1}    \left. D_y
G_0 B^k u  \right|_{y=0}, \dots, y^{2 \nu + 1}   \left. D_y
G_{m-1} B^k u  \right|_{y=0}  \right\}.
$$
 If the above assumptions (1)--(7)
 are satisfied, then $\mathfrak{U}_{\nu}$ and
$\mathfrak{U}_{\nu}'$ are called singular elliptic operators with
    $(\varepsilon, s)$-low-variation
    coefficients.
 Below, for sufficiently small $\varepsilon$, we construct regularizers for them.
 Only the operator $\mathfrak{U}_{\nu}$ is studied in detail because the operator
  $\mathfrak{U}_{\nu}'$ is studied in the same way.

    To apply the perturbation method, decompose the operator
$\mathfrak{U}_{\nu}$ into the sum of three terms:
 the
 operator $\mathfrak{U}_{\nu, 0}$ with constant coefficients, the operator
  $\mathfrak{U}_{\nu, 1}$ with a small norm, and the operator
$\mathfrak{U}_{\nu, 2}$ containing only low-order
 terms. Thus, we assign
$$
    \begin{gathered}
\mathfrak{U}_{\nu, 0} u = \left\{ A_0 \lr{0, \frac{1}{i} D_{x'},
\frac{1}{i^2} B_y} u, \sigma_{\nu} (y) G_{0, 0} \lr{0, \frac{1}{i}
D_{x'}, \frac{1}{i^2} B_y} \left. B^k_y u  \right|_{y=0}, \dots,
\right.
 \\
\left. \sigma_{\nu} (y) G_{m-1, 0} \lr{0, \frac{1}{i} D_{x'},
\frac{1}{i^2} B_y} \left. B^k_y u  \right|_{y=0} \right\},
    \end{gathered}
$$
    where
$$
 A_0 \lr{0, \frac{1}{i} D_{x'}, \frac{1}{i^2} B_y}
  = \sum\limits_{|\alpha'| + 2 \alpha_{n+1} = 2 m} a_{\alpha} (0)
    \lr{\frac{1}{i} D_{x'}}^{\alpha'} \lr{\frac{1}{i^2} B_y}^{\alpha_{n+1}}
$$
    and
$$
G_{j, 0} \lr{0, \frac{1}{i} D_{x'}, \frac{1}{i^2} B_y}  =
 \sum\limits_{|\alpha'| + 2 \alpha_{n+1} =  m_j} g_{\alpha j} (0)
    \lr{\frac{1}{i} D_{x'}}^{\alpha'} \lr{\frac{1}{i^2} B_y}^{\alpha_{n+1}},
$$
$$
j = 0, \dots, m-1.
$$
    The operator $\mathfrak{U}_{\nu, 1}$ has the form
$$
    \begin{gathered}
\mathfrak{U}_{\nu, 1} u = \left\{ A_1 \lr{x, \frac{1}{i} D_{x'},
\frac{1}{i^2} B_y} u, \sigma_{\nu} (y) G_{0, 1} \lr{x',
\frac{1}{i} D_{x'}, \frac{1}{i^2} B_y} \left. B^k_y u
\right|_{y=0},  \right.
 \\
\left. \dots, \sigma_{\nu} (y) G_{m-1, 1}  \lr{x', \frac{1}{i}
D_{x'}, \frac{1}{i^2} B_y} \left. B^k_y u  \right|_{y=0}\right\},
    \end{gathered}
$$
    where
$$
A_1 \lr{x, \frac{1}{i} D_{x'}, \frac{1}{i^2} B_y}   =
\sum\limits_{|\alpha'| + 2 \alpha_{n+1} = 2 m} (a_{\alpha} (x)-
a_{\alpha} (0))   \lr{\frac{1}{i} D_{x'}}^{\alpha'}
\lr{\frac{1}{i^2} B_y}^{\alpha_{n+1}}
$$
    and
$$
G_{j, 1} \lr{x', \frac{1}{i} D_{x'}, \frac{1}{i^2} B_y}   =
\sum\limits_{|\alpha'| + 2 \alpha_{n+1} =  m_j} (g_{\alpha j}
(x')-g_{\alpha j} (0))   \lr{\frac{1}{i} D_{x'}}^{\alpha'}
\lr{\frac{1}{i^2} B_y}^{\alpha_{n+1}},
$$
$$
j = 0, \dots, m-1.
$$
    The operator $\mathfrak{U}_{\nu, 2}$ is defined by the
    relation
$$
\mathfrak{U}_{\nu, 2} u =  \left\{ A_2 \lr{x, \frac{1}{i} D_{x'},
\frac{1}{i^2} B_y} u, \sigma_{\nu} (y) G_{0, 2} \lr{x',
\frac{1}{i} D_{x'}, \frac{1}{i^2} B_y} \left. B^k_y u
\right|_{y=0},  \right.
$$
$$
\dots, \left. \sigma_{\nu} (y) G_{m-1, 2}  \lr{x', \frac{1}{i}
D_{x'}, \frac{1}{i^2} B_y} \left. B^k_y u  \right|_{y=0} \right\},
$$
    where
$$
A_2 \lr{x, \frac{1}{i} D_{x'}, \frac{1}{i^2} B_y}   =
\sum\limits_{|\alpha'| + 2 \alpha_{n+1} \leq 2 m-1} a_{\alpha} (x)
\lr{\frac{1}{i} D_{x'}}^{\alpha'} \lr{\frac{1}{i^2}
B_y}^{\alpha_{n+1}}
$$
    and
$$
G_{j, 2} \lr{x', \frac{1}{i} D_{x'}, \frac{1}{i^2} B_y}  = \sum\limits_{|\alpha'| + 2 \alpha_{n+1} \leq m_j-1} g_{\alpha j} (x')   \lr{\frac{1}{i} D_{x'}}^{\alpha'} \lr{\frac{1}{i^2} B_y}^{\alpha_{n+1}},
$$
$$
j = 0, \dots, m-1,
$$
 and it is assumed that  $m_j \geq 1 $ in the last relation.
 If there exists a number $j$ such that $m_j=0,$ then we assign $G_{j,
2}=0.$

Hence, the following relation holds:
\begin{equation}
\mathfrak{U}_{\nu} =  \mathfrak{U}_{\nu, 0} + \mathfrak{U}_{\nu,
1}  +\mathfrak{U}_{\nu, 2}. \label{6.2.4}
\end{equation}
    In the sequel, it is assumed that relations \eqref{6.1.4} are satisfied.

    Now, note that the operator $\mathfrak{U}_{\nu, 0}$ is generated by a
    boundary-value problem with a homogeneous operator both in the equation and
  boundary-value conditions.
 In the previous section, a regularizer is constructed in the previous section;
    here, we denote it by $\mathfrak{R}_{\nu, 0}.$
    As we proved above, the operator
$\mathfrak{R}_{\nu, 0}$ belongs to $L \lr{\mathcal{H}^s_{\nu}
(E^{n+1}_{+}, E^n, m), H^{s+2m} (E^{n+1}_{+})}$ and satisfies the
relations
\begin{equation}
\mathfrak{U}_{\nu, 0}  \mathfrak{R}_{\nu, 0}= I +   T_{\mbox{r},
0}, \  \mathfrak{R}_{\nu, 0} \mathfrak{U}_{\nu, 0}  = I +
T_{\mbox{l}, 0}, \label{6.2.5}
\end{equation}
    where $I$ denotes the identity operator, while
$T_{\mbox{r}, 0}$ and $T_{\mbox{l}, 0}$ are smoothing operators,
 which means that $T_{\mbox{r}, 0} \in L \lr{H^{s+2m}_{\nu},
H^{s+2m+1}_{\nu} }$ and $T_{\mbox{l}, 0} \in L \lr{
\mathcal{H}^{s}_{\nu}, \mathcal{H}^{s+1}_{\nu} }.$

 This and expansion \eqref{6.2.4} yield relations of the kind
\begin{equation}
\mathfrak{U}_{\nu}  \mathfrak{R}_{\nu, 0}= I +   T_{\mbox{r}, 0} +
\mathfrak{U}_{\nu, 1} \mathfrak{R}_{\nu, 0}    +\mathfrak{U}_{\nu,
2} \mathfrak{R}_{\nu, 0} \label{6.2.6}
\end{equation}
    and
\begin{equation}
 \mathfrak{R}_{\nu, 0} \mathfrak{U}_{\nu} = I +   T_{\mbox{l}, 0} +
  \mathfrak{R}_{\nu, 0} \mathfrak{U}_{\nu, 1} +
   \mathfrak{R}_{\nu, 0} \mathfrak{U}_{\nu, 2}.
\label{6.2.7}
\end{equation}
    Below, we show that the operators
 $I + \mathfrak{U}_{\nu, 1}\mathfrak{R}_{\nu, 0}$ and
  $I +  \mathfrak{R}_{\nu, 0}\mathfrak{U}_{\nu, 1}$ have bounded
  inverse operators provided that  $\varepsilon$ is sufficiently small.
    Then, introducing the notation
\begin{equation}
\mathfrak{R}_{\nu, 0} \lr{I + \mathfrak{U}_{\nu, 1}
\mathfrak{R}_{\nu, 0}}^{-1}  = \mathfrak{R}_{\nu} \label{6.2.8}
\end{equation}
and taking into account that
\begin{equation}
\mathfrak{R}_{\nu, 0} \lr{I + \mathfrak{U}_{\nu, 1}
\mathfrak{R}_{\nu, 0}}^{-1}  = \lr{I + \mathfrak{R}_{\nu, 0}
\mathfrak{U}_{\nu, 1}}^{-1}  \mathfrak{R}_{\nu, 0},  \label{6.2.9}
\end{equation}
    we transform relations \eqref{6.2.6}-\eqref{6.2.7}
    to the form
\begin{equation}
\mathfrak{U}_{\nu}  \mathfrak{R}_{\nu}= I +   \lr{T_{\mbox{r}, 0}
+  \mathfrak{U}_{\nu, 2} \mathfrak{R}_{\nu, 0}}   \lr{I +
\mathfrak{U}_{\nu, 1} \mathfrak{R}_{\nu, 0}}^{-1} \equiv I +
{T_{\mbox{r}}}  \label{6.2.10}
\end{equation}
    and
\begin{equation}
\mathfrak{U}_{\nu}  \mathfrak{R}_{\nu}= I +  \lr{I +
\mathfrak{R}_{\nu, 0} \mathfrak{U}_{\nu, 1} }^{-1}
\lr{T_{\mbox{l}, 0} +   \mathfrak{R}_{\nu, 0} \mathfrak{U}_{\nu,
2}}    \equiv I +  {T_{\mbox{l}}},
    \label{6.2.11}
\end{equation}
    where $T_{\mbox{r}}$ and $T_{\mbox{l}}$ are smoothing operators.

To convert the operator $I + \mathfrak{U}_{\nu, 1}
\mathfrak{R}_{\nu, 0},$ it suffices to prove the convergence (in
the corresponding operator topology) of the Neumann series from the
    right-hand side of the relation
\begin{equation}
\lr{I + \mathfrak{U}_{\nu, 1} \mathfrak{R}_{\nu, 0} }^{-1}  =
\sum\limits_{l=0}^{\infty} (-1)^l \lr{\mathfrak{U}_{\nu, 1}
\mathfrak{R}_{\nu, 0}}^l.
    \label{6.2.12}
\end{equation}
 Let $\Phi = \{f, g_0, \dots, g_{m-1} \} \in\mathcal{H}^s_{\nu}.$
  For a while, assign  $u= \mathfrak{R}_{\nu,0} \Phi.$
  Then $u \in H^{s+2m}_{\nu} (E_{+}^{n+1})$ and
\begin{multline}
    \left\| \mathfrak{U}_{\nu, 1} \mathfrak{R}_{\nu, 0}
     \Phi \right\|^2_{\mathcal{H}^s_{\nu}}  =
      \left\| \mathfrak{U}_{\nu, 1}  u \right\|^2_{\mathcal{H}^s_{\nu}}
    \\=  \left\| A_{ 1}  u \right\|^2_{H^s_{\nu} (E_{+}^{n+1})} +
    \sum\limits_{j=0}^{m-1}
    \left\| \sigma_{\nu} G_{j, 1}  B^k_y u|_{y=0} \right\|^2_{H^{s+2m - m_j - 2 k +
     \Re \nu -1} (E^n)}.
    \label{6.2.13}
\end{multline}
    Let us estimate each term of the last relation.
    Let $c_j,$ $j=1, 2, \dots,$ denote constants independent of $\mathfrak{U}_{\nu, 1}$
    and $\Phi.$
    From Corollary \ref{cor: 3.2.1} and assumptions (2)-(3),
    we intermediately obtain that
\begin{equation}
    \left\| A_{ 1}  u \right\|^2_{H^s_{\nu} (E_{+}^{n+1})}
     \leq c_1 \varepsilon \| u \|_{ H^{2m}_{\nu} (E_{+}^{n+1})}.
    \label{6.2.14}
\end{equation}
    Applying the Leibnitz rule for the case where $s> 0,$ we find
    that
$$
    \begin{gathered}
    \left\| A_{ 1}  u \right\|^2_{H^s_{\nu} (E_{+}^{n+1})} = \sum\limits_{|\alpha'|
    + \alpha_{n+1} \leq s} \left| D_{x'}^{\alpha'} D^{\alpha_{n+1}}_y S_{\nu, e} A_1 u
     \right|^2_{L_2 (E_{+}^{n+1})}
     \\
\leq c_2 \lr{\sum\limits_{|\alpha'| + \alpha_{n+1} \leq s} \left|  D^{\alpha_{n+1}}_y S_{\nu, e} \lr{ A_1 D_{x'}^{\alpha'} u} \right|^2_{L_2 (E_{+}^{n+1})} + \sum\limits_{\substack{|\alpha'|>0 \\
|\alpha'| + \alpha_{n+1} \leq s}} \left|  D^{\alpha_{n+1}}_y
S_{\nu, e} \lr{ A_1^{\alpha'} u} \right|^2_{L_2 (E_{+}^{n+1})}},
    \end{gathered}
$$
    where
$$
A^{(\alpha')}_1 u = D_{x'}^{\alpha'} A_1 u - A_1  D_{x'}^{\alpha'} =
 \sum\limits_{\substack{|\gamma'| + |\delta'|= |\alpha'|\\
|\delta'|>0}}  \sum\limits_{|\beta'| + 2 \beta_{n+1} = 2 m} c
(\delta', \gamma') \lr{D^{\delta'}_{x'} a_{\beta} (x)}
D^{\beta'+\gamma'}_{x'} B_y^{\beta_{n+1}} u
$$
    and $c (\delta', \gamma')$ are known constants.
 The order of the operator $A^{(\alpha')}_1$ does not exceed
  $2m +|\alpha'| -1.$
  Therefore, estimating the terms according to theorem \ref{teo: 3.2.1},
  we obtain the inequalities
$$
    \begin{gathered}
\left\|  D^{\alpha_{n+1}}_y S_{\nu, e}  A_1 D_{x'}^{\alpha'} u
\right\|_{L_2 (E_{+}^{n+1})} \leq
 c_3 \, \| u  \|_{H^{s+2 m}_{\nu} (E_{+}^{n+1})} \max\limits_{l \leq |\alpha_{n+1}|
  +3N+1} \sup\limits_x
  \left| \lr{\frac{1}{y} D_y}^l (a_{\alpha} (x)-a_{\alpha} (0 )) \right|
  \\
\leq c_3 \,  \varepsilon \, \| u  \|_{H^{s+2 m}_{\nu}
(E_{+}^{n+1})}
    \end{gathered}
$$
    and
$$
\left\|  D^{\alpha_{n+1}}_y S_{\nu, e}  A_1^{\alpha'} u \right\|^2_{L_2 (E_{+}^{n+1})} \leq
 c_4 \, \| u  \|_{H^{s+2 m-1}_{\nu} (E_{+}^{n+1})} \max\limits_{l \leq 3N+ |\alpha_{n+1}| +1}
  \sup\limits_x \left| D_{x'}^{\delta'} \lr{\frac{1}{y} D_y}^l a_{\alpha} (x) \right|.
$$
    Then
\begin{equation}
\left\| A_{ 1}  u \right\|^2_{H^s_{\nu} (E_{+}^{n+1})} \leq c_5 \,
\varepsilon \, \| u \|_{ H^{2m}_{\nu} (E_{+}^{n+1})} + c_5 \, M_s
(a) \| u \|_{ H^{s+2m+1}_{\nu} (E_{+}^{n+1})} + c_6 \, M_s (a) \|
u \|_{ H^{s+2m+1}_{\nu} (E_{+}^{n+1})}, \label{6.2.15}
\end{equation}
 where
\begin{equation}
M_s (a) = \max\limits_{|\beta| \leq 3N + s +1}
\max\limits_{|\alpha|=2m} \sup\limits_{x \in E^{n+1}_{+}} \left|
D_{x'}^{\beta'} \lr{\frac{1}{y} D_y}^{\beta_{n+1}} a_{\alpha} (x)
\right| \label{6.2.16}
\end{equation}
    and a positive integer $N$ is such that  $\Re \nu < N + \dfrac{1}{2}.$

By Theorem \ref{teo: 3.3.2}, there exists a positive constant
$c_7$ such that it does not depend on the operators $G_{j, 1}$ and
function $u(x)$ and
\begin{equation}
 \left\| \sigma_{\nu} \left. G_{j, 1} B^k u \right|_{y=0}
  \right\|_{H^{s+2 m -m_j-2 k + \Re \nu -1} (E^n)}
   \leq c_7  \left\| G_{j, 1} u \right\|_{H^{s+2 m -m_j}_{\nu} (E^{n+1}_{+})}.
\label{6.2.17}
\end{equation}
 The coefficients of the operators  $G_{j, 1}$ do not depend on the last variable.
 Then, from assumptions (5)-(6),
 we obtain the estimate
$$
\| G_{j, 1} u \|_{H^0_{\nu} (E_{+}^{n+1})}  \leq c_8 \,
\varepsilon \, \| u \|_{H^{m_j}_{\nu} (E_{+}^{n+1})}.
$$
    By the Leibnitz rule, for $s > m_j -2 m$, this implies that
\begin{equation}
\left\|  G_{j, 1}  u \right\|^2_{H^{s+2m-m_j}_{\nu} (E_{+}^{n+1})}
\leq c_9 \, \varepsilon \, \| u \|_{ H^{s+2m}_{\nu} (E_{+}^{n+1})}
+ c_{10} \, \widetilde{M}_s (\widetilde{g}_j) \| u \|_{
H^{s+2m-1}_{\nu} (E_{+}^{n+1})}, \label{6.2.18}
\end{equation}
    where
\begin{equation}
\widetilde{M}_s (\widetilde{g}_j)  =  \max\limits_{|\beta'| \leq
s+2 m - m_j|} \max\limits_{|\alpha'| +\alpha_{n+1} = m_j}
\sup\limits_{x' \in E^n} \left| D^{\beta'}_{x'} g_{\alpha j} (x')
\right|. \label{6.2.19}
\end{equation}
    Since $u = \mathfrak{R}_{\nu, 0} \Phi,$ it follows from estimates \eqref{6.2.15},
    \eqref{6.2.17}, and \eqref{6.2.18} that
$$
\left\| \mathfrak{U}_{\nu, 1} \mathfrak{R}_{\nu, 0} \Phi
\right\|^2_{\mathcal{H}^s_{\nu}  (E_{+}^{n+1}, E^n, m)} \leq
c_{11} \, \varepsilon \, \left\|  \mathfrak{R}_{\nu, 0} \Phi
\right\|_{H^{s+2m}_{\nu} (E_{+}^{n+1})}+
 c_{12}  \, M_s \left\|  \mathfrak{R}_{\nu, 0}
  \Phi \right\|_{H^{s+2m-1}_{\nu} (E_{+}^{n+1})},
$$
    where $M_s = \max \left\{M_s (a), \widetilde{M}_s (\widetilde{g}_1),
\dots, \widetilde{M}_s (\widetilde{g}_{m-1}) \right\}.$
    Since $\mathfrak{R}_{\nu, 0}$ is a bounded operator, it
    follows that the last inequality takes the form
\begin{equation}
\left\| \mathfrak{U}_{\nu, 1} \mathfrak{R}_{\nu, 0} \Phi
\right\|_{\mathcal{H}^s_{\nu}}  \leq c_{13} \, \varepsilon \,
\left\| \Phi \right\|_{\mathcal{H}^{s}_{\nu}} +  c_{14} \, M_s
\left\| \Phi \right\|_{\mathcal{H}^{s-1}_{\nu}}.
 \label{6.2.20}
\end{equation}
    In this inequality change
 $\Phi$ for $ \mathfrak{U}_{\nu, 1}\mathfrak{R}_{\nu, 0} \Phi.$
 We obtain the inequality
\begin{equation}
\left\| (\mathfrak{U}_{\nu, 1} \mathfrak{R}_{\nu, 0})^2  \Phi
\right\|_{\mathcal{H}^s_{\nu}} \leq c_{13} \, \varepsilon \,
\left\| \mathfrak{U}_{\nu, 1} \mathfrak{R}_{\nu, 0} \Phi
\right\|_{\mathcal{H}^{s}_{\nu}} +  c_{14} \, M_s  \left\|
\mathfrak{U}_{\nu, 1} \mathfrak{R}_{\nu, 0} \Phi
\right\|_{\mathcal{H}^{s-1}_{\nu}}. \label{6.2.21}
\end{equation}
    To estimate the last term of the right-hand
    side, we assume (without loss of generality) that $\varepsilon \leq 1$ $M_s \geq 1.$
    Using the Erling--Nirenberg
 inequality and inequalities \eqref{6.2.14} and \eqref{6.2.17}, we
 find that
$$
    \begin{gathered}
\left\| \mathfrak{U}_{\nu, 1} \mathfrak{R}_{\nu, 0} \Phi
\right\|_{H^{s-1}_{\nu}}  \leq \left\| A_{1} \mathfrak{R}_{\nu, 0}
\Phi \right\|_{H^{s-1}_{\nu} (E^{n+1}_{+})} + c_{15}
\sum\limits_{j=0}^{m-1} \left\| G_{j,1} \mathfrak{R}_{\nu, 0} \Phi
\right\|_{H^{s+2 m -m_j-1}_{\nu} (E^{n+1}_{+})}
    \\
\leq \varepsilon_1 \left\| A_{1} \mathfrak{R}_{\nu, 0} \Phi
\right\|_{H^{s}_{\nu}  (E^{n+1}_{+})} + c (\varepsilon_1) \left\|
A_{1} \mathfrak{R}_{\nu, 0} \Phi \right\|_{H^{0}_{\nu}
(E^{n+1}_{+})}
    \\
 + \varepsilon_1 \sum\limits_{j=0}^{m-1}
 \left\| G_{j,1} \mathfrak{R}_{\nu, 0} \Phi \right\|_{H^{s+2 m -m_j}_{\nu} (E^{n+1}_{+})}
 + c (\varepsilon_1)   \sum\limits_{j=0}^{m-1} \left\| G_{j,1}
 \mathfrak{R}_{\nu, 0} \Phi \right\|_{H^{0}_{\nu} (E^{n+1}_{+})} \leq
 c_{16} \lr{\varepsilon_1 + \varepsilon  c (\varepsilon_1) } M_s \left\|
  \Phi \right\|_{H^s_{\nu}},
    \end{gathered}
$$
    where $\varepsilon_1$ is an arbitrary positive number and $c
(\varepsilon_1)$ depends on  $\varepsilon_1.$ Combining this with
inequalities \eqref{6.2.20}-\eqref{6.2.21}
 we arrive at the estimate
\begin{equation}
\left\| (\mathfrak{U}_{\nu, 1} \mathfrak{R}_{\nu, 0})^2 \Phi
\right\|_{\mathcal{H}^s_{\nu}}  \leq c_{17} M_s^2
\lr{\varepsilon_1 + \varepsilon  c (\varepsilon_1) }  \left\|
\Phi \right\|_{\mathcal{H}^{s}_{\nu}}. \label{6.2.22}
\end{equation}
    Select $\varepsilon_1$ to satisfy the relation $c_{17} M^2_s \varepsilon_1 =
\dfrac{1}{2}$ and assume that
$$
c (\varepsilon_1) c_{17} M_s^2  \varepsilon < \frac{1}{2}.
$$
    Then
$$
\left\| (\mathfrak{U}_{\nu, 1} \mathfrak{R}_{\nu, 0})^2 \Phi
\right\|_{\mathcal{H}^s_{\nu}}  \leq q  \left\|  \Phi
\right\|_{\mathcal{H}^{s}_{\nu}},
$$
where $0<q<1.$
    Then $ (\mathfrak{U}_{\nu, 1}\mathfrak{R}_{\nu, 0})^2$ is a contracting operator.
    This yields the convergence of the Neumann series \eqref{6.2.12} in the Banach space
    of linear operators bounded in $\mathcal{H}^s_{\nu}$.
    Hence, the operator
     $ (I+ \mathfrak{U}_{\nu, 1}\mathfrak{R}_{\nu, 0})^{-1}$
     is defined and continuous in the spaces $\mathcal{H}^s_{\nu}.$
    It suffices to note that, due to relations \eqref{6.2.9}--\eqref{6.2.11},
     the operator $\mathfrak{R}_{\nu}$ defined by the relation
$$
  \mathfrak{R}_{\nu} =
   \mathfrak{R}_{\nu, 0} (I+ \mathfrak{U}_{\nu,1} \mathfrak{R}_{\nu, 0})^{-1}
$$
 is a regularizer for the operator $\mathfrak{R}_{\nu}$ because, obviously, $T_{\mbox{r}}$
 and $T_{\mbox{l}}$ are smoothing operators.

Thus, the following assertion is proved.
\begin{theorem} \label{teo: 6.2.1}
 Let $\mathfrak{U}_{\nu}$ and  $\mathfrak{U'}_{\nu}$ be singular
 elliptic operators with
  $(\varepsilon, s)$-low-variation
  coefficients such that the relations
    $$
    2k+1-\Re \nu > 0
    $$
 and
    $$
    s+2m-2k-2+ \Re \nu - \max\limits_j m_j > 0
    $$
    are satisfied. Let $\varepsilon$ be sufficiently small.
    Then the operator $\mathfrak{U}_{\nu}$ $(\mathfrak{U'}_{\nu})$
    has a regularizer from the space $L
\big(\mathcal{H}^s_{\nu} (E^{n+1}_{+}, E^n, m), H^{s+2m}_{\nu}
(E^{n+1}_{+})\big)$ provided that $\Re \nu > 0$ or  $\nu=0$ $(\Re
\nu \geq 0).$
\end{theorem}

\subsection{Weight boundary-value
    problems on bounded domains}\label{sec13.3}

    Let $\Omega$ be a bounded domain of the half-space
 $E^{n+1}_{+}$ such that its boundary is smooth and all the assumptions from
 Sec. \ref{sec4} are satisfied for it.

Let  $A$ denote a $2 m$-order
 operator elliptic in the domain $\Omega$ and such that its complex-valued
 coefficients are infinitely differentiable in $\Omega$.
 Assume that, in each local coordinate system, the operator $A$ admits
 the following representation near the boundary:
\begin{equation}
A= A \lr{x, \frac{1}{i} D_{x'},  \frac{1}{i^2} B_y} =
\sum\limits_{|\alpha'|+2 \alpha_{n+1} \leq 2m} a_{\alpha} (x)
\lr{\frac{1}{i} D_{x'}}^{\alpha'} \lr{\frac{1}{i^2}
B_y}^{\alpha_{n+1}}. \label{6.3.1}
\end{equation}
    Assume that the coefficients $a_{\alpha} (x)$ are infinitely differentiable
    up to the hyperplane  $\{y=0\}$ and the condition
\begin{equation}
D_y^p a_{\alpha} (x) =0, \  y=0, \  p=1,2, \dots \label{6.3.2}
\end{equation}
    is satisfied.

In each local coordinate system, introduce $m$
 boundary-value operators
$G_j,$ $j=0, 1, \dots, m-1,$ by relations of the kind
$$
G_j \lr{x', \frac{1}{i} D_{x'},  \frac{1}{i^2} B_y} =
\sum\limits_{|\alpha'|+2 \alpha_{n+1} \leq m_j} g_{\alpha j} (x')
\lr{\frac{1}{i} D_{x'}}^{\alpha'} \lr{\frac{1}{i^2}
B_y}^{\alpha_{n+1}},
$$
    where $g_{\alpha j}$  are infinitely differentiable coefficients
of variable $x' \in E^n.$

 Assume that the polynomials
$$
A_{0} (x, \xi', \eta^2) = \sum\limits_{|\alpha'|+2 \alpha_{n+1}
\leq 2m} a_{\alpha} (x)  {\xi'}^{\alpha'} \eta^{2 \alpha_{n+1}}
$$
    and
$$
G_{j, 0} (x, \xi', \eta^2) = \sum\limits_{|\alpha'|+2 \alpha_{n+1}
\leq m_j}  g_{\alpha j} (x') {\xi'}^{\alpha'} \eta^{2
\alpha_{n+1}},
$$
    where $\xi' \in E_n,$ $j=0, \dots, m-1,$ form an elliptic collection for each
    fixed point of the boundary.

Consider a boundary-value
 problem of the kind
\begin{equation}
\begin{cases}
A u = f, & \\
\left. \sigma_{\nu} \widetilde{G}_j u \right|_{\pr \Omega} = g_j,
& j=0, \dots, m-1,
\end{cases}
\label{6.3.3}
\end{equation}
where the operators $\widetilde{G}_j$ have the form
$$
\widetilde{G}_j =  G_j \lr{x', \frac{1}{i} D_{x'},  \frac{1}{i^2} B_y} B_y^k
$$
 in each local coordinate system.

In the sequel, it is assumed that
\begin{equation}
2 k +1 - \Re \nu >0  \textrm{ and }  s+2 m - 2 k-1 - \Re \nu -
\max\limits_j m_j > 0. \label{6.3.4}
\end{equation}
    The main goal of this section (and the whole chapter) is to prove
      the  Noetherian property of the set
 boundary-value problem (another boundary-value
 problem is considered at the end of the section).
 To do that, we have to construct a regularizer for the operator
$$
\mathfrak{U}_{\nu}: u \to  \mathfrak{U} u =  \{Au, \left.
\sigma_{\nu} \, \widetilde{G}_0 \right|_{\pr \Omega}, \dots,
\left. \sigma_{\nu}\, \widetilde{G}_{m-1} \right|_{\pr \Omega}
\}.
$$
  From results of Sec. \ref{sec4}, it follows that the operator $\mathfrak{U}_{\nu}$
  continuously maps the space $H^{s+2m}_{\nu} (\Omega)$ into the  space
$$
\mathcal{H}^{s}_{\nu} (\Omega, \pr \Omega, m) = H^{s}_{\nu}
(\Omega)  \times \prod\limits_{j=0}^{m-1} H^{s+2m-2k -1 + \Re \nu - m_j} (\pr \Omega)
$$
 endowed with the direct-product
 topology.

    The regularizer is constructed locally. To do that, we need a
    special continuation of an operator $\mathfrak{U}_{\nu}$ of
    the kind
\begin{equation}
\mathfrak{U}_{\nu} u =  \{Au, \left. \sigma_{\nu} \,
\widetilde{G}_0 \right|_{y=0}, \dots, \left. \sigma_{\nu} \,
\widetilde{G}_{m-1} \right|_{y =0}  \}, \label{6.3.5}
\end{equation}
    given by the above relations in a local coordinate system,
e.\,g., in a neighborhood of a point $x^0$ of the hyperplane
$\{y=0\}.$ For definiteness, let the coefficients of operator
\eqref{6.3.5} be defined in a semicube  $K^{+}_{\delta_0} (x^0),$
$\delta_0
>0,$ where
$$
K^{+}_{\delta_0} = \{x = (x', y):x \in E^{n+1}_{+},~|x_p - x^0_p| < \delta,~p=1,\dots,n+1 \}.
$$
    For an arbitrary $\varepsilon \in (0, 1)$, introduce the
    following operator $\mathfrak{U}_{\nu}^{x^0, \varepsilon}$
    corresponding to the principal part of operator \eqref{6.3.5}:
\begin{equation}
\mathfrak{U}_{\nu, 0}^{x^0, \varepsilon} u = \{A^{x^0, \varepsilon} u, \left. \sigma_{\nu} \, \widetilde{G}_0^{x^0, \varepsilon} \right|_{y=0}, \dots, \left. \sigma_{\nu} \, \widetilde{G}_{m-1}^{x^0, \varepsilon} \right|_{y =0}  \},
\label{6.3.6}
\end{equation}
where
\begin{multline*}
A^{x^0, \varepsilon} u(x) = \hspace{-3ex}\sum\limits_{|\alpha'|+ 2
\alpha_{n+1} = 2 m }\hspace{-3ex} \varphi_{n+1} \lr{
\frac{x-x^0}{\delta_0}}
 \lr{  \lr{a_{\alpha} (\varepsilon (x- x^0) + x^0) -a_{\alpha} (x^0)}
  \lr{\frac{1}{i} D_{x'}}^{\alpha'}}  \lr{\frac{1}{i^2} B_{y}}^{\alpha_{n+1}}\hspace{-2ex} u
 {}
\\
{}+  \sum\limits_{|\alpha'|+ 2 \alpha_{n+1} = 2 m } a_{\alpha}
(x^0) \lr{\frac{1}{i} D_{x'}}^{\alpha'}  \lr{\frac{1}{i^2}
B_{y}}^{\alpha_{n+1}} u
\end{multline*}
    and
\begin{multline*}
\widetilde{G}_j^{x^0, \varepsilon} u(x) = \sum\limits_{|\alpha'|+
2 \alpha_{n+1} = m_j } (-1)^k  \varphi_{n} \lr{
\frac{x-x^0}{\delta_0}}
\\
\times \lr{  \lr{g_{\alpha j} (\varepsilon (x- x^0) + x^0)
-g_{\alpha j} (x^0)}} t \lr{\frac{1}{i} D_{x'}}^{\alpha'}
\lr{\frac{1}{i^2} B_{y}}^{\alpha_{n+1} +k} u {}
\\
{}+
 \sum\limits_{|\alpha'|+ 2 \alpha_{n+1} = m_j } (-1)^k g_{\alpha j} (x^0)
  \lr{\frac{1}{i} D_{x'}}^{\alpha'}  \lr{\frac{1}{i^2} B_{y}}^{\alpha_{n+1}+k} u,
\quad j=0, \dots, m-1.
\end{multline*}
    To simplify the notation for boundary operators, we use $x^0$
 to denote a point $(x_1^0, \dots, x_n^0) \in E^n$ and
$\varphi_p,$ $p=n, n+1,$  to denote a function of the kind
$\varphi_p (x) = \varphi (x_1)   \dots     \varphi (x_p),$ where the
 one-variable function $ \varphi (t)$ belongs to $\mathring{C}^{\infty} (E^1),$
  $ \varphi (t)=1$ for $t \leq
\dfrac{1}{2},$ and $ \varphi (t)=0$ for $t \geq 1.$

Thus, the operator $\mathfrak{U}_{\nu, 0}^{x^0, \varepsilon}$ is
defined on the whole half-space
  $E^{n+1}_{+}$ now.
    Its infinitely differentiable coefficients possess the following properties.
    From \eqref{6.3.2}, we obtain the estimates
$$
\left|  D_{x'}^{\beta'} \lr{ \frac{1}{y} D_y}^{\beta_{n+1}} \lr{
\varphi_{n+1}  \lr{{ \frac{x-x^0}{\delta_0}}} \lr{a_{\alpha}
(\varepsilon (x- x^0) + x^0) -a_{\alpha} (x^0)} } \right| \leq
c_{\beta} \, \varepsilon
$$
    and
$$
\left|  D_{x'}^{\beta'} \lr{ \varphi_{n+1} \lr{{
\frac{x'-x^0}{\delta_0}}}  \lr{g_{\alpha j} (\varepsilon (x- x^0)
+ x^0) -g_{\alpha j} (x^0)} } \right| \leq c_{\beta'} \,
\varepsilon,
$$
    where $c_{\beta}$ and $c_{\beta'} $ are positive constants
    independent of
$\varepsilon \in (0, 1)$ and $x \in \ov{E_{+}^{n+1}}.$ Hence,
according to Theorem \ref{teo: 6.2.1}, the operator
$\mathfrak{U}_{\nu, 0}^{x^0, \varepsilon}$ is a singular operator
with
    $(\varepsilon, s)$-low-variation
    coefficients.
    Therefore, if $\varepsilon := \varepsilon_0$ is sufficiently small, then it has
    a regularizer $\mathfrak{R}_{\nu, 0}^{\varepsilon} \in
L \lr{ \mathcal{H}^s_{\nu} (E^{n+1}_{+}, E^n, m), H^{s+2m}_{\nu}
(E^{n+1}_{+})},$ i.\,e., the relations
\begin{equation}
\begin{array}{l}
\mathfrak{U}_{\nu, 0}^{x^0, \varepsilon} \mathfrak{R}_{\nu}^{x^0,
\varepsilon} \Phi = \Phi + T_{\mbox{r}, \nu, 0} \Phi, \  \Phi \in
\mathcal{H}^s_{\nu} (E^{n+1}_{+}, E^n, m), \\
\mathfrak{R}_{\nu}^{x^0, \varepsilon} \mathfrak{U}_{\nu, 0}^{x^0,
\varepsilon} u = u + T_{\mbox{l}, \nu, 0} u, \  u \in
H^{s+2m}_{\nu} (E^{n+1}_{+}),
\end{array}
\label{6.3.7}
\end{equation}
    where $T_{\mbox{r}, \nu, 0}$ and $T_{\mbox{l}, \nu, 0}$ are
    smoothing operators, hold.

Introduce the following extension operator $Q_{\varepsilon}$:
$$
Q_{\varepsilon}: u \to Q_{\varepsilon} u(x) = u \lr{\frac{x-x^0}{\varepsilon}+x^0}.
$$
 The operator $Q_{\frac{1}{\varepsilon}}$ is its inverse. Also,
introduce the operator
$$
    \begin{gathered}
Q_{\varepsilon, m}: \Phi = \{f(x), g_0 (x'), \dots, g_{m-1}(x') \}
\to Q_{\varepsilon, m} \Phi
    \\
= \left\{ \varepsilon^{-2m } f_0
\lr{\frac{x-x^0}{\varepsilon}+x^0},\  \varepsilon^{-2 \nu - m_0 -
2k } g_0 \lr{\frac{x-x^0}{\varepsilon}+x^0},\ \varepsilon^{-2 \nu
- m_{m-1} - 2k } g_{m-1} \lr{\frac{x-x^0}{\varepsilon}+x^0}
\right\}.
    \end{gathered}
$$
 The operator $Q_{\frac{1}{\varepsilon}, m}$ is its inverse. If $\varepsilon>0$, then
 the operators $Q_{\varepsilon}$ and
$Q_{\varepsilon, m}$ isomorphically map the spaces $H_{\nu}^s
(E_{+}^{n+1})$ and $\mathcal{H}_{\nu}^s (E_{+}^{n+1}, E^n, m)$
(respectively) onto itself.

Apply the operator $Q_{\varepsilon_0, m}$ to the first relation of
\eqref{6.3.7} and change $\Phi$ for $Q_{\frac{1}{\varepsilon}, m}
\Phi$ in this relation.
    Then simple transformations yield a relation of the kind
\begin{equation}
\mathfrak{U}_{\nu, 0}^{x^0} \mathfrak{R}_{\nu}^{x^0} \Phi = \Phi +
T_{\mbox{r}, \nu, 0}^{x^0} \Phi, \  \Phi \in \mathcal{H}^s_{\nu}
(E^{n+1}_{+}, E^n, m), \label{6.3.8}
\end{equation}
where $T_{\mbox{r}, \nu, 0}^{x^0}  = Q_{\varepsilon_0, m}
T_{\mbox{r}, \nu, 0}^{x^0, \varepsilon_0}
Q_{\frac{1}{\varepsilon_0}, m} $ is a smoothing operator (as well
as the operator $T_{\mbox{r}, \nu, 0}^{x^0, \varepsilon_0}$). The
operator $\mathfrak{R}_{\nu}^{x^0}$ has the form
$\mathfrak{R}_{\nu}^{x^0} = Q_{\varepsilon_0}
\mathfrak{R}_{\nu}^{x^0, \varepsilon_0}
Q_{\frac{1}{\varepsilon_0}, m} $ and belongs to the same class $L
\big(\mathcal{H}^s_{\nu} (E^{n+1}_{+}, E^n, m),$ $H^{s+2 m}_{\nu}
(E^{n+1}_{+})\big)$ as the operator $\mathfrak{R}_{\nu}^{x^0,
\varepsilon_0}.$ The operator $\mathfrak{U}_{\nu, 0}^{x^0}$ has
the form
$$
    \begin{gathered}
\mathfrak{U}_{\nu, 0}^{x^0} u  = \left\{ \sum\limits_{|\alpha'|+ 2
\alpha_{n+1} = 2 m }  \lr{ \varphi_{n+1} \lr{
\frac{x-x^0}{\varepsilon_0 \delta_0}} \lr{a_{\alpha} (x) -
a_{\alpha} (x^0)} -a_{\alpha} (x^0)} \right.
    \\
\times \lr{\frac{1}{i} D_{x'}}^{\alpha'}  \lr{\frac{1}{i^2}
B_{y}}^{\alpha_{n+1}} u, \sigma_{\nu} (y) \sum\limits_{|\alpha'|+
2 \alpha_{n+1} = m_0 } (-1)^k   \left( \varphi_{n} \lr{
\frac{x-x^0}{\varepsilon_0 \delta_0}}  \right.
 \\
\times \left. \lr{g_{\alpha 0} (x) - g_{\alpha 0} (x^0)}
-g_{\alpha 0} (x^0)\right) \lr{\frac{1}{i} D_{x'}}^{\alpha'}
\left. \lr{\frac{1}{i^2} B_{y}}^{\alpha_{n+1}+k} u \right|_{y=0},
 \\
\sigma_{\nu} (y) \sum\limits_{|\alpha'|+ 2 \alpha_{n+1} = m_{m-1}
} (-1)^k   \lr{ \varphi_{n} \lr{ \frac{x-x^0}{\varepsilon_0
\delta_0}} \lr{g_{\alpha m-1} (x) - g_{\alpha m-1} (x^0)}
-g_{\alpha m-1} (x^0)}
    \\
\left. \times \lr{\frac{1}{i} D_{x'}}^{\alpha'} \left.
\lr{\frac{1}{i^2} B_{y}}^{\alpha_{n+1}+k} u \right|_{y=0}
\right\}.
    \end{gathered}
$$
    The second relation of \eqref{6.3.7} is transformed in the same way.
    Note that the coefficients of the operator $\mathfrak{U}_{\nu, 0}^{x^0}$
    coincide with  the coefficients of the principal part of the operator $\mathfrak{U}_{\nu}$
 in the semicube $K^{+}_{\frac{\varepsilon_0 \delta_0}{2}}(x^0)$ (see \eqref{6.3.5}), which is
 a neighborhood of the point $x^0.$

    Now, define the operator $\mathfrak{U}_{\nu}^{x^0}$ by the relation
    $\mathfrak{U}_{\nu}^{x^0}= \mathfrak{U}_{\nu, 0}^{x^0} +
\mathfrak{U}_{\nu, 1}^{x^0},$
    where the operator $\mathfrak{U}_{\nu, 1}^{x^0}$ has the form
$$
    \begin{gathered}
\mathfrak{U}_{\nu, 1}^{x^0} u  = \left\{ \sum\limits_{|\alpha'|+ 2
\alpha_{n+1} < 2 m }  \varphi_{n+1} \lr{ \frac{x-x^0}{\delta_0}}
a_{\alpha} (x)  \lr{\frac{1}{i} D_{x'}}^{\alpha'}
\lr{\frac{1}{i^2} B_{y}}^{\alpha_{n+1}} u, \right.
 \\
\sigma_{\nu} (y) \sum\limits_{|\alpha'|+ 2 \alpha_{n+1} < m_0 }
(-1)^k   \varphi_{n} \lr{ \frac{x'-x^0}{ \delta_0}} g_{\alpha 0}
(x)  \lr{\frac{1}{i} D_{x'}}^{\alpha'} \left. \lr{\frac{1}{i^2}
B_{y}}^{\alpha_{n+1}+k} u \right|_{y=0},
 \\
 \sigma_{\nu} (y) \sum\limits_{|\alpha'|+ 2 \alpha_{n+1} < m_{m-1} } (-1)^k
  \varphi_{n} \lr{ \frac{x'-x^0}{ \delta_0}} g_{\alpha m-1} (x)
    \lr{\frac{1}{i} D_{x'}}^{\alpha'}
\left. \left. \lr{\frac{1}{i^2} B_{y}}^{\alpha_{n+1}+k} u
\right|_{y=0}\right\}.
    \end{gathered}
$$
    Thus, the operator $\mathfrak{U}_{\nu}^{x^0}$ has a regularizer
 $\mathfrak{R}_{\nu}^{x^0} \in L \lr{\mathcal{H}^s_{\nu} (E^{n+1}_{+}, E^n, m),
  H^{s+2 m}_{\nu} (E^{n+1}_{+})}$ and smoothing (in the corresponding spaces)
  operators  $T_{\mbox{r}, \nu}^{x^0}$ and $T_{\mbox{l}, \nu}^{x^0}$
    such that
\begin{equation}
\mathfrak{U}_{\nu}^{x^0}  \mathfrak{R}_{\nu}^{x^0} = I +
T_{\mbox{r}, \nu}^{x^0} \label{6.3.9}
\end{equation}
    and
\begin{equation}
\mathfrak{R}_{\nu}^{x^0} \mathfrak{U}_{\nu}^{x^0}   = I +
T_{\mbox{l}, \nu}^{x^0}. \label{6.3.10}
\end{equation}
    There exists a semicube
     $K^{+}_{\frac{\varepsilon_0\delta_0}{2}} (x^0)$ such that the coefficients
     of the operator $\mathfrak{U}_{\nu}^{x^0}$   coincide with  the coefficients
      of operator \eqref{6.3.5} in it. The desired
      continuation of the operator $\mathfrak{U}_{\nu}$ from \eqref{6.3.5}
      is constructed.

For each boundary point $x$, one can argue in the same way.
    Thus, one can construct a covering of the boundary by images of cubes
    $K^{+}_{\delta} (x),$ $x \in \{y=0\},$ $\delta= \delta(x)>0,$
    such that the coefficients of the operators
    $\mathfrak{U}_{\nu}^{x}$ and $\mathfrak{U}_{\nu}$ coincide each other
 (in a local coordinate system) in the semicube $K^{+}_{\delta / 2} (x)$.
 From the specified covering, select a finite subcovering of the boundary and, therefore,
 of its neighborhood. The preimages of the centers of the semicubes forming this
   finite subcovering are denoted by $\widetilde{x}^q,$ $q=1, \dots, \ov{q}.$
   Their images are denoted by $x^q = \varkappa_{l_q} \widetilde{x}^q.$
 The edge of the $q$th cube is denoted by $\delta_q.$

Consider the domain $\mathsf{y}: = \Omega \setminus
\bigcup\limits_{q=1}^{\ov{q}} \varkappa_{l_q}^1 K^{+}_{\delta_q}
(x^q)$ lying strictly inside the domain  $\Omega.$
    The domain $\mathsf{y}$ is to be processed as above.
    We can  omit this arguing because $A$ is an elliptic operator  with smooth coefficients
    inside the domain $\Omega$.
    Therefore, it suffices to refer, e.\,g., to \cite{83}, where   the following fact
     is proved: there exists a finite covering of the domain $\mathsf{y},$
e.\,g., by cubes $K_{\delta_q} (x^q),$ $q= \ov{q}+1, \dots,
\ov{\ov{q}},$ located in $\Omega$ together with their closures.
 Moreover, this covering possesses the following property: for each cube
  $K_{\delta_q} (x^q),$ $q> \ov{q},$ there exists a continuation
  operator $\mbox{R}^{x_q}$ (it can be constructed according to
  the scheme explained above) preserving the
smoothness and such that the operator
$$
A^{x_q} = \sum\limits_{|\alpha| \leq 2 m} \lr{\mbox{R}^{x_q}
b_{\alpha} (x) } \lr{\frac{1}{i} D_x}^{\alpha},
$$
    where $b_{\alpha}$ are the coefficients of the operator $A,$
    defined in the domain  $\Omega,$ has a regularizer $\mathfrak{R}^{x^q} \in
L \lr{H^s (E^{n+1}),  H^{s+2m} (E^{n+1}) }.$
  This means the validity of the relations
$$
A^{x^q}  \mathfrak{R}^{x^q} = I + T^{x^q}_{\mbox{r}} \textrm{ and
} \mathfrak{R}^{x^q}  A^{x^q} = I + T^{x^q}_{\mbox{l}},
$$
    where $ T^{x^q}_{\mbox{r}} \in L \lr{H^s (E^{n+1}), H^{s+1}
(E^{n+1})}$  and $ T^{x^q}_{\mbox{l}} \in L \lr{H^{s+2m}
(E^{n+1}), H^{s+2m+1} (E^{n+1})}.$ Note that  the coefficients of
the operator  $A^{x^q}$ coincide with  the coefficients of the
operator  $A$ in the cube $K_{\delta/2} (x^q),$ $q > \ov{q}.$

    Thus, we have a covering of   $\ov{\Omega}$ by a finite set of domains
     $\varkappa_{l_q}^{-1} K_{\delta} (x^q)$ such that the diffeomorphisms $x_{l_q}$
     are treated as identity maps  for $q > \ov{q}$.
  Without loss of generality, assume that the domains $\varkappa_{l_q}^{-1} K_{\delta_q/4}
(x^q)$ form a covering of  $\ov{\Omega}$ as.
    Let $\{h_q\}_{q=1}^{\ov{\ov{q}}}$ denote the partition of unity subordinated to the last
     covering. Introduce the set of compactly supported infinitely differentiable functions
      $\{\psi_q\}_{q=1}^{\ov{\ov{q}}}$ such that
\begin{equation}
\supp \psi_q \subset \varkappa_{l_q}^{-1} K_{\delta_q / 2} (x^q)
\label{6.3.11}
    and
\end{equation}
\begin{equation}
\psi_q (x) h_q (x) \equiv h_q (x),
\label{6.3.12}
\end{equation}
    where  $q = 1, \dots, \ov{\ov{q}}.$
    Assume that there exists a neighborhood of the hyperplane $\{y=0\}$ such that
\begin{equation}
D_y \psi_q = 0
\label{6.3.13}
\end{equation}
    in each local coordinate system in this neighborhood.
    In fact, the existence of such a collection of functions
    is proved within the proof of Lemma \ref{lem: 4.1.1.}.

 Let $\Phi = \{f, g_0, \dots, g_{m-1} \} \in \mathcal{H}^s
(\Omega, \pr \Omega, m).$ Then the desired regularizer
$\mathfrak{R}_{\nu}$ of the operator $\mathfrak{U}_{\nu}$
    is defined by the relation
$$
\mathfrak{R}_{\nu} \Phi = \sum\limits_{q=1}^{\ov{q}} h_q
\varkappa_{l_q}^{-1}  \mathfrak{R}_{\nu}^{x^q} \varkappa_{l_q}
\lr{\psi_q \Phi} +  \sum\limits_{q=\ov{q}+1}^{\ov{\ov{q}}} h_q
\mathfrak{R}^{x^q} (\psi_q f),
$$
    where $\psi_q \Phi = \{\psi_q f, \psi_q|_{\pr \Omega} \, g_0,
\dots, \psi_q|_{\pr \Omega} \, g_{m-1}  \}.$
    Using Theorem  \ref{teo: 3.2.1} and the assertions on norm equivalence,
  one can easily see that $ \mathfrak{R}_{\nu} \in L
\lr{\mathcal{H}_{\nu}^{s} (\Omega, \pr \Omega, m), H^{s+2m}
(\Omega)}.$

Let us show that the operator $\mathfrak{R}_{\nu}$ is a
regularizer. Let $u \in H^{s+2m}_{\nu} (\Omega).$ Then
\begin{multline}
\mathfrak{R}_{\nu} \mathfrak{U}_{\nu} u = \sum\limits_{q=1}^{\ov{q}} h_q \varkappa_{l_q}^{-1}
  \mathfrak{R}_{\nu}^{x^q} \varkappa_{l_q} \lr{\psi_q \mathfrak{U}_{\nu}  u}
 + \sum\limits_{q= \ov{q}+ 1}^{\ov{\ov{q}}} h_q \mathfrak{R}^{x^q}  \lr{\psi_q A  u}  \\
= \sum\limits_{q=1}^{\ov{q}} h_q \varkappa_{l_q}^{-1}
\mathfrak{R}_{\nu}^{x^q} \varkappa_{l_q} \lr{\psi_q
\mathfrak{U}_{\nu}^{x^q}  u} + \sum\limits_{q= \ov{q}+
1}^{\ov{\ov{q}}} h_q \mathfrak{R}^{x^q}  \lr{\psi_q A^{x^q}  u}
\label{6.3.14}
\end{multline}
    because the coefficients of the operators $A^{x^q}$ and $G_j^{x^q}$
    coincide with the coefficients of the operators $A$ and
$G_j $ (respectively) in the preimages of the cubes $K_{\delta_q /
2} (x^q)$ and, therefore, on the supports of the functions
$\psi_q.$ Further, from the Leibnitz rule, we have the relations
\begin{equation}
\varkappa_{l_q}  \lr{ \psi_q \mathfrak{U}_{\nu}^{x^q} u }  =
\varkappa_{l_q}  \lr{ \mathfrak{U}_{\nu}^{x^q} (\psi_q u) +
\mathfrak{\hat{U}}_{\nu}^{x^q} u  },~ q = 1, \dots, \ov{q},
\label{6.3.15}
\end{equation}
    and
\begin{equation}
\psi_{q}  A^{x^q} u  = A^{x^q} (\psi_q u) + \hat{A}^{x^q} u,~ q  =
\ov{q} + 1, \dots, \ov{\ov{q}}. \label{6.3.16}
\end{equation}
 Note the   Leibnitz rule for Bessel operators has the form
$$
B (\psi v) = v B \psi + \psi B v + 2 \lr{D_y \psi} D_y v,
$$
 whence, using the induction with respect to $r$, we deduce the relation
$$
B^r (\psi v) =  \psi B^r v + \sum\limits c_{\nu} (r_1, \dots, r_5)
y^{-r_1} \lr{ D_y^{r_2} B^{r_3} \psi } D_y^{r_4} B^{r_5} v,
$$
    where the summing is implemented with respect to all nonnegative integer indices
     $r_1, \dots, r_5$ such that $r_1+r_2+2 r_3 + r_4
+2 r_5 = 2 r-1,$ $r_2+2 r_3 > 0,$ and $r_{3, 4} =0,1.$ The
constants  $c_{\nu} (r_1, \dots, r_5)$ are known (for $\nu = -
\dfrac{1}{2}$, they coincide with the binomial coefficients).
    Combining this rule and the classical  Leibnitz rule with respect to the remaining
    variables, we arrive at relations \eqref{6.3.15}-\eqref{6.3.16},
    where the orders of the operators
$\mathfrak{\hat{U}}_{\nu}$ and $\hat{A}$ are one unit less than
 the orders of the operators $\mathfrak{U}_{\nu}$ and $A$ respectively.

By virtue of properties of the function $\psi$ (see relations
 \eqref{6.3.11}--\eqref{6.3.13}), it follows from Corollaries
  \ref{cor: 3.2.1}-\ref{cor: 3.2.2} that
$$
\varkappa_{l_q}   \mathfrak{\hat{U}}_{\nu}^{x^q} \in L
\lr{H^{s+2m}_{\nu} (E^{n+1}_{+}),  \mathcal{H}^{s+2m}_{\nu}
(E^{n+1}_{+}, E^n, m)}
$$
    and
$$
\hat{A}^{x^q} \in L \lr{H^{s+2m} (E^{n+1}), H^{s+1} (E^{n+1})}.
$$
    Since $\mathfrak{R}_{\nu}^{x^q}$ are regularizers, it follows
    that
\begin{equation}
\mathfrak{R}_{\nu} \mathfrak{U}_{\nu} u =
\sum\limits_{q=1}^{\ov{q}} h_q  \psi_q u + T_{\mbox{l}, \nu} u = u
+ T_{\mbox{l}, \nu} u, \label{6.3.17}
\end{equation}
 where $T_{\mbox{l}, \nu}$ is proved to be a smoothing operator.

Thus, it is proved that the operator $\mathfrak{R}_{\nu}$ is a
left-hand regularizer. In the same way, it is proved that it is a
 right-hand regularizer as well.

\begin{theorem}\label{teo: 6.3.1}
 Let $\Re \nu > 0$ or $\nu=0$ and relations \eqref{6.3.4} be satisfied.
 Then the boundary -value
 problem set by \eqref{6.3.3} has the Noetherian property.
 If $u \in H^{s+2m}_{\nu} (\Omega)$ and
   $ \mathfrak{U}_{\nu} u  \in \mathcal{H}_{\nu}^{s+p} (\Omega, \pr \Omega, m)$
 for a positive $p,$ then $u \in H^{s+2m+p}_{\nu} (\Omega)$ and
 the estimate
 \begin{equation}
 c \, \| u \|_{ H^{s+2m+p}_{\nu} (\Omega)} \leq  \|A u \|_{ H^{s+2m+p}_{\nu} (\Omega)}
 + \sum\limits_{j=0}^{m-1} \|\left. \sigma_{\nu} \widetilde{G}_j u\right|_{\pr \Omega}
  \|_{H^{s+p+2m-m_j-2k-1 + \Re \nu} (\pr \Omega)} + \| u \|_{ H^{s+2m}_{\nu} (\Omega)}
 \label{6.3.18}
 \end{equation}
 holds, where  $c$ is a positive constant independent of $u.$
\end{theorem}

\begin{proof}
    The Noetherian property of the boundary-value
    problem is equivalent to the existence of a regularizer (see, e.\,g.,
    abstract results from \cite{71}).
    The a priori estimate given by inequality \eqref{6.3.18} holds because
    the operator $R_{\nu}$ satisfies relation \eqref{6.3.17}, where $T_{\mbox{l}, \nu} $
 is a smoothing operator. The increasing smoothness assertion follows from the same relation
 as well. Really, if $u \in H^{s+2m}_{\nu} (\Omega),$ then $T_{\mbox{l}, \nu} u \in
H^{s+2m+1}_{\nu} (\Omega).$ Since $\mathfrak{U}_{\nu} u  \in
\mathcal{H}^{s+p},$ it follows that $\mathfrak{R}_{\nu}
\mathfrak{U}_{\nu} u \in H^{s+2m+p}_{\nu} (\Omega).$
 Then \eqref{6.3.16} implies that $u \in H^{s+2m+1}_{\nu} (\Omega).$
 Arguing as above, we show that  $u \in
H^{s+2m+p}_{\nu} (\Omega),$ which completes the proof of the
theorem.
\end{proof}
  Now, consider the singular elliptic equation \eqref{6.3.1} with weight boundary-value
  conditions of the kind
\begin{equation}
\sigma_{\nu + \frac{1}{2}} \widetilde{G}'_j u|_{\pr \Omega} = g_j, \  j = 0, \dots, m-1,
\label{6.3.19}
\end{equation}
    where $\Re \nu > 0$ and the operators  $\widetilde{G}'_j$ have the form
$$
\widetilde{G}'_j = D_y G_j \lr{x', \frac{1}{i} D_{x'},  \frac{1}{i^2} B_y } B_y^k
$$
    in each local coordinate system.

 If all above assumptions are satisfied, then the following assertion holds.

\begin{theorem}\label{teo: 6.3.2}
    Let $\Re \nu \geq 0$ and relations \eqref{6.3.4} be satisfied.
Then \eqref{6.3.1},\eqref{6.3.18} is a Noetherian boundary-value
 problem. If $u \in H^{s+2m}_{\nu} (\Omega)$ and
  $\mathfrak{U}_{\nu}' u \in \mathcal{H}_{\nu}^{s+p} (\Omega, \pr
\Omega, m),$ $p>0,$ then $u \in H^{s+2m+p}_{\nu} (\Omega)$ and the
estimate
$$
c \, \| u \|_{ H^{s+2m+p}_{\nu} (\Omega)} \leq  \|A u \|_{ H^{s+p}_{\nu} (\Omega)} +
 \sum\limits_{j=0}^{m-1} \|\left. \sigma_{\nu+ \frac{1}{2}}
  \widetilde{G}'_j u\right|_{\pr \Omega} \|_{H^{s+p+2m-m_j-2k-1 + \Re \nu} (\pr \Omega)}
   + \| u \|_{ H^{s+2m}_{\nu} (\Omega)}
$$
    holds, where $c$ is a positive constant independent of $u.$
\end{theorem}

   The proof of this theorem is entirely the same as the proof of Theorem \ref{teo: 6.3.1};
   here, we omit it.

\chapter[New Boundary-Value
    Problems\\
     for the Poisson Equation with Singularities at
    Isolated Points]{New Boundary-Value
    Problems for the Poisson Equation with Singularities at
    Isolated Points}\label{ch5}

    In this chapter, new boundary-value
    problems for the Poisson equation are considered.
 Their solutions might have singularities at isolated internal points treated
 as boundary ones.
 These singularities might be power singularities (of the pole type)
 or essential (infinite-order)
  singularities.
  We introduce and study new  Frechet-type
  functions spaces in bounded domains with smooth boundaries.
  These spaces possess the following three properties.
  First, they are wider that Sobolev spaces.
  Secondly, they contain all harmonic functions with arbitrary singularities at a finite set
  of fixed internal points (without loss of generality, we consider the case where such point
  is unique). Finally, they locally coincide with Sobolev spaces
  outside the singular point. To combine these properties in one space,
  we use multidimensional transmutation operators from the second chapter.
  We introduce the notion of a nonlocal (in a way) trace at a singular point;
  it becomes nontrivial only for functions singular at the said point.
  The direct and inverse trace theorems are proved in this chapter.
    In terms of the specified trace, we classify isolated singular
    points of harmonic functions and prove the main result of the
    chapter, which is the unique solvability of the corresponding
    boundary-value problem for the Poisson equation.
    Recall that a growth of arbitrary order is admitted for the solution and right-hand
    side of the equation at the singular point.

    At singular points, it is necessary to use the new nonlocal
    boundary-value condition; in the sequel, it is called the $\sigma$-trace.
    Also, one can propose the term $K$-trace
    (in honor of Katrakhov introduced and comprehensively studied this
    boundary-value condition).

\section{Function Spaces}\label{sec14}

\subsection{Definition and embedding theorems}\label{sec14.1}

Let $\Omega$ be a bounded domain in the space $E^n$ and its
boundary $\pr \Omega$ be smooth. Let the origin $\mathbf{0}$
 belong to $\Omega.$ Denote the domain  $\Omega\setminus \mathbf{0}$ by $\Omega_0$.
 Let $R_0$ be a positive number such that $\ov{U}_{2 R_0} \subset \Omega.$

    Let $X$ denote the set of functions $\chi (r) \in
\mathring{C}^{\infty} (\ov{E_{+}^1})$ such that $\chi (r) = 1$ for
$0 \leq r \leq 1,$  $\chi (r) = 0$ for $r \geq 2,$ and $\chi (r)$
monotonously decreases for $1<r<2.$ For each $\chi,$ define the
function $\chi_R,$ $R> 0,$ by the relation $\chi_R (r) = \chi
\Big(\dfrac{r}{R}\Big).$
    Let $\mathring{T}^{\infty} (\Omega_0)$
denote the set of functions $f \in C^{\infty} \lr{\ov{\Omega}
\setminus \mathbf{0}}$ such that
\begin{equation}
\chi_R (r)  f (r, \vartheta) = \sum\limits_{k=0}^{\mathcal{K}}
\sum\limits_{l=1}^{d_k} \chi_R (r) f_{k, l} (r) Y_{k, l}
(\vartheta), \label{7.1.1}
\end{equation}
    where the functions $f_{k,l}$ are such that $r^{-k} f_{k, l} \chi_R \in
\mathring{C}^{\infty}_{\frac{n}{2} +k -1} (E_{+}^1)$ and the
positive integer $\mathcal{K}$ depends on the function $f.$
    Recall (see Chap. \ref{ch2}) that the set
$\mathring{C}^{\infty}_{\nu} (E_{+}^1)$ consists of functions
$h(r)$ of the kind $h = P_{\nu} g,$ where $ g \in
\mathring{C}^{\infty}_{\nu} (\ov{E_{+}^1})$

For each $s$ and $R,$  $s \geq 0,$ $0 <R < R_0,$ and each function
$\chi \in X,$ define the following norms on $\mathring{T}^{\infty}
(\Omega_0)$:
\begin{equation}
 \| f \|_{ s, R} =
 \left(  \sum\limits_{k=0}^{\mathcal{K}} \sum\limits_{l=1}^{d_k}  \|r^{-k}
  \chi_R  f_{k, l} \|^2_{\mathring{H}^{s}_{\frac{n}{2} +k -1} (0, 2R)} +
   \|(1 - \chi_R) f \|^2_{H^s (\Omega)} \right)^{1/2}.
\label{7.1.2}
\end{equation}
 Spaces $\mathring{H}^s_{\nu} (0, R)$ introduced in Sec. \ref{sec4.4} are used.
 It follows from results of Chap. \ref{ch2} that the relation
    \begin{equation}
 \| f \|_{ s, R}^2 = \| \mathfrak{G}_n \Delta^{\frac{s}{2} }
   (\chi_R f) \|^2_{L_2 (U_{2R})} + \|  (1-\chi_R) f \|^2_{H^s (\Omega)}.
\label{713}
\end{equation}
    holds for even $s$.

If $R \in (0, R_0)$ and a nonnegative $s$ is fixed, then system of
norms \eqref{7.1.2} defines a topology on the lineal
$\mathring{T}^{\infty} (\Omega_0)$.

Introduce the space $H^s_{loc} (\Omega_0)$ consisting of functions
$f$ such that   $ (1-\chi_R) f \in H^s (\Omega)$ for each $R \in
(0, R_0)$.
    Endow this space with the topology defined by the family of
    seminorms
$$
p_{s, R} (f) = \|  (1-\chi_R) f \|^2_{H^s (\Omega)}, \  0 <R < R_0.
$$
  The set
$H^s_{loc} (\Omega_0)$ endowed with this topology is a complete
topological vector space.
    It is obvious that $\mathring{T}^{\infty}
(\Omega_0) \subset H^s_{loc} (\Omega_0)$ and the estimate
$$
 \| f \|_{ s, R} \geq p_{s, R} (f)
$$
 holds for each $f \in \mathring{T}^{\infty}$

 Thus, this is a topological embedding as well.

Define the space $M^s (\Omega_0)$ as the closure of the space
$\mathring{T}^{\infty} (\Omega_0)$ with respect to the topology
generated by the system of norms \eqref{7.1.2}.
    Therefore, elements of the space $M^s (\Omega_0)$ are ordinary functions belonging
    to the class $H^s_{loc} (\Omega_0).$

To find the relation between the spaces $M^s (\Omega_0)$ and $H^s
(\Omega),$ introduce the notation $\mathring{T}^{\infty}_{+}
(\Omega)$ for the subset of functions $f$ from
$\mathring{T}^{\infty} (\Omega_0)$ such that their functions
$f_{k, l}$ from expansion \eqref{7.1.1} are such that  $r^{-k}
f_{k, l} \chi_R \in \mathring{C}^{\infty}_{+} (\ov{E_{+}^1}).$

 Let $f \in \mathring{T}^{\infty}_{+} (\Omega).$
 Then $\chi_R f \in \mathring{T}^{\infty}_{+}(U_{2 R})$ for each $R \in (0, R_0)$.
  From Lemma \ref{lem5.2} and results of Sec. \ref{sec4.4}, it
  follows that
$$
 \| \chi_R f \|_{H^s(U_{2, R})} = \sum\limits_{k=0}^{\mathcal{K}} \sum\limits_{l=1}^{d_k}
   \|r^{-k} f_{k, l} \chi_R \|^2_{\mathring{H}^{s}_{\frac{n}{2} +k -1, +} (0, 2R)}  \geq
  \sum\limits_{k=0}^{\mathcal{K}} \sum\limits_{l=1}^{d_k}  \|r^{-k} f_{k, l} \chi_R
  \|^2_{\mathring{H}^{s}_{\frac{n}{2} +k -1} (0, 2R)}.
$$
    Hence,
$$
\| f \|^2_{s, R} \leq c \lr{ \| \chi_R f
\|_{\mathring{H}^s(U_{2R})}^2 +  \| (1-\chi_R) f
\|_{H^s(\Omega)}^2}.
$$
  For each $R \in (0, R_0)$, the right-hand
    side is the second power of an equivalent norm of the space $H^s(\Omega).$
  Due to Lemma \ref{lem:2.1.1}, this and the fact that the set
  $\mathring{T}^{\infty}_{+} (\Omega)$ is everywhere
    dense in $H^s(\Omega)$ implies the following assertion.

\begin{theorem} \label{teo: 7.1.1}
 For any nonnegative $s,$ the embedding
$$
H^s(\Omega) \subset  M^s (\Omega_0) \subset H^s_{loc} (\Omega_0)
$$
    holds in the topological sense.
    If  $n$ is odd, then the topology induced be the left-hand
    embedding and the proper topology of the space $H^s(\Omega)$ are equivalent to each other.
  For each $R \in (0,
R_0)$, the topologies of the subspace  $\widetilde{H}^s (\Omega
\subset \ov{U}_R)$ of the space $H^s(\Omega),$ consisting of
functions from $H^s(\Omega),$ vanishing in the ball $U_R,$
 induced by the embeddings $\widetilde{H}^s (\Omega \setminus
\ov{U}_R) \subset H^s(\Omega)$ and $\widetilde{H}^s (\Omega
\setminus \ov{U}_R) \subset M^s(\Omega_0),$ are equivalent to each
other.
\end{theorem}
\begin{corollary} \label{cor: 7.1.1}
 The space $H^s(\Omega)$ is not everywhere dense in  $M^s(\Omega_0).$
\end{corollary}

Let $\chi$ and $\chi'$ be functions from the set $X.$ Let us show
that the topologies generated by them are equivalent to each
other.

 Let $f \in \mathring{T}^{\infty} (\Omega_0)$ and $0<R'<R<R_0.$
    Then
\begin{multline}
\sum\limits_{k=0}^{\mathcal{K}} \sum\limits_{l=1}^{d_k}
\|\chi_R r^{-k}   f_{k, l} \|^2_{\mathring{H}^{s}_{\frac{n}{2} +k -1} (0, 2R)}  \\
\leq 2 \sum\limits_{k=0}^{\infty} \sum\limits_{l=1}^{d_k}
\|\chi'_{R'} r^{-k}   f_{k, l} \|^2_{\mathring{H}^{s}_{\frac{n}{2}
+k -1} (0, 2R)} + 2 \sum\limits_{k=0}^{\infty}
\sum\limits_{l=1}^{d_k} \|(\chi_{R} - \chi'_{R'}) r^{-k}   f_{k,
l} \|^2_{\mathring{H}^{s}_{\frac{n}{2} +k -1} (0, 2R)}.
\label{7.1.3}
\end{multline}
    Since $\chi$ and $\chi'$ belong to the set $X,$ it follows that $\chi_R (r)
= \chi'_{R'} (r) = 1$ provided that $r< R'.$  Then   ${\chi_R (r)
- \chi'_{R'} (r) = 0}$ for  $0 \leq r \leq R'$ and for $r \geq 2
R.$ Hence, $(\chi_R (r) - \chi'_{R'} (r)) r^{-k} f_{k, l} \in
\mathring{C}_{+}^{\infty} [0, 2R)$.
 Then, due to Theorem \ref{teo:2.1.2}, the second sum of the
 right-hand side \eqref{7.1.3} is estimated by the norm
$$
\| (\chi_R - \chi'_{R'}) f \|_{\mathring{H}^s (U_{2 R_0})}  \leq c
\, \| (\chi_R (r) - \chi'_{R'}) f \|_{{H}^s (\Omega)}.
$$
    This yields the estimate
\begin{equation}
\| f \|_{ s, R}^2 \leq c  \left(  \sum\limits_{k=0}^{\infty}
\sum\limits_{l=1}^{d_k}  \|\chi'_{R'} r^{-k}   f_{k, l}
\|^2_{\mathring{H}^{s}_{\frac{n}{2} +k -1} (0, 2R')}  + \|(1 -
\chi_R) f \|^2_{H^s (\Omega)} +  \|(1 - \chi'_{R'}) f \|^2_{H^s
(\Omega)} \right). \label{7.1.4}
\end{equation}
    From the monotonicity of the functions $\chi$ and $\chi'$ and the inequality $R' < R,$
it follows that $1-  \chi'_{R'} \geq \delta > 0$ on the support of
the function $1 - \chi_R.$ Hence,  $(1 - \chi_R) / (1- \chi'_{R'})
\in C^{\infty} (\ov{\Omega}).$ Then the following estimate holds:
$$
\|(1 - \chi_R) f \|_{H^s (\Omega)}  \leq c \,  \|(1 - \chi'_{R'}) f \|_{H^s (\Omega)}.
$$
  Substituting this relation in \eqref{7.1.4}, we obtain the final
  estimate
$$
\| f \|_{ s, R}^2 \leq c  \left(  \sum\limits_{k=0}^{\infty}
\sum\limits_{l=1}^{d_k}  \|\chi'_{R'} r^{-k}   f_{k, l}
\|^2_{\mathring{H}^{s}_{\frac{n}{2} +k -1} (0, 2R')}  +  \|(1 -
\chi'_{R'}) f \|^2_{H^s (\Omega)} \right) =  c \, \| f \|_{ s,
R'}^2,
$$
    where the
 left-hand side is the norm generated by the function $\chi$ and the right-hand
 side contains the norm generated by the function $\chi'.$ Changing the places of $\chi$
 and
$\chi',$ we obtain the opposite estimate. Hence, the following
assertion is proved.

\begin{lemma}\label{lem: 7.1.1}
Topologies generated in the space $M^s (\Omega_0)$ by different
functions from the set $X,$ are equivalent to each other.
\end{lemma}
\begin{lemma}\label{lem: 7.1.2}
 The space $M^s (\Omega_0)$ is a complete denumerably normable
  topological space, i.\,e., a Frechet space.
\end{lemma}

\begin{proof}
    Consider the  denumerable set $\| \|_{s, R_0/m},$ $m=1, 2,\dots,$ of norms.
  Take an arbitrary $R \in (0, R_0).$ Then there exist nonnegative integers $m_1$ and $m_2$
  such that $\dfrac{R_0}{m_1}<R< \dfrac{R_0}{m_2}.$
  In this case, due to Lemma \ref{lem: 7.1.1}, there exist positive constants $c_1$ and $c_2$
  such that the estimate
$$
c_2 \|f \|_{s, R_0/ m_2} \leq \| f\|_{s, R} \leq  c_1 \|f \|_{s,R_0/ m}
$$
 holds for each $f \in M^s (\Omega_0)$.
\end{proof}

\begin{theorem}\label{teo: 7.1.2}
 Let $0 \leq s_1 < s_2.$ Then the space  $M^{s_1} (\Omega_0)$
 is continuously embedded in $M^{s_2} (\Omega_0).$
\end{theorem}

\begin{proof}
    Since
$$
\| g \|_{\mathring{H}^{s_1} (0, R)} = \| D^{s_1} g\|_{L_2 (0, R)}
\leq c \, (s_1 - s_2, R) \| D^{s_2} g\|_{L_2 (0, R)}=
  c\, (s_1 - s_2, R) \|  g\|_{\mathring{H}^{s_2} (0, R)}
$$
    for each function $g \in \mathring{C}^{\infty} [0, R),$ it
    follows that
$$
    \begin{gathered}
\|\chi_{R} \, r^{-k}   f_{k, l}
\|^2_{\mathring{H}^{s_1}_{\frac{n}{2} +k -1} (0, 2R)} =
\frac{\sqrt{\pi}}{2^{\frac{n}{2} +k -1} \, \Gamma \lr{\frac{n}{2}
+k}}  \| S_{\frac{n}{2} +k -1} (\chi_{R} \, r^{-k}   f_{k, l})
\|_{\mathring{H}^{s_1} (0, R)}
    \\
\leq c \, (s_1 - s_2, R) \|\chi_{R} \, r^{-k}   f_{k, l}
\|^2_{\mathring{H}^{s_2}_{\frac{n}{2} +k -1} (0, 2R)}.
    \end{gathered}
$$
    Summing these inequalities over  $k$ and $l$ and taking onto
    account that
$$
\| (1-\chi_R) f \|_{H^{s_1} (\Omega)} \leq  \| (1-\chi_R) f
\|_{H^{s_2} (\Omega)},
$$
    we obtain the inequality $ \|f\|_{s_1, R} \leq c \, \|f\|_{s_2, R},$
    where the constant does not depend on the function $f \in M^{s_2}
(\Omega_0).$
\end{proof}
 Investigate traces of functions from $M^s (\Omega_0)$ on the boundary of the domain
  $\Omega_0,$ consisting of the point $\mathbf{0}$ and the surface $\pr \Omega.$

Note that since the space $M^s$ and the space  $H^s$ have  the
same structure outside each neighborhood of the point
$\mathbf{0}$, we have the following assertion.
\begin{theorem}\label{teo: 7.1.3}
    Let $s > \dfrac{1}{2}.$ Then the map $f \to f|_{\pr
\Omega}$ continuously acts from the space $M^s (\Omega_0)$ to
$H^{s- \frac{1}{2}} (\pr \Omega)$ and there exists an operator
 continuously mapping $H^{s- \frac{1}{2}}
(\pr \Omega)$ into $M^s (\Omega_0)$ and such that $f|_{\pr
\Omega}=g,$ where $f$ is the image of the function  $g$ under this
map.
\end{theorem}

\subsection{Direct and inverse problems on
 $\sigma$-traces ($K$-traces)}\label{sec14.2}

    Define the trace at the boundary point $\mathbf{0} \in \pr\Omega_0.$
    It is easy to see that, in general, functions from $M^s(\Omega_0)$
    have no classical traces at $\mathbf{0}$ even if $s$ is arbitrarily large.
    We introduce a new (in a way, nonlocal) notion of the trace.

    Let $A(\Theta)$ be the set of functions  $\psi\lr{\vartheta}$
    defined on the sphere $\Theta$, real-analytic
    on $\Theta$, and such that the norms
\begin{equation}
\| \psi \|_h = \lr{  \sum\limits_{k=0}^{\infty}
\sum\limits_{l=1}^{d_k} \|\psi_{k, l} \|^2 h^{-2 k}
}^{\frac{1}{2}}, \label{7.2.1}
\end{equation}
    where $\psi_{k, l}$ denote the coefficients of the expansion of the function $\psi$
    with respect to spherical harmonics $Y_{k, l},$ are finite
for each $h \in (0, 1)$.

\begin{lemma}\label{lem: 7.2.1}
 The space $A(\Theta)$ is a complete denumerable normable
 topological space.
\end{lemma}

\begin{proof}
    To prove the completeness, assume that $\{ \psi^m \} \subset A(\Theta)$ is a fundamental
    sequence. Then, for each $h \in (0,1)$, we have the limit relation
     $\| \psi^{m_1} - \psi^{m_2} \|_h \to 0$ as $m_1, m_2 \to\infty.$
  Therefore, for each $h$ there exists a function $g_h(\vartheta)$ such that
   $\| g_h \|_h < \infty$ and $\| \psi^m - g_h \|_h\to 0$ as $m \to \infty.$
   From \cite[p. 496]{77}, it follows that the function $g_h$ is analytic on $\Theta.$
   Since $\| \psi \|_h \leq\| \psi \|_{h'}$ provided that $h' < h,$
   it follows from the limit relation $\|g_{h'}- \psi^m\|_{h'}  \to 0$
   that $\|g_{h'}- \psi^m \|_{h}  \to 0$ as $m \to \infty.$ Hence, all functions $g_h$
   coincide each other. Obviously, this unique function belongs to
   the space $A(\Theta)$ and is the limit of the sequence $\{ \psi^m \}$ with respect
   to each norm \eqref{7.2.1}, which completes the proof of the completeness.

    Consider the denumerable set of norms $\| \|_{1/m},$ $m=2,3,\dots$
    It is easy to see that they define a topology equivalent to the original one.
    Thus, the space $A(\Theta)$ is denumerably normable.
\end{proof}
    Note that the space $A(\Theta)$ consists of functions such that they are real-analytic
    on the sphere $\Theta$ and admit a harmonic extension to the whole space $E^n.$
    For each $\psi \in A(\Theta)$, this extension is defined by the
    relation
$$
\Psi (r, \vartheta) = \sum\limits_{k=0}^{\infty}
\sum\limits_{l=1}^{d_k} r^k \psi_{k, l} Y_{k, l} (\vartheta).
$$
    It is clear that the function $\Psi$ is harmonic in the whole space $E^n.$

    Let $f \in \mathring{T}^{\infty} (\Omega_0).$ This means that the
    following relation holds provided that $r < 2 R_0$:
\begin{equation}
f (r, \vartheta) = \sum\limits_{k=0}^{\mathcal{K}}
\sum\limits_{l=1}^{d_k} f_{k, l} (r) Y_{k, l} (\vartheta),
\label{7.2.2}
\end{equation}
    where $r^{-k} \chi_R f_{k, l} \in
\mathring{C}^{\infty}_{\frac{n}{2} +k -1} (0, 2R),$ $0<R<R_0.$
    Note that the functions $f_{k, l} (r)$ might have singularities at the point $r=0.$
  For $r < R_0$ and $n \geq 3$, define the averaging operation $\sigma$
  on the function class $\mathring{T}^{\infty}(\Omega_0)$ as follows:
\begin{equation}
\sigma f (r, \vartheta) = \sum\limits_{k=0}^{\mathcal{K}}
\sum\limits_{l=1}^{d_k} r^{n+k-2} f_{k, l} (r) Y_{k, l}
(\vartheta). \label{7.2.3}
\end{equation}
  The {\it $\sigma$-trace}
  of a function $f \in\mathring{T}^{\infty} (\Omega_0)$ is the limit
$$
\lim\limits_{r \to + 0} \sigma f (r, \vartheta) =: \sigma f|_0.
$$
    From the assumptions imposed on the function $f$, it follows from Theorem \ref{teo:1.4.1}
    that its $\sigma$-trace
    exists and is equal to the function
$$
\psi (\vartheta) = \sum\limits_{k=0}^{\mathcal{K}}
\sum\limits_{l=1}^{d_k} \psi_{k, l} (r) Y_{k, l} (\vartheta),
$$
    where
$$
\psi_{k, l} = \lim\limits_{r \to + 0}  r^{n+k-2} f_{k, l} (r).
$$
    Also, it is clear that $\psi \in A(\Theta).$

 Now, we provide an equivalent definition  of
$\sigma$-traces; it does not use the expansion with respect to
spherical harmonics. Let $n \geq 3,$ and the functions $f_{k, l}$
 from \eqref{7.2.2} be defined by the relation
$$
f_{k, l} (r) = \int\limits_{\Theta} f (r, \vartheta) Y_{k, l} (\vartheta) \,  d \vartheta.
$$
    Substituting this expression in \eqref{7.2.3}, we find that
$$
\sigma f (r, \vartheta) = \sum\limits_{k=0}^{\infty}
\sum\limits_{l=1}^{d_k} r^{n+k-2} Y_{k, l} (\vartheta)
\int\limits_{\Theta} f (r, \vartheta') Y_{k, l} (\vartheta') \,  d
\vartheta' =
 \int\limits_{\Theta} f (r, \vartheta')  \sum\limits_{k=0}^{\infty}
  \sum\limits_{l=1}^{d_k} r^{n+k-2} Y_{k, l} (\vartheta)  Y_{k, l} (\vartheta') \,
    d \vartheta'.
$$
    It is known (see, e.\,g., \cite{BE2}), that the relation
$$
\sum\limits_{k=0}^{\infty} \sum\limits_{l=1}^{d_k} r^{k} Y_{k, l}
(\vartheta)  Y_{k, l} (\vartheta') = \frac{\Gamma
\lr{\frac{n}{2}}}{2 \pi^{\frac{n}{2}}} \frac{1-r^2}{|r \vartheta -
\vartheta' |^n}=K_n (r \vartheta, \vartheta')
$$
    holds provided that $r < 1$.
 Here, the function $K_n (x, y)$ is called the {\it Poisson
 kernel} for the unit sphere  $\Theta.$ Therefore,
$$
\sigma f|_0 =  \lim\limits_{r \to + 0}  r^{n-2}
\int\limits_{\Theta} f (r, \vartheta')  K_n (r \vartheta,
\vartheta') \, d \vartheta'.
$$
    This is the desired explicit definition of the {\it $\sigma$-trace}.

Define the  $\sigma$-trace
 for the dimension $n=2$. Let $r$ and $\varphi,$ $r>0,$ $|\varphi|<\pi,$
 be the polar coordinates on the plane.
 In this case, instead of \eqref{7.2.3}, we assign
$$
\sigma f (r, \varphi) = \frac{f_0 (r)}{\ln r} +
\sum\limits_{k=1}^{\mathcal{K}} r^k \lr{f_{k, 1} (r) \cos (k
\varphi)+ f_{k, 2} (r) \sin (k \varphi)},
$$
    where
\begin{align*}
f_0 (r) &= \frac{1}{2 \pi} \int\limits_{- \pi}^{\pi}  f (r, \varphi)  \, d \varphi,
\\
f_{k,1} (r) &= \frac{1}{\pi} \int\limits_{- \pi}^{\pi}  f (r,
\varphi) \cos(k \varphi) \, d  \varphi,
\\
f_{k,2} (r) &= \frac{1}{2 \pi} \int\limits_{- \pi}^{\pi}  f (r,
\varphi) \sin (k \varphi)  \, d \varphi.
\end{align*}
    This yields the relation
$$
\sigma f (r, \varphi) = \frac{f_0 (r)}{\ln r} +  \frac{1}{\pi}
\int\limits_{- \pi}^{\pi} f (r, \varphi')
\sum\limits_{k=1}^{\infty} r^k  \cos (k( \varphi-\varphi')) \, d
\varphi'.
$$
    Since the relation
$$
\sum\limits_{k=1}^{\infty} r^k  \cos (k( \varphi-\varphi')) =
\frac{r \cos( \varphi-\varphi') -r^2 }{1-2 r \cos(
\varphi-\varphi') +r^2}
$$
 holds provided that $r<1$, it follows that the $\sigma$-trace
 for $n=2$ can be defined by the relation
$$
\sigma f|_0 =  \lim\limits_{r \to + 0}  \frac{1}{2 \pi}
\int\limits_{- \pi}^{\pi} f (r, \varphi') \lr{ \frac{2 r \cos(
\varphi-\varphi') - 2 r^2 }{1-2 r \cos( \varphi-\varphi') +r^2} +
\frac{1}{\ln r}}  \, d \varphi'
$$
    as well.

Extend the  $\sigma$-trace
    notion to functions $f \in M^s (\Omega_0).$
    For each function $f \in M^s (\Omega_0)$  there exists a sequence of functions
     $f^j \in \mathring{T}^{\infty} (\Omega_0),$ $j = 1, 2,
\dots,$ converging to $f$ as  $j \to \infty$ with respect to the
topology of the space $M^s (\Omega_0).$
 As we proved above, the $\sigma$-trace
 of each function $f^j$ exists and belongs to the space $A(\Theta).$
 If the sequence $\{\psi^j\} := \{\sigma
f^j|_0\}$ converges as $j \to \infty$ to a function $\psi \in
A(\Theta)$ with respect to the topology of the space $A(\Theta)$
and the function $\psi$ does not depend on the choice of the
sequence $\{f^j\},$ then it is called the \emph{$\sigma$-trace}
 of the function $f.$

The following direct $\sigma$-trace
 theorem holds.

\begin{theorem}\label{teo: 7.2.1}
 Let $s \geq 1$ provided that  $n \geq 3$ and $s \geq 2$ provided that   $n=2.$
 Then each $f \in M^s (\Omega_0)$ has a $\sigma$-trace
      $\sigma f|_0 \in A(\Theta)$ and the operator $f \to \sigma f|_0$ continuously maps
      the space  $ M^s (\Omega_0)$ into $A(\Theta).$
\end{theorem}

\begin{proof}
    It suffices to show that the specified operator
     continuously maps
      the space   $\mathring{T}^{\infty} (\Omega_0)$ with the topology induced by the space
       $M^s (\Omega_0)$ into the space $A(\Theta).$
 To do that, one has to prove that for each $h \in (0, 1)$ there exist  a number
  $R \in (0, R_0)$ and a positive constant $c$ such that each function
  $f \in \mathring{T}^{\infty} (\Omega_0)$ satisfies the estimate
\begin{equation}
\| \sigma f|_0 \|_h \leq c \, \| f \|_{s, R}. \label{7.2.4}
\end{equation}
    Let $f_{k, l}$ be the expansion coefficients of functions $f \in
\mathring{T}^{\infty} (\Omega_0)$ with respect to the spherical
harmonics $Y_{k, l}.$ Then the functions $r^{-k} \chi_R f_{k, l}$
 belong to the space $\mathring{C}^{\infty}_{\frac{n}{2}+k-1}(E^1_{+})$
 and, therefore, to the space $\mathring{H}^{s}_{\frac{n}{2}+k-1} (0, 2R).$

Due to Corollary \ref{cor:1.4.3}, each function $g \in
\mathring{H}^{s}_{\nu} (0, 2R)$ satisfies the estimate
$$
|\left. \sigma_{\nu} (r) g(r)\right|_{r=0} | \leq c \, (s, R)
(4R)^{\nu} (\nu+1)^{-s} \| g\|_{\mathring{H}^{s}_{\nu} (0, 2R)}
$$
    provided that $s \geq 1,$ $\nu \geq 0,$ and $s+\nu > 1.$
    Here, $\sigma_{\nu}(r) = r^{2 \nu}$ provided that $\nu > 0,$
    and $\sigma_0 (r) = \Big(\ln\dfrac{1}{r}\Big)^{-1}.$
    Assigning  $\nu =\dfrac{n}{2}+k-1,$ $ g = \chi_R r^{-k} f_{k, l}$
 in the last estimate, we obtain that
$$
|\left. \sigma_{\frac{n}{2}+k-1} (r)  r^{-k} f_{k, l}
(r)\right|_{r=0} | \leq c \, (s, n,  R) (4R)^{k} (k+1)^{-s}  \| \chi_R \, r^{-k} f_{k,
l}\|_{\mathring{H}^{s}_{\frac{n}{2}+k-1} (0, 2R)}.
$$
    Further, since
$$
\sigma f |_0 = \sum\limits_k \sum\limits_l
\sigma_{\frac{n}{2}+k-1} r^{-k} f_{k, l} (r)|_{r=0} Y_{k, l}
(\vartheta),
$$
    it follows from the previous estimate that the following inequalities hold
    provided that $h \in (0, 1)$:
$$
    \begin{gathered}
\| \sigma f |_0 \|_h^2 = \sum\limits_k \sum\limits_l |
\sigma_{\frac{n}{2}+k-1} r^{-k} f_{k, l} (r)|_{r=0} |^2 h^{-2 k}
 \\
\leq c \, (s, R, n) \sum\limits_k \sum\limits_l (4 R)^{2 k}
(k+1)^{-2 s} h^{-2 k}  \|  \chi_R \, r^{-k} f_{k, l}
\|^2_{\mathring{H}^{s}_{\frac{n}{2}+k-1} (0, 2R)}
    \\
\leq c \, (s, R, n) \sum\limits_k \sum\limits_l  \|  \chi_R \,
r^{-k} f_{k, l} \|^2_{\mathring{H}^{s}_{\frac{n}{2}+k-1} (0, 2R)},
    \end{gathered}
$$
    where $R = \dfrac{1}{4} h.$
\end{proof}
    Consider the inverse $\sigma$-trace
 theorem.

\begin{theorem}\label{teo: 7.2.2}
 For each function $\psi \in A (\Theta)$ there exists a function $f$
  such that it belongs to the spaces $M^s (\Omega_0)$ for all nonnegative $s$
  and $\psi$ is its  $\sigma$-trace.
 For each $R \in (0, R_0)$ there exists $h \in (0, 1)$ such that
\begin{equation}
\| f \|_{s, R} \leq c \, \| \psi \|_h, \label{7.2.5}
\end{equation}
    where  $c$ is a positive constant independent of the choice of the function $\psi.$
\end{theorem}

\begin{proof}
    Let $\psi \in A (\Theta)$ and $\psi_{k, l}$ be the
    coefficients of its expansion in the series with respect to
    spherical harmonics.
    If  $n \geq 3,$ then define the claimed function $f$ by the
    relation
\begin{equation}
f (r, \vartheta) = \sum\limits_{k=0}^{\infty}
\sum\limits_{l=1}^{d_k} r^{2-n - k}  \psi_{k, l}  Y_{k, l}
(\vartheta). \label{7.2.6}
\end{equation}
    If $n=2,$ then assign
\begin{equation}
f (r, \vartheta) = \psi_0 \ln r + \sum\limits_{k=1}^{\infty}
\sum\limits_{l=1}^{2} r^{ - k}  \psi_{k, l}  Y_{k, l} (\vartheta).
\label{7.2.7}
\end{equation}
    Since $\psi \in A (\Theta),$ it follows that series \eqref{7.2.6} and \eqref{7.2.7}
     converge
 to an analytic function   absolutely and uniformly outside each
 ball centered at $\mathbf{0}$.
 Moreover, the function $f$ is harmonic in $E^n$ apart from the point $\mathbf{0}.$

    Let us prove estimate \eqref{7.2.5}. First, we assume that $s$ is even.
    Then, using properties of transmutation operators (see Chap \ref{ch2}),
    we obtain that
\begin{multline}
\| \chi_R \, r^{2-n-2k} \|^2_{\mathring{H}^{s}_{\frac{n}{2}+k-1}
(0, 2R)} \leq \frac{2 \pi}{ 2^{n+2k-1}\, \Gamma^2
\lr{\frac{n}{2}+k} }
\| D^s S_{\frac{n}{2}+k-1} (\chi_R \, r^{2-n-2k}) \|^2_{L_2 (0, 2R)}  \\
= \frac{2 \pi}{ 2^{n+2k-1}\, \Gamma^2 \lr{\frac{n}{2}+k} } \|
S_{\frac{n}{2}+k-1} B^{\frac{s}{2}}_{\frac{n}{2}+k-1} (\chi_R
r^{2-n-2k}) \|^2_{L_2 (0, 2R)}, \label{7.2.8}
\end{multline}
    where $S_{\nu}$ are the transmutation operators
 from Sec. \ref{sec4.1}.

Applying the Darboux--Weinstein
 relation to \eqref{7.2.8}, we obtain the inequality
\begin{equation}
\| \chi_R \, r^{2-n-2k} \|^2_{\mathring{H}^{s}_{\frac{n}{2}+k-1}
(0, 2R)} \leq \frac{2 \pi}{ 2^{n+2k-1}\, \Gamma^2
\lr{\frac{n}{2}+k} }
 \| S_{\frac{n}{2}+k-1} (r^{2-n-2k}
B^{\frac{s}{2}}_{\frac{n}{2}+k-1} \chi_R) \|^2_{L_2 (0, 2R)}.
\label{7.2.9}
\end{equation}
    Since $\chi_R (r) =1$ provided that $0 \leq r \leq R,$ it
    follows that
$$
B_{- \nu}  \chi_R (r) = r^{2 \nu -1} \frac{\pr }{\pr r} \lr{r^{1 -
2 \nu}  \frac{\pr \chi_R}{\pr r}} =0
$$
 for the same values of $r.$
 Hence, $B_{- \nu}^{\frac{s}{2}} \chi_R(r) = 0$ provided that $r \leq R.$
 The last relation holds for $r \geq 2R$ as well because
$\chi_R (r) =0$ for $r \geq 2R.$
    This implies that the function
$B^{\frac{s}{2}}_{\frac{n}{2}+k-1} \chi_R$
    is infinitely differentiable, has a compact support, and vanishes near the origin.
 Hence, $B^{\frac{s}{2}}_{\frac{n}{2}+k-1}\chi_R \in \mathring{C}_{+}^{\infty}(\ov{E_{+}^1}).$
 Then, by virtue of Lemma \ref{lem:1.4.1}, from relation
 \eqref{7.2.9}, we conclude that
\begin{equation}
\| r^{2-n-2k}  \chi_R \|^2_{\mathring{H}^{s}_{\frac{n}{2}+k-1} (0,
2R)} \leq 2 \int\limits_{R}^{2 R}
|B^{\frac{s}{2}}_{1-\frac{n}{2}-k} \chi_R (r) |^2  r^{3-n-2k} \,
dr. \label{7.2.10}
\end{equation}
    Since $\chi (t) \in C^{\infty} [1, 2],$ we have the estimate
$$
| B^{\frac{s}{2}}_{-\nu} \chi (t) | = \left| \lr{D^2_t + \frac{1-2
\nu}{t} D_t}^{\frac{s}{2}} \chi (t) \right| \leq c \, (s) (1+
\nu)^{\frac{s}{2}},
$$
 where $c \, (s)>0$ and  $c \, (s)$ does not depend on $\nu.$
 This and \eqref{7.2.10} imply the inequality
\begin{multline}
\| r^{2-n-2k}  \chi_R \|^2_{\mathring{H}^{s}_{\frac{n}{2}+k-1} (0,
2R)} \leq
    \\
    2 R^{4 - s- n-2 k}
  \int\limits_{1}^{2} |B^{\frac{s}{2}}_{1-\frac{n}{2}-k}
\chi (t) |^2 \,  t^{3-n-2k} \, d t
\leq c (s, n, R) \, R^{- 2 k} (1+k)^{s-1}. \label{7.2.11}
\end{multline}
    Thus, we have the estimate
\begin{equation}
\sum\limits_k \sum\limits_l \| r^{2-n-2k} \psi_{k, l} \chi_R
\|^2_{\mathring{H}^{s}_{\frac{n}{2}+k-1} (0, 2R)}
 \leq c  \sum\limits_k \sum\limits_l | \psi_{k, l}|^2  R^{- 2 k} (1+k)^{s-1}.
\label{7.2.12}
\end{equation}
    Now, to estimate the expression $\| (1 - \chi_R) f \|_{H^s (\Omega)},$
 denote the diameter of the domain $\Omega$ by $\ov{R}$. Then
$\chi_{\ov{R}} (r) = 1$ on $\Omega$ and, therefore,
$$
\|  (1 - \chi_R) f \|_{H^s (\Omega)}^2 =  \| \chi_{\ov{R}} (1 -
\chi_R) f \|_{H^s (\Omega)} \leq
  \| \chi_{\ov{R}} (1 - \chi_R) f \|_{\mathring{H}^s (U_{2
\ov{R}})},
$$
    where $U_{2 \ov{R}}$ is the ball of radius ${2\ov{R}}$ centered at $\mathsf{0}$.
  Theorem \ref{teo:2.1.2} implies the relation
$$
    \begin{gathered}
 \| \chi_{\ov{R}} (1 - \chi_R) f \|_{\mathring{H}^s (U_{2 \ov{R}})}^2
  = \sum\limits_{k=0}^{\infty} \sum\limits_{l=1}^{d_k}  \| \chi_{\ov{R}} (1 - \chi_R)
   f_{k, l} \|_{\mathring{H}^{s}_{\frac{n}{2}+k-1, +} (0, 2 \ov{R})}^2
   \\
=  \sum\limits_{k=0}^{\infty} \sum\limits_{l=1}^{d_k} |\psi_{k,
l}|^2  \| \chi_{\ov{R}} (1 - \chi_R) r^{2-n-2k}
\|_{\mathring{H}^{s}_{\frac{n}{2}+k-1, +} (0, 2 \ov{R})}^2
    \\
=  \sum\limits_{k=0}^{\infty} \sum\limits_{l=1}^{d_k} |\psi_{k,
l}|^2 \int\limits_R^{2 \ov{R}} |B^{\frac{s}{2}}_{\frac{n}{2}+k-1}
(\chi_{\ov{R}} (1 - \chi_R) r^{2-n-2k}) |^2   r^{2k+n-1} \, dr
    \\
=  \sum\limits_{k=0}^{\infty} \sum\limits_{l=1}^{d_k} |\psi_{k,
l}|^2 \int\limits_R^{2 \ov{R}} |B^{\frac{s}{2}}_{1-\frac{n}{2}-k}
(\chi_{\ov{R}} (1 - \chi_R)) |^2   r^{3-n-2k} \, dr,
    \end{gathered}
$$
    (here, the Darboux--Weinstein
    relation is used again).

The following estimate is easily verified:
$$
\max\limits_{r} |B^{\frac{s}{2}}_{1-\frac{n}{2}-k}  (\chi_{\ov{R}}
(1 - \chi_R)) |  \leq c \, (1+k)^s,
$$
 where $c$ is a positive constant depending on $s,$ $n,$ $R,$ $\ov{R},$ and $\chi,$
 but independent of $k.$ Combining the last two relations, we
 obtain the inequality
\begin{multline}
 \| \chi_{\ov{R}} (1 - \chi_R) f \|_{\mathring{H}^s (U_{2 \ov{R}})}^2
  \leq c \sum\limits_{k=0}^{\infty} \sum\limits_{l=1}^{d_k} |\psi_{k, l}|^2 (1+k)^s
  \int\limits_R^{2 \ov{R}} r^{3-n-2k} \, dr
  \\
 \leq c' \, \sum\limits_{k=0}^{\infty} \sum\limits_{l=1}^{d_k}
|\psi_{k, l}|^2  (1+k)^{s-1}  R^{- 2 k}, \label{7.2.13}
\end{multline}
where $c'$ is a positive constant  independent of $\psi.$
    Now, for $h <\min\limits (1, R),$ relations \eqref{7.2.12}-\eqref{7.2.13}
    yield the estimate
$$
\| f \|^2_{s, R} =  \sum\limits_{k=0}^{\infty}
\sum\limits_{l=1}^{d_k} \| \chi_R\, r^{2-n-2k} \psi_{k, l}
\|^2_{\mathring{H}^{s}_{\frac{n}{2}+k-1} (0, 2R)} +
 \| (1 - \chi_R) f\|_{H^s (\Omega)}^2 \leq c \sum\limits_{k=0}^{\infty}
\sum\limits_{l=1}^{d_k} |\psi_{k, l}|^2 h^{-2 k} = c \, \| \psi
\|_h^2,
$$
 where $c$ is a positive constant independent of $\psi.$
 Thus, estimate \eqref{7.2.5} is proved for all even  $s.$
 Proving Theorem \ref{teo: 7.1.2}, we established the validity of the inequality
  $\| f \|_{s', R}  \leq c  \| f \|_{s, R}$
    under the assumption that $s'<s.$
 This implies the validity of estimate \eqref{7.2.5} for all nonnegative $s.$
\end{proof}
 To complete this section, we prove that the embedding of the space $M^s(\Omega_0)$ into
  $M^{s'} (\Omega_0)$ is not completely continuous  under the assumption that $s>s'$
  (unlike its continuity established by Theorem \ref{teo: 7.1.2}).

    Really, let $V$ be an arbitrary bounded set in the space $A(\Theta).$
    Due to the inverse $\sigma$-trace
    theorem, there exists a bounded domain $W$ belonging both to the space $M^s(\Omega_0)$
    and the space $M^{s'} (\Omega_0)$ such that the set of $\sigma$-traces
    of its functions coincides with $V.$
    If the specified operator of the embedding of the space  $M^s (\Omega_0)$
    into the space $M^{s'} (\Omega_0)$ is completely continuous, then  the set $W$
is relatively compact in the space $M^{s'} (\Omega_0).$ Then, due
to the direct $\sigma$-trace
    theorem, its image (i.\,e., the set of $\sigma$-traces
    of functions from $W$) is relatively compact in the space $A (\Theta).$
    Thus, each bounded set is relatively compact in the space $A (\Theta).$
   However, this is possible only in finite-dimensional
   spaces, while $A (\Theta)$ is not finite-dimensional.
 This contradiction proves the claimed assertion.

\section{New Boundary-Value
 Problems for the Poisson Equation}\label{sec15}

\subsection{Setting of boundary-value
 problems and isolated singular points of harmonic functions}\label{sec15.1}

In the domain $\Omega_0$, consider the Poisson equation
\begin{equation}
\Delta u= f(x), \  x \in \Omega_0 \label{8.1.1}
\end{equation}
 with the following boundary-value
 condition on $\pr \Omega_0$:
\begin{equation}
 u|_{\pr \Omega}  = g(x), \  x \in \pr \Omega,
\label{8.1.2}
\end{equation}
\begin{equation}
\sigma u|_{0} = \psi (\vartheta), \  \vartheta \in \Theta.
\label{8.1.3}
\end{equation}
    As above, $\sigma u |_0$ denotes the $\sigma$-trace
 of the function $u$ at the point $\mathsf{0} \in \pr \Omega_0.$ The
  boundary-value problem given by \eqref{8.1.1}--\eqref{8.1.3}
  generate an operator of the kind
$$
\Lambda: u \to  \Lambda u = \{ \Delta u, u|_{\pr \Omega}, \sigma u|_0 \}.
$$
 Equip the space
  $M^s =  M^s (\Omega_0) \times H^{s+\frac{3}{2}} (\pr \Omega) \times A (\Theta)$
  by the direct-product
  topology.

From the above results, it easily follows that the operator
$\Lambda$ continuously maps the space  $M^{s+2} (\Omega_0)$ into
$M^s$ for each nonnegative $s.$

The main result of the present chapter is as follows.

\begin{theorem}\label{teo: 8.1.1}
 The operator $\Lambda$ has an inverse operator $\Lambda^{-1}$  continuously mapping
  the space   $M^s$ onto the space $M^{s+2} (\Omega_0)$ for each even nonnegative $s.$
\end{theorem}

The proof of the theorem consists of several stages.

First, consider the following property of isolated singular points
of harmonic functions; in our opinion, it represents an
independent interest as well. It turns out that the $\sigma$-trace
    completely characterizes the behavior of the singular part of
 the harmonic function in a neighborhood of the isolated singular point.
 More exactly, this means the following.

\begin{lemma}\label{lem: 8.1.1}
    If a function $u(x)$ is harmonic in $\Omega_0,$ then it has a
$\sigma$-trace at the point $\mathsf{0},$  belonging to the space
$A(\Theta).$ Conversely, for each function $\psi \in A(\Theta)$
there exists a unique {\rm (}up to a function harmonic in a
neighborhood of $\mathsf{0})$ function $u(x)$ harmonic in
$\Omega_0$ and such that $\sigma u |_0 = \psi.$
\end{lemma}

\begin{proof}
Let  $u(x)$ be a function harmonic in the domain $\Omega_0$.
    Then the following expansion takes place for $n \geq 3$
(see, e.\,g., \cite{77}):
\begin{equation}
u (r, \vartheta) = \sum\limits_{k=0}^{\infty}
\sum\limits_{l=1}^{d_k} a_{k, l} r^k Y_{k, l} (\vartheta) +
\sum\limits_{k=0}^{\infty} \sum\limits_{l=1}^{d_k} b_{k, l}
r^{2-n-k} Y_{k, l} (\vartheta) =u_1+ u_2, \label{8.1.4}
\end{equation}
    where the former series converges absolutely and uniformly in
    each ball centered at the point $\mathsf{0}$ and located in the domain $\Omega$
    together with its closure, while the latter series
     converges absolutely and uniformly outside each ball centered at the point
$\mathsf{0}.$
 Then it is obvious that the  $\sigma$-trace
 of the function $u_1$ is equal to zero.
 If $r> 0 $, then the function  $u_2$ satisfies the relation
$$
\sigma u_2 (r, \vartheta) = r^{n-2} \int\limits_{\Theta}    u_2
(r, \vartheta') \mathcal{K}_n  u_2 (r \vartheta, \vartheta') \, d
\vartheta'=
 \sum\limits_{k=0}^{\infty} \sum\limits_{l=1}^{d_k} b_{k, l} r^{-k}
\int\limits_{\Theta}   Y_{k, l} (\vartheta') \mathcal{K}_n   (r
\vartheta, \vartheta') \, d \vartheta',
$$
    where $K_n (x, y)$ is the Poisson kernel for the sphere $\Theta.$
  The latter integral is  is the Poisson integral.
  It defines a function harmonic in the unit ball and equal to the function
   $Y_{k, l} (\vartheta)$ on the boundary of the ball, i.\,e., on the sphere $\Theta.$
   The only such function is $r^k Y_{k, l} (\vartheta).$ Hence,
$$
\sigma u_2 (r, \vartheta) =  \sum\limits_{k=0}^{\infty}
\sum\limits_{l=1}^{d_k} b_{k, l} Y_{k, l} (\vartheta).
$$
    The right-hand
    side of this relation does not depend on the variable $r.$
    Therefore, it has a limit as $r \to + 0.$ By definition, it is its $\sigma$-trace.
   Therefore,
\begin{equation}
\sigma u|_0 = \sigma u_2|_0 =  \sum\limits_{k=0}^{\infty}
\sum\limits_{l=1}^{d_k} b_{k, l} Y_{k, l} (\vartheta).
\label{8.1.5}
\end{equation}
    Further, since the function $u_2$ is harmonic in the whole space $E^n$ apart from
    the point $\mathsf{0}$ and $b_{k, l}r^{2-n-k}$ are the coefficients of the expansion
    of the function  $u_2$ in the series with respect to spherical
    harmonics, it follows that the relation
$$
b_{k, l}  r^{2-n-k} = \int\limits_{\Theta}    u_2 (r, \vartheta')
Y_{k, l} (\vartheta) d \vartheta,
$$
    holds for each positive $r.$
    By the Cauchy--Bunyakovsky
  inequality, this implies the estimate
$$
|b_{k, l}|  r^{2-n-k} \leq \lr{ \int\limits_{\Theta}    |u_2 (r,
\vartheta)|^2 \, d \vartheta}^{\frac{1}{2}}.
$$
Thus, $\sigma u_2|_0  \in A (\Theta),$ which completes the first
part of the lemma. To prove the second one, assume that $\psi \in
A (\Theta).$ Then we have the relation
$$
\psi (\vartheta) = \sum\limits_{k=0}^{\infty}
\sum\limits_{l=1}^{d_k} \psi_{k, l} Y_{k, l} (\vartheta),
$$
    where the coefficients $\psi_{k, l}$ admit the estimate $|\psi_{k,
l}| \leq c\, h^k$ such that  $h$ is any positive number and the
constant $c$ depends on $h,$ but does not depend on the parameters
$k$ and $l.$
  We have to find a function $u$ of kind \eqref{8.1.4} such that $\sigma u_0 |_0=\psi,$ while
   the
   $\sigma$-trace of the function $u$ coincides with the $\sigma$-trace
    of the function $u_2$ and is computed according to the relation \eqref{8.1.5}.
    Since the expansion with respect to spherical harmonics is unique, it follows that it is
    to be assigned $b_{k, l}= \psi_{k, l}$ in relation \eqref{8.1.4}.
    Thus, it is proved that the sought function is to be of the form
$$
u (r, \vartheta) = u_1 (r, \vartheta) + \sum\limits_{k=0}^{\infty}
\sum\limits_{l=1}^{d_k} \psi_{k, l}  r^{2-n-k}  Y_{k, l}
(\vartheta),
$$
    where $u_1$ is an arbitrary function harmonic in a neighborhood
    of the point $\mathsf{0}$. The following estimate is well known for
  spherical harmonics:
$$
| Y_{k, l} (\vartheta) | \leq c_n \, k^{\frac{n}{2}-1}, \  k = 0,
1, \dots; \  l = 1, \dots, d_k.
$$
    This implies that the function series from the previous relation
 converges absolutely and uniformly outside each ball centered at
 the point $\mathsf{0}.$
 Each term of this series is a harmonic function. Therefore, the sum
 of the series is a function harmonic outside the point $\mathsf{0}$.

For the case where $n=2$, the arguing is the same.
    The only difference is a logarithmic term arising in the expansion with respect
  to  spherical harmonics.
\end{proof}
  Note that singular points of harmonic functions can be classified in terms of
   $\sigma$-traces. If $\left. \sigma u \right|_0 = 0,$ then $\mathsf{0}$
   is a removable singular point. If the series expansion  of the
   function  $\psi = \left.\sigma u \right|_0 $ with respect
  to  spherical harmonics contains a finite amount of terms, then the function
   $u$ has a singularity of the pole type.
   Otherwise,   $u$ has an essential singularity.

\begin{lemma} \label{lem: 8.1.2}
 Let $f \in M^s (\Omega_0),$ $g \in H^{s+ \frac{3}{2}} (\pr \Omega),$
     $\psi \in A (\Theta),$ and $s \geq 0.$
 The problem \eqref{8.1.1}--\eqref{8.1.3}
 has no more than one solution in the space $M^{s+2} (\Omega_0).$
\end{lemma}

\begin{proof}
The solution of the boundary-value
 problem  $\Delta v = 0,$ $\left. v
\right|_{\pr \Omega}=0,$ $\left. \sigma u \right|_0 = 0$ has a
removable singularity at the point $\mathsf{0}$.
 Then  $v \in H^{s+2} (\Omega)$ and the homogeneous boundary-value
 problem  has only the trivial solution.
\end{proof}

\subsection{Solutions: existence and a priori estimates}\label{sec15.2}

In this section, we complete the proof of Theorem \ref{teo:
8.1.1}.

\begin{lemma} \label{lem: 8.2.1}
    Let $U_{\ov{R}, 0} = U_{\ov{R}} \setminus \mathsf{0},$ where
$U_{\ov{R}}$ is the ball of radius $\ov{R}$ in $E^n$ centered at
the point $\mathsf{0}.$ Let $s$ be nonnegative and even, $f \in
M^s (U_{\ov{R}, 0}),$ and $f$ vanish near the boundary of the ball
$U_{\ov{R}}.$ Then the exists a function $v \in M^{s+2}
(U_{\ov{R}, 0})$ such that
\begin{equation}
\Delta v= f(x), \  x \in U_{\ov{R}, 0},
\label{8.2.1}
\end{equation}
\begin{equation}
\left. \sigma v \right|_0 = 0, \label{8.2.2}
\end{equation}
 and the operator $f \to v$ continuously acts from the space
$M^s (U_{\ov{R}, 0})$ into the space $M^{s+2} (U_{\ov{R}, 0}).$
\end{lemma}

\begin{proof}
    We consider only the case where $n\geq 3.$ For  the case where  $n=2,$
    the arguing is the same. First, assume that $f \in\mathring{T}^{\infty} (U_{\ov{R}, 0}).$
    Recall that it means the validity of the expansion
\begin{equation}
    f (r, \vartheta) = \sum\limits_{k=0}^{\mathcal{K}}
     \sum\limits_{l=1}^{d_k} f_{k, l} (r) Y_{k, l} (\vartheta),
    \label{8.2.3}
\end{equation}
 where $\mathcal{K}=\mathcal{K} (f)$ is a positive integer, and $r^{-k} f_{k, l} \in
\mathring{C}^{\infty}_{\frac{n}{2}+k-1} (0, \ov{R})$ (see Chap.
\ref{ch2}). In this case, the desired solution $v$ has the form
\begin{equation}
v (r, \vartheta) = - \sum\limits_{k=0}^{\mathcal{K}}
\sum\limits_{l=1}^{d_k} Y_{k, l} (\vartheta) r^k
\int\limits_r^{\ov{R}} t^{1-2k-n} \int\limits_0^t \tau^{n+k-1}
f_{k, l} (\tau) \, d \tau dt. \label{8.2.4}
\end{equation}
For each function from $ \mathring{C}^{\infty}_{\nu} (0, \ov{R}),$
    its singularity at the origin is at most a power function of order $- 2 \nu.$
    This means that the iterated integral in \eqref{8.2.4} is equal to $O(r^{4-2k-n})$ and,
    therefore, \eqref{8.2.2} is satisfied. To verify \eqref{8.2.1}, we directly
    differentiate \eqref{8.2.4}.

Let us show that the operator $f \to v$ defined by relation
 \eqref{8.2.4} is continuous in the corresponding topologies.

    Let $v_{k, l}$ denote a function of the kind
$$
v_{k, l} = - r^k \int\limits_r^{\ov{R}} t^{1-2k-n} \int\limits_0^t
\tau^{n+k-1} f_{k, l} (\tau) \, d \tau dt.
$$
    Let $2 R < \ov{R}.$ Then
    \begin{multline*}
v_{k, l} (r) = - r^k \int\limits_r^{\ov{R}} t^{1-2k-n}
\int\limits_0^t \tau^{n+k-1} \chi_{R/4} f_{k, l} (\tau) \, d \tau
dt
 \\
 - r^k \int\limits_r^{\ov{R}} t^{1-2k-n} \int\limits_0^t
\tau^{n+k-1} (1-\chi_{R/4}) f_{k, l} (\tau) \, d \tau dt = v_{k,
l}^1+v_{k, l}^2.
    \end{multline*}
    Introduce the functions  $v^j,$ $j=1,2,$ by the relation
$$
v^j (r, \vartheta) = \sum\limits_{k=0}^{\mathcal{K}}
\sum\limits_{l=1}^{d_k} v^j_{k, l} (r) Y_{k, l} (\vartheta)
$$
    and estimate each one separately.
    First, consider the function $v^1 (r, \vartheta).$
    Due to the Leibnitz rule, we have the relation
\begin{equation}
B_{\frac{n}{2}+k-1} \lr{\chi_R \, r^{-k} v^1_{k,l}} = \chi_R
B_{\frac{n}{2}+k-1} \lr{r^{-k} v^1_{k, l}} + 2 \frac{\pr
\chi_R}{\pr r} \frac{\pr r^{-k} v^1_{k, l}}{\pr r} + r^{-k}
v^1_{k, l} B_{\frac{n}{2}+k-1} \chi_R. \label{8.2.5}
\end{equation}
    The first term of the right-hand
 side is represented by the relation
$$
\chi_R  B_{\frac{n}{2}+k-1} \lr{ r^{-k} v^1_{k,l}} = \chi_R \,
r^{-k} \chi_{R/4} f_{k, l} (r) =
 \chi_{R/4} (r) r^{-k} f_{k, l} (r)
$$
    because $\chi_R  \chi_{R/4} =  \chi_{R/4}.$ Taking into account that
     $D\chi_R (r) =0$ for $0 \leq r \leq R$ and  $\chi_{R/4} (r) = 0$
     for $r \geq \dfrac{R}{2},$ we express the second term of \eqref{8.2.5}
as follows:
\begin{equation}
2 \frac{\pr \chi_R}{\pr r} \frac{\pr (r^{-k} v^1_{k, l})}{\pr r} =
2  \frac{\pr \chi_R}{\pr r} r^{1-2k-n} \int\limits_0^{R/2}
\tau^{n+k-1} \chi_{R/4} f_{k, l} \, d \tau. \label{8.2.6}
\end{equation}
 From the same considerations, we obtain the following representation of the third
 term:
\begin{multline}
r^{-k} v^1_{k, l} B_{\frac{n}{2}+k-1} \chi_R (r)
= - \lr{B_{\frac{n}{2}+k-1} \chi_R} \int\limits_r^{\ov{R}} t^{1-2k-n}
 \int\limits_0^t \tau^{n+k-1} \chi_{R/4} f_{k, l} (\tau) \, d \tau dt  \\
= \lr{B_{\frac{n}{2}+k-1} \chi_R}
\frac{r^{2-2k-n}-\ov{R}^{2-2k-n}}{2-2k-n} \int\limits_0^{R/2}
\tau^{n+k-1} \chi_{R/4} f_{k, l} (\tau) \, d \tau. \label{8.2.7}
\end{multline}
  Taking into account that
$$
\frac{2}{r} D \chi_R +  \frac{1}{2-2k-n} B_{\frac{n}{2}+k-1}
\chi_R = \frac{-1}{n+2k-2} B_{1-\frac{n}{2}-k} \, \chi_R,
$$
    we conclude that \eqref{8.2.5}--\eqref{8.2.7}
    imply the relation
\begin{multline}
 B_{\frac{n}{2}+k-1} \lr{ \chi_R (r) r^{-k} v^1_{k, l}}
  = \chi_{R/4} r^{-k} f_{k, l}   \\
- \frac{1}{n+2k-2}  \lr{B_{1-\frac{n}{2}-k} \chi_R} r^{2-2k-n}
 \int\limits_0^{R/2} \tau^{n+k-1} \chi_{R/4} f_{k, l} (\tau) \, d \tau  \\
+ \frac{1}{n+2k-2} \ov{R}^{2-2k-n}  \lr{B_{\frac{n}{2}+k-1} \chi_R
}  \int\limits_0^{R/2} \tau^{n+k-1} \chi_{R/4} f_{k, l} (\tau) \,
d \tau. \label{8.2.8}
\end{multline}
    Now, by virtue of relation  \eqref{7.1.2}, we obtain the following estimate
    of the function $v^1$:
\begin{multline}\label{8.2.9}
\| v^1 \|^2_{s+2, R} = \sum\limits_k \sum\limits_l \| \chi_R \,
r^{-k} v^1_{k,l} \|^2_{\mathring{H}^{s+2}_{\frac{n}{2}+k-1 (0,
2R)}} + \| (1-\chi_R)  v^1 \|^2_{{H}^{s+2} (U_{\ov{R}})}
\\
\leq 3  \sum\limits_k \sum\limits_l \| \chi_{R/4} \, r^{-k}
f_{k,l} \|^2_{\mathring{H}^{s}_{\frac{n}{2}+k-1} (0, 2R)} +
3  \sum\limits_k \sum\limits_l \frac{1}{(n+2k-2)^2}
 \left(  \| r^{2-2k-n} B_{1-\frac{n}{2}-k}
  \chi_R \|^2_{\mathring{H}^{s}_{\frac{n}{2}+k-1 (0, 2R)}}   \right.  \\
+ \left. \ov{R}^{4-4k-2n} \|B_{\frac{n}{2}+k-1} \chi_R
\|^2_{\mathring{H}^{s}_{\frac{n}{2}+k-1} (0, 2R)}   \right) \left|
\int\limits_0^{R/2} \tau^{n+k-1} \chi_{R/4} f_{k, l} \, d \tau
\right|^2 + \|(1-\chi_R) v^1 \|^2_{H^{s+2} (U_{\ov{R}})}
\\
:=3 \Sigma_1 +3 \Sigma_2 + \|(1-\chi_R) v^1 \|^2_{H^{s+2}
(U_{\ov{R}})}.
\end{multline}
  Estimate each term of the right-hand side of the last relation.
  The term $\Sigma_1$ admits the obvious estimate
\begin{equation}
\Sigma_1 \leq \| f \|^2_{s, R/4}. \label{8.2.10}
\end{equation}
  Consider the integral contained in the second term.
  Introduce the notation
 $\nu = \dfrac{n}{2}+k-1,$ $\omega (r) =\chi_{\frac{R}{4}} r^{-k} f_{k, l}.$
 Then the specified integral takes the form
$$
Q_{\nu} (\omega, R) = \int\limits_0^{R/2} \tau^{2 \nu +1}
\omega(\tau) \, d \tau =  \int\limits_0^{R/2} \tau^{2 \nu +1}
P_{\nu}^{\frac{1}{2}-\nu} S_{\nu}^{\nu-\frac{1}{2}} \omega \, d
\tau,
$$
    where $P_{\nu}^{\frac{1}{2}-\nu}$ and $S_{\nu}^{\nu-\frac{1}{2}}$
    are the transmutation operators from Sec. \ref{sec4.1}.
 Since $r^{-k}f_{k, l} \in \mathring{C}^{\infty}_{\frac{n}{2}+k-1} (0, \ov{R}),$
 it follows from the definition of the class $\mathring{C}^{\infty}_{\nu} (0, \ov{R})$
 that the functions $S_{\nu}^{\nu-\frac{1}{2}} \omega  =
S_{\nu}^{\nu-\frac{1}{2}} (\chi_R \, r^{-k} f_{k, l})$ belong to
the space $\mathring{C}^{\infty} [0, R).$
 Let $\widetilde{\omega} = S_{\nu}^{\nu-\frac{1}{2}} \omega$ and
$\nu< N + \dfrac{1}{2},$ where $N$ is a positive integer.
    Then, by the definition of transmutation operators, we obtain
    that
$$
P_{\nu}^{\frac{1}{2}-\nu} \widetilde{\omega} (\tau) = \frac{(-1)^N
2^{-N} \sqrt{\pi} \, \tau^{2(N-\nu)}}{\Gamma(\nu+1) \,
\Gamma(N-\nu+\frac{1}{2})}  \lr{\frac{\pr}{\pr \tau}
\frac{1}{\tau}}^N \int\limits_{\tau}^{\infty} \tau^{2 \nu}
(t^2-\tau^2)^{N-\nu-\frac{1}{2}} t^{-2 N} \widetilde{\omega} (t)
\, d t.
$$
    Integrating this by parts, we obtain
    that
$$
    \begin{gathered}
Q_{\nu} = \int\limits_0^{R/2} \tau^{2 \nu +1} \omega(\tau) \, d
\tau
    \\
=  \frac{(-1)^N 2^{-N} \sqrt{\pi}}{\Gamma(\nu+1) \,
\Gamma(N-\nu+\frac{1}{2})} \int\limits_0^{R/2} \tau^{2 \nu +1}
\lr{\frac{\pr}{\pr \tau} \frac{1}{\tau}}^N
\int\limits_{\tau}^{\infty} \tau^{2 \nu}
(t^2-\tau^2)^{N-\nu-\frac{1}{2}} t^{-2 N} \widetilde{\omega} (t)
\, d t d \tau
    \\
= \frac{ \sqrt{\pi} \, \Gamma(N+\frac{3}{2})}{\Gamma(\nu+1) \,
\Gamma(N-\nu+\frac{1}{2})} \int\limits_0^{R/2} \tau^{2 \nu +1}
\int\limits_{\tau}^{\infty}  (t^2-\tau^2)^{N-\nu-\frac{1}{2}}
t^{-2 N} \widetilde{\omega} (t) \, d t d \tau
    \\
= \frac{ \sqrt{\pi} \, \Gamma(N+\frac{3}{2})}{\Gamma(\nu+1) \,
\Gamma(N-\nu+\frac{1}{2}) \, \Gamma(\frac{3}{2})}
\int\limits_0^{R/2}    \widetilde{\omega} (t) t^{-2 N}
\int\limits_{0}^{t} \tau^{2 \nu +1}
(t^2-\tau^2)^{N-\nu-\frac{1}{2}}  \, d \tau d t.
    \end{gathered}
$$
    The internal integral is expressed via Euler functions as follows:
$$
\int\limits_{0}^{t} \tau^{2 \nu +1}
(t^2-\tau^2)^{N-\nu-\frac{1}{2}}  \, d \tau = t^{2 N+1}
\frac{\Gamma(\nu+1) \, \Gamma(N-\nu+\frac{1}{2})}{ 2 \,
\Gamma(N+\frac{3}{2})}.
$$
    Hence, we have the relation
$$
Q_{\nu} = \int\limits_0^{R/2} t \widetilde{\omega}(t) \, d t =
\int\limits_0^{R/2} t S_{\nu}^{\nu-\frac{1}{2}}   \omega(t) \, d
t.
$$
 Further, since $S_{\nu} = I^{\frac{1}{2} - \nu}
S_{\nu}^{\nu-\frac{1}{2}},$ where $I^{\mu}$ is the Liouville
operator, it follows from the previous relation that
$$
    \begin{gathered}
Q_{\nu} = \int\limits_0^{R/2} t  I^{s+\nu-\frac{1}{2} } I^{-s}
S_{\nu} \, \omega (t)  \, d t =
 \frac{1}{\Gamma\lr{s+\nu-\frac{1}{2}}} \int\limits_0^{R/2} t
\int\limits_t^{R/2} (\tau-t)^{s+\nu-\frac{3}{2}} I^{-s} S_{\nu} \,
\omega (t)  \, d \tau d t
    \\
= \frac{1}{\Gamma\lr{s+\nu-\frac{1}{2}}} \int\limits_0^{R/2}
\lr{I^{-s} S_{\nu} \, \omega(\tau)} \int\limits_0^{\tau} t
(\tau-t)^{s+\nu-\frac{3}{2}}   \,  d t d \tau.
    \end{gathered}
$$
    Since
$$
 \int\limits_0^{\tau} t (\tau-t)^{s+\nu-\frac{3}{2}}   \,  d t
 =\tau^{s+\nu+\frac{1}{2}} \frac{\Gamma(s+\nu-\frac{1}{2})}{\Gamma(s+\nu+\frac{3}{2})},
$$
    it follows that
$$
Q_{\nu} =  \frac{1}{\Gamma(s+\nu+\frac{3}{2})} \int\limits_0^{R/2}
\tau^{s+\nu+\frac{1}{2}} I^{-s} S_{\nu} \, \omega(\tau) \, d \tau.
$$
  Now, from the  Cauchy--Bunyakovsky
  inequality, we obtain the estimate
$$
    \begin{gathered}
|Q_{\nu}| \leq  \frac{1}{\Gamma(s+\nu+\frac{3}{2})} \lr{
\int\limits_0^{R/2} \tau^{2 s+2\nu+1} \, dt}^{1/2} \|D^s S_{\nu}\,
\omega \|_{L_2 (0, R/2)}
    \\
= \frac{R^{s+\nu+1}}{2^{s+\nu+\frac{3}{2}}\sqrt{s+\nu+1} \, \Gamma
(s+\nu+\frac{3}{2})}  \|S_{\nu} \, \omega \|_{\mathring{H}^s (0,
R/2)} =
 \frac{R^{s+\nu+1} \, \Gamma(\nu+1)}{2^{s+1} \sqrt{s+\nu+1} \,
\Gamma (s+\nu+\frac{3}{2})}  \| \omega \|_{\mathring{H}^s_{\nu}
(0, R/2)}.
    \end{gathered}
$$
    Due to relation \eqref{1.4.20} for gamma-functions,
    this implies the inequality
$$
|Q_{\nu}| \leq c (s, R) R^{\nu} \nu^{-1-s} \| \omega
\|_{\mathring{H}^s_{\nu} (0, R/2)}.
$$
    Getting back to the previous notation, we obtain the following final
  form  of the estimate of the integral $Q_{\nu}$:
\begin{equation}
|Q_{\frac{n}{2}+k-1} (\chi_{R/4} \, r^{-k} f_{k, l}, R)| \leq c
(n, s, R)  \frac{R^k}{(k+1)^{1+s}} \| \chi_{R/4} \, r^{-k} f_{k,
l} \|_{\mathring{H}^s_{\frac{n}{2}+k-1} (0, R/2)}. \label{8.2.11}
\end{equation}
    It follows from results of Sec. \ref{sec4.1} that
$$
    \begin{gathered}
\| B_{\frac{n}{2}+k-1} \chi_R
\|^2_{\mathring{H}^s_{\frac{n}{2}+k-1} (0, 2 R)} \leq 2 \|
B_{\frac{n}{2}+k-1} \chi_R \|^2_{\mathring{H}^s_{\frac{n}{2}+k-1,
+} (0, 2 R)} =
 2 \int\limits_0^{2R} |B_{\frac{n}{2}+k-1}^{\frac{s}{2}+1} \chi_R
|^2 \, r^{2 (\frac{n}{2}+k-1)+1} dr
    \\
= 2 R^{n+2k-s-2} \int\limits_0^{2}
|B_{\frac{n}{2}+k-1}^{\frac{s}{2}+1} \chi(t) |^2 \, t^{n+2k-1} dt
=
 c (s, n, k) \, R^{n+2k-s-2} \leq c (s, n) \, 2^{2 k} \,
R^{n+2k-s-2} (k+1)^{s+1}.
    \end{gathered}
$$
  Proving Theorem \ref{teo: 7.2.2}, we obtained the inequality
$$
\| r^{2-2k-n} B_{1-\frac{n}{2}-k} \chi_R
\|_{\mathring{H}^s_{\frac{n}{2}+k-1} (0, 2 R)} \leq c (s, n, R) \,
R^{-2 k} (k+1)^{s+1}.
$$
    Combining the last two inequalities with \eqref{8.2.11}, we arrive at
    the following estimate of a term from relation \eqref{8.2.9}:
$$
   \begin{gathered}
\Sigma_2 \leq c \sum\limits_k \sum\limits_l R^{2k} (k+1)^{-4-2s}
\left( R^{-2 k} (k+1)^{s+1} {}\right.
  \\
\left.{} + \ov{R}^{4-4k-2n} R^{n+2k-s-2} \, 2^{2k} (k+1)^{s+1}
\right)   \| \chi_{R/4} \, r^{-k} f_{k, l}
\|^2_{\mathring{H}^s_{\frac{n}{2}+k-1} (0, R/2)}.
    \end{gathered}
$$
    Since $2 R < \ov{R},$ we obtain the following final estimate of the term $\Sigma_2$:
\begin{equation}
\Sigma_2 \leq c \sum\limits_k \sum\limits_l \| \chi_{R/4} \,
r^{-k} f_{k, l} \|^2_{\mathring{H}^s_{\frac{n}{2}+k-1} (0, R/2)}
\leq c \, \|f \|^2_{s, R/4}, \label{8.2.12}
\end{equation}
    where $c$ is a positive constant independent of $f.$

To complete the estimate of the function $v^1$, it remains to
consider the last term of \eqref{8.2.9}. Since $\chi_{\ov{R}} (r)
= 1$ in the ball $U_{\ov{R}},$ it follows that
\begin{equation}
\| (1 - \chi_R) v^1 \|_{H^{s+2} (U_{\ov{R}}) } \leq c  \|
\chi_{\ov{R}} (1 - \chi_R) v^1 \|_{\mathring{H}^{s+2} (U_{2
\ov{R}}) }. \label{8.2.13}
\end{equation}
    Similarly to relations \eqref{8.2.5}--\eqref{8.2.8},
    we have the relation
$$
   \begin{gathered}
B_{\frac{n}{2}+k-1} (\chi_{\ov{R}} (1 - \chi_R) r^{-k} v^1_{k, l})
 \\
= \frac{1}{2k+n-2}  \lr{B_{\frac{n}{2}+k-1} (\chi_{\ov{R}} (1 -
\chi_R)) \ov{R}^{2-2k-n} - B_{1-\frac{n}{2}-k} (\chi_{\ov{R}} (1 -
\chi_R) r^{2-2k-n}}   \int\limits_{0}^{R/2} \tau^{n+k-1}
\chi_{R/4} f_{k, l} \, d \tau.
   \end{gathered}
$$
    Hence,
$$
   \begin{gathered}
\| \chi_{\ov{R}} (1 - \chi_R) v^1 \|^2_{\mathring{H}^{s+2} (U_{2
\ov{R}}) } \leq \sum\limits_k \sum\limits_l
\frac{|Q_{\frac{n}{2}+k-1}|^2}{(2k+n-2)^2}  \left( \| r^{2-2k-n}
B_{1-\frac{n}{2}-k}^{\frac{s+2}{2}} (\chi_{\ov{R}} (1 - \chi_R))
\|^2_{L_2,  \frac{n}{2}+k-1}   {} \right.
 \\
\left.{} +  \ov{R}^{4-4k-2n}  \|
B_{\frac{n}{2}+k-1}^{\frac{s+2}{2}} (\chi_{\ov{R}} (1 - \chi_R))
\|_{L_2,  \frac{n}{2}+k-1}  \right).
   \end{gathered}
$$
    The norms from the last sum are already estimated (in
    particular, this is done within the proof of Theorem \ref{teo: 7.2.2}).
 This yields the inequality
$$
\| \chi_{\ov{R}} (1 - \chi_R) v^1 \|^2_{\mathring{H}^{s+2} (U_{2
\ov{R}}) } \leq c\, \sum\limits_k \sum\limits_l (k+1)^{-s-1}  \|
\chi_{R/4} \, r^{-k} f_{k, l}
\|^2_{\mathring{H}^s_{\frac{n}{2}+k-1} (0, R/2)} \leq c \, \| f
\|^2_{s, R/4}.
$$
    Therefore,
\begin{equation}
\| (1 - \chi_R) v^1 \|_{{H}^{s+2} (U_{ \ov{R}}) } \leq c \, \| f
\|^2_{s, R/4}. \label{8.2.14}
\end{equation}
Thus, the following estimate of the function $v^1$ is proved (see
\eqref{8.2.10}, \eqref{8.2.12}, and \eqref{8.2.14}):
\begin{equation}
\| v^1 \|_{s+2, R} \leq c \, \| f \|_{s, R/4}, \label{8.2.15}
\end{equation}
where $c$ is a positive constant independent of the function $f.$

Consider the function $v^2.$ It is easy to see that it belongs to
the space ${H}^{s+2} (U_{ \ov{R}})$ and satisfies the
boundary-value problem
$$
\Delta v^2 = (1 - \chi_{R/4}) f(x), \  x \in U_{ \ov{R}},
$$
$$
\left. v^2 \right|_{\pr U_{\ov{R}}} = 0.
$$
    Since  $f \in M^s (U_{ \ov{R}, 0}),$ it follows that
 $(1 -\chi_{R/4}) f \in {H}^{s} (U_{ \ov{R}}).$
 The solution uniqueness of this Dirichlet boundary-value
 problem for the Poisson equation and the following estimate of the solution
 are well known:
$$
\| v^2 \|_{{H}^{s+2} (U_{ \ov{R}}) } \leq c \, \| (1 - \chi_{R/4})
f \|_{{H}^{s} (U_{ \ov{R}}) }.
$$
Now, from the definition of norms $\| \|_{s, R}$ and from Theorem
\ref{teo: 7.1.1}, we obtain the estimate
\begin{equation}
\| v^2 \|_{s+2, R} \leq c \, \| f \|_{s, R/4}, \label{8.2.16}
\end{equation}
    where the constant does not depend on the function $f.$

Thus, an estimate for each term of the decomposition $v=v^1+v^2$
is obtained. Hence, for each $R \in \Big(0,\dfrac{\ov{R}}{2}\Big)$
 and each nonnegative $s$ there exists a positive constant $c$
 such that each $f \in \mathring{T}^{\infty}(U_{\ov{R}, 0})$
 satisfied the inequality
$$
\| v \|_{s+2, R} \leq c \, \| f \|_{s, R/4}.
$$
  To complete the proof of the lemma, we pass to the limit.
  Let
$f \in M^s (U_{\ov{R}, 0})$ and $f$ satisfy the assumption of the
lemma.
    Then there exists a sequence of functions $f^m \in \mathring{T}^{\infty}
(U_{\ov{R}, 0})$, converging to $f$ with respect to the topology
of this space. For each function $f^m$, define the function $v^m$
 by means of relation \eqref{8.2.4}. Then
\begin{equation}
\Delta v^m = f^m \to f, \  m \to \infty. \label{8.2.17}
\end{equation}
    As we proved above, the map $f^m \to v^m$ from the space $M^s (U_{\ov{R}, 0})$
    into $M^{s+2} (U_{\ov{R}, 0}).$
    Then  $\{v^m\}$ is a fundamental sequence in $M^{s+2}(U_{\ov{R}, 0}).$
    Then, since $M^{s+2}(U_{\ov{R}, 0})$ is a complete space, it follows that
    there exists a function $v \in M^{s+2} (U_{\ov{R},0})$ such that it is the limit
    of the sequence of functions  $v^m$ in the sense of the topology of this space.

    The operator $\Delta$ continuously maps the space $M^{s+2}(U_{\ov{R}, 0})$ into
 $M^{s} (U_{\ov{R}, 0}).$ Therefore, $\Delta v^m\to \Delta v,$ $m \to \infty,$
 in the sense of the space $M^s(U_{\ov{R}}).$ Then \eqref{8.2.17} implies that
  $\Delta v= f.$ Due to the direct $\sigma$-trace
  theorem, we have the relation
$$
\left. \lim\limits_{m \to \infty} \sigma v^m \right|_0 = \left.
\sigma v \right|_0,
$$
    where the limit passage is treated in the sense of the space $A(\Theta).$
    Since $\left. \sigma v^m \right|_0=0,$ it follows that
$\left. \sigma v \right|_0=0$ as well, which completes the proof
of Lemma \ref{lem: 8.2.1}.
\end{proof}

\begin{proof}
[Proof of Theorem \ref{teo: 8.1.1}]
    The uniqueness assertion for solutions of problem \eqref{8.1.1}--\eqref{8.1.3}
    is established by Lemma \ref{lem: 8.1.2}.
    The existence of as solution is proved as follows. Let  $\ov{R}$ be the diameter
    of the domain $\Omega$ and a function $u^1 \in M^{s+2} (U_{\ov{R}, 0})$
  be the solution of the boundary-value
    problem
\begin{equation}
\begin{cases}
\Delta u^1 =  \chi_{R_0}\, f(x), &  x \in U_{ \ov{R}, 0},\\
\left. \sigma u^1 \right|_{0} = 0,&
\end{cases}
\label{8.2.18}
\end{equation}
 constructed in Lemma \ref{lem: 8.1.1}.
 Let a function  $u^2$ satisfy the boundary-value
 problem
\begin{equation}
\begin{cases}
\Delta u^2 = (1- \chi_{R_0}) f(x), &  x \in \Omega,\\
\left. u^2 \right|_{\pr \Omega} = g - \left. u^1 \right|_{\pr
\Omega} - \left. u^3 \right|_{\pr \Omega},&
\end{cases}
\label{8.2.19}
\end{equation}
    where  $u^3$ is the harmonic in $E^n \setminus
\mathsf{0}$ function constructed in Theorem \ref{teo: 7.2.2} and
such that
$$
\left. \sigma u^3 \right|_{0} = \psi.
$$
    Then the function $u = u^1+u^2+u^3$ satisfies problem \eqref{8.1.1}--\eqref{8.1.3}.

The function $u$ belongs to the space $M^{s+2} (\Omega_0)$ because
$u^1 \in M^{s+2} (\Omega_0)$ due to Lemma \ref{lem: 8.2.1}, $u^3
\in M^{s+2} (\Omega_0)$  due to Theorem \ref{teo: 7.2.2}, and $u^2
\in H^{s+2} (\Omega) \subset M^{s+2} (\Omega_0)$ due to the
general theory of elliptic boundary-value
 problems (see, e.\,g., \cite{51, 62}). These results imply the continuity of the resolving
 operator $\Lambda^{-1},$   which completes the proof
of Theorem \ref{teo: 8.1.1}.
\end{proof}

    %\newpage

\chapter{Constructing of Transmutation Operators: Composition Methods}\label{ch6}

\section[Constructing of Transmutation Operators:
 General Scheme of Composition Methods]{Constructing of Transmutation Operators:\\
 General Scheme of Composition Methods}\label{sec16}

Each class of transmutation operators, investigated above, is
constructed by its own method.
    This causes a necessity to develop a general scheme to construct
     transmutation operators. In the present chapter, we propose such a scheme.
     This is the factorization (composition) method based on the representation of
transmutation operators by a composition of integral
transformations.
 The composition methods provides algorithms to construct new transmutation operators.
 Also, it includes (as special cases)
  So\-nin--Poisson--Delsarte,
    Vekua--Erd\'elyi--Lowndes,
    and Buschman--Erd\'elyi
 transmutation operators of various types,
 Sonin--Katrakhov and Poisson--Katrakhov
    unitary transmutation operators, generalized
    Erd\'elyi--Kober operators, and the classes of elliptic, hyperbolic, and
    parabolic transmutation operators introduced in \cite{Car1, Car2, Car3}.
 In the present chapter, we generalize them, introduce classes of
$B$-elliptic, $B$-hyperbolic,
 and $B$-parabolic
  transmutation operators.

The composition method to construct transmutation operators is
developed in
\cite{S6,S66,S7,S5,S46,S14,S400,S42,S38,S401,SitDis,FJSS}.
    The initial idea close to this method is applied by Katrakhov to construct
    a special class of  transmutation operators (in the present monograph, they are
    called
    Sonin--Katrakhov and Poisson--Katrakhov
 transmutation operators, see Chap. \ref{ch3}).

In \cite{FJSS}, the term {\it Integral Transforms Composition
Method} (ITCM) is proposed for the composition method.

In this section, we provide main definitions and a small part of
results that can be obtained by means of the composition method.
    Since the formulations of the results are quite cumbersome, we omit their proofs.
 In the future, it is planned to present the the composition method completely
 (including numerous examples).

The general scheme of the composition method is as follows. The
input data is a pair $A,B$ of operators of an arbitrary type and
the related generalized invertible Fourier transformations $F_A$
and $F_B$ acting by relations
\begin{equation}
\label{4301} F_A A =g(t)F_A \,  \textrm{ and } \, F_B B= g(t) F_B,
\end{equation}
    where $t$ is the dual variable and $g(t)$ is an arbitrary suitable function
    (the most obvious choice is $g(t)=-t^2,$ cf. classical integral transformations).

    % \medskip

\textit{The goal of the composition is to construct
{\rm(}formally{\rm)} an output pair of mutually inverse
transmutation operators  $P$ and $S$ acting according to the
relations}
\begin{equation}\label{4302}
S=F^{-1}_B \frac{1}{w(t)} F_A \,  \textrm{ and } \,  P=F^{-1}_A
w(t) F_B,
\end{equation}
\textit{where $w(t)$ is an arbitrary weight function. Then $P$ and
$S$ are mutually inverse transmutation operators  intertwining the
original operators $A$ and} $B$:
\begin{equation}\label{Inter}
SA=BS \,  \textrm{ and } \, PB=AP.
\end{equation}

    %   \medskip

The formal verifications is the direct substitution.
 The main difficulty is to compute the introduced composition in an explicit integral form
 and to find the corresponding domains of the operators.

The main merits of the composition method are as follows.

\begin{itemize}
\item{The composition method is simple: numerous  transmutation
operators are constructed from basic blocks (classical integral
transformations) according to a simple rule.}
 \item{All   transmutation
operators known earlier in the explicit form can be constructed by
means of the (unified) composition method.}
 \item{Using the (unified) composition method, one can construct various new
 transmutation
operators.}
 \item{The composition method provides a simple possibility to obtain inverse transformations
 in the same composition form.}
  \item{The
composition method  provides a possibility to derive estimates of
norms of direct and inverse  transmutation operators,
    using known estimates of norms for classical integral
    transformations in various function spaces.}
\item{The composition method  provides a possibility to obtain
explicit relations binding solutions of the intertwined
differential equations.}
\end{itemize}
    The composition method has a disadvantage as well: usually,
    the action of basic blocks (integral transformations) is known in classical spaces,
    but, applying them to differential equations or differential operators for particular
    problems, we need estimates and the action on functions of other classes, e.\,g.,
    classes of functions vanishing at the origin or/and infinity.
    In such cases, one can start, applying  the composition
method to obtain an explicit form of the transmutation operator,
and to extend it to desired spaces afterwards.

We emphasize that relations of kind \eqref{4301}-\eqref{4302}
    are new nor for the theory of integral transformations neither
    for applications to differential equations.
    However, the composition method is new for applications in the theory of
    transmutation operators.

Selecting the classical Fourier transformation  in other problems
for integral transformations and related differential equations,
we see that relations \eqref{4302} of the composition method lead
to the class of pseudodifferential operators with symbols $w(t)$
and $\dfrac{1}{w(t)}$ (see, e.\,g., \cite{OlRa, Agr}). Selecting
the classical Fourier transformation  and the weight function
$w(t)=(\pm it)^{-s}$ for the operators $A=B=D^2$ and $F_A=F_B=F,$
    we obtain the fractional Riemann--Liouville
integral on the whole real line. Selecting the weight function
$w(t)=|t|^{-s}$, we obtain Riesz potentials. Once
$w(t)=(1+t^2)^{-s}$ is selected, relations \eqref{4302} yield
Bessel potentials. For $w(t)=(1\pm it)^{-s}$, we obtain modified
Bessel potentials (see \cite{SKM}).

If we select the classical Hankel integral transformation and the
weight functions
\begin{equation}\label{trans}
 A=B=B_\nu,\quad  F_A=F_B=H_\nu,\quad
g(t)=-t^2,\, \textrm{ and } \, w(t)=j_\nu(st),
\end{equation}
    then we obtain the Delsarte generalized translation
operators (see \cite{Lev2, Lev3, Lev4, Mar9}); recall that the
Bessel operator is denoted by $B_\nu,$ the  Hankel transformation
 is denoted by $H_\nu,$ and the normalized (``small'') Bessel
 function is denoted by $j_\nu(\cdot)$ (see Chap. \ref{ch1}).
 In the  general case, selecting data $A=B, F_A=F_B$ and arbitrary
 weight functions $g(t)$ and $w(t)$ for the composition method, we obtain
 transmutation operators commuting with the operator $A$
 (cf. the generalized translation operator commuting with the Bessel operator).

In \cite{OZK,AbOs,Os}, the quantum oscillator expression is
selected for the operator $A$ and the quadratic Fourier
transformation (fractional Fourier transformation,
Fourier--Fresnel transformation, or Weierstrass transformation) is
selected for the corresponding integral transformation $F_A$.
    This important integral  transformation arises from the Fresnel proposal
    to change standard plane waves with linear independent variables in
    exponential functions for more general waves with quadratic independent variables in
    exponential functions (this provides a possibility to
    completely explain paradoxes with spectral lines under the     Fraunhofer diffraction).
    From the mathematical viewpoint, the operators of the quadratic Fourier
 transformation are fractional powers $F^\alpha$ of the classical  Fourier transformation,
 defined by Wiener and Weyl: they add-on
  the latter till a semigroup with respect to the parameter $\alpha.$
 In the wavelet theory (according to its tradition to treat each relation as a new one
 and endow well known objects with new terms), the quadratic Fourier
transformation is called the Gabor transformation.
 In \cite{AlKi1}, its applications to the Heisenberg group,
 quantum oscillators, and wavelets are recently obtained.
 The composition method explained above allows one to use this
 transformation to construct transmutation operators for the one-dimensional
  Schr\"odinger operator, i.\,e. quantum oscillator (see \cite{S45,S42}).
  The quadratic Hankel transformation, which is more general,  can be used for that purpose
  as well (see \cite{S60}).

    Applying the composition method instead of classical approaches, one can construct various
    explicit integral forms of fractional powers of the Bessel operator
 (see \cite{S135,S133,S127,S123,S18,S42,S700,SS,FJSS}).
 Nor semigroup-theory
 methods neither spectral methods are able to provide such constructions.

The composition method is obviously generalized to the
multidimensional case: in this case, $t$ becomes a vector, while
$g(t)$ and $w(t)$ from relations
    \eqref{4301}-\eqref{4302} become vector-functions.
 Unfortunately, only a few classes of transmutation
operators are known (or can be introduced explicitly)  in this
case. However, well-known
  classes of potentials are among them.
  For example, if the classical Fourier transformation is used
  (in the composition method) and a positive definite quadratic form is taken
  as the weight function $w(t)$ in \eqref{4302}, then we obtain Riesz elliptic
  potentials (see \cite{Riesz,SKM}); if $w(t)$ is an indefinite quadratic
  form, then we obtain Riesz hyperbolic
  potentials (see \cite{Riesz,SKM,Nogin}; if we assign $w(x,t)=(|x|^2-it)^{-\alpha/2},$
  then we obtain parabolic potentials (see \cite{SKM}).
  If, applying the composition method, we use the Hankel transformation in relations
  \eqref{4301}-\eqref{4302} and a quadratic form for the function $w(t)$,
  then we arrive at elliptic or hyperbolic Riesz $B$-potentials
  (see \cite{Lyah3,Gul1} and \cite{ShiE2} respectively).
  In these cases, the theory of generalized functions and their convolutions
  is applied, special averaging and approximation procedures are used to convert
  such potentials (see \cite{Nogin,ShiE2}), and Schwartz or Lizorkin spaces of test
  functions and dual spaces  for generalized functions are taken.

Applying the composition method, one can obtain relations binding
solutions of perturbed and unperturbed differential equations
linked by  transmutation operators. For example, applying the
composition method to equations with Bessel operators of the kind
\begin{equation}\label{GenEPD}
\sum\limits_{k=1}^n A_k\left(   \frac{\partial^2 u}{\partial
x_k^2}+\frac{\nu_k}{x_k} \frac{\partial u}{\partial x_k}\right)
\pm \lambda^2 u=0,
\end{equation}
$$
x_k>0,\qquad A_k=\const,\qquad \nu_k=\const,\qquad \lambda=\const,
$$
we find binding relations for solutions of the unperturbed
equation
\begin{equation}\label{GenWave}
\sum\limits_{k=1}^n A_k  \frac{\partial^2 v}{\partial x_k^2}\pm \lambda^2 v=0.
\end{equation}
    In particular, this class of relations includes binding relations for
solutions of equations with  Bessel operators with different
parameters; they are called translation operators with respect to
parameter (see \cite{FJSS}). A similar idea called the
\emph{subordination principle} is used in the theory of
fractional-order  equations. In fact, it consists of the
constructing of transmutation operators performing a translation
with respect to the order of the equation; in particular, this
allows one to express solutions of fractional-order
  equations via solutions of classical integer-order
  equations (see, e.\,g., \cite{Jan,EiIvKoch,Bajlekova0,Bajlekova1,FJSS}).

Thus, one can conclude that the composition method is efficient
for the constructing of transmutation operators, it  is related to
a number of other known methods and problems, all known explicit
representations for transmutation operators are obtained by means
of it, and numerous new classes of  transmutation operators can be
constructed by means of it as well.
 The algorithm of the composition method is applied in the following three stages
 (steps).

\begin{itemize}
\item Step 1.
    For the given operator pair $A,B$ and related pair  $F_A, F_B$ of generalized Fourier
transformations, determine and compute the pair $P,S$ of mutually
inverse transmutation operators by means of the main relations of
the method, given by \eqref{4301}-\eqref{4302}.

\item Step 2.
    Find exact conditions and determine function classes such that
    the transmutation operators constructed on Step 1 satisfy the intertwining properties
    given by \eqref{Inter}.

\item Step  3.
    On the corresponding function spaces, apply the transmutation operators well defined
    on Steps 1-2
 to solve problems for differential equations, e.\,g., to establish
correspondence relations between solutions of perturbed and
unperturbed differential equations.
\end{itemize}
    Further, the composition method is used to construct
    hyperbolic, elliptic, and parabolic (in the Carroll sense)
 transmutation operators, generalized Erd\'elyi--Kober
 operators, etc.
 Note that transmutation operators of the specified types, treated as operators
 representing solutions of abstract
     equations via solutions of other abstract
     equations, are considered in \cite{Lav1, BD1, BD2, BD3}.
 In particular, applying this variant of the method of transmutation
operators, one obtains unusual relations expressing solutions of
the wave equation via solutions of the heat equation and vise
versa. For the case of the abstract Bessel operator, such
relations are obtained in \cite{Glu13}.

We use the classification of transmutation operators, proposed by
Carroll; it is related to the type of the partial differential
equation such that the kernel of the said operator satisfies it.
    In the same way, we introduce and use natural generalizations
    of these notions such as $B$-elliptic,
    $B$-hyperbolic, and $B$-pa\-rabolic
 transmutation operators.
 To do that, we use the Kipriyanov definitions of the corresponding classes of singular
 partial differential equations.
    For example, in these terms, the classical
    Sonin--Poisson--Delsarte
    intertwining operators are $B$-hyperbolic.

\subsection{$B$-hyperbolic transmutation operators}\label{sec16.1}

Translation operators with respect to the parameter of the Bessel
operator, satisfying the relation
\begin{equation}\label{449}{T B_{\nu} = B_{\mu} T}
\end{equation}
    are called $B$-\emph{hyperbolic transmutation operators}.
  Such operators are sought in the factorized form
\begin{equation}\label{4410}{T_{\nu, \, \mu}^{(\varphi)} = F_{\mu}^{-1}
\lr{\varphi(t) F_{\nu}}.}\end{equation}
   If $\nu=-\dfrac{1}{2}$ or $\mu=-\dfrac{1}{2},$ then such
   operators are reduced to already studied ones.
 Assume that $\varphi (t) = t^{\alpha}$ and $T^{(\varphi)}=T^{(\alpha)}.$

\begin{theorem}
    If
    $$
    -2-2 \Re \mu < \Re \alpha < \Re(\nu-\mu),
    $$
    then the integral representation
    \begin{multline}\label{4411}
    \lr{T^{(\alpha)}_{\nu,\, \mu} f}(x) = \frac{C_1}{x^{2 \mu + \alpha +2}}
     \int\limits_0^{x} y^{2 \nu +1}\, {_2F_1}\lr{\mu+\frac{\alpha}{2}+1, \frac{\alpha}{2}+1;
      \nu+1; \frac{y^2}{x^2}} f(y)  \, dy  {} \\
    {}+ C_2 \int\limits_x^{\infty}  y^{-2\mu+2 \nu-\alpha -1}
       {_2F_1}\lr{\mu+\frac{\alpha}{2}+1, \mu-\nu+\frac{\alpha}{2}+1; \mu+1;
        \frac{x^2}{y^2}} f(y) \,dy,
    \end{multline}
    where
    $$
    C_1=   \frac{2^{ \mu-\nu+\alpha+1} \Gamma \lr{\mu+\frac{\alpha}{2}+1}}
    {\Gamma \lr{-\frac{\alpha}{2}} \Gamma \lr{\nu+1}},\qquad
    C_2 = \frac{2^{ \mu-\nu+\alpha+1} \Gamma \lr{\mu+\frac{\alpha}{2}+1}}
    {\Gamma \lr{\nu-\mu-\frac{\alpha}{2}} \Gamma \lr{\mu+1}},  \nonumber
   $$
  and ${_2F_1}$ is the Gauss hypergeometric function, holds.
\end{theorem}

Consider several special cases of operator \eqref{4411}.

\noindent
 (a) Let $\alpha=-1-2\mu+\nu,$ $\Re \nu > -1,$ and $\Re \mu > -1.$
 Then
\begin{multline*}
    \lr{T^{(-1-2\mu+\nu)}_{\nu,\, \mu} f}(x) = 2^{- \mu} \frac{1}{x^{ \mu}}
      \int\limits_x^{\infty} y^{\nu} (y^2-x^2)^{\frac{\mu}{2}}
      P_{\frac{\nu}{2}-\frac{1}{2}}^{- \mu} \lr{1-2 \frac{x^2}{y^2}}  f(y)  \, dy {} \\
    {}+2^{1- \mu} e^{i \mu \pi}  \frac{\Gamma \lr{\mu+\frac{\nu}{2}+\frac{1}{2}}}
    {\Gamma \lr{\mu-\frac{\nu}{2}+\frac{1}{2}}} \frac{1}{x^{ \mu}}
     \int\limits_0^{x} y^{\nu} (x^2-y^2)^{\frac{\mu}{2}} Q_{\frac{\nu}{2}-\frac{1}{2}}^{ \mu}
      \lr{2 \frac{x^2}{y^2}-1} f(y) \,dy.
\end{multline*}
    (b) Let $\alpha=0$ and $-1< \Re \mu < \Re \nu.$
  In the case, we obtain the remarkable operator of  ``descending'' with respect
   to the parameter, depending neither on the initial nor on the
  final values of the parameters $\nu $ and $\mu,$ but only on the value
    $\gamma = \nu -
\mu$ of the ``descending:''
$$
\lr{T^{(0)}_{\nu,\, \mu} f}(x) =
 \frac{2^{1-(\nu-\mu)}}{\Gamma(\nu-\mu)}
 \int\limits_x^{\infty} y (y^2-x^2)^{\nu-\mu-1} f(y)  \, dy.
$$
    The only difference between this integral and the fractional
    integral $I_{-, x^2}$ with respect to the function $g(x)=x^2$ is a constant
factor.
    This form of the operator is discovered by Erd\'elyi; it is a special case
    of the Erd\'elyi--Kober operators
    or Lowndes operators.

    %\medskip

\noindent
    (c) Let $\alpha=2 \nu$ and $-1< \Re (\nu + \mu) < 0.$
    Then
$$
    \begin{gathered}
\lr{T^{(2 \nu)}_{\nu,\, \mu} f}(x) = \frac{\sin (\pi \mu)}{\pi}
2^{\mu+\nu+1}  \Gamma(\nu+\mu+1) \int\limits_x^{\infty} y^{2 \nu
+1} (y^2-x^2)^{-\mu-\nu-1} f(y)  \, dy
    \\
- \frac{\sin (\pi \nu)}{\pi} 2^{\mu+\nu+1} \Gamma(\nu+\mu+1)
 \int\limits_0^x y^{2 \nu +1} (x^2-y^2)^{-\mu-\nu-1} f(y)  \, dy.
    \end{gathered}
$$
(d) If  $\mu = \nu$ and $-2 \Re \nu -2 < \Re \alpha < 0,$
    then we obtain the following family of operators commutative with $B_{\nu}$:
$$
\lr{T^{(\alpha)}_{\nu,\, \nu} f}(x) =
\frac{2^{\alpha+2}}{\sqrt{\pi}} \frac{e^{-i \frac{\pi \alpha}{2}}}
{\Gamma\lr{-\frac{1}{2}-\frac{\alpha}{2}}}\frac{1}{x^{\nu+\frac{1}{2}}}
\int\limits_0^{\infty} y^{\nu +
\frac{3}{2}}|x^2-y^2|^{-\frac{\alpha}{2}-1} Q_{\nu -
\frac{1}{2}}^{\frac{\alpha}{2}+1} \lr{\frac{x^2+y^2}{2xy}} f(y) \,
dy.
$$
(e) If $\mu = - \nu$ and $2 \Re \nu -2 < \Re \alpha < \Re 2 \nu,$
    then
$$
\lr{T^{(\alpha)}_{\nu, -\nu} f}(x) = 2^{-2\nu+\alpha+1}
\frac{\Gamma\lr{- \nu+\frac{\alpha}{2} +1}}{\Gamma\lr{2
\nu-\frac{\alpha}{2}}} \, x^{\nu} \int\limits_0^{\infty}
|y^2-x^2|^{\nu-\frac{\alpha}{2}-1} y^{\nu + 1} P_{\nu -
\frac{\alpha}{2}-1}^{\nu} \lr{\frac{x^2+y^2}{|x^2-y^2|}} f(y) \,
dy.
$$
    Another approach to the constructing of translation operators with respect to a parameter
    is as follows. Consider a Sonin-type
 operator constructed according to the above relations and such that $\varphi=\varphi_1$ and a
 Poisson-type operators such that
$\nu=\mu$  and $\varphi=\varphi_2.$ Then their composition
$T=P_{\mu} S_{\nu}$ is the desired operator.
    If the Fourier
    sine-transformation or the Fourier cosine-transformation
    are selected in both case, then we obtain \eqref{4410},
     where $\varphi=\dfrac{\varphi_2}{\varphi_1}.$
  Therefore, the problem on the boundedness of such an operator in Lebesgue spaces
  $L_{p, \gamma}(0, \infty)$   with power weights is reduced to the problem on the
  possibility to divide in the spaces of symbols. If different transformations are selected,
  then we arrive at the following factorization:
\begin{equation}\label{4412}{T^{(\varphi_1, \varphi_2)}_{\nu,\,
\mu} = F^{-1}_{\mu} \cdot \varphi_2(t)\cdot F_{\scriptsize\left\{
\begin{matrix} s
\\
c \end{matrix} \right\}} \cdot F^{-1}_{\scriptsize\left\{
\begin{matrix} c
\\
s \end{matrix} \right\}} \frac{1}{\varphi_1(t)}
F_{\nu}.}\end{equation}
 It is easy to show that the composition of the Fourier transformations is reduced to the
 so-called semiaxis Hilbert transformations
\begin{equation}
 (F_s F_c f)(x) = \int\limits_0^{\infty} \frac{x}{x^2-y^2} f(y) \, dy
   \label{4413}
 \end{equation}
    and
\begin{equation}
 (F_c F_s f)(x) = \int\limits_0^{\infty} \frac{y}{y^2-x^2} f(y)
\, dy, \label{4414}
\end{equation}
    where the integral is treated in the principal-value
    sense; these transformations multiplied by suitable constants
    are unitary in $L_2 (0, \infty).$ Thus, in the considered case, operator \eqref{4412}
    is factorized via two Fourier--Bessel
    transformations and one of the two-weight
    transformations
$$
(A_1 f)(x) = x\, \varphi_2 (x) \int\limits_0^{\infty}
\frac{f(y)}{(x^2-y^2)}  \cdot \frac{dy}{\varphi_1 (x)}~ \textrm{
and } ~(A_2 f)(x) = \varphi_2 (x) \int\limits_0^{\infty}
\frac{y}{\varphi_1 (x)} \cdot \frac{f(y)}{(y^2-x^2)} \,  dy.
$$
    In this case, the boundedness problem for \eqref{4412} is reduced to the problem on
    two-weight estimates for the semiaxis  Hilbert transformations in the corresponding
    spaces.

To construct other classes of  $B$-hyperbolic
 transmutation operators, one can use the transformation with the kernel Neumann function
  $Y_{\nu}(z)$ instead of the  Fourier--Bessel
  transformation $F_{\nu}$.

Operators of translations with respect to parameters of Bessel
 operators have important applications in the theory of singular
 differential equations (see \cite{FJSS, KarST, KarST2, KarST3, KarST4}.

    \subsection{$B$-elliptic  transmutation operators}\label{sec16.2}

 $B$-elliptic  transmutation operators satisfy the relation
\begin{equation}\label{4415}{T B_{\nu} = - D^2 T.}\end{equation}
 This unusual class of transmutation operators establishes a relation between
 solutions of
$B$-elliptic and $B$-hyperbolic
    differential equations, expressing their solutions via each
    other.

To construct such operators, in the previous factorizations, we
change the Fourier
 sine-transformation and cosine-transformation
 for the Laplace transformation or change the Fourier--Bessel
 transformation to one of the following transformations with
 the Macdonald function and Neumann function:
\begin{equation}
  (K_{\nu} f) (t) = \frac{1}{t^{\nu}} \int\limits_0^{\infty} y^{\nu+1} K_{\nu}(t y) f(y)
\, dy  \label{4416}
\end{equation}
    and
\begin{equation} (Y_{\nu} f) (t) = \frac{1}{t^{\nu}} \int\limits_0^{\infty}
y^{\nu+1} Y_{\nu}(t y) f(y) \, dy. \label{4417}
\end{equation}
    Also, it is possible to use the transformation with the Struwe
    function in the kernel.

\begin{theorem}
    Let $|\Re \nu| + \Re (\alpha+\nu) < 1.$ Then
    $$
    \lr{A_{\nu}^{\alpha} f} (x) = F_c^{-1} t^{- \alpha} K_{\nu} f
    $$
 is a $B$-elliptic
 operator satisfying relation \eqref{4415}.
  The integral representation
    $$
    \begin{gathered}
    (A^{\alpha}_{\nu} f) (x) = \frac{\pi \Gamma(1-\alpha)}{4 \sin \frac{\pi}{2}
      (1-\alpha-2 \nu)} \int\limits_0^{\infty} y^{\nu+1} (x^2+y^2)^{\frac{\alpha+\nu-1}{2}}
     \\
    \times \left[ P_{-\alpha-\nu}^{-\nu} \lr{\frac{x}{\sqrt{x^2+y^2}}} +
    P_{-\alpha-\nu}^{-\nu} \lr{ - \frac{x}{\sqrt{x^2+y^2}}}  \right] f(y) \, dy.
    \end{gathered}
    $$
    holds for it.
\end{theorem}

Define an operator satisfying relation \eqref{4415} as follows:
$$
(C^{\alpha}_{\nu} f) (x) = L (t^{-\alpha} F_{\nu} f),
$$
    where  $L$ is the Laplace transformation.

\begin{theorem}
    Let $\Re \alpha < 1.$ Then the integral
representation
    $$
    (A^{\alpha}_{\nu} f) (x) =
     \Gamma(1-\alpha) \int\limits_0^{\infty} y^{\nu+1} (x^2+y^2)^{\frac{\alpha+\nu-1}{2}}
      P_{-\alpha-\nu}^{-\nu} \lr{\frac{x}{\sqrt{x^2+y^2}}} f(y) \, dy
    $$
    holds.
\end{theorem}

For $|\Re \nu| + \Re (\alpha+\nu) < 1,$ a similar representation
 holds for the operator $L$ as well:
$$
L (t^{-\alpha} Y_{\nu} f)(x) =
 - \frac{2}{\pi} \Gamma(1-\alpha) \int\limits_0^{\infty} y^{\nu+1}
  (x^2+y^2)^{\frac{\alpha+\nu-1}{2}} Q_{-\alpha-\nu}^{-\nu}
   \lr{\frac{x}{\sqrt{x^2+y^2}}} f(y) \, dy.
$$
    Also, similar expressions are derived for a broader class of $B$-elliptic
    transmutation operators intertwining
$B_{\nu}$ and $\lr{-B_{\mu}}.$

Consider the basic above operators $A_{\nu}^{\alpha}$ and
$C_{\nu}^{\alpha}$ for $\nu = \pm \dfrac{1}{2}.$
 These operators are defined by the relations
$$
\lr{A^{\beta} f}(x) = \lr{F^{-1}_{\left\{ \begin{matrix} s \\
c \end{matrix} \right\}} t^{\beta} L} f ~\textrm{ and }~
\lr{C^{\beta}
f}(x) = \lr{L  t^{-\beta} F_{\left\{ \begin{matrix} s \\
c
\end{matrix} \right\}} } f,
$$
    where  $L$ is the Laplace transformation. They intertwine $D^2$ and $-D^2.$
    In this case, it is easier to compute the kernels of the
    integral operators directly. This leads to the relation
$$
\lr{C^{\beta} f}(x) = \sqrt{\frac{2}{\pi}} \Gamma (1-\beta)
 \int\limits_0^{\infty} \frac{ f(y)}{\lr{x^2+y^2}^{\frac{1-\beta}{2}}}
  \left\{ \begin{matrix} \sin \\
\cos \end{matrix} \right\}  \Big[(1-\beta) {\rm atan}
\frac{y}{x}\Big] dy,
$$
    where $\Re \beta < 1 + \delta$ and $\delta=\left\{ \begin{matrix} 1 \\
 0 \end{matrix} \right\}$, and
$$
\lr{A^{\beta} f}(x) = \sqrt{\frac{2}{\pi}} \Gamma (1+\beta)
 \int\limits_0^{\infty} \frac{ f(y)}{\lr{x^2+y^2}^{\frac{1+\beta}{2}}}
  \left\{ \begin{matrix} \sin \\
\cos \end{matrix} \right\}  \Big[(1+\beta) {\rm atan}
\frac{x}{y}\Big] dy,
$$
 where $\Re \beta > - \delta - 1$ and $\delta$ is defined above.
 In particular, for $\beta=0$, we obtain the following pair of transmutation operators,
  related to the half-space
  Poisson integrals:
$$
\lr{C^{0} f}(x) = \sqrt{\frac{2}{\pi}} \int\limits_0^{\infty} \frac{y f(y)}{x^2+y^2} \, dy,
\qquad
\lr{A^{0} f}(x) = \sqrt{\frac{2}{\pi}} \int\limits_0^{\infty} \frac{x f(y)}{x^2+y^2}\, dy.
$$
    These operators and the operators  $C^{\beta} $ and $A^{\beta}$ for special positive
    integer values of $\beta$ are constructed in \cite{Car1}.

\subsection{$B$-parabolic transmutation operators}\label{sec16.3}

This unusual class of transmutation operators allows one to
express solutions of parabolic equations via solutions of
hyperbolic ones and vise versa.

Introduce the integral transformations
$$
\lr{F'_c f} (x) =  \lr{F_c f} (\sqrt{x}), ~  \lr{F'_s f} (x) =  \lr{F_s f} (\sqrt{x}),
 \textrm{ and }
\lr{P f} (x) = \lr{L \varphi (t) F'_{\scriptsize\left\{ \begin{matrix} s \\
c \end{matrix} \right\}}}(x).
$$
    Then, on compactly supported functions, the operator $P$ intertwines the second and
    first derivative as follows:
$$
P D^2 f = D P f.
$$
    Thus, this operator is parabolic in the Carroll sense.

\subsection{Lowndes-type operators of translations with respect to the spectral parameter}
    \label{sec16.4}

Operators of this type arise once solutions of Helmholtz equations
    are expressed via harmonic functions.
    Their study is started by Vekua and Erd\'elyi and is continued by Lowndes.
    This is the reason to call them the
  Vekua--Erd\'elyi--Lowndes
  transmutation operators (see \cite{S66,S6,S46,S59,S125}).
  Similar operators are applied to resolve singular differential equations
  (see, e.\,g., \cite{KarST}).

    Consider the operator
$$
T_1 = F^{-1}_{\nu} \lr{ \varphi (t) F'_{\mu}},
$$
     where the following transformation with parameter $\lambda$
     is introduced:
$$
\lr{F'_{\mu} f} (t) = \frac{1}{t^{\nu}} \int\limits_0^{\infty} y^{\nu+1}
 J_{\nu}(y \sqrt{t^2+\lambda^2}) f(y) \, dy.
$$
    The operator $T_1$ satisfies the relation
\begin{equation}\label{4420}{T_1 B_{\mu} = (B_{\nu} - \lambda^2) T_1.}\end{equation}

\begin{theorem} If $-1< \Re \nu<  \Re \mu$ and
  $\varphi (x) = x^{\mu} (x^2+\lambda^2)^{-\frac{\mu}{2}},$
  then the integral representation
    $$
    \lr{T_1 f} (x) =
     \lambda^{1+\nu-\mu} \int\limits_0^{\infty} y (y^2-x^2)^{\frac{\mu-\nu-1}{2}}
      J_{\nu-\mu-1}(\lambda \sqrt{y^2-x^2}) f(y) \, dy
    $$
    holds.
\end{theorem}

A number of other operators obtained in the factorized form is
provided below.
    %\medskip

\noindent
 (a) Let $\nu=1$ and $\varphi (x) = x^{\mu-2}
(x^2+\lambda^2)^{\frac{\mu}{2}}.$

Then the operator $T_2= F_1^{-1} \varphi (t) F'_{\mu}$
 satisfying the relation
$$
T_2 B_{\mu} = (B_1 -\lambda^2)\, T_2
$$
 for $\Re \lambda > 0$  and $\Re \mu > -1$ can be represented as
 follows:
$$
(T_2 f)(x)= \frac{1}{x^2 \lambda^2} \int\limits_0^{\infty}
y^{\nu+1} J_{\mu}(\lambda y) f(y)  \, dy -
 \frac{1}{x^2 \lambda^2} \int\limits_x^{\infty} y (y^2-x^2)^{\frac{\mu}{2}}
 J_{\mu}(\lambda \sqrt{y^2-x^2}) f(y) \, dy.
$$
    (b) Let  $\varphi (x) = x^{\mu} (x^{\frac{\mu}{2}}+\lambda^2)^2,$
$\Re \lambda > 0,$ $T_3 = F^{-1}_{\nu} \varphi F'_{\mu},$ and $-1<
\Re \nu< - \Re \mu.$ Then
$$
    \begin{gathered}
(T_3 f) (x) =   \frac{2 \sin (\pi \mu)}{\pi} \lambda^{\mu+\nu+1}
\int\limits_0^x y^{2 \mu + 1} (x^2-y^2)^{-\frac{\mu+\nu+1}{2}}
K_{\mu+\nu+1}(\lambda \sqrt{x^2-y^2}) f(y) \, dy
    \\
+ \lambda^{\nu+\mu+1} \int\limits_x^{\infty} y^{2 \mu + 1} (y^2-x^2)^{-\frac{\mu+\nu+1}{2}}
  \left[ \sin (\pi \nu) Y_{\mu+\nu+1}(\lambda \sqrt{y^2-x^2}) -
 \cos (\pi \nu) J_{\nu+\mu+1}(\lambda \sqrt{y^2-x^2}) \right]  f(y) \, dy.
    \end{gathered}
$$
    (c) Let $\varphi (x) = x^{\mu-1}/(x^2+\lambda^2)
     (x^2+\frac{\lambda}{2}^2)^{\frac{\mu}{2}}.$\\
Assign $T_4 = Y^{-1}_{\nu} \varphi F'_{\mu}$ and assume that
 $\Re
\lambda
> 0$ and $-\dfrac{1}{2}< \Re \nu< 3+ \Re \mu.$
    Then the following relation holds:
$$
(T_4 f) (x) =  - \frac{\lambda^{\nu-\mu-1}}{2^{- \frac{\mu}{2}}}
\cdot \frac{K_{\nu} (\lambda x)}{x^{\nu}}  \int\limits_0^{\infty}
y^{\mu+1} J_{\mu}\lr{\frac{\lambda y}{\sqrt{2}} } f(y) \, dy.
$$
    The provided examples demonstrate that it is important to have a freedom to select
    the function $\varphi.$

Operators  $T$ satisfying the relation
\begin{equation}\label{4421}{T B_{\mu} = (B_{\mu} + \lambda^2) \,
T}\end{equation}
    are constructed in the same way.
    For example, there exists a function   $\varphi$  such that
    the operator
$$
T_5 = (F'_{\mu})^{-1} \varphi F_{\nu}
$$
    is the  Lowndes operator
$$
(T_4 f)(x) = \frac{\lambda^{\nu-\mu+1}}{x^{2 \mu}}
\int\limits_0^{\infty} y^{2 \nu +1}
(x^2-y^2)^{\frac{\mu-\nu-1}{2}} J_{\mu-\nu-1}(\lambda
\sqrt{y^2-x^2}) f(y) \, dy
$$
for it.

Operators contained in relations \eqref{4420}-\eqref{4421}
    are linear. Therefore, the said relations can be represented as
    follows:
$$
T_1 (B_{\mu} + \lambda^2) = B_{\nu} T \textrm{ and } T (B_{\nu} -
\lambda^2) = B_{\mu} T.
$$
    In all operators, the substitution $\lambda\to i \lambda$ can be
    justified. For the most general transmutation operators,
    the relation
     \begin{equation}\label{22}{T (B_{\nu}+\alpha)
= (B_{\mu} + \beta) T}\end{equation}
 is equivalent for the relation pair
$$
T(B_{\nu}+\alpha-\beta) = B_{\mu}T,\quad  (B_{\mu}+\beta-\alpha)T
=T B_{\nu}
$$
    considered above.

Operators satisfying the previous relations can be obtained
intermediately  as well.
    Lowndes-type $B$-elliptic
    operators satisfying the relation
$$
T (B_{\nu}+\lambda) = (-B_{\mu} + \beta)\, T
$$
    are considered in the same way.

Selecting the values $\nu = \mu = - \dfrac{1}{2}$ of the
parameters, we obtain operators intertwining  $D^2$ and $D^2 \pm
\lambda^2.$

Also, the composition method can be applied if the quadratic
    (fractional) Fourier or Hankel transformation is taken as the
    integral transformation.
    For this purpose, relations \eqref{Ht^2D}-\eqref{x^2HD}
    are the most interesting.
    According to the general scheme of the composition method, we can construct
    transmutation operators intertwining the differential operators $D^2$
    and
$$
\left(\sin^{2}{\frac{\alpha}{2}}L_{\nu} -
\frac{1}{2}i\sin{\alpha}\left(XD+DX\right)  -
X^{2}\cos^{2}{\frac{\alpha}{2}} \right),
$$
    where
$$
L_{\nu}= -\frac{1}{4}D^2-\frac{\nu^2-1/4}{x^2} + \frac{1}{4}x^2 - \frac{\nu+1}{2}
$$
 with arbitrary parameters $\alpha$ and $ \nu.$
 To do that, one has to assign
$$
A=\left(\sin^{2}{\frac{\alpha}{2}}L_{\nu} -
\frac{1}{2}i\sin{\alpha}\left(XD+DX\right)  -
X^{2}\cos^{2}{\frac{\alpha}{2}}\right),\quad B=D^2,
$$
$$
F(A)=H_\nu^\alpha, F(B)=F_c, \textrm{ and } g(t)=-t^2
$$
    (in our terms),
where $H_\nu^\alpha$ is the quadratic Fourier--Fresnel
 transformation and $F_c$ is the Fourier
    cosine-transformation.

Also, the composition method can be successfully applied to the
constructing of fractional powers of the  Bessel operator.
  Solutions of numerous integrodifferential equations can be obtained by means of it
  (see \cite{S46, SitDis}).

\chapter[Applications of Transmutation Operator Method to Solution
Estimating\\
 for
 Variable-Coefficient Differential Equations and Landis Problem]{Applications
  of Transmutation Operator Method\\
   to Solution Estimating for
 Variable-Coefficient Differential Equations and Landis Problem}\label{ch7}

    % \chaptermarknum{Приложения метода операторов преобразования}

\section[Applications of Transmutation Operator Method to Perturbed
 Bessel Equations\\
  with Variable Potentials]{Applications of Transmutation Operator Method\\
  to Perturbed
 Bessel Equations with Variable Potentials}\label{sec17}
%\sectionmarknum{\sП\sр\sи\sл\sо\sж\sе\sн\sи\sя\s{ }\sм\sе\sт\sо\sд\sа\s{ }\sо\sп\sе\sр\sа\sт\sо\sр\sо\sв\s{
%}\sп\sр\sе\sо\sб\sр\sа\sз\sо\sв\sа\sн\sи\sя\s{ }\sд\sл\sя\s{
%}\sв\sо\sз\sм\sу\sщ\sё\sн\sн\sо\sг\sо\s{
%}\sу\sр\sа\sв\sн\sе\sн\sи\sя\s{ }\sБ\sе\sс\sс\sе\sл\sя\s{ }\sс\s{
%}\sп\sе\sр\sе\sм\sе\sн\sн\sы\sм{
%}\sп\sо\sт\sе\sн\sц\sи\sа\sл\sо\sм}

    Consider the problem to construct an integral relation for
    solutions of the differential equation
\begin{equation}\label{4.2.1}
B_{\alpha} g(x) - q(x) g(x)=\lambda^2 g(x)
\end{equation}
 with a given asymptotic behavior, where $B_\alpha$ is the  Bessel operator;
 in this section, it is convenient to define it in the form
\begin{equation}\label{4.2.2}
B_{\alpha}g =g''(x)+\frac{2\alpha}{x}g'(x),~ \alpha>0.
\end{equation}
    Note that various form to represent constants in the
 Bessel operator are used throughout the monograph.
 This depends on the problem to be solved and is caused by easier
 expressions of the obtained results.
 No misunderstanding is expected because problems considered at different sections of
  the monograph are disjoint.

The above problem is solved by means of the transmutation operator
method. It suffices to construct a pair of mutually inverse
 transmutation operators such that the first one
\begin{equation}\label{4.2.3}
S_\alpha h(x)=h(x)+\int\limits_x^{\infty}S(x,t)h(t)\,dt
\end{equation}
    intertwines the operators $B_\alpha - q(x)$ and $B_\alpha$ by
    means of the relation
\begin{equation}\label{4.2.4}
S_{\alpha}(B_{\alpha}-q(x))h=B_{\alpha}S_{\alpha}h,
\end{equation}
 while second one is inverse to the first one, is constructed as an integral operator
 with kernel $P(x,t)$,
\begin{equation}\label{4.2.3*}
P_\alpha h(x)=h(x)+\int\limits_x^{\infty}P(x,t)h(t)\,dt,
\end{equation}
    and acts according to the relation
$$
P_\alpha B_{\alpha}h=(B_{\alpha}-q(x))P_\alpha h
$$
 on functions $h \in C^2(0, \infty).$

As a result, on solutions of Eq. \eqref{4.2.1}, the function
$S_{\alpha}u=v$ is expressed via solutions of the unperturbed
equation, which is Eq. \eqref{4.2.1} without the potential term.
    In fact, it  is expressed via Bessel functions and the function $u=P_\alpha v$
    is satisfies the original perturbed equation given by \eqref{4.2.1}.
    The integral representation of the solution, given by \eqref{4.2.3*},
    describes the kernel $P(x,t)$ explicitly.
    This technique reflects one of the main utilities of transmutation operators:
    they express complicated differential equations via simpler ones (see above).
    Note that, since linear  transmutation operators are used, it follows that the same pair
    of mutually inverse transmutation operators yields representations both for solutions
    of Eq. \eqref{4.2.1} with a spectral parameter and for solutions of the (simpler)
  case of the   homogeneous equation
$$
B_{\alpha} h(x) - q(x) h(x)=0.
$$
    If we have to represent solutions of the perturbed equation, i.\,e., Eq. \eqref{4.2.1},
    then the constructing of the direct transmutation operator can be
    skipped:
   one can immediately pass to the constructing of an integral
   representation of kind \eqref{4.2.3*} for the desired solution.

In \cite{Sta1,Sta2}, an original technique to solve such problems
is developed. Using it, one can consider singular potentials
satisfying the most exact estimate at the origin: $|q(x)| \leq c
x^{- 3/2+\varepsilon},$ $\varepsilon
> 0$ where $\alpha$ is integer.
 This technique based on the application of generalized Paley--Wiener
 theorems is broadly developed and recognized.
 The case where $q$ is continuous and $\alpha>0$ is considered in
 \cite{Soh1, Soh2, Soh3, Soh4, Volk}. Stashevskaya and Volk construct
  Povzner-type transmutation operators integrating over a finite interval, while
  Sokhin constructs
 Levin-type transmutation operators integrating over an infinite one.
    Below, we propose a new modified method combining both these
    approaches.

There many works obtaining representations of solutions of the
perturbed Bessel equation given by \eqref{4.2.1}-\eqref{4.2.2}.
    We note \cite{CFH, FH, CKT1,CKT2, Krav1, Krav2,
Krav3, Krav4, Krav5, Krav6, Krav7, Krav8, Krav9, Krav10}, where
solutions are sought in the form of special-kind
 series (also, see the critical review \cite{Hol} of a number of results about this problem).

In many problems of mathematics and physics, one has to consider
strongly singular potentials, e.\,g., potentials admitting
arbitrary power singularities at the origin.
    In the present work, we formulate results on integral representations of solutions
    of equations with such singular potentials.
 The only requirement for the potential is to be majorized by a function summable ar infinity.
 In particular, the class of admissible potentials includes the singular potential $q=x^{-2},$
 the strongly singular potential with the power singularity $q=x^{-2-\varepsilon},$
$\varepsilon > 0,$
 Ukawa-type potentials  $q=e^{-\alpha x}/x,$ Bargmann and Batman--Chadan
 potentials (see \cite{ShSa}), and a  number of others.
 No additional assumptions (such as a fast oscillation at the origin or a constant sign)
 are imposed on the function $q(x)$, which allows one to use the same method to study both
 attracting and repelling potentials.

    Special-kind transmutation operators constructed in the present book differ from the ones
    known before. Earlier, only the cases where the main integral equation for the kernel
    of the transmutation operator has the same integrations limits
    (either  $[0;a]$ or $[a;\infty]$ in both integrals).
    Here, we show that different limits can be considered.
    This provides a possibility to cover more broad class of
    potentials with singularities at the origin.
    Also, we improve the classical Levitan scheme of \cite{Lev7}, expressing the Green
    function (used in the proofs) via the Legendre function (instead of more general Gauss
    hypergeometric function) depending on less amount of parameters; this allows one to avoid
    undetermined constant in estimates of earlier papers.

Due to the volume limit for the present book, we provide only the
problem formulation and a list of main results and corollaries
(without proofs) in this section; for the detailed presentation,
see \cite{S8,S63,S46,S19,S43,S4}.

\subsection{The resolving of the main integral equation for kernels
 of transmutation operators}\label{sec17.1}

    Introduce new variables and functions as follows:
$$
\xi=\frac{t+x}{2}, ~ \eta=\frac{t-x}{2}, ~ \xi \geq \eta > 0,
$$
\begin{equation}\label{4.2.5}
K(x, t)= \left(\frac{x}{t}\right)^\alpha P(x, t),
    \textrm{ and } u(\xi, \eta)= K(\xi-\eta, \xi+\eta).
\end{equation}
    Introduce the notation $\nu=\alpha-1.$ Thus, to justify representation
    \eqref{4.2.3*} for solutions of Eq. \eqref{4.2.1}, it suffices
    to define the function $u(\xi, \eta ).$
    It is known from \cite{Soh1, Soh2, Soh3, Soh4} that if there exists a twice
    continuously differentiable solution $u(\xi, \eta)$ of the integral equation
$$
u(\xi, \eta)=-\frac{1}{2}\int\limits_\xi^{\infty}R_\nu(s, 0; \xi,
\eta) q(s) \, ds - \int\limits_\xi^{\infty} ds
\int\limits_0^{\eta} q(s+\tau) R_{\nu}(s, \tau; \xi, \eta)  u(s,
\tau) \, d \tau
$$
    under the assumption that $ 0< \tau < \eta < \xi < s,$ then the desired function
 $P(x,t)$ is defined by relations \eqref{4.2.5} via this solution $u(\xi, \eta).$
 Here, $R_{\nu}=R_{\alpha-1}$ is the Riemann  function arising under the resolving of
 a Goursat problem for the singular hyperbolic equation
$$
\frac{\partial^2 u(\xi, \eta)}{\partial \xi \partial \eta}+
\frac{4 \alpha(\alpha-1) \xi \eta}{(\xi^2-\eta^2)^2} u(\xi,
\eta)=q(\xi+\eta)u(\xi, \eta).
$$
    This function is known in its explicit form (see \cite{Soh1, Soh2, Soh3,
Soh4}): it is expressed via the  Gauss
    hypergeometric function $_2{F_1}$ as follows:
\begin{equation}\label{4.2.6}
R_\nu=\left(\frac{s^2-\eta^2}{s^2-\tau^2}\cdot
\frac{\xi^2-\tau^2}{\xi^2-\eta^2}\right)^{\nu}{_2F_1} \left(-\nu,
-\nu; 1; \frac{s^2-\xi^2}{s^2-\eta^2}\cdot
\frac{\eta^2-\tau^2}{\xi^2-\tau^2}\right).
\end{equation}
    This expression is simplified in \cite{S8}: it is proved that, in the considered case,
   the Riemann  function is expressed via the Legendre function as   follows:
\begin{equation}\label{4.2.7}
R_\nu (s, \tau, \xi, \eta)=P_\nu \left(\frac{1+A}{1-A}\right),~A
=\frac{\eta^2-\tau^2}{\xi^2-\tau^2}\cdot
\frac{s^2-\xi^2}{s^2-\eta^2}.
\end{equation}
 The main content of this section is the following  assertion.

\begin{theorem}
 Let $q(r)\in C^1 (0,\infty)$ and the conditions
    \begin{equation}\label{4.2.8}
    |q(s+\tau)|\leq |p(s)|, ~ \forall s, \forall \tau, ~ 0< \tau <s\
   \textrm{ and } \int\limits_\xi^\infty |p(t)| \, dt < \infty  \forall \xi>0
    \end{equation}
    be satisfied.
   Then there exists an integral representation of kind
   \eqref{4.2.3*} such that its kernel satisfies the estimate
    $$
    \gathered
    |P(r, t)| \leq \left(\frac{t}{r}\right)^ \alpha \frac {1}{2}
     \int\limits_{\frac{t+r}{2}}^\infty P_{\alpha-1}
     \left(\frac{y^2(t^2+r^2)-(t^2-r^2)}{2try^2}\right)|p(y)|\, dy
    \\
    \times\exp \left[ \left(\frac{t-r}{2}\right) \frac{1}{2}
     \int\limits_{\frac{t+r}{2}}^\infty P_{\alpha-1}
     \left(\frac{y^2(t^2+r^2)-(t^2-r^2)}{2try^2}\right)|p(y)|\, dy  \right].
    \endgathered
    $$
  The kernel $P(x,t)$ of the transmutation operator and the solution of Eq. \eqref{4.2.1}
  are twice   continuously differentiable on $(0,\infty)$ with respect to their independent
  variables.
\end{theorem}

The following classes of potentials satisfy Conditions
\eqref{4.2.8}. If $|q(x)|$ monotonously decreases, then one can
assign $p(x)=|q(x)|.$ If a potential has an arbitrary singularity
at the origin, increases for $0<x<M$ (e.\,g., the Coulomb
potential $q=-\dfrac{1}{x}$), and is truncated by zero at
infinity, i.\,e., $q(x)=0$ for $x>M,$ one can assign $p(x)=|q(M)|$
for $x<M$ and $p(x)=0$ for $x \geq M.$
 Potentials with the estimate $q(x+\tau) \leq c|q(x)|=|p(x)|$ satisfy Condition \eqref{4.2.8}
 as well.

In particular, the above conditions are satisfied by the following
potentials arising in applications: strongly singular potentials
with power singularities of the kind  $q(x)=x^{-2-\varepsilon},$
 various  Bargmann  potentials
$$
q_1 (x) = - \frac{e^{-ax}}{(1+\beta e^{-ax})^2}, ~ q_2 (x)=
\frac{c_2}{(1+c_3 x)^2}, \textrm{ and }  q_3 (x)=
\frac{c_4}{ch^2(c_5 x)},
$$
    and the Ukawa potentials
$$
q_4 (x) = - \frac{e^{-ax}}{x} \, \textrm{ and } \, q_5 (x) =
\int\limits_x^\infty e^{-at} \, d c(t)
$$
(see, e.\,g., \cite{ShSa}).

\begin{remark}
    In fact, the explicit form of the Riemann  function, given by \eqref{4.2.7}, is not
    required to prove the above theorem.
    Only the existence of the Riemann  function, its positivity, and a special
    monotonicity property are used.
    These facts are rather general. Therefore, the obtained results can be extended to
   a quite broad class of differential equations.
\end{remark}

The estimate from Theorem \ref{4.2.1} for general-kind
 potentials can be transformed into a rougher (though more visible) one.

\begin{theorem}
 Let assumptions of Theorem \ref{4.2.1} be satisfied.
  Then the kernel $P(x, t)$ of the transmutation operator
  satisfies the estimate
    $$
    |P(x, t)| \leq \frac{1}{2} \left(\frac{t}{x}\right)^{\alpha} P_{\alpha-1}
    \left(\frac{t^2+x^2}{2tx}\right) \int\limits_x^\infty |p(y)|\, dy\,  \exp
     \left[ \frac{1}{2} \left(\frac{t-x}{2}\right) P_{\alpha-1} \left(\frac{t^2+x^2}{2tx}
     \right) \int\limits_x^\infty |p(y)| \, dy \right].
    $$
\end{theorem}

Note that an exponential singularity as $x \to 0$ is admitted for
the kernel of the integral representation.

For the class of potentials with power singularities of the kind
\begin{equation}\label{4.2.16}
q(x)=x^{-(2\beta +1 )},~ \beta > 0,
\end{equation}
    the obtained estimates can be simplified such that the precision is preserved.
 The restriction for $\beta$ is caused by the summability condition at infinity.

\begin{theorem}
    For potentials of kind \eqref{4.2.16}, Theorem  \ref{4.2.1}
    holds with the estimate
    $$
    |P(x, t)|  \leq \left(\frac{t}{x}\right)^{\alpha}
      \frac{\Gamma(\beta)4^{\beta-1}}{(t^2-x^2)^\beta} \cdot P_{\alpha-1}^{- \beta}
       \left(\frac{t^2+x^2}{2tx}\right) \exp\left[\left(\frac{t-x}{x}\right)
        \frac{\Gamma(\beta)4^{\beta-1}}{(t^2-x^2)^\beta} P_{\alpha-1}^{- \beta}
         \left(\frac{t^2+x^2}{2tx}\right)   \right],
    $$
  where $P_\nu^\mu(\cdot)$ is the Legendre function, the value of
   $\beta$ is determined by \eqref{4.2.16}, and
    the value of $\alpha$  is determined by \eqref{4.2.2}.
\end{theorem}

 This estimate is obtained by means of a long chain of computations using the
  well-known Slater--Marichev
  theorem (see \cite{Marich1}); it helps to compute the needed integrals (in terms of
  hypergeometric functions) once they are reduced to the  Mellin convolution.

In \cite{S8}, a basic estimate of this kind is obtained for the
potential $q(x)=cx^{-2}$; for this potential, we have
$\beta=\dfrac{1}{2}.$ It follows from \cite{BE1} that the Legendre
 function $P_{\nu}^{-\frac{1}{2}}(z)$ can be expressed via elementary functions in such a case.
 Therefore, the corresponding estimate can be expressed via elementary functions as well.

Also, the obtained estimate can be simplified and expressed
  via elementary functions for potentials of the kind
$q(x)=x^{-(2 \beta + 1)},$ where the parameters are bound by the
relation $\beta=\alpha-1.$

\begin{corollary} If $\beta=\alpha-1,$ then the estimate from Theorem \ref{4.2.3}
takes the form
   \begin{multline}\label{4.2.20}
    |P(x, t)|  \leq \left(\frac{t}{x}\right)^{\beta+1} \frac{2^{\beta-2}}{\beta}
     \left[\frac{t^2+x^2}{2tx}\right]^\beta  \exp \left[ \left(\frac{t-x}{2}\right)
      \frac{2^{\beta-2}}{\beta} \left[\frac{t^2+x^2}{2tr}\right]^\beta \right]
      \\
    = \frac{1}{4 \beta} \frac{1}{x^{2 \beta +1}} (t^2+x^2)^{\beta} \exp
     \left[  \frac{2^{\beta-2}}{\beta} \left(\frac{t-x}{2}\right)
      \left(\frac{t^2+x^2}{2tx}\right)^\beta \right].
    \end{multline}
\end{corollary}

If $\alpha=0$ in relations \eqref{4.2.1}--\eqref{4.2.3*},
    then Theorem \ref{4.2.1} is reduced to the known estimates of the kernel of
    the integral representation of Jost solutions for the
Sturm--Liouville equation.

The above technique is completely extended to the problem to
construct nonclassical generalized translation operators.
    Actually, this problem is equivalent to the problem to express
    solutions of the equation
\begin{equation}\label{4.2.22}
B_{\alpha,x} u(x,y) - q(x) u(x,y)=B_{\beta,y} u(x,y)
\end{equation}
 via solutions of the unperturbed
    Euler--Poisson--Darboux
    equation (the wave equation in the regular case)
    under additional assumptions ensuring the well-posedness.
   Such representations follow just from the existence of transmutation operators.
   For the regular case (the case where $\alpha=\beta=0$), they
 are studied in \cite{Lev1, Lev2, Lev3} as corollaries from the generalized translation theory
   (see \cite{Mar9} as well).
   In \cite{Bor1,Bor2}, an interesting original technique is development
   to obtain such representations in the regular case as well.
 Results  of the  present work imply integral representations of a subclass of solutions
 of Eq. \eqref{4.2.22} in the general singular case (quite arbitrary potentials with
 singularities at the origin are admitted).
    The obtained estimates of solutions contain no undetermined constants,
 while integral equations such that the kernels of integral representations
 satisfy them are provided explicitly.

\section{Applications of Transmutation Operator Method to Landis Problem}\label{sec18}

In \cite{Lan}, the following problem is set: to prove that each
solution of the stationary Schr\"odinger equation with a bounded
potential
\begin{eqnarray}
& &  \Delta u(x) - q(x) u(x) =0,~x \in \R^n,~|x| \geq R_0 > 0,   \label{l1} \\
& & |q(x)| \leq \lambda^2,~ \lambda > 0,~u(x) \in C^2 \lr{|x| \geq
R_0}, \nonumber
\end{eqnarray}
    satisfying an estimate of the kind
$$
|u(x)| \leq \const\cdot \, e^{- (\lambda + \varepsilon)|x|},~ \varepsilon > 0,
$$
    is identically equal to zero.

In \cite{Mesh1,Mesh2}, the negative answer is given: the existence
of counterexamples of solutions (such that they are complex
functions) are proved. Also, it is proved that the answer is
positive if the estimate is strengthened is follows:
$$
|u(x)| \leq \const\cdot \, e^{- (\lambda +
\varepsilon)|x|^{4/3}},~ \varepsilon > 0;
$$
    indeed, there are no nonzero solutions satisfying it.
  The interest to these results is preserved and this direction is
  actively developed (see \cite{BuKe1,Ken1,Ken2,DKW1,Rossi}).
  The main question is the investigation of the Landis conjecture for real solutions;
  this problem is still unsolved.
  Therefore, it is proposed to call it the
\textit{Landis--Meshkov problem} and to formulate it as follows.

    %   \medskip

\textbf{Landis--Meshkov problem.} {\it Is it true that, for given
domain $D$ and positive functions $r(x)$ and $ s(x),$
    only the zero classical solution of the stationary
Schr\"odinger equation
\begin{equation}\label{LM1}
\Delta u(x) - q(x) u(x) =0,\ \ x \in \mathbb{R}^n, \ \ |q(x)|\leq r(x),
\end{equation}
 satisfies the estimate}
\begin{equation}\label{LM2}
|u(x)| \leq s(x)?
\end{equation}

    % \medskip

Meshkov results yield the negative answer for the case of complex
    solutions if  $D$ is the exterior of a disk, $q(x)=\lambda^2,$ and $s(x)=
 e^{- (\lambda + \varepsilon)|x|},$ $\varepsilon > 0,$ and
 the positive answer for the case of complex
    solutions if  $D$ is the exterior of a disk, $q(x)=\lambda^2,$ and
 $s(x)=e^{- (\lambda + \varepsilon)|x|^{4/3}},\ \varepsilon > 0.$
 For real solutions, no answers are known even in these special cases.

Below, we show that there are classes of potentials such that the
original Landis problem is solved positively for them (in spite of
the general negative answer found by Meshkov); moreover, this
positive answer is found for real solutions. To do that, we use a
special kind of transmutation operators (see \cite{S3, S71, S75}).

Further, this problem is solved for the case where the potential
depends only on one variable, i.\,e., $q(x)= q(x_i),$ where $1
\leq i \leq n;$ for  definiteness, we assume that $i=1.$
    For this case, generalize \eqref{l1} for the equation
\begin{equation}\label{l2}{\Delta u - q(x_1) u = 0,}\end{equation}
 where the potential $q(x_1)$ is bounded by an arbitrary
nondecreasing function.
 To solve the problem, we use transmutation operators reducing Eq. \eqref{l2} to the Laplace
 equation.

    %   \medskip

{\bf 1.}
    The assumptions of problem \eqref{l1} are satisfied in the
    half-space $\{x_1 \geq R_0\}$ and are invariant with respect to the substitution
     $z=x_1- R_0.$ Therefore, problem \eqref{l1} is considered in the
    half-space  $\{z \geq 0\}$ or  $\{x_1 \geq 0\}$ (to preserve the original notation
    of the spatial variable).
    We prove that the solution of problem \eqref{l1} is equal to zero  in the
    half-space $\{x_1 \geq 0\}$; then, due to the Calderon uniqueness theorem for the
    continuation (see \cite[Ch. 6, p. 14]{Miran}), such a solution is identically equal
    to zero in the whole space $\R^n.$

Let  $T\lr{\delta}$ be the set of functions satisfying the
following conditions  in the
    half-space  $\R^n_+=\{ x \in \R^n,$ $x_1 \geq 0\}$:
\begin{eqnarray}
& & u(x) \in C^2 \lr{\R^n_+}, \label{l3} \\
& & \left| u(x)\right| \leq c_1 \, e^{- \delta |x|},~\delta>0, \label{l4} \\
& &  \left| \frac{\pr u }{\pr x_1}\right| \leq c_2 \, e^{- \delta
|x|}.\label{l5}
\end{eqnarray}
    For functions from $T \lr{\lambda+\varepsilon}$, construct a transmutation operator
    of the kind
\begin{equation}\label{l6}{ S u(x) = u(x) +
\int\limits_{x_1}^{\infty} K (x_1, t) u (t,
x^1)\,dt}\end{equation}
 such that the relation
\begin{equation}\label{l7}{S \lr{\frac{\pr^2 u}{\pr x_1^2} -
q(x_1) u} = \frac{\pr^2 }{\pr x_1^2} S u,~|q(x_1)|\leq
\lambda^2}\end{equation}
  is satisfied (see \cite{S3, S71, S75}), where $(x_1, x^1):=(x_1, x_2, \dots, x_n).$
  Substituting expression \eqref{l6} in relation \eqref{l7}, we obtain
  the relations
\begin{equation}\label{l8}{ \frac{\pr^2 K}{\pr t^2} - \frac{\pr^2
K}{\pr x_1^2} = q(t) K,}\end{equation}
\begin{equation}\label{l9}{3 \frac{\pr K(x_1, x_1)}{\pr x_1} =
q(x_1),}\end{equation}
    and
\begin{equation}\label{l10}{\lim\limits_{t
\to \infty} K(x_1, t) \frac{\pr u (t, x^1)}{\pr t} -
\lim\limits_{t \to \infty} \frac{\pr  K(x_1, t)}{\pr t} u (t, x^1)
= 0.}\end{equation}
    Applying the standard substitutions $w=\dfrac{t+x_1}{2}$ and $v=\dfrac{t-x_1}{2},$
 we reduce system \eqref{l8}-\eqref{l9}
  (the fact that Condition  \eqref{l10} is satisfied on solutions of problem  \eqref{l1} is
  proved below) to the following (simpler) system:
\begin{eqnarray}
& & \frac{\pr^2 K}{\pr w \pr v} = q (w+v) K, \label{l11} \\
& & K(w, 0) =  \frac{1}{3} \int\limits_0^w q(s) \, ds. \label{l12}
\end{eqnarray}
    This system follows from the integral equation
\begin{equation}
K(w,v) = \frac{1}{3} \int\limits_0^w q(s) \, ds +
 \int\limits_0^w d\alpha  \int\limits_0^v q(\alpha+\beta) K (\alpha, \beta) \, d \beta,
 \quad |q| \leq \lambda^2,\quad w \geq v \geq 0. \label{l13}
\end{equation}
 Equation \eqref{l13} differs from the one usually used in the consideration of
  transmutation operators on infinite intervals: the integration
  domain is changed from the semiaxis $\lr{w, \infty}$ to the interval $\lr{0, w},$
  which implies the exponential growth of the kernel $K \lr{x_1, t}.$
  In the sequel, it is proved that such a kernel exists and the transmutation operator
  with this kernel, given by \eqref{l6}, is  defined on the set
   $T\lr{\lambda+\varepsilon}.$ Problem \eqref{l8}--\eqref{l10}
   can be reduced to nonequivalent integral equations because the Cauchy problem
  \eqref{l11}-\eqref{l12} is underdetermined.

Note that Eq. \eqref{l13} is to be solved in a wider region (free
of the restrictions imposed by the inequality $w \geq v$);
otherwise, the kernel under the integral signs is not defined. The
proof of the existence of a solution in this wider region is the
same as the proof provided below. Usually, this detail
 is left without attention during the proof of the existence of a solution of Eq. \eqref{l13}
 (this is noted by Borovskikh).

\begin{lemma}\label{L18.1}
    There exists a unique continuous solution of Eq. \eqref{l13},
 satisfying the inequality
\begin{equation}\label{l14}{|K (w, v)| \leq \frac{\lambda}{3}  \sqrt{\frac{w}{v}}\,
I_1 \lr{2\lambda \, \sqrt{wv}},}\end{equation}
 where $I_1 (x)$ is the modified Bessel function.
 For $q(x_1) \equiv \lambda^2,$ which is an admissible potential,  \eqref{l14}
 becomes an identity.
\end{lemma}

\begin{remark}
    In the sequel, $c$ denotes each absolute positive constant such that its
    value plays no role.
\end{remark}

\begin{proof}
    Introduce the notation
$$
K_0 (w, v) = \frac{1}{3} \int\limits_0^w q(s) \, ds
$$
    and
$$
P K (w, v) = \int\limits_0^w d\alpha  \int\limits_0^v
q(\alpha+\beta) K (\alpha+\beta) \, d \beta.
$$
    Then Eq. \eqref{l13} takes the form $K = K_0 +P K.$ We look
    for its solution represented by the Neumann series
\begin{equation}\label{l15}{K = K_0 +P K_0 + P^2 K_0 + \dots }\end{equation}
   Taking into account the condition that $|q(x_1)| \leq \lambda^2$,
    we obtain the following estimate for terms of series \eqref{l15}:
\begin{equation}\label{l16}{\left| P^n K_0 (w_0 v) \right| \leq
\frac{1}{3} \,  \lr{\lambda^2}^{n+1} \frac{w^{n+1}}{(n+1)!}
\frac{v^n}{n!},~n=0, 1, 2, \dots}\end{equation}
  This and the representation of the function $I_1 (x)$ by the
  series
$$
    I_1 (x) = \sum\limits_{k=0}^{\infty} \frac{(x/2)^{2 k+1}}{k!(k+1)!}
$$
 imply inequality \eqref{l14}.
    Estimate \eqref{l14} is exact because, for $q(x_1) \equiv
\lambda^2,$ inequalities \eqref{l16} become identities for all
nonnegative integers  $n,$ which completes the proof of the lemma.
\end{proof}

\begin{lemma}\label{L18.2}
    In terms of variables  $x_1$ and $t,$ the estimate
$$
\left| K(x_1, t) \right| \leq c \, t \,e^{\lambda t }
$$
    holds.
\end{lemma}

\begin{proof}
  Consider the inequality
$$
\left| \frac{1}{x} I_1 (x) \right| \leq c \, e^x,~x \geq 0.
$$
    To verify it, one has to consider the cases where
     (a) $x \geq 1$ and (b) $0 \leq x  \leq 1$
     and to use the known asymptotical behavior of the function $I_1 (x)$ as $x \to \infty$
     and as $x \to +0$ (see \cite{BE2}). Combining this with the
     obvious inequalities
$$
 \frac{x_1+t}{2} \leq t \textrm{ and } 2 \sqrt{w v} = \sqrt{t^2 -x_1^2} \leq t,
$$
 we deduce the assertion of the lemma from estimate \eqref{l14}.
\end{proof}
    It follows from this lemma that expression \eqref{l6} is defined on functions from
     $T \lr{\lambda+\varepsilon}.$ To prove that expression \eqref{l6}
     defines a transmutation operator on $T\lr{\lambda+\varepsilon}$ indeed, it remains
     to verify relation \eqref{l10}. The belonging  $u(x) \in T
\lr{\lambda+\varepsilon}$ and the lemma imply the inequality
$$
\lim\limits_{t \to \infty} K (x_1, t) \frac{\pr u(t, x^1)}{\pr t} = 0.
$$
    Hence, it remains to prove that if $u(x) \in T\lr{\lambda+\varepsilon},$
    then
$$
\lim\limits_{t \to \infty}  \frac{\pr K (x_1, t)}{\pr t}  u(x_1, t) = 0.
$$
    The last relation follows from the estimate
\begin{equation}\label{l17}{\left|  \frac{\pr K (x_1, t)}{\pr t}
\right| \leq c\, t\, e^{\lambda t }. }\end{equation}
  To prove inequality \eqref{l17}, one has to pass to the variables
$w$ and $v$ and, using the estimates found for the kernel $K(x_1,
t)$, estimate the derivatives $\dfrac{\pr K}{\pr w}$  and
$\dfrac{\pr K}{\pr v},$ differentiating Eq. \eqref{l13}. Since
$$
\frac{\pr K}{\pr t} = \frac{1}{2} \lr{\frac{\pr K}{\pr
w}+\frac{\pr K}{\pr v}},
$$
 we arrive at \eqref{l17}.

    %   \medskip

{\bf 2.} Let us show that each solution of problem  \eqref{l1}
    belongs to  $T \lr{\lambda+\varepsilon}$ and, therefore, operator \eqref{l6}
    is defined on such solutions. To do that, one has to check whether Condition \eqref{l5}
    is satisfied.

\begin{lemma}\label{L18.3}
 Let a function $u(x) \in C^2 \lr{|x| \geq R_0}$ satisfy problem \eqref{l1}.
 Then there exists a positive constant $c$ such that
$$
\left| \frac{\pr u}{\pr x_1} \right| \leq c \, e^{-(\lambda+\varepsilon)|x|}.
$$
\end{lemma}

\begin{proof}
    By virtue of the Schauder a priori estimates, the following inequality holds
    in the closed unit ball $B(x, 1)$ centered at  $x,$ $|x| \geq R_0+1$
(see \cite[Th. 33, II]{Miz}):
$$
u_1 \leq c \lr{{u_{1, \lambda_1}}^{\frac{1}{1+\lambda_1}} \cdot
{u_{0}}^{\frac{\lambda_1}{\lambda_1+1}} + u_0 },
$$
    where
$$
u_0= \| u(x) \|_{C^0 \lr{B(x, 1)}},~u_1 = \| u(x) \|_{C^1 \lr{B(x,
1)}},
$$
and $u_{1, \lambda_1}$ is the sum of H\"older coefficients of the
function   $u(x)$ and its first-order
  derivatives  $\dfrac{\pr u}{\pr x_i},$ $1 \leq  i \leq n.$
  This implies the inequality
\begin{equation}\label{l18}{\left|\frac{\pr u(x)}{\pr x_1}\right|
\leq c \lr{{u_{1, \lambda_1}}^{\frac{1}{1+\lambda_1}} \cdot
{u_{0}}^{\frac{\lambda_1}{\lambda_1+1}} + u_0 }.}\end{equation}
    Since all the assumptions of \cite[Assertion 33,
V]{Miz} are satisfied, it follows that the constant  $c$ in
estimate \eqref{l18} does not depend on $x.$

 Morrey results (see \cite[Th. 39, IV]{Miran}) yield the following
 estimate for  $u_{1, \lambda_1}$:
\begin{equation}\label{l19}{ u_{1, \lambda_1} \leq c \left[ \| u
\|_{L_2 \lr{B(x, 1)}}+ \|q\, u \|_{L_2 \lr{B(x,
1)}}\right],}\end{equation}
 where the constant does not depend on $x.$ By assumption, we have the inequality
  $|q(x_1)| \leq \lambda^2.$
  Hence, applying the mean-value
  theorem to \eqref{l19}, we obtain the inequality
$$
u_{1, \lambda_1}  \leq c \lr{\int\limits_{B(x, 1)} |u(y)|^2 \,
dy}^{1/2} \leq {c^1}\, e^{- \lr{\lambda+\varepsilon}|x|}.
$$
    Substituting this inequality in \eqref{l18}, we obtain that
$$
\left| \frac{\pr u}{\pr x_1} \right| \leq c \left[ \lr{ e^{-
\lr{\lambda+\varepsilon}|x|}}^{\frac{1}{1+\lambda}+\frac{\lambda}{1+\lambda}}
+
 e^{- \lr{\lambda+\varepsilon}|x|}\right] \leq c \, e^{- \lr{\lambda+\varepsilon}|x|}.
$$
  Thus, the claimed inequality is proved provided that $|x_1| \geq
R_0 +1.$ Since the set  $R_0 \leq |x| \leq R_{0}+1$ is a compactum
in  $\R^n,$ it follows that this inequality holds for  $|x| \geq
R_0$ as well, which completes the proof of Lemma \ref{L18.3}.
\end{proof}
    Again, change the coordinates as follows: $z=x_1-R_0.$
    We obtain the validity of Lemma \ref{L18.3} in the half-space
     $x_1\geq 0$ ($z$ is redenoted by $x_1$).

    %   \medskip

{\bf 3.}
    Apply the operator $S$ to Eq. \eqref{l2}.
    Identity \eqref{l7} and the fact that $S$ is commutative with the derivatives
$\dfrac{\pr^2 u}{\pr x_i^2},$ $2 \leq i \leq n,$ imply that
$$
S \lr{\Delta u - q(x_1) u} = \Delta S u = 0
$$
    in the half-space
 $\R_{+}^n$.

Denote the function  $Su$ by $v.$ It follows from \eqref{l6} and
\eqref{l13} that if $u(x) \in C^2 \lr{\R_{+}^n}$ and $q(x) \in C
\lr{\R_{+}^n},$ then $v(x) \in C^2 \lr{\R_{+}^n}.$ Let us show
that $v(x)$ exponentially decays in $\R_{+}^n$ as $|x| \to \infty$
and, therefore, is equal to zero.

\begin{lemma}
    Let $u(x) \in T \lr{\lambda+ \varepsilon}.$ Then
$$
|v|= \left| S u \right|  \leq c\, |x| \, e^{-\varepsilon |x|},~\varepsilon>0,
$$
  for each   $x \in \R_{+}^n$.
\end{lemma}

\begin{proof}
    From representation \eqref{l6} and Lemma \ref{L18.3}, we
    obtain that
$$
\left| S u \right|  \leq |u(x)| + \int\limits_{x_1}^{\infty} t \,
e^{\lambda t} c\, e^{-(\lambda+\varepsilon) \sqrt{t^2+|x^1|^2}} \,
dt \leq
 c \lr{ e^{-(\lambda + \varepsilon) |x|} +
\int\limits_{x_1}^{\infty} t \, e^{-(\lambda+\varepsilon)
\sqrt{t^2+|x^1|^2}} \, dt}.
$$
  Change the variables according to the relation $y =\sqrt{t^2 + |x^1|^2}$.
  Then integrate by parts and obtain the claimed estimate.
\end{proof}
  Thus, $v(x)=0$ in $\R_{+}^n.$ On  $T \lr{\lambda+\varepsilon}$, define the
  following operator $P$ inverse to the operator $S$:
$$
P u(x) = u(x)+ \int\limits_{x_1}^{\infty} N(x_1, t) u(t, x^1) \, dt.
$$
Then the assertions of Lemmas \ref{L18.1}--\ref{L18.3}
 hold for the kernel $N \lr{x_1,t}$ and, if  $S u \in T
\lr{\lambda+ \varepsilon},$ then
 \begin{equation}\label{l20}
    {P Su(x) = u(x).}
    \end{equation}
   Taking into account that $0 \in T \lr{\lambda+ \varepsilon}$ and applying \eqref{l20}
  to both sides of the relation  $S u = 0$ (we found that it holds in $\R_{+}^n$),
 we obtain that $u=0$ in $\R_{+}^n.$ It is proved above that this implies the relation
  $u \equiv 0$ in the whole space $\R^n.$

\begin{remark}
    The passage to the half-space
  $\R_{+}^n$ is used in the proof because expression \eqref{l6} is not defined in the
  intersection of the ball  $|x|\leq R_0$ and the infinite semicylinder
   $\{|x^1| \leq R_0,~|x_1| \leq R_0 \}.$
\end{remark}

Thus, the following assertion is proved.

\begin{theorem} \label{Theo1}
Each solution  $u(x) \in C^2 \lr{|x|>R_0}$ of the stationary
Schr\"odinger equation with a bounded potential
$$
\Delta u(x) - q(x_1) u = 0, ~x \in \R^n,~ |x| \geq R_0 > 0,
$$
$$
q(x_1) \in C \lr{|x| \geq R_0},~|q(x_1) | \leq \lambda^2,~\lambda>0,
$$
    satisfying the estimate
$$
|u(x)| \leq \const e^{- \lr{\lambda+\varepsilon}|x|},~\varepsilon>0,
$$
    is the identical zero.
\end{theorem}

{\bf 4.}
    The used technique of transmutation operators allows one to strengthen the obtained
    result. Let
$L_{2,\, loc} \lr{x_1 \geq R_0}$ denote the set of functions
$\psi$ such that  the integral $\int\limits_{R_0}^{x_1}
\psi^2(s)\,ds$ is finite for each  $x_1 \geq R_0$.
    Let $g(x)$ be a nonnegative function such that
$\int\limits_{x_1}^{\infty} t\, g(t, x^1) \, dt = p(x)$ is finite
for each  $x_1 \geq R_0$ and there exist a positive constant
$\alpha$ such that
$$
|p(x)| \leq c \cdot \exp \lr{- \alpha |x|^{\delta}},~\delta>0.
$$
    Then, arguing in the same way as in the proof of the previous theorem,
    one can prove the following assertion.

 \begin{theorem} \label{Theo2}
 Let $\psi(x_1) \in L_{2,\, loc} \lr{x_1 \geq R_0},$
  $\psi(x_1)$ be a nondecreasing function, and a function $g(x)$ satisfy the requirements
  listed above. Then each solution of the equation
$$
\Delta u(x) - q(x_1) u = 0, ~x \in \R^n,~ |x| \geq R_0 > 0,
$$
$$
|q(x_1)| \leq \psi^2 (x_1),
$$
   satisfying the inequality
$$
\psi (x_1) |u(x)| \leq \const e^{-\psi (x_1) |x|} g(x),~g(x) \geq 0,
$$
    is the identical zero.
\end{theorem}

Under the assumptions of Theorem \ref{Theo1}, one has to assign
$g(x) = e^{-\varepsilon |x|}.$ Another example of an admissible
function $g(x)$ is the function $g(x)= \exp \lr{-\varepsilon
|x|^{\delta} },$ $0< \delta <1.$
    This case is an example of problem \eqref{LM1}-\eqref{LM2}
    (the generalized Landis--Meshkov
     problem) as well.

According to the same scheme, one can consider the case of
potentials depending only on the radial variable.
 In this case, the answer to the original Landis problem is positive as well:
 after the passage to the spherical coordinates, one has to use transmutation operators
 for the Bessel perturbed operator, similar to operators considered in the previous section
 (see \cite{S3, S71, S75}).

The  generalized Landis
     problem can be considered for the case of more general differential equations
     and corresponding growth estimates for solutions. For example, it would be interesting
     to do this for the nonlinear
$p$-Laplace equation \cite{Lind, DKN} (this problem arose during
the discussion between the second author and Prof. Kon'kov at the
seminar of the Department of Differential Equations of MSU).

%\chapter*{Заключение}
%\addcontentsline{toc}{chapter}{Заключение}
%\chaptermark{\sc Заключение}

\chapter*{Acknowledgments}
\addcontentsline{toc}{chapter}{Acknowledgments}
\sectionmark{\sc Acknowledgments}

The second author is grateful to his colleagues reviewed the
manuscript and provided a number of useful remarks, corrections,
and additions: A.\,V.~Glushak, D.\,B.~Karp, V.\,V.~Kravchenko,
A.\,B.~Muravnik, and E.\,L.~Shishkina. The second author is
grateful to E.\,M.~Varfolomeev for his great labor to edit the
present book. Also, he expresses his gratitude to S.\,N.~Ushakov
    for his  qualified computer typesetting of the manuscript of the dissertation
    of V.\,V.~Katrakhov.

\newpage

\chapter*{Valery Vyacheslavovich Katrakhov\\ (1949--2010):\\a brief biographical background}
\addcontentsline{toc}{chapter}{Valery Vyacheslavovich Katrakhov: a
brief biographical background}
\sectionmark{\sc Valery Vyacheslavovich Katrakhov (1949--2010): a brief biographical background}

\epigraph{``It would be a     mission for his friends to write
down his biography; but remarkable people tracelessly disappear in
our society. We are lazy and incurious...''}{A.\,S.~Pushkin. ``A
Journey to Erzurum'' (1836).}

%\begin{figure}[H]
%\centering
%\begin{subfigure}[t]{0.45\textwidth}
%\centering
%\includegraphics[width=0.7\textwidth]{Kat1.eps}
%\subcaption{\hfil Валерий Вячеславович Катрахов\hfil\\
%\hphantom{\hspace{1.4cm}}(фото из личного дела 1972~г.)}
%\label{pic3}
%\end{subfigure}
%\qquad
%\begin{subfigure}[t]{0.45\textwidth}
%\centering
%\includegraphics[width=0.7\textwidth]{Kat2.eps}
%\subcaption{\hfil Валерий Вячеславович Катрахов\hfil\\
%\hphantom{\hspace{1.4cm}}(фото из личного дела 2006~г.)}
%\label{pic4}
%\end{subfigure}
%\end{figure}

Valery Vyacheslavovich Katrakhov was born on August 4,  1949 on
Sakhalin island. His father Vyacheslav Timofeevich Katrakhov was
an officer of the Soviet Army; afterwards, he worked at the
Voronezh Aviation Plant. His mother Klavdiya Vasil'evna Katrakhova
was  a  physician; she worked as a     polyclinic staff doctor.

\begin{figure}[H]
\centering
\begin{subfigure}[t]{0.45\textwidth}
\centering
\includegraphics[width=0.7\textwidth]{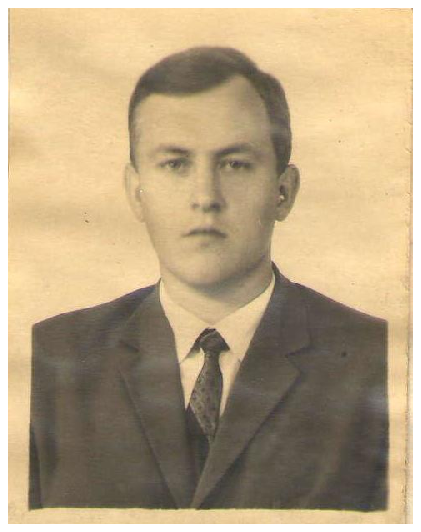}
\subcaption{\hfil Valery Vyacheslavovich Katrakhov\hfil\\
\hphantom{\hspace{1.4cm}}(foto 1972)}
\label{pic3}
\end{subfigure}
\qquad
\begin{subfigure}[t]{0.45\textwidth}
\centering
\includegraphics[width=0.7\textwidth]{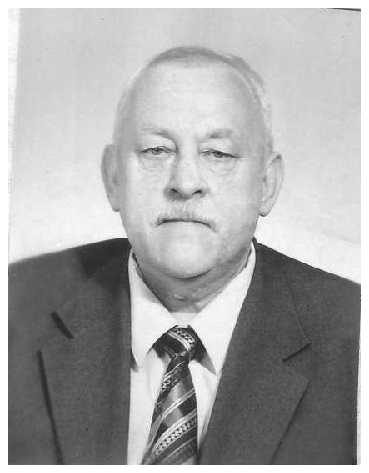}
\subcaption{\hfil Valery Vyacheslavovich Katrakhov\hfil\\
\hphantom{\hspace{1.4cm}}(foto 2006)}
\label{pic4}
\end{subfigure}
\end{figure}

In  1966, V. V. Katrakhov graduated from the mathematical class of
Voronezh school No 15. Since 1966 to 1971, he studied at the
Mathematical Faculty of Voronezh State University (VSU); he
graduated from it with honor. In 1971-1972,
 he was a postgraduate of Voronezh State University.

Further labor activity of  V. V. Katrakhov is strongly related to
the Applied Mathematics and Mechanics Faculty of Voronezh State
University.
 The faculty is founded in 1969 as a result of the separation
  of the Mathematics and Mechanics Faculty into the Mathematical Faculty
  and the Applied Mathematics and Mechanics Faculty (AMM).
 Mainly, the faculty is founded due to efforts professor Gennadiy
 Ivanovich Bykovtsev (a known scientist in mechanics): he became
 its first dean.

Since 1972,  V. V. Katrakhov worked at the Differential Equations
Department of Voronezh State University as an assistant-professor
(1972--1974), senior lecturer (1974--1980),
 and associate professor (1980--1983).
 Since 1983 to 1987, he was the head of the Computational Mathematics Department
 of the Applied Mathematics and Mechanics Faculty.
 In 1977, his work  ``Singular boundary-value
 problems and their applications to optimal control'' was awarded the Voronezh     Komsomol
 Prize in Science and Technics.
 During a number of years, he supervised the student research
 society of the Applied Mathematics and Mechanics Faculty, worked
 (pro bono)  as the deputy dean in scientific work, and was a member
 of the AMM faculty council and of the  VSU dissertation council.

Since 1987 to 2006,  V. V. Katrakhov worked at Vladivostok. In
1987--1995, he worked at the Institute of  Applied Mathematics of
Far-Eastern Division of the USSR Academy of Sciences (later,
Far-Eastern Division of the Russian Academy of Science) as a
senior scientific researcher, laboratory head, and deputy
director; in 1995--2006,
 he was a professor of the Pacific State University of Economics and a number
 of other Far-Eastern
 institutes of higher education.

In  2006, he came back to Voronezh and worked as a professor of
 the Differential Equations
Department as well as of Software and System Administration of the
Applied Mathematics and Mechanics Faculty of Voronezh State
University till 2008.

In 1974, he became a PhD in Mathematics, defending the
dissertation ``On the spectral function of singular differential
operators'' (supervised by professor Ivan Aleksandrovich
Kipriyanov) at Voronezh State University; in 1989, he became a D.
Sci. in Mathematics, defending the dissertation ``Singular
elliptic boundary-value problems. Method of transmutation
operators'' at Novosibirsk State University.

 V. V. Katrakhov achieved a number of significant results in various areas of mathematics
 such as theory of singular and degenerating differential
 equations, transmutation operators and integral transformations,
 special functions and fractional integrodifferentiating
 operators, function theory and functional analysis, spectral theory and pseudodifferential
 operators, numerical methods and optimal control, mathematical
 physics, and simulating of systems.
 To obtain the following results, V. V. Katrakhov used new ideas and methods:
 the resolving of systems of singular differential and pseudodifferential equations
 with weight boundary-value
 conditions,
    the formulation and resolving of new boundary-value
    problems for equations with essential singularities,
    the introduction of a new convolutional nonlocal boundary-value
 condition (the ``$K$-trace'')
 for solutions with singularities,
 the introducing of a new class of  denumerably normable
  Frechet spaces, based on transmutation operators,
  for problem settings and well-posedness
  proofs for elliptic equations with isolated essential singular
  points,
  the introducing and applying of new classes of transmutation operators
  and fractional integrodifferentiating,
  the discovering of new settings of boundary-value
    problems for singular differential equations in angular domains and in the
   Lobachevsky space, and the matrix method    to investigate the  Ising model.
   In his main works on the theory of differential equations, he used and developed the method
   of transmutation operators.

    Four  disciples of V. V. Katrakhov
    (N.\,I.~Golovko, A.\,B.~Muravnik, I.\,P.~Polovinkin, and S.\,M.~Sitnik) have  D. Sci.
degrees.
    About ten ones have PhD degrees.

   V. V. Katrakhov is the author of about  150 research works in leading Russian Journals
   such as  \emph{Doklady Mathematics}, \emph{Differential Equations},
    \emph{Sbornik: Mathematics},
\emph{Siberian Mathematical Journal}, and others. Also,  he  is
the author of the following seven monographs:
\begin{itemize}
\item V. V. Katrakhov and L. S. Mazelis. \emph{Continuity,
Completion, and Closure in Metric Spaces},
 Vladivostok: Far-Eastern
 Division of the Russian Academy of Science,   2000.

\item  N. I. Golovko and V. V. Katrakhov.
  \emph{Analysis of Queuing Systems Operating in Random Media},
Vladivostok:  Far-Eastern
 State Academy of Economics and Management,   2000.

\item A. A. Dmitriev, V. V. Katrakhov, and Y. N. Kharchenko.
 \emph{Root Transfer Matrices in Ising Models},
  Moscow: Nauka, 2004.

\item V. V. Katrakhov and D. E. Ryzhkov.
    \emph{Introduction in Functional-Analysis
    Methods in Dynamical Queuing Theory}, Vladivostok:  Far-Eastern
 State University, 2004.

\item V. V. Katrakhov,  N. I. Golovko,  and D. E. Ryzhkov.
   \emph{Introduction in Theory of Markov Twice Stochastic
  Queuing Systems}, Vladivostok:  Far-Eastern
 State University,  2005.

\item  N. I. Golovko and V. V. Katrakhov.
  \emph{Applications  of Queuing System Models  in Information Systems},
Vladivostok:  Far-Eastern
 State Univesity of Economics,   2008.

\item V. V. Katrakhov and S. M. Sitnik. Transmutation method and
boundary-value
 problems for singular elliptic equations,
\emph{Sovrem. Mat. Fundam. Napravl.}, Vol. 64, No 2, 211--426
(2018)

\end{itemize}
He was married. His wife Alla Anatol'evna  Katrakhova is an
instructor of mathematics.
 His daughter  Alla Valer'evna  Katrakhova is a physician.

V. V. Katrakhov died in 2010. He is buried at the Komintern
cemetery of Voronezh city.

\def\numberline{}

\end{document}